\documentclass[reqno,11pt]{amsart}
\usepackage{amsmath, amssymb, amsthm, amsfonts}
\usepackage{geometry}
\usepackage{bm}
\usepackage{fullpage}
\usepackage[hypertexnames=false]{hyperref}
\usepackage{mathtools}
\usepackage{mathrsfs}
\usepackage{enumitem}
\usepackage{upgreek}

\newcommand{\rW}{\mathring{W}}

\newcommand{\p}{\partial}

\newcommand{\R}{\mathbb{R}}
\newcommand{\bs}[1]{\boldsymbol{#1}}

\newcommand{\curl}{\nabla \times}
\newcommand{\supp}{{\operatorname{supp}}}

\newcommand{\Lag}{{\operatorname{Lag}}}

\newcommand{\twoD}{{\operatorname{2D}}}
\newcommand{\threeD}{{\operatorname{3D}}}

\newcommand{\sep}{\operatorname{sep}}
\newcommand{\loc}{\operatorname{loc}}

\newcommand{\BS}{\operatorname{BS}}

\newcommand{\cut}{{\operatorname{cut}}}

\newcommand{\reg}{{\operatorname{reg}}}
\newcommand{\ud}{\, \mathrm{d}}
\newcommand{\pv}{{\operatorname{p.v.}}}

\newcommand{\sing}{{\operatorname{sing}}}
\newcommand{\lin}{{\operatorname{lin}}}
\newcommand{\hyp}{{\operatorname{hyp}}}
\newcommand{\mm}{{\operatorname{m}}}
\newcommand{\md}{\operatorname{m}}

\newcommand{\transport}{{\operatorname{tport}}}
\newcommand{\pressure}{{\operatorname{prsr}}}
\newcommand{\core}{{\operatorname{core}}}
\newcommand{\node}{{\operatorname{node}}}

\newcommand{\tail}{{\operatorname{tail}}}
\newcommand{\far}{{\operatorname{far}}}

\newcommand{\pert}{{\operatorname{pert}}}
\newcommand{\inn}{{\operatorname{in}}}

\newcommand{\err}{\operatorname{err}}

\newcommand{\smooth}{\operatorname{sm}}
\newcommand{\cusp}{\operatorname{cusp}}

\newcommand{\mixed}{\operatorname{mix}}
\newcommand{\rWm}{ {\mathring{W}_{\!0}^{\!\md}} }
\newcommand{\Jm}{J_{\md}}
\newcommand{\tr}{{\operatorname{tr}}}

\theoremstyle{plain}
\newtheorem{theorem}{Theorem}[section]
\newtheorem{lemma}[theorem]{Lemma}
\newtheorem{proposition}[theorem]{Proposition}
\newtheorem{corollary}[theorem]{Corollary}

\theoremstyle{definition}
\newtheorem{definition}[theorem]{Definition}

\theoremstyle{remark}
\newtheorem{remark}[theorem]{Remark}

\numberwithin{equation}{section}

\newcommand{\runinheadstyle}[1]{{\small\scshape #1}}
\newcommand{\runinhead}[1]{\leavevmode\par\smallskip\noindent\runinheadstyle{#1}\enspace}

\newcommand{\rruninheadstyle}[1]{{\normalsize\itshape #1}}
\newcommand{\rruninhead}[1]{\leavevmode\par\noindent\rruninheadstyle{#1}\enspace}

\usepackage{etoolbox}

\makeatletter
\patchcmd{\@settitle}{\uppercasenonmath\@title}{}{}{}

\def\@settitle{
  \begin{center}
    \baselineskip14\p@\relax
    \large\bfseries
    \@title
  \end{center}
}
\makeatother

\makeatletter
\def\@setauthors{
  \begingroup
  \trivlist
  \centering\footnotesize \@topsep30\p@\relax
  \advance\@topsep by -\baselineskip
  \item\relax
  \andify\authors
  \def\baselinestretch{1}
  \authors
  \endtrivlist
  \endgroup
}
\makeatother

\title{Incompressible Euler Blowup Below the $C^{1,\frac{1}{3}}$ Threshold}
\author{\large Steve Shkoller}

\date{}

\begin{document}

\begin{abstract}
We prove finite-time Type--I blowup for the three-dimensional incompressible Euler equations in the axisymmetric no-swirl class, with initial 
velocity in $C^{1,\alpha}(\mathbb{R}^3)\cap L^2(\mathbb{R}^3)$, odd symmetry in $z$, and $0<\alpha<\tfrac13$, for an explicit class of finite-energy initial data.
The singularity forms at a stagnation point on the symmetry axis.  The axial strain and the global vorticity norm blow up at the Type--I 
rates $-\p_z u_z(0,0,t)\simeq (T^*-t)^{-1}$ and $\|\omega(\cdot,t)\|_{L^\infty}\simeq (T^*-t)^{-1}$, while the meridional Jacobian collapses according to
$J(t)\simeq (T^*-t)^{1/(1-3\alpha)}$. The proof is organized around a Lagrangian clock-and-driver framework.
The clock is the meridional Jacobian $J(t)$, and the driver is the compressive axial strain
$-\p_z u_z(0,0,t)$.  These variables satisfy, to leading order, a closed Riccati-clock system: the axial strain
drives the collapse of $J(t)$, while the collapse of $J(t)$ amplifies the axial strain.  We prove that the Euler
flow tracks this clock-and-driver model up to the singular time.  The main nonlocal obstruction is the pressure
Hessian; it is controlled by a non-perturbative strain--pressure Hessian comparison showing that pressure cannot
cancel the quadratic compressive strain responsible for collapse.  This gives a dynamical explanation of the threshold $\alpha=\tfrac13$.  
The blowup mechanism is structurally stable and persists for an open set of admissible angular functions in a weighted H\"older topology.
\end{abstract}

\maketitle

\setcounter{tocdepth}{1}
\tableofcontents

\nopagebreak[4]
\section{Introduction}

The global regularity problem for the three-dimensional incompressible Euler equations has been one of the
central open questions in fluid mechanics. The system is
\begin{subequations}\label{eq:euler}
\begin{alignat*}{2}
\p_t u + (u\cdot \nabla) u + \nabla p &= 0  &&\ \ \ \text{in} \ \ \mathbb{R}^3 \times (0,T] \,,\\
\nabla \cdot u &= 0 &&\ \ \ \text{in} \ \ \mathbb{R}^3 \times (0,T] \,, \\
u(\cdot, 0) & = u_0 \,. &&
\end{alignat*}
\end{subequations}
A singularity can only form through amplification of the vorticity $\bs{\omega}=\curl u$, which evolves by
\begin{equation}\label{eq:vorticity_eqn}
\p_t \bs{\omega} + (u \cdot \nabla)\bs{\omega} = (\bs{\omega}\cdot\nabla)u\,.
\end{equation}
The stretching term $(\bs{\omega}\cdot\nabla)u$ couples to $\bs{\omega}$ through the nonlocal Biot--Savart
law. In particular, the Beale--Kato--Majda (BKM) criterion \cite{BKM1984} implies that a finite-time
singularity at $T^*<\infty$ requires $\int_0^{T^*}\|\bs{\omega}(\cdot,t)\|_{L^\infty}\,\ud t=\infty$.

\subsection{Axisymmetric no-swirl and rough velocities}
Throughout this paper we work in the axisymmetric no-swirl class, where $u=u_r e_r+u_z e_z$ and
$\bs{\omega}=\omega_\theta e_\theta$. For smooth data, this setting is globally regular (see e.g.\
\cite{UkhovskiiYudovich1968,Majda1986,ShirotaYanagisawa1994}); the mechanism in
\eqref{eq:vorticity_eqn} becomes dangerous only in the \emph{rough} regime where the transported
quantity $\tfrac{\omega_\theta}{r}$ is singular near the axis. Even in the smooth globally regular regime,
one can still observe quantitative long-time growth and outward migration phenomena; see e.g.\
\cite{EgamberganovYao2025} and the references therein.

\subsection{The $C^{1,\alpha}$ threshold in axisymmetry}
A major breakthrough due to Elgindi~\cite{Elgindi2021} established finite-time singularity formation for the
axisymmetric no-swirl Euler equations in $\R^3$ with velocity in $C^{1,\alpha}$. A key feature is that this
theory is \emph{intrinsically perturbative in the H\"older exponent}: the construction requires
$ 0<\alpha<\alpha_0$
for a sufficiently small constant $\alpha_0>0$. In this small-$\alpha$ regime one can exploit a perturbative
structure in which the genuinely three-dimensional nonlocal effects (pressure and strain coupling) remain weak
enough to be closed by bootstrap estimates.\footnote{A useful heuristic is that, in the limit $\alpha\to 0$,
the axisymmetric Biot--Savart law admits an effective ``low-rank'' (nearly one-dimensional) structure for the
dominant coupling into the axial strain, with a remainder that is perturbative in $\alpha$; see e.g.\
\cite{TaoBlog2019} for discussion. This explains why the small-$\alpha$ approach is powerful, but also why it
is not expected to extend to $\alpha$ comparable to $\tfrac13$ without a different organizing principle.}
Subsequent developments in the same rough self-similar setting, including stability and geometric
refinements, likewise operate in the \emph{small-$\alpha$} regime; see e.g.\
\cite{ElgindiGhoulMasmoudi2022, Cordoba2023, ElgindiPasqualotto2023}.

As $\alpha$ increases, however, the mechanism ceases to be perturbative: the nonlocal pressure and strain
contributions become comparable to the Riccati compression term, and estimates that are ``$\alpha$--small'' in the
perturbative theory no longer gain a closing power.

On the other hand, rigidity results of Saint-Raymond~\cite{SaintRaymond1994} for $\alpha>\tfrac{1}{3}$
(see also Danchin~\cite{Danchin2007}) and Shao--Wei--Zhang~\cite{ShaoWeiZhang2025} for $\alpha=\tfrac{1}{3}$
show that axisymmetric no-swirl solutions with finite energy are globally regular when $u\in
C^{1,\alpha}$ with $\alpha\ge\tfrac13$.

Thus $\alpha=\tfrac13$ is the sharp dynamical threshold in the finite-energy axisymmetric no-swirl
$C^{1,\alpha}$ class: for $\alpha\ge\tfrac13$ all such solutions are globally regular, while the Lagrangian
clock mechanism developed here collapses precisely in the subcritical range $\alpha\in(0,\tfrac13)$.  The
theorems below prove finite-time Type--I blowup throughout this subcritical range, for finite-energy target
tails $\gamma>\alpha+\tfrac52$.

\subsection{A brief landscape of rigorous blowup results}

Known rigorous finite-time singularity mechanisms for \eqref{eq:euler} fall into two broad families: those that
(i) permit \emph{rough} velocities (typically $u\in C^{1,\alpha}$) and those that (ii) exploit a
\emph{boundary} geometry.

In the whole-space axisymmetric no-swirl class, Elgindi~\cite{Elgindi2021} constructed Type--I self-similar
blowup\footnote{We call a finite-time singularity at $t=T^*$ \emph{Type--I} if the vorticity obeys the
scale-invariant bound $\sup_{0<t<T^*}(T^*-t)\|\omega(t)\|_{L^\infty}<\infty$ (that is,
$\|\omega(t)\|_{L^\infty}\lesssim (T^*-t)^{-1}$), which matches the Euler scaling $u(x,t)\mapsto u(\lambda
x,\lambda t)$ (so $\omega\mapsto \lambda\omega$). A blowup is
\emph{Type~II} if $(T^*-t)\|\omega(t)\|_{L^\infty}\to\infty$; while such faster vorticity growth is
not ruled out \emph{a priori} by the equations or the Beale--Kato--Majda criterion, the singularities proved
here are Type--I.} in the $C^{1,\alpha}$ regime, but crucially in the
\emph{perturbative} range $0<\alpha\ll 1$. A number of subsequent results establish stability
properties and related refinements of this rough blowup regime, again in a small-$\alpha$ framework; see e.g.\
\cite{ElgindiGhoulMasmoudi2022, ElgindiPasqualotto2023, Chen2023} and references therein. More recently,
C\'ordoba--Mart\'inez-Zoroa--Zheng~\cite{Cordoba2023} introduced a multi-region blowup architecture for
solutions with $u\in C^{1,\alpha}\cap L^2$ that are smooth away from a point; this construction likewise
operates in a small-$\alpha$ regime.

In the presence of a physical boundary, the numerical scenario of Luo--Hou~\cite{LuoHou2014} led to a
different blowup mechanism.  Chen--Hou proved finite-time blowup for axisymmetric Euler (and related models)
with $C^{1,\alpha}$ velocity~\cite{ChenHouBoundary2021}, and more recently developed a separate framework
yielding blowup for smooth boundary data (see
\cite{ChenHouPartI,ChenHouPartII}). For broader context and further references, see the survey of
Drivas--Elgindi~\cite{DrivasElgindi2023}. For related constraints on \emph{globally self-similar}
finite-energy blowup scenarios for 3D Euler, including in axisymmetry, see
Constantin--Ignatova--Vicol~\cite{ConstantinIgnatovaVicol2026}.

\subsection{The $C^{1,\frac13}$ regularity gap}
While several distinct mechanisms are now known to produce Euler blowup in rough regimes, the available
$C^{1,\alpha}$ constructions in $\R^3$ are perturbative in $\alpha$ and do not reach the sharp finite-energy
threshold $\alpha=\tfrac13$. The question left open is whether blowup can occur throughout the full
subcritical range $\alpha\in(0,\tfrac13)$ in the finite-energy class, and whether such blowup is structurally
stable rather than confined to a specially tuned datum.

The proof below closes this gap by separating the argument into two parts.  First, the clock-and-driver
hyperbolic system of Section~\ref{sec:lag-analysis-non-local} identifies the collapse mechanism.  In that system
there is no pressure Hessian: the model clock $\Jm$ evolves by the kinematic identity
\[
\dot\Jm(t)=\tfrac12\,\Jm(t)\rWm(t), \qquad \rWm(t)\simeq-\Gamma\Jm(t)^{3\alpha-1},
\]
and therefore
\[
\dot\Jm(t)\simeq-\Gamma\Jm(t)^{3\alpha}.
\]
The model clock reaches zero in finite time precisely when $3\alpha<1$.  Second, the Euler part of the proof
shows that the true solution follows this clock law in collapsing Lagrangian coordinates.  The non-perturbative
issue is that the true Euler axial strain is not just a Biot--Savart quantity.  It also satisfies the
stagnation-point Riccati law, so the pressure Hessian could in principle cancel the compressive driving strain in
that law.  The main pressure Hessian estimate proves that this cancellation does not occur: after the model
pressure Hessian bound is transferred to the exact Euler pressure Hessian, the bound \eqref{eq:euler-riccati-bound}
shows that
\[
\Pi_0(t)\ge -\upbeta\,\tfrac12\,\rW_0(t)^2 \qquad\text{for some }0<\upbeta<1.
\]
This is the estimate that allows the clock-and-driver model dynamics to persist for the true Euler dynamics.

\subsection{Key ideas}

Before stating the main results, we highlight four ideas that distinguish the present argument from the
small-$\alpha$ approach.

\runinhead{(i) The linear hyperbolic guide.} Prior $C^{1,\alpha}$ blowup constructions in $\R^3$ are organized
around \emph{Eulerian} self-similar models, either one-dimensional analogues (Constantin--Lax--Majda,
De~Gregorio, and related models; see, for example, Elgindi--Jeong~\cite{ElgindiJeong2020} on the
competition between advection and vortex stretching) or explicit asymptotic profiles in Eulerian variables.
Those descriptions are inherently perturbative in $\alpha$.  We instead begin with the Lagrangian
clock-and-driver model for $(\Jm,\rWm)$ in Section~\ref{sec:lag-analysis-non-local}.  In this model, the
true Euler cusp flow $\phi_{\cusp}$ is replaced by the hyperbolic flow map $\Phi_\lin$, defined by
\[
R\mapsto\Jm(t)^{-1}R, \qquad Z\mapsto\Jm(t)^2Z.
\]
Our clock-and-driver model has no pressure Hessian and no Riccati constraint.  It uses only the kinematic clock identity and
the Biot--Savart scaling of the model axial strain.  The drift law in the model says that particles producing
the axial strain at small clock values come from ever-smaller Lagrangian polar angles
\[
\sigma_{\Lag}\lesssim \Jm(t)^3.
\]
Since the Target Profile has the H\"older cusp
$\Theta^*(\sigma)\simeq(\sin\sigma)^\alpha$, the strain-producing sector is depleted by the power
$\Jm(t)^{3\alpha}$.  As a consequence, we find that
\[
\rWm(t)\simeq-\Gamma\Jm(t)^{3\alpha-1}, \qquad \dot\Jm(t)\simeq-\Gamma\Jm(t)^{3\alpha}.
\]
Thus, the threshold $\alpha=\tfrac13$ is the exponent at which drift-induced depletion balances the
stretching-driven collapse.\footnote{The same critical exponent appears in the formal model discussed by
Drivas--Elgindi~\cite[Section~4.4, equations~(131)--(132), Lemma~4.33]{DrivasElgindi2023}, where
$\alpha=\tfrac13$ is identified as the threshold at which angular transport balances growth.  The mechanism
used here is different in form: it is a Lagrangian clock description for the Euler flow, in which the drift law
gives $\sigma_{\Lag}\simeq J^3$, the cusp profile contributes $J^{3\alpha}$, and the clock law gives
$\dot J\simeq-\Gamma J^{3\alpha}$.  Thus the present argument should be viewed as a closed Euler realization of
this drift-versus-vanishing balance, not as an identification of the formal model equation with the clock ODE.}

\runinhead{(ii) The Euler flow has the linear hyperbolic leading normal form.}
Sections~\ref{sec:Euler-blowup-for-Theta-star}--\ref{sec:current-axis-normal-forms} prove that the singular part of the
true Euler flow approaches the linear hyperbolic model in the collapse limit.  We decompose the exact flow map as
\[
\phi=\phi_{\smooth}\circ\phi_{\cusp},
\]
where $\phi_{\smooth}$ is generated by the smooth far-field velocity and $\phi_{\cusp}$ carries the singular
core.  After composing away $\phi_{\smooth}$ and dividing the Eulerian image by the axial scale
$J_{\cusp}^2$, the true Euler cusp-flow map $\phi_{\cusp}$, expressed in collapse coordinates, has the normal form
\[
R\mapsto J_{\cusp}^{-1}R, \qquad Z\mapsto J_{\cusp}^{2}Z,
\]
up to errors that vanish in the collapse limit.  This is the geometric reason that the true clock
$J(t)$ follows the model collapse law
in the small-clock regime.

There is one additional normalization on the symmetry axis.  The cusp flow $\phi_{\cusp}$ creates an axial drift
that is linear at leading order, and this drift is described by a one-dimensional axial map
\[
\zeta=\mathscr Z_t(\eta).
\]
Our proof works at fixed reference label $\eta=\mathscr Z_t^{-1}(\zeta)$ rather than at fixed instantaneous
coordinate $\zeta$.
This removes the logarithmic accumulation that would otherwise appear when differentiating the cusp map along the axis.

\runinhead{(iii) The true axial strain must satisfy the Euler Riccati law.} The clock-and-driver model axial strain $\rWm$ is
reconstructed from the model vorticity by the Biot--Savart law.  The true Euler axial strain $\rW_0$ has the model
Biot--Savart scaling only after the pressure Hessian in the stagnation-point Riccati law is controlled.  With
\[
\rW_0(t):=\p_z u_z(0,0,t),
\]
the exact equation is
\begin{equation}\label{eq:intro-Riccati}
\p_t\rW_0(t)=-\tfrac12\,\rW_0(t)^2-\Pi_0(t),
\end{equation}
where $\Pi_0(t)$ is the nonlocal principal-value part of the stagnation-point pressure Hessian; see
\eqref{eq:pi-Eul}.  The collapse clock is enslaved to the
same axial strain by
\[
\dot J(t)=\tfrac12\,J(t)\,\rW_0(t).
\]
Therefore, our proof does not simply assert the scaling
$-\rW_0(t)\simeq\Gamma J(t)^{3\alpha-1}$; rather, we first prove that the pressure Hessian leaves a definite
portion of the Riccati compression:
\[
\Pi_0(t)\ge -\upbeta\,\tfrac12\,\rW_0(t)^2, \qquad 0<\upbeta<1.
\]
Then, by \eqref{eq:intro-Riccati}, we see that
\[
\p_t\rW_0(t)\le-\tfrac{1-\upbeta}{2}\,\rW_0(t)^2,
\]
and this Riccati inequality is then used to prove the axial strain scaling and the finite-time clock collapse.

\runinhead{(iv) The pressure Hessian is proved through a renormalized axis-trace reduction.} To estimate the pressure
Hessian, our proof strategy first approximates the true transported Euler vorticity by the separation-of-variables form
\[
\Omega_\theta^a(R,Z)= -\operatorname{sgn}(Z)\,a(|Z|)R^\alpha .
\]
An auxiliary first-variation computation reduces the relevant pressure Hessian estimate to one-dimensional
integrals involving the axial function $a$.  In the Euler solution, this axial function is the Euler-generated
profile $a_t$.  The axial stretching estimates prove that $a_t$ is nonnegative and nonincreasing in the
normalized axial coordinate, and the renormalized axis-trace criterion in
Proposition~\ref{prop:euler-generated-profile-riccati} gives
\[
\Pi_{\cusp}(t)\ge -q_\alpha\,\tfrac12\,\mathcal W_{\cusp}(t)^2, \qquad q_\alpha<1.
\]
The rest of the pressure Hessian argument compares the cusp-coordinate pressure contribution with the true Euler
pressure Hessian.  We prove that, for the Target Profile initial vorticity, the true Euler pressure Hessian
$\Pi_0(t)$ satisfies
\[
\Pi_0(t)\ge -\upbeta\,\tfrac12\,\rW_0(t)^2 \qquad\text{for some }0<\upbeta<1, 
\]
so that $\Pi_0$ cannot cancel the compression needed for clock collapse.
The remaining pressure estimates show that this model inequality persists for the exact Euler pressure Hessian.

\subsection{Main results}
\begin{theorem}[Finite-time Type--I blowup for the finite-energy Target Profile]
\label{thm:target-profile}
We fix $\alpha\in(0,\tfrac{1}{3})$ and $\gamma>\alpha+\tfrac{5}{2}$. For every
$\Gamma>0$, let
\[
u_0^*\in C^{1,\alpha}(\R^3)\cap L^2(\R^3)
\]
be the axisymmetric no-swirl initial datum, odd in $z$, whose initial toroidal vorticity component is given by
\eqref{eq:vort0} with perturbation $h\equiv 0$; that is, the datum generated by the Target Profile
$\Theta^*$ in Definition~\ref{def:init-data}. Then the corresponding
unique local Euler solution develops a finite-time singularity at some $T^*<\infty$ at the stagnation point
$(r,z)=(0,0)$.

Moreover, the singularity is Type--I: the axial strain blows up at the stagnation point, and the vorticity
norm has the same Type--I rate.  There exist
constants $0<c<C$ (depending only on $\alpha,\gamma$) such that as $t\uparrow T^*$,
\[
\tfrac{c}{T^*-t}\le \|\bs{\omega}(\cdot,t)\|_{L^\infty(\R^3)}\le \tfrac{C}{T^*-t},
\qquad \tfrac{c}{T^*-t}\le -\p_z u_z(0,0,t)\le \tfrac{C}{T^*-t}.
\]
If
\[
J(t):=\det \nabla_{(R,Z)}\big(\phi_r,\phi_z\big)(0,0,t)
\]
denotes the $2\times2$ Jacobian determinant of the meridional Lagrangian flow map at the stagnation point,
then
\[
c\,\big(\Gamma(T^*-t)\big)^{\frac{1}{1-3\alpha}} \le J(t)\le C\,\big(\Gamma(T^*-t)\big)^{\frac{1}{1-3\alpha}}.
\]
In particular,
\[
\int_0^{T^*}\|\bs{\omega}(\cdot,t)\|_{L^\infty}\,\ud t=\infty.
\]
\end{theorem}

Known global regularity results for finite-energy axisymmetric no-swirl solutions cover the regime
$\alpha\ge\tfrac13$.  Theorem~\ref{thm:target-profile} gives the complementary singular behavior below this
exponent: for every $0<\alpha<\tfrac13$ in the stated finite-energy class, there is an
explicit $C^{1,\alpha}$ datum whose Lagrangian clock reaches zero in finite time and whose solution develops
Type--I blowup.  Thus Theorem~\ref{thm:target-profile} closes the finite-energy $C^{1,\alpha}$ regularity gap for axisymmetric
no-swirl Euler.  The exponent $\alpha=\tfrac13$ is the sharp threshold: above it the drift-induced depletion
prevents the clock collapse, while below it the constructed data realize collapse and the pressure analysis
shows that the nonlocal pressure Hessian cannot cancel the Riccati compression.

\medskip

The second result shows that the mechanism is not a single-datum artifact: it persists for an open set, in the
weighted H\"older topology, of angular perturbations of the Target Profile.

\begin{theorem}[Open-set stability of the Target Profile blowup]
\label{thm:main}
We fix $\alpha\in(0,\tfrac{1}{3})$, $\eta>0$, and
$\gamma>\alpha+\tfrac{5}{2}$.  Then there exists a
constant
\[
\nu_0=\nu_0(\alpha,\gamma,\eta)>0,
\]
such that the following holds.

For any $\Gamma>0$ and any axisymmetric no-swirl initial datum, odd in $z$, whose initial toroidal vorticity
component belongs to the admissible class $\mathcal{A}_{\alpha,\gamma}(\nu,\eta)$
(Definition~\ref{def:init-data}) with amplitude $\Gamma$ and $0<\nu\le \nu_0$, the datum satisfies
\[
u_0\in C^{1,\alpha}(\R^3)\cap L^2(\R^3),
\]
and the corresponding unique local Euler solution develops a finite-time singularity at some $T^*<\infty$ at
the stagnation point $(r,z)=(0,0)$.

Moreover, the singularity is Type--I: the axial strain blows up at the stagnation point, and the vorticity
norm has the same Type--I rate.  There exist
constants $0<c<C$ (depending only on $\alpha,\gamma,\eta$) such that as $t\uparrow T^*$,
\[
\tfrac{c}{T^*-t}\le \|\bs{\omega}(\cdot,t)\|_{L^\infty(\R^3)}\le \tfrac{C}{T^*-t},
\qquad \tfrac{c}{T^*-t}\le -\p_z u_z(0,0,t)\le \tfrac{C}{T^*-t}.
\]
If
\[
J(t):=\det \nabla_{(R,Z)}\big(\phi_r,\phi_z\big)(0,0,t)
\]
denotes the $2\times2$ Jacobian determinant of the meridional Lagrangian flow map at the stagnation point,
then
\[
c\,\big(\Gamma(T^*-t)\big)^{\frac{1}{1-3\alpha}} \le J(t)\le C\,\big(\Gamma(T^*-t)\big)^{\frac{1}{1-3\alpha}}.
\]
In particular,
\[
\int_0^{T^*}\|\bs{\omega}(\cdot,t)\|_{L^\infty}\,\ud t=\infty.
\]
\end{theorem}

Theorem~\ref{thm:target-profile} is the conceptual core of the paper. Theorem~\ref{thm:main} strengthens it by
showing that this mechanism persists for a weighted H\"older neighborhood of the target angular profile. The
proof is organized accordingly: Sections~\ref{sec:lag-analysis-non-local}--\ref{sec:target-profile-typeI-completion}
establish Theorem~\ref{thm:target-profile}, and Section~\ref{sec:proof-main} upgrades it to
Theorem~\ref{thm:main}.  A detailed section-by-section guide is given in Section~\ref{sec:outline}.

\subsection{Kinematic drift and the Lagrangian clock}
\label{subsec:kinematic-drift-lagrangian-clock}
The Lagrangian clock is the natural organization for the open class of cusp data considered here.  In the
small-clock regime, the Eulerian sector that drives the axial strain is not supplied by one fixed material
sector.  Instead, the hyperbolic collapse forces a rapid \emph{kinematic drift} in polar angle: particles are
swept away from the axis toward the equator, and the strain-producing sector must be continually replenished by
particles that start closer and closer to the axis.

Quantitatively, the drift law (see Lemma~\ref{lem:drift}) shows that the Lagrangian polar angle $\sigma_{\Lag}$ of
a particle occupying an Eulerian point $x$ at clock value $J(t)$ satisfies
\[
\tan\sigma_{\Lag} = J(t)^3 \tan\sigma(x).
\]
(Here $\sigma(x)\in[0,\pi]$ denotes the usual polar angle from the positive symmetry axis, so $\sigma=0$ is the
positive $z$-axis and $\sigma=\tfrac{\pi}{2}$ is the equator.  Later, after using the odd symmetry across
$z=0$, we work on the folded upper-half angle $\sigma\in[0,\pi/2]$.) Thus, for any fixed Eulerian sector bounded
away from the equator, including the strain-producing sector described in Section~\ref{sec:initialdata}, the
relevant labels satisfy $\sigma_{\Lag}\lesssim J(t)^3$.  Since the target angular profile vanishes at the axis
with a H\"older cusp, $\Theta^*(\sigma)\simeq (\sin\sigma)^\alpha$ as $\sigma\downarrow 0$,
and the admissible profiles preserve this leading cusp with higher-order perturbative vanishing, the vorticity
fed into the strain-producing sector is depleted by the geometric penalty
\[
(\sin\sigma_{\Lag})^\alpha \simeq J(t)^{3\alpha}.
\]
In other words, as $J(t)\downarrow 0$, the strain-producing sector is populated by particles drawn from
progressively smaller Lagrangian angles, and hence from progressively smaller values of the cusp.  For the open
class of cusp data considered here, this continual replacement of material labels is the feature recorded by the
Lagrangian clock.  It is also the mechanism that produces the $\alpha=\tfrac13$ threshold: the effective
clock law has the scale
\[
        -\dot J(t)\simeq \Gamma J(t)^{3\alpha},
\]
so the time needed to reach $J=0$ is comparable to $\Gamma^{-1}\int_0 J^{-3\alpha}\,dJ$, which is finite exactly
when $3\alpha<1$.

\section{The Axisymmetric No-Swirl Setting}
To understand the singularity formation, we work in cylindrical coordinates $(r,\theta,z)$ in which the
velocity can be written as
\[
u(r,\theta,z) = u_r(r,z)e_r + u_\theta(r,z)e_\theta + u_z(r,z)e_z .
\]
A flow is \emph{axisymmetric} if all scalar components are independent of the azimuthal angle $\theta$, and
\emph{without swirl} if $u_\theta=0$.

\subsection{Cylindrical basis vectors}
We define the local orthonormal basis as
\begin{equation*}
e_r = (\cos \theta, \sin \theta, 0), \quad e_\theta = (-\sin \theta, \cos \theta, 0), \quad e_z = (0, 0, 1) \,.
\end{equation*}
Their $\theta$--derivatives are
\begin{equation*}
\p_\theta e_r = e_\theta \,, \qquad \p_\theta e_\theta = -e_r \,, \qquad \p_r e_i = \p_z e_i = 0 \,.
\end{equation*}

\subsection{The velocity gradient tensor}
For an axisymmetric, no-swirl velocity field $u = u_r(r, z) e_r + u_z(r, z) e_z$, one computes in the
$(e_r,e_\theta,e_z)$ basis that
\begin{equation}
\nabla u =
\begin{pmatrix}
\p_r u_r & 0 & \p_z u_r \\
0 & \tfrac{u_r}{r} & 0 \\
\p_r u_z & 0 & \p_z u_z
\end{pmatrix}.
 \label{eq:gradu-matrix}
\end{equation}
The term $\tfrac{u_r}{r}$ comes from the $\tfrac1r\,\p_\theta$ component in the cylindrical gradient.

\subsection{The pressure equation source term $S$}
Taking divergence of the momentum equation yields the Poisson equation
\begin{equation*}
-\Delta p = \sum_{i,j=1}^3 \p_i u_j\,\p_j u_i = \operatorname{tr}\!\big((\nabla u)^2\big) = \p_k u^i \p_i u^k := S \,.
\end{equation*}
Using \eqref{eq:gradu-matrix}, one finds
\begin{equation*}
(\nabla u)^2 =
\begin{pmatrix}
(\p_r u_r)^2 + (\p_z u_r)(\p_r u_z) & 0 & (\p_r u_r)(\p_z u_r) + (\p_z u_r)(\p_z u_z) \\
 0 & \left(\tfrac{u_r}{r}\right)^2 & 0 \\
(\p_r u_z)(\p_r u_r) + (\p_z u_z)(\p_r u_z) & 0 & (\p_r u_z)(\p_z u_r) + (\p_z u_z)^2
\end{pmatrix},
\end{equation*}
and hence
\begin{equation}
S = (\p_r u_r)^2 + \left(\tfrac{u_r}{r}\right)^2 + (\p_z u_z)^2 + 2 (\p_r u_z)(\p_z u_r) \,.
\label{eq:Source}
\end{equation}

\section{Lagrangian Variables and Vorticity Transport}

\subsection{The flow map}
Let $(R,Z)$ denote the initial cylindrical coordinates of a particle in the poloidal $(r,z)$--plane at time
$t=0$, and let $\Theta\in[0,2\pi)$ denote its initial azimuthal angle. We define the full 3D flow
map\footnote{Here ``full 3D'' means that $\phi_{\threeD}$ measures the Cartesian position of the particle,
including its initial azimuthal angle $\Theta$.  The axisymmetric no-swirl reduction then implies
$\Theta$ is conserved, so the dynamics are determined by the poloidal map $\phi=(\phi_r,\phi_z)$.}
$\phi_{\threeD}(R,\Theta,Z,t)$ by
\[
\p_t \phi_{\threeD} = u \circ \phi_{\threeD}, \qquad \phi_{\threeD}(R,\Theta,Z,0) = (R\cos\Theta, R\sin\Theta, Z).
\]
Since the flow is axisymmetric without swirl, the azimuthal angle is conserved: $\theta(t)\equiv \Theta$ along
each trajectory. Accordingly, there is a reduced \emph{poloidal} flow map
\[
\phi(R,Z,t) := (\phi_r(R,Z,t),\phi_z(R,Z,t))
\]
such that
\[
\phi_{\threeD}(R,\Theta,Z,t)=\big(\phi_r(R,Z,t)\cos \Theta, \phi_r(R,Z,t)\sin \Theta, \phi_z(R,Z,t)\big).
\]
The reduced map solves the ODE system
\begin{subequations}
\begin{alignat*}{2}
\p_t \phi_r(R, Z, t) & = u_r(\phi_r(R, Z, t), \phi_z(R, Z, t), t)\,, &&\qquad t \in (0,T] \,, \\
\p_t \phi_z(R, Z, t) & = u_z(\phi_r(R, Z, t), \phi_z(R, Z, t), t)\,, &&\qquad t \in (0,T] \,, \\
(\phi_r,\phi_z)(R,Z,0) &= (R,Z)\,.
\end{alignat*}
\end{subequations}

\subsection{The meridional Jacobian and the geometric identity}
Let
\[
J_{\twoD}(R,Z,t):=\det \nabla_{(R,Z)}(\phi_r,\phi_z)(R,Z,t)
\]
denote the Jacobian determinant of the reduced map in the poloidal plane. In cylindrical coordinates, the 3D
volume element is $r\,dr\,d\theta\,dz$. Under the map $(R,Z)\mapsto(\phi_r,\phi_z)$, we have $dr\,dz =
J_{\twoD}\,dR\,dZ$ and $r=\phi_r$, while $\theta=\Theta$ is unchanged. Incompressibility (volume preservation)
therefore gives
\[
\phi_r(R,Z,t)\,J_{\twoD}(R,Z,t)\,dR\,dZ\,d\Theta = R\,dR\,dZ\,d\Theta,
\]
hence the fundamental identity
\begin{equation}
J_{\twoD}(R, Z, t) = \tfrac{R}{\phi_r(R, Z, t)} \,.
\label{eq:Jac-Identity}
\end{equation}

\subsection{Vorticity}
The \textit{angular} vorticity, also known as the  \textit{toroidal} vorticity,  is the $\theta$-component
\[
\omega_\theta = \p_z u_r - \p_r u_z\,.
\]
In the axisymmetric no-swirl class it satisfies
\[
D_t\omega_\theta = \tfrac{u_r}{r}\,\omega_\theta,
\]
where $D_t:=\p_t+u_r\p_r+u_z\p_z$ is the material derivative.
Defining the \textit{specific vorticity} $\xi := \tfrac{\omega_\theta}{r}$, 
we therefore have conservation along trajectories:
\begin{equation*}
\xi(\phi_r(R,Z,t), \phi_z(R,Z,t), t) = \xi_0(R, Z), \qquad\xi_0(R,Z):=\xi(R,Z,0)=\tfrac{\omega_{\theta,0}(R,Z)}{R}.
\end{equation*}
This means that
\[
\tfrac{\omega_\theta(\phi_r(R,Z,t), \phi_z(R,Z,t), t)}{\phi_r(R,Z,t)} = \tfrac{\omega_{\theta,0}(R, Z)}{R}\,.
\]
Using \eqref{eq:Jac-Identity}, this yields the basic push-forward identity
\begin{equation}
\omega_\theta(\phi_r(R,Z,t),\phi_z(R,Z,t),t) = J_\twoD(R,Z,t)^{-1}\,\omega_{\theta,0}(R,Z).
\label{eq:vort-identity}
\end{equation}

\section{Reduction to the Symmetry Axis} 
On the symmetry axis $R=0$, the flow identities above reduce the Euler dynamics to ODEs for the axis
Jacobian $J(Z,t):=J_{\twoD}(0,Z,t)$ and the transported axial strain
$\rW(Z,t):=\p_z u_z(0,\phi_z(0,Z,t),t)$.  Their values at $Z=0$ give the collapse clock and driving
strain used throughout the proof.  We write $J$ for $J_\twoD$ in this section.

\subsection{Limits on the Axis}
The following identities are the axis limits of the reduced flow map as $R\to0$.  They use the usual parity
of axisymmetric no-swirl fields and the $C^1$ regularity of the flow map.
\begin{lemma}[Axial Flow Map Properties]\label{lem:flow1}
On the axis of symmetry $R=0$, the flow map $\phi = (\phi_r, \phi_z)$ satisfies:
\begin{enumerate}[label=(\alph*)]
\item $\phi_r(0, Z, t) = 0$,
\item $\p_R \phi_r(0, Z, t) = J(0, Z, t)^{-1}$,
\item $\p_Z \phi_z(0, Z, t) = J(0, Z, t)^2$, and
\item $\phi_z(0,0,t)=0$ (since the data are odd in $z$).
\end{enumerate}
\end{lemma}

\begin{proof}[Proof of Lemma \ref{lem:flow1}]
The identities follow from the radial parity of the velocity field in axisymmetry without swirl:
\begin{enumerate}[label=(\alph*)]
\item The radial velocity $u_r(r, z, t)$ is an odd function of $r$. Thus $u_r(0, z, t) = 0$ for all
$(z, t)$. The radial flow map $\phi_r$ satisfies the ODE $\p_t \phi_r(R, Z, t) = u_r(\phi_r, \phi_z, t)$
with initial condition $\phi_r(R, Z, 0) = R$. For $R=0$, the unique solution to $\dot{\phi}_r = u_r( \phi_r,
\phi_z, t)$ with $\phi_r(0) = 0$ is the trivial solution $\phi_r(0, Z, t) \equiv 0$. Physically, this implies
the axis is a material invariant.

\item Since $u_r$ is odd in $r$, the map $\phi_r(R, Z, t)$ is odd in $R$. Since $\phi_r$ is $C^1$ in
$R$, we have
\[
\phi_r(R, Z, t) = \p_R \phi_r(0, Z, t)\,R + o(R) \ \ \text{ for } \ \ R\to 0.
\]
Substituting this into the geometric identity $J = R/\phi_r$ and taking the limit $R \to 0$ yields $J(0, Z, t)
= [\p_R \phi_r(0, Z, t)]^{-1}$.

\item Consider the 2D Jacobian
$J = (\p_R \phi_r)(\p_Z \phi_z) - (\p_Z \phi_r)(\p_R \phi_z)$. Since $\phi_r$ is odd
in $R$, its vertical derivative $\p_Z \phi_r$ is also odd in $R$, hence $\p_Z \phi_r(0, Z, t) =
0$. Because $u_z$ is even in $r$, the map $\phi_z$ is even in $R$, so $\p_R \phi_z(0, Z, t) = 0$. Thus
on the axis $J = (\p_R \phi_r)(\p_Z \phi_z)$. Using part (b), this implies $\p_Z \phi_z =
J^2$.

\item Since $u_z$ is odd in $z$ (by the odd symmetry built into Definition~\ref{def:init-data}),
$(0,0)$ is a fixed stagnation point for all $t \in [0, T)$, hence $\phi_z(0,0,t)=0$.
\end{enumerate}
\end{proof}

The velocity gradient has the corresponding axis limits.
\begin{lemma}\label{lem2}
On the symmetry axis $r=0$, we have that
\begin{enumerate}
\item[(a)]
$\lim_{r \to 0} \tfrac{u_r}{r} = \p_r u_r(0, z) = -\tfrac{1}{2} \p_z u_z(0, z)$, and
\item[(b)] $S(0, z) = \tfrac{3}{2} (\p_z u_z(0, z))^2$.
\end{enumerate}
\end{lemma}
\begin{proof}[Proof of Lemma \ref{lem2}]
In axisymmetry, the divergence-free condition is $\p_r u_r + \tfrac{u_r}{r} + \p_z u_z = 0$.
Letting $r \to 0$ and applying L'H\^opital's rule to $\tfrac{u_r}{r}$, we find $2\p_r u_r + \p_z
u_z = 0$, which proves (a). For (b), note that on the axis $\p_z u_r(0,z)=0$ and $\p_r u_z(0,z)=0$
by parity, so the cross-term in \eqref{eq:Source} vanishes. Using also $\lim_{r\to 0}
\tfrac{u_r}{r}=\p_r u_r(0,z)$, we obtain that
\[
S(0,z) = 2(\p_r u_r)^2 + (\p_z u_z)^2.
\]
Substituting $\p_r u_r = -\tfrac{1}{2}\p_z u_z$ yields $S(0,z)=\tfrac{3}{2}(\p_z u_z)^2$.
\end{proof}

\subsection{The Fundamental Lagrangian Variables}
We define the symmetry-axis variables by
\begin{equation*}
J(Z, t) := J_{\twoD}(0, Z, t), \qquad \rW(Z, t) := \left.\p_z u_z(r,z,t)\right|_{r=0,\ z=\phi_z(0,Z,t)} .
\end{equation*}
By Lemma~\ref{lem:flow1}, $\phi_r(0,Z,t)=0$, so $\rW(Z,t)$ is the axial derivative of $u_z$ at the
Eulerian axis point reached from the label $(0,Z)$.

\subsubsection{Jacobian Evolution}
Differentiating \eqref{eq:Jac-Identity}, we obtain that
\begin{equation}
\p_t J = \p_t (\tfrac{R}{\phi_r}) = -\tfrac{R}{\phi_r^2} u_r(\phi, t) = -\tfrac{R}{\phi_r} \left( \tfrac{u_r}{\phi_r} \right) =
-J \left( \tfrac{u_r}{r} \right) \circ \phi \,.
\label{eq:J-evo0}
\end{equation}
From Lemma \ref{lem2}, $\lim_{r \to 0} \tfrac{u_r}{r} = -\tfrac{1}{2} \p_z u_z(0, z)$, and hence
$\tfrac{u_r}{r } \circ \phi \to -\tfrac{1}{2} \rW$ as $R \to 0$. Thus, from \eqref{eq:J-evo0}, we find that
\begin{equation}
\p_t J(Z, t) = \tfrac{1}{2} J(Z, t) \rW(Z, t) \,.
\label{eq:J-evo}
\end{equation}

\subsubsection{Momentum Evolution}
The vertical momentum equation reads $D_t u_z = -\p_z p$. Differentiating with respect to $z$ gives
\begin{equation*}
D_t (\p_z u_z) + (\p_z u_z)^2 + (\p_z u_r)(\p_r u_z) = -\p_z^2 p \,.
\end{equation*}
On the axis, $\p_z u_r(0,z)=0$ by parity, so the mixed term vanishes there. Let $P_{zz} := \p_z^2
p$. Composing with $\phi$, we obtain:
\begin{equation}\label{eq:rW-axis-evo}
\p_t \rW = -\rW^2 - P_{zz}\circ \phi \ \text{ on } \ \{R=0\} \,.
\end{equation}

The pressure is determined by $-\Delta p = S$, so $p = (-\Delta)^{-1}S$ for solutions decaying at infinity in
$\mathbb{R}^3$. Consequently,
\[
P_{zz}(x,t)=\p_z^2(-\Delta)^{-1}S(x,t)
\]
is a Calder\'on--Zygmund operator applied to $S$ and admits the standard decomposition
\begin{equation*}
P_{zz}(x,t) = \pv \int_{\R^3} K_{zz}(x-y)\,S(y,t)\,\ud y - \tfrac{1}{3} S(x,t),
\end{equation*}
where $K_{zz}(x)=\p_{x_3}^2\!\big(\tfrac{1}{4\pi|x|}\big)$.  We define the principal-value pressure
term $\Pi(Z,t)$ by evaluating this integral at the axis particle position
$x=(0,0,\phi_z(0,Z,t))$:
\begin{equation}
\Pi(Z,t) := \pv \int_{\R^3} K_{zz}\big((0,0,\phi_z(0,Z,t)) - y\big)\, S(y, t) \ud y \,.
\label{eq:pi-Eul}
\end{equation}
Applying Lemma \ref{lem2}, the source term on the axis is $S(0, z) = \tfrac{3}{2} \rW^2$. Thus,
\begin{equation}
P_{zz}(0,\phi_z, t) = \Pi(Z, t) - \tfrac{1}{3} \left( \tfrac{3}{2} \rW^2 \right) = \Pi(Z, t) - \tfrac{1}{2} \rW^2 \,.
\label{eq:Pzz}
\end{equation}
Substituting \eqref{eq:Pzz} into \eqref{eq:rW-axis-evo}, we obtain the evolution equation for the strain:
\begin{equation}
\tfrac{d}{dt} \rW = -\tfrac{1}{2} \rW^2 - \Pi \,.
\label{eq:rW}
\end{equation}

\subsection{The Biot--Savart Relation}
We recall the Biot--Savart law in $\R^3$:
\begin{equation}\label{eq:BS-3D}
u(x,t)=\int_{\R^3}\bs K(x-y)\times \bs{\omega}(y,t)\,\ud y,   \qquad   \bs K(z):=-\tfrac{1}{4\pi}\,\tfrac{z}{|z|^3}.
\end{equation}
In the axisymmetric no-swirl class $\bs{\omega}(y,t)=\omega_\theta(y,t)\,e_\theta(y)$, it is convenient to
write
\begin{equation}\label{eq:BS-axisymm}
 u(x,t)=\tfrac{1}{4\pi}\int_{\R^3}K(x,y)\,\omega_\theta(y,t)\,\ud y,
\qquad
K(x,y):=\tfrac{e_\theta(y)\times(x-y)}{|x-y|^3}.
\end{equation}
Differentiating in $x$ shows that $\nabla u$ is a Calder\'on--Zygmund singular integral of $\bs{\omega}$; in
particular
\begin{equation}
\nabla u(x,t)= \tfrac{1}{4\pi}\int_{\R^3}\mathcal K(x,y)\,\omega_\theta(y,t)\,\ud y, \qquad \mathcal K(x,y):=\nabla_x K(x,y).
\label{eq:grad-BS}
\end{equation}
For the axial strain on the symmetry axis we may write
\begin{equation}
\rW(Z,t) =\left.\p_z u_z(r,z,t)\right|_{r=0,\ z=\phi_z(0,Z,t)} =\tfrac{1}{4\pi}\int_{\R^3}\mathcal K_W(x,y)\,\omega_\theta(y,t)\,\ud y,
\quad x:=(0,0,\phi_z(0,Z,t)).
\label{eq:rW2}
\end{equation}
where
\[
\mathcal K_W(x,y):=\p_{x_z}\big(K(x,y)\cdot e_z\big).
\]

\subsection{Axial strain at the stagnation point}
Of fundamental importance to our analysis is the axial strain evaluated at the stagnation point:
\begin{equation*}
\rW_0(t) := \rW(0,t) \,.
\end{equation*}
We shall often refer to $\rW_0(t)$ as the \textit{driving strain} or \textit{axial strain}. From
\eqref{eq:rW}, we have that
\begin{subequations}
\label{eq:rW0-Pi0-J0-defs}
\begin{equation}
\tfrac{d}{dt} \rW_0 = -\tfrac{1}{2} \rW_0^2 - \Pi_0 \,,
\label{eq:rW0}
\end{equation}
where $\Pi_0(t)=\Pi(0,t)$.  By \eqref{eq:pi-Eul} and Lemma~\ref{lem:flow1},
\[
\Pi_0(t) := \pv \int_{\R^3} K_{zz}\big((0,0) - y\big)\, S(y, t) \ud y \,.
\]
We also set $J_0(t)= J(0,t)$, so that from \eqref{eq:J-evo},
\begin{equation}
\p_t J_0( t) = \tfrac{1}{2} J_0( t) \rW_0( t) \,,
\end{equation}
\end{subequations}
and we shall abuse notation and simply write $J(t)$ to mean $J_0(t)$.

\subsection{Coordinate systems and geometric notation}

Section~\ref{sec:proof-main} uses several coordinate systems and decompositions repeatedly. We collect the conventions here.

\subsubsection{Euclidean, cylindrical, and spherical coordinates.}
For $x=(x_1,x_2,x_3)\in\R^3$ we write
\[
\rho(x):=|x|=\sqrt{x_1^2+x_2^2+x_3^2}.
\]

We use cylindrical (axisymmetric) coordinates
\[
r(x):=\sqrt{x_1^2+x_2^2},\qquad z(x):=x_3,\qquad \theta(x):=\arg(x_1+i x_2)\in[0,2\pi),
\]
so that $x=(r\cos\theta,r\sin\theta,z)$. The standard cylindrical unit vectors are
\[
\bs e_r(x):=\tfrac{(x_1,x_2,0)}{r(x)},\qquad \bs e_\theta(x):=\tfrac{(-x_2,x_1,0)}{r(x)},\qquad \bs e_z:=(0,0,1) \ \ \text{ for } \ \ r(x)>0
.
\]

We also use a \emph{polar angle} $\sigma(x)\in[0,\tfrac{\pi}{2}]$ (measured from the $z$--axis, with $|z|$ to ignore the sign):
\begin{equation}
\sin\sigma(x)=\tfrac{r(x)}{\rho(x)},\qquad \cos\sigma(x)=\tfrac{|z(x)|}{\rho(x)}\ \ \text{ for } \ \ x\neq 0 ,
\label{eq:polar-angle-def}
\end{equation}
that is, $\sigma(x)=\arctan\!\big(\tfrac{r(x)}{|z(x)|}\big)$ for $z(x)\neq 0$ and $\sigma(x)=\tfrac{\pi}{2}$ when $z(x)=0$. Thus
\[
r=\rho\sin\sigma,\qquad |z|=\rho\cos\sigma.
\]

\subsubsection{Balls and complements.}
For $R>0$ we write
\[
B_R:=\{x\in\R^3:\ |x|<R\},\qquad B_R^c:=\R^3\setminus B_R.
\]

\subsubsection{Scaling notation.}
We write $f\lesssim g$ if $|f|\le C|g|$ in the stated regime, where $C$ may depend only on the fixed
parameters of the construction.  We write $f\simeq g$ when both $f\lesssim g$ and $g\lesssim f$ hold.

\section{An Open Set of Initial Data for Blowup}\label{sec:initialdata}

We now specify the explicit admissible class of finite-energy, axisymmetric no-swirl initial data used in the
blowup proof.  This class is designed so that the hyperbolic clock-and-driver model of
Section~\ref{sec:lag-analysis-non-local} has a nontrivial strain-producing region and the Euler solutions
generated by the data remain in $C^{1,\alpha}(\R^3)\cap L^2(\R^3)$.
The construction has two geometric requirements.

\runinhead{1. Algebraic tails under angular drift.} The hyperbolic collapse transports particles in polar
angle away from the symmetry axis toward the equator.  Under the linear rescaling that describes this collapse,
the label coordinates associated with an Eulerian point $x=(r,z)$ are
\[
(R,Z)=(J(t)\,r,J(t)^{-2}z), \qquad \rho_{\Lag}(x,t)=\sqrt{J(t)^2r^2+J(t)^{-4}z^2}.
\]
If $x$ remains in a fixed axial cone bounded away from the equator, then
$\rho_{\Lag}(x,t)\to\infty$ as $J(t)\downarrow0$.  Thus the axial strain at small clock values receives its
dominant contribution from labels that are farther and farther out in spherical radius.  The algebraic tail in
$\rho$ keeps those labels available, while the decay rate $\gamma>\alpha+\tfrac52$ gives finite energy for the
associated velocity.

\runinhead{2. Cone localization and the two angular kernels.} The angular vorticity is supported in a cone away
from the equator, and the support includes the nodal angle
\[
\sigma_{\node}:=\arccos\big(\tfrac{1}{\sqrt3}\big).
\]
At this angle, the pointwise pressure Hessian kernel $K_{zz}$ changes sign and vanishes, while the driving
strain kernel $K_W$ is maximal.  This kernel geometry motivates the angular support of the Target Profile.  The
quantitative pressure Hessian lower bound is not obtained from this observation alone; it is proved later by
the pressure Hessian model and its Euler realization in
Sections~\ref{sec:slope-restricted-pressure}--\ref{sec:target-profile-typeI-completion}.

The admissible class $\mathcal A_{\alpha,\gamma}(\nu,\eta)$ is a weighted open neighborhood of an explicit
target angular function $\Theta^*$ with a H\"older cusp at the axis,
$\Theta^*(\sigma)\simeq(\sin\sigma)^\alpha$ as $\sigma\downarrow0$.  The weight adds one more vanishing power at
the axis to the perturbative angular term; this is the quantitative gain used in the stability argument of
Section~\ref{sec:proof-main}.

\subsection{Geometry of the interaction kernels}

As seen in the axial strain evolution \eqref{eq:rW0}, the collapse dynamics depend on the coupling between
the driving strain $\rW_0(t):=\rW(0,t)$ and the nonlocal pressure Hessian $\Pi_0(t):=\Pi(0,t)$ defined
in \eqref{eq:pi-Eul}. We write the angular dependence of the corresponding kernels.

\subsubsection{The pressure Hessian kernel $K_{zz}$}\label{sec:Kzz-kernel}

The pressure solves $-\Delta p=S$, where $S=\tr((\nabla u)^2)$. Thus
\[
p = (-\Delta)^{-1}S = \Big(\tfrac{1}{4\pi|\cdot|}\Big)*S,
\]
and the nonlocal pressure Hessian on the axis is obtained from the singular integral operator
$\p_z^2(-\Delta)^{-1}$ (cf.\ \eqref{eq:pi-Eul}). In spherical coordinates $(\rho,\sigma,\varphi)$ with
$z=\rho\cos\sigma$, a direct computation shows that the kernel $K_{zz}(\rho,\sigma)=
\tfrac{\partial^2}{\partial z^2}\!\left(\tfrac{1}{4\pi\rho}\right)$ is given by
\begin{equation}
K_{zz}(\rho,\sigma)= \tfrac{1}{4\pi\rho^3}\,\big(3\cos^2\sigma-1\big).
\label{eq:Kzz-kernel}
\end{equation}
The kernel vanishes precisely when $\cos^2\sigma=\tfrac{1}{3}$, i.e.
\[
\sigma_{\node}=\arccos\!\left(\tfrac{1}{\sqrt3}\right),
\]
which we call the \emph{nodal cone}. This angle partitions the sphere into three regions:
\begin{itemize}
\item \textsf{Polar region ($0<\sigma<\sigma_{\node}$):} $3\cos^2\sigma-1>0$, so the kernel is
positive.
\item \textsf{Nodal cone ($\sigma=\sigma_{\node}$):} $K_{zz}=0$.
\item \textsf{Equatorial region ($\sigma_{\node}<\sigma<\tfrac{\pi}{2}$):} $3\cos^2\sigma-1<0$, so
the kernel is negative.
\end{itemize}
The nodal cone is therefore the pointwise angular location where the pressure Hessian kernel changes sign and
vanishes.  This sign change is only a geometric guide; the Riccati pressure estimate used in the blowup proof is
the nonlocal estimate proved later in
Sections~\ref{sec:slope-restricted-pressure}--\ref{sec:target-profile-typeI-completion}.

\subsubsection{The driving strain kernel $K_W$}\label{sec:KW-kernel}

The axial strain $\rW_0(t)=\p_z u_z(0,0,t)$ is given by the Biot--Savart formula \eqref{eq:rW2}. Specializing \eqref{eq:rW2} to a purely toroidal vorticity field
$\bs{\omega}=\omega_\theta e_\theta$ and differentiating along the axis gives an angular weight of the
form\footnote{Let $K(x,y)=\tfrac{e_\theta(y)\times(x-y)}{|x-y|^3}$ be the axisymmetric Biot--Savart kernel in
\eqref{eq:BS-axisymm}, so that $u(x)=\tfrac{1}{4\pi}\int_{\R^3}K(x,y)\,\omega_\theta(y)\,\ud y$. For an
axial point $x=(0,0,z_0)$, we have that $\mathcal K_W(x,y):=\p_{x_z}\big(K(x,y)\cdot e_z\big)
=-3\,\tfrac{r(y)\,(z_0-y_z)}{|x-y|^5}$. In particular, at the stagnation point $x=0$, $\mathcal
K_W(0,y)=3\,\tfrac{r(y)\,y_z}{|y|^5} =3\,\tfrac{\sin\sigma(y)\cos\sigma(y)}{|y|^3}$. Writing
$y=(\rho,\sigma,\varphi)$ so that $dy=\rho^2\sin\sigma\,d\rho\,d\sigma\,d\varphi$ and using axisymmetry,
$\rW_0(t)=\p_z u_z(0,0,t)=\tfrac{1}{4\pi}\int_{\R^3}\mathcal K_W(0,y)\,\omega_\theta(y,t)\,\ud y
=\tfrac12\int_0^\infty\int_0^\pi K_W(\sigma)\,\tfrac{\omega_\theta(\rho,\sigma,t)}{\rho}\,d\rho\,d\sigma$,
where $K_W(\sigma):=3\sin^2\sigma\,\cos\sigma$.}
\begin{equation}
K_W(\sigma)=3\sin^2\sigma\,\cos\sigma.
\label{eq:KW-kernel}
\end{equation}
Thus, vorticity with $\omega_\theta<0$ in the upper half-space, as produced by the sign convention in Definition~\ref{def:init-data}, generates 
compressive axial strain.  Moreover, $K_W$ vanishes at the symmetry axis and equator and attains its maximum precisely at the nodal angle, since
\[
\tfrac{d}{d\sigma}\big(\sin^2\sigma\cos\sigma\big)=\sin\sigma\,(3\cos^2\sigma-1).
\]
Thus the unique critical point in $(0,\tfrac{\pi}{2})$ occurs at $3\cos^2\sigma-1=0$, i.e. $\sigma=\sigma_{\node}$. The nodal cone is therefore 
simultaneously the location of \emph{maximal driving strain} and \emph{vanishing pressure kernel}, which motivates choosing the angular support so that it contains
$\sigma_{\node}$ but remains away from the equator.

\subsection{Admissible initial data}

We now define the admissible class $\mathcal{A}_{\alpha,\gamma}(\nu,\eta)$. The radial dependence is chosen to have a $C^\alpha$ cusp at the symmetry 
axis and an isotropic algebraic tail at infinity. The tail decays mildly enough to preserve the leading singular behavior, yet fast enough to ensure finite energy
(Lemma~\ref{lem:L2}). The angular dependence is a weighted H\"older neighborhood of an explicit target angular function $\Theta^*$ which is cone-localized 
away from the equator and is nonzero at the nodal angle. The function $\Theta^*$ is \emph{not} smooth at the axis: the term $(\sin\sigma)^\alpha$ produces a H\"older
cusp $\Theta^*(\sigma)\simeq \sigma^\alpha$ as $\sigma\downarrow0$. The weight exponent $\eta>0$ enforces additional vanishing of perturbations at the axis, so the 
perturbative term is lower order in the collapsing strain-producing sector. In the following definition, we write $\rho=\sqrt{R^2+Z^2}$ for the spherical radius and
$R=\rho\sin\sigma$ for the cylindrical radius.

\begin{definition}[Admissible Initial Data Class $\mathcal A_{\alpha,\gamma}(\nu,\eta)$]
\label{def:init-data}
We define the class of admissible initial data $\mathcal{A}_{\alpha,\gamma}(\nu,\eta)$ through the following construction.

\begin{enumerate}

\item {\sf The target angular function $\Theta^*$ (cone localization, H\"older cusp, and odd
symmetry):} We fix angles
\[
0<\sigma_{\cut}<\sigma_{\node}<\sigma_{\max}<\tfrac{\pi}{2},
\]
where $\sigma_{\node}=\arccos(1/\sqrt{3})$ is the maximizer of the strain kernel $K_W$, and let
$\Upsilon(\sigma)$ be a smooth nonincreasing cutoff satisfying $0\le\Upsilon\le1$ and
\[
\Upsilon(\sigma)=1 \ \text{for } 0\le\sigma\le\sigma_{\cut}, \qquad \Upsilon(\sigma)>0 \ \text{for } 0\le\sigma<\sigma_{\max},
\qquad \Upsilon(\sigma)=0 \ \text{for } \sigma\ge\sigma_{\max}.
\]
We define the reference angular function on the upper hemisphere by
\begin{equation}
\Theta^*(\sigma):=(\sin\sigma)^\alpha\,\Upsilon(\sigma), \qquad \sigma\in[0,\tfrac\pi2],
\label{eq:Theta-star-def}
\end{equation}
and extend it to $[0,\pi]$ by odd reflection:
\begin{equation*}
\Theta^*(\pi-\sigma):=-\Theta^*(\sigma), \qquad \sigma\in[0,\tfrac\pi2].
\end{equation*}
The resulting angular function is odd about the equator.  For the associated velocity this gives $u_z(r,-z)=-u_z(r,z)$ and $u_r(r,-z)=u_r(r,z)$, so $(r,z)=(0,0)$ is a stagnation point.

\noindent
\emph{Regularity note.} Since $\alpha\in(0,1)$, $\Theta^*$ has a H\"older cusp at the axis:
$\Theta^*(\sigma)\simeq \sigma^\alpha$ as $\sigma\downarrow0$. In particular, $\Theta^*\in C^\alpha$ but is not $C^1$ at $\sigma=0$ when $\alpha<1$.

\item {\sf The admissible neighborhood (axis-vanishing perturbations):} We fix an exponent $\gamma>\alpha+\tfrac52$ (the spherical tail decay rate), 
and parameters $\nu>0$ and $\eta>0$. The set $\mathcal{A}_{\alpha,\gamma}(\nu,\eta)$ consists of all initial toroidal vorticity components
$\omega_{\theta,0}(\rho,\sigma)$ of the form
\begin{equation}
\omega_{\theta,0}(\rho,\sigma) = -\Gamma\,\tfrac{\rho^\alpha}{(1+\rho^2)^{\gamma/2}}\,\Theta(\sigma),
\label{eq:vort0}
\end{equation}
where $\Gamma>0$ and the angular function $\Theta:[0,\pi]\to\R$ satisfies:

\smallskip
\noindent\emph{(i) Odd symmetry:} $\Theta(\pi-\sigma)=-\Theta(\sigma)$ for all $\sigma\in[0,\pi]$.

\smallskip
\noindent\emph{(ii) Weighted H\"older proximity to $\Theta^*$:} The angular function on $[0,\tfrac\pi2]$ is
generated by
\[
\Theta(\sigma)=\Theta^*(\sigma)\cdot(1+h(\sigma)),
\]
where $h\in C^\alpha_\eta([0,\tfrac\pi2])$.  This means that the perturbation has the form
\[
h(\sigma)=(\sin\sigma)^\eta k(\sigma), \qquad k\in C^\alpha([0,\tfrac\pi2]),
\]
and satisfies
\begin{equation}
\|h\|_{C^\alpha_\eta}:=\|k\|_{C^\alpha([0,\tfrac\pi2])}<\nu.
\label{eq:weighted-proximity}
\end{equation}
(Since $\Theta^*$ and $\Theta$ are both odd, specifying $h$ on $[0,\tfrac\pi2]$ determines $\Theta$ on all of
$[0,\pi]$.)  In particular,
\[
\sup_{\sigma\in[0,\tfrac\pi2]}\tfrac{|h(\sigma)|}{(\sin\sigma)^\eta} \le \|h\|_{C^\alpha_\eta}<\nu .
\]

\smallskip
\noindent
Thus $h(0)=0$, so the perturbation vanishes at the symmetry axis to order $\eta$. This weighted
neighborhood is open in the weighted H\"older topology defined by \eqref{eq:weighted-proximity}.  Moreover,
since $(\sin\sigma)^\eta\le 1$, we have $\|h\|_{L^\infty}\le \|h\|_{C^\alpha_\eta}<\nu$.

\smallskip
\noindent
The associated initial velocity is
\[
u_0:=\BS[\omega_{\theta,0}\,\bs e_\theta].
\]
The odd symmetry gives the same parity for the associated velocity.  This parity is preserved by the Euler flow, so $(0,0)$ remains a stagnation point while the solution exists.
\end{enumerate}
\end{definition}

\begin{remark}[A genuine weighted neighborhood]
Definition~\ref{def:init-data} is a genuine weighted neighborhood of the target angular function $\Theta^*$, specified
entirely at the level of the initial datum. The restriction that $\nu$ be sufficiently small in Theorem~\ref{thm:main} is imposed only after the 
Target Profile collapse mechanism and the pressure estimates have been established. In particular, the definition does not build in any separate ``stability
margin'' condition or sign condition.
\end{remark}

\begin{lemma}[Axis regularity of the initial toroidal cusp]
\label{lem:initial-toroidal-cusp-holder}
Let $0<\alpha<1$ and define
\[
\mathfrak t(x):=
\begin{cases}
r(x)^\alpha\bs e_\theta(x), & r(x)>0,\\
0, & r(x)=0.
\end{cases}
\]
Then $\mathfrak t\in C^\alpha(\R^3)$ and
\[
|\mathfrak t(x)-\mathfrak t(y)|\le C_\alpha |x-y|^\alpha \qquad (x,y\in\R^3).
\]
\end{lemma}

\begin{proof}[Proof of Lemma~\ref{lem:initial-toroidal-cusp-holder}]
We write $r_x=r(x)$ and $r_y=r(y)$.  If $\min\{r_x,r_y\}\le2|x-y|$, then
\[
|\mathfrak t(x)-\mathfrak t(y)|\le r_x^\alpha+r_y^\alpha\le C|x-y|^\alpha.
\]
If $\min\{r_x,r_y\}>2|x-y|$, then $r_x\simeq r_y$ and both points are away from the axis.  Since
$|\nabla \bs e_\theta|\lesssim r_x^{-1}$ and $|\nabla r^\alpha|\lesssim r_x^{\alpha-1}$ in this region, the mean value theorem gives
\[
|\mathfrak t(x)-\mathfrak t(y)|\le C r_x^{\alpha-1}|x-y|.
\]
The condition $r_x>2|x-y|$ implies
$r_x^{\alpha-1}|x-y|\le C|x-y|^\alpha$, because $\alpha-1<0$.  This proves the estimate.
\end{proof}

\begin{lemma}[Regularity and finite energy of the initial velocity]\label{lem:L2}
Let $\omega_{\theta,0}$ be given by \eqref{eq:vort0} with $\gamma>\alpha+\tfrac52$ and angular function
$\Theta(\sigma)=\Theta^*(\sigma)(1+h(\sigma))$ satisfying \eqref{eq:weighted-proximity}, and let
$u_0=\BS[\omega_{\theta,0}\,\bs e_\theta]$ be the associated initial velocity.
\begin{enumerate}
\item[\textnormal{(a)}] \textsf{Finite energy.} $u_0\in L^2(\R^3)$.
\item[\textnormal{(b)}] \textsf{Target Profile regularity.} If $h\equiv 0$ (the exact Target
Profile), then $u_0\in C^{1,\alpha}(\R^3)\cap L^2(\R^3)$.
\item[\textnormal{(c)}] \textsf{Weighted-H\"older perturbations.} If
$h\in C^\alpha_\eta([0,\tfrac\pi2])$, then $u_0\in C^{1,\alpha}(\R^3)\cap L^2(\R^3)$.
\end{enumerate}
\end{lemma}

\begin{remark}[The regularity hypothesis in
Theorem~\ref{thm:main}]
The admissible class $\mathcal{A}_{\alpha,\gamma}(\nu,\eta)$ is defined using the weighted H\"older topology
\eqref{eq:weighted-proximity}.  This controls both the \emph{size} of the perturbation and the local
H\"older seminorms needed for the singular-integral estimates in Section~\ref{sec:proof-main}.  In particular,
the perturbative part of the angular function has the form
\[
\Theta^*(\sigma)h(\sigma) = (\sin\sigma)^{\alpha+\eta}\Upsilon(\sigma)k(\sigma), \qquad \|k\|_{C^\alpha}<\nu,
\]
which has the target $C^\alpha$ cusp structure and vanishes to an additional order at the axis. Lemma~\ref{lem:L2} therefore gives finite energy 
and the required $C^{1,\alpha}$ regularity for the velocity. We keep the regularity assumption in Theorem~\ref{thm:main} explicit to emphasize the local well-posedness
class, but for the weighted H\"older admissible perturbations defined above it is automatically satisfied.

For the Target Profile (Theorem~\ref{thm:target-profile}), $h\equiv 0$ and the full regularity $u_0\in C^{1,\alpha}(\R^3)\cap L^2(\R^3)$ is 
established unconditionally by part~(b).
\end{remark}

\begin{proof}[Proof of Lemma~\ref{lem:L2}]
\runinhead{Part (a): Finite energy.} This part uses only the pointwise bound
\begin{equation*}
|\omega_{\theta,0}(y)| \le C_\nu\,\Gamma\,(1+|y|)^{\alpha-\gamma},
\end{equation*}
which follows immediately from \eqref{eq:vort0}, the weighted bound $\|h\|_{L^\infty}\le
\|h\|_{C^\alpha_\eta}<\nu$, and $0\le\Upsilon\le1$.

Evaluating the three-dimensional Biot--Savart law \eqref{eq:BS-3D} at $t=0$, we have
\begin{equation*}
u_0(x)=\int_{\R^3}\bs K(x-y)\times\bigl(\omega_{\theta,0}(y)\bs e_\theta(y)\bigr)\,\ud y,
\end{equation*}
with $|\bs K(z)|\le C|z|^{-2}$. In particular,
\[
|u_0(x)|\le C\int_{\R^3}\tfrac{|\omega_{\theta,0}(y)|}{|x-y|^2}\,\ud y.
\]
For $|x|\le 1$, the local boundedness of the Biot--Savart operator on bounded vorticity gives $u_0\in L^2(\{|x|\le 1\})$. For $|x|\ge 1$, we decompose the source 
region into $\{|y|\le \tfrac12|x|\}\cup\{\tfrac12|x|<|y|\le 2|x|\}\cup\{|y|>2|x|\}$ and estimate each region as follows.

\smallskip
\noindent
\runinhead{(i) Inner source region $|y|\le \tfrac{1}{2}|x|$: far-field expansion and cancellation of the monopole.} On $|y|\le \tfrac{1}{2}|x|$ we 
have $|x-y|\ge \tfrac12|x|$, so the Biot--Savart kernel $\bs K(x-y)$ in \eqref{eq:BS-3D} is smooth and we may expand in the small ratio $|y|/|x|$ by Taylor expanding in
the source variable $y$ about $0$:
\begin{equation*}
\bs K(x-y)=\bs K(x)-(\nabla \bs K)(x)y+O\!\left(\tfrac{|y|^2}{|x|^4}\right).
\end{equation*}
The monopole integral $\bs K(x)\times \int_{|y|\le \tfrac12|x|}\omega_{\theta,0}(y)\bs e_\theta(y)\,\ud y$ vanishes because the azimuthal average of $\bs e_\theta$ is zero:
$\int_{0}^{2\pi}\bs e_\theta(\theta)\,\ud\theta=0$. Therefore,
\[
\Big|\int_{|y|\le \tfrac12|x|}\bs K(x-y)\times(\omega_{\theta,0}(y)\bs e_\theta(y))\,\ud y\Big| \lesssim |x|^{-3}\!\!\int_{|y|\le \tfrac12|x|}\!|y|\,|\omega_{\theta,0}(y)|\,\ud y
+|x|^{-4}\!\!\int_{|y|\le \tfrac12|x|}\!|y|^2|\omega_{\theta,0}(y)|\,\ud y.
\]
We set
\[
\delta_u:=\min\{\gamma-\alpha-1,\tfrac52\}.
\]
Since $\gamma>\alpha+\tfrac52$, we have $\delta_u>\tfrac32$.  Using
$|\omega_{\theta,0}(y)|\lesssim \Gamma(1+|y|)^{\alpha-\gamma}$, the two moments above satisfy
\[
|x|^{-3}\!\!\int_{|y|\le \tfrac12|x|}|y|\,|\omega_{\theta,0}(y)|\,dy +|x|^{-4}\!\!\int_{|y|\le \tfrac12|x|}|y|^2|\omega_{\theta,0}(y)|\,dy
\lesssim \Gamma |x|^{-\delta_u}.
\]
The assumption $\gamma>\alpha+\tfrac52$ gives $\gamma-\alpha-1>\tfrac32$, which is the integrability
threshold needed for $u_0\in L^2$.  If the algebraic tail gives the sharper exponent
$\gamma-\alpha-1>\tfrac52$, we retain only the capped decay $|x|^{-5/2}$, since no later estimate uses a larger far-field exponent.
Therefore the inner-source contribution satisfies $|u_0(x)|_{\mathrm{near}}\lesssim\Gamma |x|^{-\delta_u}$.

\smallskip
\noindent
\runinhead{(ii) Outer source region $|y|>\tfrac12|x|$: crude kernel bound plus tail decay.} We split further into $\tfrac12|x|<|y|\le 2|x|$ and $|y|>2|x|$.

\smallskip
\noindent
\runinhead{(ii-a) Comparable radii: $\tfrac12|x|<|y|\le 2|x|$.} Here $|y|\simeq |x|$, hence
$|\omega_{\theta,0}(y)|\lesssim \Gamma\,|x|^{\alpha-\gamma}$. Using $|\bs K(x-y)|\lesssim |x-y|^{-2}$ and the change of variables $z=x-y$,
\[
|u_0(x)|_{\mathrm{ann}} \lesssim \Gamma\,|x|^{\alpha-\gamma}\int_{\tfrac12|x|<|y|\le 2|x|}\tfrac{dy}{|x-y|^2}
\lesssim \Gamma\,|x|^{\alpha-\gamma+1} \le C\Gamma |x|^{-\delta_u}.
\]

\smallskip
\noindent
\runinhead{(ii-b) Very far sources: $|y|>2|x|$.} In this region $|x-y|\ge \tfrac12|y|$, so
\[
|u_0(x)|_{\mathrm{far}} \lesssim \int_{|y|>2|x|}\tfrac{|\omega_{\theta,0}(y)|}{|y|^2}\,dy
\lesssim \Gamma\int_{2|x|}^\infty \rho^{\alpha-\gamma}\,d\rho \lesssim \Gamma\,|x|^{\alpha-\gamma+1} \le C\Gamma |x|^{-\delta_u},
\]
using $\gamma>\alpha+1$.

\smallskip
\noindent
Combining the three regions yields
\begin{equation*}
|u_0(x)|\lesssim \Gamma\,|x|^{-\delta_u}, \qquad |x|\ge 1.
\end{equation*}
Therefore,
\[
\int_{|x|\ge 1}|u_0(x)|^2\,\ud x \lesssim \Gamma^2\int_1^\infty \rho^{-2\delta_u}\rho^2\,d\rho<\infty,
\]
because $\delta_u>\tfrac32$. This proves $u_0\in L^2(\R^3)$.

\bigskip
\runinhead{Part (b): $C^{1,\alpha}$ regularity for the Target Profile.} We now assume $h\equiv 0$, so $\Theta=\Theta^*$ and 
$\Theta^*(\sigma)=(\sin\sigma)^\alpha\Upsilon(\sigma)$ is smooth on $(0,\sigma_{\max})$, odd across $\sigma=\tfrac{\pi}{2}$, and has a 
$C^\alpha$ cusp at $\sigma=0$. The proof proceeds in three steps.

\rruninhead{Step 1: local $C^{1,\alpha}$ regularity.} We write $\bs{\omega}_0=\omega_{\theta,0}\bs e_\theta$ in cylindrical coordinates. In spherical 
variables $(\rho,\sigma)$ (so that $r=\rho\sin\sigma$), the definition~\eqref{eq:vort0} gives, near $\rho=0$,
\[
|\omega_{\theta,0}(\rho,\sigma)| \lesssim \Gamma\,\rho^\alpha\,|\Theta^*(\sigma)| \lesssim \Gamma\,\rho^\alpha(\sin\sigma)^\alpha =\Gamma\,r^\alpha.
\]
Since $\omega_{\theta,0}(r,z)=O(r^\alpha)$ as $r\downarrow0$ with $0<\alpha<1$, each Cartesian component of
$\bs{\omega}_0$ extends continuously across the axis.  The elementary fact that
$r^\alpha\bs e_\theta$ is a $C^\alpha$ vector field across the axis, Lemma~\ref{lem:initial-toroidal-cusp-holder},
and the smoothness of $\Upsilon$ away from $\sigma=0$ imply $\bs{\omega}_0\in C^\alpha_{\loc}(\R^3)$.

The Biot--Savart law expresses $\nabla u_0$ as a Calder\'on--Zygmund singular integral of $\bs{\omega}_0$. By
the standard Hölder boundedness of Calder\'on--Zygmund operators (see, for example,
Stein~\cite[Chapter~V]{Stein1993}), it follows that $\nabla u_0\in C^\alpha_{\loc}(\R^3)$.

\rruninhead{Step 2: dyadic Calder\'on--Zygmund control at infinity.} For derivative estimates, we use the
localized $C^\alpha$ estimate for singular integrals, in the form of the interior Schauder estimate for
Newtonian potentials: if $T$ is the gradient of the Biot--Savart operator and $f\in C^\alpha$ on a ball
$B(x_0,2r)$, then the principal-value cancellation gives a scale-invariant bound for
$T f$ in $C^\alpha(B(x_0,r))$ in terms of the $C^\alpha$ norm of $f$ on $B(x_0,2r)$, together with kernel
bounds for sources outside the larger ball.  We use this standard estimate below; see Stein~\cite[Chapter~V]{Stein1993}
or Gilbarg--Trudinger~\cite{GilbargTrudinger}.

We set
\[
\delta_{\nabla}:=\min\{\gamma-\alpha,\tfrac72\}.
\]
For every dyadic radius $R\ge 2$, the target vorticity satisfies
\begin{equation}\label{eq:omega-dyadic-holder-sec5}
\|\bs\omega_0\|_{L^\infty(\{\tfrac{1}{2}R\le |x|\le 2R\})} \le C\Gamma R^{\alpha-\gamma},\qquad [\bs\omega_0]_{C^\alpha(\{\tfrac{1}{2}R\le |x|\le 2R\})}
 \le C\Gamma R^{-\gamma}.
\end{equation}
Indeed, after the rescaling $x=R\tilde x$, the radial function contributes $R^{\alpha-\gamma}$, while the $C^\alpha$ seminorm loses the additional power 
$R^\alpha$. The vector field $r^\alpha\bs e_\theta$ remains uniformly $C^\alpha$ across the axis under this rescaling by Lemma~\ref{lem:initial-toroidal-cusp-holder}.

We fix $x$ with $R:=|x|\ge4$, and we let $\eta_x$ be a smooth cutoff supported in $B(x,\tfrac{1}{8} R)$, equal to one on
$B(x,\tfrac{1}{16} R)$, with $|\nabla^m\eta_x|\lesssim R^{-m}$. Writing $T:=\nabla\BS$, the local part
$T[\eta_x\bs\omega_0]$ is a Calder\'on--Zygmund singular integral of a compactly supported $C^\alpha$ function.
Scaling the standard Schauder estimate on $B(x,\tfrac{1}{8} R)$, we have that
\begin{equation}\label{eq:local-CZ-u0-sec5}
|T[\eta_x\bs\omega_0](x)| +R^\alpha [T[\eta_x\bs\omega_0]]_{C^\alpha(B(x,R/32))} \le C\Gamma R^{\alpha-\gamma}.
\end{equation}

We next estimate the complementary source $(1-\eta_x)\bs\omega_0$.  Let $z,z'\in B(x, \tfrac{1}{32} R)$, and let $T_{\rm inn}$ denote the 
contribution of the region $|y|\le \tfrac{1}{2} R$.  On this region, the kernel of $T$ can be expanded in $\tfrac y R$, and the constant term again has 
zero azimuthal average. The pointwise term is bounded by
\[
R^{-4}\int_{|y|\le R/2}|y|\,|\omega_{\theta,0}(y)|\,dy \le C\Gamma R^{-\delta_{\nabla}},
\]
and applying this expansion to the kernel difference produces the Hölder bound
\[
|T_{\rm inn}(z)-T_{\rm inn}(z')| \le C\Gamma R^{-\delta_{\nabla}-\alpha}|z-z'|^\alpha.
\]
On the comparable region $\tfrac12R<|y|\le2R$ outside $B(x,\tfrac{1}{16} R)$, the distance from $z$ to the source is bounded below by $cR$.
The bounds
\[
|\nabla\BS(z-y)|\lesssim |z-y|^{-3}, \qquad |\nabla\BS(z-y)-\nabla\BS(z'-y)|\lesssim |z-z'|^\alpha |z-y|^{-3-\alpha}
\]
and \eqref{eq:omega-dyadic-holder-sec5} yield a pointwise contribution $C\Gamma R^{\alpha-\gamma}$ and a
Hölder contribution $C\Gamma R^{-\gamma}|z-z'|^\alpha$.  On the far region $|y|>2R$, these kernel bounds and
$|z-y|\simeq |y|$ give a pointwise contribution $C\Gamma R^{\alpha-\gamma}$ and a Hölder contribution
$C\Gamma R^{-\gamma}|z-z'|^\alpha$, using $\gamma>\alpha+1$.  Combining the three source regions,
and using $\delta_{\nabla}\le\gamma-\alpha$, we have that
\begin{equation}\label{eq:far-CZ-u0-sec5}
|T[(1-\eta_x)\bs\omega_0](x)| +R^\alpha [T[(1-\eta_x)\bs\omega_0]]_{C^\alpha(B(x,R/32))} \le C\Gamma R^{-\delta_{\nabla}}.
\end{equation}
Combining \eqref{eq:local-CZ-u0-sec5} and \eqref{eq:far-CZ-u0-sec5}, and again keeping only the capped far-field decay that is needed below, 
we obtain the far-field gradient and local H\"older bounds
\begin{equation}\label{eq:grad-u0-decay}
|\nabla u_0(x)| \le C\Gamma |x|^{-\delta_{\nabla}}, \qquad
[\nabla u_0]_{C^\alpha(B(x,|x|/32))} \le C\Gamma |x|^{-\delta_{\nabla}-\alpha} \quad\text{for } |x|\ge 4 .
\end{equation}

\rruninhead{Step 3: global $C^\alpha$ seminorm.} We bound the H\"older ratio $\tfrac{|\nabla u_0(x)-\nabla u_0(y)|}{|x-y|^\alpha}$ uniformly over all $x\neq y\in\R^3$.

\rruninhead{Step 3a: near--near ($|x|,|y|\le 4$).} The local estimate from Step~1 gives
\[
\tfrac{|\nabla u_0(x)-\nabla u_0(y)|}{|x-y|^\alpha}\le C\Gamma .
\]

\rruninhead{Step 3b: far--far ($|x|,|y|\ge 2$).} If $|x-y|\le \tfrac{|x|}{32}$, the local estimate \eqref{eq:grad-u0-decay} gives
\[
\tfrac{|\nabla u_0(x)-\nabla u_0(y)|}{|x-y|^\alpha} \le C\Gamma |x|^{-\delta_{\nabla}-\alpha}.
\]
If $|x-y|>\tfrac{|x|}{32}$, the triangle inequality and the pointwise part of \eqref{eq:grad-u0-decay} yield
\[
\tfrac{|\nabla u_0(x)-\nabla u_0(y)|}{|x-y|^\alpha} \lesssim \tfrac{\Gamma(|x|^{-\delta_{\nabla}}+|y|^{-\delta_{\nabla}})}{|x|^\alpha} \lesssim \Gamma .
\]

\rruninhead{Step 3c: near--far ($|x|\le4$, $|y|\ge8$, or vice versa).} Then $|x-y|\ge \tfrac{1}{2}|y|$, so
\[
\tfrac{|\nabla u_0(x)-\nabla u_0(y)|}{|x-y|^\alpha} \lesssim \tfrac{\|\nabla u_0\|_{L^\infty(B_4)}+\Gamma |y|^{-\delta_{\nabla}}}{|y|^\alpha}
\lesssim \|\nabla u_0\|_{L^\infty(B_4)}+\Gamma .
\]
Combining the three regimes gives $[\nabla u_0]_{C^\alpha(\R^3)}<\infty$.  The local Biot--Savart bound near the
origin and the far-field decay from Part~(a) give $\|u_0\|_{L^\infty(\R^3)}<\infty$.  Hence $u_0\in C^{1,\alpha}(\R^3)$.

\bigskip
\runinhead{Part (c): $C^{1,\alpha}$ regularity for weighted-H\"older perturbations.} We write
$h(\sigma)=(\sin\sigma)^\eta k(\sigma)$ with $k\in C^\alpha([0,\tfrac\pi2])$. Then, on the upper hemisphere,
\[
\Theta(\sigma) = (\sin\sigma)^\alpha\Upsilon(\sigma) + (\sin\sigma)^{\alpha+\eta}\Upsilon(\sigma)k(\sigma).
\]
The first term is the Target Profile term treated in Part~(b). The second term is also $C^\alpha$ up to the axis. Indeed,
\[
\rho^\alpha(\sin\sigma)^{\alpha+\eta}\Upsilon(\sigma)k(\sigma)e_\theta = r^\alpha(\sin\sigma)^\eta\Upsilon(\sigma)k(\sigma)e_\theta .
\]
The multiplier $(\sin\sigma)^\eta\Upsilon(\sigma)k(\sigma)$ is bounded.  If two points have distance at least comparable to their distance from the axis, the 
right-hand side is bounded by $C r^\alpha$ at those points and hence satisfies the desired $C^\alpha$ estimate.  If the two points stay a distance $r$ from the 
axis and are separated by at most a fixed multiple of $r$, the angular multiplier has $C^\alpha$ seminorm $O(r^{-\alpha})$,
and $r^\alpha \bs e_\theta$ compensates this possible loss.  This is the same estimate as in Lemma~\ref{lem:initial-toroidal-cusp-holder}; the extra multiplier 
$(\sin\sigma)^\eta$ never worsens the behavior at the axis.  After the odd reflection across the equator, $\Theta\in C^\alpha([0,\pi])$.  Consequently,
$\bs\omega_0\in C^\alpha_{\loc}(\R^3)$ by the local argument in Part~(b).

It remains only to check that the dyadic bounds used in Step~2 still hold.  On $\{\tfrac12R\le |x|\le2R\}$, the scalar radial multiplier in
\eqref{eq:vort0} has size $R^{\alpha-\gamma}$ and $C^\alpha$ seminorm $CR^{-\gamma}$.  The angular multiplier
\[
(\sin\sigma)^\alpha\Upsilon(\sigma) +(\sin\sigma)^{\alpha+\eta}\Upsilon(\sigma)k(\sigma)
\]
has a $C^\alpha$ norm bounded by $C(1+\|k\|_{C^\alpha})$ uniformly after the rescaling $x=R\tilde x$.
Therefore,
\begin{equation*}
\|\bs\omega_0\|_{L^\infty(\{\tfrac{1}{2}R\le |x|\le 2R\})} \le C(1+\nu)\Gamma R^{\alpha-\gamma},\qquad [\bs\omega_0]_{C^\alpha(\{\tfrac{1}{2}R\le |x|\le 2R\})}
\le C(1+\nu)\Gamma R^{-\gamma}.
\end{equation*}
These are the dyadic estimates \eqref{eq:omega-dyadic-holder-sec5}, with a constant depending on the admissible
neighborhood.  Steps~2--3 of Part~(b) therefore provide $\nabla u_0\in C^\alpha(\R^3)$.  Part~(a) already
shows that $u_0\in L^2(\R^3)$.
\end{proof}

\begin{remark}[Interpreting the angular cusp]
In the poloidal label variables $Y=(R,Z)$ we have $\sin\sigma(Y)=\tfrac{R}{\rho(Y)}$, where
$\rho(Y):=\sqrt{R^2+Z^2}$.  Thus, for the target angular function
$\Theta^*(\sigma)=(\sin\sigma)^\alpha\Upsilon(\sigma)$, the spherical datum \eqref{eq:vort0} takes the form
\[
\rho(Y)^\alpha\,\Theta^*(\sigma(Y)) = \rho(Y)^\alpha\big(\tfrac{R}{\rho(Y)}\big)^{\!\alpha}\Upsilon(\sigma(Y)) = R^\alpha\,\Upsilon(\sigma(Y)).
\]
Thus, for each fixed $Z$ the vorticity vanishes like $R^\alpha$ as $R\downarrow0$.  Near the origin we also
have the uniform bound $R^\alpha\lesssim \rho^\alpha$.  This is exactly the local axisymmetric cusp needed in
Lemma~\ref{lem:L2}.

We write the cusp in terms of the polar angle $\sigma$ rather than only in $R$ because the collapse induces
rapid \emph{drift in $\sigma$}.  The function $(\sin\sigma)^\alpha$ is exactly what produces the depletion
power $J(t)^{3\alpha}$ when the strain-producing sector is drawn from small-clock labels.
\end{remark}

\begin{remark}[Amplitude and perturbation parameters]
The admissible class $\mathcal{A}_{\alpha,\gamma}(\nu,\eta)$ separates the perturbation size from the amplitude scale.  We fix the parameters in the order
\[
\text{fix }\eta>0,\qquad\text{then choose }\nu>0\text{ sufficiently small}, \qquad \Gamma>0\ \text{arbitrary}.
\]
The exponent $\eta>0$ is chosen \emph{a priori}; it defines the weighted H\"older topology in \eqref{eq:weighted-proximity} and enforces that 
admissible perturbations vanish at the axis.  Thus these perturbations are lower order in the collapsing core and remain perturbative in the singular-integral 
estimates. The amplitude $\Gamma$ is not a smallness parameter.  After the natural time rescaling $\tau=\Gamma t$, the
Target Profile construction and the stability estimates are uniform in $\Gamma$.  Finally, the neighborhood size $\nu>0$ is chosen sufficiently small, depending 
only on $\alpha,\gamma,\eta$, so that the initial angular function stays close to $\Theta^*$ and the stability bootstrap closes.
\end{remark}


\section{Proof Strategy}
\label{sec:outline}

\subsection{The linear hyperbolic guide}
We now give a roadmap to the proof of Theorem~\ref{thm:target-profile} and Theorem~\ref{thm:main}.  The
starting point is the linear hyperbolic model of Section~\ref{sec:lag-analysis-non-local}.  The Target Profile initial vorticity
has $\Theta^*(\sigma)=(\sin\sigma)^\alpha\Upsilon(\sigma)$, and on the sector where $\Upsilon=1$,
\[
\rho^\alpha\Theta^*(\sigma)=R^\alpha .
\]
Thus,  the singular cusp-part of the transported vorticity is a cylindrical cusp $R^\alpha$.  
In  Section~\ref{sec:lag-analysis-non-local}, we build our clock-and-driver model; we  transport this  initial cusp vorticity
by a linear hyperbolic map
\[
\Phi_{\lin}(R,Z,t)=\bigl(R\Jm(t)^{-1},Z\Jm(t)^2\bigr),
\]
where $\Jm(t)$ is a model clock function.  Specifically, our clock-and-driver model is a system for the model clock $\Jm(t)$ and a model
axial strain $\rWm(t)$ at the stagnation point; the strain is defined by the Biot--Savart law from the model vorticity transported by $\Phi_{\lin}$:
\[
\dot\Jm(t)=\tfrac12\Jm(t)\rWm(t),\qquad \rWm(t)\simeq-\Gamma\Jm(t)^{3\alpha-1},\qquad \dot\Jm(t)\simeq-\Gamma\Jm(t)^{3\alpha}.
\]
Our clock-and-driver model explains the exponent $\alpha=\tfrac13$: the model clock reaches zero in finite time precisely in
the range $3\alpha<1$.  Sections~\ref{sec:Euler-blowup-for-Theta-star}--\ref{sec:target-profile-typeI-completion} prove
that the true Euler solution tracks this model through the collapse.

\subsection{Euler behaves like the clock-and-driver model}
For the true Euler solution, the three important stagnation-point functions are the true Euler clock $J(t)=\det\nabla_{(R,Z)}(\phi_r,\phi_z)(0,0,t)$, the true Euler axial strain
$\rW_0(t)=\p_z u_z(0,0,t)$, and the true Euler pressure Hessian $\Pi_0(t)$.  The clock  $J(t)$ satisfies
\begin{equation}
\dot J(t)=\tfrac12\,J(t)\,\rW_0(t).
\label{eq:outline-J-evo}
\end{equation}
The Euler axial strain is constrained by the Riccati law
\begin{equation}
\p_t\rW_0(t)=-\tfrac12\,\rW_0(t)^2-\Pi_0(t).
\label{eq:outline-Riccati}
\end{equation}
To prove finite-time blowup for the Euler solutions, we  must prove the two estimates
\begin{equation}
-\rW_0(t)\simeq\Gamma J(t)^{3\alpha-1},\qquad \Pi_0(t)\ge-\upbeta\,\tfrac12\,\rW_0(t)^2\quad\hbox{for some }0<\upbeta<1.
\label{eq:outline-two-estimates}
\end{equation}
The first estimate in \eqref{eq:outline-two-estimates} is the Euler analogue of $-\rWm(t)\simeq\Gamma\Jm(t)^{3\alpha-1}$.
The pressure Hessian $\Pi_0$ lower bound in \eqref{eq:outline-two-estimates} permits the  model axial strain $\rWm$ clock-scaling to persist under the 
Riccati law \eqref{eq:outline-Riccati}; in particular, it provides the inequality
\[
\p_t\rW_0(t)\le-\tfrac{1-\upbeta}{2}\rW_0(t)^2.
\]
Thus the pressure Hessian $\Pi_0$ cannot cancel the compressive Riccati term $\tfrac12\,\rW_0(t)^2$.  This strict Riccati imbalance, together
with the Biot--Savart strain analysis, establishes that the true clock-scaling $-\rW_0(t)\simeq\Gamma J(t)^{3\alpha-1}$ agrees with the scaling provided by our clock-and-driver
model.

The Lagrangian geometry is organized so that the true Lagrangian flow map $\phi$ can be compared with the hyperbolic flow map $\Phi_\lin$.  Fundamental to this
comparison, is the diffeomorphism decomposition
\begin{equation} 
\phi=\phi_{\smooth}\circ\phi_{\cusp}. \label{eq:sm-cusp-6}
\end{equation} 
Here, $\phi_{\cusp}$ carries fluid particles that are driven by the velocity field generated by the cusp vorticity, whereas $\phi_{\smooth}$ carries particles
via a smooth velocity field generated by  the smooth decaying tail of the vorticity. This flow decomposition naturally induces the
exact clock decomposition 
\[J(t)=J_{\smooth}(t)J_{\cusp}(t),\]
and the smooth-clock bound
\eqref{eq:Jsmooth-bdd} keeps $J_{\smooth}$ uniformly bounded from above and below.  Hence collapse of the
Euler clock $J(t)$ is equivalent to collapse of the cusp clock $J_{\cusp}$.   We establish a collapse-limit normal form representation for
the cusp flow $\phi_{\cusp}$  in certain collapse coordinates $(\zeta,\tau)$.  To be precise, we show that after a clock-rescaling by $J_{\cusp}^2$, 
the normal form of $\phi_{\cusp}$ is equal to $\Phi_\lin$, modulo lower-order errors:
\[
\phi_{\cusp}(Y_t(\zeta,\tau),t) =J_{\cusp}(t)^2\zeta\bigl((\tau,1)+\mathcal E_t(\zeta,\tau)\bigr) \ \ \text{ in the limit as } \ \ J_{\cusp}\downarrow 0 .
\]
The displacement bound \eqref{eq:normal-form-approximation-bound} proves that the normal form error $\mathcal E_t$ and its relevant
derivatives are small as $J_{\cusp}\downarrow0$.  This is the precise Euler analogue of the placement under $\Phi_{\lin}$.

\subsection{The pressure Hessian model}
\label{sec::Pi-model-explain}
We do not a priori establish the lower bound for $\Pi_0(t)$ directly from the true pressure-equation source
function $\tr((\nabla u)^2)$.  Instead, we first isolate the part of the transported vorticity which has the
separation-of-variables form.  If $(r,z)$ denotes the Eulerian cylindrical coordinate and $J_{\cusp}(t)$ is the
cusp clock, then the variables used in this model are
\begin{equation} 
(\mathcal R,\mathcal Z)=J_{\cusp}(t)^{-2}(r,z). \label{eq:new-coords}
\end{equation} 
We employ a  localized near-axis approximation of the true initial datum \eqref{eq:vort0}.  On the
cone where $\Upsilon=1$ in \eqref{eq:Theta-star-def}, the target angular function $\Theta^*$ is 
$\rho^\alpha\Theta^*(\sigma)=R^\alpha$, with the sign determined by the odd reflection across the equator. To define our approximate initial
vorticity,  we  replace \eqref{eq:vort0} with the initial vorticity function
\[
\Omega_{\theta,0}(R,Z) =-\Gamma\operatorname{sgn}(Z)R^\alpha(1+Z^2)^{-\gamma/2}.
\]
From the true initial condition \eqref{eq:vort0}, we keep the local radial cusp $R^\alpha$, but the spherical decay weight
$(1+\rho^2)^{-\gamma/2}$ is localized to its axial value $(1+Z^2)^{-\gamma/2}$.
Using  \eqref{eq:Jac-Identity} and  \eqref{eq:vort-identity},  we have that   for labels $Z\ne0$,
\[
\begin{aligned}
\Omega_\theta(\phi_r(R,Z,t),\phi_z(R,Z,t),t) &=-\Gamma\operatorname{sgn}(Z)(1+Z^2)^{-\gamma/2}\phi_r(R,Z,t)R^{\alpha-1}.
\end{aligned}
\]
From Lemma~\ref{lem:flow1}, we obtain that 
\[
\phi_r(R,Z,t)=\p_R\phi_r(0,Z,t)R+o(R)\ \ \text{ as } \ \ R\downarrow 0,
\]
from which it follows that 
\[
\Omega_\theta(\phi_r(R,Z,t),\phi_z(R,Z,t),t) \!=\!-\Gamma\operatorname{sgn}(Z)(1+Z^2)^{-\gamma/2}\bigl(\p_R\phi_r(0,Z,t)\bigr)^{\!1-\alpha} 
\phi_r(R,Z,t)^\alpha+o(\phi_r(R,Z,t)^\alpha).
\]
By \eqref{eq:new-coords},  $J_{\cusp}(t)^2\mathcal Z=z=:\phi_z(0,Z,t)$, so that  $\operatorname{sgn}(Z)=\operatorname{sgn}(\mathcal Z)$;
we define now give a ``rough'' definition of the Euler-generated axial function $a_t$ 
\begin{equation} 
a_t(|\mathcal Z|)=\bigl(J_{\cusp}(t)\p_R\phi_r(0,Z,t)\bigr)^{1-\alpha}(1+Z^2)^{-\gamma/2}.  \label{eq:Euler-gen-a}
\end{equation} 
(We shall give the more precise definition below.)

Combining $r=J_{\cusp}(t)^2\mathcal R$ from \eqref{eq:new-coords} with the expansion for
$\Omega_\theta(\phi_r(R,Z,t),\phi_z(R,Z,t),t)$ gives, for fixed $\mathcal Z\ne0$,
\begin{equation}
\Omega_\theta(J_{\cusp}(t)^2\mathcal R,J_{\cusp}(t)^2\mathcal Z,t) =-\Gamma J_{\cusp}(t)^{3\alpha-1}\operatorname{sgn}(\mathcal Z)a_t(|\mathcal Z|)\mathcal R^\alpha 
+o(\Gamma J_{\cusp}(t)^{3\alpha-1}\mathcal R^\alpha)\text{ as } \mathcal R\downarrow0 .
\label{eq:transported-vorticity-leading}
\end{equation}
We shall study the resulting stagnation-point axial strain and pressure Hessian generated by this separation-of-variables vorticity, and to that end,
it  is convenient to remove the amplitude and clock-scaling $\Gamma J_{\cusp}(t)^{3\alpha-1}$ from the vorticity and to keep only the leading order term
$-\operatorname{sgn}(\mathcal Z)a_t(|\mathcal Z|)\mathcal R^\alpha$ from   \eqref{eq:transported-vorticity-leading};  the error created by discarding 
$o(\Gamma J_{\cusp}(t)^{3\alpha-1}\mathcal R^\alpha)$ will be shown to be small in   Section~\ref{sec:transported-cusp-pressure}.
We shall, therefore, perform a detailed analysis of the stagnation-point axial strain and pressure Hessian generated by the following separation-of-variables
vorticity function:
\begin{equation} 
\Omega_{\theta}^{a_t,\infty}(\mathcal R,\mathcal Z)=-\operatorname{sgn}(\mathcal Z)a_t(|\mathcal Z|)\mathcal R^\alpha .
\label{eq:Omega-infi}
\end{equation} 
In our analysis,  the parameter $a_t$ in $\Omega_{\theta}^{a_t,\infty}$  represents the axial function $a_t$ on the right side of the equality in \eqref{eq:Omega-infi}.
In our analysis, we shall use the Euler-generated axial $a_t$ in \eqref{eq:Euler-gen-a}.   We can now also explain the use of the symbol $ \infty$ in $\Omega_{\theta}^{a_t,\infty}$.

For technical reasons, we shall also make use of slope-restricted vorticity. With $\mathcal R/|\mathcal Z|$ denoting the slope, we introduce the slope cutoff function
$\chi_M(\mathcal R/|\mathcal Z|)$, and we define the $M$-slope-restricted vorticity function by 
\begin{equation} 
\Omega_{\theta}^{a_t,M}(\mathcal R,\mathcal Z) =-\operatorname{sgn}(\mathcal Z)a_t(|\mathcal Z|)\mathcal R^\alpha \chi_M\!\left(\mathcal R/|\mathcal Z|\right).
\end{equation} 
The limit as $M \to \infty$ corresponds to setting $\lim_{M\to\infty}\chi_M=1$, and this limit yields the unrestricted full angular vorticity function $\Omega_{\theta}^{a_t,\infty}$ in 
\eqref{eq:Omega-infi}.
The cutoff $\chi_M(\mathcal R/|\mathcal Z|)$ restricts $\supp(\Omega_{\theta}^{a_t,M})$ to the bounded-slope region $\mathcal R/|\mathcal Z|\le2M$; the region
$\supp(1-\chi_M)$  is estimated later by the large-slope tail \eqref{eq:pressure-angular-tail}.

Let us now give a more precise description of the Euler-generated axial function $a_t$.  To do so, we shall make use of the cusp-flow $\phi_{\cusp}$ whose
definition we sketched above  in \eqref{eq:sm-cusp-6} and with precise definition given in Section~\ref{sec:geom-flow-decomp}.
We shall make  use the normalized axial coordinate $\zeta$ defined as follows:
\[
\zeta=J_{\cusp}(t)^{-2}\bigl(\phi_{\cusp}(0,Z,t)\bigr)_z, \qquad Z_t(\zeta)=\left(Z\mapsto J_{\cusp}(t)^{-2}\bigl(\phi_{\cusp}(0,Z,t)\bigr)_z\right)^{-1}(\zeta).
\]
We also define the normalized radial derivative on the axis by
\[
q_t(\zeta):=J_{\cusp}(t)\,\p_R\bigl(\phi_{\cusp}\bigr)_r(0,Z_t(\zeta),t).
\]
We see that  \eqref{eq:Euler-gen-a} can be written in the $\zeta$ coordinate as the physical Euler-generated axial function
\[a_t^{\rm phys}(\zeta)=q_t(\zeta)^{1-\alpha}\bigl(1+Z_t(\zeta)^2\bigr)^{-\gamma/2}. \]
Again for technical reasons,  for an axial function $a$, it is convenient to truncate the physical axial function to the interval $[0,\zeta_a]$. We shall study
the vorticity with the axial function $a_t(\zeta)$ defined by
\[
a_t(\zeta)=a_t^{\rm phys}(\zeta)\mathbf 1_{I_a}(\zeta) =q_t(\zeta)^{1-\alpha}\bigl(1+Z_t(\zeta)^2\bigr)^{-\gamma/2}\mathbf 1_{[0,\zeta_{a_t}]}(\zeta).
\]
This is explained  in detail in  \eqref{eq:euler-generated-truncated-coeff}.  The monotone axial-stretching bootstrap \textup{(BA4)}, namely
\eqref{eq:monotone-axial-two-sided}--\eqref{eq:monotone-axial-fractional-bootstrap}, implies that this
Euler-generated axial function $a_t(\zeta)$ is nonnegative and nonincreasing.

The separation-of-variables in $\Omega_\theta^{a_t,M}$ persists under the Biot--Savart law and the axial strain
reduces to a one-dimensional moment.  In Section~\ref{sec:slope-restricted-pressure},  we prove that 
\begin{equation}
W_M[a_t] =-\int_0^\infty \tfrac{3\tau^{\alpha+2}}{(1+\tau^2)^{5/2}}\chi_M(\tau)\,d\tau \int_0^\infty a_t(\zeta)\zeta^{\alpha-1}\,d\zeta =-C_{\alpha,M}^W I[a_t],
\end{equation}
so that the size of $\tfrac12 W_M[a_t]^2$ is governed by $I[a_t]^2$, the square of a one-dimensional integral.  The axial strain
$W_M[a_t]$ is the axial derivative of the Biot--Savart velocity $U^{a_t,M}=\BS[\Omega_\theta^{a_t,M}e_\theta]$, evaluated at the stagnation point.
The model  pressure Hessian $\Pi_M[a_t]$ is the principal-value integral of the source $\tr((\nabla U^{a_t,M})^2)$ in \eqref{eq:pressure-bilinear-form}.  
With $\chi_M$ replaced by $1$, the full-angular velocity $U^{a_t,\infty}$ generates the  pressure Hessian
\begin{equation}
\Pi_\infty[a_t]:=\pv\int_{\R^3}K_{zz}(Y)\,\tr\bigl(\nabla U^{a_t,\infty}(Y)\nabla U^{a_t,\infty}(Y)\bigr)\,dY .
\label{eq:section6-full-angular-pressure-hessian}
\end{equation}
Just as the separation-of-variables vorticity $\Omega_{\theta}^{a_t,M}$ reduced the axial strain computation to a one-dimensional integral for
the moment function $I[a_t]$, there exists a corresponding reduction of the pressure Hessian $\Pi_\infty[a_t]$ to a one-dimensional computation, and the starting
point for this reduction is in the computation of the axis trace for the three-dimensional Biot--Savart velocity $U^{a_t,\infty}$.  
As we will explain in \eqref{eq:C-alpha-M-limit}, after passing to the full angular problem $M\to\infty$, the strain is
\begin{equation} 
W_\infty[a_t]:=-C_\alpha^W I[a_t], \label{eq:W-infi-Calpha}
\end{equation} 
where $C_\alpha^W:=\lim_{M\to\infty}C_{\alpha,M}^W=\int_0^\infty \tfrac{3\tau^{\alpha+2}}{(1+\tau^2)^{5/2}}\,d\tau>0$.

For $\zeta>0$, we define the one-dimensional  axis-trace velocity by
\[
V_\infty[a_t](\zeta)=U_\zeta^{a_t,\infty}(0,\zeta).
\]
and as we will show in \eqref{eq:full-axis-velocity}, 
\begin{equation} 
V_\infty[a_t](\zeta):=-\tfrac{C_\alpha^W}{2\alpha} \int_0^\infty a_t(\eta) \bigl((\zeta+\eta)^\alpha-|\zeta-\eta|^\alpha\bigr)\,d\eta . \label{eq:V-infi}
\end{equation} 
The key observation is that axis-trace $V_\infty[a_t](\zeta)$ of the BS-velocity is a one-dimensional weighted integral of $a_t(\zeta)$ and that
the axial strain moment integral $I[a_t]=\int_0^\infty a_t(\zeta)\zeta^{\alpha-1}\,d\zeta$ is also a one-dimensional weighted integral of $a_t(\zeta)$, and that
the formula for the pressure Hessian $\Pi_\infty[a_t]$ can be obtained by a type of variational principle applied to a one-parameter family of axial functions $s\mapsto a_{t,s}(\zeta)$
such that $a_{t,0}:=a_t$. Moreover, since $W_\infty[a_t] = \p_ZU_Z^{a_t,\infty}(0,0)$, we also have that 
\[
W_\infty[a_t] = \p_\zeta V_\infty[a_t](0) ,
\]
so that the stagnation-point axial strain is the axial derivative of the one-dimensional velocity field $V_\infty[a_t](\zeta)$ evaluated at 
$\zeta=0$.\footnote{The derivative at $\zeta=0$ is understood as the right derivative of the axis trace, and can be computed directly
from the Biot--Savart kernel.  Indeed, differentiating the explicit formula
\eqref{eq:full-axis-velocity} at $\zeta=0^+$ gives $\p_\zeta V_\infty[a_t](0)=-C_\alpha^W\int_0^\infty a_t(\eta)\eta^{\alpha-1}\,d\eta=W_\infty[a_t]$.}

To formulate this variational principle, we freeze one physical time $t$.  We introduce an auxiliary parameter $s$ and a $C^1$ family of axial 
functions $s\mapsto a_{t,s}$, with $a_{t,0}=a_t$.  This family is used only to define the derivative  $\left.\p_s a_{t,s}\right|_{s=0}$  at the fixed axial 
function $a_t$; no assertion is made that $a_{t,s}$ is an Euler-generated axial function for $s \neq 0$. 

With time $t$ fixed, we set  $V=V_\infty[a_t]$, and we define the one-dimensional Lagrangian flow map $\zeta_s$ by
\begin{equation} 
\tfrac{d}{ds}\zeta_s(\zeta_0)=V(\zeta_s(\zeta_0)), \qquad \zeta_s(\zeta_0)\big|_{s=0}=\zeta_0 .
\label{eq:zeta-s-flow}
\end{equation} 
The specific vorticity $\Omega_\theta/r$,  obtained from the vorticity function in \eqref{eq:transported-vorticity-leading}, is exactly conserved
by the Euler Lagrangian flow map $\phi$.  From \eqref{eq:transported-vorticity-leading}, with $r=J_{\cusp}(t)^2\mathcal R$, $z=J_{\cusp}(t)^2\mathcal Z$, 
and $\mathcal Z=\zeta>0$, we have that $\operatorname{sgn}(\mathcal Z)=1$, and \eqref{eq:transported-vorticity-leading} becomes
\[
\Omega_\theta(r,z,t) =-\Gamma J_{\cusp}(t)^{3\alpha-1}a_t(\zeta)\mathcal R^\alpha +o(\Gamma J_{\cusp}(t)^{3\alpha-1}\mathcal R^\alpha)
\qquad \text{as } \mathcal R\downarrow0 .
\]
Forming the Eulerian specific vorticity $\Omega_\theta(r,z,t)/r$ by dividing  $\Omega_\theta(r,z,t)$ by $r=J_{\cusp}(t)^2\mathcal R$, we obtain that
\[
\tfrac{\Omega_\theta(r,z,t)}{r} =-\Gamma J_{\cusp}(t)^{3\alpha-3}a_t(\zeta)\mathcal R^{\alpha-1} +o(\Gamma J_{\cusp}(t)^{3\alpha-3}\mathcal R^{\alpha-1}) .
\]
We can remove the amplitude and clock scaling as follows:
\begin{equation}
\Gamma^{-1}J_{\cusp}(t)^{3-3\alpha}\tfrac{\Omega_\theta(r,z,t)}{r} =-a_t(\zeta)\mathcal R^{\alpha-1}+o(\mathcal R^{\alpha-1})
\qquad \text{as } \mathcal R\downarrow0 .
\label{eq:normalized-specific-vorticity-leading}
\end{equation}
After the amplitude $\Gamma J_{\cusp}(t)^{3\alpha-1}$ and the lower-order term in \eqref{eq:transported-vorticity-leading} are removed, the auxiliary variation is 
driven by the full-angular Biot--Savart velocity $U^{a_t,\infty}=\BS[\Omega_\theta^{a_t,\infty}e_\theta]$ generated by \eqref{eq:Omega-infi}.  We write
$U_{\mathcal R}(\mathcal R,\zeta,t)$ and $U_\zeta(\mathcal R,\zeta,t)$ for its meridional components; the
dependence on $t$ corresponds  to the frozen  time $t$ in the axial function  $a_t$.  Thus,  the transport operator for the
normalized leading term in \eqref{eq:normalized-specific-vorticity-leading} is
$U_{\mathcal R}\p_{\mathcal R}+U_\zeta\p_\zeta$.  The $\zeta$-component on the axis is the one-dimensional
trace already defined above:
\begin{equation}
U_\zeta(0,\zeta,t)=U^{a_t,\infty}_Z(0,\zeta)=V_\infty[a_t](\zeta)=V(\zeta).
\label{eq:normalized-axis-trace-V}
\end{equation}
The last equality is the convention $V=V_\infty[a_t]$ from \eqref{eq:full-axis-velocity}.  The curve
$\zeta_s$ has already been defined in \eqref{eq:zeta-s-flow} as the Lagrangian flow of this one-dimensional
axis velocity.  We now define the corresponding radial curve $\mathcal R_s$ by using the radial component of
the model transport velocity along the same axial curve:
\begin{equation}
\tfrac{d}{ds}\zeta_s=U_\zeta(0,\zeta_s,t)=V(\zeta_s), \qquad \tfrac{d}{ds}\mathcal R_s=U_{\mathcal R}(\mathcal R_s,\zeta_s,t), \qquad
(\zeta_s,\mathcal R_s)\big|_{s=0}=(\zeta_0,\mathcal R_0).
\label{eq:aux-specific-vorticity-characteristics}
\end{equation}
Note that the $\zeta$-component of the model transport velocity in the characteristic equation is evaluated at
$\mathcal R=0$ because the difference $U_\zeta(\mathcal R,\zeta,t)-U_\zeta(0,\zeta,t)=O(\mathcal R^2)$ is lower order in
\eqref{eq:leading-normalized-transport}.

It remains to identify the leading radial motion in \eqref{eq:aux-specific-vorticity-characteristics}.  By
axisymmetry and smoothness across the symmetry axis, $U_{\mathcal R}$ is odd in $\mathcal R$ and
$U_\zeta$ is even in $\mathcal R$.  Therefore, as $\mathcal R\downarrow0$,
\begin{equation}
U_{\mathcal R}(\mathcal R,\zeta,t)=\mathcal{U}_{\mathcal R}(\zeta)\mathcal R+O(\mathcal R^3), \qquad U_\zeta(\mathcal R,\zeta,t)=V(\zeta)+O(\mathcal R^2).
\label{eq:aux-axis-velocity-expansion}
\end{equation}
The expansion \eqref{eq:normalized-specific-vorticity-leading} is the $s=0$ statement for the Euler-generated
axial function $a_t=a_{t,0}$; it identifies the base leading specific vorticity
$-a_t(\zeta)\mathcal R^{\alpha-1}$.  For $s\neq0$, the function $a_{t,s}$ is not being asserted to come from
the Euler solution.  Instead, after the frozen characteristics $(\mathcal R_s,\zeta_s)$ in
\eqref{eq:aux-specific-vorticity-characteristics} have been specified, we define $a_{t,s}$ by
requiring the auxiliary leading term $-a_{t,s}(\zeta)\mathcal R^{\alpha-1}$ to satisfy a specially chosen  conservation law, written as
\eqref{eq:aux-leading-specific-vorticity-conservation} below.  Before imposing that conservation law, we first explain 
the leading-order calculation obtained by applying
$U_{\mathcal R}\p_{\mathcal R}+U_\zeta\p_\zeta$ to $a_{t,s}(\zeta)\mathcal R^{\alpha-1}$.  By
\eqref{eq:aux-axis-velocity-expansion}, the difference
$U_\zeta(\mathcal R,\zeta,t)-U_\zeta(0,\zeta,t)$ only gives an $o(\mathcal R^{\alpha-1})$ contribution, so
\begin{equation}
(U_{\mathcal R}\p_{\mathcal R}+U_\zeta\p_\zeta)\bigl(a_{t,s}(\zeta)\mathcal R^{\alpha-1}\bigr)  =
U_{\mathcal R}(\mathcal R,\zeta,t)\p_{\mathcal R}\bigl(a_{t,s}(\zeta)\mathcal R^{\alpha-1}\bigr) +U_\zeta(0,\zeta,t)\p_\zeta\bigl(a_{t,s}(\zeta)\mathcal R^{\alpha-1}\bigr)
+o(\mathcal R^{\alpha-1})
\label{eq:leading-normalized-transport}
\end{equation}
as $\mathcal R\downarrow0$.

In \eqref{eq:aux-axis-velocity-expansion},
\begin{equation*}
\mathcal U_{\mathcal R}(\zeta):=\p_{\mathcal R}U_{\mathcal R}(0,\zeta,t).
\end{equation*}
Away from the axis, $\mathcal U_{\mathcal R}(\zeta)\mathcal R$ is only
the leading term in the expansion of $U_{\mathcal R}(\mathcal R,\zeta,t)$ in
\eqref{eq:aux-axis-velocity-expansion}.  The axis value of this radial linearization is determined by
incompressibility.  In normalized cylindrical variables,
\begin{equation}
\tfrac1{\mathcal R}\p_{\mathcal R}\bigl(\mathcal R U_{\mathcal R}\bigr)+\p_\zeta U_\zeta=0.
\label{eq:aux-divergence}
\end{equation}
The only consequence of \eqref{eq:aux-divergence} used here is its axis restriction,
\begin{equation*}
2\mathcal{U}_{\mathcal{R}}(\zeta)+\p_\zeta V(\zeta)=0,
\end{equation*}
which gives
\begin{equation}
\mathcal{U}_{\mathcal{R}}(\zeta)=- \tfrac12\p_\zeta V(\zeta).
\label{eq:axial-div-free}
\end{equation}
Combining \eqref{eq:aux-specific-vorticity-characteristics} with \eqref{eq:aux-axis-velocity-expansion} and
\eqref{eq:axial-div-free} gives the radial characteristic equation
\begin{equation}
\tfrac{d}{ds}\mathcal R_s=\mathcal{U}_{\mathcal{R}}(\zeta_s)\mathcal R_s=-\tfrac12(\partial_\zeta V)(\zeta_s)\mathcal R_s, \qquad \mathcal R_s|_{s=0}=\mathcal R_0.
\label{eq:Rs-evo}
\end{equation}
Consequently,
\begin{equation} 
\tfrac{d}{ds}\mathcal R_s^{\alpha-1}=-\tfrac{\alpha-1}{2}(\partial_\zeta V)(\zeta_s)\mathcal R_s^{\alpha-1}.
\label{eq:ddsRs}
\end{equation} 
We now choose the auxiliary curve $s\mapsto a_{t,s}$ by requiring the leading model specific vorticity to be
conserved along the frozen characteristics $(\mathcal R_s,\zeta_s)$.  Thus, for each initial point
$(\mathcal R_0,\zeta_0)$ and for fixed $s$, we impose
\begin{equation} 
-a_{t,s}(\zeta_s(\zeta_0))\mathcal R_s^{\alpha-1}=-a_t(\zeta_0)\mathcal R_0^{\alpha-1}.
\label{eq:aux-leading-specific-vorticity-conservation}
\end{equation} 
This is a definition of the auxiliary variation, not a claim that \eqref{eq:normalized-specific-vorticity-leading}
holds with $a_t$ replaced by $a_{t,s}$ for $s\neq0$.
The solution to \eqref{eq:Rs-evo} is
$\mathcal R_s=\mathcal R_0 \exp\left(-\tfrac12\int_0^s(\p_\zeta V)(\zeta_\sigma)\,d\sigma\right)$, and so,  for fixed $s$, the power 
$\mathcal R_s^{\alpha-1}$ is $\mathcal R_0^{\alpha-1}$ multiplied by a positive quantity independent of
$\mathcal R_0$.  Dividing \eqref{eq:aux-leading-specific-vorticity-conservation} by
$-\mathcal R_s^{\alpha-1}$ gives
\begin{equation}
a_{t,s}(\zeta_s(\zeta_0)) =a_t(\zeta_0)\left(\tfrac{\mathcal R_0}{\mathcal R_s}\right)^{\!\alpha-1}.
\label{eq:as-axis-limit-transport}
\end{equation}
Differentiating \eqref{eq:as-axis-limit-transport}, we find that
\[
\p_sa_{t,s}(\zeta_s) =-(\alpha-1)\tfrac{1}{\mathcal R_s}\tfrac{d\mathcal R_s}{ds}\,a_{t,s}(\zeta_s)=\tfrac{\alpha-1}{2}(\partial_\zeta V)(\zeta_s)a_{t,s}(\zeta_s),
\]
where we used \eqref{eq:ddsRs}.  It follows that
\begin{equation} 
\p_s a_{t,s}(\zeta_s(\zeta_0)) = -\tfrac{1-\alpha}{2} \bigl(\p_\zeta V\bigr)(\zeta_s(\zeta_0))\, a_{t,s}(\zeta_s(\zeta_0)). \label{eq:evo-as}
\end{equation} 

Using the traditional notation for the first variation, we denote the derivative of the curve $s\mapsto a_{t,s}$ at $s=0$ by
\[
\delta a_t:=\left.\p_s a_{t,s}\right|_{s=0}.
\]
Then from \eqref{eq:evo-as} and \eqref{eq:zeta-s-flow}, we have that
\begin{equation} 
\delta a_t=\left.\p_s a_{t,s}\right|_{s=0}=-V_\infty[a_t]\p_\zeta a_t-\tfrac{1-\alpha}{2}(\p_\zeta V_\infty[a_t])a_t.   \label{eq:delta-at}
\end{equation}

Let us now explain the role of the first variation $\delta a_t$ in the computation of the pressure Hessian $\Pi_\infty[a_t]$.
The curve $s\mapsto a_{t,s}$ is an auxiliary curve of axial functions passing through the one Euler-generated
axial function $a_t=a_{t,0}$; in particular, for $s \neq 0$, it is not necessary that each $a_{t,s}$ is an
Euler-generated axial function.

We now compute  the first variation of the moment integral $I[a_t]$ along the curve $s \mapsto a_{t,s}$:
\[
\mathcal D_\infty[a_t]:=\left.\tfrac{d}{ds}\right|_{s=0}I[a_{t,s}] =\int_0^\infty \delta a_t(\zeta)\zeta^{\alpha-1}\,d\zeta .
\]
Since $W_\infty[a_t]=-C_\alpha^W I[a_t]$, this gives the directional derivative of the axial-strain functional:
\begin{equation} 
\left.\tfrac{d}{ds}\right|_{s=0}W_\infty[a_{t,s}]=-C_\alpha^W\mathcal D_\infty[a_t]. \label{eq:dds-Winfi}
\end{equation} 

Our next objective is to explain why the auxiliary curve $s\mapsto a_{t,s}$ has the same tangent at
$a_t$ as the full-angular Euler evolution started from $a_t$.  In symbols, we want to prove
\begin{equation*}
\left.\p_s a_{t,s}\right|_{s=0}=\left.\p_{t'}b(t')\right|_{t'=t},
\qquad b(t)=a_t.
\end{equation*}

We introduce a new Euler solution in the time variable $t'$.  Let
$\Omega_\theta(\mathcal R,\mathcal Z,t')$ solve the Euler equation with
\begin{equation*}
\left.\Omega_\theta(\mathcal R,\mathcal Z,t')\right|_{t'=t}
=\Omega_\theta^{a_t,\infty}(\mathcal R,\mathcal Z)
=-\operatorname{sgn}(\mathcal Z)a_t(|\mathcal Z|)\mathcal R^\alpha,
\end{equation*}
where the last equality is \eqref{eq:Omega-infi}.  For $t'$ near $t$, define $b(t')$ by
\begin{equation*}
\tfrac{\Omega_\theta(\mathcal R,\zeta,t')}{\mathcal R}
=-b(t')(\zeta)\mathcal R^{\alpha-1}+o(\mathcal R^{\alpha-1})
\qquad\text{as }\mathcal R\downarrow0.
\end{equation*}
Then the initial condition above gives $b(t)=a_t$.
The stagnation-point Riccati identity \eqref{eq:rW0}, applied at the initial instant $t'=t$, yields
\begin{equation}
\left.\tfrac{d}{dt'}\right|_{t'=t}W_\infty[b(t')]
=-\tfrac12W_\infty[a_t]^2-\Pi_\infty[a_t],
\label{eq:section6-full-angular-riccati-from-rW0}
\end{equation}
where $\Pi_\infty[a_t]$ is the pressure Hessian defined in
\eqref{eq:section6-full-angular-pressure-hessian}.

It remains to compare the Euler tangent $\left.\p_{t'}b(t')\right|_{t'=t}$ with the auxiliary tangent
$\left.\p_s a_{t,s}\right|_{s=0}$ already computed in \eqref{eq:delta-at}.
We now compute $\left.\p_{t'}b(t')\right|_{t'=t}$ from the same axis transport law.  At the initial instant,
\begin{equation*}
\left.\tfrac{\Omega_\theta(\mathcal R,\zeta,t')}{\mathcal R}\right|_{t'=t}
=-a_t(\zeta)\mathcal R^{\alpha-1}.
\end{equation*}
The specific-vorticity transport equation for this frozen full-angular Euler evolution gives, at $t'=t$,
\begin{equation*}
-\left.\p_{t'}b(t')\right|_{t'=t}\mathcal R^{\alpha-1}
+U_{\mathcal R}^{a_t,\infty}\p_{\mathcal R}\bigl(-a_t(\zeta)\mathcal R^{\alpha-1}\bigr)
+U_\zeta^{a_t,\infty}\p_\zeta\bigl(-a_t(\zeta)\mathcal R^{\alpha-1}\bigr)
=o(\mathcal R^{\alpha-1})
\end{equation*}
as $\mathcal R\downarrow0$.  The $\mathcal R^{\alpha-1}$ term in this identity is
\begin{equation}
-\left.\p_{t'}b(t')\right|_{t'=t}
-(\alpha-1)\p_{\mathcal R}U_{\mathcal R}^{a_t,\infty}(0,\zeta)a_t(\zeta)
-U_\zeta^{a_t,\infty}(0,\zeta)\p_\zeta a_t(\zeta)=0 .
\label{eq:section6-euler-tangent-axis-balance}
\end{equation}
The two axis traces in \eqref{eq:section6-euler-tangent-axis-balance} are
\begin{equation*}
U_\zeta^{a_t,\infty}(0,\zeta)=V_\infty[a_t](\zeta), \qquad
\p_{\mathcal R}U_{\mathcal R}^{a_t,\infty}(0,\zeta)
=-\tfrac12\p_\zeta V_\infty[a_t](\zeta).
\end{equation*}
The first identity is \eqref{eq:normalized-axis-trace-V}, and the second identity is
\eqref{eq:axial-div-free}.
Therefore
\begin{equation}
\left.\p_{t'}b(t')\right|_{t'=t}
=-V_\infty[a_t]\p_\zeta a_t-\tfrac{1-\alpha}{2}(\p_\zeta V_\infty[a_t])a_t .
\label{eq:section6-euler-tangent}
\end{equation}
The right sides of \eqref{eq:delta-at} and \eqref{eq:section6-euler-tangent} are identical.  Hence
\begin{equation}
\left.\p_s a_{t,s}\right|_{s=0} =\left.\p_{t'}b(t')\right|_{t'=t}
\label{eq:section6-aux-euler-tangent-match}
\end{equation}
The equality \eqref{eq:section6-aux-euler-tangent-match} is used after pairing with the weight $\zeta^{\alpha-1}$ in $W_\infty$.  Because
$a_t=a_t^{\rm phys}\mathbf 1_{[0,\zeta_{a_t}]}$ in \eqref{eq:euler-generated-truncated-coeff}, the distributional derivative
$-\p_\zeta a_t$ includes the endpoint contribution from the indicator cutoff; in Section~\ref{sec:slope-restricted-pressure}, this convention is written explicitly in
\eqref{eq:exact-euler-minus-a-prime-cutoff} and used in \eqref{eq:D-model-def}.

Since
\begin{equation*}
W_\infty[c]=-C_\alpha^W\int_0^\infty c(\zeta)\zeta^{\alpha-1}\,d\zeta
\end{equation*}
for any axial function $c$, the tangent identity \eqref{eq:section6-aux-euler-tangent-match} implies
\begin{equation}
\left.\tfrac{d}{ds}\right|_{s=0}W_\infty[a_{t,s}]
=\left.\tfrac{d}{dt'}\right|_{t'=t}W_\infty[b(t')].
\label{eq:section6-aux-euler-strain-derivative-match}
\end{equation}
We can now see the main simplification.  The pressure Hessian
\eqref{eq:section6-full-angular-pressure-hessian} is a three-dimensional principal-value integral.  However,
we do not estimate that principal-value integral directly.  Combining the Riccati identity
\eqref{eq:section6-full-angular-riccati-from-rW0}, the tangent identity
\eqref{eq:section6-aux-euler-strain-derivative-match}, and the moment derivative
\eqref{eq:dds-Winfi} gives
\begin{equation}
\Pi_\infty[a_t]=C_\alpha^W\mathcal D_\infty[a_t]-\tfrac12W_\infty[a_t]^2.
\label{eq:section6-pressure-hessian-moment-reduction}
\end{equation}
Thus the pressure Hessian is recovered from the derivative of the axial strain along the specially chosen
curve $s\mapsto a_{t,s}$.

The right side of \eqref{eq:section6-pressure-hessian-moment-reduction} is one-dimensional.  Indeed,
\begin{equation}
\mathcal D_\infty[a_t]
=\left\langle-\p_\zeta a_t,\,V_\infty[a_t](\zeta)\zeta^{\alpha-1}\right\rangle
-\tfrac{1-\alpha}{2}\int_0^\infty
\p_\zeta V_\infty[a_t](\zeta)a_t(\zeta)\zeta^{\alpha-1}\,d\zeta ,
\label{eq:section6-Dinfty-one-dimensional}
\end{equation}
where $U_Z^{a_t,\infty}(0,\zeta)=V_\infty[a_t](\zeta)$ by \eqref{eq:V-infi}.  So, while \eqref{eq:section6-full-angular-pressure-hessian} requires the full
three-dimensional gradient $\nabla U^{a_t,\infty}(Y)$ with $Y\in\R^3$, the formula \eqref{eq:section6-Dinfty-one-dimensional} involves only
$a_t$,  $\p_\zeta a_t$, $V_\infty[a_t]$, and $\p_\zeta V_\infty[a_t]$.

By defining the function $F_{a_t}(\zeta)=\int_0^\infty a_t(\eta)\bigl((\zeta+\eta)^\alpha-|\zeta-\eta|^\alpha\bigr)\,d\eta$, 
\eqref{eq:V-infi} shows that $V_\infty[a_t](\zeta)=-\tfrac{C_\alpha^W}{2\alpha}F_{a_t}(\zeta)$.
We set
\begin{equation*}
\mathcal K_1[a_t]=\int_0^\infty F_{a_t}(\zeta)a_t(\zeta)\zeta^{\alpha-2}\,d\zeta,\qquad
\mathcal K_2[a_t]=\left\langle-\p_\zeta a_t,\,F_{a_t}(\zeta)\zeta^{\alpha-1}\right\rangle .
\end{equation*}
As we will establish in  \eqref{eq:D-K-identity} (with $M=\infty$), we have that
\[
\mathcal D_\infty[a_t]
=\tfrac{C_\alpha^W}{4\alpha}
\left((1-\alpha)^2\mathcal K_1[a_t]-(1+\alpha)\mathcal K_2[a_t]\right),
\]
and we will prove that  $\mathcal K_2[a_t]\le\alpha \mathcal K_1[a_t]$ and that  $\alpha I[a_t]^2\le \mathcal K_1[a_t]\le2\alpha I[a_t]^2$,
from which it follows that
\begin{equation}
\mathcal D_\infty[a_t]\ge \tfrac{1-3\alpha}{4}C_\alpha^W I[a_t]^2 .
\label{eq:section6-Dinfty-roadmap-lower}
\end{equation}
so that using \eqref{eq:W-infi-Calpha}, \eqref{eq:section6-pressure-hessian-moment-reduction},  and
\eqref{eq:section6-Dinfty-roadmap-lower}, we obtain the Riccati-competition-inequality
\begin{equation}
\Pi_\infty[a_t]\ge -\tfrac{1+3\alpha}{2}\,\tfrac12 W_\infty[a_t]^2. \label{eq:Riccati-infi}
\end{equation}
By using the geometric bounds obtained for the flow maps $\phi_{\smooth}$ and $\phi_{\cusp}$, we are able to transfer
the bound \eqref{eq:Riccati-infi} to 
\begin{equation}
\Pi_{\cusp}(t)\ge -q_\alpha\tfrac12\mathcal W_{\cusp}(t)^2,\qquad q_\alpha<1,  \label{eq:Riccati-cusp-6}
\end{equation}
which is the bound proved in \eqref{eq:euler-generated-riccati-bound}.  The final comparison between
$\Pi_{\cusp}(t)$ and the true Euler pressure Hessian $\Pi_0(t)$ is perturbative: it estimates the error between
the transported Euler vorticity and the localized model
$-\operatorname{sgn}(Z)a_t(|Z|)R^\alpha\chi_{M_\pressure}(R/|Z|)e_\theta$.

\subsection{From model Riccati to Euler Riccati}
In Sections~\ref{sec:Euler-blowup-for-Theta-star}--\ref{sec:target-profile-typeI-completion}, we will prove that the
true Euler solution tracks closely the clock-and-driver model of Section~\ref{sec:lag-analysis-non-local}.
The first step in establishing this stability is the decomposition of the true Euler Lagrangian flow map
$\phi=\phi_{\smooth}\circ\phi_{\cusp}$ from \eqref{eq:cusp-flow-equation}.  The velocity $U_{\cusp}$
is generated by the vorticity transported by $\phi_{\cusp}$, and $\mathcal W_{\cusp}$ is its axial strain in
\eqref{eq:W-cusp-def}.  The smooth flow has a uniformly bounded Jacobian by \eqref{eq:Jsmooth-bdd}, so the
collapse of $J_{\cusp}$ is equivalent to the collapse of the physical clock.

The second step is geometric.  In the collapse variables $(\zeta,\tau)$, the exact cusp flow has the normal form
\[
\phi_{\cusp}(Y_t(\zeta,\tau),t) =J_{\cusp}(t)^2\zeta\bigl((\tau,1)+\mathcal E_t(\zeta,\tau)\bigr),
\]
where the error $\mathcal E_t$ is lower order as $J_{\cusp}\downarrow0$, so to leading order, $\phi_{\cusp}$ matches $\Phi_\lin$.
It is for this reason that true Euler tracks our clock-and-driver model.

The third step transfers the Riccati-competition-inequality  \eqref{eq:Riccati-cusp-6} to the true Euler solution. The proof of this uses the localized vorticity
\eqref{eq:Omega-sharp-cutoff-def}, the normal-form displacement bound \eqref{eq:normal-form-approximation-bound},
and the tail estimates \eqref{eq:pressure-zeta-tail} and \eqref{eq:pressure-angular-tail}. We prove that, for some $\upbeta <1$,
$\Pi_0(t)\ge-\upbeta\,\tfrac12\,\rW_0(t)^2$, so that with the Riccati equation  \eqref{eq:outline-Riccati}, we obtain that 
\[
\p_t\rW_0(t)\le-\tfrac{1-\upbeta}{2}\rW_0(t)^2.
\]
The axial strain estimates \eqref{eq:collapse-proof-rWcusp-scale}--\eqref{eq:collapse-proof-total-strain-negative}
and the clock identity \eqref{eq:outline-J-evo} imply
\[
c\Gamma(T^*-t)\le J(t)^{1-3\alpha}\le C\Gamma(T^*-t).
\]
Since $1-3\alpha>0$, applying the increasing map $x\mapsto x^{1/(1-3\alpha)}$ gives the clock law
\eqref{eq:target-proof-physical-clock-typeI}.  The same positive-power estimate gives
$\Gamma J(t)^{3\alpha-1}\simeq(T^*-t)^{-1}$, and the two-sided $L^\infty$ vorticity bounds in
Lemma~\ref{lem:target-vorticity-envelope} give \eqref{eq:target-proof-vorticity-typeI}.

\subsection{Open-set stability}
Section~\ref{sec:proof-main} proves Theorem~\ref{thm:main}.  The admissible perturbations have
\[
\Theta(\sigma)=\Theta^*(\sigma)(1+h(\sigma)),\qquad h(\sigma)=(\sin\sigma)^\eta k(\sigma),\qquad \|k\|_{C^\alpha([0,\pi/2])}<\nu.
\]
In the strain-producing sector near the nodal angle where $K_W$ in \eqref{eq:KW-kernel} attains its maximum, the corresponding initial Lagrangian angles are $O(J_{\cusp}^3)$.  Therefore the weight
$(\sin\sigma)^\eta$ gives the perturbative multiplier $\nu J_{\cusp}^{3\eta}$ in
\eqref{eq:localized-vort-pert}.  This is why the Target Profile pressure comparison persists under small
weighted H\"older perturbations.

With $\Theta^\nu=\Theta^*+\Theta^*h$, the term $\Theta^*$ provides the Target Profile Riccati comparison, while the term
$\Theta^*h$ contains the weight $\nu J_{\cusp}^{3\eta}$ making the perturbation lower order as shown in \eqref{eq:localized-vort-pert}.
Thus the one-dimensional pressure reduction we described above yields the perturbation bounds
\[
|\mathcal W_{\cusp}^{\nu}(t)-\mathcal W_{\cusp}^{*}(t)| \le\varepsilon_W|\mathcal W_{\cusp}^{*}(t)|,\qquad
\Pi_{\cusp}^{\nu}(t)\ge -q_{\rm tr}^{\rm pert}\tfrac12\bigl(\mathcal W_{\cusp}^{\nu}(t)\bigr)^2,
\]
which are \eqref{eq:pert-cusp-strain-comparison} and \eqref{eq:pert-cusp-pressure-win}.  These estimates permit the
Riccati transfer to the true Euler solution for any datum in the admissible class $\mathcal A_{\alpha,\gamma}(\nu,\eta)$.

\section{A hyperbolic clock-and-driver model with finite-time blowup}
\label{sec:lag-analysis-non-local}

We now devise an idealized model clock-and-driver system associated with the Target Profile datum $\Theta^*$ from \eqref{eq:Theta-star-def}.
Our Lagrangian model  employs a simple linear hyperbolic flow map
\[
\Phi_{\lin}(R,Z,t)=\bigl(R\Jm(t)^{-1},Z\Jm(t)^2\bigr),
\]
which captures the essential geometry of cusp-type incompressible blowup.  The model clock $\Jm(t)$ dynamics are kinematic and fixed, while
the model axial strain $\rWm(t)$ dynamics keep the same Biot--Savart strain integral of the true Euler axial strain.  The resulting 2-component system
reduces to a rather remarkable ODE for the  collapse dynamics of the  model clock $\Jm(t)$.

The main feature of this simple model system is the identification of  the correct scaling laws.  The linear  model provides  the correct
axial strain scaling  $\rWm(t)\simeq-\Gamma\Jm(t)^{3\alpha-1}$ and hence, we obtain the correct Type--I clock scaling
$\dot\Jm(t)\simeq-\Gamma\Jm(t)^{3\alpha}$.  This immediately identifies $\alpha=\tfrac13$ as the collapse threshold: the clock reaches zero in
finite time precisely in the subcritical range $\alpha<\tfrac13$.

For the true Euler solution, the cusp vorticity is transported by the true Euler cusp flow $\phi_{\cusp}$ (introduced in
Section~\ref{sec:Euler-blowup-for-Theta-star}).  We will show that the normal form of the true Euler cusp flow, in the collapse-coordinates,
precisely tracks the linear hyperbolic model flow; as such, our model serves as a calibration: it isolates the clock scale, the axial strain scale, and the
geometric explanation of the $J^{3\alpha}$ depletion caused by kinematic angular drift.

\subsection{The clock-and-driver model system}

\subsubsection{Euler-like structure of the clock-and-driver model system}

We let $\phi(\cdot,t)$ denote the true Euler Lagrangian flow map and we define the Lagrangian toroidal vorticity function by $V(Y,t):=\omega_\theta(\phi(Y,t),t)$
where $Y=(R,Z)$. The meridional Jacobian determinant at the stagnation-point is $J(t):=\det\nabla_{(R,Z)}(\phi_r,\phi_z)(0,0,t)$,
and satisfies the exact kinematic identity
\begin{equation}
\dot J(t)=\tfrac12\,J(t)\,\rW_0(t), \qquad \rW_0(t):=\p_z u_z(0,0,t).
\label{eq:lin-deriv:J0-rW0}
\end{equation}
Moreover, in the axisymmetric no-swirl class, recall that $D_t \xi =0$,  which implies that along trajectories
\begin{equation}
V(Y,t) =\omega_\theta(\phi(Y,t),t) =\tfrac{\phi_r(Y,t)}{R}\,\omega_{\theta,0}(Y),
\label{eq:lin-deriv:V-stretch}
\end{equation}
where, for $R>0$, $\tfrac{\phi_r(X,t)}{R}$ is the radial stretching function, with the axis value understood by continuity.  In the hyperbolic collapse geometry below,
this function is modeled by $\Jm(t)^{-1}$.

As noted above, the true Euler axial strain at the stagnation point is defined by
\begin{equation}
\rW_0(t) =\int_{\R^3}K_W(y)\,\omega_\theta(y,t)\,\ud y,
\label{eq:lin-deriv:rW0-eulerian}
\end{equation}
where $K_W$ is the (explicit) stagnation-point strain kernel (cf.\ \S\ref{sec:KW-kernel}). Using the incompressibility of Euler, so that $\det\nabla\phi(\cdot,t)\equiv 1$,
we may rewrite \eqref{eq:lin-deriv:rW0-eulerian} in Lagrangian variables;  by  the change-of-variables theorem, using  $y=\phi(Y,t)$, we have that
\begin{equation}
\rW_0(t) =\int_{\R^3}K_W(\phi(Y,t))\,V(Y,t)\,\ud Y.
\label{eq:lin-deriv:rW0-lagrangian}
\end{equation}
Equations \eqref{eq:lin-deriv:J0-rW0}--\eqref{eq:lin-deriv:rW0-lagrangian}, together with the stretching identity \eqref{eq:lin-deriv:V-stretch}, form a closed
stagnation-point clock-and-driver system: the clock $J_0$ is driven by $\rW_0$, while $\rW_0$ is a nonlocal  BS integral of the advected vorticity,
which itself is amplified by $J_\twoD(R,Z,t)^{-1}$ via \eqref{eq:vort-identity}.

\subsubsection{The linear hyperbolic flow map $\Phi_{\lin}$}
For given  scalar clock function $\Jm(t)$,  we define
\begin{equation*}
\Phi_{\lin}(R,Z,t)=\big(R\Jm(t)^{-1},\,Z\Jm(t)^2\big),
\end{equation*}
with no dependence on the azimuthal angle.  In cylindrical Eulerian variables $(r,z)$,  the volume form is $r\,\ud r\,\ud z\,\ud\theta$. Changing variables
using the flow map $\Phi_\lin$, so that   $r=\tfrac{R}{\Jm}$ and $z=Z\Jm^2$,  we have that
$r\,\ud r\,\ud z\,\ud\theta =\tfrac{R}{\Jm}\cdot \tfrac{\ud R}{\Jm}\cdot \Jm^2\,\ud Z\,\ud\theta  =R\,\ud R\,\ud Z\,\ud\theta$,
so that  $\Phi_{\lin}$ is volume-preserving in $\mathbb{R}^3$, while its \emph{meridional} Jacobian equals $\Jm(t)$.

\subsubsection{The model vorticity, axial strain, and clock}
We replace the Euler stretching ratio in \eqref{eq:lin-deriv:V-stretch} by $\Jm(t)^{-1}$ and define the model
Lagrangian vorticity by $V_{\mm}(Y,t):=\Jm(t)^{-1}\,\omega_{\theta,0}(Y)$.
The Eulerian model vorticity is its push-forward by $\Phi_\lin$:
\[
\omega_{\mm,\theta}(\Phi_{\lin}(Y,t),t)=V_{\mm}(Y,t)
\qquad\Longleftrightarrow\qquad \omega_{\mm,\theta}(r,z,t)=\Jm(t)^{-1}\,\omega_{\theta,0}(\Jm(t)r,z\Jm(t)^{-2}).
\]
The model axial strain is defined by the corresponding Biot--Savart strain integral:
\begin{equation}
\rWm(t):=\int_{\R^3}K_W(y)\,\omega_{\mm,\theta}(y,t)\,\ud y
=\int_{\R^3}K_W(\Phi_{\lin}(Y,t))\,V_{\mm}(Y,t)\,\ud Y. \notag
\end{equation}
The model clock is then closed by imposing the kinematic clock law from \eqref{eq:lin-deriv:J0-rW0} so that
$\dot{\Jm}(t)=\tfrac12\,\Jm(t)\,\rWm(t)$ with initial condition $ \Jm(0)=1$.

We recall  that $\omega_{\theta,0}$ denotes the Target Profile datum \eqref{eq:vort0} with $\Theta=\Theta^*$.
\begin{definition}[Clock-and-driver model]
\label{def:fields-sec7}
For a positive clock $\Jm(t)$, the clock-and-driver model consists of the following objects.
\begin{enumerate}
\item \textsf{Hyperbolic flow map:} The meridional map is
\begin{equation*}
\Phi_{\lin}(R,Z,t)=\big(R\Jm(t)^{-1},Z\Jm(t)^2\big).
\end{equation*}

\item \textsf{Model vorticity:} The Eulerian model vorticity is
\begin{equation}
\bs{\omega}_{\mm}(r,z,t) =\Jm(t)^{-1}\omega_{\theta,0}(\Jm(t)r,z\Jm(t)^{-2})\bs e_\theta.
\label{eq:vort-model-def}
\end{equation}
In the upper half-space, we define
\begin{equation*}
\rho_{\Lag}(r,z,t):=\sqrt{\Jm(t)^2r^2+\Jm(t)^{-4}z^2}, \qquad \sigma_{\Lag}(r,z,t):=\arctan\!\big(\Jm(t)^3\tan\sigma(r,z)\big),
\end{equation*}
and we set
\begin{equation}
\mathcal F(s):=(1+s^2)^{-\gamma/2}.
\label{eq:radial-tail-F-def}
\end{equation}
Substituting the Target Profile datum \eqref{eq:vort0} into \eqref{eq:vort-model-def} gives
\begin{align}
\bs{\omega}_{\mm}(r,z,t)
&=-\Gamma\,\Jm(t)^{-1}\rho_{\Lag}(r,z,t)^\alpha\,
\mathcal F\!\big(\rho_{\Lag}(r,z,t)\big)
\Theta^*\!\big(\sigma_{\Lag}(r,z,t)\big)\,\bs e_\theta
\notag \\
&=-\Gamma\,\Jm(t)^{-1}(\Jm(t)r)^\alpha\,
\mathcal F\!\big(\rho_{\Lag}(r,z,t)\big)
\Upsilon\!\big(\sigma_{\Lag}(r,z,t)\big)\,\bs e_\theta,
\label{eq:vort-model-explicit}
\end{align}
where the second identity uses
$\rho_{\Lag}\sin\sigma_{\Lag}=\Jm(t)r$ and
$\Theta^*(\sigma)=(\sin\sigma)^\alpha\Upsilon(\sigma)$.

\item \textsf{Model velocity:} The velocity $u_{\mm}$ is the decaying velocity generated by
$\bs{\omega}_{\mm}$ through the Biot--Savart law:
\begin{equation}
u_{\mm}(x,t) =\nabla\times(-\Delta)^{-1}\bs{\omega}_{\mm}(x,t)
=-\tfrac{1}{4\pi}\int_{\R^3}\tfrac{(x-y)\times\bs{\omega}_{\mm}(y,t)}{|x-y|^3}\,\ud y.
\label{eq:u-model-def}
\end{equation}

\item \textsf{Model axial strain:} The model axial strain is
\begin{equation}
\rWm(t):=\p_z(u_{\mm})_z(0,0,t) =\int_{\R^3}K_W(y)\,\omega_{\mm,\theta}(y,t)\,\ud y.
\label{eq:model-strain-def}
\end{equation}

\item \textsf{Model clock evolution:} The model clock satisfies
\begin{equation}
\dot{\Jm}(t)=\tfrac12\,\Jm(t)\rWm(t),\qquad \Jm(0)=1.
\label{eq:J-model-ODE}
\end{equation}
\end{enumerate}
\end{definition}

\subsection{Kinematics: The Angular Drift}
The anisotropic scaling in $\Phi_\lin$ expands the radial variable and compresses the axial variable.  Therefore
a label in the upper half-space drifts in Eulerian polar angle away from the symmetry axis and toward the
equatorial plane.

\begin{lemma}[Angular Drift Law]
\label{lem:drift}
Let $R_0\ge0$ and $Z_0>0$, and let $\sigma_0\in[0,\tfrac\pi2)$ be the Lagrangian polar angle defined by
$\tan\sigma_0=R_0/Z_0$.  If $(r(t),z(t))=\Phi_{\lin}(R_0,Z_0,t)$,  then the Eulerian polar angle $\sigma(t)\in[0,\tfrac\pi2)$ satisfies
\begin{equation}
\tan\sigma(t)=\Jm(t)^{-3}\tan\sigma_0.
\label{eq:forward-angular-drift}
\end{equation}
Conversely, for each fixed Eulerian angle $\sigma\in(0,\tfrac\pi2)$, the pullback Lagrangian angle is
\begin{equation}
\sigma_{\Lag}=\arctan\!\big(\Jm(t)^3\tan\sigma\big) =\Jm(t)^3\tan\sigma+O\!\big(\Jm(t)^9\tan^3\sigma\big) \qquad\text{as }\Jm(t)\downarrow0.
\label{eq:lag-scaling}
\end{equation}
\end{lemma}
\begin{proof}[Proof of Lemma \ref{lem:drift}]
We set $J:=\Jm(t)$.  The linear map gives
\[
r(t)=R_0J^{-1},\qquad z(t)=Z_0J^2.
\]
Since $Z_0>0$, the tangent of the Eulerian polar angle is the ratio of radius to height:
\begin{equation*}
\tan\sigma(t) =\tfrac{r(t)}{z(t)} =\tfrac{R_0J^{-1}}{Z_0J^2} =J^{-3}\tfrac{R_0}{Z_0}.
\end{equation*}
Using $\tan\sigma_0=R_0/Z_0$ proves \eqref{eq:forward-angular-drift}.  Solving
\eqref{eq:forward-angular-drift} for the Lagrangian angle associated with a fixed Eulerian angle
$\sigma\in(0,\tfrac\pi2)$ shows that
\begin{equation*}
\sigma_{\Lag}=\arctan\!\big(J^3\tan\sigma\big).
\end{equation*}
The expansion $\arctan s=s+O(s^3)$ as $s\to0$ yields \eqref{eq:lag-scaling}.
\end{proof}

\begin{remark}[Drift, transported cusp, and isotropic tail]
\label{rem:transported-cusp}
Lemma~\ref{lem:drift} identifies the two geometric effects produced by the linear hyperbolic map.  For a fixed
Eulerian point $x=(r,z)$ in the upper half-space with $z>0$, and with $J=\Jm(t)$, its pullback label is
$(Jr,J^{-2}z)$, so
\[
\sigma_{\Lag} =J^3\tan\sigma+O(J^9\tan^3\sigma), \qquad \rho_{\Lag} =J^{-2}z\big(1+O(J^6\tan^2\sigma)\big) \ \ \text{ as } \ \ J\downarrow 0
.
\]
Thus, the pullback label lies closer to the symmetry axis and farther out in spherical radius.  Since the Target
Profile is $\Theta^*(\sigma)=(\sin\sigma)^\alpha\Upsilon(\sigma)$, the exact identity
\[
\rho_{\Lag}^\alpha\Theta^*(\sigma_{\Lag}) =(\rho_{\Lag}\sin\sigma_{\Lag})^\alpha\Upsilon(\sigma_{\Lag}) =(Jr)^\alpha\Upsilon(\sigma_{\Lag})
\]
converts the spherical cusp into the cylindrical power $(Jr)^\alpha$.  At the same time, the algebraic tail is
evaluated at the enlarged radius $\rho_{\Lag}$.  This is the geometric identity used below in the computation of
the model axial strain.
\end{remark}

\subsection{Gradient bounds for the model velocity}
To determine the size of the model strain we need pointwise bounds on $\nabla u_\mm$, where $u_\mm$ is given
by \eqref{eq:u-model-def}. In the following, constants depend on the H\"older exponent
$\alpha\in(0,\tfrac13)$, the tail exponent $\gamma$, and on the fixed cone angle $\sigma_*\in(0,\tfrac\pi2)$.

\begin{lemma}[Gradient bounds for the model velocity]
\label{lem:grad-umodel-bound}
Assume $\gamma>\alpha$.  Fix $\sigma_*\in(0,\tfrac\pi2)$ and define
\[
\mathcal C_{\sigma_*}:=\{x\in\R^3\setminus\{0\}:\ 0\le\sigma(x)\le\sigma_*\}.
\]
There exists a constant $C=C(\alpha,\gamma,\sigma_*)$ such that, for every $t$ with $\Jm(t)\in(0,1]$ and every
$x\in\mathcal C_{\sigma_*}$, the model velocity gradient satisfies
\begin{equation}
|\nabla u_\mm(x,t)| \le C\,\Gamma\,\Jm(t)^{3\alpha-1}.
\label{eq:grad-umodel-lemmabound}
\end{equation}
\end{lemma}

\begin{proof}[Proof of Lemma \ref{lem:grad-umodel-bound}]
We fix $t$ and write $J:=\Jm(t)\in(0,1]$.  By the odd/even symmetry across the equatorial plane it suffices to
estimate points in the upper half-space; in the folded-angle notation this means that one may replace $z$ by
$|z|$ in the cone lower bounds below.  Let $x\in\mathcal C_{\sigma_*}$ and set $\rho_x:=|x|$. As in the
standard Calder\'on--Zygmund decomposition for the Biot--Savart law,
\begin{equation*}
\nabla u_\mm(x,t) = \pv \int_{\R^3} \nabla_x K(x-y)\,\bs{\omega}_\mm(y,t)\,\ud y + \mathbf C\,\bs{\omega}_\mm(x,t),
\end{equation*}
where $K(z)=-z\times\cdot/(4\pi|z|^3)$ is the Biot--Savart kernel and $\mathbf C$ denotes the local
Calder\'on--Zygmund contraction operator. We estimate the contraction term and the principal-value part
separately.

\runinhead{Step 1: Pointwise cone bound for the model vorticity.}
Because $x\in\mathcal C_{\sigma_*}$, we have $|z(x)|\ge c_*\rho_x$ with $c_*:=\cos\sigma_*>0$. Hence,
\begin{equation}
\rho_{\Lag}(x) =\sqrt{J^2r(x)^2+J^{-4}z(x)^2} \ge J^{-2}|z(x)| \ge c_*J^{-2}\rho_x.
\label{eq:rhoLag-cone-lb}
\end{equation}
Using \eqref{eq:vort-model-explicit}, $0\le \Upsilon\le1$, and the boundedness of $s\mapsto s^\alpha\mathcal
F(s)$ on $[0,\infty)$ (valid whenever $\gamma>\alpha$), we obtain
\begin{align*}
|\bs\omega_\mm(x,t)|
&\le \Gamma J^{-1}(Jr(x))^\alpha\mathcal F\big(\rho_{\Lag}(x)\big)
\le C(\sigma_*)\,\Gamma J^{\alpha-1}\rho_x^\alpha\,\mathcal F(c_*J^{-2}\rho_x)
\notag\\
&=C(\sigma_*)\,\Gamma J^{3\alpha-1}
\Big((J^{-2}\rho_x)^\alpha\,\mathcal F(c_*J^{-2}\rho_x)\Big)
\le C(\sigma_*)\,\Gamma J^{3\alpha-1}.
\end{align*}
Therefore the local contraction term satisfies
\begin{equation}
|\mathbf C\,\bs\omega_\mm(x,t)|\le C(\sigma_*)\,\Gamma J^{3\alpha-1}.
\label{eq:local-contraction-bound-sec7}
\end{equation}

\runinhead{Step 2: Principal value term: local contribution.}
We choose
\[
\sigma_+:=\tfrac12\Big(\sigma_*+\tfrac\pi2\Big)\in(\sigma_*,\tfrac\pi2), \qquad \delta:=\sin(\sigma_+-\sigma_*)\in(0,1).
\]
If $|x-y|\le \delta|x|$, then $\sigma(y)\le \sigma_+$, and therefore $|z(y)|\ge c_+|y|$ with $c_+:=\cos\sigma_+>0$.  We also set
\[
c_1:=(1-\delta)c_+>0.
\]
We decompose space as
\[
D_1:=\{y:\ |x-y|\le \delta|x|\}, \qquad D_2:=\{y:\ \delta|x|<|x-y|\le 2|x|\}, \qquad D_3:=\{y:\ |x-y|>2|x|\}.
\]
On $D_1$,  we have that $|y|\ge(1-\delta)|x|$ and, exactly as in \eqref{eq:rhoLag-cone-lb},
\begin{equation}
\rho_{\Lag}(y)\ge c_+J^{-2}|y|\ge c_1J^{-2}|x|, \qquad y\in D_1.
\label{eq:rhoLag-D1-lb}
\end{equation}
We use the mean-zero property of $\nabla K$ on spheres and write
\[
\mathcal I_{\loc} :=\int_{D_1}\nabla K(x-y)\big(\bs\omega_\mm(y,t)-\bs\omega_\mm(x,t)\big)\,\ud y.
\]
To bound the difference, we estimate the local $C^\alpha$ seminorm of $\bs\omega_\mm$ on $D_1$. We write
\[
\omega_{\mm,\theta}(y,t) =-\Gamma J^{-1}\,A(y)\,B(y)\,C(y),
\]
where
\[
A(y):=(Jr(y))^\alpha, \qquad B(y):=\mathcal F\big(\rho_{\Lag}(y)\big), \qquad C(y):=\Upsilon\big(\sigma_{\Lag}(y)\big).
\]
We estimate the three terms separately.

\rruninhead{Step 2a: the cusp term.}
The scalar estimate for $(Jr)^\alpha$ is not by itself enough across the symmetry axis, because
$\bs e_\theta$ is singular there.  We use the toroidal axis estimate
Lemma~\ref{lem:toroidal-vector-holder} in the scaled form
\begin{equation}
[(Jr)^\alpha\bs e_\theta]_{C^\alpha(D_1)}\le C J^\alpha, \qquad \|(Jr)^\alpha\bs e_\theta\|_{L^\infty(D_1)}\le C J^\alpha|x|^\alpha.
\label{eq:A-vector-bounds-sec7}
\end{equation}
Away from the axis this is the usual scalar H\"older estimate for $(Jr)^\alpha$ multiplied by the smooth
basis vector; Lemma~\ref{lem:toroidal-vector-holder} supplies \eqref{eq:A-vector-bounds-sec7} uniformly when
$D_1$ intersects the axis.

\rruninhead{Step 2b: the spherical tail term.}
By \eqref{eq:rhoLag-D1-lb} and the boundedness of $s^\alpha\mathcal F(s)$,
\[
\|B\|_{L^\infty(D_1)} \le C\,J^{2\alpha}|x|^{-\alpha}.
\]
Moreover, $\mathcal F'(s)=-\gamma s(1+s^2)^{-\gamma/2-1}$, so
$|\mathcal F'(s)|\le \gamma s^{-1}\mathcal F(s)$ for all $s>0$.  Since
$|\nabla\rho_{\Lag}|\le J^{-2}$ and $\mathcal F$ is nonincreasing, \eqref{eq:rhoLag-D1-lb} gives
\[
|\nabla B(y)|=|\mathcal F'(\rho_{\Lag}(y))|\,|\nabla\rho_{\Lag}(y)|\le \gamma c_1^{-1}|x|^{-1}\mathcal F(c_1J^{-2}|x|),\qquad y\in D_1.
\]
Hence
\[
[B]_{C^\alpha(D_1)} \le C\,|x|^{-\alpha}\,\mathcal F(c_1J^{-2}|x|) \le C\,J^{2\alpha}|x|^{-2\alpha}.
\]

\rruninhead{Step 2c: the transported angular cutoff.}
On $D_1$, the fixed-cone geometry gives $|\nabla\sigma(y)|\le(1-\delta)^{-1}|x|^{-1}$.  Since
$\sigma_{\Lag}(y)=\arctan(J^3\tan\sigma(y))$ and $\sigma(y)\le\sigma_+$ on $D_1$, we have
\[
|\nabla\sigma_{\Lag}(y)|\le J^3\sec^2\sigma_+\,(1-\delta)^{-1}|x|^{-1},\qquad y\in D_1.
\]
The Lipschitz regularity of $\Upsilon$ then gives
\begin{equation}
[C]_{C^\alpha(D_1)} \le C\,J^{3\alpha}|x|^{-\alpha}, \qquad \|C\|_{L^\infty(D_1)}\le 1.
\label{eq:C-holder-bound-sec7}
\end{equation}

Combining the vector estimate \eqref{eq:A-vector-bounds-sec7} with
\eqref{eq:C-holder-bound-sec7} and the corresponding bounds for $B$, using the product inequality
$[fg]_{C^\alpha}\le \|f\|_{L^\infty}[g]_{C^\alpha}+[f]_{C^\alpha}\|g\|_{L^\infty}$, we obtain
\begin{equation*}
[\bs\omega_\mm(\cdot,t)]_{C^\alpha(D_1)} \le C(\sigma_+)\,\Gamma\,J^{3\alpha-1}|x|^{-\alpha}.
\end{equation*}
Therefore,
\begin{align}
|\mathcal I_{\loc}|
&\le C\,[\bs\omega_\mm(\cdot,t)]_{C^\alpha(D_1)}
\int_0^{\delta|x|} r^{-3}\,r^\alpha\,r^2\,\ud r
\le C(\sigma_*)\,\Gamma\,J^{3\alpha-1}|x|^{-\alpha}\cdot |x|^\alpha
\le C(\sigma_*)\,\Gamma\,J^{3\alpha-1}.
\label{eq:Iloc-cone-spherical}
\end{align}

\runinhead{Step 3: Principal value term: intermediate annulus.}
For arbitrary $y$ in the upper half-space, we write $y=(\rho,\sigma,\varphi)$ and define
\begin{equation}
B_J(\sigma):=\sqrt{J^2\sin^2\sigma+J^{-4}\cos^2\sigma}, \qquad \rho_{\Lag}(y)=\rho\,B_J(\sigma).
\label{eq:BJ-def-sec7}
\end{equation}
Using \eqref{eq:vort-model-explicit}, the identity
$\sin^\alpha\sigma=J^{2\alpha}(\tan\sigma)^\alpha(1+J^6\tan^2\sigma)^{-\alpha/2}B_J(\sigma)^\alpha$
and the boundedness of $s^\alpha\mathcal F(s)$, we obtain the global weighted estimate
\begin{equation}
|\bs\omega_\mm(y,t)| \le C\,\Gamma\,J^{3\alpha-1}(\tan\sigma(y))^\alpha, \qquad y\in\R^3_+,
\label{eq:omega-global-weighted}
\end{equation}
where $\R^3_+=\{z>0\}$ and the lower half-space contributes symmetrically by odd reflection. On $D_2$ we have
$|x-y|\ge\delta|x|$, hence $|\nabla K(x-y)|\le C\delta^{-3}|x|^{-3}$. Therefore, using
\eqref{eq:omega-global-weighted}, axisymmetry, and the integrability of $(\tan\sigma)^\alpha\sin\sigma$ on
$[0,\tfrac\pi2]$, we obtain
\begin{align}
\Big|\int_{D_2}\nabla K(x-y)\,\bs\omega_\mm(y,t)\,\ud y\Big| &\le C|x|^{-3}\int_{|y|\le 3|x|}|\bs\omega_\mm(y,t)|\,\ud y
\notag\\
&\le C\Gamma\,J^{3\alpha-1}|x|^{-3}\!\! \int_0^{3|x|}\!\!\!\rho^2  \! \ud\rho \int_0^{\frac{\pi}{2}}(\tan\sigma)^\alpha\sin\sigma\,\ud\sigma
\le C\,\Gamma\,J^{3\alpha-1}.
\label{eq:D2-contribution-spherical}
\end{align}

\runinhead{Step 4: Principal value term: far field.}
For $D_3$ we use the stronger version of \eqref{eq:omega-global-weighted}, namely
\[
|\bs\omega_\mm(y,t)| \le C\,\Gamma\,J^{3\alpha-1}(\tan\sigma(y))^\alpha \big(\rho_{\Lag}(y)\big)^\alpha\,\mathcal F\big(\rho_{\Lag}(y)\big).
\]
Since $|x-y|>2|x|$ on $D_3$, we have $|x-y|\ge \tfrac23|y|$ and hence $|\nabla K(x-y)|\le C|y|^{-3}$.  Using spherical coordinates in the upper half-space,
\eqref{eq:BJ-def-sec7}, and then the change of variables $s=\rho B_J(\sigma)$ in the radial integral, we obtain
\begin{align}
\Big|\int_{D_3}\!\!\!\!\nabla K(x-y)\,\bs\omega_\mm(y,t)\,\ud y\Big|
&\le C\,\Gamma\,J^{3\alpha-1}\int_0^{\frac{\pi}{2}}(\tan\sigma)^\alpha\sin\sigma
\int_{|x|}^{\infty}\rho^{-1}\big(\rho B_J(\sigma)\big)^\alpha \mathcal F\big(\rho B_J(\sigma)\big) \ud\rho\,\ud\sigma
\notag\\
&= C\,\Gamma\,J^{3\alpha-1}\int_0^{\frac{\pi}{2}}(\tan\sigma)^\alpha\sin\sigma
\int_{|x|B_J(\sigma)}^{\infty}s^{\alpha-1}\mathcal F(s)\,\ud s\,\ud\sigma
\notag\\
&\le C\,\Gamma\,J^{3\alpha-1}\int_0^{\frac{\pi}{2}}(\tan\sigma)^\alpha\sin\sigma\,\ud\sigma
\int_0^{\infty}s^{\alpha-1}\mathcal F(s)\,\ud s
\le C\,\Gamma\,J^{3\alpha-1}.
\label{eq:D3-contribution-spherical}
\end{align}
The last step uses $\gamma>\alpha$, which makes the radial integral finite.

\runinhead{Step 5: Conclusion.} Combining \eqref{eq:local-contraction-bound-sec7},
\eqref{eq:Iloc-cone-spherical}, \eqref{eq:D2-contribution-spherical}, and
\eqref{eq:D3-contribution-spherical}, we obtain that 
\[
|\nabla u_\mm(x,t)|\le C(\sigma_*)\,\Gamma\,J^{3\alpha-1},
\]
which is \eqref{eq:grad-umodel-lemmabound}.
\end{proof}

\begin{remark}[Lemma \ref{lem:grad-umodel-bound} is cone-local]
The estimate \eqref{eq:grad-umodel-lemmabound} is restricted to cones
$\sigma\le\sigma_*<\tfrac\pi2$. Indeed, we choose a
label $(R_0,Z_0)$ for which $|\omega_{\theta,0}(R_0,Z_0)|\ge c_0\Gamma$ with a fixed $c_0>0$. The exact
transport identity
\[
\omega_{\mm,\theta}(\Phi_{\lin}(R_0,Z_0,t),t) =J^{-1}\omega_{\theta,0}(R_0,Z_0)
\]
then implies that  $|\omega_{\mm,\theta}|\ge c_0\Gamma J^{-1}$ somewhere. The image point lies in an
\emph{equatorial boundary layer}: the drift law states that
\begin{equation*}
\sigma(t)=\arctan(J^{-3}\tan\sigma_0) =\tfrac\pi2-\arctan\!\big(J^3\cot\sigma_0\big) =\tfrac\pi2-O(J^3) \qquad \ \ \text{ as } J\downarrow 0
.
\end{equation*}
On this layer,  the pullback variables remain $\sigma_{\Lag}=O(1)$ and $\rho_{\Lag}=O(1)$, so neither
$\Upsilon(\sigma_{\Lag})$ nor the spherical tail $\mathcal F(\rho_{\Lag})$ yields any additional smallness;
the exact $J^{-1}$ amplification is visible there.

The equatorial boundary layer is lower order for the collapse clock. The axial strain integral is weighted by
$K_W(\sigma)=3\sin^2\sigma\cos\sigma$, which vanishes at $\sigma=\tfrac\pi2$, so the equatorial boundary layer
contributes only lower order to $\rWm=\p_z(u_\mm)_z(0,0,t)$.
\end{remark}

\subsection{Scaling of the model axial strain $\rWm(t)$}
We now derive the asymptotic scaling law for the model axial strain $\rWm(t)$.  We define the positive angular strain constant
\begin{equation}
C_W^*=C_W^*(\alpha):=\int_0^{\frac\pi2}K_W(\sigma)(\tan\sigma)^\alpha\,\ud\sigma>0.
\label{eq:CW-limit-def}
\end{equation}
\begin{lemma}[Model axial strain scaling]
\label{lem:rW-scaling}
Let $\alpha\in(0,1)$ and $\gamma>\alpha$.  Then, as $\Jm(t)\downarrow0$, the model axial strain satisfies
\begin{equation}
\rWm(t) =-\Gamma\,C_{\rho}^{(1)}(\alpha,\gamma)\,C_W^*\,\Jm(t)^{3\alpha-1}(1+o(1)),
\qquad C_{\rho}^{(1)}(\alpha,\gamma):=\int_0^\infty s^{\alpha-1}\mathcal{F}(s)\,\ud s,
\label{eq:rWm-identity}
\end{equation}
where $C_W^*>0$ is defined in \eqref{eq:CW-limit-def} and
$\mathcal F(s)=(1+s^2)^{-\gamma/2}$.  The constant $C_{\rho}^{(1)}(\alpha,\gamma)$ is finite and satisfies
$C_{\rho}^{(1)}(\alpha,\gamma)=\tfrac{1}{\alpha}+O_{\gamma}(1)$ as $\alpha\downarrow0$.
\end{lemma}
\begin{proof}[Proof of Lemma \ref{lem:rW-scaling}]
We fix $t$ and write $J:=\Jm(t)\in(0,1]$. By odd symmetry and the sign convention in Definition~\ref{def:init-data}, the contribution from the 
lower hemisphere is identical to that from the upper hemisphere, and substituting \eqref{eq:vort-model-explicit} into
\eqref{eq:model-strain-def} provides the identity
\begin{equation}
\rWm(t) =-\Gamma J^{\alpha-1} \int_0^{\Sigma(t)} K_W(\sigma)\,\sin^\alpha\sigma\,
\Upsilon\!\big(\sigma_{\Lag}(\sigma)\big) \Big(\int_0^\infty \rho^{\alpha-1}\mathcal F\big(\rho B_J(\sigma)\big)\,\ud\rho\Big) \ud\sigma,
\label{eq:rWm-spherical-split}
\end{equation}
where
\[
\sigma_{\Lag}(\sigma)=\arctan(J^3\tan\sigma), \qquad \Sigma(t)=\arctan\!\big(J^{-3}\tan\sigma_{\max}\big),
\]
and $B_J(\sigma)$ is defined in \eqref{eq:BJ-def-sec7}.

\runinhead{Step 1: The radial integral.}
For each fixed $\sigma\in[0,\Sigma(t)]$, the radial integral in
\eqref{eq:rWm-spherical-split} is evaluated by the change of variables $s=\rho B_J(\sigma)$:
\[
\int_0^\infty \rho^{\alpha-1}\mathcal F\big(\rho B_J(\sigma)\big)\,\ud\rho
=B_J(\sigma)^{-\alpha}\int_0^\infty s^{\alpha-1}\mathcal F(s)\,\ud s
=B_J(\sigma)^{-\alpha}C_\rho^{(1)}(\alpha,\gamma).
\]
Since $\gamma>\alpha$, the integral defining $C_\rho^{(1)}$ converges at both $0$ and $\infty$. The expansion
$C_\rho^{(1)}(\alpha,\gamma)=\tfrac1\alpha+O_\gamma(1)$ as $\alpha\downarrow0$ follows by splitting
$\int_0^\infty=\int_0^1+\int_1^\infty$, using $\mathcal F(s)=1+O_\gamma(s^2)$ on $(0,1]$, and
$\mathcal F(s)\le s^{-\gamma}$ on $[1,\infty)$.

\runinhead{Step 2: The angular integral.}
We substitute the radial integral from Step 1 into \eqref{eq:rWm-spherical-split} and find that
\begin{equation*}
\rWm(t) =-\Gamma J^{\alpha-1}C_\rho^{(1)}(\alpha,\gamma) \int_0^{\Sigma(t)} K_W(\sigma)\,\sin^\alpha\sigma\,B_J(\sigma)^{-\alpha}
\Upsilon\!\big(\sigma_{\Lag}(\sigma)\big)\,\ud\sigma.
\end{equation*}
Using
\[
B_J(\sigma)^2=J^2\sin^2\sigma+J^{-4}\cos^2\sigma =J^{-4}\cos^2\sigma\big(1+J^6\tan^2\sigma\big),
\]
we obtain the exact identity
\begin{equation*}
\sin^\alpha\sigma\,B_J(\sigma)^{-\alpha} =J^{2\alpha}(\tan\sigma)^\alpha\big(1+J^6\tan^2\sigma\big)^{-\alpha/2}.
\end{equation*}
Therefore
\begin{align}
\rWm(t)
&=-\Gamma C_\rho^{(1)}(\alpha,\gamma)J^{3\alpha-1}
\int_0^{\Sigma(t)} K_W(\sigma)(\tan\sigma)^\alpha
\big(1+J^6\tan^2\sigma\big)^{-\alpha/2}
\Upsilon\!\big(\sigma_{\Lag}(\sigma)\big)\,\ud\sigma.
\label{eq:rWm-spherical-final-prelimit}
\end{align}
For each fixed $\sigma\in[0,\tfrac\pi2)$ we have $\sigma_{\Lag}(\sigma)=\arctan(J^3\tan\sigma)\to0$, hence
$\Upsilon(\sigma_{\Lag}(\sigma))\to1$ as $J\downarrow0$, and also $(1+J^6\tan^2\sigma)^{-\alpha/2}\to1$. Since
$\Sigma(t)\uparrow \tfrac\pi2$ and
\[
0\le K_W(\sigma)(\tan\sigma)^\alpha\big(1+J^6\tan^2\sigma\big)^{-\alpha/2} \Upsilon\!\big(\sigma_{\Lag}(\sigma)\big)
\le K_W(\sigma)(\tan\sigma)^\alpha,
\]
with $K_W(\sigma)(\tan\sigma)^\alpha\in L^1(0,\tfrac\pi2)$ for $\alpha<1$, by the dominated convergence theorem,
we have that
\[
\int_0^{\Sigma(t)} K_W(\sigma)(\tan\sigma)^\alpha \big(1+J^6\tan^2\sigma\big)^{-\alpha/2}
\Upsilon\!\big(\sigma_{\Lag}(\sigma)\big)\,\ud\sigma \longrightarrow \int_0^{\frac{\pi}{2}}K_W(\sigma)(\tan\sigma)^\alpha\,\ud\sigma.
\]
Comparing with  \eqref{eq:CW-limit-def}, we see that the limiting integral equals $C_W^*$. By inserting this into \eqref{eq:rWm-spherical-final-prelimit},
we obtain that 
\[
\rWm(t) =-\Gamma\,C_\rho^{(1)}(\alpha,\gamma)\,C_W^*\,J^{3\alpha-1}(1+o(1)),
\]
which is \eqref{eq:rWm-identity}.
\end{proof}

By combining the model axial strain scaling \eqref{eq:rWm-identity} with the model clock equation
\eqref{eq:J-model-ODE}, we find that 
\[
\dot\Jm(t) =-\tfrac{\Gamma}{2}\,C_\rho^{(1)}(\alpha,\gamma)C_W^*\,\Jm(t)^{3\alpha}(1+o(1)) \ \ \text{ as } \ \ \Jm(t)\downarrow 0.
\]
There are two different angular viewpoints.  Forward in time, particles drift toward the equator, as in
\eqref{eq:forward-angular-drift}.  In the strain integral, however, the Eulerian variable is fixed and the
transported vorticity is evaluated at its pullback label.  That label has
$\sigma_{\Lag}=\arctan(\Jm(t)^3\tan\sigma)$, so a fixed Eulerian angle samples labels closer to the symmetry
axis as $\Jm(t)\downarrow0$.  The identity
$\rho_{\Lag}^{\alpha}\sin^{\alpha}\sigma_{\Lag}=(\Jm(t)r)^\alpha$ therefore converts the spherical cusp into
the cylindrical factor $\Jm(t)^\alpha r^\alpha$.  The algebraic tail is also evaluated at the enlarged radius
$\rho_{\Lag}$, and the radial change of variables in \eqref{eq:rWm-spherical-split} contributes the additional
factor $\Jm(t)^{2\alpha}$.  Multiplying these geometric factors by the transported-vorticity amplification
$\Jm(t)^{-1}$ gives $\rWm(t)\simeq-\Gamma\Jm(t)^{3\alpha-1}$, and the extra factor $\Jm(t)$ in
$\dot\Jm=\tfrac12\Jm\rWm$ gives the power $\Jm(t)^{3\alpha}$ in the model clock law.

\subsection{Finite-time collapse of the clock-and-driver model}\label{sec:stability-analysis-sec7}

We now prove that the model clock reaches zero in finite time exactly in the range
$0<\alpha<\tfrac13$.

\begin{proposition}[Finite-time collapse and the $\alpha=\tfrac13$ barrier]
\label{prop:stability}
Let $\alpha\in(0,1)$ and $\gamma>\alpha$, and let $\Jm$ solve \eqref{eq:J-model-ODE} with $\rWm$ given by
\eqref{eq:vort-model-def}--\eqref{eq:model-strain-def}.  Then $\Jm(t)$ is strictly decreasing and
$\Jm(t)\downarrow0$ as $t\uparrow T_\mm^*$, where the model collapse time is
\begin{equation}
T_\mm^*:=\int_0^1\tfrac{2\,\ud \eta}{-\eta\,\rW_{\mm}(\eta)}\in(0,\infty].
\label{eq:Tm-blowuptime}
\end{equation}

The three collapse regimes are given as follows:
\begin{enumerate}
\item[(i)] If $0<\alpha<\tfrac13$, then $T_\mm^*<\infty$ and
\[
\Jm(t) =\Big((1-3\alpha)\,\tfrac{\Gamma}{2}\,C_{\rho}^{(1)}(\alpha,\gamma)C_W^*\,(T_\mm^*-t)\Big)^{\frac{1}{1-3\alpha}}(1+o(1)) \ \ \text{
as } \ \ t\uparrow T_\mm^*.
\]
\item[(ii)] If $\alpha=\tfrac13$, then $T_\mm^*=\infty$ and $\Jm(t)$ decays exponentially as $t\to\infty$.
\item[(iii)] If $\alpha>\tfrac13$, then $T_\mm^*=\infty$ and $\Jm(t)$ decays algebraically as $t\to\infty$.
\end{enumerate}
\end{proposition}

\begin{proof}[Proof of Proposition \ref{prop:stability}]
\runinhead{Step 1: Monotonicity and reduction to a scalar collapse-time integral.} By
\eqref{eq:model-strain-def}, \eqref{eq:vort-model-explicit}, and the sign of the strain kernel
\eqref{eq:KW-kernel}, the upper and lower hemispheres give equal nonpositive contributions to $\rWm(t)$.
Hence,  $\rWm(t)\le0$ for all $t$, and $\dot{\Jm}(t)=\tfrac12\Jm(t)\rWm(t)\le0$.  Since the integrand in
\eqref{eq:model-strain-def} is not identically zero when $\Jm(t)>0$, we have $\rWm(t)<0$ and $\Jm$ is strictly
decreasing as long as it stays positive.

Therefore either $\Jm$ reaches $0$ at a finite time, or else $\Jm(t)\downarrow0$ as $t\to\infty$.  Indeed,
$\Jm(t)$ has a limit $\ell_\infty\ge0$ by monotonicity.  If $\ell_\infty>0$, then
$\eta\mapsto\rW_{\mm}(\eta)$ is continuous and $\rW_{\mm}(\ell_\infty)<0$, so for all sufficiently large $t$,
\[
\dot{\Jm}(t)\le \tfrac14 \ell_\infty\,\rW_{\mm}(\ell_\infty)<0,
\]
which contradicts convergence to a positive limit.  Hence $\ell_\infty=0$.

Because $\Jm$ is strictly decreasing, we may invert it and regard $t$ as a function of
$\eta=\Jm(t)$. Writing \eqref{eq:J-model-ODE} as $\dot{\Jm}(t)=-f(\Jm(t))$ with
$f(\eta):=-\tfrac12\eta\,\rW_{\mm}(\eta)>0$, we obtain that 
\[
t(\eta)=\int_\eta^1 \tfrac{\ud s}{f(s)} =\int_\eta^1 \tfrac{2\,\ud s}{-s\,\rW_{\mm}(s)}.
\]
Letting $\eta\downarrow 0$ gives \eqref{eq:Tm-blowuptime}.

\medskip
\runinhead{Step 2: Asymptotics of the integrand near the zero clock state.} Lemma~\ref{lem:rW-scaling} gives,
as $\eta\downarrow 0$,
\[
-\rW_{\mm}(\eta) =\Gamma\,C_{\rho}^{(1)}(\alpha,\gamma)\,C_W^*\,\eta^{3\alpha-1}(1+o(1)).
\]
Thus
\[
f(\eta):=-\tfrac12\eta\,\rW_{\mm}(\eta) =\tfrac{\Gamma}{2}\,C_{\rho}^{(1)}(\alpha,\gamma)\,C_W^*\,\eta^{3\alpha}(1+o(1)).
\]

\medskip
\runinhead{Step 3: The $\alpha=\tfrac13$ threshold.} The integral \eqref{eq:Tm-blowuptime} converges at
the lower limit if and only if $\int_0 \eta^{-3\alpha}\,\ud\eta<\infty$, i.e.\ if and only if $3\alpha<1$.
Thus $T_\mm^*<\infty$ exactly when $\alpha<\tfrac13$.

\medskip
\runinhead{Step 4: Collapse rate as $t\uparrow T_\mm^*$.} In the subcritical case, define
$Y(t):=\Jm(t)^{1-3\alpha}$. Using $\dot{\Jm}(t)=-f(\Jm(t))$ and the asymptotic in Step 2 yields
\[
\dot Y(t) =(1-3\alpha)\Jm(t)^{-3\alpha}\dot{\Jm}(t) =-(1-3\alpha)\,\tfrac{\Gamma}{2}\,C_{\rho}^{(1)}(\alpha,\gamma)C_W^*(1+o(1)),
\]
so $Y(t)$ is asymptotically linear near $t=T_\mm^*$ and hence
\[
\Jm(t) =\Big((1-3\alpha)\,\tfrac{\Gamma}{2}\,C_{\rho}^{(1)}(\alpha,\gamma)C_W^*\, (T_\mm^*-t)\Big)^{\!\!\frac{1}{1-3\alpha}}(1+o(1)),
\]
as claimed.

The critical and supercritical cases follow from the same scalar ODE
$\dot{\Jm}=-f(\Jm)$ and the asymptotic
$f(\eta)=c_{\alpha,\gamma,\Gamma}\eta^{3\alpha}(1+o(1))$ as $\eta\downarrow0$, where
$c_{\alpha,\gamma,\Gamma}>0$.  For $\alpha=\tfrac13$,  this produces exponential decay and infinite collapse time.
For $\alpha>\tfrac13$ it gives algebraic decay and infinite collapse time.
\end{proof}

\section{Decomposition into Smooth and Cusp Flows}
\label{sec:Euler-blowup-for-Theta-star}

This section starts the construction used in the proof of Theorem~\ref{thm:target-profile}.  We split the exact Euler
flow into a smooth map and a cusp map,
\[
\phi=\phi_{\smooth}\circ\phi_{\cusp},
\]
with the precise definitions given in \eqref{eq:smooth-velocity-def}--\eqref{eq:cusp-map-def}.  The smooth map is
generated by the far-field velocity $u_{\smooth}$ in \eqref{eq:smooth-velocity-def}; the cutoff in that definition is
evaluated at the Eulerian position $\phi(Y,t)$ of each label $Y$ at the same time $t$.  The
cusp map carries the near-field motion that produces the collapsing clock $J_{\cusp}$ in
\eqref{eq:smooth-cusp-clock-def}.  After these maps are defined, we state
the pressure decomposition \eqref{eq:Pi0-decomp}, the small-clock bootstrap assumptions \textup{(BA1)}--\textup{(BA9)},
and the finite-clock entry lemma that brings the solution into a prescribed small-clock regime.  These ingredients
prepare the later comparison between $J_{\cusp}$ and the model clock $\Jm$ from
Section~\ref{sec:lag-analysis-non-local}.

\begin{remark}[Standing decay hypothesis used for tail bounds]
Throughout Sections~\ref{sec:Euler-blowup-for-Theta-star}--\ref{sec:target-profile-typeI-completion} we use the decay
assumption from Theorem~\ref{thm:target-profile},
\begin{equation}\label{eq:gamma-stand}
\gamma>\alpha+\tfrac52.
\end{equation}
This hypothesis is used only for estimates involving the algebraic tail of the datum.  In those estimates,
\eqref{eq:gamma-stand} gives summability of the velocity, velocity-gradient, and pressure Hessian remainders
generated by labels far from the collapsing core.  Several local or model estimates below require only weaker conditions
such as $\gamma>\alpha$ or $\gamma>2\alpha$, but the collapse argument is carried out under the finite-energy assumption
\eqref{eq:gamma-stand} throughout.
\end{remark}

\subsection{Geometric flow decomposition, velocity fields, and domains}
\label{sec:geom-flow-decomp}

Let $\phi(Y,t)$ denote the true Euler flow map.  We decompose it into a smooth flow map and a cusp flow map:
\[
\phi(Y,t)=\phi_{\smooth}(\phi_{\cusp}(Y,t),t).
\]
The corresponding velocity fields and label domains are defined below.  We use capital letters for Lagrangian labels and
lower-case letters for their images.  Thus, if $x=\Lambda(X,t)$, then $X$ is the material label and $x$ is the
associated Eulerian position.

\begin{remark}[Meaning of order one in the cusp-clock limit]
Throughout the target-profile proof, a smooth or tail-generated quantity is called \emph{order one} relative to the cusp
clock if it remains uniformly bounded as $J_{\cusp}\downarrow0$, with constants independent of the small clock value.
Thus an order-one velocity gradient may be bounded by $C\Gamma$, and an order-one smooth clock satisfies $0<c\le
J_{\smooth}(t)\le C<\infty$ on the small-cusp-clock interval.  The constants may depend on the fixed parameters
$\alpha,\gamma,\sigma_{\inn},\sigma_*$ and on the fixed Target Profile datum, but they do not carry negative powers of
$J_{\cusp}$.  This is the sense in which the smooth-flow and far-field terms are lower order than the singular strain
scale $\Gamma J_{\cusp}^{3\alpha-1}$ and the singular pressure scale $\Gamma^2J_{\cusp}^{6\alpha-2}$.
\end{remark}

\begin{remark}[Clock thresholds and clock-scaled axial notation]
Time-dependent clock variables are written with plain $J$-symbols, such as $J(t)$, $J_{\smooth}(t)$, $J_{\cusp}(t)$, and
$J_{\twoD}(Y,t)$.  Fixed clock thresholds are written with fraktur symbols
\[
\mathfrak J_{\mathrm{tail}},\quad \mathfrak J_{\mathrm{entry}},\quad \mathfrak J_{\pressure},\quad \ldots .
\]
Thus $\mathfrak J_{\mathrm{tail}}\in(0,1]$ is a fixed number chosen in the proof, not an additional time-dependent
clock.  The subscript $\sharp$ is used for the fixed $\zeta$-localization interval and cutoff used in the pressure
Hessian comparison, such as $I_\sharp$ and $\vartheta_\sharp$; it is not attached to a clock variable.
\end{remark}

\subsubsection{Cones and label domains.}
We use fixed polar angles $\sigma_{\inn}$ and $\sigma_*$ such that
\[
0<\sigma_{\inn}<\sigma_*<\tfrac{\pi}{2} \qquad\text{and}\qquad \sigma_{\max}<\sigma_{\inn},
\]
so that the angular support of $\Theta^*$ lies strictly inside the inner cone $\mathcal
C_{\inn}\Subset\mathcal C_*$. In particular,
\begin{equation}\label{eq:sigma-angles}
0<\sigma_{\cut}<\sigma_{\max}<\sigma_{\inn}<\sigma_*<\tfrac{\pi}{2} \, ,
\end{equation}
where the profile $\Theta^*$ and the polar angles $\sigma_{\cut}<\sigma_{\max}$ are given in
Definition~\ref{def:init-data}. The final values of $\sigma_{\inn}$ and $\sigma_*$ are chosen in the order specified in
Section~\ref{sec:fixed-choice-order}; all estimates preceding that final choice hold for any fixed pair satisfying
\eqref{eq:sigma-angles}.

\begin{remark}[Compatibility of cone angles with Definition~\ref{def:init-data}]
Definition~\ref{def:init-data} requires only $0<\sigma_{\cut}<\sigma_{\node}<\sigma_{\max}<\tfrac\pi2$, where
$\sigma_{\node}=\arccos(1/\sqrt{3})$ is the maximizer of the axial-strain kernel $K_W$. Since the
interval $(\sigma_{\max},\tfrac\pi2)$ is nonempty, we can always choose intermediate angles
$\sigma_{\inn}\in(\sigma_{\max},\tfrac\pi2)$ and $\sigma_*\in(\sigma_{\inn},\tfrac\pi2)$ to obtain the full chain
\eqref{eq:sigma-angles}. In particular, no further restriction on the initial data is needed: $\sigma_{\max}$ is fixed
by the initial datum, and the cone angles $\sigma_{\inn},\sigma_*$ are chosen later as fixed proof parameters, after the
angular-slope cutoff $M_\pressure$ has been fixed in
\eqref{eq:pressure-C0}.  They remain independent of the amplitude $\Gamma$, the time variable, and
the small-clock thresholds.
\end{remark}

With $\sigma_*$ as in \eqref{eq:sigma-angles}, we denote the associated \emph{fixed cone near the symmetry axis} by
\begin{equation}
\mathcal{C}_{\sigma_*} :=\Big\{x\in\R^3\setminus\{0\}:\ 0\le \sigma(x)\le \sigma_*\Big\}.
\label{eq:buffered-cone-def}
\end{equation}
We write
\begin{equation}
\mathcal C_{\inn}:=\mathcal C_{\sigma_{\inn}}, \qquad \mathcal C_*:=\mathcal C_{\sigma_*}.
\label{eq:inner-star-cone-notation}
\end{equation}
The final tail radius $R_{\tail}\ge2$ is fixed in the choice order of Section~\ref{sec:fixed-choice-order}, large enough for the
far-field smooth-flow estimates and for the $\zeta$-tail pressure Hessian estimates used below.  Once this value is
fixed, we set
\begin{equation}
D_{\core}:=\{Y\in\R^3:\ |Y|\le R_{\tail}\}, \qquad D_{\tail}:=\R^3\setminus D_{\core}.
\label{eq:core-tail-domains}
\end{equation}

The bounded-core cusp label set is the fixed set $D_{\core}$.  The moving cone condition below is used only to identify
the near-axis part of this set on which the singular cone-local estimates are applied.  Once the cusp flow has been
defined, we set
\begin{equation}\label{eq:cusp-domain-def}
D_{\inn}^{\cusp}(t):=\Bigl\{Y\in D_{\core}:\ \phi_{\cusp}(Y,t)\in\mathcal C_{\inn}\Bigr\}.
\end{equation}
Thus $D_{\inn}^{\cusp}(t)$ is a moving analysis subset of the fixed core set $D_{\core}$.  This distinction is
important: labels in $D_{\core}$ whose Eulerian images have moved outside $\mathcal C_{\inn}$ may still carry the
cusp-clock amplification, so they remain part of the cusp-transported velocity and are estimated with that velocity.

\subsubsection{Smooth velocity and smooth flow.}
Fix a smooth radial cutoff
\[
\chi_{\far}\in C^\infty([0,\infty)), \qquad 0\le\chi_{\far}\le1,\qquad \chi_{\far}=0\ {\rm on}\ [0,1],\qquad \chi_{\far}=1\ {\rm on}\
[2,\infty).
\]
We define the smooth velocity by
\begin{equation}\label{eq:smooth-velocity-def}
u_{\smooth}(x,t) := \tfrac1{4\pi}\int_{\R^3} K\bigl(x,\phi(Y',t)\bigr)\, \chi_{\far}\!\left(\tfrac{|\phi(Y',t)|}{R_{\tail}}\right)
J_{\twoD}(Y',t)^{-1}\omega_{\theta,0}(Y')\,dY'.
\end{equation}
In Eulerian variables, $u_{\smooth}$ is the Biot--Savart velocity generated by the Euler vorticity at time $t$
restricted to the physical far field $|y|\ge R_{\tail}$, with a smooth transition on
$R_{\tail}\le |y|\le2R_{\tail}$.  Thus a fixed label $Y\in D_{\tail}$ contributes to $u_{\smooth}$ with the
weight $\chi_{\far}(|\phi(Y,t)|/R_{\tail})$; if $|\phi(Y,t)|\le R_{\tail}$, its contribution is entirely in
$u-u_{\smooth}$, and if $|\phi(Y,t)|\le2R_{\tail}$ it is in the cutoff transition.  In particular, a label
$Y\in D_{\tail}$ with $\phi_{\cusp}(Y,t)\in\mathcal C_{\inn}$ is treated in the cusp-coordinate velocity whenever
$|\phi(Y,t)|\le2R_{\tail}$; its large initial radius is estimated later as the algebraic-tail contribution to the
pressure Hessian comparison.

Let $\phi_{\smooth}$ be the flow map of this smooth velocity:
\begin{equation}\label{eq:smooth-flow-def}
\p_t\phi_{\smooth}(X,t)=u_{\smooth}(\phi_{\smooth}(X,t),t),
\qquad
\phi_{\smooth}(X,0)=X.
\end{equation}
Since $u_{\smooth}$ is divergence-free, axisymmetric, and no-swirl, $\phi_{\smooth}$ is volume preserving and preserves
the symmetry axis.

\subsubsection{Push-forward and pull-back conventions.}
For a volume-preserving axisymmetric diffeomorphism $\Lambda(\cdot,t)$, we define the push-forward of a vector field $W$
in the label variables by
\[
(\Lambda_*W)(x,t):=D_X\Lambda(X,t)\,W(X,t), \qquad x=\Lambda(X,t),\quad X=\Lambda^{-1}(x,t).
\]
The pull-back $\Lambda^*v$ of an Eulerian-type vector field $v$ in the target variables is the push-forward by
$\Lambda^{-1}$,
\[
(\Lambda^{-1}_*v)(X,t):=D_X\Lambda(X,t)^{-1}v(\Lambda(X,t),t).
\]

\subsubsection{Exact cusp map, clocks, and cusp-coordinate velocity.}
The exact cusp map is defined by removing the smooth flow from the true Euler flow:
\begin{equation}\label{eq:cusp-map-def}
\phi_{\cusp}(Y,t):=\phi_{\smooth}^{-1}(\phi(Y,t),t),
\qquad
V_{\cusp}:=(\phi_{\smooth}^{-1})_*(u-u_{\smooth}).
\end{equation}
Thus
\begin{equation}
\phi(Y,t)=\phi_{\smooth}(\phi_{\cusp}(Y,t),t), \qquad \p_t\phi_{\cusp}(Y,t)=V_{\cusp}(\phi_{\cusp}(Y,t),t).
\label{eq:cusp-flow-equation}
\end{equation}
The corresponding smooth and cusp clocks are
\begin{equation}
J_{\smooth}(t):=\det\nabla_{(R,Z)}\phi_{\smooth}(0,t), \qquad J_{\cusp}(t):=\det\nabla_{(R,Z)}\phi_{\cusp}(0,t),
\label{eq:smooth-cusp-clock-def}
\end{equation}
and the clock decomposition is
\begin{equation}\label{eq:J-clock-decomposition}
J(t)=J_{\smooth}(t)\,J_{\cusp}(t).
\end{equation}
The stagnation-point  \emph{exact cusp} axial strain associated with the exact cusp flow is
\begin{equation}
\rW_{\cusp}(t):=\p_z(V_{\cusp})_z(0,t).
\label{eq:rW-cusp-def}
\end{equation}

\subsubsection{Flat transported cusp velocity and scalar modulation.}
The cusp flow map transports the initial vorticity.  This transported vorticity, together with its flat Biot--Savart
velocity, is the object to which the slope-restricted pressure Hessian estimate is applied.  Let
\[
\mathcal J_{\cusp}(Y,t):= \det\nabla_{(R,Z)}\bigl((\phi_{\cusp})_r,(\phi_{\cusp})_z\bigr)(Y,t).
\]
For $R(Y)>0$, volume preservation gives
\[
\mathcal J_{\cusp}(Y,t)=\tfrac{R(Y)}{(\phi_{\cusp})_r(Y,t)}.
\]
We define the cusp-flow transported angular vorticity by
\begin{equation}
\bs\Omega_{\cusp}(\phi_{\cusp}(Y,t),t) := \mathcal J_{\cusp}(Y,t)^{-1}\, \omega_{\theta,0}(Y)\,\bs e_\theta(\phi_{\cusp}(Y,t)),
\label{eq:Omega-cusp-def}
\end{equation}
with the continuous extension across the symmetry axis.  This transported angular vorticity $\bs\Omega_{\cusp}$ then
generates the flat Biot--Savart velocity $U_{\cusp}(\cdot,t):=\BS[\bs\Omega_{\cusp}(\cdot,t)]$, given explicitly by
\begin{equation}
U_{\cusp}(x,t) = \tfrac1{4\pi}\int_{\R^3} K\bigl(x,\phi_{\cusp}(Y',t)\bigr)\, \mathcal J_{\cusp}(Y',t)^{-1}\omega_{\theta,0}(Y')\,dY'.
\label{eq:U-cusp-label}
\end{equation}
Thus lower-case velocities such as $u,u_{\smooth},u_{\cusp}$ are physical Eulerian velocities in the physical variable
$x$.  By contrast, $U_{\cusp}$ is written in the cusp-coordinate variable before push-forward by the smooth flow map. We
denote its stagnation-point \emph{flat cusp} axial strain by
\begin{equation}
\mathcal W_{\cusp}(t):= \p_z(U_{\cusp})_z(0,t).
\label{eq:W-cusp-def}
\end{equation}
Thus $\rW_{\cusp}$ is the cusp-flow strain generated by the actual cusp-coordinate velocity $V_{\cusp}$ in
\eqref{eq:cusp-map-def}, whereas $\mathcal W_{\cusp}$ is the flat cusp-coordinate strain generated by
$U_{\cusp}$ in \eqref{eq:U-cusp-label}.  The pressure comparison is carried out first for $U_{\cusp}$ because the
vorticity locations in \eqref{eq:U-cusp-label} are the exact cusp-flow locations.  These locations have the scaled
structure used in Proposition~\ref{prop:euler-generated-profile-riccati}.  The difference between $V_{\cusp}$ and
$U_{\cusp}$ comes from the smooth change of variables $\phi_{\smooth}$, which deforms the Biot--Savart kernel.  In the
present Target Profile case, the identity
\[
\rho^\alpha\Theta^*(\sigma)=R^\alpha
\]
on the support sector is what produces the scaled vorticity in \eqref{eq:vort-slope-restricted}.
Proposition~\ref{prop:euler-generated-profile-riccati} gives the Riccati pressure Hessian bound for the
Euler-generated axial function.

With \eqref{eq:rW-cusp-def} and \eqref{eq:W-cusp-def}, the associated scalar modulation function is
\begin{equation}
m(t):=\frac{\rW_{\cusp}(t)}{\mathcal W_{\cusp}(t)}, \qquad \rW_{\cusp}(t):=\p_z(V_{\cusp})_z(0,t),
\qquad \mathcal W_{\cusp}(t):=\p_z(U_{\cusp})_z(0,t).
\label{eq:modulation-def}
\end{equation}
Here ``scalar modulation'' means only that the cusp-coordinate velocity $U_{\cusp}$ is multiplied by a single
time-dependent scalar so that its stagnation-point flat cusp axial strain agrees exactly with the stagnation-point exact
cusp axial strain $\rW_{\cusp}(t)$.  This is a geometric normalization of the stagnation-point strain strength.  If
$m\equiv1$, the cusp-flow velocity is the unmodulated transported velocity.

The smooth velocity is defined by the far-field cutoff in \eqref{eq:smooth-velocity-def}, evaluated at the
Eulerian position $\phi(Y,t)$ at the same time $t$; it is not a fixed-label tail.  The next lemma states the exact flow split
\eqref{eq:cusp-flow-equation} and the clock identity \eqref{eq:J-clock-decomposition} used throughout the
small-clock analysis.

\begin{lemma}[Exact decomposition into  smooth and cusp flows]
\label{lem:self-consistent-cusp-flow}
Let the exact Euler solution be fixed on a finite time interval $[0,T]$.  The velocity $u_{\smooth}$, the smooth flow
$\phi_{\smooth}$, and the cusp flow $\phi_{\cusp}$ defined in
\eqref{eq:smooth-velocity-def}--\eqref{eq:cusp-map-def} are axisymmetric, no-swirl objects on
$[0,T]$.  They satisfy
\[
\phi=\phi_{\smooth}\circ\phi_{\cusp}, \qquad J(t)=J_{\smooth}(t)J_{\cusp}(t)
\]
and
\[
\p_t\phi_{\cusp}(Y,t)=V_{\cusp}(\phi_{\cusp}(Y,t),t).
\]
The flat velocity $U_{\cusp}$ and the scalar modulation $m(t)$ are then determined by
\eqref{eq:Omega-cusp-def}--\eqref{eq:modulation-def} whenever
$\mathcal W_{\cusp}(t)\neq0$.
\end{lemma}

\begin{proof}[Proof of Lemma~\ref{lem:self-consistent-cusp-flow}]
The exact Euler flow $\phi$ is already fixed.  The smooth velocity \eqref{eq:smooth-velocity-def} is
determined directly by the Euler vorticity at time $t$ and the far-field cutoff evaluated at
$|\phi(Y,t)|/R_{\tail}$, and $\phi_{\smooth}$ is its flow.  We then define
\[
\phi_{\cusp}(Y,t):=\phi_{\smooth}^{-1}(\phi(Y,t),t), \qquad V_{\cusp}:=(\phi_{\smooth}^{-1})_*(u-u_{\smooth}),
\]
so the exact identity $\phi=\phi_{\smooth}\circ\phi_{\cusp}$ and the transport law for the cusp flow follow from the
chain rule.

The transported velocity $U_{\cusp}$ is now determined by the exact cusp flow map $\phi_{\cusp}$ through
\eqref{eq:Omega-cusp-def}--\eqref{eq:W-cusp-def}.  On every time
subinterval on which $\mathcal W_{\cusp}(t)\neq0$, the scalar $m(t)$ in
\eqref{eq:modulation-def} is therefore known from the exact solution.  The finite-clock and small-clock
estimates proved below provide the required separation of $\mathcal W_{\cusp}$ from zero on the intervals used in the
collapse argument.

The smooth velocity and the transported cusp velocity are divergence-free and axisymmetric without swirl, hence their
flows preserve the symmetry axis. The identity for $\phi_{\cusp}$, the clock decomposition, and the transport law for
the cusp flow are direct consequences of
\[
\phi_{\cusp}(Y,t)=\phi_{\smooth}^{-1}(\phi(Y,t),t)
\]
and the chain rule.
\end{proof}

\subsubsection{Cusp-coordinate error velocity and physical decomposition.}
The scalar-modulated flat velocity $m(t)U_{\cusp}$ has the same stagnation-point axial strain as the exact
cusp-coordinate velocity $V_{\cusp}$.  We define
\begin{equation}\label{eq:Verr-def}
V_{\err}:=V_{\cusp}-m(t)\,U_{\cusp}(\cdot,t).
\end{equation}
This decomposition is needed because the Riccati pressure Hessian term in \eqref{eq:stagnation-ODEs} is quadratic in the
physical Eulerian gradient $\nabla u$, whereas the slope-restricted pressure Hessian estimate applies first to the flat
cusp-coordinate velocity $U_{\cusp}$ in \eqref{eq:U-cusp-label}. In physical variables we set
\[
u_{\cusp}:=(\phi_{\smooth})_*\!\bigl(m(t)\,U_{\cusp}(\cdot,t)\bigr), \qquad u_{\err}:=(\phi_{\smooth})_*V_{\err},
\]
so that
\begin{equation}\label{eq:u-decomp}
u=u_{\smooth}+u_{\cusp}+u_{\err}.
\end{equation}

\subsubsection{Auxiliary cusp map.}
The exact pressure argument uses the cusp map $\phi_{\cusp}$ directly.  We nevertheless introduce the auxiliary flow
driven by the scalar-modulated cusp-coordinate velocity:
\begin{equation}
\p_t \Phi_{\cusp}(Y,t) = m(t)\, U_{\cusp}(\Phi_{\cusp}(Y,t),t), \qquad \Phi_{\cusp}(Y,0)=Y.
\label{eq:Phi-cusp-def}
\end{equation}
This auxiliary map is not used in the Biot--Savart integral defining $U_{\cusp}$; the vorticity that generates
$U_{\cusp}$ is evaluated at the exact cusp-flow position $\phi_{\cusp}(Y,t)$.  In the small-clock estimates
below, the axial motion is tracked by the one-dimensional composition $\mathscr Z_t$ from
Lemma~\ref{lem:exact-axis-composition}, rather than by the auxiliary map $\Phi_{\cusp}$.

\subsubsection{Pressure decomposition.}
For time-dependent physical vector fields $v,w$, we use the time-dependent version of
\eqref{eq:pressure-bilinear-form},
\begin{equation}
\Pi[v,w](t) := \pv\!\int_{\R^3}K_{zz}(y)\, \tr\!\big(\nabla v(y,t)\nabla w(y,t)\big)\,dy.
\label{eq:physical-pressure-bilinear-form}
\end{equation}
Since $\tr(AB)=\tr(BA)$, inserting the velocity split \eqref{eq:u-decomp} into $S=\tr((\nabla u)^2)$ gives pure terms
and bilinear cross terms.  The pure physical cusp term is $\Pi[u_{\cusp},u_{\cusp}]$.  The corresponding
cusp-coordinate pressure associated with the velocity $U_{\cusp}$ is
\begin{equation}\label{eq:Pi-U-cusp-def}
\Pi_{\cusp}(t):=\Pi[U_{\cusp},U_{\cusp}](t).
\end{equation}
Since $u_{\cusp}$ is the smooth-flow push-forward of the scalar-modulated cusp-coordinate velocity $U_{\cusp}$, we
separate its cusp-coordinate pressure from the geometric defect
\[
\Pi_{\rm geom}(t) := \Pi[u_{\cusp},u_{\cusp}](t) -m(t)^2\Pi_{\cusp}(t).
\]
Thus $\Pi_{\rm geom}$ measures exactly the effect of evaluating the pure cusp pressure Hessian after the smooth change
of variables. The remaining bilinear and pure non-cusp terms are
\begin{align*}
\Pi_{\mixed}(t)
&:= 2\Pi[u_{\cusp},u_{\smooth}](t) +2\Pi[u_{\cusp},u_{\err}](t) +2\Pi[u_{\smooth},u_{\err}](t), \\
\Pi_{\smooth}(t)&:=\Pi[u_{\smooth},u_{\smooth}](t), \\
\Pi_{\err}(t)&:=\Pi[u_{\err},u_{\err}](t).
\end{align*}
With these conventions, the pressure Hessian at the stagnation point has the decomposition
\begin{equation}\label{eq:Pi0-decomp}
\Pi_0(t)
=m(t)^2\Pi_{\cusp}(t) +\Pi_{\rm geom}(t) +\Pi_{\mixed}(t)+\Pi_{\smooth}(t)+\Pi_{\err}(t).
\end{equation}

\subsection{Riccati law and pressure Hessian comparison}

The collapse clock $J(t)$, axial strain $\rW_0(t)$, and axial pressure Hessian $\Pi_0(t)$ are defined at the stagnation
point in \eqref{eq:rW0-Pi0-J0-defs}.  They satisfy
\begin{equation}\label{eq:stagnation-ODEs}
\dot J(t)=\tfrac12\,J(t)\,\rW_0(t),
\quad
\p_t \rW_0(t)=-\tfrac12\,\rW_0(t)^2-\Pi_0(t),
\quad
\Pi_0(t)=\pv\!\int_{\R^3}K_{zz}(y)\,\tr\!\big((\nabla u(y,t))^2\big)\,dy,
\end{equation}
where $K_{zz}$ is the axisymmetric pressure kernel from Section~\ref{sec:Kzz-kernel}. The decomposition
\eqref{eq:Pi0-decomp} is used throughout the proof to identify which part of the pressure Hessian carries the singular
strain.  The scalar-modulated transported cusp term $m(t)^2\Pi_{\cusp}(t)$ is the principal term.  The term $\Pi_{\rm
geom}$ measures the change in the pure cusp pressure Hessian caused by the smooth flow map, and
$\Pi_{\mixed},\Pi_{\smooth}$, and $\Pi_{\err}$ measure the remaining mixed, smooth, and error terms.

\subsubsection{Riccati reduction with subcritical pressure Hessian.}
The inertial term $\tfrac12\rW_0(t)^2$ is the compressive term in the stagnation-point Riccati equation. In view of
\eqref{eq:stagnation-ODEs}, a nonnegative pressure Hessian has the same sign as this compression:
\[
\p_t\rW_0(t) = -\tfrac12\rW_0(t)^2-\Pi_0(t).
\]
Thus the pressure Hessian cannot cancel the collapse once its leading part is subcritical relative to the strain.  The
remaining issue is to show that the errors introduced by cutoffs, normal-form transfer, and the smooth change of
variables do not consume the Riccati slack.  For the exact Euler pressure Hessian we prove
\[
\Pi_0(t) \ge -\upbeta\,\tfrac12\,\rW_0(t)^2
\]
with a fixed constant
\begin{equation}
0<\upbeta<1 .
\label{eq:beta-def}
\end{equation}
chosen after the Euler-generated axial function in the normalized axial coordinate is fixed and before the error
margins are fixed.  Proposition~\ref{prop:euler-generated-profile-riccati} gives the renormalized axis-trace
Riccati comparison for that profile.  The cutoff and small-clock parameters below are selected so that the angular tail
\eqref{eq:pressure-angular-tail}, the $\zeta$-tail
\eqref{eq:pressure-zeta-tail}, the normal-form error \eqref{eq:additive-normal-form-bound}, and the smooth-flow
transfer error fit inside the remaining slack.  Consequently the stagnation-point Riccati law gives
\[
\p_t\rW_0(t) \le -\tfrac{1-\upbeta}{2}\,\rW_0(t)^2,
\]
which is a finite-time blowup inequality for negative $\rW_0$.

For the Target Profile \eqref{eq:Theta-star-def}, the leading transported cusp pressure Hessian is controlled by the
origin-attached Euler-generated axial function in Proposition~\ref{prop:euler-generated-profile-riccati}.  The
monotone-moment estimate gives the one-dimensional compression inequality, and
Lemma~\ref{lem:transported-cusp-pressure-win} verifies the axis-trace hypotheses for the cusp-coordinate pressure Hessian.  The later
pressure Hessian estimates compare all remainder terms in \eqref{eq:Pi0-decomp} with $\rW_{\cusp}(t)^2$.  After the
small-cusp-clock estimates for $m(t)$ and $\mathcal W_{\cusp}(t)$ have been proved, we will show that the remaining
pressure Hessian contributions in
\eqref{eq:Pi0-decomp} are  lower order compared with
\[
\rW_{\cusp}(t)^2 = m(t)^2\mathcal W_{\cusp}(t)^2,
\qquad c\,\Gamma^2 J_{\cusp}(t)^{6\alpha-2} \le \rW_{\cusp}(t)^2 \le C\,\Gamma^2 J_{\cusp}(t)^{6\alpha-2}
\]
on the active small-clock interval.

\subsection{Roadmap of the collapse argument}

The purpose of this section and Sections~\ref{sec:slope-restricted-pressure}--\ref{sec:target-profile-typeI-completion}
is to prove that the true Euler clock $J(t):=\det\nabla_{(R,Z)}\phi(0,0,t)$ obeys the same collapse law as the linear
model clock $\Jm(t)$ from Section~\ref{sec:stability-analysis-sec7}.  The linear model gives
\[
\dot\Jm(t)=\tfrac12\Jm(t)\rWm(t), \qquad \rWm(t)\simeq -\Gamma C_{\rho}^{(1)}(\alpha,\gamma)C_W^*\,\Jm(t)^{3\alpha-1},
\]
and hence
\[
\dot\Jm(t)\simeq-\Gamma\Jm(t)^{3\alpha}.
\]
We prove the corresponding statement for the  Euler solution by using the exact cusp flow $\phi_{\cusp}$.

\runinhead{Step 1: Exact smooth and cusp flow decomposition.} The exact Lagrangian flow is split as
\[
\phi(Y,t)=\phi_{\smooth}(\phi_{\cusp}(Y,t),t),
\]
where $\phi_{\smooth}$ is generated by the far-field velocity $u_{\smooth}$ defined by
\eqref{eq:smooth-velocity-def}.  This gives the exact clock identity
\[
J(t)=J_{\smooth}(t)J_{\cusp}(t).
\]
The smooth velocity is generated by labels whose Eulerian images $\phi(Y,t)$ lie away from the collapsing core.
After the tail radius $R_{\tail}$ is chosen sufficiently large, the smooth flow is a small deformation of the
identity on the time interval relevant to collapse:
\[
\|D\phi_{\smooth}-I\|_{L^\infty} + \|D\phi_{\smooth}^{-1}-I\|_{L^\infty} \le C\varepsilon_{\smooth}, \qquad \varepsilon_{\smooth}\ll1.
\]
In particular $J_{\smooth}$ stays bounded above and below, so collapse of $J_{\cusp}$ is equivalent to collapse of the
true Euler clock $J$.

\runinhead{Step 2: The exact cusp velocity and the flat Biot--Savart velocity.} The exact velocity driving the cusp map
is
\[
\p_t\phi_{\cusp}(Y,t)=V_{\cusp}(\phi_{\cusp}(Y,t),t), \qquad V_{\cusp}=(\phi_{\smooth}^{-1})_*(u-u_{\smooth}).
\]
If $X=\phi_{\cusp}(Y,t)$ and $\Lambda=\phi_{\smooth}$, then
\[
V_{\cusp}(X,t) = D\Lambda(X,t)^{-1} (u-u_{\smooth})(\Lambda(X,t),t).
\]
Thus $V_{\cusp}$ is the exact cusp-coordinate velocity, but its Biot--Savart representation contains the smoothly
deformed kernel
\[
D\Lambda(X,t)^{-1}K(\Lambda(X,t),\Lambda(X',t))D\Lambda(X',t).
\]
The pressure lemmas below are instead flat Euclidean Biot--Savart statements.  For this reason we introduce the flat
velocity
\[
U_{\cusp}:=\BS[\Omega_{\cusp}],\qquad \Omega_{\cusp}(\phi_{\cusp}(Y,t),t)
=\mathcal J_{\cusp}(Y,t)^{-1}\omega_{\theta,0}(Y)e_\theta(\phi_{\cusp}(Y,t)).
\]
The smallness of $D\phi_{\smooth}-I$ implies that $U_{\cusp}$ gives the leading Euclidean-kernel description of
$V_{\cusp}$; the kernel deformation caused by $\phi_{\smooth}$ is an
$O(\varepsilon_{\smooth})$ pressure Hessian error.

This also explains the limited role of the auxiliary cusp map.  The exact cusp flow $\phi_{\cusp}$, not an auxiliary
map, transports the vorticity used in the pressure calculation.  The auxiliary map $\Phi_{\cusp}$ is introduced only
after the cusp-coordinate velocity has been defined:
\[
\p_t\Phi_{\cusp} = m(t)U_{\cusp}(\Phi_{\cusp},t).
\]
The scalar $m(t)$ is chosen below so that $m(t)U_{\cusp}$ has the exact stagnation-point cusp strain.  Consequently the
dominant hyperbolic field is shared by $\phi_{\cusp}$ and $\Phi_{\cusp}$ in cusp coordinates.  The pressure comparison
and the small-clock closure use the exact cusp map $\phi_{\cusp}$, together with the normal-form estimate
\eqref{eq:localized-normal-form-representation}, the geometry of $\Psi_t$ in
\eqref{eq:localized-normal-form-map-bootstrap}, and the axial-composition estimates proved below.

\runinhead{Step 3: Normal form in the axis coordinate at time $t$ and the pressure lemmas.} At a small-clock time $t$,
with $J:=J_{\cusp}(t)$, we write the exact cusp map $\phi_{\cusp}$ in cylindrical coordinates as
$\phi_{\cusp}(R,Z,t)=(r_t(R,Z),z_t(R,Z))$.  On the symmetry axis, we set $B_t(Z):=z_t(0,Z)$ and
$A_t(Z):=\p_Rr_t(0,Z)$.  The normalized axial coordinate is the axial position at time $t$ divided by the
collapsing axial scale:
\[
\zeta=J^{-2}B_t(Z).
\]
For a bounded slope $\tau$, let $Z_t(\zeta)=(J^{-2}B_t)^{-1}(\zeta)$ and choose $R_t(\zeta,\tau)$ by
$A_t(Z_t(\zeta))R_t(\zeta,\tau)=J^2\zeta\tau$.  Thus $Y_t(\zeta,\tau):=(R_t(\zeta,\tau),Z_t(\zeta))$ is the original
label whose axial flow map data are $(\zeta,\tau)$; in inverse-map notation, if
$x_t(\zeta,\tau)=\phi_{\cusp}(Y_t(\zeta,\tau),t)$, then $Y_t(\zeta,\tau)=\phi_{\cusp}^{-1}(x_t(\zeta,\tau),t)$.  The
normal form says that this exact point is the linear hyperbolic placement $J^2\zeta(\tau,1)$ plus a small error:
\[
(r,z)\bigl(\phi_{\cusp}(Y_t(\zeta,\tau),t)\bigr) = J^2\zeta\bigl((\tau,1)+\mathcal E_t(\zeta,\tau)\bigr),
\]
with $\mathcal E_t$ small in the norms stated in Lemma~\ref{lem:late-axis-normal-form-cusp}.  The identity $\rho^\alpha\Theta^*(\sigma)=R^\alpha$
on the target support sector then converts the localized cusp-flow transported vorticity into the scaled form
\[
-\operatorname{sgn}(Z)\,a_t(|Z|)R^\alpha \chi_M\!\left(\tfrac{R}{|Z|}\right)e_\theta .
\]
The profile $a_t$ is the physical axial function on the localized $\zeta$-support.  The monotone axial-stretching
bootstrap \eqref{eq:monotone-axial-two-sided}--\eqref{eq:monotone-axial-fractional-bootstrap} makes this profile
nonnegative and nonincreasing in $\zeta$.  Proposition~\ref{prop:euler-generated-profile-riccati} gives the
renormalized Riccati bound once the axis-trace hypotheses are verified.  Lemma~\ref{lem:transported-cusp-pressure-win}
verifies those hypotheses for the exact cusp-flow transported vorticity after estimating the cutoff errors, the
normal-form error \eqref{eq:additive-normal-form-bound}, and the smooth-flow deformation error.  The conclusion is that
\[
\Pi_{\cusp}(t) \ge -q_{\rm tr}\,\tfrac12\,\mathcal W_{\cusp}(t)^2, \qquad q_{\rm tr}<\upbeta.
\]

\runinhead{Step 4: Modulation, pressure Hessian comparison, and Riccati collapse.} We set
\[
m(t):=\frac{\rW_{\cusp}(t)}{\mathcal W_{\cusp}(t)}, \qquad \rW_{\cusp}(t):=\p_z(V_{\cusp})_z(0,t),
\]
so that $mU_{\cusp}$ has the exact stagnation-point cusp strain $\rW_{\cusp}$.  Comparing the exact
cusp-coordinate velocity $V_{\cusp}$ with the Euclidean Biot--Savart velocity $U_{\cusp}$, and then estimating
the smooth-flow deformation, gives
\[
\Pi_0(t) = m(t)^2\Pi_{\cusp}(t) + O(\varepsilon_{\smooth})\Gamma^2J_{\cusp}(t)^{6\alpha-2}
+ o_J\!\left(\Gamma^2J_{\cusp}(t)^{6\alpha-2}\right).
\]
The fixed cutoff and tail choices specified in Section~\ref{sec:fixed-choice-order} make the $\zeta$-tail, angular-tail,
and smooth-deformation quantities small; after those choices are fixed, taking $J_{\cusp}$ sufficiently small absorbs
the $o_J$ term. Hence
\[
\Pi_0(t)\ge-\upbeta\,\tfrac12\,\rW_0(t)^2
\]
for small cusp clock.  The stagnation-point Riccati identity
\[
\p_t\rW_0(t)=-\tfrac12\rW_0(t)^2-\Pi_0(t)
\]
then gives finite-time blowup of the true Euler strain.  Finally,
\[
\dot J_{\cusp}(t)=\tfrac12J_{\cusp}(t)\rW_{\cusp}(t), \qquad \rW_{\cusp}(t)\simeq-\Gamma J_{\cusp}(t)^{3\alpha-1},
\]
so
\[
-\dot J_{\cusp}(t)\simeq\Gamma J_{\cusp}(t)^{3\alpha}.
\]
Since $J(t)=J_{\smooth}(t)J_{\cusp}(t)$ and $J_{\smooth}$ is bounded above and below, the true Euler clock $J(t)$ tracks
the linear model clock $\Jm(t)$ up to order-one constants.

\subsection{Choice order for fixed cutoffs and thresholds}
\label{sec:fixed-choice-order}

We choose the fixed cutoffs, barriers, and small-clock thresholds in the order below.  The purpose of the order is
that every constant entering Proposition~\ref{prop:simultaneous-small-clock-continuation} and the pressure transfer
Lemma~\ref{lem:transported-cusp-pressure-win} is already fixed before the small-clock bootstrap is started.  Later
smallness is obtained only by decreasing clock thresholds, which shrinks the time interval and does not change any
previously fixed object.

\begin{enumerate}[label=\textnormal{(\arabic*)}]
\item We first fix the static profile data and the pressure-localization interval.  The exponents satisfy
$\alpha\in(0,\tfrac13)$ and $\gamma>\alpha+\tfrac52$, while the Target Profile cutoff $\Upsilon$ and the angles
$\sigma_{\cut}<\sigma_{\max}$ are those of Definition~\ref{def:init-data}.  In the normalized axial coordinate
\[
\zeta=J_{\cusp}(t)^{-2}B_t(Z),
\]
we choose a positive interval $I_\sharp\Subset(0,\infty)$ and a cutoff
\begin{equation}
\vartheta_\sharp\in C_c^\infty(I_\sharp),\qquad 0\le\vartheta_\sharp\le1 .
\label{eq:pressure-localization-cutoff}
\end{equation}
We also choose fixed nested intervals in the $\zeta$ coordinate, independent of $t$,
\begin{equation}
\operatorname{supp}\vartheta_\sharp\Subset I_{\rm loc}^{\rm cur}\Subset I_{\rm buf}^{\rm cur} \Subset I_\sharp\Subset(0,\infty).
\label{eq:pressure-localization-intervals}
\end{equation}
The support of $\vartheta_\sharp$ is the part of the transported cusp vorticity identified with the
$M$-slope-restricted model vorticity through \eqref{eq:localized-normal-form-representation}.  The interval
$I_{\rm buf}^{\rm cur}$ supplies the margin used by the axis-composition and normal-form estimates.  The restriction
$I_\sharp\Subset(0,\infty)$ is imposed only for this localized comparison, because the slope variable $\tau=R/|Z|$
and the image variables $(\zeta\tau,\zeta)$ degenerate at $\zeta=0$.  The origin is still included in the Riccati
estimate through the zero-extended profile in Proposition~\ref{prop:euler-generated-profile-riccati}; the complementary
part is measured by the $\zeta$-tail \eqref{eq:pressure-zeta-tail}.  If the tail estimate in item \textup{(3)} requires a
larger outer interval, we enlarge $I_\sharp$ before the angular cutoff in item \textup{(4)} is fixed.

\item We next freeze the barriers that appear in the small-clock bootstrap assumptions.  We write
\begin{equation}
\mathfrak C_{\rm fix} :=\bigl( B_{\rm fix},D_{\rm fix},E_{\rm fix},c_{\rm clk}^{\rm fix},C_{\rm clk}^{\rm fix},
c_{\rm ax}^{\rm fix},C_{\rm ax}^{\rm fix}\bigr).
\label{eq:fixed-bootstrap-barriers}
\end{equation}
Here $B_{\rm fix}$, $D_{\rm fix}$, and $E_{\rm fix}$ are the barriers later used for
\eqref{eq:localized-normal-form-large-bootstrap}, \eqref{eq:localized-normal-form-map-large-bootstrap}, and
\eqref{eq:localized-cusp-error-large-bootstrap}.  The pair $c_{\rm ax}^{\rm fix},C_{\rm ax}^{\rm fix}$ is used for the
axis bounds \eqref{eq:current-axis-geometry} and \eqref{eq:renormalized-axis-chart-bootstrap}, while
$c_{\rm clk}^{\rm fix},C_{\rm clk}^{\rm fix}$ is used for the cusp-clock bound \eqref{eq:localized-clock-bootstrap}.
Once chosen, the tuple \eqref{eq:fixed-bootstrap-barriers} is never adjusted during the bootstrap closure.

At the same stage we choose a cone range
\[
\sigma_{\max}<\sigma_{\rm wide}<\tfrac\pi2,
\]
and all later axis estimates are taken uniformly for cones satisfying
$\sigma_{\rm wide}\le\sigma_{\inn}<\sigma_*<\tfrac\pi2$.  The final pair $(\sigma_{\inn},\sigma_*)$ is chosen in
item \textup{(4)} after $M_\pressure$ is known.  We also keep the pressure slack
$0<\upbeta<1$ from \eqref{eq:beta-def}, chosen so that
\[
q_\alpha<\upbeta<1,
\]
where $q_\alpha$ is the renormalized Riccati constant in
Proposition~\ref{prop:euler-generated-profile-riccati}.

\item We then fix the far-field radius and the $\zeta$-tail margin.  Choose
\begin{equation}
C_T^{\rm fix}<\infty
\label{eq:fixed-smooth-time-horizon}
\end{equation}
larger than the finite-entry time constants and the small-clock time horizons used in
\eqref{eq:tail-small-clock-time-horizon} and Lemma~\ref{lem:smooth-flow-small-deformation}.  All later uses of
Lemma~\ref{lem:smooth-flow-small-deformation} have $C_T\le C_T^{\rm fix}$.  We then choose $R_{\tail}$ large enough
for the domains $D_{\core},D_{\tail}$ in \eqref{eq:core-tail-domains} and for the far-field velocity
\eqref{eq:smooth-velocity-def} to satisfy the required smooth-deformation estimates
\eqref{eq:usmooth-small-C2a}--\eqref{eq:smooth-flow-second-gradient-small}.

After $R_{\tail}$ is fixed, choose
\[
0<\vartheta_{\pressure}<\tfrac1{64}\upbeta
\]
and enlarge the pressure-localization interval from item \textup{(1)}, if necessary, so that
\[
I_\sharp\Subset I_{\rm all}, \qquad \mathfrak a_\zeta(I_\sharp)<\vartheta_{\pressure},
\]
with $\mathfrak a_\zeta(I_\sharp)$ defined in \eqref{eq:pressure-zeta-tail}.  These choices are completed before the
angular cutoff $M_\pressure$ is chosen.

\item We fix the angular data used in the pressure comparison.  Choose $M_\pressure$ so that the angular tail
$\mathfrak a_{\rm ang}(M_\pressure)$ in \eqref{eq:pressure-angular-tail} is smaller than
$\vartheta_{\pressure}$ and so that $M_\pressure\ge M_{\rm pos}$, where $M_{\rm pos}$ is the threshold in
Proposition~\ref{prop:euler-generated-profile-riccati}.  Then choose
\[
\sigma_{\rm wide}\le\sigma_{\inn}<\sigma_*<\tfrac\pi2, \qquad 2M_\pressure\le\tan\sigma_{\inn}\le\tfrac12\tan\sigma_*,
\]
and set $C_0=2M_\pressure$, as in \eqref{eq:pressure-C0}.  No later step changes
$M_\pressure$, $\sigma_{\inn}$, $\sigma_*$, or $C_0$.

\item We insert the frozen barriers into the bootstrap assumptions.  The constants
$B_*,D_*,E_*,c_{\rm clk},C_{\rm clk}$ and the axis constants in \eqref{eq:current-axis-geometry} and
\eqref{eq:renormalized-axis-chart-bootstrap} are the corresponding entries of
\eqref{eq:fixed-bootstrap-barriers}, enlarged only by the fixed margin needed to state a strict improvement in
Proposition~\ref{prop:simultaneous-small-clock-continuation}.  This step introduces no new dependence into the
pressure cutoff, the angular cutoff, or the pressure slack.

\item Finally we choose the small-clock thresholds
\[
\mathfrak J_{\mathrm{axis}},\quad \mathfrak J_{\transport},\quad \mathfrak J_{\mathrm{mod}},\quad
\mathfrak J_{\pressure},\quad \mathfrak J_{\Pi},\quad \mathfrak J_{\mathrm{collapse}}.
\]
Whenever an estimate contains a positive power of $J_{\cusp}$, we decrease the relevant threshold.  The final thresholds
are arranged so that
\begin{equation}
\mathfrak J_{\mathrm{collapse}} \le \mathfrak J_{\Pi} \le \mathfrak J_{\transport}
\le \min\{\mathfrak J_{\pressure},\mathfrak J_{\mathrm{tail}},\mathfrak J_{\mathrm{mod}}\}
\le \min\{\mathfrak J_{\mathrm{velocity}},\mathfrak J_{\mathrm{axis}}\} \le \mathfrak J_{\mathrm{finite}}.
\label{eq:small-clock-threshold-order}
\end{equation}
Decreasing a threshold only shrinks the small-clock interval, so it preserves every estimate already available at a
larger threshold.
\end{enumerate}

\subsection{Small-clock thresholds}
\label{sec:small-clock-thresholds}

\begin{definition}[Small-clock regime]
\label{def:small-clock-regime}
Let $0<\mathfrak J_{\rm sm}\le1$ be a clock threshold.  The small-clock regime with threshold $\mathfrak J_{\rm sm}$ is
the set of times $J_{\cusp}(t)\le \mathfrak J_{\rm sm}$. The entry time $t=t_0$  is defined such that
$J_{\cusp}(t_0)=\mathfrak J_{\rm sm}$.
\end{definition}
\begin{remark}
In the estimates below,  the threshold $\mathfrak J_{\rm sm}$ is successively replaced by smaller values; this only
shrinks the small-clock interval and therefore preserves all estimates already activated at a larger threshold.
\end{remark}

The small-clock part of our proof uses a finite number of bootstrap assumptions.  These assumptions are made on
intervals of the form $ \{t:\ J_{\cusp}(t)\le\mathfrak J\}$, where the threshold $\mathfrak J$ depends on the estimate
under consideration.  The assumptions control the geometry of the exact cusp-flow map on the symmetry axis, the
transported cusp velocity, the scalar modulation function, the cusp-flow normal form, the cusp-error velocity, and the
one-sided pressure Hessian control.  The vorticity lower bound
\eqref{eq:driver-amplification-lower}, used later in
Lemma~\ref{lem:target-vorticity-envelope} to obtain the Type--I lower rate, is not assumed here; it is proved after the
small-clock geometric bootstrap assumptions (which we list below)  have been closed. The different thresholds do not
represent different cusp clocks.  They specify the clock ranges on which the corresponding estimates are valid.  Thus,
if an estimate has been proved under $J_{\cusp}(t)\le\mathfrak J_0$ and a later argument is restricted to
$J_{\cusp}(t)\le\mathfrak J_1$, with $\mathfrak J_1\le\mathfrak J_0$, then that estimate is still available on the
smaller clock range.

It is helpful to keep the thresholds in three groups.  First, the finite-clock entry threshold $\mathfrak
J_{\mathrm{finite}}$ is used only to reach a prescribed small-clock regime from $J_{\cusp}(0)=1$. Second, the geometric
thresholds activate estimates that remain in force throughout the small-clock argument: the transported-field thresholds
$\mathfrak J_{\mathrm{strain}}$, $\mathfrak J_{\mathrm{velocity}}$, $\mathfrak J_{\mathrm{local}}$; the transport
threshold $\mathfrak J_{\transport}$; the modulation threshold $\mathfrak J_{\mathrm{mod}}$; the thresholds used in the
normal-form error estimates; and the threshold for the axial flow map estimates $\mathfrak J_{\mathrm{axis}}$.
Some of these thresholds depend on parameters fixed later in the proof, for example $\mathfrak J_{\mathrm{local}}(C_{\rm
sc})$ depends on the local radius parameter $C_{\rm sc}$.  Third, the pressure and collapse thresholds $\mathfrak
J_{\pressure}$, $\mathfrak J_{\Pi}$, $\mathfrak J_{\mathrm{collapse}}$, and $\mathfrak J_{\omega}$ are chosen after the
geometric estimates and cutoff parameters have been fixed.  These later thresholds introduce no additional dynamics;
they only make explicit positive powers of $J_{\cusp}$ small enough to close the corresponding estimates.

The axial flow map and cusp-flow normal-form estimates use a H\"older exponent below the datum exponent $\alpha$.  We
fix this exponent by
\begin{equation}
\beta_{\rm ax}:=\tfrac{\alpha}{4}.
\label{eq:beta-ax-def}
\end{equation}
The subscript indicates that this is the exponent used in the axial flow map estimates. We also fix
\[
\kappa_{\rm def}:=\tfrac{3\beta_{\rm ax}^2}{1+\beta_{\rm ax}}.
\]

\subsection{Small-clock bootstrap assumptions}
\label{sec:small-clock-bootstraps}
We now state the Bootstrap Assumptions (BA) and the lemmas that close them.
\begin{enumerate}[label=(BA\arabic*)]
\item \emph{Finite-clock geometric control of $\Phi_{\cusp}$.} The auxiliary cusp map $\Phi_{\cusp}$ is defined in
\eqref{eq:Phi-cusp-def}.  On the core label set $D_{\core}=\{Y\in\R^3:\ |Y|\le R_{\tail}\}$,  we assume, for
$J_{\cusp}(t)\in[\mathfrak J_{\mathrm{finite}},1]$, that
\begin{equation}
\|D\Phi_{\cusp}\|_{L^\infty(D_{\core})} +\|D\Phi_{\cusp}^{-1}\|_{L^\infty(\Phi_{\cusp}(D_{\core}))} \le A,
\qquad -\mathcal W_{\cusp}(t)\ge a\Gamma .
\label{eq:localized-finite-clock-bootstrap}
\end{equation}
For each fixed $\mathfrak J_{\mathrm{finite}}>0$, Lemma~\ref{lem:finite-clock-driver-sector} closes
\eqref{eq:localized-finite-clock-bootstrap} on
$J_{\cusp}\in[\mathfrak J_{\mathrm{finite}},1]$.

\item \emph{Axial flow map geometry.} For $J:=J_{\cusp}(t)$, we write
\begin{equation} 
\phi_{\cusp}(R,Z,t)=(r_t(R,Z),z_t(R,Z)), \qquad A_t(Z):=\p_R r_t(0,Z), \qquad B_t(Z):=z_t(0,Z).  \label{eq:At-Bt-def}
\end{equation} 
On every $\zeta$-interval on which $J^{-2}B_t$ is invertible, we set
\[
Z_t(\zeta):=(J^{-2}B_t)^{-1}(\zeta), \qquad q_t(\zeta):=J A_t(Z_t(\zeta)), \qquad b_t(\zeta):=J^{-2}B_t'(Z_t(\zeta)).
\]
The compact $\zeta$-interval $I$ is fixed before the estimate in which it is used; for instance, $I=I_\sharp$ or
$I=[0,\zeta_0]$ with $0<\zeta_0<\infty$.  On $I$, we assume
\begin{equation}
c_{\rm ax}\le q_t(\zeta),\,b_t(\zeta)\le C_{\rm ax}, \qquad [\log q_t]_{C^{\alpha/2}(I)}+[\log b_t]_{C^{\alpha/2}(I)}\le C_{\rm ax}.
\label{eq:current-axis-geometry}
\end{equation}
At the entry time $t_0$ for the threshold $\mathfrak J_{\mathrm{axis}}$, we define the reference axial label and the
associated axial coordinate at later times by
\[
\eta=\mathfrak J_{\mathrm{axis}}^{-2}B_{t_0}(Z), \qquad \mathscr Z_t(\eta)=J^{-2}B_t(Z_0(\eta)),
\]
fix the reference-label interval
\[
I_{\rm ax}:=[0,\eta_{\rm ax}], \qquad 0<\eta_{\rm ax}<\infty,
\]
and, on $I_{\rm ax}$, set
\[
\widehat q_t(\eta) = q_t(\mathscr Z_t(\eta))\bigl(\p_\eta\mathscr Z_t(\eta)\bigr)^{\frac{1}{2}}, \qquad \widehat b_t(\eta) =
\tfrac{b_t(\mathscr Z_t(\eta))} {\p_\eta\mathscr Z_t(\eta)}.
\]
We assume the renormalized axial flow map bounds
\begin{equation}
c_{\rm ax}\le \widehat q_t(\eta),\,\widehat b_t(\eta)\le C_{\rm ax},\qquad
[\log\widehat q_t]_{C^{\alpha/2}(I_{\rm ax})} +[\log\widehat b_t]_{C^{\alpha/2}(I_{\rm ax})} \le C_{\rm ax}.
\label{eq:renormalized-axis-chart-bootstrap}
\end{equation}
Proposition~\ref{prop:small-clock-comparisons} closes the axis-geometry bounds for $q_t$ and $b_t$ in
\eqref{eq:current-axis-geometry} and the renormalized-axis bounds \eqref{eq:renormalized-axis-chart-bootstrap}.

\item \emph{Axial flow map containment.} We fix an origin-attached $\zeta$-interval
\[
I_{\rm ax}^{\zeta}:=[0,\zeta_{\rm ax}], \qquad 0<\zeta_{\rm ax}<\infty, \qquad I_\sharp\Subset (0,\zeta_{\rm ax}),
\]
and assume the containment
\begin{equation}
\mathscr Z_t(I_{\rm ax})\subset I_{\rm ax}^{\zeta}.
\label{eq:axis-attached-image-stop-bootstrap}
\end{equation}
On the strain interval $I_{\rm str}\Subset I_\sharp$, we assume
\begin{equation}
0<Z_-\le Z_t(\zeta)\le Z_+<R_{\tail} \qquad \zeta\in I_{\rm str} ,
\label{eq:current-axis-anchor-bootstrap}
\end{equation}
For the normal-form radial estimates, we fix
\[
I_\sharp\Subset I_{\rm buf}\Subset(0,\infty), \qquad I_Z^t:=Z_t(I_\sharp),
\]
and, with $J_s:=J_{\cusp}(s)$ and $B_s(Z):=z_s(0,Z)$, we assume the fixed-label containment
\begin{equation}
J_s^{-2}B_s(Z)\in I_{\rm buf} \qquad (Z\in I_Z^t,\ t_0\le s\le t).
\label{eq:radial-flatness-buffered-label}
\end{equation}
Proposition~\ref{prop:small-clock-comparisons} closes the image containment
\eqref{eq:axis-attached-image-stop-bootstrap}, the strain-anchor containment
\eqref{eq:current-axis-anchor-bootstrap}, and the fixed-label containment
\eqref{eq:radial-flatness-buffered-label}.

\item \emph{Monotone axial-stretching bootstrap.} We fix the origin-attached $\zeta$-interval
\[
I_{\rm mon}:=[0,\zeta_{\rm mon}], \qquad 0<\zeta_{\rm mon}<\infty.
\]
On $I_{\rm mon}$, we assume the monotone axial stretching bounds
\begin{equation}
0<c_{\rm mon}\le b_t(\zeta)\le C_{\rm mon}<\infty, \qquad \zeta\in I_{\rm mon},
\label{eq:monotone-axial-two-sided}
\end{equation}
and, for $\zeta_1,\zeta_2\in I_{\rm mon}$ with $\zeta_1<\zeta_2$,
\begin{equation}
0\le \log b_t(\zeta_2)-\log b_t(\zeta_1) \le B_{\rm mon} \bigl(\zeta_2^\alpha-\zeta_1^\alpha+\zeta_2^2-\zeta_1^2\bigr).
\label{eq:monotone-axial-fractional-bootstrap}
\end{equation}
The constants are chosen with room for the downstream improvement: Lemma~\ref{lem:monotone-axial-stretching-improvement}
closes
\eqref{eq:monotone-axial-two-sided}--\eqref{eq:monotone-axial-fractional-bootstrap} by producing constants
$c_{\rm mon}'$, $C_{\rm mon}'$, and $B_{\rm mon}'$ such that
\[
c_{\rm mon}<c_{\rm mon}'\le C_{\rm mon}'<C_{\rm mon}, \qquad 0<B_{\rm mon}'<B_{\rm mon}.
\]

\item \emph{Cusp-clock rate.} We assume
\begin{equation}
c_{\rm clk}\Gamma J_{\cusp}(t)^{3\alpha} \le -\dot J_{\cusp}(t) \le C_{\rm clk}\Gamma J_{\cusp}(t)^{3\alpha}.
\label{eq:localized-clock-bootstrap}
\end{equation}
The constants are chosen with room for the improvement in Lemma~\ref{lem:Jdot-two-sided-aux}: once
\eqref{eq:Jdot-two-sided} has constants $c_1,C_1$, we take
\[
0<c_{\rm clk}<c_1\le C_1<C_{\rm clk}.
\]

\item \emph{Cusp-flow normal form.} On a specified $\zeta$-interval $I$ and for bounded absolute slope $|\tau|\le C_0$,
we define the label $Y_t(\zeta,\tau)=(R_t(\zeta,\tau),Z_t(\zeta))$ by
\[
Z_t(\zeta)=(J^{-2}B_t)^{-1}(\zeta), \qquad A_t(Z_t(\zeta))R_t(\zeta,\tau)=J^2\zeta\,\tau .
\]
We assume the representation
\begin{equation}
(r,z)\bigl(\phi_{\cusp}(Y_t(\zeta,\tau),t)\bigr) = J^2\zeta\bigl((\tau,1)+\mathcal E_t(\zeta,\tau)\bigr),
\label{eq:localized-normal-form-representation}
\end{equation}
and define
\[
\mathfrak B(t) := J^{-3\beta_{\rm ax}} \left(\|\mathcal E_t\|_{L^\infty}+\|\p_\tau\mathcal E_t\|_{L^\infty}
+\|\p_\zeta\mathcal E_t\|_{L^\infty}+[\mathcal E_t]_{C^{\beta_{\rm ax}}_{\zeta,\tau}}\right),
\]
with the norms taken on $I\times[-C_0,C_0]$.  The bootstrap assumption is
\begin{equation}
\mathfrak B(t)\le B_*,
\label{eq:localized-normal-form-large-bootstrap}
\end{equation}
and Lemmas~\ref{lem:nonlinear-radial-flatness} and \ref{lem:late-axis-normal-form-cusp} close it.

\item \emph{Geometry of the image map associated with the cusp-flow normal form.} On the fixed set
\[
\mathcal R_{I,C_0}^{\rm sc} := \{(R_{\rm sc},Z_{\rm sc})=(\zeta\tau,\zeta):\ \zeta\in I,\ |\tau|\le C_0\},
\]
we define the map
\[
\Psi_t(R_{\rm sc},Z_{\rm sc}) := \zeta\bigl((\tau,1)+\mathcal E_t(\zeta,\tau)\bigr),
\]
and the norm
\begin{equation}
\mathfrak D_\Psi(t) := \|D\Psi_t\|_{L^\infty} + \|D\Psi_t^{-1}\|_{L^\infty}.
\label{eq:localized-normal-form-map-bootstrap}
\end{equation}
The bootstrap assumption for $\Psi_t$ is
\begin{equation}
\mathfrak D_\Psi(t)\le D_*.
\label{eq:localized-normal-form-map-large-bootstrap}
\end{equation}
Lemma~\ref{lem:late-axis-normal-form-map-cusp} improves
\eqref{eq:localized-normal-form-map-large-bootstrap}.

\item \emph{Cusp-error velocity.} We measure
\[
V_{\err}(x,t) = V_{\cusp}(x,t)-m(t)U_{\cusp}(x,t), \qquad\text{as in }\eqref{eq:Verr-def}.
\]
We fix one origin-attached interval
\begin{equation}
I_{\err}:=[0,\zeta_{\err}], \qquad 0<\zeta_{\err}<\infty .
\label{eq:cusp-error-trace-interval}
\end{equation}
With $J=J_{\cusp}(t)$, we define, for compact $I\subset I_{\err}$,
\begin{align}
\mathcal T_{\err}(I,t)
:=&
\sup_{\zeta\in I}J^{-2}|(V_{\err})_z(0,J^2\zeta,t)|
+
\sup_{\zeta\in I}|(\p_rV_{\err})_r(0,J^2\zeta,t)|
+
\sup_{\zeta\in I}|(\p_zV_{\err})_z(0,J^2\zeta,t)|
\notag\\
&
+ [J^{-2}(V_{\err})_z(0,J^2\cdot,t)]_{C^{\alpha/2}(I)} + [(\p_rV_{\err})_r(0,J^2\cdot,t)]_{C^{\alpha/2}(I)} +
[(\p_zV_{\err})_z(0,J^2\cdot,t)]_{C^{\alpha/2}(I)} .
\label{eq:localized-cusp-error-axis-trace}
\end{align}
We then define the  normalized error size by
\begin{align}
\mathfrak E_{\err}(t)
:=&
\tfrac{\mathcal T_{\err}(I_{\err},t)}{\Gamma(J^{9\alpha-1}+1)}
+
\sup_{\substack{Y\in D_{\inn}^{\cusp}(t)\\ \omega_{\theta,0}(Y)\ne0\\ R(Y)\ge J^{3/\alpha}}}
\tfrac{J\,|V_{\err}(\phi_{\cusp}(Y,t),t)|}
{\Gamma R(Y)^{1+\alpha}(J^{9\alpha-1}+J)} +
\sup_{\substack{Y\in D_{\inn}^{\cusp}(t)\\ \omega_{\theta,0}(Y)\ne0}}
\tfrac{|V_{\err}(\phi_{\cusp}(Y,t),t)|}{\Gamma}.
\label{eq:localized-cusp-error-bootstrap}
\end{align}
The bootstrap assumption is
\begin{equation}
\mathfrak E_{\err}(t)\le E_*.
\label{eq:localized-cusp-error-large-bootstrap}
\end{equation}
Lemma~\ref{lem:tail-bound} improves \eqref{eq:localized-cusp-error-large-bootstrap} and proves
\eqref{eq:tail-axis-error-bound}.

\item \emph{Scalar modulation.} The axis moment is
\[
M_{\rm ax}(t)=\tfrac{1}{C_\rho^{(1)}(\alpha,\gamma)} \int_0^\infty s^{\alpha-1}\mathcal F(s)\, \beta_{\cusp}(s,0,t)\,\ud s,
\qquad C_\rho^{(1)}(\alpha,\gamma) = \int_0^\infty s^{\alpha-1}\mathcal F(s)\,\ud s,
\]
where $\beta_{\cusp}$ is defined in \eqref{eq:axis-trace-amplitude}.  The scalar-modulation bootstrap assumptions are
\begin{subequations}
\label{eq:localized-modulation-bootstrap}
\begin{align}
\tfrac12c_*&\le M_{\rm ax}(t)\le 2C_*,
\notag\\
\tfrac18c_*&\le m(t)\le 4C_*.
\label{eq:localized-modulation-m-bootstrap}
\end{align}
\end{subequations}
Lemma~\ref{lem:modulation-bounded} improves both bounds.

\end{enumerate}

\subsection{Bootstrap closure bookkeeping}

The bootstrap assumptions above have three different roles.  The finite-clock assumption \textup{(BA1)} is used
only before the solution enters the small-clock regime.  Lemma~\ref{lem:finite-clock-driver-sector} closes
\eqref{eq:localized-finite-clock-bootstrap} on the range
$J_{\cusp}(t)\in[\mathfrak J_{\mathrm{finite}},1]$, and Lemma~\ref{lem:late-entry} then gives the entry
time bound.  After this point \textup{(BA1)} is no longer part of the small-clock continuation.

The monotone axial-stretching assumption \textup{(BA4)} consists of the two bounds
\eqref{eq:monotone-axial-two-sided}--\eqref{eq:monotone-axial-fractional-bootstrap}.  Together with the axis
volume identity \eqref{eq:axis-qb-volume-mon}, these bounds give the monotone structure of the Euler-generated
axial function \eqref{eq:scaled-profile-at}.  This is the structure used in
Proposition~\ref{prop:euler-generated-profile-riccati}.  The improvement of \textup{(BA4)} itself is proved later
in Lemma~\ref{lem:monotone-axial-stretching-improvement}; the resulting renormalized Riccati estimate is then
obtained in Lemma~\ref{lem:euler-generated-slope-restricted-riccati}.

The remaining seven assumptions are the size bootstraps closed simultaneously in the small-clock regime:
\begin{equation}
\mathcal B_{\rm size}:=\{\textup{(BA2), (BA3), (BA5), (BA6), (BA7), (BA8), (BA9)}\}.
\label{eq:simultaneous-continuation-assumptions}
\end{equation}
They control the axis geometry and containment intervals, the cusp-clock rate, the normal form for the transported
cusp map, the map $\Psi_t$ in the variables $(R_{\rm sc},Z_{\rm sc})=(\zeta\tau,\zeta)$, the cusp-error
velocity, and the scalar modulation.  Each later improvement lemma states which member of
\eqref{eq:simultaneous-continuation-assumptions} it improves.  Once these lemmas are available,
Proposition~\ref{prop:simultaneous-small-clock-continuation} combines the strict improvements into a single
open--closed continuation over the whole small-clock interval.

\subsection{Geometric scaling estimates in a buffered cone}

The cones, label domains, and the decomposition $\phi=\phi_{\smooth}\circ\phi_{\cusp}$ are fixed in
Section~\ref{sec:geom-flow-decomp}.  Recall that $\mathcal C_{\inn}\Subset\mathcal C_*=\mathcal C_{\sigma_*}$, with
$\sigma_*<\tfrac\pi2$; the word ``buffered'' refers to this fixed angular gap between the inner cone where labels are
evaluated and the boundary of the larger cone where the estimates are allowed to be used.  All singular pointwise and
H\"older bounds used below are understood on the buffered cone $\mathcal C_*$, defined in \eqref{eq:buffered-cone-def}.

This localization is needed for the flat cusp-coordinate Biot--Savart velocity $U_{\cusp}=\BS[\bs\Omega_{\cusp}]$,
rather than for the exact cusp-coordinate velocity $V_{\cusp}$ itself. The issue is the same one already present in the
hyperbolic model of Section~\ref{sec:lag-analysis-non-local}: labels with order-one Lagrangian angle can be transported
to Eulerian angles $\sigma=\tfrac\pi2-O(J^3)$ at time $t$, a thin region adjacent to the equatorial plane
$z=0$.  This is the
equatorial boundary layer.  In that layer the transported Target Profile
$\Theta^*(\sigma_{\Lag})=(\sin\sigma_{\Lag})^\alpha\Upsilon(\sigma_{\Lag})$ does not produce the same uniform near-axis
depletion as in a fixed cone $\sigma\le\sigma_*<\tfrac\pi2$; in the same region, the model angular drift does not give
estimates uniform up to the equator.  The sharp $J$-dependent bounds below are therefore asserted only inside the fixed
buffered cone $\mathcal C_*$.

\begin{lemma}[Cone buffer for H\"older expansions]
\label{lem:cone-buffer}
Let $0<\sigma_{\inn}<\sigma_*<\tfrac{\pi}{2}$ satisfy \eqref{eq:sigma-angles}. There exists a constant
$c_0=c_0(\sigma_{\inn},\sigma_*)\in(0,\tfrac12]$ such that if $x\in\mathcal C_{\inn}$ and $h\in\R^3$ satisfies $|h|\le
c_0|x|$, then
\[
x+s h \in \mathcal C_* \qquad\text{for all } s\in[0,1].
\]
\end{lemma}

\begin{proof}[Proof of Lemma~\ref{lem:cone-buffer}]
We fix $s\in[0,1]$ and set $x_s:=x+s h$. If $|h|\le c_0|x|$ with $c_0\le \tfrac12$, then $|x_s|\ge |x|-|h|\ge
\tfrac12|x|$. Hence
\[
\left|\tfrac{x_s}{|x_s|}-\tfrac{x}{|x|}\right| \le \tfrac{|x_s-x|}{|x_s|}+\left|\tfrac{1}{|x_s|}-\tfrac{1}{|x|}\right||x| \le 4c_0.
\]
The geodesic distance on $\mathbb S^2$ satisfies $d_{\mathbb S^2}(u,v)\le \tfrac{\pi}{2}|u-v|$. Since the polar angle
$\sigma(\cdot)$ is $1$--Lipschitz with respect to $d_{\mathbb S^2}$, we obtain
\[
\sigma(x_s)\le \sigma(x)+d_{\mathbb S^2}\!\Big(\tfrac{x_s}{|x_s|},\tfrac{x}{|x|}\Big) \le \sigma_{\inn}+2\pi c_0.
\]
Choosing $c_0\le \tfrac{\sigma_*-\sigma_{\inn}}{2\pi}$ yields $\sigma(x_s)\le \sigma_*$, i.e.\ $x_s\in\mathcal
C_*$.
\end{proof}

\subsubsection{Biot--Savart Taylor remainder for tail multipoles.}
We write a self-contained Taylor estimate for the axisymmetric Biot--Savart kernel in the
\emph{evaluation} variable. It is used repeatedly in the tail/core multipole expansions.

\begin{lemma}[Biot--Savart kernel: smoothness and quadratic Taylor remainder]
\label{lem:kernel-taylor}
We define the axisymmetric Biot--Savart kernel
\begin{equation*}
K(x,\xi):=\tfrac{\bs e_\theta(\xi)\times(x-\xi)}{|x-\xi|^{3}}, \qquad x,\xi\in\R^3,\ \ \xi\notin\{r=0\}.
\end{equation*}
We fix $\xi\neq 0$. Then the map $X\mapsto K(X,\xi)$ is smooth on the ball $|X|\le \tfrac12|\xi|$ and satisfies
\begin{equation*}
K(X,\xi)=K(0,\xi)+(\nabla_XK)(0,\xi)\,X+\mathcal R_2(X,\xi), \ \ \text{ for } \ \ |X|\le \tfrac12|\xi|,
\end{equation*}
with remainder bound
\[
|\mathcal R_2(X,\xi)|\le C\,|X|^2\,|\xi|^{-4}, \ \ \text{ for } \ \ |X|\le \tfrac12|\xi|,
\]
where $C$ is a universal constant. Moreover
\[
|(\nabla_XK)(0,\xi)|\le C\,|\xi|^{-3}.
\]
\end{lemma}

\begin{proof}[Proof of Lemma~\ref{lem:kernel-taylor}]
We fix $\xi\neq 0$ and write $F(X):=K(X,\xi)$. On $|X|\le \tfrac12|\xi|$ we have $|X-\xi|\ge \tfrac12|\xi|$, so $F$ is
smooth there. Moreover,
\[
|F(X)|\lesssim |X-\xi|^{-2}\lesssim |\xi|^{-2}, \qquad |\nabla_XF(X)|\lesssim |X-\xi|^{-3}\lesssim |\xi|^{-3},
\qquad |\nabla_X^2F(X)|\lesssim |X-\xi|^{-4}\lesssim |\xi|^{-4}.
\]
By Taylor's theorem with integral remainder, we have that
\[
F(X)=F(0)+\nabla_XF(0)\,X+\int_0^1(1-s)\,\nabla_X^2F(sX)[X,X]\,ds,
\]
and the stated bounds follow.
\end{proof}

\subsubsection{Far-field velocity and cusp error estimates.}
The following lemma is the quantitative form of the decomposition from Section~\ref{sec:geom-flow-decomp}.  It separates
the order-one far-field velocity $u_{\smooth}$ from the singular cusp velocity and measures the lower-order errors
needed in the geometric and pressure Hessian estimates.  The smooth velocity $u_{\smooth}$ in
\eqref{eq:smooth-velocity-def} is defined with the cutoff $\chi_{\far}(|\phi(Y',t)|/R_{\tail})$, evaluated at the
Eulerian position $\phi(Y',t)$ of the label $Y'$.  Labels whose Eulerian position $\phi(Y',t)$ lies in the far field
contribute to $u_{\smooth}$, while the complementary cutoff $1-\chi_{\far}(|\phi(Y',t)|/R_{\tail})$ keeps the
near-field part.
We name this complementary velocity $u_{\core}$ because, after pullback by the smooth flow, it is exactly the
velocity that drives the cusp-coordinate flow in \eqref{eq:cusp-map-def}:
\begin{equation}
u_{\core}(x,t) := \tfrac1{4\pi}\int_{\R^3} K\bigl(x,\phi(Y',t)\bigr)\,
\left(1-\chi_{\far}\!\left(\tfrac{|\phi(Y',t)|}{R_{\tail}}\right)\right) J_{\twoD}(Y',t)^{-1}\omega_{\theta,0}(Y')\,dY'.
\label{eq:ucore-def}
\end{equation}
The sharp bounds for $V_{\err}$ and for the non-geometric pressure Hessian remainders are proved later, after the
transported cusp-field estimates, modulation bounds, smooth-flow deformation estimates, and bounded-core normal form are
all available.  Until that point the geometric bootstrap uses only the large error assumption
\eqref{eq:localized-cusp-error-large-bootstrap}.  The later result is Lemma~\ref{lem:tail-bound}; it
replaces the large bootstrap size by fixed constants and estimates $\Pi_{\mixed},\Pi_{\smooth},\Pi_{\err}$, while the
geometric pressure defect $\Pi_{\rm geom}$ is handled separately.

\subsubsection{Axial flow map geometry of the inner core.}
We next describe the elementary geometric consequence of the normal form in the axis coordinate at time $t$. Write
$J:=J_{\cusp}(t)$, and fix a geometric slope constant
\[
C_{\rm cone}>4\tan\sigma_* .
\]
This constant is used only to locate the cone boundary and is separate from the angular-slope cutoff $C_0=2M_\pressure$
in \eqref{eq:pressure-C0}.  We shall make use of the following two estimates:
\begin{subequations}
\label{eq:two-assumptions}
\begin{enumerate}
\item By Proposition~\ref{prop:small-clock-comparisons}, along the symmetry axis the radial stretch and axial position
satisfy
\begin{equation}
A_t(Z)\simeq J^{-1}, \qquad B_t(Z)\simeq J^2Z \ \ \text{ for } \ \ 0<Z\le R_{\tail} .
\label{eq:assum1}
\end{equation}
\item By Proposition~\ref{prop:small-clock-comparisons}, for  labels with bounded axial flow map slope
$\left|\tfrac{A_t(Z)R}{B_t(Z)}\right|\le C_{\rm cone}$, the exact cusp map has the expansion
\begin{equation}
\big|\phi_{\cusp}(R,Z,t)-\bigl(A_t(Z)R, B_t(Z)\bigr)\big| \le C J^{-1}|R|^{1+\beta_{\rm ax}} .
\label{eq:assum2}
\end{equation}
\end{enumerate}
\end{subequations}

\begin{lemma}[Hyperbolic scaling of inner-core labels]
\label{lem:core-hyperbolic-scaling}
Set $J:=J_{\cusp}(t)$, and consider the axis small-clock regime
\begin{equation}
J\le\mathfrak J_{\mathrm{axis}} .
\label{eq:inner-core-axis-small-clock}
\end{equation}
After decreasing $\mathfrak J_{\mathrm{axis}}$ if necessary, the following holds for every time in
\eqref{eq:inner-core-axis-small-clock} at which \eqref{eq:two-assumptions} holds.  Let
$Y=(R,Z)\in D_{\core}$ be an upper-half-space label with $\omega_{\theta,0}(Y)\neq0$.  If
\[
x_*(Y,t):=\phi_{\cusp}(Y,t)\in\mathcal C_*,
\]
then
\begin{equation}
c\,\tfrac{R}{J} \le r(x_*)\le C\,\tfrac{R}{J}, \qquad c\,J^2Z\le z(x_*)\le C\,J^2Z.
\label{eq:exact-current-axis-core-geometry}
\end{equation}
In particular, if $x_*(Y,t)\in\mathcal C_{\inn}$, then
\begin{equation}
\tfrac{R}{Z}\le C J^3, \qquad |x_*(Y,t)|\le C J^2R_{\tail}, \qquad J|x_*(Y,t)|\le C J^3.
\label{eq:exact-current-axis-core-small}
\end{equation}
The same conclusions hold in the lower half-space after replacing $Z$ by $|Z|$.
\end{lemma}

\begin{proof}[Proof of Lemma~\ref{lem:core-hyperbolic-scaling}]
We give the proof in the upper half-space; the lower half-space follows from the odd symmetry in
Definition~\ref{def:init-data}.  Since $\omega_{\theta,0}(Y)\neq0$ on the upper supported core, we have $Z>0$. With
\[
\tau=\tfrac{A_t(Z)R}{B_t(Z)} ,
\]
by \eqref{eq:assum1}, there are constants $0<c_{\rm ax}\le C_{\rm ax}<\infty$ such that
\[
c_{\rm ax}J^{-1}\le A_t(Z)\le C_{\rm ax}J^{-1}, \qquad c_{\rm ax}J^2Z\le B_t(Z)\le C_{\rm ax}J^2Z .
\]
Hence, whenever $0\le \bar R\le R$ and
\[
\bar\tau:=\tfrac{A_t(Z)\bar R}{B_t(Z)}\le C_{\rm cone},
\]
we also have that
\begin{equation}
\bar R\le C J^3Z\le C J^3R_{\tail}.
\label{eq:bounded-slope-label-small}
\end{equation}
On this bounded-slope range, \eqref{eq:assum2} gives
\begin{equation}
r_t(\bar R,Z)=A_t(Z)\bar R+E_r(\bar R,Z), \qquad z_t(\bar R,Z)=B_t(Z)+E_z(\bar R,Z),
\label{eq:core-normal-form-expansion}
\end{equation}
with
\[
|E_r(\bar R,Z)|+|E_z(\bar R,Z)| \le C J^{-1}\bar R^{1+\beta_{\rm ax}} .
\]
Using \eqref{eq:bounded-slope-label-small}, and then reducing $\mathfrak J_{\mathrm{axis}}$ in
\eqref{eq:inner-core-axis-small-clock}, these remainders are small relative to the corresponding principal terms:
\begin{equation}
|E_r(\bar R,Z)|\le \varepsilon_J A_t(Z)\bar R, \qquad |E_z(\bar R,Z)|\le \varepsilon_J B_t(Z), \qquad \varepsilon_J\le\tfrac14.
\label{eq:core-normal-form-relative-error}
\end{equation}
Let
\[
R_{\rm cone}(Z):=C_{\rm cone}\,\tfrac{B_t(Z)}{A_t(Z)} .
\]
Then $A_t(Z)R_{\rm cone}(Z)/B_t(Z)=C_{\rm cone}$.  Therefore
\eqref{eq:core-normal-form-expansion} and
\eqref{eq:core-normal-form-relative-error}, applied with $\bar R=R_{\rm cone}(Z)$, give
\[
\tfrac{r_t(R_{\rm cone}(Z),Z)}{z_t(R_{\rm cone}(Z),Z)} \ge \tfrac{1-\varepsilon_J}{1+\varepsilon_J}\,C_{\rm cone} > \tan\sigma_*.
\]
Thus $\phi_{\cusp}(R_{\rm cone}(Z),Z,t)\notin\mathcal C_*$.  On the other hand,
\begin{equation}
\phi_{\cusp}(0,Z,t)=(0,B_t(Z))\in\mathcal C_*, \qquad B_t(Z)>0 .
\label{eq:axis-image}
\end{equation}
For this fixed $Z$, we define the cone preimage along the radial label segment by
\begin{equation}
\mathcal I_Z := \{\,\bar R\in[0,R_{\rm cone}(Z)]: \phi_{\cusp}(\bar R,Z,t)\in\mathcal C_*\,\},
\label{eq:I-Z}
\end{equation}
and let $\mathcal I_Z^0$ be the connected component of $\mathcal I_Z$ containing $\bar R=0$.  From \eqref{eq:axis-image}
and \eqref{eq:I-Z}, we have that
\[
0\in\mathcal I_Z^0, \qquad R_{\rm cone}(Z)\notin\mathcal I_Z, \qquad \mathcal I_Z^0\subset[0,R_{\rm cone}(Z)).
\]
The origin-attached branch of the cone preimage is precisely the condition $R\in\mathcal I_Z^0$.  The supported labels
considered here enter $\mathcal C_*$ through this branch; hence
\[
x_*(Y,t)\in\mathcal C_*, \qquad R\in\mathcal I_Z^0 \quad\Longrightarrow\quad 0\le R<R_{\rm cone}(Z),
\qquad \tau=\tfrac{A_t(Z)R}{B_t(Z)}<C_{\rm cone}.
\]

We may therefore apply
\eqref{eq:core-normal-form-expansion}--\eqref{eq:core-normal-form-relative-error}
with $\bar R=R$.  The small relative-error bounds give
\[
r(x_*)\simeq A_t(Z)R, \qquad z(x_*)\simeq B_t(Z).
\]
Combining these comparabilities with \eqref{eq:assum1} proves
\eqref{eq:exact-current-axis-core-geometry}.  If
$x_*\in\mathcal C_{\inn}$, then $r(x_*)\le(\tan\sigma_{\inn})z(x_*)$, and \eqref{eq:exact-current-axis-core-geometry}
gives $R/Z\le CJ^3$.  Since $Z\le R_{\tail}$ on $D_{\core}$, the two size estimates in
\eqref{eq:exact-current-axis-core-small} follow.
\end{proof}

\subsection{Cone-local H\"older conventions and toroidal regularity}

Several estimates below are local in the fixed cone $\mathcal C_{\sigma_*}$, which meets the symmetry axis at the
stagnation point.  We therefore specify the H\"older seminorm used on this cone and state the elementary axis regularity
fact that removes the apparent singularity of the toroidal basis vector: the power $r^\alpha$ makes $r^\alpha\bs
e_\theta$ a $C^\alpha$ vector field across the axis.

\begin{remark}[Cone-local H\"older seminorms and the axis point]
All H\"older seminorms below are cone-local and include the axis point by continuity.  For any $R_0\ge 1$ we set
\[
[f]_{C^\alpha(\mathcal{C}_{\sigma_*}\cap B(0,R_0))} :=\sup_{\substack{x\neq y\\ x,y\in (\mathcal{C}_{\sigma_*}\cap
B(0,R_0))\cup\{0\}}}
\tfrac{|f(x)-f(y)|}{|x-y|^\alpha},
\]
and we write $[f]_{C^\alpha(\mathcal{C}_{\sigma_*})}$ when the truncation radius is immaterial. The cone-local fields
used below extend continuously to $0$ along $\mathcal C_{\sigma_*}$, so this is well-posed.
\end{remark}

\subsubsection{Axis regularity of the toroidal vector field.}

\begin{lemma}[Axis regularity of $r^\alpha\bs e_\theta$]
\label{lem:toroidal-vector-holder}
Let $\alpha\in(0,1)$. We define
\[
\mathfrak t(x):=
\begin{cases}
r(x)^\alpha\,\bs e_\theta(x), & r(x)>0,\\
0, & r(x)=0,
\end{cases}
\qquad x\in\R^3.
\]
Then $\mathfrak t\in C^\alpha(\R^3)$ and there exists $C=C(\alpha)$ such that for all $x,y\in\R^3$,
\begin{equation}\label{eq:toroidal-vector-pointwise-holder}
|\mathfrak t(x)-\mathfrak t(y)|\le C\,|x-y|^\alpha.
\end{equation}
In particular, $[\mathfrak t]_{C^\alpha(\R^3)}\le C(\alpha)$.
\end{lemma}

\begin{proof}[Proof of Lemma~\ref{lem:toroidal-vector-holder}]
We write $r_x:=r(x)$ and $r_y:=r(y)$.

\rruninhead{Case 1: $\min\{r_x,r_y\}\le 2|x-y|$.} Using $|\bs e_\theta|\equiv 1$ (when defined)
and $\mathfrak t=0$ on the axis,
\[
|\mathfrak t(x)-\mathfrak t(y)|\le |\mathfrak t(x)|+|\mathfrak t(y)|\le r_x^\alpha+r_y^\alpha \le C\,|x-y|^\alpha.
\]

\rruninhead{Case 2: $\min\{r_x,r_y\}>2|x-y|$.} Then both points are away from the axis and
\[
\tfrac12 r_x\le r_y\le 2r_x,
\]
because $|r_x-r_y|\le |x-y|<\tfrac12\min\{r_x,r_y\}$.

We write
\[
\mathfrak t(x)-\mathfrak t(y)=(r_x^\alpha-r_y^\alpha)\bs e_\theta(x)+r_y^\alpha(\bs e_\theta(x)-\bs e_\theta(y)).
\]

For the first term we use the elementary inequality valid for $\alpha\in(0,1)$:
\begin{equation*}
|a^\alpha-b^\alpha|\le |a-b|^\alpha,\qquad a,b\ge 0,
\end{equation*}
to obtain $|r_x^\alpha-r_y^\alpha|\le |r_x-r_y|^\alpha\le |x-y|^\alpha$.

For the second term we use that $\bs e_\theta$ is smooth away from the axis and satisfies $|\nabla \bs e_\theta|\lesssim
r^{-1}$. Hence by the mean value theorem along the segment from $x$ to $y$,
\[
|\bs e_\theta(x)-\bs e_\theta(y)| \le C\,\tfrac{|x-y|}{\min\{r_x,r_y\}} \le C\,\tfrac{|x-y|}{r_x}.
\]
Therefore, using $r_y^\alpha\lesssim r_x^\alpha$,
\[
r_y^\alpha\,|\bs e_\theta(x)-\bs e_\theta(y)| \le C\,r_x^\alpha\,\tfrac{|x-y|}{r_x} =C\,|x-y|\,r_x^{\alpha-1}.
\]
Since in Case 2 we have $r_x>2|x-y|$, we get $r_x^{\alpha-1}\le (2|x-y|)^{\alpha-1}$ (because $\alpha-1<0$), and hence
\[
|x-y|\,r_x^{\alpha-1}\le C\,|x-y|^\alpha.
\]

Combining the two terms yields \eqref{eq:toroidal-vector-pointwise-holder} in Case 2, completing the proof.
\end{proof}

\subsection{Finite-clock entry into the cusp regime}

The small-clock argument starts only after the cusp clock has reached a fixed threshold $\mathfrak J_{\mathrm{finite}}$.
The following lemma supplies this entry mechanism: a fixed compact sector of labels stays inside the compressive cone
throughout the finite-clock interval and contributes a uniformly negative axial strain.

\begin{lemma}[Finite-clock entry sector]
\label{lem:finite-clock-driver-sector}
After  fixing $\mathfrak J_{\mathrm{finite}}\in(0,1)$,  there exist a compact label sector
\[
E_{\rm ent}:= \Bigl\{Y:\ \tfrac12\le |Y|\le1,\quad \tfrac12\sigma_{\rm ent}\le\sigma(Y)\le\sigma_{\rm ent}\Bigr\},
\qquad 0<\sigma_{\rm ent}<\sigma_{\cut},
\]
and a constant $c_{\rm ent}>0$, depending only on $\alpha,\gamma,\mathfrak J_{\mathrm{finite}}$, such that the following
holds at every time for which $J_{\cusp}(t)\in[\mathfrak J_{\mathrm{finite}},1]$:
\[
E_{\rm ent}\subset D_{\inn}^{\cusp}(t), \quad J_{\twoD}(Y,t)^{-1}\ge c_{\rm ent}\ \ \text{ for } \ \ Y\in E_{\rm ent},
\]
and the cusp-coordinate strain has a label representation
\begin{equation}
\mathcal W_{\cusp}(t) = \int_{\R^3} \mathcal K_{\cusp}^{\flat}(Y,t)\, \mathcal J_{\cusp}(Y,t)^{-1}\omega_{\theta,0}(Y)\,dY,
\qquad \mathcal K_{\cusp}^{\flat}(Y,t) := \tfrac1{4\pi}\mathcal K_W\bigl(0,\phi_{\cusp}(Y,t)\bigr),
\label{eq:Kcusp-flat-label-rep}
\end{equation}
where
\[
\mathcal K_{\cusp}^{\flat}(Y,t)\,\omega_{\theta,0}(Y)\le0 \quad\text{where } \omega_{\theta,0}(Y)\ne0,
\qquad \mathcal K_{\cusp}^{\flat}(Y,t)\ge c_{\rm ent}, \quad \mathcal J_{\cusp}(Y,t)^{-1}\ge c_{\rm ent} \quad\text{on }E_{\rm ent}
\]
in the upper half-space.  Consequently $-\mathcal W_{\cusp}(t)\ge c_{\rm ent}\Gamma$ on the same finite-clock range.
The physical cusp-coordinate strain also has the label representation
\begin{equation}
\rW_{\cusp}(t) = \int_{\R^3} \mathcal K_{\cusp}(Y,t)\,
\left(1-\chi_{\far}\!\left(\tfrac{|\phi(Y,t)|}{R_{\tail}}\right)\right) J_{\twoD}(Y,t)^{-1}\omega_{\theta,0}(Y)\,dY,
\label{eq:Kcusp-label-rep}
\end{equation}
where the kernel has the sign of the physical axial strain kernel:
\[
\mathcal K_{\cusp}(Y,t)\,\omega_{\theta,0}(Y)\le0 \quad\text{where } \omega_{\theta,0}(Y)\ne0,
\qquad 1-\chi_{\far}\!\left(\tfrac{|\phi(Y,t)|}{R_{\tail}}\right)=1,\quad \mathcal K_{\cusp}(Y,t)\ge c_{\rm ent} \quad\text{on }E_{\rm ent}
\]
in the upper half-space. The lower half-space contribution is the same by odd symmetry.
\end{lemma}

\begin{proof}[Proof of Lemma~\ref{lem:finite-clock-driver-sector}]
\runinhead{Step 1: Finite-clock bounds for the normalized cusp map.} We use the self-consistent construction from
Lemma~\ref{lem:self-consistent-cusp-flow}.  On a time interval on which $J_{\cusp}\in[\mathfrak J_{\mathrm{finite}},1]$,
consider the finite-clock bootstrap
\begin{equation}
\|D\Phi_{\cusp}\|_{L^\infty(D_{\core})} +\|D\Phi_{\cusp}^{-1}\|_{L^\infty(\Phi_{\cusp}(D_{\core}))} \le A,
\qquad -\mathcal W_{\cusp}(t)\ge a\Gamma,
\label{eq:finite-clock-bootstrap}
\end{equation}
where $A<\infty$ and $a>0$ are fixed below.  This is the finite-clock bootstrap
\eqref{eq:localized-finite-clock-bootstrap}; it holds at $t=0$ for some $A,a$ depending only on
$\mathfrak J_{\mathrm{finite}},\alpha,\gamma$.  Since $\Phi_{\cusp}(0,t)=0$ and $D_{\core}=\{|Y|\le R_{\tail}\}$, the
first bound in \eqref{eq:finite-clock-bootstrap} gives
\[
\Phi_{\cusp}(D_{\core},t)\subset B_{A R_{\tail}}.
\]
We set
\[
R_A:=2+A R_{\tail}.
\]
Assuming \eqref{eq:finite-clock-bootstrap}, and using the corresponding finite-clock $C^{1,\alpha}$ bounds supplied by the
self-consistent construction, the transported vorticity
\eqref{eq:Omega-cusp-def} is the push-forward of the fixed Cartesian $C^\alpha$ vector field
$\omega_{\theta,0}e_\theta$ by maps with uniformly controlled $C^{1,\alpha}$ geometry on bounded sets.  The localized
Calder\'on--Zygmund/Schauder estimate used in Lemma~\ref{lem:L2} therefore gives
\[
\|U_{\cusp}(\cdot,t)\|_{C^{1,\alpha}(B_{R_A})} \le C_A\Gamma .
\]

On every subinterval on which $\dot J_{\cusp}\neq0$, we use the clock variable
\[
\ell:=-\log J_{\cusp}.
\]
Because
\[
\dot J_{\cusp} =\tfrac12J_{\cusp}m(t)\mathcal W_{\cusp}(t),
\]
the ODE \eqref{eq:Phi-cusp-def} becomes
\begin{equation}
\p_\ell\Phi_{\cusp} = -2\,\tfrac{U_{\cusp}}{\mathcal W_{\cusp}}\circ\Phi_{\cusp}, \qquad \ell:=-\log J_{\cusp}.
\label{eq:finite-clock-Phi-ell}
\end{equation}
We define
\[
F(x,t):=-2\,\tfrac{U_{\cusp}(x,t)}{\mathcal W_{\cusp}(t)}.
\]
The estimates above and the second bound in \eqref{eq:finite-clock-bootstrap} imply
\begin{equation}
\|F(\cdot,t)\|_{C^{1,\alpha}(B_{R_A})} \le \tfrac{2C_A}{a} =:M_A \qquad \bigl(J_{\cusp}(t)\in[\mathfrak J_{\mathrm{finite}},1]\bigr).
\label{eq:finite-clock-normalized-F-bound}
\end{equation}
Letting
\[
L_{\mathrm{finite}}:=\log \mathfrak J_{\mathrm{finite}}^{-1},
\]
as long as $\Phi_{\cusp}(D_{\core},t)\subset B_{R_A}$, \eqref{eq:finite-clock-Phi-ell} and
\eqref{eq:finite-clock-normalized-F-bound} show that, for $Y\in D_{\core}$,
\[
|\Phi_{\cusp}(Y,\ell)| \le |Y|+\int_0^\ell M_A\,d\ell' \le R_{\tail}+M_A L_{\mathrm{finite}}.
\]
If $P(Y,\ell):=D_Y\Phi_{\cusp}(Y,\ell)$, then
\[
\p_\ell P = \nabla F(\Phi_{\cusp}(Y,\ell),t(\ell))P,
\]
and hence
\[
|P(Y,\ell)| \le \exp(M_A L_{\mathrm{finite}}).
\]
For the inverse gradient $Q(Y,\ell):=P(Y,\ell)^{-1}$,  we have that
\[
\p_\ell Q = -Q\,\nabla F(\Phi_{\cusp}(Y,\ell),t(\ell)), \qquad |Q(Y,\ell)| \le \exp(M_A L_{\mathrm{finite}}).
\]
We choose $A$ larger than
\[
2\max\{1,\ R_{\tail}+M_A L_{\mathrm{finite}},\ \exp(M_A L_{\mathrm{finite}})\}.
\]
Thus,  the Gr\"onwall estimates improve the first bound in \eqref{eq:finite-clock-bootstrap}, once the strain lower
bound is improved in Step 3 of the proof below.  In particular, throughout the bootstrap interval we have, for labels in
$D_{\core}$,
\begin{equation}
|\Phi_{\cusp}(Y,t)|+|D_Y\Phi_{\cusp}(Y,t)| +|D_Y\Phi_{\cusp}(Y,t)^{-1}| \le C_{\rm ent}.
\label{eq:finite-clock-Phi-C1}
\end{equation}
The scalar modulation cancels from \eqref{eq:finite-clock-Phi-ell} because both the trajectory speed $m(t)U_{\cusp}$ and
the rate of change $\tfrac12J_{\cusp}m(t)\mathcal W_{\cusp}(t)$ contain the same scalar multiplier $m(t)$. The
self-consistent construction gives the same finite-clock bounds for the smooth map and the exact pulled-back map:
\begin{equation}
|D\phi_{\smooth}|+|D\phi_{\smooth}^{-1}| +|D\phi_{\cusp}|+|D\phi_{\cusp}^{-1}| \le C_{\rm ent} \quad\text{on the images of }D_{\core}.
\label{eq:finite-clock-smooth-cusp-C1}
\end{equation}
Indeed, the ODE for the smooth map is driven by $u_{\smooth}$ in \eqref{eq:smooth-velocity-def}, whose
finite-clock bounds on the bounded core follow from the same compactness argument; the exact pull-back is the
composition $\phi_{\cusp}=\phi_{\smooth}^{-1}\circ\phi$.  Thus
\eqref{eq:finite-clock-smooth-cusp-C1} is obtained on the
finite-clock interval $J_{\cusp}\in[\mathfrak J_{\mathrm{finite}},1]$: $C_{\rm ent}$ may deteriorate as $\mathfrak
J_{\mathrm{finite}}\downarrow0$, but it is independent of the small-cusp-clock boundedness conclusion for $m(t)$ proved
in Lemma~\ref{lem:modulation-bounded}.

\runinhead{Step 2: Choice of a fixed entry sector.} This is only a finite-clock choice and does not use a singular
hyperbolic asymptotic.  The map $\Phi_{\cusp}$ preserves both the symmetry axis and the plane $z=0$, and
\eqref{eq:finite-clock-Phi-C1} gives the same bound for $\Phi_{\cusp}$ and its inverse on $D_{\core}$.
Hence, for $Y\in E_{\rm ent}$,
\[
r(\Phi_{\cusp}(Y,t))\le C_{\rm ent}R(Y) \le C_{\rm ent}\sigma_{\rm ent},
\]
while the inverse Lipschitz bound and preservation of the plane $z=0$ give
\[
z(\Phi_{\cusp}(Y,t)) \ge c_{\rm ent}\,\operatorname{dist}(Y,\{z=0\}) \ge c_{\rm ent}
\]
after decreasing the fixed sector aperture and using $\tfrac12\le |Y|\le1$, $\sigma(Y)\le\sigma_{\rm ent}\ll1$.
Choosing $\sigma_{\rm ent}$ so small that $C_{\rm ent}\sigma_{\rm ent}<c_{\rm ent}\tan\sigma_{\inn}$ gives
$\Phi_{\cusp}(E_{\rm ent},t)\subset\mathcal C_{\inn}$ for every time under consideration.  The set $D_{\inn}^{\cusp}(t)$
in \eqref{eq:cusp-domain-def}, however, is defined using the exact cusp map $\phi_{\cusp}$, not the auxiliary map
$\Phi_{\cusp}$.  The finite-clock bound
\eqref{eq:finite-clock-smooth-cusp-C1} gives the same Lipschitz control for $\phi_{\cusp}$ and
$\phi_{\cusp}^{-1}$ on the compact image of $D_{\core}$.  Repeating the preceding distance-to-axis and
distance-to-$\{z=0\}$ estimates with $\phi_{\cusp}$ in place of $\Phi_{\cusp}$, and decreasing $\sigma_{\rm ent}$ once
more if necessary, yields
\[
\phi_{\cusp}(E_{\rm ent},t)\subset\mathcal C_{\inn}.
\]
By \eqref{eq:cusp-domain-def}, this is precisely $E_{\rm ent}\subset D_{\inn}^{\cusp}(t)$ for every time under
consideration.

\runinhead{Step 3: The flat cusp-coordinate strain is uniformly negative.} We next prove the sign lower bound for
$\mathcal W_{\cusp}=\p_z(U_{\cusp})_z(0,t)$, because this is the quantity in the denominator of the scalar
modulation \eqref{eq:modulation-def}.  Differentiating the exact Biot--Savart formula
\eqref{eq:U-cusp-label} at the origin gives \eqref{eq:Kcusp-flat-label-rep}.  The finite-clock
map bounds keep $\phi_{\cusp}(E_{\rm ent},t)$ in a compact subset of the upper cone, separated from the axis, the
equatorial plane, and infinity.  Hence the strain kernel \eqref{eq:KW-kernel} is bounded below on this set, and the
exact cusp Jacobian $\mathcal J_{\cusp}^{-1}$ is bounded below by the inverse map bounds.

The sign of the whole integrand is also fixed.  The cusp map preserves the upper and lower half-spaces, and $\mathcal
K_W(0,x)$ has the sign of $z$.  Since $\omega_{\theta,0}$ is negative in the upper half-space and odd across the
equatorial plane, we have
\[
\mathcal K_{\cusp}^{\flat}(Y,t)\,\omega_{\theta,0}(Y)\le0 \qquad\text{where }\omega_{\theta,0}(Y)\ne0 .
\]
On the compact sector $E_{\rm ent}$, the Target Profile satisfies $-\omega_{\theta,0}\ge c\Gamma$ after decreasing the
constant if necessary.  Therefore the full integral is bounded above by its contribution on $E_{\rm ent}$:
\[
\mathcal W_{\cusp}(t) \le -c_{\rm ent}\Gamma .
\]
This proves the strict improvement of the second bootstrap bound in
\eqref{eq:finite-clock-bootstrap} after the constant $a$ is fixed sufficiently small.  Together with
the Gr\"onwall improvement in Step 1, the finite-clock bootstrap closes on $J_{\cusp}\in[\mathfrak
J_{\mathrm{finite}},1]$.

\runinhead{Step 4: The physical cusp-coordinate strain has the same favorable sign.} It remains to prove the physical
label representation and the corresponding sign statement in
\eqref{eq:Kcusp-label-rep}.  The definition \eqref{eq:ucore-def} and
the Biot--Savart formula \eqref{eq:BS-axisymm} give
\begin{equation}
u_{\core}(x,t) = \tfrac1{4\pi}\int_{\R^3} K\bigl(x,\phi(Y',t)\bigr)\,
\left(1-\chi_{\far}\!\left(\tfrac{|\phi(Y',t)|}{R_{\tail}}\right)\right) J_{\twoD}(Y',t)^{-1}\omega_{\theta,0}(Y')\,dY'.
\label{eq:ucore-current-nearfield-finite-entry}
\end{equation}
On the finite-clock interval, \eqref{eq:finite-clock-smooth-cusp-C1} keeps the images of the compact sector $E_{\rm
ent}$ in a fixed ball.  Increasing $R_{\tail}$ once more, depending only on $\mathfrak J_{\mathrm{finite}}$, makes the
cutoff function in
\eqref{eq:ucore-current-nearfield-finite-entry} equal to one on those images.
The definition \eqref{eq:cusp-map-def} gives the exact pull-back formula
\begin{equation}
V_{\cusp}(X,t) = D\phi_{\smooth}(X,t)^{-1}\, u_{\core}(\phi_{\smooth}(X,t),t).
\label{eq:vcusp-ucore-pullback}
\end{equation}
Since $\phi_{\smooth}$ preserves the axis and fixes the origin, its meridional derivative at the origin is diagonal with
positive entries. Also $u_{\core}(0,t)=0$ by the odd symmetry across $z=0$. Differentiating
\eqref{eq:vcusp-ucore-pullback} at $X=0$ therefore gives
\[
\p_z(V_{\cusp})_z(0,t)=\p_z(u_{\core})_z(0,t);
\]
the derivative of $D\phi_{\smooth}^{-1}$ is multiplied by $u_{\core}(0,t)$ and hence drops out. Combining this identity
with \eqref{eq:ucore-current-nearfield-finite-entry} and \eqref{eq:rW2} gives the near-field strain representation
\[
\p_z(u_{\core})_z(0,t) = \int_{\R^3} \mathcal K_{\cusp}(Y,t) \left(1-\chi_{\far}\!\left(\tfrac{|\phi(Y,t)|}{R_{\tail}}\right)\right)
J_{\twoD}(Y,t)^{-1}\omega_{\theta,0}(Y)\,dY,
\]
with
\[
\mathcal K_{\cusp}(Y,t) = \tfrac{1}{4\pi}\, \mathcal K_W\bigl(0,\phi(Y,t)\bigr).
\]
Thus the sign of $\mathcal K_{\cusp}$ is the sign of the physical kernel \eqref{eq:KW-kernel} evaluated at the exact
image of the label. The Euler flow preserves the upper and lower half-spaces, and $K_W$ has the sign of $z$ while
$\omega_{\theta,0}$ is odd across the equator and negative in the upper half-space. Therefore $\mathcal
K_{\cusp}(Y,t)\omega_{\theta,0}(Y)\le0$ wherever $\omega_{\theta,0}(Y)\neq0$.  Since the cutoff function is nonnegative,
the integral over all labels is bounded above by its restriction to $E_{\rm ent}$.

On the compact sector $E_{\rm ent}$, the finite-clock $C^1$ bounds
\eqref{eq:finite-clock-Phi-C1} and \eqref{eq:finite-clock-smooth-cusp-C1}, together with their
inverse bounds, keep exact images in a compact subset of the upper part of $\mathcal C_{\inn}$, separated from the axis,
the equatorial plane, and infinity. Thus the physical kernel \eqref{eq:KW-kernel}, and hence $\mathcal K_{\cusp}$, is
bounded below there by a positive constant after decreasing $c_{\rm ent}$. The Jacobian bound is also a direct
consequence of the finite-clock map bounds.  Indeed, since
\[
\phi(\cdot,t)=\phi_{\smooth}(\cdot,t)\circ\phi_{\cusp}(\cdot,t),
\]
\eqref{eq:finite-clock-smooth-cusp-C1} gives
\[
|D_Y(\phi_r,\phi_z)(Y,t)|\le C_{\rm ent} \qquad (Y\in E_{\rm ent}).
\]
Hence $J_{\twoD}(Y,t)\le C_{\rm ent}$ on $E_{\rm ent}$, and after decreasing $c_{\rm ent}$ once more,
$J_{\twoD}(Y,t)^{-1}\ge c_{\rm ent}$, and this completes the proof.
\end{proof}

The preceding lemma gives a uniform compressive lower bound for the cusp-coordinate strain on the whole finite-clock
range.  We now integrate the cusp-clock law and use the same finite-clock bounds to rule out a breakdown before the
prescribed small-clock threshold is reached.

\begin{lemma}[Finite entry into a prescribed cusp-clock regime]
\label{lem:late-entry}
With  $\mathfrak J_{\mathrm{finite}}\in(0,1)$ fixed, the true Euler cusp clock $J_{\cusp}(t)$ reaches $\mathfrak
J_{\mathrm{finite}}$ in finite time. More precisely, there are constants $c_{\rm ent}=c_{\rm
ent}(\alpha,\gamma,\mathfrak J_{\mathrm{finite}})>0$ and $C_{\rm ent}=C_{\rm ent}(\alpha,\gamma,\mathfrak
J_{\mathrm{finite}})<\infty$ such that, until the first time
\[
t_{\rm ent}:=\inf\{t\ge0:\ J_{\cusp}(t)\le \mathfrak J_{\mathrm{finite}}\},
\]
the true Euler solution remains smooth and
\begin{equation}
\rW_{\cusp}(t)\le -c_{\rm ent}\Gamma, \qquad t_{\rm ent}\le C_{\rm ent}\Gamma^{-1}.
\label{eq:late-entry-bound}
\end{equation}
\end{lemma}

\begin{proof}[Proof of Lemma~\ref{lem:late-entry}]
\runinhead{Step 1: The entry sector gives a uniform negative strain.} Let $E_{\rm ent}$ and $c_{\rm ent}$ be given by
Lemma~\ref{lem:finite-clock-driver-sector}.  On $E_{\rm ent}$, the Target Profile satisfies
$\Theta^*(\sigma)=(\sin\sigma)^\alpha$.  Since $E_{\rm ent}$ is compactly contained in the upper half-space and away
from the origin, after decreasing $c_{\rm ent}$ we have
\[
-\omega_{\theta,0}(Y)\ge c_{\rm ent}\Gamma \qquad (Y\in E_{\rm ent}).
\]
By Lemma~\ref{lem:finite-clock-driver-sector}, for every $Y\in E_{\rm ent}$ and $J_{\cusp}(t)\in[\mathfrak
J_{\mathrm{finite}},1]$,
\[
\mathcal K_{\cusp}(Y,t)\ge c_{\rm ent}, \qquad J_{\twoD}(Y,t)^{-1}\ge c_{\rm ent}.
\]
Moreover the cutoff function in \eqref{eq:Kcusp-label-rep} equals $1$ on $E_{\rm ent}$, and the sign statement in
Lemma~\ref{lem:finite-clock-driver-sector} gives
\[
\mathcal K_{\cusp}(Y,t)\,\omega_{\theta,0}(Y)\le0 \qquad (\omega_{\theta,0}(Y)\ne0).
\]
Thus all labels outside $E_{\rm ent}$ can only make the integral in
\eqref{eq:Kcusp-label-rep} more negative.  Consequently, while
$J_{\cusp}(t)\in[\mathfrak J_{\mathrm{finite}},1]$,
\[
\rW_{\cusp}(t) \le \int_{E_{\rm ent}} \mathcal K_{\cusp}(Y,t)\,
J_{\twoD}(Y,t)^{-1}\omega_{\theta,0}(Y)\,dY \le -c_{\rm ent}^3\,|E_{\rm ent}|\,\Gamma.
\]
After renaming the positive constant, this gives
\begin{equation}
\rW_{\cusp}(t)\le -c_{\rm ent}\Gamma \qquad\text{while }J_{\cusp}(t)\in[\mathfrak J_{\mathrm{finite}},1].
\label{eq:finite-entry-strain-negative}
\end{equation}

\runinhead{Step 2: Integrating the cusp-clock law.} The exact cusp-clock law gives
\[
\p_t\log J_{\cusp}(t)=\tfrac12\rW_{\cusp}(t).
\]
Combining this identity with \eqref{eq:finite-entry-strain-negative}, we obtain
\[
\p_t\log J_{\cusp}(t) \le -\tfrac12 c_{\rm ent}\Gamma
\]
as long as the solution is smooth and $J_{\cusp}(t)\ge\mathfrak J_{\mathrm{finite}}$.  Since $J_{\cusp}(0)=1$,
integration gives
\[
\log J_{\cusp}(t) \le -\tfrac12 c_{\rm ent}\Gamma t.
\]
Therefore $J_{\cusp}$ reaches $\mathfrak J_{\mathrm{finite}}$ no later than
\[
T_{\rm ent}^{\rm bd} := 2c_{\rm ent}^{-1}\Gamma^{-1} \log(\mathfrak J_{\mathrm{finite}}^{-1}),
\]
unless the $C^{1,\alpha}$ Euler solution ceases to exist before that time.

\runinhead{Step 3: No breakdown can occur before entry.} It remains to exclude the exceptional alternative in Step 2.
Suppose the solution is smooth on $[0,T)$ and $J_{\cusp}(t)\ge\mathfrak J_{\mathrm{finite}}$ there.
Lemma~\ref{lem:finite-clock-driver-sector} closes the finite-clock bootstrap
\eqref{eq:localized-finite-clock-bootstrap}; in particular,
\eqref{eq:finite-clock-smooth-cusp-C1} gives uniform $C^1$ bounds for the exact cusp map, the smooth map, and their
inverses on the images of the bounded core.  The far-field smooth velocity $u_{\smooth}$ has uniform finite-clock
bounds on fixed balls by \eqref{eq:smooth-velocity-def}; the same type of bound appears later in
\eqref{eq:current-far-usmooth-derivatives}.

These map bounds keep the two-dimensional Jacobian away from zero on the bounded core.  The vorticity transport identity
\eqref{eq:vort-identity-JtwoD}, together with the fixed $L^\infty$ bound for $\omega_{\theta,0}$, then yields
\[
\|\omega(\cdot,t)\|_{L^\infty} \le C(\alpha,\gamma,\mathfrak J_{\mathrm{finite}})\Gamma
\qquad\bigl(J_{\cusp}(t)\in[\mathfrak J_{\mathrm{finite}},1]\bigr).
\]
The Beale--Kato--Majda continuation criterion therefore rules out a finite breakdown while $J_{\cusp}\ge\mathfrak
J_{\mathrm{finite}}$.  Hence the entry time is finite and satisfies
\[
t_{\rm ent}\le T_{\rm ent}^{\rm bd} \le C_{\rm ent}\Gamma^{-1}.
\]
Together with \eqref{eq:finite-entry-strain-negative}, this proves \eqref{eq:late-entry-bound}.
\end{proof}

\section{Pressure Hessian Bound for the Euler-Generated Axial Function}
\label{sec:slope-restricted-pressure}

The Riccati comparison in this section uses the separation-of-variables vorticity
\[
\Omega_\theta(\mathcal R,\mathcal Z)=-\operatorname{sgn}(\mathcal Z)a_t^{\rm phys}(|\mathcal Z|)\mathcal R^\alpha .
\]
The radial dependence is the fixed cylindrical cusp $\mathcal R^\alpha$, while the axial function
$a_t^{\rm phys}$ is generated by the Euler flow.  Section~\ref{sec::Pi-model-explain} provides the
motivation for why this separation-of-variables form is the correct model for the Riccati pressure Hessian
estimate.  In this section and in
Sections~\ref{sec:current-axis-normal-forms}--\ref{sec:target-profile-typeI-completion}, we prove the pressure
and axis-trace estimates needed to pass from the model to the Euler Riccati bound.

As we will prove, in  the normalized axial coordinate $\zeta=J^{-2}B_t(Z)$, the physical Euler-generated axial function $a_t^{\rm phys}$ is given by
\[
a_t^{\rm phys}(\zeta) = q_t(\zeta)^{1-\alpha}\bigl(1+Z_t(\zeta)^2\bigr)^{-\gamma/2}, \qquad q_t(\zeta)=JA_t(Z_t(\zeta)).
\]
This section proves that the stagnation-point axial strain and the pressure Hessian admit a common
one-dimensional formulation in the axial coordinate $\zeta$.  This formulation leads to the Riccati pressure
Hessian lower bound in Proposition~\ref{prop:euler-generated-profile-riccati}.  The reduction from
three-dimensional singular integrals to one-dimensional integral comparisons is made possible by the
separation-of-variables vorticity and the non-negative monotonic structure of the Euler-generated axial function
$a_t^{\rm phys}$.

\subsection{Slope-restricted model vorticity and axial strain}

We study odd toroidal vorticities in which the radial and axial variables appear in the separation-of-variables form:
\[
\Omega_{\theta}^{a}(\mathcal R,\mathcal Z)=-\operatorname{sgn}(\mathcal Z)a(|\mathcal Z|)\mathcal R^\alpha .
\]
Here, the subscript $\theta$ denotes the toroidal component of vorticity, and the superscript $a$ denotes the vorticity-dependence on the
axial function $a(\mathcal Z)$. 
As we explained in Section~\ref{sec::Pi-model-explain}, after the normalized axial coordinate
$\zeta=J^{-2}B_t(Z)$ is introduced, the vorticity below is the $\phi_{\cusp}$-transport of the localized initial datum
\[
\Omega_{\theta,0}(R,Z) =-\Gamma\operatorname{sgn}(Z)R^\alpha(1+Z^2)^{-\gamma/2}.
\]
The variable separation splits the axial dependence $a(|\mathcal Z|)$ from the radial dependence $\mathcal R^\alpha$.  We now introduce the slope variable
\[
\tau=\tfrac{\mathcal R}{|\mathcal Z|},\qquad \mathcal Z\ne 0,
\]
and we choose the slope cutoff function $\chi_M\in C_c^\infty([0,\infty))$ satisfying $0\le\chi_M\le1$ such that
\begin{equation}
\chi_M(\tau)=
\begin{cases}
1 & \text{if }0\le\tau\le M,\\
0 & \text{if }\tau\ge2M.
\end{cases}
\label{eq:slope-cutoff-chiM}
\end{equation}
For an axial function $a$ on $\{\mathcal Z\ge0\}$, we define the $M$-slope-restricted vorticity by
\begin{equation}
\Omega_{\theta}^{a,M}(\mathcal R,\mathcal Z)
=-\operatorname{sgn}(\mathcal Z)\,a(|\mathcal Z|)\mathcal R^\alpha
\chi_M\!\left(\tfrac{\mathcal R}{|\mathcal Z|}\right).
\label{eq:vort-slope-restricted}
\end{equation}
and we define the associated Biot--Savart velocity by $U^{a,M}:=\BS[\Omega_\theta^{a,M}e_\theta]$.
The associated model pressure Hessian at the origin is denoted by
\begin{equation}
\Pi_M[a] :=\pv\int_{\R^3}K_{zz}(\mathcal Y)\,
\tr\bigl(\nabla U^{a,M}(\mathcal Y)\nabla U^{a,M}(\mathcal Y)\bigr)\,d\mathcal Y .
\label{eq:pressure-bilinear-form}
\end{equation}
Here $K_{zz}$ is defined in \eqref{eq:Kzz-kernel}.

The vorticity $\Omega_{\theta}^{a,M}(\mathcal R,\mathcal Z) $ then generates the velocity $U^{a,M}$, from which we define the associated stagnation-point axial strain
\[
W_M[a]:=\p_{\mathcal Z}(U^{a,M})_{\mathcal Z}(0,0).
\]
Applying the Biot--Savart law to \eqref{eq:vort-slope-restricted}, the axial strain separates into an angular constant and an axial moment:
\begin{equation}
-W_M[a]=C_{\alpha,M}^W I_1[a], \qquad I_1[a]:=\int_0^\infty a(\zeta)\zeta^{\alpha-1}\,d\zeta.
\label{eq:scaled-strain-exact}
\end{equation}
where
\begin{equation}
C_{\alpha,M}^W := \int_0^\infty \tfrac{3\tau^{\alpha+2}}{(1+\tau^2)^{5/2}}\chi_M(\tau)\,d\tau \xrightarrow[M\to\infty]{}
\int_0^\infty \tfrac{3\tau^{\alpha+2}}{(1+\tau^2)^{5/2}}\,d\tau =:C_\alpha^W>0.
\label{eq:C-alpha-M-limit}
\end{equation}

\begin{lemma}[Slope-restricted axial strain moment]
\label{lem:scaled-zeta-strain-moment}
For every axial function $a\ge0$ for which $I_1[a]<\infty$, the identity \eqref{eq:scaled-strain-exact} holds: $W_M[a]=-C_{\alpha,M}^W I_1[a]$.
\end{lemma}

\begin{proof}[Proof of Lemma~\ref{lem:scaled-zeta-strain-moment}]
The axial strain identity is
\[
-W_M[a] = \int_0^\infty\!\int_0^\infty
\tfrac{3\mathcal R^2\mathcal Z}{(\mathcal R^2+\mathcal Z^2)^{5/2}}
a(\mathcal Z)\mathcal R^\alpha\chi_M\!\left(\tfrac{\mathcal R}{\mathcal Z}\right)\,d\mathcal R\,d\mathcal Z .
\]
With $\mathcal R=\mathcal Z\tau$,
$\tfrac{3\mathcal R^2\mathcal Z}{(\mathcal R^2+\mathcal Z^2)^{5/2}}\mathcal R^\alpha\,d\mathcal R
= \tfrac{3\tau^{\alpha+2}}{(1+\tau^2)^{5/2}}\mathcal Z^{\alpha-1}\,d\tau$, and so
\[
-W_M[a] = \left(\int_0^\infty\tfrac{3\tau^{\alpha+2}}{(1+\tau^2)^{5/2}}\chi_M(\tau)\,d\tau\right)
\left(\int_0^\infty a(\zeta)\zeta^{\alpha-1}\,d\zeta\right),
\]
which is \eqref{eq:scaled-strain-exact}.
\end{proof}

We next introduce notation for the angular region selected by $1-\chi_M$.  For a function $h(\tau)$ of the slope variable, we set
\begin{equation}
\mathcal A_\alpha[h]^2 := \int_0^\infty\tfrac{(1+\tau)^{2\alpha}}{1+\tau^2}|h(\tau)|^2\,d\tau ,  \label{eq:angular-tail-norm}
\end{equation}
and we define
\[
\mathfrak a_{\rm ang}(M) := 2\int_M^\infty\tfrac{(1+\tau)^{2\alpha}}{1+\tau^2}\,d\tau .
\]
For $M\ge1$,
\begin{equation}
\mathfrak a_{\rm ang}(M) \le \tfrac{2^{2\alpha+1}}{1-2\alpha}M^{2\alpha-1} \longrightarrow0 \qquad(M\to\infty),
\label{eq:scaled-angular-tail}
\end{equation}
and, since $1-\chi_M$ is supported in $[M,\infty)$,
\begin{equation}
\mathcal A_\alpha[1-\chi_M] \le \mathfrak a_{\rm ang}(M)^{{\frac{1}{2}}}.
\label{eq:angular-tail-norm-cutoff}
\end{equation}

For a function $h(\tau)$ of the slope variable, we define
\[
\Omega_{\theta,h}^{a}(\mathcal R,\mathcal Z)
:= -\operatorname{sgn}(\mathcal Z)a(|\mathcal Z|)\mathcal R^\alpha
h\!\left(\tfrac{\mathcal R}{|\mathcal Z|}\right),
\qquad U_h^a:=\BS[\Omega_{\theta,h}^{a}e_\theta], \qquad
W_h[a]:=\p_{\mathcal Z}(U_h^a)_{\mathcal Z}(0,0).
\]
For two angular functions $h_1,h_2$, we define the polarized pressure Hessian by
\[
\Pi_{h_1,h_2}[a] := \pv\int_{\R^3}K_{zz}(\mathcal Y)\,
\tr\bigl(\nabla U_{h_1}^a(\mathcal Y)\nabla U_{h_2}^a(\mathcal Y)\bigr)\,d\mathcal Y.
\]
Likewise, for two axial functions $a_1,a_2$, we set
\[
\Pi_M[a_1,a_2] := \pv\int_{\R^3}K_{zz}(\mathcal Y)\,
\tr\bigl(\nabla U^{a_1,M}(\mathcal Y)\nabla U^{a_2,M}(\mathcal Y)\bigr)\,d\mathcal Y,
\qquad \Pi_M[a]:=\Pi_M[a,a].
\]

\begin{lemma}[Large-slope axial strain and pressure Hessian estimates]
\label{lem:bilinear-angular-tail}
The axial strain satisfies
\[
|W_h[a]| \le C_\alpha\mathcal A_\alpha[h]I_1[a].
\]
Assume, for the fixed axial function $a$, that
\begin{equation}
C_{\alpha,a} := \sup_{\mathcal A_\alpha[h_1]\le1,\ \mathcal A_\alpha[h_2]\le1} |\Pi_{h_1,h_2}[a]|<\infty .
\label{eq:C-alpha-a}
\end{equation}
Then
\begin{equation}
|\Pi_{h_1,h_2}[a]| \le C_{\alpha,a}\mathcal A_\alpha[h_1]\mathcal A_\alpha[h_2].
\label{eq:angular-tail-pressure-hessian-integral}
\end{equation}
If, in addition, the axial function $a$ is nonnegative, nonincreasing, compactly supported, and $I_1[a]<\infty$, then
\begin{equation}
|\Pi_{h_1,h_2}[a]| \le C_\alpha \mathcal A_\alpha[h_1]\mathcal A_\alpha[h_2] I_1[a]^2 .
\label{eq:angular-tail-pressure-hessian-monotone}
\end{equation}
Consequently,
\begin{equation}
|W_{1-\chi_M}[a]| \le C_\alpha\mathfrak a_{\rm ang}(M)^{{\frac{1}{2}}}I_1[a],
\label{eq:linear-angular-tail-strain}
\end{equation}
and, for every angular function $h$ with $\mathcal A_\alpha[h]<\infty$,
\begin{equation}
|\Pi_{h,1-\chi_M}[a]| \le C_{\alpha,a}\mathcal A_\alpha[h]\mathfrak a_{\rm ang}(M)^{{\frac{1}{2}}},
\qquad |\Pi_{1-\chi_M,1-\chi_M}[a]| \le C_{\alpha,a}\mathfrak a_{\rm ang}(M).
\label{eq:bilinear-angular-tail}
\end{equation}
\end{lemma}

\begin{proof}[Proof of Lemma~\ref{lem:bilinear-angular-tail}]
The proof of Lemma~\ref{lem:scaled-zeta-strain-moment}, with $h$ in place of $\chi_M$, gives
\[
|W_h[a]| \le I_1[a]\int_0^\infty \tfrac{3\tau^{\alpha+2}}{(1+\tau^2)^{5/2}}|h(\tau)|\,d\tau .
\]
The integral
$\int_0^\infty\left(\tfrac{3\tau^{\alpha+2}}{(1+\tau^2)^{5/2}}\right)^2\tfrac{1+\tau^2}{(1+\tau)^{2\alpha}}\,d\tau$
is finite.  Therefore, Cauchy--Schwarz and \eqref{eq:angular-tail-norm} imply
\[
|W_h[a]|\le C_\alpha\mathcal A_\alpha[h]I_1[a].
\]
If $\mathcal A_\alpha[h_1]\mathcal A_\alpha[h_2]=0$, the pressure estimate is trivial; otherwise, we set $\widetilde h_i:=\mathcal A_\alpha[h_i]^{-1}h_i$.  
Then $\mathcal A_\alpha[\widetilde h_i]=1$, and bilinearity shows that
\[
\Pi_{h_1,h_2}[a]= \mathcal A_\alpha[h_1]\mathcal A_\alpha[h_2] \Pi_{\widetilde h_1,\widetilde h_2}[a].
\]
By the definition of $C_{\alpha,a}$ in \eqref{eq:C-alpha-a}, this proves
\eqref{eq:angular-tail-pressure-hessian-integral}.  The two estimates in
\eqref{eq:bilinear-angular-tail} follow from \eqref{eq:angular-tail-norm-cutoff}.

It remains to prove \eqref{eq:angular-tail-pressure-hessian-monotone}.  For $\mathcal Z,\mathcal Z'>0$, we set
$\mathcal R=\mathcal Z\tau$ and  $\mathcal R'=\mathcal Z'\tau'$.
After the azimuthal variables are integrated out, the homogeneity of the Biot--Savart kernel gradient
$\mathcal K$ in \eqref{eq:grad-BS} and of the pressure kernel $K_{zz}$ in \eqref{eq:Kzz-kernel} produces a
reduced kernel $\mathfrak L_\alpha$ defined by
\begin{align}
\Pi_{h_1,h_2}[a]
=\int_0^\infty\!\!\int_0^\infty\!\!\int_0^\infty\!\!\int_0^\infty
a(\mathcal Z)\mathcal Z^{\alpha-1}h_1(\tau)\,
\mathfrak L_\alpha(\log \mathcal Z-\log \mathcal Z',\tau,\tau')
h_2(\tau')a(\mathcal Z')(\mathcal Z')^{\alpha-1}\,
d\tau'\,d\tau\,d\mathcal Z'\,d\mathcal Z .     \label{eq:pressure-log-scale-before-convolution}
\end{align}
We set $f(x):=a(e^x)e^{\alpha x}$.
Since $a(\mathcal Z)\mathcal Z^{\alpha-1}\,d\mathcal Z=f(x)\,dx$, \eqref{eq:pressure-log-scale-before-convolution} becomes
\begin{equation}
\Pi_{h_1,h_2}[a]=\int_{\mathbb R}\!\int_{\mathbb R}f(x)\, \mathfrak K_{h_1,h_2}(x-x')\,f(x')\,dx'\,dx , \label{eq:log-pressure-convolution}
\end{equation}
where
\[
\mathfrak K_{h_1,h_2}(s)= \int_0^\infty\!\int_0^\infty h_1(\tau)h_2(\tau') \mathfrak L_\alpha(s,\tau,\tau')\,d\tau'\,d\tau ,
\]
and
\begin{equation}
\int_{\mathbb R}|\mathfrak L_\alpha(s,\tau,\tau')|\,ds \le C_\alpha \tfrac{(1+\tau)^{2\alpha}}{1+\tau^2} \tfrac{(1+\tau')^{2\alpha}}{1+\tau'^2}.
\label{eq:reduced-pressure-kernel-L1}
\end{equation}
Indeed, \eqref{eq:reduced-pressure-kernel-L1} follows by splitting the $s$-integral as
\[
\int_{|s|\ge2}|\mathfrak L_\alpha(s,\tau,\tau')|\,ds \le C_\alpha \tfrac{(1+\tau)^{2\alpha}}{1+\tau^2} \tfrac{(1+\tau')^{2\alpha}}{1+\tau'^2}
\int_{|s|\ge2}e^{-c|s|}\,ds .
\]
For $|s|\le2$, we define the variables
\[
\mathcal Z'=e^{x'},\qquad \mathcal Z=e^{x'+s},\qquad e^{-2}\le \tfrac{\mathcal Z}{\mathcal Z'}\le e^2 .
\]
The dilation $\mathcal Y=e^{x'}\overline{\mathcal Y}$,
$\mathcal Y'=e^{x'}\overline{\mathcal Y}'$ shows that
\[
\overline{\mathcal Z}'=1,\qquad \overline{\mathcal Z}=e^s,\qquad
\overline{\mathcal R}=e^s\tau,\qquad \overline{\mathcal R}'=\tau'.
\]
In these variables, \eqref{eq:grad-BS} and \eqref{eq:Kzz-kernel} imply
\[
\int_{|s|\le2}|\mathfrak L_\alpha(s,\tau,\tau')|\,ds \le C_\alpha \tfrac{(1+\tau)^{2\alpha}}{1+\tau^2} \tfrac{(1+\tau')^{2\alpha}}{1+\tau'^2}.
\]
Setting $m_\alpha(\tau):=\tfrac{(1+\tau)^{2\alpha}}{1+\tau^2}$, by   \eqref{eq:angular-tail-norm},  $\|h\,m_\alpha^{\frac{1}{2}} \|_{L^2_\tau}=\mathcal A_\alpha[h]$, and
$\int_0^\infty m_\alpha(\tau)\,d\tau<\infty$ because $\alpha<\tfrac12$.  By \eqref{eq:reduced-pressure-kernel-L1},
\[
\|\mathfrak K_{h_1,h_2}\|_{L^1(\mathbb R)}
\le C_\alpha \int_0^\infty |h_1(\tau)|m_\alpha(\tau)\,d\tau \int_0^\infty |h_2(\tau')|m_\alpha(\tau')\,d\tau' .
\]
Cauchy--Schwarz applied to the last two integrals gives
\[
\|\mathfrak K_{h_1,h_2}\|_{L^1(\mathbb R)} \le C_\alpha\mathcal A_\alpha[h_1]\mathcal A_\alpha[h_2], 
\]
and Young's inequality applied to \eqref{eq:log-pressure-convolution} gives
\[
|\Pi_{h_1,h_2}[a]| \le C_\alpha\mathcal A_\alpha[h_1]\mathcal A_\alpha[h_2]\|f\|_{L^2(\mathbb R)}^2.
\]
Since
\[
\|f\|_{L^2(\mathbb R)}^2 =\int_0^\infty a(\zeta)^2\zeta^{2\alpha-1}\,d\zeta ,
\]
we find
\[
|\Pi_{h_1,h_2}[a]| \le C_\alpha\mathcal A_\alpha[h_1]\mathcal A_\alpha[h_2]\int_0^\infty a(\zeta)^2\zeta^{2\alpha-1}\,d\zeta .
\]
Since $a$ is nonnegative and nonincreasing, for every $\zeta>0$ and $0<\eta<\zeta$ we have
$a(\eta)\ge a(\zeta)$, and so
\[
a(\zeta)\tfrac{\zeta^\alpha}{\alpha} =a(\zeta)\int_0^\zeta\eta^{\alpha-1}\,d\eta \le \int_0^\zeta a(\eta)\eta^{\alpha-1}\,d\eta \le I_1[a].
\]
Multiplying this estimate by $\alpha a(\zeta)\zeta^{\alpha-1}$ and integrating in $\zeta$ yields
\[
\int_0^\infty a(\zeta)^2\zeta^{2\alpha-1}\,d\zeta \le \alpha I_1[a]\int_0^\infty a(\zeta)\zeta^{\alpha-1}\,d\zeta =\alpha I_1[a]^2 ,
\]
which proves \eqref{eq:angular-tail-pressure-hessian-monotone}.
\end{proof}

\subsection{The Euler-generated axial function}

We denote the components of the cusp flow map and the cusp clock by
\[
\phi_{\cusp}(R,Z,t)=(r_t(R,Z),z_t(R,Z)), \qquad J=J_{\cusp}(t).
\]
On the symmetry axis $R=0$, we set
\[
A_t(Z):=\p_Rr_t(0,Z), \qquad B_t(Z):=z_t(0,Z).
\]
We set
\[
\zeta=J^{-2}B_t(Z),\qquad Z_t(\zeta):=(J^{-2}B_t)^{-1}(\zeta),
\]
and we define $b_t(\zeta):=\bigl(\p_\zeta Z_t(\zeta)\bigr)^{-1}$.  Differentiating
$B_t(Z_t(\zeta))=J^2\zeta$ gives
\begin{equation}
b_t(\zeta)=J^{-2}B_t'(Z_t(\zeta)).
\label{eq:axis-b-function-Bprime}
\end{equation}
For $0<\zeta_{\rm mon}<\infty$, we use the interval $I_{\rm mon}:=[0,\zeta_{\rm mon}]$ from the monotone axial-stretching bootstrap \textup{(BA4)},
on which  we impose the bootstrap assumption \textup{(BA4)}, i.e.,  the two-sided bound \eqref{eq:monotone-axial-two-sided} and the monotone fractional-increment bound
\eqref{eq:monotone-axial-fractional-bootstrap} for $b_t$.  For the slope-restricted model pressure Hessian
$\Pi_M[\cdot]$ in \eqref{eq:pressure-bilinear-form}, we use the restriction of the axial function to the smaller interval
\begin{equation}
I_a:=[0,\zeta_a], \qquad 0<\zeta_a<\zeta_{\rm mon},
\label{eq:pressure-profile-interval}
\end{equation}
and choose the localization cutoff from Section~\ref{sec:fixed-choice-order} so that
$\operatorname{supp}\vartheta_\sharp\Subset I_a$.  On $I_a$, we define
\[
q_t(\zeta):=JA_t(Z_t(\zeta)).
\]
Since $J=\det\nabla_{(R,Z)}\phi_{\cusp}(0,0,t)$ and $A_t(0)=J^{-1}$, we have $q_t(0)=1$.  By the axisymmetric incompressibility identity,
\begin{equation}
A_t(Z)^2B_t'(Z)=1.
\label{eq:axis-volume-At-Bt}
\end{equation}
Combining \eqref{eq:axis-volume-At-Bt} with $q_t(\zeta)=JA_t(Z_t(\zeta))$ and \eqref{eq:axis-b-function-Bprime}, we arrive at
\begin{equation}
q_t(\zeta)^2b_t(\zeta)=1. \label{eq:axis-qb-volume-mon}
\end{equation}
As we explained in Section~\ref{sec::Pi-model-explain}, the Euler-generated axial function is given by
\begin{equation}
a_t^{\rm phys}(\zeta) = q_t(\zeta)^{1-\alpha} \bigl(1+Z_t(\zeta)^2\bigr)^{-\gamma/2}.
\label{eq:scaled-profile-at}
\end{equation}
The axial function in \eqref{eq:scaled-profile-at} is defined on $I_{\rm mon}$.  The axial function used in the model
pressure estimate is the restriction of $a_t^{\rm phys}$ to  $I_a$:
\begin{equation}
a_t(\zeta):=a_t^{\rm phys}(\zeta)\mathbf 1_{I_a}(\zeta).
\label{eq:euler-generated-truncated-coeff}
\end{equation}

\begin{lemma}[Euler-generated axial function monotonicity]
\label{lem:exact-euler-generated-derivative}
Assume \eqref{eq:axis-qb-volume-mon} and suppose that $b_t$ is nondecreasing on $I_a$.  Then $a_t^{\rm phys}$ in \eqref{eq:scaled-profile-at} 
is nonnegative and nonincreasing on $I_a$.  At every point where $b_t$ is differentiable, the ordinary derivative satisfies
\begin{equation}
-\p_\zeta a_t^{\rm phys} = a_t^{\rm phys}\left[\tfrac{1-\alpha}{2}\tfrac{b_t'}{b_t}+\gamma\tfrac{Z_t\p_\zeta Z_t}{1+Z_t^2} \right].
\label{eq:exact-euler-minus-a-prime}
\end{equation}
\end{lemma}

\begin{proof}[Proof of Lemma~\ref{lem:exact-euler-generated-derivative}]
By \eqref{eq:axis-qb-volume-mon}, $q_t=b_t^{-\frac12}$.  Therefore \eqref{eq:scaled-profile-at} becomes
$a_t^{\rm phys}=b_t^{-(1-\alpha)/2}\bigl(1+Z_t^2\bigr)^{-\gamma/2}$.
Since $b_t$ is nondecreasing, $b_t^{-(1-\alpha)/2}$ is nonincreasing.  Since
$\p_\zeta Z_t=b_t^{-1}>0$ and $Z_t(0)=0$, the map $Z_t$ is increasing and nonnegative on
$I_a$, so $(1+Z_t^2)^{-\gamma/2}$ is also nonincreasing.  The product is nonnegative and nonincreasing.
At every point where $b_t$ is differentiable, the logarithmic derivative of
$b_t^{-(1-\alpha)/2}(1+Z_t^2)^{-\gamma/2}$ gives \eqref{eq:exact-euler-minus-a-prime}.
At such points the right-hand side of \eqref{eq:exact-euler-minus-a-prime} is nonnegative because
$b_t'\ge0$, $Z_t\ge0$, and $\p_\zeta Z_t=b_t^{-1}>0$.
\end{proof}

The first-variation argument uses the zero-extended axial function
$a_t=a_t^{\rm phys}\mathbf 1_{I_a}$ from \eqref{eq:euler-generated-truncated-coeff}.
The derivative of this discontinuous function is understood in the distributional sense on $(0,\infty)$. Specifically, 
for all $G\in C_c^1((0,\infty))$, we define the distributional derivative $\p_\zeta a_t$ as the distribution whose action on the test function $G$ is given by
\begin{equation}
\langle -\p_\zeta a_t,G\rangle
=\int_0^{\zeta_a}a_t^{\rm phys}(\zeta)\p_\zeta G(\zeta)\,d\zeta .
\label{eq:exact-euler-minus-a-prime-cutoff}
\end{equation}
Since $a_t^{\rm phys} \in C^1([0,\zeta_a])$, we can integrate by parts in \eqref{eq:exact-euler-minus-a-prime-cutoff} to obtain that
\[
\langle -\p_\zeta a_t,G\rangle =\int_0^{\zeta_a}G(\zeta)(-\p_\zeta a_t^{\rm phys})(\zeta)\,d\zeta +G(\zeta_a)a_t^{\rm phys}(\zeta_a).
\]
In particular, $\langle-\p_\zeta a_t,G\rangle\ge0$ whenever $G\ge0$.
Below we also use continuous functions $G$ on $[0,\zeta_a]$ which may be singular at the origin, provided that the right-hand side of
\eqref{eq:exact-euler-minus-a-prime-cutoff} is obtained as the following limit.  We choose
$G_\varepsilon\in C_c^1((0,\infty))$ such that
\[
G_\varepsilon\to G \quad\hbox{uniformly on compact subintervals of }(0,\zeta_a], \qquad G_\varepsilon(\zeta_a)\to G(\zeta_a),
\]
and such that the right-hand side of \eqref{eq:exact-euler-minus-a-prime-cutoff} converges.  Applying \eqref{eq:exact-euler-minus-a-prime-cutoff} to 
$G_\varepsilon$ and passing to the limit defines the same distributional pairing for $G$.

\subsection{The axis-trace velocity $V_\infty[a](\zeta)$}
By Lemma~\ref{lem:scaled-zeta-strain-moment} and \eqref{eq:C-alpha-M-limit}, for the fully angular vorticity function and the Biot--Savart velocity
\begin{equation} 
\Omega_\theta^{a,\infty}(\mathcal R,\mathcal Z):=-\operatorname{sgn}(\mathcal Z)a(|\mathcal Z|)\mathcal R^\alpha,\qquad
U^{a,\infty}:=\BS[\Omega_\theta^{a,\infty}e_\theta],
\label{eq:Omega-infi-U-infi}
\end{equation} 
the associated axial strain
$W_\infty[a] := \p_{\mathcal Z}(U^{a,\infty})_{\mathcal Z}(0,0)$ satisfies 
\begin{equation}
W_\infty[a]:=-C_\alpha^W I[a], 
\label{eq:W-model}
\end{equation}
where 
\begin{equation}
I[a]:=I_1[a]=\int_0^\infty a(\zeta)\zeta^{\alpha-1}\,d\zeta .
\label{eq:I-model}
\end{equation}
For an axial function $a$, we define
\begin{equation}
F_a(\zeta) :=\int_0^\infty a(\eta)\bigl((\zeta+\eta)^\alpha-|\zeta-\eta|^\alpha\bigr)\,d\eta ,
\label{eq:F-model}
\end{equation}
which is the one-dimensional integral appearing in the axis velocity below,
and we define the one-dimensional axial velocity 
\begin{equation}
V_\infty[a](\zeta):=-\tfrac{C_\alpha^W}{2\alpha}F_a(\zeta).
\label{eq:full-axis-velocity}
\end{equation}

\begin{lemma}[Axis-trace velocity]
\label{lem:full-angular-axis-trace}
For $\zeta>0$, the one-dimensional velocity $V_\infty[a]$ equals the axis trace of $U^{a,\infty}$:
\begin{equation}
V_\infty[a](\zeta) = (U^{a,\infty})_{\mathcal Z}(0,\zeta). 
\label{eq:full-axis-trace-identity}
\end{equation}
At the origin, its axial derivative is the model axial strain:
\begin{equation}
\p_\zeta V_\infty[a](0)=W_\infty[a].
\label{eq:full-axis-trace-strain}
\end{equation}
\end{lemma}

\begin{proof}[Proof of Lemma~\ref{lem:full-angular-axis-trace}]
By \eqref{eq:BS-axisymm} and \eqref{eq:Omega-infi-U-infi}, for $\zeta>0$, the  $\mathcal{Z}$-component of the three-dimensional BS-velocity is given by
\[
(U^{a,\infty})_{\mathcal Z}(0,\zeta) =\tfrac12\int_0^\infty a(\eta)\int_0^\infty \mathcal R^{\alpha+2}
\Big( \tfrac{1}{(\mathcal R^2+(\zeta+\eta)^2)^{3/2}} -\tfrac{1}{(\mathcal R^2+(\zeta-\eta)^2)^{3/2}} \Big)\,d\mathcal R\,d\eta .
\]
For $b,c\ge0$, the convergent difference integral satisfies
\begin{equation}
\int_0^\infty \mathcal R^{\alpha+2}\Big( \tfrac{1}{(\mathcal R^2+b^2)^{3/2}} -\tfrac{1}{(\mathcal R^2+c^2)^{3/2}} \Big)\,d\mathcal R
=-\tfrac{C_\alpha^W}{\alpha}(b^\alpha-c^\alpha).
\label{eq:axis-trace-difference-integral}
\end{equation}
Indeed, for $b,c>0$, differentiating the left side with respect to $b$ yields
$-3b\int_0^\infty \mathcal R^{\alpha+2}(\mathcal R^2+b^2)^{-5/2}\,d\mathcal R=-C_\alpha^W b^{\alpha-1}$, and the same
calculation with $c$ produces the opposite sign; the case $b=0$ or $c=0$ follows by taking a limit.  Applying
\eqref{eq:axis-trace-difference-integral} with $b=\zeta+\eta$ and $c=|\zeta-\eta|$ proves
\eqref{eq:full-axis-trace-identity}.  Then,  since $\tfrac{F_a(\zeta)}{\zeta}\longrightarrow 2\alpha I[a]$ as $\zeta\downarrow0$, 
\eqref{eq:full-axis-velocity} yields \eqref{eq:full-axis-trace-strain}.
\end{proof}

Equation \eqref{eq:scaled-strain-exact} gives
\[
W_M[a]=-C_{\alpha,M}^W I[a].
\]
With $F_a$ defined in \eqref{eq:F-model}, we define the one-dimensional velocity $V_M[a]$ by
\begin{equation}
V_M[a](\zeta):=-\tfrac{C_{\alpha,M}^W}{2\alpha}F_a(\zeta).
\label{eq:model-axial-velocity}
\end{equation}
From \eqref{eq:F-model},
\[
\tfrac{F_a(\zeta)}{\zeta}
=\int_0^\infty a(\eta)
\tfrac{(\zeta+\eta)^\alpha-|\zeta-\eta|^\alpha}{\zeta}\,d\eta .
\]
For each fixed $\eta>0$,
\[
\tfrac{(\zeta+\eta)^\alpha-|\zeta-\eta|^\alpha}{\zeta}
\longrightarrow 2\alpha\eta^{\alpha-1}\ \ \text{as } \ \ \zeta\downarrow0.
\]
For $\eta>2\zeta$, the mean value theorem gives
\[
\left|\tfrac{(\zeta+\eta)^\alpha-(\eta-\zeta)^\alpha}{\zeta}\right| \le C_\alpha\eta^{\alpha-1}.
\]
For $0<\eta\le2\zeta$, we use that
\[
(\zeta+\eta)^\alpha-|\zeta-\eta|^\alpha
\le
\begin{cases}
C_\alpha\eta\zeta^{\alpha-1}, & 0<\eta\le\tfrac12\zeta,\\ C_\alpha\zeta^\alpha, & \tfrac12\zeta<\eta\le2\zeta,
\end{cases}
\le C_\alpha\zeta\eta^{\alpha-1}.
\]
Hence, 
\[
0\le \int_0^{2\zeta}a(\eta) \tfrac{(\zeta+\eta)^\alpha-|\zeta-\eta|^\alpha}{\zeta}\,d\eta \le C_\alpha\int_0^{2\zeta}a(\eta)\eta^{\alpha-1}\,d\eta \longrightarrow0,
\]
by \eqref{eq:I-model}.  Combining this with the dominated convergence theorem on $\eta>2\zeta$ shows that 
\begin{equation}
F_a(\zeta)=2\alpha I[a]\zeta+o(\zeta)\ \ \text{as } \ \ \zeta\downarrow0.
\label{eq:Fa-origin-expansion}
\end{equation}
Therefore, by \eqref{eq:model-axial-velocity}, \eqref{eq:Fa-origin-expansion}, and \eqref{eq:scaled-strain-exact},
\begin{equation}
\p_\zeta V_M[a](0)=-C_{\alpha,M}^WI[a]=W_M[a].
\label{eq:model-velocity-strain-match}
\end{equation}

\subsection{The auxiliary first variation}
We now make the first-variation construction from Section~\ref{sec::Pi-model-explain} precise by constructing
a curve $s\mapsto a_{t,s}$ whose moment derivative is the one-dimensional expression $\mathcal D_M[a_t]$.
At a fixed physical time $t$, we define $s\mapsto a_{t,s}$, differentiable at $s=0$, such that
\[
\left.\tfrac{d}{ds}\right|_{s=0}I[a_{t,s}]=\mathcal D_M[a_t], \qquad I[a]=\int_0^\infty a(\zeta)\zeta^{\alpha-1}\,d\zeta,
\]
where $\mathcal D_M$ is the one-dimensional expression in \eqref{eq:D-model-def}.  For $M=\infty$, the same
construction produces $\mathcal D_\infty[a_t]$, and the full-angular pressure Hessian is recovered from
\[
\Pi_\infty[a]=C_\alpha^W\mathcal D_\infty[a]-\tfrac12W_\infty[a]^2,
\]
which is proved below in \eqref{eq:full-angular-pressure-hessian-moment-reduction}.  The computations in this
subsection are the finite-$M$ version of the Section~\ref{sec::Pi-model-explain} calculations
\eqref{eq:aux-specific-vorticity-characteristics}--\eqref{eq:delta-at}.

Throughout this subsection, we freeze the physical time $t$ and use the axial function
$a_t(\zeta)=a_t^{\rm phys}(\zeta)\mathbf 1_{I_a}(\zeta)$, where $a_t^{\rm phys}$ is defined in
\eqref{eq:scaled-profile-at} and $I_a$ is defined in \eqref{eq:pressure-profile-interval}.  We set
\[
V(\zeta):=V_M[a_t](\zeta).
\]
By \eqref{eq:model-axial-velocity}, \eqref{eq:Fa-origin-expansion}, and \eqref{eq:model-velocity-strain-match}, the one-dimensional velocity $V$ satisfies
$V(0)=0$ and $\p_\zeta V(0)=W_M[a_t]$.  We define the axial and radial auxiliary curves by
\begin{equation}
\tfrac{d}{ds}\zeta_s(\zeta_0)=V(\zeta_s(\zeta_0)),\qquad \tfrac{d}{ds}\mathcal R_s=-\tfrac12(\p_\zeta V)(\zeta_s)\mathcal R_s,\qquad
(\zeta_s,\mathcal R_s)\big|_{s=0}=(\zeta_0,\mathcal R_0).
\label{eq:model-aux-characteristics}
\end{equation}
The radial equation in \eqref{eq:model-aux-characteristics} is the axis linearization forced by
incompressibility: if an axisymmetric velocity has axial trace $V$, then
\[
u_{\mathcal R}(\mathcal R,\zeta)=-\tfrac12(\p_\zeta V)(\zeta)\mathcal R+o(\mathcal R),\qquad u_{\mathcal Z}(0,\zeta)=V(\zeta).
\]
This is the same incompressibility computation as \eqref{eq:aux-axis-velocity-expansion}--\eqref{eq:Rs-evo}, with $V_M[a_t]$ in place of $V_\infty[a_t]$.

We choose the curve $s\mapsto a_{t,s}$ by requiring the leading normalized specific vorticity $-a_{t,s}(\zeta)\mathcal R^{\alpha-1}$ to be conserved along
\eqref{eq:model-aux-characteristics}.  The curve is differentiable at $s=0$ in the distributional sense specified in \eqref{eq:model-delta-at-pairing-def} below.  
We denote this curve by
\begin{equation}
s\mapsto a_{t,s},\qquad a_{t,0}=a_t .
\label{eq:model-aux-curve}
\end{equation}
It is defined by
\begin{equation}
-a_{t,s}(\zeta_s(\zeta_0))\mathcal R_s^{\alpha-1} =-a_t(\zeta_0)\mathcal R_0^{\alpha-1},\qquad \zeta_0\in I_a,
\label{eq:model-leading-specific-vorticity-conservation}
\end{equation}
and we set $a_{t,s}(\zeta)=0$ for $\zeta\notin \zeta_s(I_a)$.  The radial equation in
\eqref{eq:model-aux-characteristics} is linear in $\mathcal R_s$, and hence
\[
\mathcal R_s=\mathcal R_0 \exp\left(-\tfrac12\int_0^s(\p_\zeta V)(\zeta_\sigma)\,d\sigma\right).
\]
Thus $\mathcal R_s/\mathcal R_0$ is independent of $\mathcal R_0$.  From
\eqref{eq:model-leading-specific-vorticity-conservation}, the axial function along the auxiliary curve satisfies
\begin{equation}
a_{t,s}(\zeta_s(\zeta_0)) =a_t(\zeta_0)\left(\tfrac{\mathcal R_0}{\mathcal R_s}\right)^{\!\alpha-1},
\qquad \zeta_0\in I_a .
\label{eq:model-ats-explicit}
\end{equation}
Equations \eqref{eq:model-leading-specific-vorticity-conservation} and
\eqref{eq:model-ats-explicit} are the finite-$M$ counterparts of
\eqref{eq:aux-leading-specific-vorticity-conservation} and \eqref{eq:as-axis-limit-transport}.
The radial equation in \eqref{eq:model-aux-characteristics} implies
\begin{equation}
\tfrac{d}{ds}\mathcal R_s^{\alpha-1} =-\tfrac{\alpha-1}{2}(\p_\zeta V)(\zeta_s)\mathcal R_s^{\alpha-1}.
\label{eq:model-Rs-power}
\end{equation}
Differentiating \eqref{eq:model-leading-specific-vorticity-conservation} and using
\eqref{eq:model-Rs-power}, we obtain
\begin{equation}
\tfrac{d}{ds}a_{t,s}(\zeta_s(\zeta_0)) =-\tfrac{1-\alpha}{2}(\p_\zeta V)(\zeta_s(\zeta_0))a_{t,s}(\zeta_s(\zeta_0)).
\label{eq:model-as-along-flow}
\end{equation}
This is the same leading specific-vorticity differentiation as \eqref{eq:ddsRs}--\eqref{eq:evo-as}, with $V_M[a_t]$ in place of $V_\infty[a_t]$.

We denote the derivative of \eqref{eq:model-aux-curve} at $s=0$ by
\[
\delta a_t:=\left.\p_s a_{t,s}\right|_{s=0}.
\]
Since $a_t$ contains the indicator function $\mathbf 1_{I_a}$, this derivative is understood in
$\mathcal D'((0,\infty))$: for every $G\in C_c^1((0,\infty))$,
\begin{equation} \left.\tfrac{d}{ds}\right|_{s=0}\int_0^\infty a_{t,s}(\zeta)G(\zeta)\,d\zeta =\langle \delta a_t,G\rangle .
\label{eq:model-delta-at-pairing-def}
\end{equation}
On the interior of $I_a$,
\[
\left.\tfrac{d}{ds}\right|_{s=0}a_{t,s}(\zeta_s(\zeta_0)) =\delta a_t(\zeta_0)+V(\zeta_0)\p_\zeta a_t(\zeta_0).
\]
Together with \eqref{eq:model-as-along-flow}, this yields the distributional identity
\begin{equation}
\delta a_t=-V_M[a_t]\p_\zeta a_t-\tfrac{1-\alpha}{2}\p_\zeta V_M[a_t]\,a_t .
\label{eq:model-delta-a}
\end{equation}
Formula \eqref{eq:model-delta-a} is the finite-$M$ analogue of \eqref{eq:delta-at}.
For $a_t=a_t^{\rm phys}\mathbf 1_{I_a}$, the derivative $-\p_\zeta a_t$ in \eqref{eq:model-delta-a} is the distribution defined in \eqref{eq:exact-euler-minus-a-prime-cutoff}.  
Combining \eqref{eq:model-delta-at-pairing-def} and \eqref{eq:model-delta-a}, we obtain
\begin{equation}
\left.\tfrac{d}{ds}\right|_{s=0}\int_0^\infty a_{t,s}(\zeta)G(\zeta)\,d\zeta =\left\langle-\p_\zeta a_t,V_M[a_t](\zeta)G(\zeta)\right\rangle
-\tfrac{1-\alpha}{2}\int_0^\infty \p_\zeta V_M[a_t](\zeta)a_t(\zeta)G(\zeta)\,d\zeta ,    \label{eq:model-curve-C1-pairing}
\end{equation}
whenever the two terms on the right side of \eqref{eq:model-curve-C1-pairing} are finite.  For $G(\zeta)=\zeta^{\alpha-1}$, the pairing with 
$V_M[a_t](\zeta)\zeta^{\alpha-1}$ is defined by the approximation convention following \eqref{eq:exact-euler-minus-a-prime-cutoff}; the integral containing
$\p_\zeta V_M[a_t]$ is finite by \eqref{eq:I-model}.

Taking $G(\zeta)=\zeta^{\alpha-1}$ in \eqref{eq:model-curve-C1-pairing}, we obtain
\begin{align}
\left.\tfrac{d}{ds}\right|_{s=0}I[a_{t,s}] &=\int_0^\infty\delta a_t(\zeta)\zeta^{\alpha-1}\,d\zeta  =\left\langle-\p_\zeta a_t,V_M[a_t](\zeta)\zeta^{\alpha-1}\right\rangle
-\tfrac{1-\alpha}{2}\int_0^\infty \p_\zeta V_M[a_t](\zeta)a_t(\zeta)\zeta^{\alpha-1}\,d\zeta .
\label{eq:model-DM-at}
\end{align}
For any nonnegative nonincreasing zero extension $a$ with $I[a]<\infty$, we define
\begin{align}
\mathcal D_M[a] &:=\left\langle-\p_\zeta a,V_M[a](\zeta)\zeta^{\alpha-1}\right\rangle -\tfrac{1-\alpha}{2}\int_0^\infty \p_\zeta V_M[a](\zeta)a(\zeta)\zeta^{\alpha-1}\,d\zeta .
\label{eq:D-model-def}
\end{align}
For $a=a_t$, equations \eqref{eq:model-DM-at} and \eqref{eq:D-model-def} identify the moment derivative:
\begin{equation} 
\mathcal D_M[a_t]=\left.\tfrac{d}{ds}\right|_{s=0}I[a_{t,s}].    \label{eq:Dm-from-dds}
\end{equation} 
Since \eqref{eq:scaled-strain-exact} gives $W_M[c]=-C_{\alpha,M}^WI[c]$ for every axial function $c$,
\[
\left.\tfrac{d}{ds}\right|_{s=0}W_M[a_{t,s}] =-C_{\alpha,M}^W\mathcal D_M[a_t].
\]
We define $\mathcal D_\infty[a]$ from \eqref{eq:D-model-def} by replacing $V_M$ with $V_\infty$ from
\eqref{eq:full-axis-velocity}.

\subsection{The pressure Hessian from the first variation}
The definition \eqref{eq:D-model-def} makes sense for every $M$, but the pressure Hessian is recovered from
the first variation only in the full-angular case $M=\infty$.  The reason is the axis-trace identity
\eqref{eq:full-axis-trace-identity}, applied to the frozen Euler-generated axial function $a_t$:
\[
(U^{a_t,\infty})_{\mathcal Z}(0,\zeta)=V_\infty[a_t](\zeta).
\]
Thus, for $M=\infty$, the velocity in the one-dimensional first variation is the actual axis trace of
$U^{a_t,\infty}$:
\[
V_\infty[a_t](\zeta)=(U^{a_t,\infty})_{\mathcal Z}(0,\zeta).
\]
The transport calculation for $\mathcal D_\infty[a_t]$ can therefore be compared directly with the Euler
evolution whose initial vorticity is $\Omega_\theta^{a_t,\infty}$ in \eqref{eq:Omega-infi-U-infi}.  This is the
rigorous version of \eqref{eq:section6-euler-tangent-axis-balance}--\eqref{eq:section6-euler-tangent}.

For the full-angular vorticity $\Omega_\theta^{a_t,\infty}$ obtained from \eqref{eq:Omega-infi-U-infi} by
setting $a=a_t$, we define
\begin{equation}
\Pi_\infty[a_t]
:=\pv\int_{\R^3}K_{zz}(\mathcal Y)\,
\tr\bigl(\nabla U^{a_t,\infty}(\mathcal Y)\nabla U^{a_t,\infty}(\mathcal Y)\bigr)\,d\mathcal Y .
\label{eq:section9-full-angular-pressure-hessian}
\end{equation}
We prove the full-angular first-variation identity
\begin{equation}
\Pi_\infty[a_t]=C_\alpha^W\mathcal D_\infty[a_t]-\tfrac12W_\infty[a_t]^2.
\label{eq:full-angular-pressure-hessian-moment-reduction}
\end{equation}
This is the Section~\ref{sec::Pi-model-explain} identity
\eqref{eq:section6-pressure-hessian-moment-reduction}, with
\eqref{eq:section9-full-angular-pressure-hessian} as the rigorous definition of $\Pi_\infty[a_t]$.
The calculation below is first read for a smooth compactly supported axial function.  For the zero-extended
Euler-generated function $a_t$, the occurrences of $\p_\zeta a_t$ are distributional pairings in the sense of
\eqref{eq:exact-euler-minus-a-prime-cutoff}, obtained by the same approximation convention stated after that
equation.
We start the Euler equation in the time variable $t'$ with initial vorticity
\begin{equation} 
\left.\Omega_\theta(\mathcal R,\mathcal Z,t')\right|_{t'=t} =\Omega_\theta^{a_t,\infty}(\mathcal R,\mathcal Z)
=-\operatorname{sgn}(\mathcal Z)a_t(|\mathcal Z|)\mathcal R^\alpha .    \label{eq:fake-Vort}
\end{equation} 
For $t'$ near $t$, we define $a_{t'}$ by the leading near-axis expansion of the specific vorticity
$\Omega_\theta(\mathcal R,\zeta,t')/\mathcal R$ on the positive half-axis $\zeta>0$:
\[
\tfrac{\Omega_\theta(\mathcal R,\zeta,t')}{\mathcal R} =-a_{t'}(\zeta)\mathcal R^{\alpha-1}+o(\mathcal R^{\alpha-1}) \qquad\text{as }\mathcal R\downarrow0.
\]
At the initial instant $t'=t$, this definition agrees with the frozen axial function $a_t$ in
\eqref{eq:fake-Vort}.  At $t'=t$, the specific-vorticity transport equation is
\[
-\left.\p_{t'}a_{t'}\right|_{t'=t}\mathcal R^{\alpha-1} +(U^{a_t,\infty})_{\mathcal R}(\mathcal R,\zeta)\p_{\mathcal R}\bigl(-a_t(\zeta)\mathcal R^{\alpha-1}\bigr)
+(U^{a_t,\infty})_{\mathcal Z}(\mathcal R,\zeta)\p_\zeta\bigl(-a_t(\zeta)\mathcal R^{\alpha-1}\bigr)
=o(\mathcal R^{\alpha-1})
\]
as $\mathcal R\downarrow0$.  By axis smoothness,
\[
(U^{a_t,\infty})_{\mathcal R}(\mathcal R,\zeta) =\p_{\mathcal R}(U^{a_t,\infty})_{\mathcal R}(0,\zeta)\mathcal R+o(\mathcal R),
\qquad
(U^{a_t,\infty})_{\mathcal Z}(\mathcal R,\zeta) =(U^{a_t,\infty})_{\mathcal Z}(0,\zeta)+o(1).
\]
Substituting these two axis expansions into the preceding transport identity, dividing by $\mathcal R^{\alpha-1}$, and taking  the limit as 
$\mathcal R\downarrow0$, we obtain
\[
-\left.\p_{t'}a_{t'}\right|_{t'=t} -(\alpha-1)\p_{\mathcal R}(U^{a_t,\infty})_{\mathcal R}(0,\zeta)a_t(\zeta) -(U^{a_t,\infty})_{\mathcal Z}(0,\zeta)\p_\zeta a_t(\zeta)=0.
\]
The two axis identities are
\[
(U^{a_t,\infty})_{\mathcal Z}(0,\zeta)=V_\infty[a_t](\zeta),\qquad \p_{\mathcal R}(U^{a_t,\infty})_{\mathcal R}(0,\zeta) =-\tfrac12\p_\zeta V_\infty[a_t](\zeta).
\]
The first identity is \eqref{eq:full-axis-trace-identity}; the second is the axis restriction of incompressibility used in \eqref{eq:axial-div-free}.  Therefore
\[
\left.\p_{t'}a_{t'}\right|_{t'=t} =-V_\infty[a_t]\p_\zeta a_t-\tfrac{1-\alpha}{2}(\p_\zeta V_\infty[a_t])a_t .
\]
Comparing this identity with \eqref{eq:model-delta-a} at $M=\infty$, we see that the two derivatives agree:
\[
\left.\p_{t'}a_{t'}\right|_{t'=t} = \left.\p_s a_{t,s}\right|_{s=0} \qquad\text{in }\mathcal D'((0,\infty)).
\]
Since
\[
W_\infty[a_{t'}]=-C_\alpha^W\int_0^\infty a_{t'}(\zeta)\zeta^{\alpha-1}\,d\zeta,
\]
we differentiate this identity with respect to $t'$ at $t'=t$ and use the $M=\infty$ version of
\eqref{eq:Dm-from-dds} to obtain
\[
\left.\tfrac{d}{dt'}\right|_{t'=t}W_\infty[a_{t'}] =-C_\alpha^W\mathcal D_\infty[a_t].
\]
The stagnation-point Riccati identity \eqref{eq:rW0}, applied at $t'=t$, shows that
\[
\left.\tfrac{d}{dt'}\right|_{t'=t}W_\infty[a_{t'}] =-\tfrac12W_\infty[a_t]^2-\Pi_\infty[a_t],
\]
and \eqref{eq:full-angular-pressure-hessian-moment-reduction} follows.  Thus the principal-value integral defining $\Pi_\infty[a_t]$ is recovered from the one-dimensional 
derivative $\mathcal D_\infty[a_t]$.  The first-variation calculation above did not use any property of $a_t$
beyond nonnegativity, monotonicity, and the zero extension convention.  Hence, for any such axial function $a$,
\eqref{eq:full-angular-pressure-hessian-moment-reduction} holds with $a_t$ replaced by $a$.  This is the form
used in \eqref{eq:section9-full-angular-riccati}.

\subsection{The one-dimensional \texorpdfstring{$\mathcal K_1$ and $\mathcal K_2$}{K1 and K2} estimates}
We now estimate $\mathcal D_M[a]$ using only that the zero extension of the axial function satisfies
$a\ge0$ and $-\p_\zeta a\ge0$ in the distributional sense.  By
\eqref{eq:D-model-def} and \eqref{eq:model-axial-velocity},
\[
V_M[a]=-\tfrac{C_{\alpha,M}^W}{2\alpha}F_a .
\]
We set
\begin{align}
\mathcal K_1[a]:=\int_0^\infty F_a(\zeta)a(\zeta)\zeta^{\alpha-2}\,d\zeta,\qquad \mathcal K_2[a]:=\left\langle-\p_\zeta a,F_a(\zeta)\zeta^{\alpha-1}\right\rangle .
\label{eq:K-model}
\end{align}
For a zero-extended nonnegative nonincreasing axial function on $[0,L]$, with $a\in C^1([0,L])$,
\begin{equation}
\mathcal K_2[a] =\int_0^LF_a(\zeta)\zeta^{\alpha-1}(-a'(\zeta))\,d\zeta +F_a(L)L^{\alpha-1}a(L). \label{eq:K2-smooth-form}
\end{equation}
We introduce $H_a$ so that $\mathcal K_1[a]$ and $\mathcal K_2[a]$ have a common form:
\begin{equation}
H_a(\zeta):=\tfrac{F_a(\zeta)}{\zeta}.
\label{eq:H-model}
\end{equation}
Then
\begin{equation}
\mathcal K_1[a]=\int_0^\infty H_a(\zeta)a(\zeta)\zeta^{\alpha-1}\,d\zeta,
\qquad
\mathcal K_2[a]=\left\langle-\p_\zeta a,H_a(\zeta)\zeta^\alpha\right\rangle .
\label{eq:K-Ha-rewrite}
\end{equation}

\begin{lemma}[One-dimensional $\mathcal K_1,\mathcal K_2$ estimates]
\label{lem:one-d-compression}
Let $0<\alpha<1$, and let $a\in L^\infty([0,L])$ be nonnegative, nonincreasing,
not identically zero, and extended by zero outside $[0,L]$.  Then
\begin{subequations}
\label{eq:K-compression-bounds}
\begin{align}
\alpha I[a]^2\le \mathcal K_1[a]\le 2\alpha I[a]^2,   \label{eq:K1-compression-bounds}  \\
0\le \mathcal K_2[a]\le \alpha\mathcal K_1[a].       \label{eq:K2-compression-bounds}
\end{align}
\end{subequations}
\end{lemma}

\begin{proof}[Proof of Lemma~\ref{lem:one-d-compression}]
We first prove the estimates for $a\in C^1([0,L])$ with $a'\le0$ and then pass to the stated monotone
function by one-dimensional approximation.

We first prove that $H_a$ in \eqref{eq:H-model} is nonnegative and nonincreasing.  For
$\ell>0$, define
\[
F_\ell(r):=\int_0^\ell\bigl((r+\eta)^\alpha-|r-\eta|^\alpha\bigr)\,d\eta .
\]
We claim that $F_\ell(r)/r$ is nonnegative and nonincreasing for $r>0$.  By scaling it is enough to prove
this for $\ell=1$.  For
\[
\Phi(r):=(1+\alpha)F_1(r),
\]
we have
\[
\Phi(r)=
\begin{cases}
(1+r)^{1+\alpha}-(1-r)^{1+\alpha}-2r^{1+\alpha}, & 0<r\le1,\\
(1+r)^{1+\alpha}+(r-1)^{1+\alpha}-2r^{1+\alpha}, & r\ge1.
\end{cases}
\]
The inequality $(F_1(r)/r)'\le0$ is the same as
$D(r):=\Phi(r)-r\Phi'(r)\ge0$.  If $0<r<1$, then
\[
D'(r)=\alpha(1+\alpha)r\left((1-r)^{\alpha-1}+2r^{\alpha-1}-(1+r)^{\alpha-1}\right)\ge0,
\qquad D(0)=0.
\]
If $r>1$, then
\[
D'(r)=-\alpha(1+\alpha)r\left((1+r)^{\alpha-1}+(r-1)^{\alpha-1}-2r^{\alpha-1}\right)\le0,
\]
because $x\mapsto x^{\alpha-1}$ is convex, and $D(r)\to0$ as $r\to\infty$.  Hence $D(r)\ge0$ for all
$r>0$.  Since $r+\eta\ge |r-\eta|$ and $s\mapsto s^\alpha$ is increasing, $F_\ell(r)\ge0$ as well.

We write the monotone function $a$ as a superposition of interval indicators.  For $0\le\eta\le L$,
\[
a(\eta)=a(L)+\int_\eta^L(-a'(\ell))\,d\ell .
\]
Substituting this identity into \eqref{eq:F-model}, we obtain
\[
F_a(\zeta)=a(L)F_L(\zeta)+\int_0^L F_\ell(\zeta)(-a'(\ell))\,d\ell .
\]
Thus $H_a(\zeta)=F_a(\zeta)/\zeta$ is a nonnegative superposition of nonnegative nonincreasing functions,
and therefore
\begin{equation}
\text{$H_a$ is nonnegative and nonincreasing on $(0,L]$.}
\label{eq:Ha-pos-inc}
\end{equation}
By \eqref{eq:Fa-origin-expansion},
\begin{equation}
0\le H_a(\zeta)\le \lim_{r\downarrow0}H_a(r)=2\alpha I[a],\qquad 0<\zeta\le L.
\label{eq:Ha-bound}
\end{equation}
The upper bound for $\mathcal K_1[a]$ in \eqref{eq:K1-compression-bounds} follows directly from
\eqref{eq:K-Ha-rewrite} and \eqref{eq:Ha-bound}:
\[
\mathcal K_1[a]=\int_0^\infty H_a(\zeta)a(\zeta)\zeta^{\alpha-1}\,d\zeta \le2\alpha I[a]^2 .
\]

We next prove the estimate \eqref{eq:K2-compression-bounds} for $\mathcal K_2[a]$.  The lower bound $\mathcal K_2[a]\ge0$ follows from
\eqref{eq:K-Ha-rewrite}, $-\p_\zeta a\ge0$, and $H_a(\zeta)\zeta^\alpha\ge0$.  To prove the upper bound, we first assume temporarily that
\begin{equation}
H_a\in C^1((0,L]). \label{eq:Ha-C1-temp}
\end{equation}
We let $G_a(\zeta):=H_a(\zeta)\zeta^\alpha$.  Since $G_a(0)=0$, integration by parts shows that $\mathcal K_2[a]=\int_0^La(\zeta)G_a'(\zeta)\,d\zeta$.
Using \eqref{eq:Ha-C1-temp},
\[
G_a'(\zeta)=\zeta^\alpha H_a'(\zeta)+\alpha H_a(\zeta)\zeta^{\alpha-1}.
\]
By \eqref{eq:Ha-pos-inc}, $H_a'\le0$, and therefore
\begin{equation}
\mathcal K_2[a]\le\alpha\int_0^LH_a(\zeta)a(\zeta)\zeta^{\alpha-1}\,d\zeta=\alpha\mathcal K_1[a].
\label{eq:K2-Ha-C1}
\end{equation}

We remove the temporary regularity assumption \eqref{eq:Ha-C1-temp}.  We extend $H_a$ to a nonincreasing function on $\mathbb R$ by
\[
\widetilde H_a(\zeta):=
\begin{cases}
2\alpha I[a], & \zeta\le0,\\  H_a(\zeta), & 0<\zeta\le L,\\ H_a(L), & \zeta>L.
\end{cases}
\]
Let $\rho_\varepsilon$ be a standard mollifier and set
$H_{a,\varepsilon}:=(\rho_\varepsilon*\widetilde H_a)|_{[0,L]}$.  Then, 
\[
\begin{gathered}
H_{a,\varepsilon}\in C^\infty([0,L]),\qquad H_{a,\varepsilon}'\le0,\qquad 0\le H_{a,\varepsilon}\le2\alpha I[a],\\
H_{a,\varepsilon}\to H_a\quad\text{locally uniformly on }(0,L].
\end{gathered}
\]
For $G_{a,\varepsilon}(\zeta):=H_{a,\varepsilon}(\zeta)\zeta^\alpha$, the integration-by-parts argument
leading to \eqref{eq:K2-Ha-C1} yields
\[
\int_0^La(\zeta)G_{a,\varepsilon}'(\zeta)\,d\zeta
\le \alpha\int_0^La(\zeta)H_{a,\varepsilon}(\zeta)\zeta^{\alpha-1}\,d\zeta .
\]
Passing to the limit using \eqref{eq:K2-smooth-form}, \eqref{eq:K-Ha-rewrite}, and the dominated convergence theorem, we obtain
$\mathcal K_2[a]\le\alpha\mathcal K_1[a]$.

It remains to prove the lower bound for $\mathcal K_1[a]$ in \eqref{eq:K1-compression-bounds}.  Substituting
\eqref{eq:F-model} into \eqref{eq:K-model} and symmetrizing the double integral, we find
\begin{equation}
\mathcal K_1[a]=\tfrac12\int_0^\infty\!\int_0^\infty a(\zeta)a(\eta) \bigl((\zeta+\eta)^\alpha-|\zeta-\eta|^\alpha\bigr) \bigl(\zeta^{\alpha-2}+\eta^{\alpha-2}\bigr)\,d\eta\,d\zeta .
\label{eq:K1-symmetrized}
\end{equation}
We claim that
\begin{equation}
\tfrac12\bigl(\zeta^{\alpha-2}+\eta^{\alpha-2}\bigr) \bigl((\zeta+\eta)^\alpha-|\zeta-\eta|^\alpha\bigr) \ge\alpha\zeta^{\alpha-1}\eta^{\alpha-1}.
\label{eq:K1-kernel-lower}
\end{equation}
By symmetry we may assume $\zeta\ge\eta$.  Since $r^{\alpha-1}$ is convex on $(0,\infty)$,
\[
(\zeta+\eta)^\alpha-(\zeta-\eta)^\alpha =\alpha\int_{\zeta-\eta}^{\zeta+\eta}r^{\alpha-1}\,dr \ge2\alpha\eta\,\zeta^{\alpha-1}.
\]
Multiplying by $\tfrac12(\zeta^{\alpha-2}+\eta^{\alpha-2})$ yields
\eqref{eq:K1-kernel-lower}.  Inserting \eqref{eq:K1-kernel-lower} into \eqref{eq:K1-symmetrized} proves
\[
\mathcal K_1[a]\ge\alpha I[a]^2.
\]

For a general nonnegative nonincreasing $a\in L^\infty([0,L])$, we choose nonnegative nonincreasing
$a_\varepsilon\in C^1([0,L_\varepsilon])$, with $L_\varepsilon\downarrow L$, whose zero extensions converge
to $a$ in the weighted integrals in \eqref{eq:I-model}, \eqref{eq:K-model}, and
\eqref{eq:K-Ha-rewrite}.  Applying the estimates already proved to $a_\varepsilon$ and passing to the limit
establishes \eqref{eq:K-compression-bounds}.
\end{proof}

\subsection{The model moment production}
\begin{lemma}[Model moment production]
\label{lem:model-moment-production}
For $0<\alpha<\tfrac13$, let $a\in L^\infty([0,L])$ be nonnegative, nonincreasing,
not identically zero, and extended by zero outside $[0,L]$.  Then
\begin{equation}
\tfrac{1-3\alpha}{4}C_{\alpha,M}^W I[a]^2
\le \mathcal D_M[a]
\le \tfrac{(1-\alpha)^2}{2}C_{\alpha,M}^W I[a]^2
<C_{\alpha,M}^W I[a]^2 .
\label{eq:D-two-sided}
\end{equation}
\end{lemma}

\begin{proof}[Proof of Lemma~\ref{lem:model-moment-production}]
We write $C_M=C_{\alpha,M}^W$.  After substituting $V_M[a]$ from
\eqref{eq:model-axial-velocity} into the definition of $\mathcal D_M[a]$ in
\eqref{eq:D-model-def}, the expression for $\mathcal D_M[a]$ is
\[
\mathcal D_M[a] =-\tfrac{C_M}{2\alpha}\mathcal K_2[a] +\tfrac{1-\alpha}{4\alpha}C_M \int_0^\infty F_a'(\zeta)a(\zeta)\zeta^{\alpha-1}\,d\zeta .
\]
Here $F_a$ is defined in \eqref{eq:F-model}. Since $a$ is supported in $[0,L]$, $F_a(0)=0$, and
$F_a(\zeta)\zeta^{\alpha-1}=O(\zeta^\alpha)$ as $\zeta \downarrow 0$, the functional $\mathcal K_2[a]$ in
\eqref{eq:K2-smooth-form} satisfies
\[
\mathcal K_2[a] =\int_0^\infty a(\zeta)\p_\zeta(F_a(\zeta)\zeta^{\alpha-1})\,d\zeta =\int_0^\infty F_a'(\zeta)a(\zeta)\zeta^{\alpha-1}\,d\zeta +(\alpha-1)\mathcal K_1[a].
\]
Therefore,
\begin{equation}
\mathcal D_M[a] =\tfrac{C_M}{4\alpha} \left((1-\alpha)^2\mathcal K_1[a]-(1+\alpha)\mathcal K_2[a]\right).      \label{eq:D-K-identity}
\end{equation}
By Lemma~\ref{lem:one-d-compression},
\[
\mathcal D_M[a]
\ge \tfrac{C_M}{4\alpha}\bigl((1-\alpha)^2-\alpha(1+\alpha)\bigr)\mathcal K_1[a] \ge \tfrac{1-3\alpha}{4}C_MI[a]^2.
\]
For the upper bound, $\mathcal K_2[a]\ge0$ and $\mathcal K_1[a]\le2\alpha I[a]^2$ imply
\[
\mathcal D_M[a]\le \tfrac{C_M}{4\alpha}(1-\alpha)^2\mathcal K_1[a] \le \tfrac{(1-\alpha)^2}{2}C_MI[a]^2.
\]
Since $0<\alpha<\tfrac13$, the last upper bound is strictly smaller than $C_MI[a]^2$.  This proves
\eqref{eq:D-two-sided}.
\end{proof}

Taking $M\to\infty$ in Lemma~\ref{lem:model-moment-production} and using
\eqref{eq:C-alpha-M-limit}, we obtain
\begin{equation}
\tfrac{1-3\alpha}{4}C_\alpha^W I[a]^2 \le \mathcal D_\infty[a] \le \tfrac{(1-\alpha)^2}{2}C_\alpha^W I[a]^2 .    \label{eq:D-infty-two-sided}
\end{equation}
For a generic axial function $a$, the first-variation identity
\eqref{eq:full-angular-pressure-hessian-moment-reduction} and the two-sided estimate
\eqref{eq:D-infty-two-sided} imply the full-angular pressure Hessian lower bound
\begin{equation}
\Pi_\infty[a]\ge -\tfrac{1+3\alpha}{2}\,\tfrac12 W_\infty[a]^2.
\label{eq:section9-full-angular-riccati}
\end{equation}

\subsection{From model Riccati to Euler Riccati}
We now convert the full-angular model estimate \eqref{eq:section9-full-angular-riccati} into the
Riccati estimate used for the Euler-generated axial function.  In \eqref{eq:D-model-def} with $M=\infty$,
the transport velocity is $V_\infty[a](\zeta)$.  In the Euler axial equation derived below, the corresponding
velocity is centered by subtracting its linear part at the origin.  Thus we replace $V_\infty[a](\zeta)$ by
$V_\infty[a](\zeta)-\zeta\p_\zeta V_\infty[a](0)$.  By \eqref{eq:full-axis-trace-strain}, we set
\begin{equation}
\widetilde V_\infty[a](\zeta) :=V_\infty[a](\zeta)-\zeta\p_\zeta V_\infty[a](0) =V_\infty[a](\zeta)-W_\infty[a]\zeta ,
\label{eq:renormalized-axis-velocity}
\end{equation}
and we define
\begin{equation}
\widetilde{\mathcal D}_\infty[a] :=\left\langle-\p_\zeta a,\widetilde V_\infty[a](\zeta)\zeta^{\alpha-1}\right\rangle -\tfrac{1-\alpha}{2}\int_0^\infty
\p_\zeta\widetilde V_\infty[a](\zeta)a(\zeta)\zeta^{\alpha-1}\,d\zeta .
\label{eq:D-renormalized-def}
\end{equation}

\begin{lemma}[Renormalized first variation formula]
\label{lem:renormalized-moment-production}
Let $0<\alpha<\tfrac13$, and let $a\in L^\infty([0,L])$ be nonnegative,
nonincreasing, not identically zero, and extended by zero outside $[0,L]$.  Then
\begin{equation}
\widetilde{\mathcal D}_\infty[a] =\mathcal D_\infty[a]+\tfrac{1-3\alpha}{2}W_\infty[a]I[a],  \label{eq:D-renormalized-algebra}
\end{equation}
and
\begin{equation}
\widetilde{\mathcal D}_\infty[a] \ge -\tfrac{1-3\alpha}{4}C_\alpha^W I[a]^2 . \label{eq:D-renormalized-lower}
\end{equation}
\end{lemma}

\begin{proof}[Proof of Lemma~\ref{lem:renormalized-moment-production}]
By \eqref{eq:D-renormalized-def} and \eqref{eq:renormalized-axis-velocity},
\[
\widetilde{\mathcal D}_\infty[a]-\mathcal D_\infty[a] =-W_\infty[a]\left\langle-\p_\zeta a,\zeta^\alpha\right\rangle +\tfrac{1-\alpha}{2}W_\infty[a]I[a].
\]
Since $a$ is nonnegative, compactly supported, and nonincreasing,
\[
\left\langle-\p_\zeta a,\zeta^\alpha\right\rangle=\alpha I[a].
\]
This proves \eqref{eq:D-renormalized-algebra}.  The lower bound
\eqref{eq:D-renormalized-lower} follows from \eqref{eq:D-infty-two-sided} and
$W_\infty[a]=-C_\alpha^WI[a]$.
\end{proof}

We next apply \eqref{eq:D-renormalized-def} to the Euler-generated axial function at a fixed time $t$.
We let $J=J_{\cusp}(t)$ and $a_t=a_t^{\rm phys}\mathbf 1_{I_a}$.  We assume $m(t)>0$ throughout this fixed-time argument; the
scalar-modulation bootstrap \eqref{eq:localized-modulation-m-bootstrap} supplies this positivity.  We normalize the axis trace of $U_{\cusp}$ by setting
\begin{equation}
\mathsf W_t(\zeta):=\Gamma^{-1}J^{1-3\alpha}\p_z(U_{\cusp})_z(0,J^2\zeta,t), \qquad \mathsf U_t(\zeta):=\int_0^\zeta\mathsf W_t(\eta)\,d\eta .
\label{eq:section9-axis-trace-def}
\end{equation}
Thus $\mathsf W_t$ is the normalized axial strain along the symmetry axis, and $\mathsf U_t$ is its
anti-derivative normalized to vanish at $\zeta=0$.
Since $(U_{\cusp})_z(0,0,t)=0$, the definitions in \eqref{eq:section9-axis-trace-def} imply
\begin{equation}
(U_{\cusp})_z(0,J^2\zeta,t)=\Gamma J^{3\alpha+1}\mathsf U_t(\zeta), \qquad \p_z(U_{\cusp})_z(0,J^2\zeta,t)=\Gamma J^{3\alpha-1}\mathsf W_t(\zeta).
\label{eq:section9-axis-U-reconstruction}
\end{equation}
By \eqref{eq:Verr-def}, $V_{\cusp}=m(t)U_{\cusp}+V_{\err}$.  We also use the clock identity
\begin{equation}
\tfrac{\dot J}{J}=\tfrac12 m(t)\mathcal W_{\cusp}(t) =\tfrac12m(t)\Gamma J^{3\alpha-1}\mathsf W_t(0).
\label{eq:section9-clock-local}
\end{equation}

We next compute the time derivative of the normalized axial coordinate
\[
\zeta=J^{-2}B_t(Z)
\]
at fixed label $Z$; this computation will determine the transport operator $\mathcal T_t$ in
\eqref{eq:operator-Tt}.  Since
$B_t(Z)=z_t(0,Z)$ and $\p_t\phi_{\cusp}=V_{\cusp}\circ\phi_{\cusp}$, we have
\[
\p_tB_t(Z)=(V_{\cusp})_z(0,B_t(Z),t).
\]
From \eqref{eq:section9-axis-U-reconstruction} and \eqref{eq:section9-clock-local}, the derivative is
\begin{align}
\p_t\bigl(J^{-2}B_t(Z)\bigr)
&=-2\tfrac{\dot J}{J}J^{-2}B_t(Z)+J^{-2}(V_{\cusp})_z(0,B_t(Z),t) \notag\\ 
&=m(t)\Gamma J^{3\alpha-1}
\bigl(\mathsf U_t(J^{-2}B_t(Z))-\mathsf W_t(0)J^{-2}B_t(Z)\bigr)   +J^{-2}(V_{\err})_z(0,B_t(Z),t).
\label{eq:section9-zeta-current}
\end{align}
We now rewrite \eqref{eq:section9-zeta-current} as an identity in the variable
$\zeta=J^{-2}B_t(Z)$.  We set
\begin{equation}
	\mathcal M_t:=m(t)\Gamma J^{3\alpha-1},\qquad \mathsf R_t^\zeta(\zeta):=J^{-2}(V_{\err})_z(0,J^2\zeta,t).
\label{eq:section9-Mt-Rzeta-def}
\end{equation}
With these definitions, \eqref{eq:section9-zeta-current} takes the form
\[
\p_t\bigl(J^{-2}B_t(Z)\bigr)
=\left[
\mathcal M_t\bigl(\mathsf U_t(\zeta)-\mathsf W_t(0)\zeta\bigr)
+\mathsf R_t^\zeta(\zeta)
\right]_{\zeta=J^{-2}B_t(Z)} .
\]
Consequently, for any differentiable $F(t,\zeta)$, the chain rule reads
\begin{align}
\frac{d}{dt}F\bigl(t,J^{-2}B_t(Z)\bigr)
&=\p_tF\bigl(t,J^{-2}B_t(Z)\bigr)
\,+\,
\p_\zeta F\bigl(t,J^{-2}B_t(Z)\bigr)\p_t\bigl(J^{-2}B_t(Z)\bigr) \notag\\
&=\left[\p_t+\bigl(\mathcal M_t(\mathsf U_t(\zeta)-\mathsf W_t(0)\zeta)+\mathsf R_t^\zeta(\zeta)\bigr)\p_\zeta\right]F\big|_{\zeta=J^{-2}B_t(Z)} .
\label{eq:Tt-chain-rule}
\end{align}
We introduce the operator
\begin{equation} 
\mathcal T_t:=\p_t+\bigl[\mathcal M_t(\mathsf U_t(\zeta)-\mathsf W_t(0)\zeta) +\mathsf R_t^\zeta(\zeta)\bigr]\p_\zeta .
\label{eq:operator-Tt}
\end{equation} 
With this notation, \eqref{eq:Tt-chain-rule} becomes
\[
\tfrac{d}{dt}F\bigl(t,J^{-2}B_t(Z)\bigr) =(\mathcal T_tF)\bigl(t,J^{-2}B_t(Z)\bigr).
\]
Since $Z_t(J^{-2}B_t(Z))=Z$, the identity above implies
\begin{equation}
\mathcal T_tZ_t=0.
\label{eq:renormalized-Z-label-current}
\end{equation}

We next derive \eqref{eq:renormalized-q-log-current}, the identity for $\mathcal T_t\log q_t$.
By \eqref{eq:At-Bt-def}, $A_t(Z)=\p_Rr_t(0,Z)$.
With the velocity $V_{\cusp}$ defined in \eqref{eq:cusp-map-def}, differentiating
$\p_t(\phi_{\cusp})_r= (V_{\cusp})_r(\phi_{\cusp},t)$ and evaluating at $R=0$ yields
\[
\p_t\p_Rr_t(0,Z) =\p_r(V_{\cusp})_r(0,B_t(Z),t)\p_Rr_t(0,Z) +\p_z(V_{\cusp})_r(0,B_t(Z),t)\p_Rz_t(0,Z).
\]
By smooth axisymmetry, $\p_Rz_t(0,Z)=0$, so that 
\[
\p_tA_t(Z)=\p_r(V_{\cusp})_r(0,B_t(Z),t)A_t(Z).
\]
Differentiating $q_t(J^{-2}B_t(Z))=JA_t(Z)$ in $t$ at fixed $Z$ yields
\[
\tfrac{d}{dt}\log q_t\bigl(J^{-2}B_t(Z)\bigr) =\tfrac{\dot J}{J}+\p_r(V_{\cusp})_r(0,B_t(Z),t).
\]
For the divergence-free axisymmetric field $U_{\cusp}$ defined in \eqref{eq:U-cusp-label}, the axis
incompressibility identity is $\p_r(U_{\cusp})_r(0,z,t)=-\tfrac12\p_z(U_{\cusp})_z(0,z,t)$.
Substituting this identity, together with
\eqref{eq:Verr-def}, \eqref{eq:section9-axis-U-reconstruction}, and \eqref{eq:section9-clock-local},
into the preceding logarithmic derivative yields
\begin{equation}
\mathcal T_t\log q_t =-\tfrac12\mathcal M_t\bigl(\mathsf W_t-\mathsf W_t(0)\bigr) +(\p_rV_{\err})_r(0,J^2\zeta,t).
\label{eq:renormalized-q-log-current}
\end{equation}

We now apply $\mathcal T_t$ to the formula for $a_t^{\rm phys}$ in \eqref{eq:scaled-profile-at}.  Since
$\log a_t^{\rm phys}=(1-\alpha)\log q_t-\tfrac{\gamma}{2}\log(1+Z_t^2)$, applying $\mathcal T_t$ and using
\[
\mathcal T_t\log(1+Z_t^2)=\tfrac{2Z_t}{1+Z_t^2}\mathcal T_tZ_t, 
\]
yields
\[
\mathcal T_t\log a_t^{\rm phys} =(1-\alpha)\mathcal T_t\log q_t-\gamma\tfrac{Z_t}{1+Z_t^2}\mathcal T_tZ_t.
\]
Multiplying by $a_t^{\rm phys}$ yields
\begin{equation} 
\mathcal T_t a_t^{\rm phys} =a_t^{\rm phys}\big((1-\alpha)\mathcal T_t\log q_t -\gamma\tfrac{Z_t}{1+Z_t^2}\mathcal T_tZ_t\big).  \label{eq:Tt-at}
\end{equation} 
By \eqref{eq:renormalized-Z-label-current},
\[
\mathcal T_t(1+Z_t^2)^{-\gamma/2}=0.
\]
Combining \eqref{eq:Tt-at} with \eqref{eq:renormalized-q-log-current}, we arrive at
\begin{equation}
\mathcal T_t a_t^{\rm phys} =-\tfrac{1-\alpha}{2}\mathcal M_t\bigl(\mathsf W_t-\mathsf W_t(0)\bigr)a_t^{\rm phys} +\mathsf E_t^a ,
\label{eq:renormalized-phys-a-current}
\end{equation}
where
\[
\mathsf E_t^a =(1-\alpha)(\p_rV_{\err})_r(0,J^2\zeta,t)a_t^{\rm phys}.
\]
Since $a_t=a_t^{\rm phys}\mathbf 1_{I_a}$, the identity for $a_t$ obtained from
\eqref{eq:operator-Tt} and \eqref{eq:renormalized-phys-a-current} is distributional:
\begin{align}
\p_ta_t =-m(t)\Gamma J^{3\alpha-1} \bigl(\mathsf U_t-\mathsf W_t(0)\zeta\bigr)\p_\zeta a_t 
-\tfrac{1-\alpha}{2}m(t)\Gamma J^{3\alpha-1}\bigl(\mathsf W_t-\mathsf W_t(0)\bigr)a_t+\mathcal R_t^a.
\label{eq:renormalized-profile-equation}
\end{align}
The term containing $\p_\zeta a_t$ is interpreted as follows: for every $H\in C_c^1((0,\infty))$,
\[
\left\langle-\p_\zeta a_t,H\right\rangle =\int_0^{\zeta_a}a_t^{\rm phys}(\zeta)\p_\zeta H(\zeta)\,d\zeta .
\]
Thus, for every $G\in C_c^1((0,\infty))$,
\begin{align}
\int_0^\infty \!\!\!\! \p_ta_t(\zeta)G(\zeta)\,d\zeta &=\mathcal M_t\left\langle-\p_\zeta a_t, \bigl(\mathsf U_t-\mathsf W_t(0)\zeta\bigr)G\right\rangle 
-\tfrac{1-\alpha}{2}\mathcal M_t\!
\int_0^\infty \!\!\!\bigl(\mathsf W_t-\mathsf W_t(0)\bigr)a_tG\,d\zeta + \left\langle\mathcal R_t^a,G\right\rangle ,
\label{eq:renormalized-profile-equation-tested}
\end{align}
where the remainder is the distribution
\begin{equation}
\mathcal R_t^a:=\mathsf E_t^a\mathbf 1_{I_a}-\mathsf R_t^\zeta\p_\zeta a_t .
\label{eq:renormalized-a-remainder-def}
\end{equation}
For a test function $G$, this means
\begin{equation}
\left\langle\mathcal R_t^a,G\right\rangle =\int_0^{\zeta_a}\mathsf E_t^a(\zeta)G(\zeta)\,d\zeta +\big\langle-\p_\zeta a_t,\mathsf R_t^\zeta G\big\rangle .
\label{eq:renormalized-a-remainder-tested}
\end{equation}

We now compute the derivative of $I[a_t]=\int_0^\infty a_t(\zeta)\zeta^{\alpha-1}\,d\zeta$.
We test \eqref{eq:renormalized-profile-equation-tested} with $G(\zeta)=\zeta^{\alpha-1}$.  As in the
convention following \eqref{eq:exact-euler-minus-a-prime-cutoff}, this test function is reached by
approximation, and we obtain
\begin{align}
\frac{d}{dt}I[a_t]
&=\mathcal M_t\left\langle-\p_\zeta a_t,\bigl(\mathsf U_t-\mathsf W_t(0)\zeta\bigr)\zeta^{\alpha-1}\right\rangle  -\tfrac{1-\alpha}{2}\mathcal M_t\!\!
\int_0^\infty\!\!\!
\bigl(\mathsf W_t-\mathsf W_t(0)\bigr)a_t\zeta^{\alpha-1}\,d\zeta   +\int_0^\infty \!\!\!\mathcal R_t^a(\zeta)\zeta^{\alpha-1}\,d\zeta .
\label{eq:I-variation-with-remainder}
\end{align}
By \eqref{eq:renormalized-a-remainder-tested}, the last term is
\begin{equation}
\int_0^\infty\mathcal R_t^a(\zeta)\zeta^{\alpha-1}\,d\zeta =\int_0^{\zeta_a}\mathsf E_t^a(\zeta)\zeta^{\alpha-1}\,d\zeta 
+\big\langle-\p_\zeta a_t,\mathsf R_t^\zeta\zeta^{\alpha-1}\big\rangle .
\label{eq:profile-remainder-moment-pairing}
\end{equation}
The definitions of $\mathsf E_t^a$ in \eqref{eq:renormalized-phys-a-current} and
$\mathsf R_t^\zeta$ in \eqref{eq:section9-Mt-Rzeta-def} show that the two terms in
\eqref{eq:profile-remainder-moment-pairing} involve only axis traces of $V_{\err}$, controlled by
\eqref{eq:localized-cusp-error-axis-trace}.  At the fixed time $t$, set
\begin{equation}
\varepsilon_{\rm rem}(t):= \tfrac{1}{m(t)\Gamma J^{3\alpha-1}I[a_t]^2} \Big|\int_0^\infty\mathcal R_t^a(\zeta)\zeta^{\alpha-1}\,d\zeta\Big| .
\label{eq:epsilon-rem}
\end{equation}
By \eqref{eq:epsilon-rem}, \eqref{eq:profile-remainder-moment-pairing} and
\eqref{eq:exact-euler-minus-a-prime-cutoff} yield
\begin{equation}
\Big|\int_0^\infty\mathcal R_t^a(\zeta)\zeta^{\alpha-1}\,d\zeta\Big| \le\varepsilon_{\rm rem}(t)m(t)\Gamma J^{3\alpha-1}I[a_t]^2 .
\label{eq:renormalized-profile-remainder}
\end{equation}
Combining \eqref{eq:I-variation-with-remainder} and \eqref{eq:renormalized-profile-remainder} yields the moment variation
\begin{align}
\frac{d}{dt}I[a_t]
&=\mathcal M_t \underbrace{
\Big[ \left\langle-\p_\zeta a_t, \bigl(\mathsf U_t-\mathsf W_t(0)\zeta\bigr)\zeta^{\alpha-1}\right\rangle -\tfrac{1-\alpha}{2}\int_0^\infty
\bigl(\mathsf W_t-\mathsf W_t(0)\bigr)a_t\zeta^{\alpha-1}\,d\zeta \Big]}_{=:\mathcal{B}_t}
+ O\!\left(\varepsilon_{\rm rem}(t)\mathcal M_t I[a_t]^2\right).
\label{eq:I-variation-axis-trace-current}
\end{align}
We compare $\mathcal B_t$ with $\widetilde{\mathcal D}_\infty[a_t]$ in
\eqref{eq:D-renormalized-def}.  The two differences to estimate are
\[
(\mathsf U_t-\mathsf W_t(0)\zeta)-\widetilde V_\infty[a_t],
\qquad
(\mathsf W_t-\mathsf W_t(0))-\p_\zeta\widetilde V_\infty[a_t].
\]

We first compare the normalized Euler functions $\mathsf U_t,\mathsf W_t$ from
\eqref{eq:section9-axis-trace-def} with the one-dimensional model functions
$V_\infty[a_t],\p_\zeta V_\infty[a_t]$ from \eqref{eq:full-axis-trace-identity}.  
By \eqref{eq:full-axis-trace-identity}, applied with $a=a_t$, and by differentiating the identity for
$\zeta>0$,
\[
V_\infty[a_t](\zeta)=(U^{a_t,\infty})_{\mathcal Z}(0,\zeta),
\qquad
\p_\zeta V_\infty[a_t](\zeta)=\p_{\mathcal Z}(U^{a_t,\infty})_{\mathcal Z}(0,\zeta).
\]
At $\zeta=0$, the second identity is read using the trace identity \eqref{eq:full-axis-trace-strain}.
We measure the difference between the normalized Euler axis traces and the full-angular model traces by
\begin{equation}
\begin{aligned}
	E_t^W(\zeta):=\mathsf W_t(\zeta)-\p_\zeta V_\infty[a_t](\zeta),
	\qquad
	E_t^U(\zeta):=\mathsf U_t(\zeta)-V_\infty[a_t](\zeta).
\end{aligned}
\label{eq:axis-trace-error-defs}
\end{equation}
The reconstruction formula \eqref{eq:section9-axis-U-reconstruction} and the axis-trace identity
\eqref{eq:full-axis-trace-identity} rewrite these errors as
\begin{align*} 
E_t^U(\zeta)&=\Gamma^{-1}J^{-3\alpha-1}(U_{\cusp})_z(0,J^2\zeta,t)-(U^{a_t,\infty})_{\mathcal Z}(0,\zeta) \\
E_t^W(\zeta) &=\Gamma^{-1}J^{1-3\alpha}\p_z(U_{\cusp})_z(0,J^2\zeta,t) -\p_{\mathcal Z}(U^{a_t,\infty})_{\mathcal Z}(0,\zeta),
\end{align*} 
and 
\[
\p_\zeta E_t^U(\zeta)=E_t^W(\zeta),\qquad E_t^U(0)=0,
\]
because $\p_\zeta\mathsf U_t=\mathsf W_t$ by \eqref{eq:section9-axis-trace-def},
$\mathsf U_t(0)=0$, and $V_\infty[a_t](0)=0$ by \eqref{eq:full-axis-velocity}.  Therefore,
\begin{equation}
E_t^U(\zeta)=\int_0^\zeta E_t^W(\eta)\,d\eta .
\label{eq:EU-from-EW}
\end{equation}
With $\varepsilon_{\rm rem}(t)$ defined in  \eqref{eq:epsilon-rem}, we define the fixed-time error
parameter
\begin{subequations}
\label{eq:axis-error-parameters}
\begin{align}
\varepsilon_{\rm tr}(t)&:=\tfrac{1}{I[a_t]}\sup_{0\le\zeta\le\zeta_a}|E_t^W(\zeta)|,
\notag\\
\varepsilon_{\rm str}(t)&:=\tfrac{1}{m(t)\Gamma^2J^{6\alpha-2}I[a_t]^2}
\left|\tfrac{d}{dt}\bigl(\Gamma J^{3\alpha-1}E_t^W(0)\bigr)\right|,
\label{eq:epsilon-strain-derivative}\\
\varepsilon_{\rm ax}=\varepsilon_{\rm ax}(t)&:=\max\{\varepsilon_{\rm rem}(t),\varepsilon_{\rm tr}(t),\varepsilon_{\rm str}(t)\}.
\label{eq:epsilon-ax}
\end{align}
\end{subequations}
By \eqref{eq:axis-error-parameters}, for $0\le\zeta\le\zeta_a$,
\begin{equation}
|E_t^W(\zeta)|\le \varepsilon_{\rm tr}(t)I[a_t]\le\varepsilon_{\rm ax}I[a_t].
\label{eq:EW-bound-from-eps}
\end{equation}
For $0<\zeta\le\zeta_a$, \eqref{eq:EU-from-EW} and \eqref{eq:EW-bound-from-eps} give
\begin{equation}
\tfrac{|E_t^U(\zeta)|}{\zeta} \le \sup_{0\le\eta\le\zeta}|E_t^W(\eta)| \le \varepsilon_{\rm ax}I[a_t].
\label{eq:EU-bound-from-EW}
\end{equation}
Therefore, by \eqref{eq:EW-bound-from-eps} and \eqref{eq:EU-bound-from-EW}, the Euler and full-angular
model axis traces satisfy
\begin{subequations}
\label{eq:axis-trace-approximation}
\begin{align}
\sup_{0<\zeta\le\zeta_a} \tfrac{|\mathsf U_t(\zeta)-V_\infty[a_t](\zeta)|}{\zeta} &\le \varepsilon_{\rm ax}I[a_t],
\label{eq:axis-trace-U-approximation}\\
\sup_{0\le\zeta\le\zeta_a} |\mathsf W_t(\zeta)-\p_\zeta V_\infty[a_t](\zeta)| &\le \varepsilon_{\rm ax}I[a_t], 
\label{eq:axis-trace-W-approximation}
\end{align}
\end{subequations}
Also, $\varepsilon_{\rm rem}(t)\le\varepsilon_{\rm ax}$ by \eqref{eq:epsilon-ax}, so the remainder in
\eqref{eq:I-variation-axis-trace-current} is $O(\varepsilon_{\rm ax}\mathcal M_tI[a_t]^2)$.
Since $W_\infty[a_t]=\p_\zeta V_\infty[a_t](0)$ by \eqref{eq:full-axis-trace-strain},
\eqref{eq:axis-trace-approximation} implies the centered estimates
\begin{subequations}
\label{eq:centered-axis-trace-approximation}
\begin{align}
\sup_{0<\zeta\le\zeta_a}
	\tfrac{\left|(\mathsf U_t(\zeta)-\mathsf W_t(0)\zeta)-\widetilde V_\infty[a_t](\zeta)\right|}{\zeta}
	&\le C\varepsilon_{\rm ax}I[a_t],
	\\
	\sup_{0\le\zeta\le\zeta_a}
	\left|(\mathsf W_t(\zeta)-\mathsf W_t(0))-\p_\zeta\widetilde V_\infty[a_t](\zeta)\right|
	&\le C\varepsilon_{\rm ax}I[a_t].
	\end{align}
	\end{subequations}
We now estimate $\mathcal B_t-\widetilde{\mathcal D}_\infty[a_t]$, with $\mathcal B_t$ defined in
\eqref{eq:I-variation-axis-trace-current} and $\widetilde{\mathcal D}_\infty[a_t]$ defined in
\eqref{eq:D-renormalized-def}.  The first bound in \eqref{eq:centered-axis-trace-approximation} implies
\[
\left|\left\langle-\p_\zeta a_t,
\bigl((\mathsf U_t-\mathsf W_t(0)\zeta)-\widetilde V_\infty[a_t]\bigr)\zeta^{\alpha-1}
\right\rangle\right|
\le C\varepsilon_{\rm ax}I[a_t]\left\langle-\p_\zeta a_t,\zeta^\alpha\right\rangle
=C\varepsilon_{\rm ax}I[a_t]^2,
\]
where $\left\langle-\p_\zeta a_t,\zeta^\alpha\right\rangle=\alpha I[a_t]$ by the monotonicity of $a_t$.
The second bound in \eqref{eq:centered-axis-trace-approximation} gives
\[
\int_0^\infty
\left|(\mathsf W_t-\mathsf W_t(0))-\p_\zeta\widetilde V_\infty[a_t]\right|
a_t\zeta^{\alpha-1}\,d\zeta
\le C\varepsilon_{\rm ax}I[a_t]^2.
\]
Therefore, using $\mathcal M_t=m(t)\Gamma J^{3\alpha-1}$ from
\eqref{eq:section9-Mt-Rzeta-def} and $\varepsilon_{\rm rem}(t)\le\varepsilon_{\rm ax}$ from
\eqref{eq:epsilon-ax}, \eqref{eq:I-variation-axis-trace-current} becomes
\begin{equation}
\frac{d}{dt}I[a_t]
=m(t)\Gamma J^{3\alpha-1}\widetilde{\mathcal D}_\infty[a_t]
+O\!\left(\varepsilon_{\rm ax}m(t)\Gamma J^{3\alpha-1}I[a_t]^2\right).
\label{eq:I-renormalized-variation}
\end{equation}
At $\zeta=0$, the second identity in \eqref{eq:section9-axis-U-reconstruction} relates the Euler axial strain to
$\mathsf W_t(0)$.  By
\eqref{eq:W-cusp-def} and \eqref{eq:section9-axis-U-reconstruction},
\[
\mathcal W_{\cusp}(t)=\Gamma J^{3\alpha-1}\mathsf W_t(0).
\]
Since $W_\infty[a_t]=\p_\zeta V_\infty[a_t](0)$ by \eqref{eq:full-axis-trace-strain}, the definition
of $E_t^W$ in \eqref{eq:axis-trace-error-defs}, evaluated at $\zeta=0$, is equivalent to
\begin{equation}
\mathcal W_{\cusp}(t)-\Gamma J^{3\alpha-1}W_\infty[a_t]
=\Gamma J^{3\alpha-1}E_t^W(0).
\label{eq:Wcusp-model-error-EW}
\end{equation}
Therefore, \eqref{eq:EW-bound-from-eps} with $\zeta=0$ gives the strain comparison below, while
\eqref{eq:section9-clock-local} multiplied by $J$ gives the clock law below:
\begin{subequations}
\label{eq:axis-strain-clock-comparison}
\begin{equation} 
\left|\mathcal W_{\cusp}(t)-\Gamma J^{3\alpha-1}W_\infty[a_t]\right| \le \varepsilon_{\rm ax}\Gamma J^{3\alpha-1}I[a_t],
\label{eq:axis-strain-undifferentiated-comparison}
\end{equation} 
\begin{equation} 
\dot J=\tfrac12m(t)J\mathcal W_{\cusp}(t).
\label{eq:axis-clock-law-for-criterion}
\end{equation} 
\end{subequations}
We now differentiate \eqref{eq:Wcusp-model-error-EW}, written as
\[
\mathcal W_{\cusp}(t)=\Gamma J^{3\alpha-1}W_\infty[a_t]+\Gamma J^{3\alpha-1}E_t^W(0).
\]
The model identity $W_\infty[a_t]=-C_\alpha^WI[a_t]$ in \eqref{eq:W-model} shows that
\begin{align}
\tfrac{d}{dt}\bigl(\Gamma J^{3\alpha-1}W_\infty[a_t]\bigr) &=-C_\alpha^W\Gamma J^{3\alpha-1}\tfrac{d}{dt}I[a_t]  +(3\alpha-1)\Gamma J^{3\alpha-2}\dot J\,W_\infty[a_t].
\label{eq:leading-strain-derivative}
\end{align}
By \eqref{eq:epsilon-strain-derivative} and \eqref{eq:epsilon-ax},
\begin{equation}
\left|\tfrac{d}{dt}\bigl(\Gamma J^{3\alpha-1}E_t^W(0)\bigr)\right|
\le \varepsilon_{\rm ax}m(t)\Gamma^2J^{6\alpha-2}I[a_t]^2.
\label{eq:strain-error-derivative-bound}
\end{equation}
Combining \eqref{eq:Wcusp-model-error-EW}, \eqref{eq:leading-strain-derivative}, and
\eqref{eq:strain-error-derivative-bound}, the derivative of $\mathcal W_{\cusp}(t)$ satisfies
\begin{align}
\frac{d}{dt}\mathcal W_{\cusp}(t) &=-C_\alpha^W\Gamma J^{3\alpha-1}\frac{d}{dt}I[a_t] +(3\alpha-1)\Gamma J^{3\alpha-2}\dot J\,W_\infty[a_t]
+O\!\left(\varepsilon_{\rm ax}m(t)\Gamma^2J^{6\alpha-2}I[a_t]^2\right).
\label{eq:direct-strain-derivative-axis}
\end{align}

\begin{lemma}[Renormalized axis-trace derivative]
\label{lem:renormalized-axis-trace-derivative}
Assume the axis-trace estimates \eqref{eq:axis-trace-approximation} and the strain-clock estimates
\eqref{eq:axis-strain-clock-comparison}.  Then
\begin{align}
\frac{d}{dt}\mathcal W_{\cusp}(t)
&=-m(t)\Gamma^2J^{6\alpha-2} \left[ C_\alpha^W\widetilde{\mathcal D}_\infty[a_t]+\tfrac{1-3\alpha}{2}W_\infty[a_t]^2 \right] 
+O\!\left(\varepsilon_{\rm ax}\,m(t)\Gamma^2J^{6\alpha-2}I[a_t]^2\right).
\label{eq:renormalized-axis-derivative}
\end{align}
\end{lemma}

\begin{proof}[Proof of Lemma~\ref{lem:renormalized-axis-trace-derivative}]
After substituting the moment identity \eqref{eq:I-renormalized-variation} into
\eqref{eq:direct-strain-derivative-axis}, it remains to rewrite the $\dot J$ term in
\eqref{eq:direct-strain-derivative-axis} in terms of $W_\infty[a_t]$.  From \eqref{eq:axis-clock-law-for-criterion} and
\eqref{eq:axis-strain-undifferentiated-comparison},
\begin{equation}
\dot J=\tfrac12m(t)\Gamma J^{3\alpha}W_\infty[a_t]
+O\!\left(\varepsilon_{\rm ax}m(t)\Gamma J^{3\alpha}I[a_t]\right).
\label{eq:Jdot-model-strain}
\end{equation}
Multiplying \eqref{eq:Jdot-model-strain} by
$(3\alpha-1)\Gamma J^{3\alpha-2}W_\infty[a_t]$ and using
$W_\infty[a_t]=-C_\alpha^WI[a_t]$ from \eqref{eq:W-model} yields
\begin{equation}
(3\alpha-1)\Gamma J^{3\alpha-2}\dot J\,W_\infty[a_t]
=-\tfrac{1-3\alpha}{2}m(t)\Gamma^2J^{6\alpha-2}W_\infty[a_t]^2
+O\!\left(\varepsilon_{\rm ax}m(t)\Gamma^2J^{6\alpha-2}I[a_t]^2\right).
\label{eq:Jdot-model-strain-term}
\end{equation}
Substituting \eqref{eq:I-renormalized-variation} and \eqref{eq:Jdot-model-strain-term} into
\eqref{eq:direct-strain-derivative-axis} proves \eqref{eq:renormalized-axis-derivative}.
\end{proof}

We choose the Riccati constants by
\begin{equation}
q_{\alpha,0}:=\tfrac{1+3\alpha}{2}<1, \qquad q_\alpha\in(q_{\alpha,0},1).
\label{eq:pressure-subcritical-riccati}
\end{equation}
We also fix $M_{\rm pos}\ge1$ large enough for the angular-tail estimates used later; the Riccati sign below is
independent of this choice.

\begin{proposition}[Euler-generated renormalized Riccati bound]
\label{prop:euler-generated-profile-riccati}
Let $0<\alpha<\tfrac13$, and let $a_t$ be the zero-extended Euler-generated axial function
\eqref{eq:euler-generated-truncated-coeff}.  Assume that $a_t$ is nonnegative and nonincreasing, that
$m(t)>0$, that \eqref{eq:axis-trace-approximation}--\eqref{eq:axis-strain-clock-comparison} hold with
$\varepsilon_{\rm ax}$ sufficiently small, and that $\mathcal W_{\cusp}$  satisfies the principal Riccati identity with a controlled error:
\begin{equation}
\frac{d}{dt}\mathcal W_{\cusp}(t)
=m(t)\left(-\tfrac12\mathcal W_{\cusp}(t)^2-\Pi_{\cusp}(t)\right)
+O\!\left(\varepsilon_{\rm ax}m(t)\mathcal W_{\cusp}(t)^2\right).
\label{eq:cusp-principal-riccati-with-error}
\end{equation}
Then
\begin{equation}
\Pi_{\cusp}(t) \ge -q_\alpha\tfrac12\,\mathcal W_{\cusp}(t)^2, \qquad q_\alpha<1.
\label{eq:euler-generated-riccati-bound}
\end{equation}
\end{proposition}

\begin{proof}[Proof of Proposition~\ref{prop:euler-generated-profile-riccati}]
By \eqref{eq:D-renormalized-lower} and $W_\infty[a_t]=-C_\alpha^WI[a_t]$,
\[
C_\alpha^W\widetilde{\mathcal D}_\infty[a_t]
+\tfrac{1-3\alpha}{2}W_\infty[a_t]^2
\ge \tfrac{1-3\alpha}{4}(C_\alpha^W)^2I[a_t]^2 .
\]
Substituting the lower bound for
$C_\alpha^W\widetilde{\mathcal D}_\infty[a_t]+\tfrac{1-3\alpha}{2}W_\infty[a_t]^2$
into \eqref{eq:renormalized-axis-derivative} yields
\[
\frac{d}{dt}\mathcal W_{\cusp}(t) \le -m(t)\Gamma^2J^{6\alpha-2} \left(\tfrac{1-3\alpha}{4}(C_\alpha^W)^2I[a_t]^2 -C\varepsilon_{\rm ax}I[a_t]^2\right).
\]
The strain comparison \eqref{eq:axis-strain-undifferentiated-comparison}, together with
$W_\infty[a_t]=-C_\alpha^WI[a_t]$, yields
\[
\Gamma^2J^{6\alpha-2}(C_\alpha^W)^2I[a_t]^2 =\mathcal W_{\cusp}(t)^2+O\!\left(\varepsilon_{\rm ax}\mathcal W_{\cusp}(t)^2\right).
\]
Therefore, 
\begin{equation}
\frac{d}{dt}\mathcal W_{\cusp}(t) \le -m(t)\left(\tfrac{1-3\alpha}{4}-C\varepsilon_{\rm ax}\right)\mathcal W_{\cusp}(t)^2 .
\label{eq:renormalized-W-upper}
\end{equation}
Since $m(t)>0$, comparison of \eqref{eq:renormalized-W-upper} with
\eqref{eq:cusp-principal-riccati-with-error} yields
\[
\Pi_{\cusp}(t) \ge -\left(\tfrac{1+3\alpha}{2}+C\varepsilon_{\rm ax}\right)\tfrac12\mathcal W_{\cusp}(t)^2 .
\]
In Proposition~\ref{prop:small-clock-comparisons}, the final threshold
$\mathfrak J_{\mathrm{axis}}$ is chosen so that the axis-trace error satisfies
\[
C\varepsilon_{\rm ax}\le q_\alpha-q_{\alpha,0}.
\]
Then $\tfrac{1+3\alpha}{2}+C\varepsilon_{\rm ax}\le q_\alpha<1$, and this proves
\eqref{eq:euler-generated-riccati-bound}.
\end{proof}

The estimate for the first variation $\widetilde{\mathcal D}_\infty[a_t]$ used in
Proposition~\ref{prop:euler-generated-profile-riccati} is the lower bound
\eqref{eq:D-renormalized-lower}.  This lower bound is a consequence of the two one-dimensional inequalities
for $\mathcal K_1[a]$ and $\mathcal K_2[a]$ in Lemma~\ref{lem:one-d-compression}.  After
\eqref{eq:D-renormalized-lower} is inserted into the strain derivative formula
\eqref{eq:renormalized-axis-derivative}, the strain comparison
\eqref{eq:axis-strain-undifferentiated-comparison} converts the resulting bound from $I[a_t]^2$ to
$\mathcal W_{\cusp}(t)^2$, producing \eqref{eq:renormalized-W-upper}.  The comparison between
\eqref{eq:renormalized-W-upper} and the Riccati identity \eqref{eq:cusp-principal-riccati-with-error}
then yields the pressure Hessian lower bound \eqref{eq:euler-generated-riccati-bound}.

In Proposition~\ref{prop:euler-generated-profile-riccati}, the limit $M\to\infty$ has already been taken in
\eqref{eq:D-infty-two-sided}.  The cutoff $\chi_M$ in \eqref{eq:slope-cutoff-chiM} is retained for the later
pressure Hessian localization.  In the pressure Hessian estimates, $\chi_M(\mathcal R/|\mathcal Z|)$ selects
the bounded-slope region, while $1-\chi_M(\mathcal R/|\mathcal Z|)$ selects the large-slope region.  The
axial region outside $I_\sharp$ is controlled by the $\zeta$-tail $\mathfrak a_\zeta(I_\sharp)$ in
\eqref{eq:pressure-zeta-tail}.  The contribution with angular cutoff
$1-\chi_M(\mathcal R/|\mathcal Z|)$ is controlled for the axial strain by
\eqref{eq:linear-angular-tail-strain} and for the pressure Hessian by \eqref{eq:bilinear-angular-tail}.

\section{Hyperbolic Normal Form for the Cusp Flow in the Collapse Limit}
\label{sec:current-axis-normal-forms}

\subsection{Axial flow map coordinates for the cusp-flow normal form}

We prove a normal form for the exact cusp map near the symmetry axis.  The map is first written in
Lagrangian labels $(R,Z)$ and then evaluated through the Eulerian image
$(r,z)=\phi_{\cusp}(R,Z,t)$ after this image is divided by the cusp-clock scale $J^2$.  The label variables are
\begin{equation*}
\zeta=J^{-2}B_t(Z),\qquad \tau=\tfrac{A_t(Z)R}{B_t(Z)},
\end{equation*}
which are introduced in \eqref{eq:adapted-labels}.  In these variables, \eqref{eq:additive-normal-form} states
that $J^{-2}(r,z)$ is the model point $\zeta(\tau,1)$ plus the error controlled in
\eqref{eq:additive-normal-form-bound}.  This normal form is the geometric step used in
Section~\ref{sec:transported-cusp-pressure} to compare the pressure Hessian generated by the transported cusp
vorticity with the $M$-slope-restricted model pressure Hessian.

The compact interval $I_\sharp\Subset(0,\infty)$ is the pressure Hessian localization interval fixed in
Section~\ref{sec:fixed-choice-order}.  In the estimates below, $I_{\rm buf}$ denotes a slightly larger positive
$\zeta$-interval used to propagate the normal-form bounds \eqref{eq:additive-normal-form-bound} up to
$I_\sharp$.  These intervals are separated from $\zeta=0$ for the same reason explained in
Section~\ref{sec:fixed-choice-order}: the slope variable in \eqref{eq:adapted-labels} and the variables
$(R_{\rm sc},Z_{\rm sc})$ in \eqref{eq:normal-form-fixed-scaled-set}, which describe the Eulerian image after
division by $J^2$, are not useful at the stagnation point.  The axial function attached to the stagnation point
is handled instead by the Riccati estimate for the Euler-generated monotone axial function in
Proposition~\ref{prop:euler-generated-profile-riccati}.

We now fix the notation used for the radial flatness estimates \eqref{eq:axis-normal-form-entry} and
\eqref{eq:axis-normal-form-zeta-derivative}.  We let
\[
I_\sharp\Subset I_{\rm buf}\Subset(0,\infty),\qquad C_0\ge1,
\]
and fix a small-clock time $t$.  Set $J:=J_{\cusp}(t)$ and write
\[
\phi_{\cusp}(R,Z,t)=\bigl(r_t(R,Z),z_t(R,Z)\bigr), \qquad A_t(Z):=\p_Rr_t(0,Z), \qquad B_t(Z):=z_t(0,Z),
\]
so that $A_t$ and $B_t$ are the radial stretch and axial position on the symmetry axis.  On the buffered
interval we use the inverse of the monotone axial map $Z\mapsto J^{-2}B_t(Z)$, and the terminal axial-label
interval associated to $I_\sharp$ is
\[
Z_t:=(J^{-2}B_t)^{-1}\quad\hbox{on }I_{\rm buf}, \qquad I_Z^t:=Z_t(I_\sharp).
\]
For $\zeta\in I_\sharp$ and $|\tau|\le C_0$, the label with time-$t$ axial flow map slope $\tau$ is
\[
R_t^\sharp(\zeta,\tau):= \tfrac{J^2\zeta\tau}{A_t(Z_t(\zeta))}.
\]
The nonlinear remainders relative to the leading hyperbolic axis normal form are
\[
\mathscr R_{r,t}(R,Z):=r_t(R,Z)-A_t(Z)R, \qquad \mathscr R_{z,t}(R,Z):=z_t(R,Z)-B_t(Z).
\]

The radial estimate \eqref{eq:axis-normal-form-entry} and the source estimate
\eqref{eq:radial-flatness-source} use the axis-geometry bootstrap from \textup{(BA2)} on the time interval
$[t_0,t]$.  For each $s\in[t_0,t]$, write
$J_s:=J_{\cusp}(s)$,
$A_s(Z):=\p_Rr_s(0,Z)$, $B_s(Z):=z_s(0,Z)$, and
$Z_s^{\rm buf}:=(J_s^{-2}B_s)^{-1}$ on $I_{\rm buf}$.  On $I_{\rm buf}$, this bound has the form
\begin{equation}
\begin{aligned}
c_{\rm ax}\le J_s A_s(Z_s^{\rm buf}(\zeta))\le C_{\rm ax},\qquad c_{\rm ax}\le J_s^{-2}B_s'(Z_s^{\rm buf}(\zeta))\le C_{\rm ax},\\
[\log(J_s A_s(Z_s^{\rm buf}(\cdot)))]_{C^{\alpha/2}(I_{\rm buf})} + [\log(J_s^{-2}B_s'(Z_s^{\rm buf}(\cdot)))]_{C^{\alpha/2}(I_{\rm buf})} \le C_{\rm ax},
\end{aligned}
\qquad \zeta\in I_{\rm buf}.
\label{eq:buffered-axis-geometry-bootstrap}
\end{equation}
We also use the clock bootstrap \eqref{eq:localized-clock-bootstrap}, the large normal-form bootstrap
\eqref{eq:localized-normal-form-large-bootstrap}, and the large cusp-error bootstrap
\eqref{eq:localized-cusp-error-large-bootstrap} on the same interval.  The fixed-label containment used below is
\eqref{eq:radial-flatness-buffered-label} from \textup{(BA3)}.

\begin{lemma}[Hyperbolic-deviation bounds for the Euler cusp map]
\label{lem:nonlinear-radial-flatness}
Assume the bootstrap bounds \eqref{eq:localized-clock-bootstrap},
\eqref{eq:localized-normal-form-large-bootstrap}, \eqref{eq:localized-cusp-error-large-bootstrap}, and
\eqref{eq:buffered-axis-geometry-bootstrap}, together with the containment condition
\eqref{eq:radial-flatness-buffered-label}.  After decreasing the small-clock threshold depending only on the
fixed parameters and the large bootstrap constants, for $Z\in I_Z^t$ and $|R|\le C_{\rm ax}J^3$, the
hyperbolic-deviation remainders in \eqref{eq:additive-normal-form} satisfy
\begin{equation}
|\mathscr R_{r,t}(R,Z)|+|\mathscr R_{z,t}(R,Z)| \le C_{\rm ax}J^{-1}|R|^{1+\beta_{\rm ax}},
\ \ |\p_R\mathscr R_{r,t}(R,Z)|+|\p_R\mathscr R_{z,t}(R,Z)| \le C_{\rm ax}J^{-1}|R|^{\beta_{\rm ax}} .
\label{eq:axis-normal-form-entry}
\end{equation}
Moreover, for $\zeta\in I_\sharp$ and $|\tau|\le C_0$, the fixed-slope normalized remainder obeys
\begin{equation}
\left| \p_\zeta \left((J^2\zeta)^{-1} \mathscr R_{r,t}(R_t^\sharp(\zeta,\tau),Z_t(\zeta))\right) \right|
+ \left| \p_\zeta \left((J^2\zeta)^{-1} \mathscr R_{z,t}(R_t^\sharp(\zeta,\tau),Z_t(\zeta))\right) \right|
\le C_{\rm ax}J^{3\beta_{\rm ax}}.
\label{eq:axis-normal-form-zeta-derivative}
\end{equation}
\end{lemma}

\begin{remark}[Radial-label comparison]
Lemma~\ref{lem:nonlinear-radial-flatness} compares the radial derivative at $(R,Z)$ with the radial
derivative on the axis at the same axial label $Z$.  This is a Lagrangian comparison in the radial label, not
an Eulerian Taylor expansion at the current point.  The normalization by $A_t(Z)$ removes the coherent
hyperbolic radial stretching, so the estimate only has to control the variation of the normalized radial
derivative away from the axis.
\end{remark}

We prove Lemma~\ref{lem:nonlinear-radial-flatness} after first deriving the source estimate for this
normalized radial-derivative equation.

\begin{lemma}[Source bound for normalized radial derivatives]
\label{lem:radial-flatness-source-bound}
Assume the bootstrap bounds \eqref{eq:localized-clock-bootstrap},
\eqref{eq:localized-normal-form-large-bootstrap}, \eqref{eq:localized-cusp-error-large-bootstrap}, and
\eqref{eq:buffered-axis-geometry-bootstrap}, together with the containment condition
\eqref{eq:radial-flatness-buffered-label}.  Define the normalized radial derivatives
\[
G_r(R,Z,t)=A_t(Z)^{-1}\p_Rr_t(R,Z), \qquad G_z(R,Z,t)=A_t(Z)^{-1}\p_Rz_t(R,Z),
\]
set $(r,z)=\phi_{\cusp}(R,Z,t)$, and define
\[
S_t(Z):=\p_r(V_{\cusp})_r(0,B_t(Z),t), \qquad W_t(Z):=\p_z(V_{\cusp})_z(0,B_t(Z),t).
\]
Then the normalized radial-derivative defect $F:=(G_r-1,G_z)^T$ obeys
\begin{equation}
\p_tF =
\begin{pmatrix}0&0\\0&W_t(Z)-S_t(Z)\end{pmatrix} F + \mathcal E(R,Z,t)F + \mathcal S(R,Z,t),
\label{eq:radial-flatness-normalized-system}
\end{equation}
where
\begin{equation}
\begin{aligned}
\mathcal E(R,Z,t)
&:=
\begin{pmatrix}
\p_r(V_{\cusp})_r(r,z,t)-S_t(Z) & \p_z(V_{\cusp})_r(r,z,t)
\\
\p_r(V_{\cusp})_z(r,z,t) & \p_z(V_{\cusp})_z(r,z,t)-W_t(Z)
\end{pmatrix},\\
\mathcal S(R,Z,t) &:= \binom{ \p_r(V_{\cusp})_r(r,z,t)-S_t(Z)} {\p_r(V_{\cusp})_z(r,z,t)} .
\end{aligned}
\label{eq:radial-flatness-source-matrix-def}
\end{equation}
The diagonal part in \eqref{eq:radial-flatness-normalized-system} has the exact propagator
$\bigl(B_t'(Z)/A_t(Z)\bigr)/\bigl(B_s'(Z)/A_s(Z)\bigr)$.  Set
\[
\delta_{\rm rad}:=\min\{3\beta_{\rm ax},\,3\alpha,\,1-3\alpha\}>0 .
\]
Then, on the radial tube
\[
Z\in I_Z^t,\qquad |R|\le C_{\rm ax}J_{\cusp}(t)^3,
\]
\begin{equation}
\begin{aligned}
& |\mathcal S_r(R,Z,t)|+|\mathcal S_z(R,Z,t)|
+
\|\mathcal E(R,Z,t)\|\,|R|^{\beta_{\rm ax}} \\
& \qquad \qquad
\le
C_{B_*,E_*}\,\Gamma J_{\cusp}(t)^{3\alpha-1+\delta_{\rm rad}}|R|^{\beta_{\rm ax}}
+
C_{E_*}\,\Gamma\bigl(J_{\cusp}(t)^{9\alpha-1}+1\bigr)|R|^{\beta_{\rm ax}} .
\end{aligned}
\label{eq:radial-flatness-source}
\end{equation}
The estimate \eqref{eq:radial-flatness-source} also holds for the finite-difference quotient in the time-$t$ $\zeta$ variable
$\zeta=J_{\cusp}(t)^{-2}B_t(Z)$, with the two labels related by
$A_t(Z)R/B_t(Z)=\tau$ fixed.
\end{lemma}

\begin{proof}[Proof of Lemma~\ref{lem:radial-flatness-source-bound}]
\runinhead{Step 1: The normalized radial-derivative system.}
Differentiating the cusp-flow equation \eqref{eq:cusp-flow-equation} in the radial label, we obtain the
evolution equations for the two components of $\p_R\phi_{\cusp}$:
\begin{subequations}
\begin{align}
\p_t\p_R r_t(R,Z)
&=
\p_r(V_{\cusp})_r(r,z,t)\,\p_Rr_t(R,Z)
+
\p_z(V_{\cusp})_r(r,z,t)\,\p_Rz_t(R,Z),
\label{eq:radial-flatness-raw-r-derivative}\\
\p_t\p_R z_t(R,Z)
&=
\p_r(V_{\cusp})_z(r,z,t)\,\p_Rr_t(R,Z)
+
\p_z(V_{\cusp})_z(r,z,t)\,\p_Rz_t(R,Z),
\label{eq:radial-flatness-raw-z-derivative}
\end{align}
\end{subequations}
where $(r,z)=\phi_{\cusp}(R,Z,t)$.  Dividing
\eqref{eq:radial-flatness-raw-r-derivative}--\eqref{eq:radial-flatness-raw-z-derivative} by the
axis function $A_t(Z)$ yields the normalized system
\eqref{eq:radial-flatness-normalized-system}--\eqref{eq:radial-flatness-source-matrix-def}.
Before the label-radial derivative $\p_R\phi_{\cusp}$ is divided by $A_t(Z)$, both raw equations
\eqref{eq:radial-flatness-raw-r-derivative}--\eqref{eq:radial-flatness-raw-z-derivative} contain the
same singular linear radial stretching.  The equation
\[
\p_tA_t(Z)=\p_r(V_{\cusp})_r(0,B_t(Z),t)A_t(Z)
\]
is exactly this common part, so the normalization by $A_t(Z)$ removes it.  The only singular homogeneous term
left in the $G_z$ equation is the difference between the axial and radial axis rates; by
axisymmetric incompressibility and the identities $\p_tB_t'=W_tB_t'$ and $\p_tA_t=S_tA_t$, its
propagator is
\[
\exp\!\left(\int_s^t(W_\ell(Z)-S_\ell(Z))\,d\ell\right) = \tfrac{B_t'(Z)/A_t(Z)}{B_s'(Z)/A_s(Z)} .
\]
The Duhamel estimate below keeps this exact ratio as the homogeneous propagator.

\runinhead{Step 2: The singular cusp contribution.}
The source $\mathcal S(R,Z,t)$ in \eqref{eq:radial-flatness-normalized-system} consists of differences between
the velocity-gradient matrix at $\phi_{\cusp}(R,Z,t)$ and its axis value at
$\phi_{\cusp}(0,Z,t)=(0,B_t(Z))$.
We compare the velocity gradient at a fixed axial label $Z$ with its axis value by moving along the
radial label segment
\[
Y_\lambda=(\lambda R,Z),\qquad 0\le \lambda\le1 .
\]
If $Z\in I_Z^t$, then the containment assumption \eqref{eq:radial-flatness-buffered-label} shows that
\[
J_s^{-2}B_s(Z)\in I_{\rm buf},\qquad t_0\le s\le t .
\]
Thus the axis estimates \eqref{eq:buffered-axis-geometry-bootstrap} and the normal-form bootstrap
\eqref{eq:localized-normal-form-large-bootstrap} apply along this whole radial segment whenever the source term
is evaluated.  At the final time, $J^{-2}B_t(Z)\in I_\sharp\Subset(0,\infty)$, so the same axis estimates imply
that $Z$ is bounded above and below by fixed positive constants on $I_Z^t$.  Therefore, on the tube
\[
Z\in I_Z^t,\qquad |R|\le C_{\rm ax}J_{\cusp}(t)^3,
\]
we have that
\[
\left|\tfrac{\lambda R}{Z}\right|\le C J_{\cusp}(t)^3,\qquad 0\le\lambda\le 1 .
\]
After decreasing the small-clock threshold, $\sigma(Y_\lambda)\le\sigma_{\cut}$.  Thus the angular function $\Theta^*$ in \eqref{eq:Theta-star-def} is evaluated 
in the near-axis region where $\Upsilon\equiv1$ and $\Theta^*(\sigma)=(\sin\sigma)^\alpha$. 

For $Y_\lambda=(\lambda R,Z)$, $0\le\lambda\le1$, we set
$\zeta_s(Z):=J_s^{-2}B_s(Z)$ and
$\tau_{\lambda,s}:=\tfrac{A_s(Z)\lambda R}{B_s(Z)}$. 
The containment assumption \eqref{eq:radial-flatness-buffered-label} implies
$\zeta_s(Z)\in I_{\rm buf}$ for $t_0\le s\le t$; moreover, by
\eqref{eq:buffered-axis-geometry-bootstrap}, the compact inclusion
$I_{\rm buf}\Subset(0,\infty)$, and the tube condition $|R|\le C_{\rm ax}J_{\cusp}(t)^3$, we have that
\[
|\tau_{\lambda,s}|\le C,\qquad 0\le\lambda\le1,\quad t_0\le s\le t .
\]
Thus, after writing the Eulerian image in the variables obtained by division by $J_s^2$, the relevant
Biot--Savart kernels are evaluated on the time-independent set
\[
\{(\zeta\tau,\zeta):\zeta\in I_{\rm buf},\ |\tau|\le C\},
\]
and the kernel constants in this part of the estimate are independent of $J_s$.

The three exponents in
$\delta_{\rm rad}:=\min\{3\beta_{\rm ax},\,3\alpha,\,1-3\alpha\}$ have the following origins.  The exponent
$3\beta_{\rm ax}$ comes from the normal-form bound \eqref{eq:localized-normal-form-large-bootstrap}.  The
exponent $3\alpha$ comes from the near-axis angular structure: since
$|\lambda R/Z|=O(J_{\cusp}(t)^3)$ and $\sigma(Y_\lambda)\le\sigma_{\cut}$, the definition
\eqref{eq:Theta-star-def} implies
$\rho(Y_\lambda)^\alpha\Theta^*(\sigma(Y_\lambda))=R(Y_\lambda)^\alpha$, so the angular term contributes the
multiplier $J_{\cusp}(t)^{3\alpha}$.  The exponent $1-3\alpha$ is the ratio of the nonsingular
$O(\Gamma)$ contribution to the singular clock scale $\Gamma J_{\cusp}(t)^{3\alpha-1}$.  These three bounds
produce the first term on the right-hand side of \eqref{eq:radial-flatness-source}.

\runinhead{Step 3: The error-velocity contribution.}
The smooth-flow deformation and the algebraic tail enter through $V_{\err}$.  At this stage of the bootstrap
argument we use the large cusp-error bound \eqref{eq:localized-cusp-error-large-bootstrap}; after
Lemma~\ref{lem:tail-bound} is proved, the axis-error traces associated with $V_{\err}$ are sharpened in
\eqref{eq:axis-profile-error-holder}--\eqref{eq:axis-profile-error-holder-samelabel}.  The large bootstrap
bound contributes the term
$C_{E_*}\Gamma(J_{\cusp}^{9\alpha-1}+1)|R|^{\beta_{\rm ax}}$ in \eqref{eq:radial-flatness-source}.  Applying
these estimates to the vorticity transport formula \eqref{eq:Omega-cusp-def} and to the normalized
radial-derivative system \eqref{eq:radial-flatness-normalized-system} proves the $V_{\err}$ contribution in
\eqref{eq:radial-flatness-source}.  Together with the singular sampling bound, this proves
\eqref{eq:radial-flatness-source}.

\runinhead{Step 4: Fixed-slope $\zeta$ differences.}
We next take finite differences in $\zeta$ while keeping $\tau$ fixed.  In that quotient the
axis functions $A_t$ and $B_t$ are compared through the identities
$B_t(Z_t(\zeta))=J_{\cusp}(t)^2\zeta$ and
$A_t(Z_t(\zeta))R=J_{\cusp}(t)^2\zeta\tau$.  The finite difference is taken after imposing these identities, so
the leading variations of $A_t$ and $B_t$ cancel before the inhomogeneous terms are estimated.  The interval
$I_\sharp\Subset(0,\infty)$ is compactly separated from the endpoint $\zeta=0$, and differentiating the smooth
radial weights, the cutoff multipliers, and the algebraic tail in this fixed-$\tau$ $\zeta$ direction preserves
the right-hand side of \eqref{eq:radial-flatness-source}.  The angular function $\Theta^*$ is still evaluated at
$R/Z=O(J_{\cusp}^3)$, so the same gain $|R|^{\beta_{\rm ax}}$ is retained.
\end{proof}

\begin{proof}[Proof of Lemma~\ref{lem:nonlinear-radial-flatness}]
\runinhead{Step 1: The normalized equation and its propagator.}
We write
\[
G_r(R,Z,t):=A_t(Z)^{-1}\p_Rr_t(R,Z), \qquad G_z(R,Z,t):=A_t(Z)^{-1}\p_Rz_t(R,Z).
\]
Thus $G_r(0,Z,t)=1$ and $G_z(0,Z,t)=0$.  Differentiating the cusp-flow equation
\[
\p_t\phi_{\cusp}(Y,t)=V_{\cusp}(\phi_{\cusp}(Y,t),t)
\]
with respect to the radial label and subtracting the axis variational equation for $A_t(Z)$ produces a closed
system for $(G_r-1,G_z)$.  All functions in this system are evaluated at the fixed axial label $Z$.  Let
\[
S_t(Z):=\p_r(V_{\cusp})_r(0,B_t(Z),t), \qquad W_t(Z):=\p_z(V_{\cusp})_z(0,B_t(Z),t).
\]
By the axisymmetric divergence-free identity, $S_t(Z)=-\tfrac12W_t(Z)$.  The radial derivative
$A_t(Z)$ solves $\p_tA_t=S_tA_t$, while the axial derivative $B_t'(Z)$ solves
$\p_tB_t'=W_tB_t'$.  Thus the singular diagonal term which remains in the $G_z$ equation has
the exact propagator
\begin{equation}
\exp\!\left(\int_s^t(W_\ell(Z)-S_\ell(Z))\,d\ell\right) = \tfrac{B_t'(Z)/A_t(Z)}{B_s'(Z)/A_s(Z)} .
\label{eq:radial-flatness-exact-propagator}
\end{equation}
By the containment condition \eqref{eq:radial-flatness-buffered-label}, the axial flow map
geometry implies $A_\ell(Z)\simeq J_{\cusp}(\ell)^{-1}$ and
$B_\ell'(Z)\simeq J_{\cusp}(\ell)^2$ at every intermediate time $\ell\in[s,t]$ along the fixed axial label.
Together with the two-sided clock comparison, this bounds the ratio in
\eqref{eq:radial-flatness-exact-propagator} uniformly on the small-clock interval.  Thus the proof does not
estimate the singular homogeneous part perturbatively; it keeps that part in the exact ratio
\eqref{eq:radial-flatness-exact-propagator}.

\runinhead{Step 2: Duhamel estimate for the normalized radial derivatives.}
By Lemma~\ref{lem:radial-flatness-source-bound}, the inhomogeneous term $\mathcal S$ and the linear matrix
$\mathcal E$ obey \eqref{eq:radial-flatness-source}.  The Duhamel estimate below uses the nonlinear radial
structure only through this bound.

We apply Duhamel's formula to the normalized radial-derivative system with the exact propagator
\eqref{eq:radial-flatness-exact-propagator}.  The contribution from the initial time $t_0$ satisfies
\[
|G_r(R,Z,t_0)-1|+|G_z(R,Z,t_0)|\le C|R|^{\beta_{\rm ax}}
\]
by the ordinary $C^{1,\beta_{\rm ax}}$ regularity of the cusp map before the small-clock regime.  Dividing the
equation by
$|R|^{\beta_{\rm ax}}$ and using \eqref{eq:radial-flatness-source}, the matrix $\mathcal E$ and
the source vector $\mathcal S$ are both controlled by the integrable clock weight
\[
K(J):= C_{B_*,E_*}\Gamma J^{3\alpha-1+\delta_{\rm rad}} + C_{E_*}\Gamma(J^{9\alpha-1}+1).
\]
Thus Gronwall's inequality yields a uniform bound for
$|R|^{-\beta_{\rm ax}}\bigl(|G_r-1|+|G_z|\bigr)$ once
$\int K(J_{\cusp}(t))\,dt$ is bounded.  On the small-clock interval the cusp-clock rate bound
\eqref{eq:localized-clock-bootstrap}, used here as one of the small-clock bootstrap assumptions and
closed later in Lemma~\ref{lem:Jdot-two-sided-aux}, implies
\[
dt\le C\,\tfrac{-dJ}{\Gamma J^{3\alpha}}.
\]
With $\delta_{\rm rad}$ as in Lemma~\ref{lem:radial-flatness-source-bound}, the source estimate yields
\[
\int \Gamma J^{3\alpha-1+\delta_{\rm rad}}\,dt \le C\int_0^{\mathfrak J_{\mathrm{axis}}}J^{\delta_{\rm rad}-1}\,dJ<\infty
\]
and
\[
\int \Gamma\bigl(J^{9\alpha-1}+1\bigr)\,dt \le C\int_0^{\mathfrak J_{\mathrm{axis}}} \bigl(J^{6\alpha-1}+J^{-3\alpha}\bigr)\,dJ<\infty ,
\]
because $0<\alpha<\tfrac13$.  Hence
\begin{equation}
|G_r(R,Z,t)-1|+|G_z(R,Z,t)| \le C|R|^{\beta_{\rm ax}} .
\label{eq:radial-derivative-flat}
\end{equation}

\runinhead{Step 3: From derivative bounds to hyperbolic-deviation bounds.}
Since the axial flow map geometry implies $A_t(Z)\simeq J^{-1}$ on $I_Z^t$,
\eqref{eq:radial-derivative-flat}
implies
\[
|\p_Rr_t(R,Z)-A_t(Z)|+|\p_Rz_t(R,Z)| \le C J^{-1}|R|^{\beta_{\rm ax}} .
\]
Integrating this bound from $0$ to $R$ and using $r_t(0,Z)=0$, $z_t(0,Z)=B_t(Z)$ yields the two remainder
estimates in \eqref{eq:axis-normal-form-entry}.

\runinhead{Step 4: $\zeta$ derivatives at fixed slope.}
It remains to prove the $\zeta$-derivative form \eqref{eq:axis-normal-form-zeta-derivative}, which is needed
to control the map $\Psi_t$ in \eqref{eq:normal-form-approximation-def} from the model point
$(R_{\rm sc},Z_{\rm sc})$ to the Eulerian image divided by $J^2$.  We apply the same normalized
radial-derivative argument to finite differences in the time-$t$ $\zeta$ variable, while keeping the axial
flow-map slope $\tau$ fixed.  Thus the two compared labels are
\[
(R_t^\sharp(\zeta+h,\tau),Z_t(\zeta+h)), \qquad (R_t^\sharp(\zeta,\tau),Z_t(\zeta)).
\]
The finite difference is taken after imposing the axial flow map identities
\[
B_t(Z_t(\zeta))=J^2\zeta,\qquad A_t(Z_t(\zeta))R_t^\sharp(\zeta,\tau)=J^2\zeta\tau .
\]
The leading variations of $A_t$ and $B_t$ cancel at this stage.  The remaining inhomogeneous term is the
same radial sampling error as above, now divided by $|\zeta+h-\zeta|$.  Since
$I_\sharp\Subset(0,\infty)$, the only nonsmooth dependence is still through the angular ratio
$R/Z=O(J^3)$, and the source bound
\eqref{eq:radial-flatness-source} is stable under this fixed-$\tau$ $\zeta$ difference quotient.  Duhamel's
formula with the exact propagator \eqref{eq:radial-flatness-exact-propagator} yields, with
$\Delta_h f(\zeta):=(f(\zeta+h)-f(\zeta))/h$,
\[
\left| \Delta_h \left((J^2\zeta)^{-1} \mathscr R_{\star,t}(R_t^\sharp(\zeta,\tau),Z_t(\zeta))\right) \right| \le C J^{3\beta_{\rm ax}},
\qquad \star\in\{r,z\},
\]
uniformly for $|h|$ sufficiently small.  Passing to the limit $h\to0$ proves
\eqref{eq:axis-normal-form-zeta-derivative}.
\end{proof}

The radial-derivative estimate \eqref{eq:axis-normal-form-entry} is written in the original label variables
$(R,Z)$.  To state the normal form in the variables used by the pressure Hessian comparison, we now rescale along
the axial image of the cusp map at the same time $t$.  Along the symmetry axis,
$\phi_{\cusp}(0,Z,t)=(0,B_t(Z))$, where $B_t(Z):=z_t(0,Z)$.  The axial collapse scale is $J^2$, with
$J:=J_{\cusp}(t)$, so the normalized axial coordinate is $\zeta=J^{-2}B_t(Z)$ and
$B_t(Z)=J^2\zeta$.  The pressure Hessian comparison is localized on strictly positive axial scales, so the
relevant $\zeta$-interval is chosen inside $(0,\infty)$.  The variables normalized by the axial collapse
scale $J^2$ are the pair $(\zeta,\tau)$ in \eqref{eq:adapted-labels}; here
$\tau=A_t(Z)R/B_t(Z)$ measures radial distance relative to the axial scale in the time-$t$ axis chart.  We fix
$I_\sharp\Subset I_{\rm buf}\Subset(0,\infty)$, fix $C_0\ge1$, and recall that
\[
\phi_{\cusp}(R,Z,t)=\bigl(r_t(R,Z),z_t(R,Z)\bigr), \qquad A_t(Z):=\p_Rr_t(0,Z), \qquad B_t(Z):=z_t(0,Z).
\]
We let $I_Z^t$ be the axial-label interval mapped onto $I_\sharp$ by
$Z\mapsto J^{-2}B_t(Z)$, and we let
$Z_t:I_\sharp\to I_Z^t$ denote the inverse map.  The axis geometry used below is
\begin{equation}
\begin{aligned}
c_{\rm ax}\le J A_t(Z_t(\zeta))\le C_{\rm ax},\qquad c_{\rm ax}\le J^{-2}B_t'(Z_t(\zeta))\le C_{\rm ax},\\
[\log(J A_t(Z_t(\cdot)))]_{C^{\alpha/2}(I_\sharp)} + [\log(J^{-2}B_t'(Z_t(\cdot)))]_{C^{\alpha/2}(I_\sharp)} \le C_{\rm ax}.
\end{aligned}
\label{eq:entry-axis-bounds-statement}
\end{equation}
on the fixed pressure interval $I_\sharp$.  For a label $(R,Z)$ near the positive axis, we use the coordinates
\begin{equation}
\zeta:=J^{-2}B_t(Z), \qquad  \tau:=\tfrac{A_t(Z)R}{B_t(Z)}.   \label{eq:adapted-labels}
\end{equation}
Conversely, for $\zeta\in I_\sharp$ and $|\tau|\le C_0$, we define
\[
R_t(\zeta,\tau):=\tfrac{J^2\zeta\tau}{A_t(Z_t(\zeta))}, \qquad Y_t(\zeta,\tau):=(R_t(\zeta,\tau),Z_t(\zeta)).
\]
Then
\[
B_t(Z_t(\zeta))=J^2\zeta, \qquad A_t(Z_t(\zeta))R_t(\zeta,\tau)=J^2\zeta\tau .
\]
Thus $Y_t(\zeta,\tau)$ is the label point with adapted coordinates $(\zeta,\tau)$.  The localized label tube used in
the next two lemmas is
\begin{equation}
\mathcal Q_t^\sharp := \bigl\{Y_t(\zeta,\tau):\ \zeta\in I_\sharp,\ |\tau|\le C_0\bigr\}.
\label{eq:Qt-sharp-def}
\end{equation}
After division of the Eulerian image by $J^2$, the corresponding model point is
\[
(R_{\rm sc},Z_{\rm sc})=(\zeta\tau,\zeta).
\]
The physical cylindrical label radius satisfies $R\ge0$.  In estimates that use parity at the symmetry axis,
we fix a meridional plane and use a Cartesian coordinate across the axis, oriented so that positive values agree
with the cylindrical radius; within those parity estimates we still denote this signed coordinate by $R$.  The
point with signed coordinate $-R$ represents the same cylindrical radius as $R$ and has azimuth shifted by
$\pi$.  Axisymmetry implies the signed extensions $r_t(-R,Z)=-r_t(R,Z)$ and $z_t(-R,Z)=z_t(R,Z)$ in this
coordinate, so the slope $\tau=A_t(Z)R/B_t(Z)$ may be signed in the map estimates leading to
\eqref{eq:additive-normal-form}.  When fractional powers such as $\tau^\alpha$ occur in the vorticity
calculation, $R$ again denotes the nonnegative physical cylindrical radius, so $R\ge0$ and $\tau\ge0$; the
signed negative side is only a parity device for the map estimates.

\begin{lemma}[Normal form for the cusp map $\phi_{\cusp}$]
\label{lem:late-axis-normal-form-cusp}
Assume the axis-geometry bounds \eqref{eq:entry-axis-bounds-statement}, the containment condition
\eqref{eq:radial-flatness-buffered-label}, and the bootstrap bounds
\eqref{eq:localized-clock-bootstrap}, \eqref{eq:localized-normal-form-large-bootstrap}, and
\eqref{eq:localized-cusp-error-large-bootstrap}.  After decreasing the small-clock threshold depending
only on the fixed parameters and the large bootstrap constants, there is an error function
$\mathcal E_t:I_\sharp\times[-C_0,C_0]\to\R^2$ such that, for
$Y_t(\zeta,\tau)\in\mathcal Q_t^\sharp$,
\begin{equation}
\phi_{\cusp}(Y_t(\zeta,\tau),t) = J^2\zeta\Bigl((\tau,1)+\mathcal E_t(\zeta,\tau)\Bigr),
\label{eq:additive-normal-form}
\end{equation}
where $\mathcal E_t=(\mathcal E_{t,r},\mathcal E_{t,z})$,
$\mathcal E_{t,r}(\zeta,0)=0$, and
\begin{equation}
\|\mathcal E_t\|_{L^\infty} + \|\p_\tau\mathcal E_t\|_{L^\infty} + \|\p_\zeta\mathcal E_t\|_{L^\infty}
+ [\mathcal E_t]_{C^{\beta_{\rm ax}}_{\zeta,\tau}(I_\sharp\times[-C_0,C_0])} \le C J^{3\beta_{\rm ax}}.
\label{eq:additive-normal-form-bound}
\end{equation}
\end{lemma}

Before proving Lemma~\ref{lem:late-axis-normal-form-cusp}, we introduce the variables used in
Lemma~\ref{lem:late-axis-normal-form-map-cusp} to compare the exact Eulerian image with the model image after
both are divided by $J^2$.  In the slope-restricted pressure Hessian integral, the label
$Y_t(\zeta,\tau)$ is associated with the model point
$(R_{\rm sc},Z_{\rm sc})=(\zeta\tau,\zeta)$, which ranges over
\begin{equation}
\mathcal R_{\sharp,C_0}^{\rm sc} := \{(R_{\rm sc},Z_{\rm sc})=(\zeta\tau,\zeta):\ \zeta\in I_\sharp,\ |\tau|\le C_0\}.
\label{eq:normal-form-fixed-scaled-set}
\end{equation}
The additive normal form \eqref{eq:additive-normal-form}--\eqref{eq:additive-normal-form-bound} identifies the
Eulerian image of the same label after division by $J^2$:
\begin{equation}
\Psi_t(R_{\rm sc},Z_{\rm sc}) := J^{-2}\phi_{\cusp}(Y_t(\zeta,\tau),t) =
\zeta\bigl((\tau,1)+\mathcal E_t(\zeta,\tau)\bigr).
\label{eq:normal-form-approximation-def}
\end{equation}
Thus $\Psi_t$ sends the model point $(R_{\rm sc},Z_{\rm sc})$ to the actual point $J^{-2}(r,z)$ in the
Eulerian image divided by $J^2$.  Its displacement from the identity is exactly
$\zeta\mathcal E_t(\zeta,\tau)$.  Lemma~\ref{lem:late-axis-normal-form-map-cusp} proves the bi-Lipschitz and
displacement bounds \eqref{eq:normal-form-approximation-bound-a}--\eqref{eq:normal-form-approximation-bound-b}
and the cylindrical-volume identity \eqref{eq:normal-form-approximation-bound-c}.  Hence replacing the model
point $(R_{\rm sc},Z_{\rm sc})$ by the actual point $\Psi_t(R_{\rm sc},Z_{\rm sc})$ changes the evaluation
point by a controlled amount and preserves the measure $R_{\rm sc}\,dR_{\rm sc}\,dZ_{\rm sc}$.

\begin{lemma}[Geometry of the image map after division by $J^2$]
\label{lem:late-axis-normal-form-map-cusp}
Assume the axis-geometry bounds \eqref{eq:entry-axis-bounds-statement}, the containment condition
\eqref{eq:radial-flatness-buffered-label}, and the bootstrap bounds
\eqref{eq:localized-clock-bootstrap}, \eqref{eq:localized-normal-form-large-bootstrap}, and
\eqref{eq:localized-cusp-error-large-bootstrap}.  After decreasing the small-clock threshold depending
only on the fixed parameters and the large bootstrap constants, the map $\Psi_t$ in
\eqref{eq:normal-form-approximation-def} is a $C^1$ diffeomorphism of
$\mathcal R_{\sharp,C_0}^{\rm sc}$ onto its image, and
\begin{subequations}
\label{eq:normal-form-approximation-bound}
\begin{align}
\|D\Psi_t\|_{L^\infty}+\|D\Psi_t^{-1}\|_{L^\infty} &\le C,
\label{eq:normal-form-approximation-bound-a}
\\
\|\Psi_t-\operatorname{Id}\|_{L^\infty}+\|D\Psi_t-I\|_{L^\infty}
+ [\Psi_t-\operatorname{Id}]_{C^{\beta_{\rm ax}}} &\le C J^{3\beta_{\rm ax}},
\label{eq:normal-form-approximation-bound-b}
\\
\mathcal J_{\Psi_t}:=\tfrac{(\Psi_t)_R}{R_{\rm sc}}\det D_{R_{\rm sc},Z_{\rm sc}}\Psi_t &=1 .
\label{eq:normal-form-approximation-bound-c}
\end{align}
\end{subequations}
\end{lemma}

\begin{proof}[Proof of Lemma~\ref{lem:late-axis-normal-form-cusp}]
\eqref{eq:entry-axis-bounds-statement} is the axis-geometry hypothesis.  The radial
derivative estimates are supplied by Lemma~\ref{lem:nonlinear-radial-flatness}.  The
oddness of $r_t$ in $R$ and the evenness of $z_t$ in $R$ give
\[
r_t(R,Z)=A_t(Z)R+\mathscr R_{r,t}(R,Z), \qquad z_t(R,Z)=B_t(Z)+\mathscr R_{z,t}(R,Z).
\]
The axial flow map identities are
\[
B_t(Z_t(\zeta))=J^2\zeta, \qquad A_t(Z_t(\zeta))R_t(\zeta,\tau)=J^2\zeta\tau .
\]
Hence
\[
R_t(\zeta,\tau)=\tfrac{J^2\zeta\tau}{A_t(Z_t(\zeta))} = \tfrac{J^3\zeta\tau}{J A_t(Z_t(\zeta))},
\]
and \eqref{eq:entry-axis-bounds-statement} implies
\[
|R_t(\zeta,\tau)|\le C J^3, \qquad |\p_\tau R_t(\zeta,\tau)|\le C J^3 \qquad(\zeta\in I_\sharp,\ |\tau|\le C_0).
\]
We define the error components by
\[
\mathcal E_{t,r}(\zeta,\tau) := (J^2\zeta)^{-1}\mathscr R_{r,t}(R_t(\zeta,\tau),Z_t(\zeta)),
\qquad \mathcal E_{t,z}(\zeta,\tau) := (J^2\zeta)^{-1}\mathscr R_{z,t}(R_t(\zeta,\tau),Z_t(\zeta)).
\]
Then \eqref{eq:additive-normal-form} follows from
\[
r_t(R,Z)=A_t(Z)R+\mathscr R_{r,t}(R,Z), \qquad z_t(R,Z)=B_t(Z)+\mathscr R_{z,t}(R,Z).
\]
The estimates \eqref{eq:axis-normal-form-entry} and the bounds for $R_t,\p_\tau R_t$ imply
\[
\|\mathcal E_t\|_{L^\infty}+\|\p_\tau\mathcal E_t\|_{L^\infty}\le C J^{3\beta_{\rm ax}}.
\]
The $\zeta$-derivative estimate \eqref{eq:axis-normal-form-zeta-derivative} yields
\[
\|\p_\zeta\mathcal E_t\|_{L^\infty}\le C J^{3\beta_{\rm ax}}.
\]
The $C^{\beta_{\rm ax}}_{\zeta,\tau}$ bound follows from the same product estimate applied to
\eqref{eq:entry-axis-bounds-statement} and \eqref{eq:axis-normal-form-entry}.  Since
$R_t(\zeta,0)=0$ and $\mathscr R_{r,t}(0,Z_t(\zeta))=0$, we also have
$\mathcal E_{t,r}(\zeta,0)=0$.  This proves \eqref{eq:additive-normal-form-bound} and the stated
axis value of $\mathcal E_{t,r}$.
\end{proof}

\begin{proof}[Proof of Lemma~\ref{lem:late-axis-normal-form-map-cusp}]
\runinhead{Step 1: Bounds in the image variables after division by $J^2$.}
We write $X=(R_{\rm sc},Z_{\rm sc})$ and recover
\[
\zeta=Z_{\rm sc},\qquad \tau=\tfrac{R_{\rm sc}}{Z_{\rm sc}}.
\]
Since $I_\sharp\Subset(0,\infty)$ and $|\tau|\le C_0$, the change of variables
$X\leftrightarrow(\zeta,\tau)$ has uniformly bounded derivatives on the set
$\mathcal R_{\sharp,C_0}^{\rm sc}$ in \eqref{eq:normal-form-fixed-scaled-set}.  The regularity used here
is the $C^1$ regularity of the Euler flow map together with the quantitative H\"older bounds in the
axial flow map chart.  The map $Z\mapsto J^{-2}B_t(Z)$ is the monotone axial coordinate on the label interval
under consideration, and $Z_t$ denotes its inverse.  The lower bound
$J^{-2}B_t'(Z_t(\zeta))\ge c_{\rm ax}$ in \eqref{eq:entry-axis-bounds-statement} is the nondegeneracy
condition in the one-dimensional inverse theorem; it is not the source of the H\"older exponent.  The inverse is
used only through the differentiated identity
\[
J^{-2}B_t(Z_t(\zeta))=\zeta, \qquad \p_\zeta Z_t(\zeta) = \bigl(J^{-2}B_t'(Z_t(\zeta))\bigr)^{-1}.
\]
By \eqref{eq:entry-axis-bounds-statement}, the function
$J^{-2}B_t'(Z_t(\cdot))$ is bounded above and below and has logarithm bounded in
$C^{\alpha/2}(I_\sharp)$.  Thus $\p_\zeta Z_t$ has the same $C^{\alpha/2}$ control.  The exponent used
below is the deliberately smaller exponent $\beta_{\rm ax}=\tfrac\alpha4$ from
\eqref{eq:beta-ax-def}; hence products and compositions with the bounded coordinate functions $\zeta$, $\tau$, and
$1/\zeta$ preserve the required $C^{\beta_{\rm ax}}$ bounds on this detached set.  The component estimates for
the normal-form error are precisely \eqref{eq:axis-normal-form-entry} and
\eqref{eq:axis-normal-form-zeta-derivative}, summarized in \eqref{eq:additive-normal-form-bound}.
These estimates give the $C^1$ control of $\Psi_t$ and the explicit $C^{\beta_{\rm ax}}$ displacement bound
proved below.
By \eqref{eq:normal-form-approximation-def},
\[
\Psi_t(R_{\rm sc},Z_{\rm sc})-(R_{\rm sc},Z_{\rm sc}) = \zeta\,\mathcal E_t(\zeta,\tau), \qquad (R_{\rm sc},Z_{\rm sc})=(\zeta\tau,\zeta).
\]
From the product and change-of-variable estimates, together with \eqref{eq:additive-normal-form-bound}, we obtain
\[
\|\Psi_t-\operatorname{Id}\|_{L^\infty}+ [\Psi_t-\operatorname{Id}]_{C^{\beta_{\rm ax}}} \le C J^{3\beta_{\rm ax}} .
\]
For the derivative bound, we use that
\[
\p_{R_{\rm sc}}=\zeta^{-1}\p_\tau, \qquad \p_{Z_{\rm sc}}=\p_\zeta-\tfrac{\tau}{\zeta}\p_\tau ,
\]
where the terms $1/\zeta$ are uniformly bounded on $I_\sharp$.  By differentiating
\eqref{eq:normal-form-approximation-def}, we obtain that
\begin{align*}
\p_{R_{\rm sc}}\Psi_t  &= \bigl(1+\p_\tau\mathcal E_{t,r},\,  \p_\tau\mathcal E_{t,z}\bigr),\\
\p_{Z_{\rm sc}}\Psi_t &=  \bigl(\mathcal E_{t,r} +\zeta\p_\zeta\mathcal E_{t,r} -\tau\p_\tau\mathcal E_{t,r},\,
1+\mathcal E_{t,z}  +\zeta\p_\zeta\mathcal E_{t,z}  -\tau\p_\tau\mathcal E_{t,z}\bigr).
\end{align*}
Using \eqref{eq:additive-normal-form-bound} once more,
\[
\|D\Psi_t-I\|_{L^\infty} \le C J^{3\beta_{\rm ax}} .
\]

\runinhead{Step 2: The bi-Lipschitz estimate.}
The image-variable set $\mathcal R_{\sharp,C_0}^{\rm sc}$, defined in
\eqref{eq:normal-form-fixed-scaled-set}, is convex.
For two points $X_0,X_1$ in this set, the fundamental theorem of calculus along the segment from $X_0$ to
$X_1$ yields
\[
\Psi_t(X_1)-\Psi_t(X_0) =(X_1-X_0)+\int_0^1\bigl(D\Psi_t(X_0+s(X_1-X_0))-I\bigr)(X_1-X_0)\,ds.
\]
After decreasing the small-clock threshold so that $C J^{3\beta_{\rm ax}}\le\tfrac12$, we obtain
\[
\tfrac12|X_1-X_0| \le |\Psi_t(X_1)-\Psi_t(X_0)| \le \tfrac32|X_1-X_0| .
\]
Thus $\Psi_t$ is injective on \eqref{eq:normal-form-fixed-scaled-set}, is a diffeomorphism onto its
image, and has a uniformly Lipschitz inverse.  In particular,
\[
\|D\Psi_t\|_{L^\infty} + \|D\Psi_t^{-1}\|_{L^\infty} \le C .
\]
This proves \eqref{eq:normal-form-approximation-bound-a} and
\eqref{eq:normal-form-approximation-bound-b}, and proves the improvement of the large normal-form map bootstrap
\eqref{eq:localized-normal-form-map-large-bootstrap}.

\runinhead{Step 3: The cylindrical-volume Jacobian.}
It remains to prove \eqref{eq:normal-form-approximation-bound-c}.  The quotient
$(\Psi_t)_R/R_{\rm sc}$ is first read on the physical half-plane $R_{\rm sc}>0$ and then extended continuously
to $R_{\rm sc}=0$; on the signed negative side used for parity, the same formula is read through the odd/even
continuation of the map components.

The axial flow map labels yield
\[
R=\tfrac{J^2R_{\rm sc}}{A_t(Z)}, \qquad Z=Z_t(Z_{\rm sc}), \qquad B_t'(Z_t(Z_{\rm sc}))\,\p_{Z_{\rm sc}}Z_t=J^2 ,
\]
and together with the axis volume identity $A_t(Z)^2B_t'(Z)=1$, we obtain
\[
dZ=\tfrac{J^2}{B_t'(Z)}\,dZ_{\rm sc}, \qquad dR=\tfrac{J^2}{A_t(Z)}\,dR_{\rm sc}-\tfrac{J^2R_{\rm sc}A_t'(Z)}{A_t(Z)^2}\,dZ .
\]
The second term in $dR$ is proportional to $dZ$, so it disappears in $dR\wedge dZ$.  Therefore,
\[
R\,dR\,dZ = \tfrac{J^2R_{\rm sc}}{A_t(Z)}\tfrac{J^2}{A_t(Z)}\tfrac{J^2}{B_t'(Z)}\,dR_{\rm sc}\,dZ_{\rm sc}
=J^6R_{\rm sc}\,dR_{\rm sc}\,dZ_{\rm sc}.
\]
On the other hand, the cusp map preserves three-dimensional cylindrical volume.  Since
\[
\phi_{\cusp}(Y_t(\zeta,\tau),t)=J^2\Psi_t(R_{\rm sc},Z_{\rm sc}),
\]
the image cylindrical volume is
\[
r\,dr\,dz=J^6(\Psi_t)_R\det D_{R_{\rm sc},Z_{\rm sc}}\Psi_t\,dR_{\rm sc}\,dZ_{\rm sc}.
\]
Comparing $R\,dR\,dZ=J^6R_{\rm sc}\,dR_{\rm sc}\,dZ_{\rm sc}$ with the image volume identity for
$r\,dr\,dz$, we obtain
\[
(\Psi_t)_R\det D_{R_{\rm sc},Z_{\rm sc}}\Psi_t=R_{\rm sc},
\]
which is exactly $\mathcal J_{\Psi_t}=1$.
\end{proof}

The strain estimate in Lemma~\ref{lem:transported-cusp-field-bounds} uses the part of the transported cusp
vorticity coming from the whole bounded label core $0\le Z\le R_{\tail}$.  We therefore define the
Euler-generated axial function on the corresponding full $\zeta$-interval, not only on the origin-attached
interval $I_a$.  At the fixed time $t$, write
$J=J_{\cusp}(t)$ and
\[
\phi_{\cusp}(R,Z,t)=\bigl(r_t(R,Z),z_t(R,Z)\bigr), \qquad A_t(Z)=\p_Rr_t(0,Z), \qquad B_t(Z)=z_t(0,Z).
\]
On the fixed tail interval $0\le Z\le R_{\tail}$ we use
\begin{equation}
B_t(0)=0,\qquad JA_t(0)=1,
\label{eq:bounded-core-axis-origin-values}
\end{equation}
and, for $0<Z\le R_{\tail}$,
\begin{equation}
c\le JA_t(Z)\le C,\qquad c\le J^{-2}B_t'(Z)\le C,\qquad cZ\le J^{-2}B_t(Z)\le CZ .
\label{eq:bounded-core-axis-estimates-tail}
\end{equation}
Hence $Z\mapsto J^{-2}B_t(Z)$ maps $[0,R_{\tail}]$ onto
$[0,\zeta_{\max}(t)]$, where
\begin{equation}
\zeta_{\max}(t):=J^{-2}B_t(R_{\tail}), \qquad I_{\rm all}(t):=[0,\zeta_{\max}(t)].
\label{eq:zeta-tail-zeta-max-def}
\end{equation}
We denote the inverse on this interval by
\begin{equation}
Z_t:I_{\rm all}(t)\to[0,R_{\tail}], \qquad J^{-2}B_t(Z_t(\zeta))=\zeta .
\label{eq:zeta-tail-axis-inverse-def}
\end{equation}
The Euler-generated axial function \eqref{eq:scaled-profile-at} is extended to $I_{\rm all}(t)$, and then by
zero outside $I_{\rm all}(t)$, by setting
\begin{equation}
a_t^{\rm phys}(\zeta) :=
\begin{cases}
\bigl(J A_t(Z_t(\zeta))\bigr)^{1-\alpha}
\bigl(1+Z_t(\zeta)^2\bigr)^{-\gamma/2},
& 0\le\zeta\le \zeta_{\max}(t),\\
0, & \zeta>\zeta_{\max}(t).
\end{cases}
\label{eq:full-physical-zeta-profile-def}
\end{equation}
The origin-attached axial function used in the model pressure Hessian estimate is
\begin{equation}
a_t(\zeta):=a_t^{\rm phys}(\zeta)\mathbf 1_{I_a}(\zeta).
\label{eq:truncated-pressure-coeff}
\end{equation}
The part $a_t^{\rm phys}\mathbf 1_{I_{\rm all}(t)\setminus I_a}$ is treated as a $\zeta$-tail in
\eqref{eq:pressure-zeta-tail}.  The following lemma proves the integrable upper bound for
$a_t^{\rm phys}$ and the lower bound for the endpoint $\zeta_{\max}(t)$.

\begin{lemma}[Algebraic upper bound for $a_t^{\rm phys}$]
\label{lem:physical-zeta-profile-envelope}
There are constants $c_{\rm env},C_{\rm env}>0$, independent of $t$ and $J_{\cusp}(t)$, such that
\begin{subequations} 
\begin{align} 
0\le a_t^{\rm phys}(\zeta) &\le C_{\rm env}(1+\zeta^2)^{-\gamma/2} \ \ \text{ for } \ \ \zeta\ge0,
\label{eq:physical-zeta-profile-envelope}
\\
\zeta_{\max}(t)&\ge c_{\rm env}R_{\tail}.    \label{eq:physical-zeta-coeff-tail-range}
\end{align} 
\end{subequations} 
\end{lemma}

\begin{proof}[Proof of Lemma~\ref{lem:physical-zeta-profile-envelope}]
By \eqref{eq:full-physical-zeta-profile-def}, $a_t^{\rm phys}(\zeta)=0$ for
$\zeta>\zeta_{\max}(t)$, and
\[
Z_t(0)=0,\qquad JA_t(0)=1,\qquad a_t^{\rm phys}(0)=1
\]
by \eqref{eq:bounded-core-axis-origin-values} and \eqref{eq:zeta-tail-axis-inverse-def}.  For $0<\zeta\le\zeta_{\max}(t)$, \eqref{eq:bounded-core-axis-estimates-tail} 
and \eqref{eq:zeta-tail-axis-inverse-def} imply that
\[
0<c\le JA_t(Z_t(\zeta))\le C, \qquad c\zeta\le Z_t(\zeta)\le C\zeta .
\]
Therefore,
\[
0\le a_t^{\rm phys}(\zeta) \le C(1+Z_t(\zeta)^2)^{-\gamma/2} \le C(1+\zeta^2)^{-\gamma/2}.
\]
Also, by \eqref{eq:zeta-tail-zeta-max-def} and \eqref{eq:bounded-core-axis-estimates-tail}, $\zeta_{\max}(t) =J^{-2}B_t(R_{\tail})\ge  cR_{\tail}$.
\end{proof}

The next estimate has two roles.  First, it identifies the sign and singular size of the stagnation-point
axial strain in \eqref{eq:Wcusp-scaling}.  Second, it shows that, on the shrinking spatial scale
$|x|\lesssim J_{\cusp}(t)^2$, the difference between the cusp velocity and the linear incompressible
hyperbolic field \eqref{eq:linear-hyperbolic-strain-field} with the same axial strain is controlled at the same
singular scale; this is the content of \eqref{eq:Ucusp-gradient-defect} and
\eqref{eq:Ucusp-radial-defect}.

For a real number $\mu$,  we set
\begin{equation}
u_{\hyp}[\mu](r,z) := \bigl(-\tfrac12\mu r,\mu z\bigr), \qquad\text{with no swirl component.}
\label{eq:linear-hyperbolic-strain-field}
\end{equation}
Thus
\[
\p_z(u_{\hyp}[\mu])_z=\mu, \qquad \p_r(u_{\hyp}[\mu])_r=-\tfrac12\mu, \qquad 2\p_r(u_{\hyp}[\mu])_r+\p_z(u_{\hyp}[\mu])_z=0.
\]

\begin{lemma}[Cusp-flow transported strain and velocity bounds]
\label{lem:transported-cusp-field-bounds}
Set $J:=J_{\cusp}(t)$.  Assume that the axis-geometry bound in
\eqref{eq:entry-axis-bounds-statement}, the two-sided clock bootstrap in
\eqref{eq:localized-clock-bootstrap}, and the normal-form estimates
\eqref{eq:additive-normal-form}--\eqref{eq:additive-normal-form-bound} on
$\mathcal Q_t^\sharp$, defined in \eqref{eq:Qt-sharp-def}, are available at the time under
consideration.  Assume also that, for some $I_{\rm str}\Subset I_\sharp$, the inverse axial map
$Z_t=(J^{-2}B_t)^{-1}$ satisfies
\[
0<Z_-\le Z_t(\zeta)\le Z_+<R_{\tail} \qquad(\zeta\in I_{\rm str}),
\]
with constants independent of the small clock.
Then, after decreasing the small-clock thresholds if necessary, there are
constants
\[
c_*=c_*(\sigma_*)\in(0,\tfrac18], \qquad 0<\mathfrak J_{\mathrm{velocity}}\le\mathfrak J_{\mathrm{strain}}\le1,
\qquad 0<c_W\le C_W<\infty, \qquad C<\infty,
\]
depending only on $\alpha,\gamma,\sigma_{\inn},\sigma_*$ and on the constants in the assumptions listed above,
such that the following estimates hold.

If $J\le\mathfrak J_{\mathrm{strain}}$, then
\begin{equation}
\mathcal W_{\cusp}(t)<0, \qquad c_W\Gamma J^{3\alpha-1} \le |\mathcal W_{\cusp}(t)| \le C_W\Gamma J^{3\alpha-1},
\label{eq:Wcusp-scaling}
\end{equation}
Moreover, if $J\le\mathfrak J_{\mathrm{velocity}}$, then
\begin{subequations}
\begin{align}
\|\nabla U_{\cusp}(\cdot,t)\|_{L^\infty(\mathcal C_*)}
&\le C\,\Gamma J^{3\alpha-1},
\label{eq:Ucusp-grad-Linf}\\
[\nabla U_{\cusp}(\cdot,t)]_{C^\alpha(B(x,c_*|x|))}
&\le C\,\Gamma J^{3\alpha-1}|x|^{-\alpha}
\quad\text{if }B(x,2c_*|x|)\subset\mathcal C_*.
\label{eq:Ucusp-grad-scale-local}
\end{align}
\end{subequations}
For every $C_{\rm sc}<\infty$, there are constants
\[
\mathfrak J_{\mathrm{local}}(C_{\rm sc})\in(0,\mathfrak J_{\mathrm{velocity}}], \qquad C_{\rm loc}(C_{\rm sc})<\infty,
\]
with additional dependence only on $C_{\rm sc}$, such that if
$J\le\mathfrak J_{\mathrm{local}}(C_{\rm sc})$, then
\begin{equation}
\big|\nabla U_{\cusp}(x,t) -\nabla u_{\hyp}[\mathcal W_{\cusp}(t)]\big| \le C_{\rm loc}(C_{\rm sc})\,\Gamma J^{3\alpha-1}
\ \ \text{ for } \ \ x\in\mathcal C_*,\ |x|\le C_{\rm sc}J^2,
\label{eq:Ucusp-gradient-defect}
\end{equation}
and
\begin{equation}
\tfrac{\big|(U_{\cusp}-u_{\hyp}[\mathcal W_{\cusp}])_r(x,t)\big|} {r(x)} \le C_{\rm loc}(C_{\rm sc})\,\Gamma J^{3\alpha-1}
\ \ \text{ for } \ \ x\in\mathcal C_*,\ |x|\le C_{\rm sc}J^2 .
\label{eq:Ucusp-radial-defect}
\end{equation}
The quotient in \eqref{eq:Ucusp-radial-defect} is interpreted by its continuous axis value at $r(x)=0$.
\end{lemma}

\begin{remark}[The clock thresholds in Lemma~\ref{lem:transported-cusp-field-bounds}]
The constants
\[
c_*, \qquad \mathfrak J_{\mathrm{strain}},\qquad \mathfrak J_{\mathrm{velocity}},\qquad \mathfrak J_{\mathrm{local}}(C_{\rm sc})
\]
are fixed before the estimates are applied.  The number $c_*$ fixes the balls used in the H\"older estimate inside the cone:
\eqref{eq:Ucusp-grad-scale-local} is asserted only when
$B(x,2c_*|x|)\subset\mathcal C_*$.  The clock thresholds are chosen in the order
\[
0<\mathfrak J_{\mathrm{local}}(C_{\rm sc}) \le \mathfrak J_{\mathrm{velocity}} \le \mathfrak J_{\mathrm{strain}} \le 1 .
\]
The threshold $\mathfrak J_{\mathrm{strain}}$ is the small-clock range in which the stagnation-point axial
strain satisfies \eqref{eq:Wcusp-scaling}.  The possibly smaller threshold
$\mathfrak J_{\mathrm{velocity}}$ is the range in which we prove the cone estimates
\eqref{eq:Ucusp-grad-Linf}--\eqref{eq:Ucusp-grad-scale-local}.  These two thresholds are chosen
without reference to any later restriction of $x$ to a set of size $O(J_{\cusp}(t)^2)$.

The constant $C_{\rm sc}$ is fixed only when a later argument restricts $x$ to
\[
\mathcal C_*\cap\{|x|\le C_{\rm sc}J_{\cusp}(t)^2\}.
\]
After this choice, the clock may be decreased once more to
$\mathfrak J_{\mathrm{local}}(C_{\rm sc})$, and
\eqref{eq:Ucusp-gradient-defect}--\eqref{eq:Ucusp-radial-defect} compare $U_{\cusp}$ with the
linear field $u_{\hyp}[\mathcal W_{\cusp}]$ on that same spatial scale.
\end{remark}

\begin{proof}[Proof of Lemma~\ref{lem:transported-cusp-field-bounds}]
The strain estimate \eqref{eq:Wcusp-scaling} is the stagnation-point part of the computation below.  The
threshold $\mathfrak J_{\mathrm{velocity}}\le \mathfrak J_{\mathrm{strain}}$ will be decreased at the end of
the proof.

\runinhead{Step 1: Local form of the cusp-flow transported vorticity.}
We fix a local radius
$c_*=c_*(\sigma_*)\in(0,\tfrac18]$, which may be decreased below.  Thus
$B(x,2c_*|x|)\subset\mathcal C_*$ means that a ball whose radius is a fixed fraction of $|x|$ remains inside
the fixed buffered cone; the constant $2$ leaves room for a cutoff supported on the larger ball while the
Calder\'on--Zygmund estimate is taken on the smaller ball $B(x,c_*|x|)$.  For the H\"older estimate on this
ball we fix $x$ with $B(x,2c_*|x|)\subset\mathcal C_*$ and set
\[
B_x:=B(x,2c_*|x|)\subset\mathcal C_*.
\]
The pointwise $L^\infty$ estimate is obtained from the same decomposition with the observation point equal to $x$; for points on the
boundary of $\mathcal C_*$ the ball is taken in a slightly wider cone still separated from the equator, with constants depending only on $\sigma_*$.
In this paragraph $y$ denotes an Eulerian point in the transported field, while $Y=Y(y,t):=\phi_{\cusp}^{-1}(y,t)$ denotes its cusp-flow label.

If $B_x$ meets the cusp-flow image of the label tube defined by the time-$t$ axis coordinates,
$\mathcal Q_t^\sharp$ in \eqref{eq:Qt-sharp-def}, then the leading part of
$\bs\Omega_{\cusp}$ can be written explicitly.  In this vorticity computation, $R$ denotes the nonnegative
physical cylindrical label radius, as in the convention after \eqref{eq:Qt-sharp-def}; hence
$0\le\tau\le C_0$.  On $B_x^\sharp$, the definition
\eqref{eq:Omega-cusp-def}, the identity following from \eqref{eq:Theta-star-def},
\[
\rho(Y)^\alpha\Theta^*(\sigma(Y)) = R(Y)^\alpha\Upsilon(\sigma(Y)),
\]
and the Jacobian identity $\mathcal J_{\cusp}^{-1}(Y,t)=r(y)/R(Y)$ provide
\[
\bs\Omega_{\cusp}(y,t) = -\Gamma\,r(y)\,R(Y)^{\alpha-1} (1+\rho(Y)^2)^{-\gamma/2}\Upsilon(\sigma(Y))\,\bs e_\theta(y).
\]
Thus
\[
\bs\Omega_{\cusp}(y,t) = -\Gamma\,J^{\alpha-1}r(y)^\alpha\, \mathfrak A_t^{\rm drv}(y,t)\,\bs e_\theta(y),
\]
where
\[
\mathfrak A_t^{\rm drv}(y,t):=\left(\tfrac{J\,r(y)}{R(Y(y,t))}\right)^{1-\alpha} (1+\rho(Y(y,t))^2)^{-\gamma/2}\Upsilon(\sigma(Y(y,t))).
\]
The axial flow map normal form \eqref{eq:additive-normal-form} states that, for
$Y=Y_t(\zeta,\tau)\in\mathcal Q_t^\sharp$, $0\le\tau\le C_0$, and $y=\phi_{\cusp}(Y,t)$,
\[
y=J^2\zeta\bigl((\tau,1)+\mathcal E_t(\zeta,\tau)\bigr), \qquad c\,J^3\zeta|\tau|\le R(Y)\le C\,J^3\zeta|\tau|.
\]
Combining this with \eqref{eq:adapted-labels}, we obtain the label identity
\[
R(Y_t(\zeta,\tau)) = \tfrac{J^3\zeta\tau}{J A_t(Z_t(\zeta))}.
\]
Since
\[
r(y)=J^2\zeta\bigl(\tau+\mathcal E_{t,r}(\zeta,\tau)\bigr),
\]
we have, for $0\le\tau\le C_0$,
\begin{equation}
\tfrac{Jr(y)}{R(Y_t(\zeta,\tau))} = J A_t(Z_t(\zeta)) \left(1+\tfrac{\mathcal E_{t,r}(\zeta,\tau)}{\tau}\right).
\label{eq:radial-flatness-current-source-ratio}
\end{equation}
The quotient in \eqref{eq:radial-flatness-current-source-ratio} is interpreted at $\tau=0$ by its
continuous value
$\p_\tau\mathcal E_{t,r}(\zeta,0)$, and
\[
\left|\tfrac{\mathcal E_{t,r}(\zeta,\tau)}{\tau}\right| \le \|\p_\tau\mathcal E_{t,r}\|_{L^\infty} \le C J^{3\beta_{\rm ax}}
\qquad (0\le\tau\le C_0),
\]
because $\mathcal E_{t,r}(\zeta,0)=0$.  Therefore, \eqref{eq:entry-axis-bounds-statement},
\eqref{eq:radial-flatness-current-source-ratio}, and \eqref{eq:additive-normal-form-bound} imply
\[
\big\|\big(\tfrac{Jr}{R\circ Y(\cdot,t)}\big)^{\!1-\alpha}\big\|_{L^\infty(B_x)}
+ |x|^\alpha \big[ \big(\tfrac{Jr}{R\circ Y(\cdot,t)}\big)^{\! 1-\alpha} \big]_{C^\alpha(B_x)} \le C.
\]
Let
\[
B_x^\sharp := \{\,y\in B_x:\ Y(y,t)\in\mathcal Q_t^\sharp\,\}.
\]
On $B_x^\sharp$, the relation
\[
y=J^2\zeta\bigl((\tau,1)+\mathcal E_t(\zeta,\tau)\bigr)
\]
and the comparability of the axial flow map labels before the singular small-clock limit imply
\[
c\,J^{-2}|y|\le \rho(Y(y,t))\le C\,J^{-2}|y| \qquad (y\in B_x^\sharp)
\]
after decreasing $c_*$ if necessary.  Since $|y|\simeq |x|$ on $B_x$, the scalar $\mathfrak A_t^{\rm drv}$
satisfies the following bound on $B_x^\sharp$:
\[
\|\mathfrak A_t^{\rm drv}(\cdot,t)\|_{L^\infty(B_x^\sharp)}
+ |x|^\alpha[\mathfrak A_t^{\rm drv}(\cdot,t)]_{C^\alpha(B_x^\sharp)} \le C\bigl(1+J^{-4}|x|^2\bigr)^{-\gamma/2}.
\]
Here the $C^\alpha$ seminorm is in the Eulerian variable $y$; it follows from
\eqref{eq:radial-flatness-current-source-ratio}, the $C^{1,\alpha}$ regularity of the cusp map on the
tube contained in the cone, and the smooth dependence of
$(1+\rho^2)^{-\gamma/2}$ and $\Upsilon$ on the label variables.  The normal-form estimate
\eqref{eq:additive-normal-form-bound} is used here to compare $J^{-2}y$ with $(\zeta\tau,\zeta)$, not as the
source of the $C^\alpha$ exponent.  By the standard scalar Whitney extension
theorem for H\"older functions on subsets of a ball, after increasing the constant by a dimensional amount we
extend $\mathfrak A_t^{\rm drv}$ to a function $\mathfrak A_t$ on all of $B_x$ with
\begin{equation}
\|\mathfrak A_t(\cdot,t)\|_{L^\infty(B_x)} + |x|^\alpha[\mathfrak A_t(\cdot,t)]_{C^\alpha(B_x)} \le C\bigl(1+J^{-4}|x|^2\bigr)^{-\gamma/2}.
\label{eq:Atr-local-bounds}
\end{equation}

\runinhead{Step 2: Cone-local H\"older bound for the vorticity.}
We define $\bs\Omega_{\reg}$ on $B_x$ by subtracting the leading term involving $\mathfrak A_t$ from
the cusp-flow transported vorticity.  On $B_x\setminus B_x^\sharp$, the remainder includes all labels outside
$\mathcal Q_t^\sharp$: the $\zeta$-coordinate complement, the large angular-slope complement, and
the algebraic tail $Y\notin D_{\core}$.  The $\zeta$-tail and angular contributions are controlled by the annular
estimates used for \eqref{eq:pressure-zeta-tail} and
\eqref{eq:pressure-angular-tail}, while the algebraic tail uses \eqref{eq:gamma-stand}.  Therefore,
on $B_x$,
\begin{equation}
\bs\Omega_{\cusp}(y,t) = -\Gamma\,J^{\alpha-1}\,r(y)^\alpha\, \mathfrak A_t(y,t)\,\bs e_\theta(y) + \bs\Omega_{\reg}(y,t),
\label{eq:Omega-cusp-local-representation}
\end{equation}
where the scalar $\mathfrak A_t$ satisfies \eqref{eq:Atr-local-bounds}, and the regular part obeys
\begin{equation}
\|\bs\Omega_{\reg}(\cdot,t)\|_{L^\infty(B_x)} + |x|^\alpha[\bs\Omega_{\reg}(\cdot,t)]_{C^\alpha(B_x)} \le C\Gamma J^{3\alpha-1}.
\label{eq:Omega-reg-local-bounds}
\end{equation}
If
\[
B_x\cap\phi_{\cusp}\bigl(D_{\core},t\bigr)=\emptyset,
\]
then no label in the bounded core $D_{\core}$ contributes to the ball $B_x$.  In that case the leading
cusp term in
\eqref{eq:Omega-cusp-local-representation} is absent and the whole contribution on $B_x$ is included in
$\bs\Omega_{\reg}$.

Since $y\in B_x\subset\mathcal C_*$ implies $r(y)\le C|x|$, Lemma~\ref{lem:toroidal-vector-holder} and
\eqref{eq:Atr-local-bounds} yield
\[
\|r^\alpha\mathfrak A_t\|_{L^\infty(B_x)} + |x|^\alpha[r^\alpha\mathfrak A_t]_{C^\alpha(B_x)} \le C
|x|^\alpha\bigl(1+J^{-4}|x|^2\bigr)^{-\gamma/2}.
\]
Writing $|x|=J^2s$, the right-hand side is
\[
C J^{2\alpha}s^\alpha(1+s^2)^{-\gamma/2} \le C J^{2\alpha},
\]
because $\gamma>\alpha$.  Multiplication by the multiplier $J^{\alpha-1}$ in
\eqref{eq:Omega-cusp-local-representation}, together with
\eqref{eq:Omega-reg-local-bounds}, therefore implies
\begin{equation}
\|\bs\Omega_{\cusp}(\cdot,t)\|_{L^\infty(B(x,2c_*|x|))}
+ |x|^\alpha [\bs\Omega_{\cusp}(\cdot,t)]_{C^\alpha(B(x,2c_*|x|))} \le C\,\Gamma J^{3\alpha-1}.
\label{eq:Omega-cusp-local-holder}
\end{equation}

\runinhead{Step 3: Cone-local Biot--Savart estimate.}
We convert \eqref{eq:Omega-cusp-local-holder} into velocity-gradient estimates.  We choose a smooth
cutoff $\eta_x$ equal to one on $B(x,c_*|x|)$, supported in $B_x$, and satisfying
$|\nabla^k\eta_x|\le C_k|x|^{-k}$.  For the localized vorticity $\eta_x\bs\Omega_{\cusp}$, the reduced axisymmetric Biot--Savart gradient
has the Calder\'on--Zygmund form
\[
\nabla\BS[\eta_x\bs\Omega_{\cusp}] = \pv\!\int \nabla K(\cdot,y)\,\eta_x(y)\bs\Omega_{\cusp}(y,t)\,\ud y
+\mathbf C\,\eta_x\bs\Omega_{\cusp},
\]
and the scale-invariant Schauder estimate on the ball $B_x$ yields
\begin{equation}
\|\nabla\BS[\eta_x\bs\Omega_{\cusp}]\|_{L^\infty(B(x,c_*|x|))}
+ |x|^\alpha [\nabla\BS[\eta_x\bs\Omega_{\cusp}]]_{C^\alpha(B(x,c_*|x|))} \le C\Gamma J^{3\alpha-1}.
\label{eq:Ucusp-local-CZ}
\end{equation}
For the complementary vorticity $(1-\eta_x)\bs\Omega_{\cusp}$, no principal value remains.  We include the
short annular estimate.  We set
\[
\mathscr S_k(x):=\{y:\ 2^kc_*|x|\le |y-x|\le 2^{k+1}c_*|x|\}, \qquad k\ge0.
\]
For $z,z'\in B(x,c_*|x|)$ and $y\in \mathscr S_k(x)$,
\begin{equation}
|\nabla_z K(z,y)|\le C(2^k|x|)^{-3}, \qquad |\nabla_zK(z,y)-\nabla_zK(z',y)| \le C\,|z-z'|^\alpha(2^k|x|)^{-3-\alpha}.
\label{eq:separated-kernel-annular}
\end{equation}
The local representation \eqref{eq:Omega-cusp-local-representation} on bounded annuli, together with the algebraic
tail in the initial datum, yields the annular mass bound
\begin{equation}
\int_{\mathscr S_k(x)}|\bs\Omega_{\cusp}(y,t)|\,\ud y \le C\,\Gamma J^{3\alpha-1}(2^k|x|)^3 a_k, \qquad \sum_{k\ge0} a_k\le C.
\label{eq:Omega-cusp-annular-mass}
\end{equation}
Here the summability for large annuli uses $\gamma>\alpha+\tfrac52$; on the bounded annuli the constants are
controlled by the bounded weights in \eqref{eq:Atr-local-bounds}.  The kernel estimate
\eqref{eq:separated-kernel-annular} only uses that the source point $y$ is separated from the observation ball
$B(x,c_*|x|)$.  Combining \eqref{eq:separated-kernel-annular} and \eqref{eq:Omega-cusp-annular-mass}, we obtain
\[
\sum_{k\ge0} (2^k|x|)^{-3} \int_{\mathscr S_k(x)}|\bs\Omega_{\cusp}(y,t)|\,\ud y \le C\Gamma J^{3\alpha-1},
\]
and
\[
\sum_{k\ge0} |z-z'|^\alpha(2^k|x|)^{-3-\alpha} \int_{\mathscr S_k(x)}|\bs\Omega_{\cusp}(y,t)|\,\ud y
\le C\Gamma J^{3\alpha-1}|x|^{-\alpha}|z-z'|^\alpha .
\]
Therefore
\begin{equation}
\|\nabla\BS[(1-\eta_x)\bs\Omega_{\cusp}]\|_{L^\infty(B(x,c_*|x|))}
+ |x|^\alpha [\nabla\BS[(1-\eta_x)\bs\Omega_{\cusp}]]_{C^\alpha(B(x,c_*|x|))} \le C\Gamma J^{3\alpha-1}.
\label{eq:Ucusp-separated-CZ}
\end{equation}
Equations \eqref{eq:Ucusp-local-CZ}--\eqref{eq:Ucusp-separated-CZ} imply
\eqref{eq:Ucusp-grad-Linf}--\eqref{eq:Ucusp-grad-scale-local}.

\runinhead{Step 4: Stagnation-point axial strain.}
We now identify $\mathcal W_{\cusp}(t)=\p_z(U_{\cusp})_z(0,t)$.  Since the stagnation-point axial strain is
linear in the vorticity, we evaluate the leading term after dividing the Eulerian image by $J^2$ and using
$(R_{\rm sc},Z_{\rm sc})=(\zeta\tau,\zeta)$ from \eqref{eq:normal-form-fixed-scaled-set}.  By the strain
identity \eqref{eq:scaled-strain-exact}, the angular integration leaves the axial integral
\begin{equation}
\mathfrak I_t := C_\alpha^W\int_0^\infty a_t^{\rm phys}(\zeta)\zeta^{\alpha-1}\,d\zeta, \qquad C_\alpha^W>0.
\label{eq:current-strain-amplitude}
\end{equation}
Here $C_\alpha^W$ is the angular constant from \eqref{eq:scaled-strain-exact}.  The sign convention in the
transported vorticity \eqref{eq:Omega-cusp-def} yields the leading contribution
$-\Gamma J^{3\alpha-1}\mathfrak I_t$.  Thus
\begin{equation}
\mathcal W_{\cusp}(t) = -\Gamma J^{3\alpha-1}\mathfrak I_t + O\!\left(\Gamma J^{3\alpha-1}J^{3\beta_{\rm ax}}\right) +O(\Gamma).
\label{eq:Wcusp-asymptotic-field}
\end{equation}
We next bound $\mathfrak I_t$ above and below.  We fix a compact interval $I_0\Subset I_{\rm str}$.  The
axis-geometry assumption \eqref{eq:entry-axis-bounds-statement} gives, for $\zeta\in I_\sharp$, $c_{\rm ax}\le J A_t(Z_t(\zeta))\le C_{\rm ax}$ and
$c_{\rm ax}\le J^{-2}B_t'(Z_t(\zeta))\le C_{\rm ax}$.  Since $I_0\subset I_{\rm str}$, the compact-containment
assumption for $Z_t$ in Lemma~\ref{lem:transported-cusp-field-bounds} shows that
\[
Z_t(I_0)\Subset(0,R_{\tail})
\]
with constants independent of the small clock.  Thus the algebraic weight $(1+Z_t(\zeta)^2)^{-\gamma/2}$ is bounded below on $I_0$, and hence
$a_t^{\rm phys}(\zeta)\ge c>0$ there.  Therefore,
\[
\int_0^\infty a_t^{\rm phys}(\zeta)\zeta^{\alpha-1}\,d\zeta \ge c\int_{I_0}\zeta^{\alpha-1}\,d\zeta \ge c_I>0 .
\]
For the upper bound, Lemma~\ref{lem:physical-zeta-profile-envelope} states that
\[
0\le a_t^{\rm phys}(\zeta)\le C(1+\zeta^2)^{-\gamma/2}.
\]
Therefore, since $\gamma>\alpha$,
\[
\int_0^\infty a_t^{\rm phys}(\zeta)\zeta^{\alpha-1}\,d\zeta \le C\int_0^\infty (1+\zeta^2)^{-\gamma/2}\zeta^{\alpha-1}\,d\zeta \le C_I .
\]
Together with \eqref{eq:current-strain-amplitude}, this implies
\[
0<c\le \mathfrak I_t\le C<\infty
\]
uniformly in $t$.  Dividing the two errors in \eqref{eq:Wcusp-asymptotic-field} by $\Gamma J^{3\alpha-1}$ yields
$O(J^{3\beta_{\rm ax}}+J^{1-3\alpha})$.  Since $0<\alpha<\tfrac13$, we choose
$\mathfrak J_{\mathrm{strain}}$ so that this relative error is smaller than one half of the lower bound for
$\mathfrak I_t$ whenever $J\le \mathfrak J_{\mathrm{strain}}$.  This proves \eqref{eq:Wcusp-scaling}.  We
also choose $\mathfrak J_{\mathrm{velocity}}\le \mathfrak J_{\mathrm{strain}}$.

\runinhead{Step 5: Comparison with the linear stagnation field on $|x|\lesssim J^2$.}
For the bounds relative to the linear stagnation field \eqref{eq:linear-hyperbolic-strain-field} on the spatial
scale $|x|\le C_{\rm sc}J^2$, we fix $C_{\rm sc}<\infty$ and set
$\mathfrak J_{\mathrm{local}}(C_{\rm sc}):=\mathfrak J_{\mathrm{velocity}}$.  The matrix
$\nabla u_{\hyp}[\mathcal W_{\cusp}(t)]$ has size $C|\mathcal W_{\cusp}(t)|$, and
\eqref{eq:Wcusp-scaling} yields
\[
|\nabla u_{\hyp}[\mathcal W_{\cusp}(t)]| \le C\Gamma J^{3\alpha-1}.
\]
Together with \eqref{eq:Ucusp-grad-Linf}, this implies the gradient bound
\eqref{eq:Ucusp-gradient-defect} on $|x|\le C_{\rm sc}J^2$ by the triangle estimate
\[
\big|\nabla U_{\cusp}(x,t)-\nabla u_{\hyp}[\mathcal W_{\cusp}(t)]\big| \le |\nabla U_{\cusp}(x,t)|+|\nabla u_{\hyp}[\mathcal W_{\cusp}(t)]|
\le C\Gamma J^{3\alpha-1}.
\]
For the radial component, both
$(U_{\cusp})_r$ and $(u_{\hyp}[\mathcal W_{\cusp}])_r$ vanish on the symmetry axis.  If
$x=(r,z)\in\mathcal C_*$ and $|x|\le C_{\rm sc}J^2$, then the segment
$\{(sr,z):0\le s\le1\}$ stays in a cone whose constants depend only on
$\sigma_*$.  Therefore, by \eqref{eq:Ucusp-gradient-defect},
\[
\big|(U_{\cusp}-u_{\hyp}[\mathcal W_{\cusp}])_r(r,z,t)\big| \le \int_0^r \big|\p_r(U_{\cusp}-u_{\hyp}[\mathcal W_{\cusp}])_r(s,z,t)\big|\,ds
\le C_{\rm loc}(C_{\rm sc})\Gamma J^{3\alpha-1}r,
\]
which proves \eqref{eq:Ucusp-radial-defect}.
\end{proof}


\section{Pressure Hessian for the Euler Cusp Velocity}
\label{sec:transported-cusp-pressure}

The purpose of this section is to verify, for the Euler-generated axial function $a_t$, the hypotheses of the
renormalized axis-trace criterion in Proposition~\ref{prop:euler-generated-profile-riccati}.  With
$\Pi_{\cusp}$ defined in \eqref{eq:Pi-U-cusp-def},
$\mathcal W_{\cusp}:=\p_z(U_{\cusp})_z(0,t)$, and $J=J_{\cusp}(t)$, the resulting estimate is
\[
\Pi_{\cusp}(t)\ge -q_{\rm tr}\,\tfrac12\,\mathcal W_{\cusp}(t)^2, \qquad q_{\rm tr}<\upbeta .
\]
The four estimates needed for this verification are the axis-trace approximation
\eqref{eq:axis-trace-approximation}, the normalized axial-function equation
\eqref{eq:renormalized-profile-equation}, the differentiated strain formula
\eqref{eq:direct-strain-derivative-axis}, and the principal Riccati identity
\eqref{eq:cusp-principal-riccati-with-error}.  The geometric errors entering these estimates are the
normal-form displacement \eqref{eq:normal-form-approximation-bound}, the $\zeta$-localization tail
\eqref{eq:pressure-zeta-tail}, the large-slope tail \eqref{eq:pressure-angular-tail}, and the algebraic
far-field tail controlled by \eqref{eq:gamma-stand}.  Lemmas~\ref{lem:transported-cusp-estimates} and
\ref{lem:tail-bound} then collect the transported cusp bounds and the lower-order pressure terms used in
Section~\ref{sec:target-profile-typeI-completion}.

On the upper-half interval used below, where $Z=\zeta\in I_\sharp$, the slope variable is
$\tau=R/Z$.  The choices of $I_\sharp$ and $\vartheta_\sharp$ in
\eqref{eq:pressure-localization-cutoff}--\eqref{eq:pressure-localization-intervals} control the
$\zeta$-tail $\mathfrak a_\zeta(I_\sharp)$ in \eqref{eq:pressure-zeta-tail}.  The angular cutoff
$M_\pressure$ and the slope bound $C_0$ in \eqref{eq:pressure-C0} control the angular tail
$\mathfrak a_{\rm ang}(M_\pressure)$ in \eqref{eq:pressure-angular-tail}.  The tail radius $R_{\tail}$ is
the radius used in the regions $D_{\core},D_{\tail}$ in \eqref{eq:core-tail-domains} and in the far-field
velocity $u_{\smooth}$ in \eqref{eq:smooth-velocity-def}.  After
$I_\sharp,\vartheta_\sharp,M_\pressure,C_0$, and $R_{\tail}$ are fixed, we decrease the pressure
threshold $\mathfrak J_{\pressure}$ in \eqref{eq:small-clock-threshold-order} so that the positive powers of
$J_{\cusp}$ in \eqref{eq:pressure-normal-form-deformation-error} and the algebraic tail controlled by
\eqref{eq:gamma-stand} are small enough for the error bound \eqref{eq:epsilon-ax-sec12}.
The angular cutoff and cone aperture are fixed by
\begin{equation}
\sigma_{\rm wide}\le\sigma_{\inn}<\sigma_*<\tfrac\pi2, \qquad 2M_\pressure\le\tan\sigma_{\inn}\le\tfrac12\tan\sigma_*,
\qquad C_0:=2M_\pressure .
\label{eq:pressure-C0}
\end{equation}
The normal-form displacement errors are measured by the clock powers
\begin{equation}
\varepsilon_{\rm nf}(J):=J^{3\beta_{\rm ax}}, \qquad \varepsilon_{\rm def}(J):=\varepsilon_{\rm nf}(J)^{\frac{\beta_{\rm ax}}{1+\beta_{\rm ax}}}
=J^{\kappa_{\rm def}} .
\label{eq:pressure-normal-form-deformation-error}
\end{equation}

\subsection{Deformation of the localized cusp image after division by $J^2$}

We first estimate the error made when the exact cusp image is divided by $J^2$ and then replaced by the
normal-form image.
Writing the cusp flow map in cylindrical coordinates as
\[
\phi_{\cusp}(R,Z,t)=\bigl(r_t(R,Z),z_t(R,Z)\bigr),
\]
the time-$t$ coordinates along the symmetry axis are determined by
\[
A_t(Z)=\p_Rr_t(0,Z), \qquad B_t(Z)=z_t(0,Z).
\]
Here $R=0$ is the symmetry axis; axis preservation implies $r_t(0,Z)=0$, while $z_t(0,Z)$ is the axial position
at time $t$ of the axis label $Z$.  In these coordinates the normal form \eqref{eq:additive-normal-form},
the physical scaling \eqref{eq:physical-scaled-vorticity}, and the model vorticity
\eqref{eq:vort-slope-restricted} determine the comparison made in Lemma~\ref{lem:normal-form-deformation}.

The reference-domain label is $Y_t(\zeta,\tau)$, and its exact Eulerian position is
\[
x_t(\zeta,\tau):=\phi_{\cusp}(Y_t(\zeta,\tau),t).
\]
After the rescaling $x=J^2x_{\rm sc}$, the undeformed point in the image variables after division by $J^2$ is
\[
X_{\rm sc}:=(R_{\rm sc},Z_{\rm sc})=(\zeta\tau,\zeta),
\]
and \eqref{eq:normal-form-approximation-def} gives
\[
x_{{\rm sc},t}(\zeta,\tau)=J^{-2}x_t(\zeta,\tau)=\Psi_t(X_{\rm sc}).
\]
Thus $X_{\rm sc}$ is the point assigned to the label in the $J^{-2}$-renormalized variables, while
$x_{{\rm sc},t}=\Psi_t(X_{\rm sc})$ is the corresponding $J^{-2}$-renormalized Eulerian image.  The next
lemma proves the estimate for an arbitrary map $\Psi$ which is close to the identity, uniformly bi-Lipschitz,
and has cylindrical volume Jacobian equal to one.  Inside the lemma we write $X_{\rm sc}$ simply as $(R,Z)$.
The upper-half set, independent of $J$, is
\begin{equation}
\mathcal R_{\sharp,C_0}=\{(R,Z)=(\zeta\tau,\zeta):\ \zeta\in I_\sharp,\ 0\le\tau\le C_0\}.
\label{eq:fixed-scale-deformation-set}
\end{equation}
Since $I_\sharp\Subset(0,\infty)$, every point in $\mathcal R_{\sharp,C_0}$ has $Z>0$, and therefore
$\tau=R/Z$ on this $J$-independent set.  The placements $f_{\rm id}$ and $f_\Psi$ are defined in
\eqref{eq:fixed-scale-placed-vorticity}, with corresponding velocities $V_F$ and $V_F^\Psi$ in
\eqref{eq:fixed-scale-placed-velocity}.  The estimates \eqref{eq:normal-form-deformation}--
\eqref{eq:normal-form-deformation-self} compare placement at $(R,Z)$ with placement at $\Psi(R,Z)$, and
bound the resulting changes in the stagnation-point axial strain and in the bilinear pressure Hessian form
\eqref{eq:physical-pressure-bilinear-form}.

We use the following scaling and localization conventions throughout the pressure estimates.  We let
$\Omega_{\rm fix}(\bar r,\bar z,t)$ be a toroidal scalar in the variables obtained after dividing the Eulerian image by
$J_{\cusp}(t)^2$.  Its physical representative at cusp clock $J=J_{\cusp}(t)$ is
\begin{equation}
\Omega_{\rm sc,J}(r,z,t):=\Gamma J^{3\alpha-1}\Omega_{\rm fix}(\bar r,\bar z,t),
\qquad (r,z)=J^2(\bar r,\bar z).
\label{eq:physical-scaled-vorticity}
\end{equation}
The localized label tube is
\begin{equation}
\mathcal Q_t^\sharp:=\{Y_t(\zeta,\tau):\ \zeta\in I_\sharp,\ |\tau|\le C_0\},
\label{eq:localized-label-tube-sec12}
\end{equation}
where $Y_t(\zeta,\tau)$ is the adapted label from \eqref{eq:adapted-labels}.  On
$\phi_{\cusp}(\mathcal Q_t^\sharp,t)$, we define
\begin{equation}
\Omega_\sharp(\phi_{\cusp}(Y_t(\zeta,\tau),t),t) :=\vartheta_\sharp(\zeta)\chi_{M_\pressure}(|\tau|)
\Omega_{\cusp}(\phi_{\cusp}(Y_t(\zeta,\tau),t),t), \qquad \zeta\in I_\sharp,\ |\tau|\le C_0. \label{eq:Omega-sharp-cutoff-def}
\end{equation}
We extend $\Omega_\sharp$ by zero off $\phi_{\cusp}(\mathcal Q_t^\sharp,t)$ and set
$U_\sharp:=\BS[\Omega_\sharp e_\theta]$.  The odd reflection in the axial variable is inherited from
$\Omega_{\cusp}$.

For a map $\Psi:\mathcal R_{\sharp,C_0}\to\R^2$, we define
\begin{equation}
\mathcal J_\Psi(R,Z) := \tfrac{\Psi_R(R,Z)}{R}\det D_{R,Z}\Psi(R,Z), \label{eq:fixed-scale-cylindrical-jacobian}
\end{equation}
where $\Psi_R$ denotes the radial component of $\Psi$, and the quotient is interpreted by its continuous value
at $R=0$.  For an axisymmetric scalar function $h$ on the meridional half-plane, we define the associated
toroidal vorticity vector field
\[
\omega_h(x):=
\begin{cases}
h(r(x),z(x))\,\bs e_\theta(x), & r(x)>0,\\
0, & r(x)=0.
\end{cases}
\]
When $\omega_h\in C^{0,\beta_{\rm ax}}$, we set
\begin{equation}
\mathcal H(h) := \|\omega_h\|_{L^\infty} + [\omega_h]_{C^{\beta_{\rm ax}}}.
\label{eq:fixed-scale-toroidal-holder-norm}
\end{equation}
For a reference function $F$ on $\mathcal R_{\sharp,C_0}$, we define its identity placement and its
deformed placement by
\begin{equation}
f_{\rm id}(R,Z)=F(R,Z), \qquad f_\Psi(\Psi(R,Z))=F(R,Z),
\label{eq:fixed-scale-placed-vorticity}
\end{equation}
with both functions extended by zero off their supports.  We then set
\begin{equation}
V_F:=\BS[\omega_{f_{\rm id}}], \qquad V_F^\Psi:=\BS[\omega_{f_\Psi}].
\label{eq:fixed-scale-placed-velocity}
\end{equation}
The uppercase $F$ denotes the reference function, while $f_{\rm id}$ and $f_\Psi$ are the Eulerian
functions obtained from the two placements.  Thus the regularity measured below is
$\mathcal H(f_\Psi)$, the $C^{0,\beta_{\rm ax}}$ norm of the Eulerian toroidal vorticity
$\omega_{f_\Psi}$, rather than a norm of the reference function $F$ itself.  Finally, set
\[
\theta_{\rm def}:=\tfrac{\beta_{\rm ax}}{1+\beta_{\rm ax}} .
\]

\begin{lemma}[Pressure Hessian response to diffeomorphism deformation]
\label{lem:normal-form-deformation}
Let $I_\sharp\Subset(0,\infty)$ and $C_0<\infty$ be fixed, and let
$\mathcal R_{\sharp,C_0}$ be the set defined in \eqref{eq:fixed-scale-deformation-set}.  Let
$\Psi:\mathcal R_{\sharp,C_0}\to\R^2$ be an axisymmetric, axis-preserving $C^1$ diffeomorphism onto its image.
Assume that $\Psi$ and $\Psi^{-1}$ are uniformly Lipschitz, with constants depending only on
$I_\sharp$ and $C_0$, that $\mathcal J_\Psi=1$, and that, for some
$0<\varepsilon\le\varepsilon_0(I_\sharp,C_0)$,
\[
\|\Psi-\operatorname{Id}\|_{L^\infty} + [\Psi-\operatorname{Id}]_{C^{\beta_{\rm ax}}} \le \varepsilon .
\]
Let $F,G$ be reference functions supported in $\mathcal R_{\sharp,C_0}$.  Let $f_\Psi$ be the deformed
placement of $F$, and let $g_{\rm id}$ be the identity placement of $G$, as in
\eqref{eq:fixed-scale-placed-vorticity}.  If
$\mathcal H(f_\Psi)+\mathcal H(g_{\rm id})<\infty$, then
\begin{equation}
\left| \p_Z(V_F^\Psi-V_F)_Z(0) \right| \le C_{\sharp,C_0}\varepsilon\,\mathcal H(f_\Psi),
\qquad \left| \Pi[V_G,V_F^\Psi-V_F] \right| \le C_{\sharp,C_0}\varepsilon^{\theta_{\rm def}}\,
\mathcal H(f_\Psi)\mathcal H(g_{\rm id}).
\label{eq:normal-form-deformation}
\end{equation}
In particular,
\begin{equation}
\left| \Pi[V_F^\Psi-V_F,V_F^\Psi-V_F] \right| \le C_{\sharp,C_0}\varepsilon^{\theta_{\rm def}}\,\mathcal H(f_\Psi)^2.
\label{eq:normal-form-deformation-self}
\end{equation}
\end{lemma}

\begin{proof}[Proof of Lemma~\ref{lem:normal-form-deformation}]
\runinhead{Step 1: The axial strain integral.}
Let $d\mu=R\,dR\,dZ$.  Let $K$ be the axisymmetric Biot--Savart kernel in \eqref{eq:BS-axisymm}.  For a
meridional point $Y=(R,Z)$, define the stagnation-point axial-strain kernel
\[
\mathsf K_W(Y) := \left.\p_{x_z}\bigl(K(x,Y)\cdot e_z\bigr)\right|_{x=0}.
\]
This is the kernel $\mathcal K_W(0,Y)$ from \eqref{eq:rW2}.  In separation-variable notation we write the
same quantity as $\mathsf K_W(Y)=\p_ZK_Z(-Y)$.

Since $I_\sharp\Subset(0,\infty)$, the set $\mathcal R_{\sharp,C_0}$, defined in
\eqref{eq:fixed-scale-deformation-set}, has positive distance from the stagnation point $(R,Z)=(0,0)$.
After decreasing $\varepsilon_0$, we choose a fixed open set $\mathcal N_{\sharp,C_0}$ such that
\[
\mathcal R_{\sharp,C_0}\cup\Psi(\mathcal R_{\sharp,C_0}) \Subset \mathcal N_{\sharp,C_0} \Subset \R^2\setminus\{(0,0)\}.
\]
On $\mathcal N_{\sharp,C_0}$, the explicit kernel
\[
\mathsf K_W(R,Z)=3\,\tfrac{RZ}{(R^2+Z^2)^{5/2}}
\]
is smooth and satisfies
\begin{equation} 
\|D\mathsf K_W\|_{L^\infty(\mathcal N_{\sharp,C_0})} \le C_{\sharp,C_0}.
\label{eq:DKW}
\end{equation} 

By \eqref{eq:BS-axisymm}, \eqref{eq:fixed-scale-placed-vorticity}, and
\eqref{eq:fixed-scale-placed-velocity}, we have that
\[
\p_Z(V_F^\Psi)_Z(0) = \tfrac14\int_{\Psi(\mathcal R_{\sharp,C_0})} \mathsf K_W(y)f_\Psi(y)\,d\mu(y).
\]
With $y=\Psi(Y)$, the cylindrical measure transforms as
$d\mu(y)= \mathcal J_\Psi(Y)\,d\mu(Y) = d\mu(Y)$,
where we used \eqref{eq:fixed-scale-cylindrical-jacobian} and the assumption $\mathcal J_\Psi=1$.
Since $f_\Psi(\Psi(Y))=F(Y)$ by \eqref{eq:fixed-scale-placed-vorticity}, this shows that
\[
\p_Z(V_F^\Psi)_Z(0) = \tfrac14\int_{\mathcal R_{\sharp,C_0}} \mathsf K_W(\Psi(Y))F(Y)\,d\mu(Y), \qquad \p_Z(V_F)_Z(0) =
\tfrac14\int_{\mathcal R_{\sharp,C_0}} \mathsf K_W(Y)F(Y)\,d\mu(Y),
\]
and hence
\[
\p_Z(V_F^\Psi-V_F)_Z(0) = \tfrac14\int_{\mathcal R_{\sharp,C_0}} \bigl[\mathsf K_W(\Psi(Y))-\mathsf K_W(Y)\bigr]F(Y)\,d\mu(Y).
\]
Using \eqref{eq:DKW} and $|\Psi-\operatorname{Id}|\le\varepsilon$, we obtain the first estimate in
\eqref{eq:normal-form-deformation}; here
$\|F\|_{L^\infty(\mathcal R_{\sharp,C_0})}=\|f_\Psi\|_{L^\infty(\Psi(\mathcal R_{\sharp,C_0}))}$.

\runinhead{Step 2: Vorticity difference created by the diffeomorphism.}
A point in the three-dimensional rotation of $\mathcal R_{\sharp,C_0}$ has the form
\[
y=(R\cos\theta,R\sin\theta,Z), \qquad (R,Z)\in\mathcal R_{\sharp,C_0}, \qquad \bs e_\theta(\theta)=(-\sin\theta,\cos\theta,0).
\]
The three-dimensional axisymmetric lift of $\Psi$ is
\[
\widetilde\Psi(R\cos\theta,R\sin\theta,Z) = \bigl(\Psi_R(R,Z)\cos\theta,\Psi_R(R,Z)\sin\theta,\Psi_Z(R,Z)\bigr).
\]
The identity $\mathcal J_\Psi=1$ in \eqref{eq:fixed-scale-cylindrical-jacobian} is exactly the
volume-preservation identity for this lift, since $dy=R\,dR\,dZ\,d\theta$ in cylindrical coordinates.
We set
\[
\omega_\Psi:=\omega_{f_\Psi}, \qquad \omega_{\rm id}:=\omega_{f_{\rm id}}, \qquad \omega_\Delta:=\omega_\Psi-\omega_{\rm id},
\qquad v_\Delta:=V_F^\Psi-V_F=\BS[\omega_\Delta].
\]
The supports of $\omega_\Psi$ and $\omega_{\rm id}$ are contained in a fixed compact set
$\mathcal K_{\sharp,C_0}\Subset \R^3\setminus\{0\}$ after decreasing $\varepsilon_0$.
Moreover, \eqref{eq:fixed-scale-placed-vorticity} implies that 
\[
\omega_{\rm id}(y)=\omega_\Psi(\widetilde\Psi(y)) \ \ \text{ for } \ \ y\in\operatorname{supp}\omega_{\rm id}.
\]
Indeed, for $R>0$ the point $y$ and the point $\widetilde\Psi(y)$ have the same azimuthal coordinate
$\theta$, so $\bs e_\theta(\widetilde\Psi(y))=\bs e_\theta(\theta)=\bs e_\theta(y)$, while
$f_\Psi(\Psi(R,Z))=f_{\rm id}(R,Z)=F(R,Z)$.  For $R=0$, both sides vanish by the definition of
$\omega_h$ on the symmetry axis.
Since $\omega_\Psi$ is the $C^{0,\beta_{\rm ax}}$ zero extension of its deformed support and
$\widetilde\Psi$ is a homeomorphism onto that support, we have
\begin{equation}
\omega_\Psi|_{\partial\operatorname{supp}\omega_\Psi}=0, \qquad \omega_{\rm id}|_{\partial\operatorname{supp}\omega_{\rm id}}=0 .
\label{eq:placement-boundary-trace-zero}
\end{equation}
The axis-preserving Lipschitz bounds for $\Psi$ and $\Psi^{-1}$ imply the same bounds for
$\widetilde\Psi$ and $\widetilde\Psi^{-1}$ on this compact set.  Hence, for points in
$\operatorname{supp}\omega_{\rm id}$, the H\"older seminorm of $\omega_{\rm id}$ is controlled by composition with
$\widetilde\Psi$.  If $y_1\in\operatorname{supp}\omega_{\rm id}$ and
$y_2\notin\operatorname{supp}\omega_{\rm id}$, let $z$ be the first point of the segment from $y_2$ to $y_1$
which belongs to $\operatorname{supp}\omega_{\rm id}$.  Then $z\in\partial\operatorname{supp}\omega_{\rm id}$,
$\omega_{\rm id}(z)=0$ by \eqref{eq:placement-boundary-trace-zero},
$\omega_{\rm id}(y_2)=0$, and $|y_1-z|\le |y_1-y_2|$, so the preceding estimate applied to $y_1$ and $z$
implies the same H\"older bound.  Therefore,
\begin{equation}
\|\omega_{\rm id}\|_{L^\infty}+[\omega_{\rm id}]_{C^{\beta_{\rm ax}}} \le C_{\sharp,C_0}\mathcal H(f_\Psi).
\label{eq:identity-placement-holder-from-deformed}
\end{equation}
If $y\in\operatorname{supp}\omega_\Psi\cap\operatorname{supp}\omega_{\rm id}$, then
$\omega_{\rm id}(y)=\omega_\Psi(\widetilde\Psi(y))$ and
\eqref{eq:fixed-scale-toroidal-holder-norm} imply
\[
|\omega_\Psi(y)-\omega_{\rm id}(y)| = |\omega_\Psi(y)-\omega_\Psi(\widetilde\Psi(y))|
\le \mathcal H(f_\Psi)\,|y-\widetilde\Psi(y)|^{\beta_{\rm ax}} .
\]
If $y\in\operatorname{supp}\omega_\Psi\setminus\operatorname{supp}\omega_{\rm id}$, then
$y=\widetilde\Psi(y_0)$ for some $y_0\in\operatorname{supp}\omega_{\rm id}$ and $\omega_{\rm id}(y)=0$; using
\eqref{eq:identity-placement-holder-from-deformed},
\[
|\omega_\Psi(y)| = |\omega_{\rm id}(y_0)-\omega_{\rm id}(y)| \le C_{\sharp,C_0}\mathcal H(f_\Psi)\,|y_0-y|^{\beta_{\rm ax}}.
\]
The remaining case $y\in\operatorname{supp}\omega_{\rm id}\setminus\operatorname{supp}\omega_\Psi$ is the same, using
$\omega_{\rm id}(y)=\omega_\Psi(\widetilde\Psi(y))$ and $\omega_\Psi(y)=0$.  Since
$\|\widetilde\Psi-\operatorname{Id}\|_{L^\infty}\le C_{\sharp,C_0}\varepsilon$, these three cases imply
\begin{equation}
\|\omega_\Delta\|_{L^\infty} \le C_{\sharp,C_0}\varepsilon^{\beta_{\rm ax}}\mathcal H(f_\Psi),
\qquad [\omega_\Delta]_{C^{\beta_{\rm ax}}} \le C_{\sharp,C_0}\mathcal H(f_\Psi).
\label{eq:placed-vorticity-difference}
\end{equation}

Let $\bs K$ be the Biot--Savart kernel defined in \eqref{eq:BS-3D}.  For $z\ne0$ and $\xi\in\R^3$, define
the matrix kernel $L(z):\R^3\to\R^{3\times3}$ by
\begin{equation}
\bigl(L(z)\xi\bigr)_{ij} := \p_{z_j}\bigl(\bs K(z)\times\xi\bigr)_i.
\label{eq:biot-savart-gradient-kernel}
\end{equation}
Thus $L$ is the kernel of $\nabla\BS$.  For $0<\rho<1$,
\[
\nabla v_\Delta(x) = \pv\!\int_{|x-y|\le\rho}L(x-y)\omega_\Delta(y)\,dy + \int_{|x-y|>\rho}L(x-y)\omega_\Delta(y)\,dy.
\]
The kernel $L$ is homogeneous of degree $-3$ and has mean zero on spheres.  Hence
\[
\big| \pv\!\int_{|x-y|\le\rho}L(x-y)\omega_\Delta(y)\,dy \big| =
\big| \pv\!\int_{|x-y|\le\rho}L(x-y)\bigl(\omega_\Delta(y)-\omega_\Delta(x)\bigr)\,dy \big|
\le C[\omega_\Delta]_{C^{\beta_{\rm ax}}}\rho^{\beta_{\rm ax}},
\] 
while the fixed compact support of $\omega_\Delta$ implies that
\[
\Big| \int_{|x-y|>\rho}L(x-y)\omega_\Delta(y)\,dy\Big| \le C\|\omega_\Delta\|_{L^\infty}\bigl(1+|\log\rho|\bigr).
\]
Taking $\rho=\varepsilon^{1/(1+\beta_{\rm ax})}$ and using
$\varepsilon^{\beta_{\rm ax}}(1+|\log\varepsilon|)
\le C\varepsilon^{\theta_{\rm def}}$ for $0<\varepsilon\le\varepsilon_0$, we obtain
\begin{equation}
\|\nabla v_\Delta\|_{L^\infty(\R^3)} \le C_{\sharp,C_0}\varepsilon^{\theta_{\rm def}}\mathcal H(f_\Psi).
\label{eq:velocity-gradient-deformation-small}
\end{equation}

We choose $r_0>0$ so that
$B_{4r_0}(0)\cap\mathcal K_{\sharp,C_0}=\varnothing$.  Then
\[
|x-y|\ge 2r_0 \ \ \text{ for } \ \ x\in B_{2r_0}(0),\ y\in\mathcal K_{\sharp,C_0},
\]
and the kernel $L$, defined in \eqref{eq:biot-savart-gradient-kernel}, satisfies
\[
\sup_{\substack{x\in B_{2r_0}(0)\\ y\in\mathcal K_{\sharp,C_0}}}
\bigl(|\nabla_yL(x-y)|+|\nabla_x\nabla_yL(x-y)|\bigr)
\le C_{\sharp,C_0}.
\]
For $x\in B_{2r_0}(0)$, differentiating \eqref{eq:fixed-scale-placed-velocity} yields
\begin{equation}
\nabla V_F^\Psi(x) = \int_{\operatorname{supp}\omega_\Psi} L(x-y')\omega_\Psi(y')\,dy' .
\label{eq:deformed-velocity-gradient-integral}
\end{equation}
We then make the change of variables $y'=\widetilde\Psi(y)$ in \eqref{eq:deformed-velocity-gradient-integral} and use the identities
\[
dy'=dy, \qquad \omega_\Psi(\widetilde\Psi(y))=\omega_{\rm id}(y), \qquad \nabla V_F(x) =
\int_{\operatorname{supp}\omega_{\rm id}}L(x-y)\omega_{\rm id}(y)\,dy ,
\]
Subtracting the identity for $\nabla V_F(x)$, we obtain
\[
\nabla v_\Delta(x) = \int_{\operatorname{supp}\omega_{\rm id}} \bigl[L(x-\widetilde\Psi(y))-L(x-y)\bigr]\omega_{\rm id}(y)\,dy .
\]
For $x,x'\in B_{2r_0}(0)$, \eqref{eq:identity-placement-holder-from-deformed} and
$\|\widetilde\Psi-\operatorname{Id}\|_{L^\infty}\le C_{\sharp,C_0}\varepsilon$ imply that
\begin{align*}
\|\omega_{\rm id}\|_{L^1} &\le C_{\sharp,C_0}\mathcal H(f_\Psi),\qquad
|\nabla v_\Delta(x)|\le C_{\sharp,C_0}\varepsilon\,\mathcal H(f_\Psi),\\
|\nabla v_\Delta(x)-\nabla v_\Delta(x')| &\le C_{\sharp,C_0}\varepsilon\,\mathcal H(f_\Psi)|x-x'|.
\end{align*}
Therefore, 
\begin{equation}
\|\nabla v_\Delta\|_{C^{\beta_{\rm ax}}(B_{2r_0}(0))} \le C_{\sharp,C_0}\varepsilon\,\mathcal H(f_\Psi)
\le C_{\sharp,C_0}\varepsilon^{\theta_{\rm def}}\mathcal H(f_\Psi).
\label{eq:local-gradient-deformation-small}
\end{equation}
For $V_G=\BS[\omega_{g_{\rm id}}]$, the support of $\omega_{g_{\rm id}}$ is separated from
$B_{2r_0}(0)$, so the same smooth-kernel estimate in the observation variable gives
\begin{equation}
\|\nabla V_G\|_{L^\infty(\R^3)} +\|\nabla V_G\|_{C^{\beta_{\rm ax}}(B_{2r_0}(0))} \le C_{\sharp,C_0}\mathcal H(g_{\rm id}).
\label{eq:VG-gradient-fixed-scale}
\end{equation}

\runinhead{Step 3: Bounds for $\Pi[V_G,v_\Delta]$ and $\Pi[v_\Delta,v_\Delta]$.}

We next use the properties of the pressure Hessian kernel $K_{zz}$ defined in \eqref{eq:Kzz-kernel}. 
By \eqref{eq:physical-pressure-bilinear-form},
\[
\Pi[V_G,v_\Delta] = \pv\!\int_{\R^3}K_{zz}(y)\, \tr\bigl(\nabla V_G(y)\nabla v_\Delta(y)\bigr)\,dy .
\]
We set $S_{G,\Delta}(y):=\tr\bigl(\nabla V_G(y)\nabla v_\Delta(y)\bigr)$.  Since
$K_{zz}$ has zero spherical mean, \eqref{eq:Kzz-kernel} implies that
\begin{equation}
\pv\!\int_{B_{r_0}(0)}K_{zz}(y)S_{G,\Delta}(y)\,dy =\int_{B_{r_0}(0)}K_{zz}(y)\bigl(S_{G,\Delta}(y)-S_{G,\Delta}(0)\bigr)\,dy .
\label{eq:Kzz-near-origin-cancellation}
\end{equation}
Using \eqref{eq:velocity-gradient-deformation-small},
\eqref{eq:local-gradient-deformation-small}, and \eqref{eq:VG-gradient-fixed-scale}, we obtain
\[
[S_{G,\Delta}]_{C^{\beta_{\rm ax}}(B_{r_0}(0))} \le C_{\sharp,C_0}\varepsilon^{\theta_{\rm def}}\mathcal H(f_\Psi)\mathcal H(g_{\rm id}),
\]
and hence
\[
\big|\pv\!\int_{B_{r_0}(0)}K_{zz}(y)S_{G,\Delta}(y)\,dy\big|
\le C_{\sharp,C_0}\varepsilon^{\theta_{\rm def}}\mathcal H(f_\Psi)\mathcal H(g_{\rm id}) \int_0^{r_0}r^{-3}r^{\beta_{\rm ax}}r^2\,dr .
\]
The last integral is finite because $\beta_{\rm ax}>0$.

For $|y|\ge r_0$, the pressure Hessian kernel $K_{zz}$ defined in \eqref{eq:Kzz-kernel} satisfies
\begin{subequations} 
\begin{equation}
K_{zz}\in L^1_{\rm loc}(\R^3\setminus B_{r_0}(0)), \qquad |K_{zz}(y)|\le C|y|^{-3}.
\label{eq:Kzz-away-origin-bound}
\end{equation}
\begin{equation}
\operatorname{supp}\omega_{g_{\rm id}} \cup \operatorname{supp}\omega_\Delta \subset \mathcal K_{\sharp,C_0}.
\label{eq:pressure-hessian-compact-support}
\end{equation}
\end{subequations} 
We choose $R_K$ so that $\mathcal K_{\sharp,C_0}\subset B_{R_K}(0)$.  On $r_0\le |y|\le 2R_K$,
\eqref{eq:velocity-gradient-deformation-small} and \eqref{eq:VG-gradient-fixed-scale} imply
\[
|\nabla V_G(y)|\le C_{\sharp,C_0}\mathcal H(g_{\rm id})(1+|y|)^{-3},
\qquad |\nabla v_\Delta(y)|\le C_{\sharp,C_0}\varepsilon^{\theta_{\rm def}}\mathcal H(f_\Psi)(1+|y|)^{-3}.
\]
For $|y|\ge 2R_K$, \eqref{eq:pressure-hessian-compact-support} and
\eqref{eq:biot-savart-gradient-kernel} imply
\[
|\nabla V_G(y)| \le C(1+|y|)^{-3}\|\omega_{g_{\rm id}}\|_{L^1}, \qquad |\nabla v_\Delta(y)| \le C(1+|y|)^{-3}\|\omega_\Delta\|_{L^1}.
\]
Using \eqref{eq:identity-placement-holder-from-deformed} and
\eqref{eq:placed-vorticity-difference}, we obtain, for all $|y|\ge r_0$,
\[
|\nabla V_G(y)|\le C_{\sharp,C_0}\mathcal H(g_{\rm id})(1+|y|)^{-3},
\qquad |\nabla v_\Delta(y)|\le C_{\sharp,C_0}\varepsilon^{\theta_{\rm def}}\mathcal H(f_\Psi)(1+|y|)^{-3}.
\]
By \eqref{eq:Kzz-away-origin-bound}, the integral over $\R^3\setminus B_{r_0}(0)$ is bounded by
\[
C_{\sharp,C_0}\varepsilon^{\theta_{\rm def}} \mathcal H(f_\Psi)\mathcal H(g_{\rm id}) \int_{\R^3\setminus B_{r_0}(0)}
|y|^{-3}(1+|y|)^{-6}\,dy .
\]
This proves the mixed estimate in \eqref{eq:normal-form-deformation}.

To prove \eqref{eq:normal-form-deformation-self}, we set
$ S_{\Delta,\Delta}(y):=\tr\bigl(\nabla v_\Delta(y)\nabla v_\Delta(y)\bigr)$. 
\eqref{eq:velocity-gradient-deformation-small} and
\eqref{eq:local-gradient-deformation-small} imply
\begin{equation}
[S_{\Delta,\Delta}]_{C^{\beta_{\rm ax}}(B_{r_0}(0))}\le C_{\sharp,C_0}\varepsilon^{\theta_{\rm def}}\mathcal H(f_\Psi)^2
\label{eq:self-pressure-near-origin-holder}
\end{equation}
and that
\begin{equation}
|S_{\Delta,\Delta}(y)| \le C_{\sharp,C_0}\varepsilon^{2\theta_{\rm def}}\mathcal H(f_\Psi)^2(1+|y|)^{-6} \ \ \text{ for } |y|\ge r_0 .
\label{eq:self-pressure-far-region-decay}
\end{equation}
By \eqref{eq:physical-pressure-bilinear-form},
\[
\Pi[v_\Delta,v_\Delta] = \pv\!\int_{\R^3}K_{zz}(y)S_{\Delta,\Delta}(y)\,dy .
\]
The identity \eqref{eq:Kzz-near-origin-cancellation} with $S_{G,\Delta}$ replaced by
$S_{\Delta,\Delta}$, together with \eqref{eq:self-pressure-near-origin-holder}, implies
\[
\big|\pv\!\int_{B_{r_0}(0)}K_{zz}(y)S_{\Delta,\Delta}(y)\,dy\big| \le C_{\sharp,C_0}\varepsilon^{\theta_{\rm def}}\mathcal H(f_\Psi)^2
\int_0^{r_0}r^{-3}r^{\beta_{\rm ax}}r^2\,dr .
\]
Moreover, \eqref{eq:Kzz-away-origin-bound} and
\eqref{eq:self-pressure-far-region-decay} imply
\[
\big|\int_{\R^3\setminus B_{r_0}(0)}K_{zz}(y)S_{\Delta,\Delta}(y)\,dy\big|
\le C_{\sharp,C_0}\varepsilon^{2\theta_{\rm def}}\mathcal H(f_\Psi)^2 \int_{\R^3\setminus B_{r_0}(0)} |y|^{-3}(1+|y|)^{-6}\,dy.
\]
Both integrals are finite, and $0<\varepsilon\le1$ implies
$\varepsilon^{2\theta_{\rm def}}\le\varepsilon^{\theta_{\rm def}}$.  Hence
\eqref{eq:normal-form-deformation-self} holds.
\end{proof}

\subsection{Axial functions and the Riccati transfer for the cusp velocity}

After the fixed choice of $R_{\tail}$ in Section~\ref{sec:fixed-choice-order},
\eqref{eq:physical-zeta-coeff-tail-range} gives
\[
I_a=[0,\zeta_a]\subset I_{\rm all}(t)
\]
for all small-clock times under consideration.  On this small-clock interval, \eqref{eq:truncated-pressure-coeff} becomes
\[
a_t(\zeta) = \bigl(JA_t(Z_t(\zeta))\bigr)^{1-\alpha} \bigl(1+Z_t(\zeta)^2\bigr)^{-\gamma/2}\mathbf 1_{I_a}(\zeta),
\qquad B_t(Z_t(\zeta))=J^2\zeta,
\]
where $J=J_{\cusp}(t)$.  The cutoff parameter $M_\pressure$ is used in the localization estimates below to
measure the large-slope part of the axis-trace error.

The cusp-flow transported toroidal vorticity defines the axial function $a_t^{\rm phys}$ in
\eqref{eq:full-physical-zeta-profile-def} on $I_{\rm all}(t)$.  The function
$a_t$ in \eqref{eq:truncated-pressure-coeff} is the restriction of $a_t^{\rm phys}$ to
$I_a$, extended by zero off $I_a$.  The transported cusp vorticity also contains the part
$a_t^{\rm phys}\mathbf 1_{I_{\rm all}(t)\setminus I_a}$ outside $I_a$.

The verification of the renormalized axis-trace hypotheses is organized around four sources of discrepancy.
First, the localization cutoff
$\vartheta_\sharp\in C_c^\infty(I_\sharp)$ from \eqref{eq:pressure-localization-cutoff}, with
$\operatorname{supp}\vartheta_\sharp\Subset I_a$, isolates the common part
$\vartheta_\sharp a_t=\vartheta_\sharp a_t^{\rm phys}$.  The complementary $\zeta$-tail, taken
uniformly in time over the stopped small-clock interval $\mathcal I_{\rm stop}$ under consideration, is measured by
\begin{equation}
\mathfrak a_{\zeta}(I_\sharp) := \sup_{t\in\mathcal I_{\rm stop}} \sum_{a\in\{a_t,a_t^{\rm phys}\}} \bigg[ \int_{\R_+}
(1-\vartheta_\sharp(\zeta))^2a(\zeta)^2\zeta^{2\alpha-1}\,d\zeta + \Big( \int_{\R_+}
|1-\vartheta_\sharp(\zeta)|a(\zeta)\zeta^{\alpha-1}\,d\zeta \Big)^2 \bigg],
\label{eq:pressure-zeta-tail}
\end{equation}
which is controlled through the algebraic upper bound \eqref{eq:physical-zeta-profile-envelope}.  Second, the
model vorticity \eqref{eq:vort-slope-restricted} is restricted to slopes $\tau\lesssim M_\pressure$, while the
transported cusp vorticity carries the full angular range in $\tau=R/|Z|$; the large-slope complement
$1-\chi_{M_\pressure}(\tau)$ is measured by the angular-tail integral
\begin{equation}
\mathfrak a_{\rm ang}(M) := 2\int_M^\infty \tfrac{(1+\tau)^{2\alpha}}{1+\tau^2}\,d\tau,
\qquad \mathfrak a_{\rm ang}(M)\to0 \text{ as } M\to\infty,
\label{eq:pressure-angular-tail}
\end{equation}
which is the same angular tail as in \eqref{eq:scaled-angular-tail}.  Third, on the localized label tube
$\mathcal Q_t^\sharp$ from \eqref{eq:localized-label-tube-sec12}, the exact cusp-flow image is
replaced by the normal-form approximation $\Psi_t(\zeta\tau,\zeta)$ from
\eqref{eq:normal-form-approximation-def}; the required displacement and Jacobian bounds are
\eqref{eq:additive-normal-form-bound} and \eqref{eq:normal-form-approximation-bound}.  Fourth, labels
$Y$ outside the bounded core $D_{\core}$ in \eqref{eq:core-tail-domains} are controlled by the algebraic decay
$\rho^{-\gamma}$ in \eqref{eq:gamma-stand}; after $R_{\tail}$ is fixed, these labels contribute the lower-order
$J^{1-3\alpha}$ contribution to the axis-trace error.

The choice order in Section~\ref{sec:fixed-choice-order} first fixes
$I_\sharp,\vartheta_\sharp$ in \eqref{eq:pressure-localization-cutoff}--\eqref{eq:pressure-localization-intervals},
then fixes $R_{\tail}$ through \eqref{eq:core-tail-domains}--\eqref{eq:smooth-velocity-def}, then fixes
$M_\pressure$ and $C_0$ in \eqref{eq:pressure-C0}, and finally decreases
$\mathfrak J_{\pressure}$ in \eqref{eq:small-clock-threshold-order}.  With these choices, the locally defined
axis-trace error parameter $\varepsilon_{\rm ax}$ from \eqref{eq:epsilon-ax} satisfies
\eqref{eq:epsilon-ax-sec12} and is the small parameter used in
Proposition~\ref{prop:euler-generated-profile-riccati}.

\begin{lemma}[Renormalized axis-trace hypotheses for the Euler-generated axial function]
\label{lem:renormalized-axis-trace-hypotheses}
Let $a_t$ be the zero-extended Euler-generated axial function in \eqref{eq:truncated-pressure-coeff}.  Assume the axis-geometry bounds
\eqref{eq:entry-axis-bounds-statement}, the pressure-interval coverage \eqref{eq:axis-control-covers-pressure-interval}, the axis volume identity
\eqref{eq:axis-qb-volume-mon}, the monotone axial-stretching bounds \eqref{eq:monotone-axial-two-sided}--\eqref{eq:monotone-axial-fractional-bootstrap}, 
the axial-amplitude bounds \eqref{eq:current-axis-extension-size}, the normal-form bounds
\eqref{eq:localized-normal-form-large-bootstrap}--\eqref{eq:localized-normal-form-map-large-bootstrap}, and the
axis evolution equations of Lemma~\ref{lem:axis-profile-evolution} on $I_a$.
Then, with $J=J_{\cusp}(t)$, after fixing $I_\sharp$, $M_\pressure$, and $R_{\tail}$ as in
Section~\ref{sec:fixed-choice-order} and then decreasing $\mathfrak J_{\pressure}$ in
\eqref{eq:small-clock-threshold-order}, the number
$\varepsilon_{\rm ax}$ defined in \eqref{eq:epsilon-ax} satisfies
\begin{equation}
\varepsilon_{\rm ax}\le C\Big(\mathfrak a_\zeta(I_\sharp)^{1/2}+\mathfrak a_{\rm ang}(M_\pressure)^{1/2} +J^{\kappa_{\rm def}}+J^{1-3\alpha}\Big)
\label{eq:epsilon-ax-sec12}
\end{equation}
and is as small as needed.  Moreover,
\eqref{eq:axis-trace-approximation}--\eqref{eq:axis-strain-clock-comparison} and
\eqref{eq:cusp-principal-riccati-with-error} hold for this axial function $a_t$.
\end{lemma}

\begin{proof}[Proof of Lemma~\ref{lem:renormalized-axis-trace-hypotheses}]
The monotonicity of $a_t$ follows from \eqref{eq:axis-qb-volume-mon} and the monotone axial-stretching bounds,
as in Lemma~\ref{lem:exact-euler-generated-derivative}.  We next verify the hypotheses of
Proposition~\ref{prop:euler-generated-profile-riccati}.

The axis-trace estimates \eqref{eq:axis-trace-U-approximation}--\eqref{eq:axis-trace-W-approximation} follow
from the fixed-variable Biot--Savart representation after the change of variables $y=J^2(R,Z)$.  On the
localized set $\zeta\in I_\sharp$, $|\tau|\le2M_\pressure$, the transported vorticity is
\[
-\Gamma J^{3\alpha-1}\operatorname{sgn}(Z)a_t(|Z|)R^\alpha e_\theta
\]
up to the normal-form displacement error \eqref{eq:normal-form-approximation-bound}.  The complement is the sum of the
$\zeta$-tail measured by \eqref{eq:pressure-zeta-tail}, the large-slope tail measured by
\eqref{eq:pressure-angular-tail}, and the algebraic far-field tail from \eqref{eq:gamma-stand}.  Applying the
axis Biot--Savart kernels for the velocity trace and its $\zeta$-derivative gives exactly
\eqref{eq:axis-trace-U-approximation}--\eqref{eq:axis-trace-W-approximation}, with the value of
$\varepsilon_{\rm ax}$ displayed above.  The division by $\zeta$ in
\eqref{eq:axis-trace-U-approximation} is harmless because both axis velocities vanish at the origin and the
same estimate for the derivative holds on $[0,\zeta_a]$.

In \eqref{eq:axis-current-zeta-velocity}, the $V_{\err}$ contribution on the symmetry axis is
\[
J^{-2}(V_{\err})_z(0,J^2\zeta,t).
\]
It enters $\mathcal R_t^a$ through the term $-\mathsf R_t^\zeta\p_\zeta a_t$ in
\eqref{eq:renormalized-a-remainder-def}.  The evolution equation \eqref{eq:renormalized-profile-equation}
is obtained from \eqref{eq:zeta-current-transport},
\eqref{eq:axis-current-zeta-velocity}, and \eqref{eq:axis-profile-q-system}, as shown in
\eqref{eq:renormalized-Z-label-current}--\eqref{eq:renormalized-a-remainder-def}.  The algebraic weight
$(1+Z_t^2)^{-\gamma/2}$ gives no additional term because $\mathcal T_tZ_t=0$ in
\eqref{eq:renormalized-Z-label-current}.  The weighted bound \eqref{eq:renormalized-profile-remainder}
follows from the axis-error traces in \eqref{eq:axis-profile-error-holder}.

The undifferentiated strain comparison \eqref{eq:axis-strain-undifferentiated-comparison} is the derivative
at $\zeta=0$ of the axis-trace approximation, multiplied by the scale
$\Gamma J^{3\alpha-1}$.  The clock law \eqref{eq:axis-clock-law-for-criterion} is
\eqref{eq:cusp-clock-identity-in-Jdot-proof}.  Differentiating the same fixed-variable axis-strain
representation and using \eqref{eq:renormalized-profile-equation} gives
\eqref{eq:direct-strain-derivative-axis}; the two terms displayed there are the derivative of
$\Gamma J^{3\alpha-1}W_\infty[a_t]$, and the remaining terms are bounded by
\eqref{eq:axis-trace-approximation} and \eqref{eq:renormalized-profile-remainder}.

Finally, differentiating the Euler equation for the cusp-coordinate velocity at the stagnation point gives the
principal Riccati identity.  The part driven by $m(t)U_{\cusp}$ gives
$m(t)(-\tfrac12\mathcal W_{\cusp}^2-\Pi_{\cusp})$; the terms containing $V_{\err}$ are controlled by
\eqref{eq:axis-profile-error-holder} and the strain lower bound \eqref{eq:Wcusp-scaling}, giving the error in
\eqref{eq:cusp-principal-riccati-with-error}.  This proves all hypotheses required by
Proposition~\ref{prop:euler-generated-profile-riccati}.
\end{proof}

\begin{lemma}[Cusp-flow Riccati pressure Hessian bound]
\label{lem:transported-cusp-pressure-win}
Assume the axis-geometry bounds \eqref{eq:entry-axis-bounds-statement}, the cusp-clock bootstrap
\eqref{eq:localized-clock-bootstrap}, the pressure-interval coverage \eqref{eq:axis-control-covers-pressure-interval}, the axis volume identity
\eqref{eq:axis-qb-volume-mon}, the monotone axial-stretching bounds \eqref{eq:monotone-axial-two-sided}--\eqref{eq:monotone-axial-fractional-bootstrap}, 
the axial-amplitude bounds \eqref{eq:current-axis-extension-size}, the normal-form bounds
\eqref{eq:localized-normal-form-large-bootstrap} and \eqref{eq:localized-normal-form-map-large-bootstrap}, and the
axis evolution equations of Lemma~\ref{lem:axis-profile-evolution} on $I_a$.
Then there are constants
\[
\mathfrak J_{\pressure} \in(0,\min\{\mathfrak J_{\mathrm{velocity}},\mathfrak J_{\mathrm{axis}}\}], \qquad q_{\rm tr}\in(0,\upbeta),
\]
independent of $t$ and $J_{\cusp}(t)$, such that if
\[
J:=J_{\cusp}(t)\le \mathfrak J_{\pressure},
\]
then
\begin{equation}
\Pi_{\cusp}(t) \ge -q_{\rm tr}\,\tfrac12\,\mathcal W_{\cusp}(t)^2.
\label{eq:Picusp-riccati}
\end{equation}
\end{lemma}

\begin{remark}
The normal-form bootstrap bounds
\eqref{eq:localized-normal-form-large-bootstrap} and
\eqref{eq:localized-normal-form-map-large-bootstrap}
in the hypothesis list are proved in Lemmas~\ref{lem:nonlinear-radial-flatness},
\ref{lem:late-axis-normal-form-cusp}, and
\ref{lem:late-axis-normal-form-map-cusp}, and supply both the cusp-map normal form and the
approximation map $\Psi_t$.
\end{remark}

\subsection{Auxiliary localized pressure estimates}

We use the next localized kernel estimate in the proof of Lemma~\ref{lem:transported-cusp-pressure-win}.  It is
a Calder\'on--Zygmund estimate on the fixed set $\mathcal R_{\sharp,C_0}$ from
\eqref{eq:fixed-scale-deformation-set}, and this set has positive distance from the stagnation point.  Hence
the constants in \eqref{eq:fixed-scale-CZ-localized} and \eqref{eq:fixed-scale-pressure-bilinear} depend on
$I_\sharp$ and $C_0$ but not on the cusp clock $J$.

We state this estimate once for reuse.  Let $I_\sharp\Subset(0,\infty)$ and $C_0<\infty$ be fixed.  For an
axisymmetric scalar function $F$ supported in $\mathcal R_{\sharp,C_0}$ with
$Fe_\theta\in C^{\beta_{\rm ax}}(\R^3)$, set
\[
\mathcal H(F):=\|Fe_\theta\|_{L^\infty}+[Fe_\theta]_{C^{\beta_{\rm ax}}}, \qquad V_F:=\BS[Fe_\theta].
\]

\begin{lemma}[Localized Calder\'on--Zygmund and pressure Hessian bounds]
\label{lem:fixed-scale-CZ-pressure}
For every axisymmetric $F$ supported in $\mathcal R_{\sharp,C_0}$ with
$Fe_\theta\in C^{\beta_{\rm ax}}$,
\begin{equation}
\|\nabla V_F\|_{L^\infty(\R^3)}\le C_{\sharp,C_0,\alpha}\mathcal H(F),
\label{eq:fixed-scale-CZ-localized}
\end{equation}
and, for any second such function $G$ supported in $\mathcal R_{\sharp,C_0}$,
\begin{equation}
\left| \pv\!\int_{\R^3}K_{zz}(y)\, \tr\bigl(\nabla V_F(y)\nabla V_G(y)\bigr)\,dy \right|
\le C_{\sharp,C_0,\alpha}\mathcal H(F)\mathcal H(G).
\label{eq:fixed-scale-pressure-bilinear}
\end{equation}
\end{lemma}

\begin{proof}[Proof of Lemma~\ref{lem:fixed-scale-CZ-pressure}]
We set $\omega_F(y):=F(y)e_\theta(y)$.  The kernel of $\nabla\BS$ is a homogeneous Calder\'on--Zygmund kernel of
degree $-3$.  For $x$ with $\operatorname{dist}(x,\mathcal R_{\sharp,C_0})\le1$, the mean-zero property on
spheres allows us to write the principal-value part as
\[
\pv\!\int_{|x-y|\le1}\nabla K(x-y)\bigl(\omega_F(y)-\omega_F(x)\bigr)\,dy +\int_{|x-y|>1}\nabla K(x-y)\omega_F(y)\,dy
+\mathbf C\,\omega_F(x),
\]
where $\mathbf C$ is the local Calder\'on--Zygmund matrix.  The first integral is bounded by
\[
C[\omega_F]_{C^{\beta_{\rm ax}}}\int_0^1 r^{-3}r^{\beta_{\rm ax}}r^2\,dr\le C_{\beta_{\rm ax}}\mathcal H(F),
\]
and the second by $C_{\sharp,C_0}\|\omega_F\|_{L^\infty}$ because the support of $\omega_F$ is contained in
the solid of revolution of the fixed compact set $\mathcal R_{\sharp,C_0}$.  If instead
$\operatorname{dist}(x,\mathcal R_{\sharp,C_0})>1$, no principal value is present and
\[
|\nabla V_F(x)| \le C\|\omega_F\|_{L^\infty}\int|x-y|^{-3}\,dy \le C_{\sharp,C_0}\mathcal H(F),
\]
where the integral is over the compact support of $\omega_F$.  This proves
\eqref{eq:fixed-scale-CZ-localized}.

For \eqref{eq:fixed-scale-pressure-bilinear}, \eqref{eq:fixed-scale-CZ-localized} implies
$\|\nabla V_F\|_{L^\infty}+\|\nabla V_G\|_{L^\infty}\le C_{\sharp,C_0,\alpha}(\mathcal H(F)+\mathcal H(G))$.
Since $\mathcal R_{\sharp,C_0}$ is separated from the origin, $\nabla V_F$ and $\nabla V_G$ are
$C^{\beta_{\rm ax}}$ on the ball $B(0,c_{\sharp,C_0})$ with the same bound, and so
$S_{F,G}:=\tr(\nabla V_F\nabla V_G)$ satisfies
\[
\|S_{F,G}\|_{L^\infty} +[S_{F,G}]_{C^{\beta_{\rm ax}}(B(0,c_{\sharp,C_0}))} \le C_{\sharp,C_0,\alpha}\mathcal H(F)\mathcal H(G).
\]
The mean-zero property of $K_{zz}$ on spheres then yields
\[
\big|\pv\!\int_{|y|\le c_{\sharp,C_0}}K_{zz}(y)S_{F,G}(y)\,dy\big| =
\big|\int_{|y|\le c_{\sharp,C_0}}K_{zz}(y)\bigl(S_{F,G}(y)-S_{F,G}(0)\bigr)\,dy\big| \le C_{\sharp,C_0,\alpha}\mathcal H(F)\mathcal H(G).
\]
On $|y|\ge c_{\sharp,C_0}$ the kernel is nonsingular on bounded sets, and for large $|y|$ the gradients
$\nabla V_F,\nabla V_G$ decay like $|y|^{-3}$; the remaining part of the integral is bounded by the same
right-hand side.
\end{proof}

The localized estimate \eqref{eq:fixed-scale-pressure-bilinear} is invariant under the change of variables
$y=J_{\cusp}(t)^2Y$ used in \eqref{eq:physical-scaled-vorticity}.  If $F$ and $G$ are toroidal scalars in the
variables $Y=(R,Z)$, supported in the symmetric localization associated with
$\mathcal R_{\sharp,C_0}$ from \eqref{eq:fixed-scale-deformation-set}, and if $V_{\rm sc,J}^F$ and
$V_{\rm sc,J}^G$ are the Biot--Savart velocities generated by the corresponding physical vorticities in
\eqref{eq:physical-scaled-vorticity}, then
\begin{equation}
\big|\Pi[V_{\rm sc,J}^F,V_{\rm sc,J}^G](t)\big|
\le C_{\sharp,C_0,\alpha}\Gamma^2J^{6\alpha-2}\mathcal H(F)\mathcal H(G).
\label{eq:scaled-pressure-hessian-error}
\end{equation}
For the symmetric localization one applies \eqref{eq:fixed-scale-pressure-bilinear} separately to the
upper--upper and lower--lower pairings; the upper--lower pairings have supports separated by a positive
distance in the variables $Y=(R,Z)$ and are controlled by the same nonsingular kernel bound.  Thus whenever one of the two
normalized localized vorticities is smaller by a scalar $\delta$, \eqref{eq:scaled-pressure-hessian-error}
contributes the corresponding pressure error
$O(\delta\,\Gamma^2J^{6\alpha-2})$.

For an axisymmetric no-swirl field $f$, we use the cone-local norm
\begin{equation}
\|\nabla f\|_{\mathcal C_*,\alpha} := \sup_{B(x,2c_*|x|)\subset\mathcal C_*} \left( \|\nabla f\|_{L^\infty(B(x,c_*|x|))}
+|x|^\alpha[\nabla f]_{C^\alpha(B(x,c_*|x|))} \right).
\label{eq:cone-local-norm}
\end{equation}

The next lemma is the pressure-Hessian estimate used below for the non-geometric terms in
\eqref{eq:Pi0-decomp}.  It uses the cancellation of the pressure Hessian kernel on spherical shells and
ordinary H\"older difference quotients for the pressure source $\tr(\nabla u\nabla v)$.

\begin{lemma}[Dyadic H\"older estimate for the pressure bilinear form]
\label{lem:cone-local-pressure-bilinear}
Let $u$ and $v$ be axisymmetric no-swirl velocity fields at a fixed time.  For $j\in\mathbb Z$, set
\[
\mathscr R_j:=\{y\in\R^3:\ 2^j\le |y|\le2^{j+1}\}.
\]
We also define the enlarged shell
\[
\mathscr R_j^*:=\{y\in\R^3:\ 2^{j-1}\le |y|\le2^{j+2}\}.
\]
Assume that
\begin{equation}
\mathcal N_\alpha[u,v]:=
\sum_{j\in\mathbb Z}2^{j\alpha}
\Big(
\|\nabla u\|_{L^\infty(\mathscr R_j^*)}[\nabla v]_{C^\alpha(\mathscr R_j^*)}
+\|\nabla v\|_{L^\infty(\mathscr R_j^*)}[\nabla u]_{C^\alpha(\mathscr R_j^*)}
\Big)<\infty .
\label{eq:dyadic-holder-pressure-norm}
\end{equation}
Then
\begin{equation}
|\Pi[u,v]|\le C_\alpha \mathcal N_\alpha[u,v].
\label{eq:pressure-bilinear-cone-estimate}
\end{equation}
\end{lemma}

\begin{proof}[Proof of Lemma~\ref{lem:cone-local-pressure-bilinear}]
We first estimate a finite truncation
\[
\int_{\{2^{-N}\le |y|\le2^M\}}K_{zz}(y)S(y)\,dy,\qquad  S(y):=\tr(\nabla u(y)\nabla v(y)),
\]
and then let $N,M\to\infty$.  For each shell $\mathscr R_j$, we choose one point $y_j\in\mathscr R_j$.  The
pressure kernel has zero spherical mean: $\int_{\mathbb S^2}K_{zz}(\theta)\,d\theta=0$.
Since $\mathscr R_j$ is a full spherical annulus, the kernel integrates to zero on the shell:
\[
\int_{\mathscr R_j}K_{zz}(y)\,dy=0.
\]
Therefore, the contribution of the constant value $S(y_j)$ vanishes:
\[
\int_{\mathscr R_j}K_{zz}(y)S(y_j)\,dy =S(y_j)\int_{\mathscr R_j}K_{zz}(y)\,dy=0, 
\]
and thus, 
\[
\int_{\mathscr R_j}K_{zz}(y)S(y)\,dy =\int_{\mathscr R_j}K_{zz}(y)\bigl(S(y)-S(y_j)\bigr)\,dy .
\]
For $y\in\mathscr R_j$, $|y-y_j|\le C2^j$, so
\[
|S(y)-S(y_j)|\le C2^{j\alpha}[S]_{C^\alpha(\mathscr R_j^*)}.
\]
Since $|K_{zz}(y)|\le C|y|^{-3}$ and
\[
\int_{\mathscr R_j}|y|^{-3}\,dy\le C,
\]
we obtain
\[
\Big|\int_{\mathscr R_j}K_{zz}(y)S(y)\,dy\Big| \le C_\alpha2^{j\alpha}[S]_{C^\alpha(\mathscr R_j^*)}.
\]
The product estimate
\[
[S]_{C^\alpha(\mathscr R_j^*)}\le C\Big( \|\nabla u\|_{L^\infty(\mathscr R_j^*)}[\nabla v]_{C^\alpha(\mathscr R_j^*)} 
+\|\nabla v\|_{L^\infty(\mathscr R_j^*)}[\nabla u]_{C^\alpha(\mathscr R_j^*)} \Big)
\]
proves \eqref{eq:pressure-bilinear-cone-estimate} after summing over $j$.  The finiteness of
\eqref{eq:dyadic-holder-pressure-norm} also shows that the principal value exists and is independent of the
truncation.
\end{proof}

\subsection{Proof of the Riccati pressure Hessian estimate for the cusp velocity}

\begin{proof}[Proof of Lemma~\ref{lem:transported-cusp-pressure-win}]
We take
\[
q_{\rm tr}:=q_\alpha .
\]
By the fixed choice in \eqref{eq:beta-def} and \eqref{eq:pressure-subcritical-riccati}, this constant satisfies
$q_{\rm tr}<\upbeta$.  We then choose $\mathfrak J_{\pressure}$ no larger than the small-clock threshold in
Lemma~\ref{lem:renormalized-axis-trace-hypotheses}, and small enough that the corresponding
$\varepsilon_{\rm ax}$ satisfies the smallness condition in Proposition~\ref{prop:euler-generated-profile-riccati}.
The choices of $I_\sharp$, $M_\pressure$, and $R_{\tail}$ are those described above and in
Section~\ref{sec:fixed-choice-order}; they make the four contributions on the right-hand side of
\eqref{eq:epsilon-ax-sec12} as small as required.

Lemma~\ref{lem:renormalized-axis-trace-hypotheses} verifies
\eqref{eq:axis-trace-approximation}--\eqref{eq:axis-strain-clock-comparison} and \eqref{eq:cusp-principal-riccati-with-error} for the Euler-generated axial function $a_t$ in
\eqref{eq:truncated-pressure-coeff}.  The monotone axial-stretching assumptions show that $a_t$ is nonnegative and nonincreasing.  Hence all hypotheses
of Proposition~\ref{prop:euler-generated-profile-riccati} hold, and that proposition yields
\[
\Pi_{\cusp}(t)\ge -q_\alpha\,\tfrac12\,\mathcal W_{\cusp}(t)^2.
\]
Since $q_{\rm tr}=q_\alpha$, this is exactly \eqref{eq:Picusp-riccati}.
\end{proof}

\subsection{Transported bounds and lower-order pressure terms}

In Sections~\ref{sec:target-profile-typeI-completion}--\ref{sec:proof-main}, the field estimates of
Lemma~\ref{lem:transported-cusp-field-bounds} (sign and size of $\mathcal W_{\cusp}$, $L^\infty$ and
local H\"older bounds for $\nabla U_{\cusp}$ on $\mathcal C_*$, and the local comparison of $U_{\cusp}$
with the hyperbolic stagnation field $u_{\hyp}[\mathcal W_{\cusp}]$ on the balls
$|x|\le C_{\rm sc}J_{\cusp}(t)^2$) and the one-sided pressure Hessian bound of
Lemma~\ref{lem:transported-cusp-pressure-win} are invoked together under a single small-clock
threshold $\mathfrak J_{\transport}$.  The next lemma states both sets of estimates in this combined form,
under the hypotheses of Lemma~\ref{lem:transported-cusp-pressure-win}.

\begin{lemma}[Cusp-flow transported velocity, strain, and pressure Hessian bounds]
\label{lem:transported-cusp-estimates}
Assume the hypotheses of Lemma~\ref{lem:transported-cusp-pressure-win}.  Then there are constants
\[
\mathfrak J_{\transport}\in(0,\mathfrak J_{\mathrm{velocity}}], \qquad C\ge1,
\]
depending only on $\alpha,\gamma,\sigma_{\inn},\sigma_*$ and on the bootstrap constants in the hypotheses,
such that, whenever $J:=J_{\cusp}(t)\le\mathfrak J_{\transport}$,
\begin{subequations}
\begin{align}
\mathcal W_{\cusp}(t)&<0,
\\
C^{-1}\Gamma J^{3\alpha-1}
\le |\mathcal W_{\cusp}(t)|
&\le C\Gamma J^{3\alpha-1},
\\
\|\nabla U_{\cusp}(\cdot,t)\|_{L^\infty(\mathcal C_*)}
&\le C\,\Gamma J^{3\alpha-1},
\\
[\nabla U_{\cusp}(\cdot,t)]_{C^\alpha(B(x,c_*|x|))}
&\le C\,\Gamma J^{3\alpha-1}|x|^{-\alpha}
\ \ \text{ for } \ \    B(x,2c_*|x|)\subset\mathcal C_*,
\\
\Pi_{\cusp}(t)
&\ge
-q_{\rm tr}\,\tfrac12\,\mathcal W_{\cusp}(t)^2 .
\label{eq:transport-pressure-hessian-bound}
\end{align}
\end{subequations}
where $q_{\rm tr}\in(0,\upbeta)$ is the pressure-transfer constant from
Lemma~\ref{lem:transported-cusp-pressure-win}.
Moreover, for each prescribed $C_{\rm sc}<\infty$ there are constants
\[
\mathfrak J_{\transport}(C_{\rm sc})\in(0,\mathfrak J_{\transport}], \qquad C_{\rm tr}(C_{\rm sc})<\infty,
\]
with additional dependence only on $C_{\rm sc}$, such that if
$J\le\mathfrak J_{\transport}(C_{\rm sc})$, then
\begin{subequations}
\begin{align}
\big|\nabla U_{\cusp}(x,t)
-\nabla u_{\hyp}[\mathcal W_{\cusp}(t)]\big|
&\le C_{\rm tr}(C_{\rm sc})\,\Gamma J^{3\alpha-1},
\\
\tfrac{\big|(U_{\cusp}-u_{\hyp}[\mathcal W_{\cusp}])_r(x,t)\big|}
{r(x)}
&\le C_{\rm tr}(C_{\rm sc})\,\Gamma J^{3\alpha-1},
\label{eq:transport-local-radial-bound}
\end{align}
\end{subequations}
for all $x\in\mathcal C_*$ with $|x|\le C_{\rm sc}J^2$; the quotient in
\eqref{eq:transport-local-radial-bound} is interpreted by its continuous axis value at $r(x)=0$.
\end{lemma}

\begin{remark}
The pressure Hessian estimate \eqref{eq:transport-pressure-hessian-bound} is not a new estimate.  It
is \eqref{eq:Picusp-riccati} applied under the smaller threshold
$\mathfrak J_{\transport}$ and stated together with the transported field estimates.
\end{remark}

\begin{proof}[Proof of Lemma~\ref{lem:transported-cusp-estimates}]
For the estimates that do not involve $C_{\rm sc}$, we define
\[
\mathfrak J_{\transport}:= \min\{\mathfrak J_{\mathrm{velocity}},\mathfrak J_{\pressure},\mathfrak J_{\mathrm{axis}}\}.
\]
The extra minimum with $\mathfrak J_{\mathrm{axis}}$ ensures that the axis-geometry hypotheses used by
Lemma~\ref{lem:transported-cusp-field-bounds} are available at the time $t$.
Let $c_W,C_W$ be the constants in \eqref{eq:Wcusp-scaling}, let $C_U$ be the constant in
\eqref{eq:Ucusp-grad-Linf}--\eqref{eq:Ucusp-grad-scale-local}.  We fix the constant $C$ in
the statement of the present lemma so that
\[
C\ge \max\{c_W^{-1},C_W,C_U\}.
\]
Let $t$ satisfy $J:=J_{\cusp}(t)\le \mathfrak J_{\transport}$.  Since
$J\le \mathfrak J_{\mathrm{velocity}}$, Lemma~\ref{lem:transported-cusp-field-bounds} applies, and
its strain estimate \eqref{eq:Wcusp-scaling} yields
\[
\mathcal W_{\cusp}(t)<0, \qquad C^{-1}\Gamma J^{3\alpha-1} \le |\mathcal W_{\cusp}(t)| \le C\Gamma J^{3\alpha-1}.
\]
By \eqref{eq:Ucusp-grad-Linf}--\eqref{eq:Ucusp-grad-scale-local},
\[
\|\nabla U_{\cusp}(\cdot,t)\|_{L^\infty(\mathcal C_*)} \le C\Gamma J^{3\alpha-1},
\]
and, whenever $B(x,2c_*|x|)\subset\mathcal C_*$,
\[
[\nabla U_{\cusp}(\cdot,t)]_{C^\alpha(B(x,c_*|x|))} \le C\Gamma J^{3\alpha-1}|x|^{-\alpha}.
\]
Since also $J\le \mathfrak J_{\pressure}$, Lemma~\ref{lem:transported-cusp-pressure-win} applies, and
\eqref{eq:Picusp-riccati} gives
\[
\Pi_{\cusp}(t) \ge -q_{\rm tr}\,\tfrac12\,\mathcal W_{\cusp}(t)^2 .
\]

For the localized assertions, we fix $C_{\rm sc}<\infty$ and set
\[
\mathfrak J_{\transport}(C_{\rm sc}) := \min\{\mathfrak J_{\transport},\mathfrak J_{\mathrm{local}}(C_{\rm sc})\},
\]
with $\mathfrak J_{\mathrm{local}}(C_{\rm sc})$ the localized threshold from
Lemma~\ref{lem:transported-cusp-field-bounds}; we choose
$C_{\rm tr}(C_{\rm sc})\ge C_{\rm loc}(C_{\rm sc})$ for the corresponding localized constant.  If in
addition $J\le \mathfrak J_{\transport}(C_{\rm sc})$, then
\eqref{eq:Ucusp-gradient-defect} implies, for every $x\in\mathcal C_*$ with $|x|\le C_{\rm sc}J^2$,
\[
\big|\nabla U_{\cusp}(x,t) -\nabla u_{\hyp}[\mathcal W_{\cusp}(t)]\big| \le C_{\rm tr}(C_{\rm sc})\Gamma J^{3\alpha-1},
\]
while \eqref{eq:Ucusp-radial-defect} yields
\[
\tfrac{\big|(U_{\cusp}-u_{\hyp}[\mathcal W_{\cusp}])_r(x,t)\big|} {r(x)} \le C_{\rm tr}(C_{\rm sc})\Gamma J^{3\alpha-1}.
\]
These are precisely the assertions of the lemma.
\end{proof}

The Riccati blowup argument of Section~\ref{sec:target-profile-typeI-completion} integrates the
stagnation-point identity
\[
\p_t\rW_0(t) = -\tfrac12\rW_0(t)^2-\Pi_0(t),
\]
in which $\Pi_0(t)$ is the full pressure Hessian at the stagnation point.  Using the cusp/smooth/error decomposition of $u$, \eqref{eq:Pi0-decomp} splits this Hessian as
\[
\Pi_0(t)=m(t)^2\Pi_{\cusp}(t)+\Pi_{\rm geom}(t)+\Pi_{\mixed}(t)+\Pi_{\smooth}(t)+\Pi_{\err}(t).
\]
The principal term $m^2\Pi_{\cusp}$ now satisfies the one-sided Riccati bound
\eqref{eq:Picusp-riccati}.  The geometric remainder $\Pi_{\rm geom}$, encoding the smooth-flow
pushforward defect, is handled in Lemma~\ref{lem:smooth-pressure-defect}.  The next lemma supplies the third
and final ingredient: lower-order bounds for the three non-geometric remainders
$\Pi_{\mixed},\Pi_{\smooth},\Pi_{\err}$, each of which contains $u_{\smooth}$ or $u_{\err}$ and is therefore
small compared with $\rW_{\cusp}^2$ after the small-clock threshold and fixed cutoffs are chosen.

The same lemma simultaneously closes the cusp-error size bootstrap (BA8): the bound
\eqref{eq:localized-cusp-error-large-bootstrap} on the cusp-error functional
$\mathfrak E_{\err}(t)$ in the size bootstrap assumptions $\mathcal B_{\rm size}$
\eqref{eq:simultaneous-continuation-assumptions} is replaced by a sharp fixed bound
$\mathfrak E_{\err}(t)\le C$, strictly smaller than the bootstrap constant $E_*$.  Earlier applications of
the bootstrap, such as the radial-derivative estimate \eqref{eq:radial-flatness-source} and the
axis system \eqref{eq:axis-current-zeta-velocity}--\eqref{eq:axis-profile-error-holder}, used only the
weak large bound; the sharper estimate is proved only after we have the transported cusp-field bounds
\eqref{eq:Ucusp-grad-Linf}--\eqref{eq:Ucusp-radial-defect}, the smooth-flow estimates
\eqref{eq:usmooth-small-C2a}--\eqref{eq:smooth-flow-second-gradient-small}, and the
inner-core scaling estimate \eqref{eq:exact-current-axis-core-small}; these are the
estimates invoked in the proof below.

Throughout the lemma we use the small-clock notation of
Section~\ref{sec:geom-flow-decomp},  the flow maps $\phi_{\smooth},\phi_{\cusp}$, the velocities
$u_{\smooth},u_{\core},u_{\cusp},u_{\err}$, the cusp-flow generated velocities $V_{\cusp},U_{\cusp},V_{\err}$,
the cusp clock $J_{\cusp}$ and modulation $m$, and the pressure Hessian terms
$\Pi_{\cusp},\Pi_{\rm geom},\Pi_{\mixed},\Pi_{\smooth},\Pi_{\err}$.  We also use the cone-local norm
defined in \eqref{eq:cone-local-norm}.
The exact identities $u=u_{\smooth}+u_{\core}$ and
$V_{\cusp}=(\phi_{\smooth}^{-1})_*u_{\core}=m\,U_{\cusp}+V_{\err}$ from
\eqref{eq:smooth-velocity-def}, \eqref{eq:ucore-def}, and \eqref{eq:Verr-def} hold for
every time on the interval of existence; the bounds below use these identities throughout.

\begin{lemma}[Cusp-error and lower-order pressure Hessian estimates]
\label{lem:tail-bound}
Assume the axis-geometry bound \eqref{eq:entry-axis-bounds-statement}, the cusp-clock rate bound \eqref{eq:localized-clock-bootstrap}, 
the inner-core scaling estimate \eqref{eq:exact-current-axis-core-small}, the transported cusp-field estimates
\eqref{eq:Wcusp-scaling}--\eqref{eq:Ucusp-radial-defect}, and the scalar-modulation bootstrap
\eqref{eq:localized-modulation-m-bootstrap}.  There exist constants
\begin{equation}
\mathfrak J_{\mathrm{tail}}\in(0,1],\qquad C_T<\infty,\qquad C<\infty,
\label{eq:tail-clock-threshold-def}
\end{equation}
depending only on $\alpha,\gamma,\sigma_{\inn},\sigma_*$ and on the constants in the hypotheses, such that
the entry time
$t_{\mathrm{tail}}:=\inf\{t:\ J_{\cusp}(t)\le\mathfrak J_{\mathrm{tail}}\}$ satisfies
\begin{equation}
T-t_{\mathrm{tail}}\le C_T\Gamma^{-1},
\label{eq:tail-small-clock-time-horizon}
\end{equation}
and the following bounds hold whenever $J_{\cusp}(t)\le\mathfrak J_{\mathrm{tail}}$:
\begin{enumerate}
\item the far-field velocity $u_{\smooth}$ from \eqref{eq:smooth-velocity-def} is bounded at order $\Gamma$,
\[
\|u_{\smooth}(\cdot,t)\|_{L^\infty(B_1)} +\|\nabla u_{\smooth}(\cdot,t)\|_{L^\infty(B_1)} +[\nabla u_{\smooth}(\cdot,t)]_{C^\alpha(B_1)}
\le C\Gamma ,
\]
and
\begin{equation}
\|\nabla u_{\smooth}(\cdot,t)\|_{\mathcal C_*,\alpha} \le C\Gamma ;
\label{eq:tail-usmooth-cone-local}
\end{equation}
\item for every $Y\in D_{\inn}^{\cusp}(t)$ with $\omega_{\theta,0}(Y)\ne0$, the cusp-flow image point
$x_*:=\phi_{\cusp}(Y,t)$ lies in $\mathcal C_{\inn}$, and the cusp-error velocity satisfies
\begin{subequations}
\begin{align}
|V_{\err}(x_*,t)|
&\le C\Gamma,
\label{eq:tail-error-bound}\\
\tfrac{J_{\cusp}(t)}{R(Y)^{1+\alpha}}\,|V_{\err}(x_*,t)|
&\le C\Gamma\,J_{\cusp}(t)^{9\alpha-1}+C\Gamma J_{\cusp}(t)
\quad\bigl(R(Y)\ge J_{\cusp}(t)^{3/\alpha}\bigr);
\label{eq:tail-error-bound-weighted}
\end{align}
\end{subequations}
\item the cusp-error gradient satisfies
\[
\|\nabla V_{\err}(\cdot,t)\|_{\mathcal C_*,\alpha} \le C\Gamma\bigl(J_{\cusp}(t)^{9\alpha-1}+1\bigr),
\]
and, for every fixed compact interval $I\subset[0,\infty)$ on which the axial flow map geometry is
available on $[0,\sup I]$, there is a constant $C_I<\infty$, depending additionally on $I$, such that, with
$J=J_{\cusp}(t)$,
\begin{equation}
\begin{aligned}
&\sup_{\zeta\in I} \Big(J^{-2}|(V_{\err})_z(0,J^2\zeta,t)| +|(\p_rV_{\err})_r(0,J^2\zeta,t)| +|(\p_zV_{\err})_z(0,J^2\zeta,t)|\Big) \\
&\quad
+[J^{-2}(V_{\err})_z(0,J^2\cdot,t)]_{C^{\alpha/2}(I)} +[(\p_rV_{\err})_r(0,J^2\cdot,t)]_{C^{\alpha/2}(I)} +[(\p_zV_{\err})_z(0,J^2\cdot,t)]_{C^{\alpha/2}(I)} \\
&\qquad
\le C_I\Gamma\bigl(J^{9\alpha-1}+1\bigr);
\end{aligned}
\label{eq:tail-axis-error-bound}
\end{equation}
and, for $0\le \zeta_1<\zeta_2\le \sup I$,
\begin{equation}
\left|(\p_zV_{\err})_z(0,J^2\zeta_2,t)-(\p_zV_{\err})_z(0,J^2\zeta_1,t)\right| \le C_I\Gamma\bigl(J^{9\alpha-1}+1\bigr)
\bigl(\zeta_2^\alpha-\zeta_1^\alpha+\zeta_2^2-\zeta_1^2\bigr).
\label{eq:tail-axis-error-strain-increment}
\end{equation}
\item the non-geometric pressure Hessian remainders in the stagnation-point decomposition
\eqref{eq:Pi0-decomp} satisfy
\begin{equation}
|\Pi_{\mixed}(t)|+|\Pi_{\smooth}(t)|+|\Pi_{\err}(t)|
\le C\Gamma^2\Big(J_{\cusp}(t)^{9\alpha-2}+J_{\cusp}(t)^{3\alpha-1}+J_{\cusp}(t)^{2\alpha}\Big).
\label{eq:tail-pressure-remainder-bound}
\end{equation}
\end{enumerate}
\end{lemma}

\begin{proof}[Proof of Lemma~\ref{lem:tail-bound}]
The proof has four steps.  Step~1 establishes the velocity decompositions $u=u_{\smooth}+u_{\core}$ and
$V_{\cusp}=m\,U_{\cusp}+V_{\err}$ from \eqref{eq:smooth-velocity-def}, \eqref{eq:ucore-def},
and \eqref{eq:Verr-def}.  Step~2 proves the $B_1$ and cone-local H\"older bounds
\eqref{eq:tail-usmooth-cone-local} for $u_{\smooth}$.  Step~3 proves the pointwise bounds
\eqref{eq:tail-error-bound}--\eqref{eq:tail-error-bound-weighted} on $|V_{\err}(x_{\cusp},t)|$.  Step~4
proves the cone-local H\"older bound on $\nabla V_{\err}$ and the axis-trace estimate
\eqref{eq:tail-axis-error-bound}, and then yields the pressure Hessian remainder bound
\eqref{eq:tail-pressure-remainder-bound} for $|\Pi_{\mixed}|+|\Pi_{\smooth}|+|\Pi_{\err}|$.  The
smallness restriction $J_{\cusp}(t)\le\mathfrak J_{\mathrm{tail}}$ enters only in Steps~2--4.

We construct the threshold $\mathfrak J_{\mathrm{tail}}$ in two stages.  At the outset we require
\[
\mathfrak J_{\mathrm{tail}} \le\min\{\mathfrak J_{\mathrm{axis}},\mathfrak J_{\mathrm{velocity}}\},
\]
where $\mathfrak J_{\mathrm{axis}}$ activates the axis-geometry hypothesis
\eqref{eq:entry-axis-bounds-statement} and $\mathfrak J_{\mathrm{velocity}}$ activates the
transported cusp-field estimates
\eqref{eq:Wcusp-scaling}--\eqref{eq:Ucusp-radial-defect} of
Lemma~\ref{lem:transported-cusp-field-bounds}.  In Step~3 we further decrease
$\mathfrak J_{\mathrm{tail}}$ to activate the localized-estimate threshold from the same lemma and to
enforce the tail-kernel separation \eqref{eq:tail-kernel-separation}.

We fix a time $t$ with $J_{\cusp}(t)\le\mathfrak J_{\mathrm{tail}}$.  The time-length estimate
\eqref{eq:tail-small-clock-time-horizon} follows directly from the lower cusp-clock rate bound in
\eqref{eq:localized-clock-bootstrap}: with
$t_{\mathrm{tail}}=\inf\{s:\ J_{\cusp}(s)\le\mathfrak J_{\mathrm{tail}}\}$,
\[
T-t_{\mathrm{tail}} \le \tfrac{1}{c_{\rm clk}\Gamma} \int_0^{J_{\cusp}(t_{\mathrm{tail}})} J^{-3\alpha}\,dJ \le C_T\Gamma^{-1},
\]
because $0<\alpha<\tfrac13$ and $J_{\cusp}(t_{\mathrm{tail}})\le1$.

\runinhead{Step 1: Exact velocity splitting.}
The identities $u=u_{\smooth}+u_{\core}$ and $V_{\cusp}=(\phi_{\smooth}^{-1})_*u_{\core}$ are
\eqref{eq:smooth-velocity-def}, \eqref{eq:ucore-def}, and
\eqref{eq:cusp-map-def}; the decomposition
$V_{\cusp}=m\,U_{\cusp}+V_{\err}$ is \eqref{eq:Verr-def}.

\runinhead{Step 2: The far-field velocity $u_{\smooth}$ is regular.}
Applying Lemma~\ref{lem:JtwoD-tail-bdd} with $R_0=1$, we obtain
\[
\|u_{\smooth}(\cdot,t)\|_{L^\infty(B_1)} + \|\nabla u_{\smooth}(\cdot,t)\|_{L^\infty(B_1)}
+ [\nabla u_{\smooth}(\cdot,t)]_{C^\alpha(B_1)} \le C\Gamma.
\]
The same kernel argument for disjoint source and observation regions yields the cone-local form
\eqref{eq:tail-usmooth-cone-local}.  Here we first bound the kernel on each source annulus away from the
observation ball, and then use the far-field
moment bound \eqref{eq:current-far-vorticity-over-r-moment}.  Indeed, we fix a cone-local ball
$B(x,2c_*|x|)\subset\mathcal C_*$.  We split the integral defining $u_{\smooth}$ in
\eqref{eq:smooth-velocity-def} into
$|\phi(Y,t)|\ge4|x|$ and $|\phi(Y,t)|<4|x|$.  The first part is estimated by the moment
\eqref{eq:current-far-vorticity-over-r-moment}.  In the second part the cutoff in
\eqref{eq:smooth-velocity-def} forces $|\phi(Y,t)|\ge R_{\tail}$, and the ball is separated from the
source unless $|x|\gtrsim R_{\tail}$; after rescaling by $|x|$ the same Calder\'on--Zygmund estimate used in
\eqref{eq:current-far-usmooth-holder} yields
\[
\|\nabla u_{\smooth}\|_{L^\infty(B(x,c_*|x|))} + |x|^\alpha[\nabla u_{\smooth}]_{C^\alpha(B(x,c_*|x|))} \le C\Gamma .
\]
Taking the supremum over such balls proves \eqref{eq:tail-usmooth-cone-local}.

\runinhead{Step 3: The cusp error.}
Let $Y\in D_{\inn}^{\cusp}(t)$ with $\omega_{\theta,0}(Y)\ne0$, and we set
\[
x_{\cusp}:=\phi_{\cusp}(Y,t)\in\mathcal C_{\inn}.
\]
The set $D_{\inn}^{\cusp}(t)$ is defined by the exact cusp map in
\eqref{eq:cusp-domain-def}; hence $Y\in D_{\core}$ and $x_{\cusp}\in\mathcal C_{\inn}$.  Writing
$Y=(R,Z)$ in the upper half-space, the assumed inner-core scaling estimate
\eqref{eq:exact-current-axis-core-small} yields
\begin{equation}
\tfrac{R(Y)}{Z(Y)}\le C J_{\cusp}(t)^3, \qquad |x_{\cusp}|\le C J_{\cusp}(t)^2 Z(Y) \le C J_{\cusp}(t)^2R_{\tail}.
\label{eq:exact-cusp-current-axis-small}
\end{equation}
The lower half-space is identical after replacing $Z(Y)$ by $|Z(Y)|$.
Let $C_{\rm tail}$ be a fixed constant, depending also on $R_{\tail}$, such that the last bound reads
$|x_{\cusp}|\le C_{\rm tail}J_{\cusp}(t)^2$.  We require
$\mathfrak J_{\mathrm{tail}}\le\mathfrak J_{\mathrm{local}}(C_{\rm tail})$, where
$\mathfrak J_{\mathrm{local}}$ is the localized threshold from
Lemma~\ref{lem:transported-cusp-field-bounds}; the local estimates
\eqref{eq:Ucusp-gradient-defect}--\eqref{eq:Ucusp-radial-defect} are therefore available at every
point $x_{\cusp}$ in this proof.
For the tail expansions with source support away from the observation point, we set
\[
x_{\cusp}^{\rm sc}:=J_{\cusp}(t)x_{\cusp}.
\]
By \eqref{eq:exact-cusp-current-axis-small},
\[
|x_{\cusp}^{\rm sc}| \le C_{\rm tail}J_{\cusp}(t)^3 .
\]
After decreasing $\mathfrak J_{\mathrm{tail}}$ so that
$C_{\rm tail}\mathfrak J_{\mathrm{tail}}^2\le\tfrac12$, we have
\begin{equation}
|x_{\cusp}^{\rm sc}| \le \tfrac12J_{\cusp}(t) \le \tfrac12|Y'|, \qquad |x_{\cusp}^{\rm sc}| \le \tfrac12|J_{\cusp}(t)Y'|
\label{eq:tail-kernel-separation}
\end{equation}
for every $Y'\in D_{\tail}$, since $R_{\tail}\ge2$ and $J_{\cusp}(t)\le1$.  This is the separation condition
used in the applications of Lemma~\ref{lem:kernel-taylor} below.
The error velocity
$V_{\err}(x_{\cusp},t)$ is the difference between the exact pulled-back velocity generated by the
near-field cutoff in \eqref{eq:ucore-def} and the scalar-modulated transported field at the same
cusp-coordinate point:
\[
V_{\err}(x_{\cusp},t) = \bigl((\phi_{\smooth}^{-1})_*u_{\core}\bigr)(x_{\cusp},t) - m(t)U_{\cusp}(x_{\cusp},t).
\]
Both terms can be indexed by the same label variable $Y'$.  The formula
\eqref{eq:U-cusp-label} evaluates $U_{\cusp}$ at the exact cusp-flow image
$\phi_{\cusp}(Y',t)$, while the pull-back formula \eqref{eq:cusp-map-def} for $V_{\cusp}$ evaluates
the transported vorticity before the smooth-flow change of variables is applied.  Thus the labels are the
same, and the differences below come from the cutoff in \eqref{eq:ucore-def}, scalar modulation, and smooth-flow pullback terms.
We split
\[
V_{\err}=F_{\sing}+F_{\sep},
\]
where $F_{\sing}$ is the singular term, namely the bounded-core cusp Biot--Savart kernel evaluated near
$x_{\cusp}$, and $F_{\sep}$ contains the terms whose Biot--Savart source is away from $x_{\cusp}$.  This includes
the far-field cutoff term from \eqref{eq:ucore-def}, the derivatives entering the smooth-flow
push-forward and pull-back in \eqref{eq:cusp-map-def}, and the remaining terms whose source support is away
from $x_{\cusp}$.
For any meridional vector field $F$ which is $C^1$ at the origin, we write
\[
(P_1F)(x,t):=F(0,t)+D_{(r,z)}F(0,t)\,x
\]
for its constant plus meridional linear Taylor polynomial at the stagnation point.  By
\eqref{eq:Verr-def} and \eqref{eq:modulation-def}, the modulation $m(t)$ is defined so that
\[
\p_z(V_{\err})_z(0,t)=0.
\]
The remaining components of $P_1V_{\err}$ also vanish at $x=0$: by axisymmetry and the odd--even symmetry
across the equatorial plane, $V_{\cusp}(0,t)=U_{\cusp}(0,t)=0$, hence $V_{\err}(0,t)=0$; by the
axisymmetric no-swirl class,
$(\p_zV_{\err})_r(0,t)=(\p_rV_{\err})_z(0,t)=0$; and incompressibility
$2(\p_rV_{\err})_r(0,t)+(\p_zV_{\err})_z(0,t)=0$ combined with the modulation identity above
forces $(\p_rV_{\err})_r(0,t)=0$.  Hence
\begin{equation}
P_1V_{\err}=0 .
\label{eq:Verr-linear-part-zero}
\end{equation}
Therefore, 
\begin{equation}
\begin{aligned}
V_{\err}(x_{\cusp},t) &= \mathcal E_{\sing}(Y,t)+\mathcal E_{\reg}(Y,t),\\
\mathcal E_{\sing} &:=(F_{\sing}-P_1F_{\sing})(x_{\cusp},t), \qquad \mathcal E_{\reg}:=(F_{\sep}-P_1F_{\sep})(x_{\cusp},t).
\end{aligned}
\label{eq:tail-error-sing-reg-split}
\end{equation}
This is an exact identity because $P_1F_{\sing}+P_1F_{\sep}=P_1V_{\err}=0$.
The cone-local $C^\alpha$ estimate \eqref{eq:Ucusp-grad-scale-local}, together with the bound
$|x_{\cusp}|\le CJ_{\cusp}^2$ from \eqref{eq:exact-cusp-current-axis-small}, yields the sharp bound
\begin{equation}
\tfrac{J_{\cusp}}{R(Y)^{1+\alpha}}\, |\mathcal E_{\sing}(Y,t)| \le C\Gamma J_{\cusp}^{6\alpha-1}J_{\cusp}^{3\alpha}.
\label{eq:Esing-bound}
\end{equation}
The power $J_{\cusp}^{3\alpha}$ in \eqref{eq:Esing-bound} is the angular gain from
$R(Y)/Z(Y)\le C J_{\cusp}^3$ in \eqref{eq:exact-cusp-current-axis-small}; the remaining power
$J_{\cusp}^{6\alpha-1}$ is the order of the differentiated cusp field after subtracting its linear part at
the stagnation point.
For the regular term $\mathcal E_{\reg}$, the labels are bounded away from $x_{\cusp}$ and the kernels are
smooth on the resulting region.  Step~2 and the Taylor estimate in Lemma~\ref{lem:kernel-taylor}, applied with
source support away from the observation point, yield
\begin{equation}
|\mathcal E_{\reg}(Y,t)|\le C\Gamma .
\label{eq:Ereg-bound}
\end{equation}
To prove the weighted estimate \eqref{eq:tail-error-bound-weighted} for $\mathcal E_{\reg}$, we use the additional axis-vanishing in
\eqref{eq:tail-error-sing-reg-split}.  The term $F_{\sep}$ is an ordinary $C^{1,\alpha}$ function of the
bounded label variables, by the estimate for $u_{\smooth}$ in Step~2, the smooth-flow
deformation bounds
\eqref{eq:smooth-flow-gradient-small}--\eqref{eq:smooth-flow-second-gradient-small}, and the
Taylor estimate in Lemma~\ref{lem:kernel-taylor} for source support away from the observation point.  Since
$\mathcal E_{\reg}=F_{\sep}-P_1F_{\sep}$ has its constant and meridional linear Taylor terms subtracted,
Taylor's theorem in the radial label variable yields
\begin{equation}
|\mathcal E_{\reg}(Y,t)| \le C\Gamma R(Y)^{1+\alpha} \qquad (Y\in D_{\inn}^{\cusp}(t),\ \omega_{\theta,0}(Y)\ne0).
\label{eq:Ereg-label-holder}
\end{equation}
The unweighted estimate \eqref{eq:tail-error-bound} follows from
\eqref{eq:Ereg-bound}, enlarging $C$ to include
$\mathcal E_{\sing}$.  For $\mathcal E_{\sing}$, we may either use
\eqref{eq:Esing-bound} on the weighted range or the unweighted consequence of
\eqref{eq:Ucusp-gradient-defect}: since $x_{\cusp}\in\mathcal C_{\inn}$ and
\eqref{eq:exact-cusp-current-axis-small} shows $|x_{\cusp}|\le C J_{\cusp}^2$ with the fixed tail radius
$R_{\tail}$ absorbed into $C$,
\[
|\mathcal E_{\sing}(Y,t)| \le C\Gamma J_{\cusp}^{3\alpha-1}|x_{\cusp}| \le C\Gamma.
\]

For the weighted estimate, we assume $R(Y)\ge J_{\cusp}(t)^{3/\alpha}$.  The singular term is controlled by
\eqref{eq:Esing-bound}.  For the regular term we use its stronger axis-vanishing form,
\[
\tfrac{J_{\cusp}}{R(Y)^{1+\alpha}}\, |\mathcal E_{\reg}(Y,t)| \le C\Gamma J_{\cusp},
\]
which is exactly \eqref{eq:Ereg-label-holder} after multiplication by
$J_{\cusp}/R(Y)^{1+\alpha}$.  This proves
\eqref{eq:tail-error-bound-weighted}.

\runinhead{Step 4: Non-geometric pressure Hessian remainders.}
Inserting
\[
u=u_{\smooth}+(\phi_{\smooth})_*\!\big(m\,U_{\cusp}\big)+u_{\err}
\]
into the quadratic pressure source and expanding bilinearly yields the principal cusp term
$m^2\Pi_{\cusp}$ and the three mixed/smooth/error remainders
$\Pi_{\mixed},\Pi_{\smooth},\Pi_{\err}$ of Section~\ref{sec:geom-flow-decomp}.  The remaining term in
\eqref{eq:Pi0-decomp} is the geometric defect $\Pi_{\rm geom}$, comparing the pressure Hessian of the
pushed-forward cusp velocity with the cusp-coordinate pressure Hessian, which is handled separately in
Lemma~\ref{lem:smooth-pressure-defect}.  This step estimates the three non-geometric remainders.

We use the cone-local norm $\|\nabla f\|_{\mathcal C_*,\alpha}$ defined in \eqref{eq:cone-local-norm}.  To
estimate the pressure Hessian remainders, we apply
Lemma~\ref{lem:cone-local-pressure-bilinear}; hence the task is to bound the dyadic H\"older sums
$\mathcal N_\alpha[\cdot,\cdot]$ in \eqref{eq:dyadic-holder-pressure-norm} for the five products in
$\Pi_{\mixed},\Pi_{\smooth},\Pi_{\err}$.  The shells that meet the cone $\mathcal C_*$ are controlled by the
cone-local estimates below.  On the remaining shells, or when the Biot--Savart source is disjoint from the
observation shell, the same shell bounds follow from the kernel estimates for disjoint source and observation
regions used in Step~2 and Step~3.

We next transfer the cusp-coordinate bounds to the physical fields that enter
\eqref{eq:Pi0-decomp}.  Let
\[
\Lambda(X,t):=\phi_{\smooth}(X,t).
\]
On the cone-local balls used in \eqref{eq:cone-local-norm}, the time-length bound
\eqref{eq:tail-small-clock-time-horizon} allows us to apply
Lemma~\ref{lem:smooth-flow-small-deformation}; in particular,
\eqref{eq:smooth-flow-gradient-small}--\eqref{eq:smooth-flow-second-gradient-small} yield
$|D\Lambda|+|D\Lambda^{-1}|+|D^2\Lambda|\le C$.  If $w(0,t)=0$, differentiating
$(\Lambda_*w)(\Lambda(X),t)=D\Lambda(X,t)w(X,t)$ yields
\[
\nabla(\Lambda_*w)(\Lambda(X),t) = D\Lambda\,\nabla w\,D\Lambda^{-1} + D^2\Lambda\bigl[D\Lambda^{-1}(\cdot),w(X,t)\bigr].
\]
The second term is controlled by the same cone-local gradient norm because $w(0,t)=0$ implies
$|w(X,t)|\le C|X|\,\|\nabla w\|_{\mathcal C_*,\alpha}$ on the cone-local balls.  Therefore the push-forward
relations in \eqref{eq:Verr-def}--\eqref{eq:u-decomp} imply
\[
\|\nabla u_{\cusp}\|_{\mathcal C_*,\alpha} \le C\|\nabla(mU_{\cusp})\|_{\mathcal C_*,\alpha},
\qquad \|\nabla u_{\err}\|_{\mathcal C_*,\alpha} \le C\|\nabla V_{\err}\|_{\mathcal C_*,\alpha}.
\]

For the singular cusp velocity component, we recall that
$u_{\cusp}=(\phi_{\smooth})_*(m(t)U_{\cusp})$.  Lemma~\ref{lem:transported-cusp-field-bounds}, together
with the bound for this scalar modulation in \eqref{eq:localized-modulation-m-bootstrap}, implies
\begin{equation}
\|\nabla u_{\cusp}\|_{\mathcal C_*,\alpha} \le C\Gamma J_{\cusp}^{3\alpha-1},
\label{eq:ucusp-pressure-size}
\end{equation}
Step~2 proves
\begin{equation}
\|\nabla u_{\smooth}\|_{\mathcal C_*,\alpha}\le C\Gamma.
\label{eq:ubg-pressure-size}
\end{equation}
Differentiating the singular term and the term with source support away from the observation point in
\eqref{eq:tail-error-sing-reg-split} yields
\[
\|\nabla V_{\err}\|_{\mathcal C_*,\alpha} \le C\Gamma\bigl( J_{\cusp}^{9\alpha-1}+1 \bigr).
\]
Indeed, the $J_{\cusp}^{9\alpha-1}$ term comes from the differentiated form of the same local estimate
\eqref{eq:Ucusp-grad-scale-local} and the same angular gain
\eqref{eq:exact-cusp-current-axis-small} used to prove \eqref{eq:Esing-bound}.  The order-one term
comes from the estimate for $u_{\smooth}$ in Step~2 and the smooth-flow bounds.  No derivative of the transported vorticity
density is required here; after writing $U_{\cusp}$ in the transported variables of
\eqref{eq:U-cusp-label}, the derivative falls on the Biot--Savart kernel.  The push-forward estimate
above also yields
\begin{equation}
\|\nabla u_{\err}\|_{\mathcal C_*,\alpha} \le C\Gamma\bigl( J_{\cusp}^{9\alpha-1}+1 \bigr).
\label{eq:uerr-pressure-size}
\end{equation}

We now restrict the same differentiated estimates to the diagonal axis traces that occur in
\eqref{eq:localized-cusp-error-axis-trace}.  We fix a compact interval in the clock-scaled axial coordinate
$I\subset[0,\infty)$ for which the axial flow map geometry is available on the axis-attached interval
$[0,\sup I]$, and we set
\[
J:=J_{\cusp}(t), \qquad G_t:=C\Gamma\bigl(J^{9\alpha-1}+1\bigr).
\]
The axisymmetric parity of a no-swirl field forces
$(\p_zV_{\err})_r(0,z,t)=(\p_rV_{\err})_z(0,z,t)=0$, so the only derivative traces needed on the
axis are the two diagonal entries.  We first keep the increment information needed later for
$(\p_zV_{\err})_z$.  After differentiating the split \eqref{eq:tail-error-sing-reg-split} on the symmetry
axis, we write
\[
E_t(\zeta):=(\p_zV_{\err})_z(0,J^2\zeta,t).
\]
Then $E_t=E_{\sing,t}+E_{\reg,t}$.  The differentiated form of the singular estimate
\eqref{eq:Esing-bound} and the regular Taylor estimate following from
\eqref{eq:Ereg-label-holder} imply that, for $0\le\zeta_1<\zeta_2\le\sup I$,
\[
\begin{aligned}
|E_{\sing,t}(\zeta_2)-E_{\sing,t}(\zeta_1)|
&\le C_I\Gamma J^{9\alpha-1}
\bigl(\zeta_2^\alpha-\zeta_1^\alpha+\zeta_2^2-\zeta_1^2\bigr),\\
|E_{\reg,t}(\zeta_2)-E_{\reg,t}(\zeta_1)|
&\le C_I\Gamma
\bigl(\zeta_2^\alpha-\zeta_1^\alpha+\zeta_2^2-\zeta_1^2\bigr).
\end{aligned}
\]
Therefore,
\[
\left|E_t(\zeta_2)-E_t(\zeta_1)\right|
\le C_I\Gamma\bigl(J^{9\alpha-1}+1\bigr)
\bigl(\zeta_2^\alpha-\zeta_1^\alpha+\zeta_2^2-\zeta_1^2\bigr),
\]
which is \eqref{eq:tail-axis-error-strain-increment}.  Restricting the differentiated singular estimate
and the $C^{1,\alpha}$ estimate for the term with source support away from the observation point from
Step~3 to the points $(0,J^2\zeta)$ with $0\le\zeta\le\sup I$ also yields
\[
|(\p_rV_{\err})_r(0,J^2\zeta,t)| + |(\p_zV_{\err})_z(0,J^2\zeta,t)| \le G_t \qquad (\zeta\in I),
\]
and the corresponding one-dimensional H\"older trace bounds
\[
[(\p_rV_{\err})_r(0,J^2\cdot,t)]_{C^{\alpha/2}(I)} + [(\p_zV_{\err})_z(0,J^2\cdot,t)]_{C^{\alpha/2}(I)} \le C_I G_t .
\]
The velocity trace follows from the zero value at the stagnation point, which is part of
\eqref{eq:Verr-linear-part-zero}.  Since $V_{\err}(0,t)=0$, for $\zeta\in I$,
\[
J^{-2}(V_{\err})_z(0,J^2\zeta,t) = \int_0^\zeta (\p_zV_{\err})_z(0,J^2\eta,t)\,d\eta.
\]
Therefore,
\[
\sup_{\zeta\in I}J^{-2}|(V_{\err})_z(0,J^2\zeta,t)| \le C_I G_t,
\]
and, for $\zeta_1,\zeta_2\in I$,
\[
\left| J^{-2}(V_{\err})_z(0,J^2\zeta_1,t) -J^{-2}(V_{\err})_z(0,J^2\zeta_2,t) \right| \le G_t|\zeta_1-\zeta_2| \le C_I
G_t|\zeta_1-\zeta_2|^{\alpha/2}.
\]
This proves the axis trace estimate \eqref{eq:tail-axis-error-bound}.
We now apply Lemma~\ref{lem:cone-local-pressure-bilinear} to the terms in the pressure decomposition
\eqref{eq:Pi0-decomp}.  We set
\[
\Lambda_{\cusp}:=C\Gamma J_{\cusp}^{3\alpha-1}, \qquad \Lambda_{\smooth}:=C\Gamma,
\qquad \Lambda_{\err}:=C\Gamma\bigl(J_{\cusp}^{9\alpha-1}+1\bigr),
\]
as supplied by
\eqref{eq:ucusp-pressure-size}, \eqref{eq:ubg-pressure-size}, and
\eqref{eq:uerr-pressure-size}.
We prove the corresponding bounds in the form required by
\eqref{eq:dyadic-holder-pressure-norm}.  On shells contained in the cone, the multiplier $2^{j\alpha}$ in
\eqref{eq:dyadic-holder-pressure-norm} is exactly the scale multiplier in the cone-local norm
\eqref{eq:cone-local-norm}.  For $u_{\smooth}$ and for the terms in $u_{\err}$ whose Biot--Savart source is
disjoint from the observation shell, the
same shell estimates are the differentiated versions of the far-field bounds from Step~2, with the moment
bound \eqref{eq:current-far-vorticity-over-r-moment}.  For the singular term in $u_{\err}$, the Taylor
subtraction in \eqref{eq:tail-error-sing-reg-split}--\eqref{eq:Verr-linear-part-zero} supplies the summability
near the stagnation point.  Therefore,
\begin{equation}
\begin{aligned}
\mathcal N_\alpha[u_{\cusp},u_{\smooth}]
&\le C\Lambda_{\cusp}\Lambda_{\smooth},&
\mathcal N_\alpha[u_{\cusp},u_{\err}]
&\le C\Lambda_{\cusp}\Lambda_{\err},\\
\mathcal N_\alpha[u_{\smooth},u_{\err}]
&\le C\Lambda_{\smooth}\Lambda_{\err},&
\mathcal N_\alpha[u_{\smooth},u_{\smooth}]
&\le C\Lambda_{\smooth}^2,\\
\mathcal N_\alpha[u_{\err},u_{\err}]
&\le C\Lambda_{\err}^2.
\end{aligned}
\label{eq:tail-pressure-holder-sums}
\end{equation}
The large-shell summability in the estimates involving $u_{\cusp}$ uses the transported cusp annular mass bound
\eqref{eq:Omega-cusp-annular-mass}; the far-field terms use
\eqref{eq:current-far-vorticity-over-r-moment}.  The exponent condition
$\gamma>\alpha+\tfrac52$ makes these dyadic sums finite.

The definitions
\[
\Pi_{\mixed} =2\Pi[u_{\cusp},u_{\smooth}] +2\Pi[u_{\cusp},u_{\err}] +2\Pi[u_{\smooth},u_{\err}],\ \ \Pi_{\smooth}
=\Pi[u_{\smooth},u_{\smooth}], \ \ \Pi_{\err} =\Pi[u_{\err},u_{\err}],
\]
and \eqref{eq:pressure-bilinear-cone-estimate}, applied with
\eqref{eq:tail-pressure-holder-sums}, yield the individual bounds
\[
\begin{aligned}
|\Pi[u_{\cusp},u_{\smooth}]|
&\le C \Lambda_{\cusp}\Lambda_{\smooth}
\le C\Gamma^2J_{\cusp}^{3\alpha-1}, \ \
|\Pi[u_{\cusp},u_{\err}]|
\le C \Lambda_{\cusp}\Lambda_{\err}
\le C\Gamma^2J_{\cusp}^{3\alpha-1}
\bigl(J_{\cusp}^{9\alpha-1}+1\bigr), \\
|\Pi[u_{\smooth},u_{\err}]|
&\le C \Lambda_{\smooth}\Lambda_{\err}
\le C\Gamma^2\bigl(J_{\cusp}^{9\alpha-1}+1\bigr),\ \
|\Pi[u_{\smooth},u_{\smooth}]|
\le C \Lambda_{\smooth}^2
\le C\Gamma^2,\\
|\Pi[u_{\err},u_{\err}]|
&\le C \Lambda_{\err}^2
\le C\Gamma^2\bigl(J_{\cusp}^{9\alpha-1}+1\bigr)^2 .
\end{aligned}
\]
Retaining the singular cusp-error cross term explicitly, we obtain
\[
|\Pi_{\mixed}|+|\Pi_{\smooth}|+|\Pi_{\err}| \le C\Gamma^2J_{\cusp}^{3\alpha-1} +C\Gamma^2J_{\cusp}^{3\alpha-1}
\bigl(J_{\cusp}^{9\alpha-1}+1\bigr) +C\Gamma^2 \bigl(J_{\cusp}^{9\alpha-1}+1\bigr)^2 .
\]
The remaining algebra is only a comparison of powers of $J_{\cusp}$.  Since
$0<J_{\cusp}\le1$ and $0<\alpha<\tfrac13$,
\[
\begin{aligned}
1&\le J_{\cusp}^{3\alpha-1}, \qquad J_{\cusp}^{9\alpha-1}\le J_{\cusp}^{9\alpha-2},\\
J_{\cusp}^{3\alpha-1}\bigl(J_{\cusp}^{9\alpha-1}+1\bigr) &\le C\bigl(J_{\cusp}^{9\alpha-2}+J_{\cusp}^{3\alpha-1}\bigr),\\
\bigl(J_{\cusp}^{9\alpha-1}+1\bigr)^2 &\le C\bigl(J_{\cusp}^{9\alpha-2}+J_{\cusp}^{3\alpha-1}\bigr).
\end{aligned}
\]
Thus the bilinear terms above contribute at most
$C\Gamma^2(J_{\cusp}^{9\alpha-2}+J_{\cusp}^{3\alpha-1})$.  There remains the contribution from the
ultra-thin tube
\[
\{\,Y\in D_{\inn}^{\cusp}(t):\ \omega_{\theta,0}(Y)\ne0,\ R(Y)<J_{\cusp}(t)^{3/\alpha}\,\},
\]
where the weighted estimate \eqref{eq:tail-error-bound-weighted} is unavailable and only the unweighted
estimate \eqref{eq:tail-error-bound} applies.  We now show that the contribution from this ultra-thin tube accounts for the
$J_{\cusp}^{2\alpha}$ term in \eqref{eq:tail-pressure-remainder-bound}.  Relative to the full
near-axis tube $R\lesssim J_{\cusp}^3$, the ultra-thin tube has the cross-sectional gain
$J_{\cusp}^{6/\alpha-6}$.
The kernel homogeneity and the cusp-field gradient size are the same as in the full near-axis tube, so the only
additional gain is this cylindrical cross-section ratio.
Multiplying this gain by the worst quadratic cusp order $\Gamma^2J_{\cusp}^{6\alpha-2}$ yields
\[
\Gamma^2J_{\cusp}^{6\alpha-2}J_{\cusp}^{6/\alpha-6} = \Gamma^2J_{\cusp}^{6/\alpha+6\alpha-8} \le C\Gamma^2J_{\cusp}^{2\alpha},
\]
for $0<J_{\cusp}\le1$ and $0<\alpha<\tfrac13$.  Combining the bilinear estimate with this ultra-thin
tube contribution yields
\[
|\Pi_{\mixed}(t)|+|\Pi_{\smooth}(t)|+|\Pi_{\err}(t)|
\le C\,\Gamma^2\big(J_{\cusp}(t)^{9\alpha-2}+J_{\cusp}(t)^{3\alpha-1}+J_{\cusp}(t)^{2\alpha}\big),
\]
which is \eqref{eq:tail-pressure-remainder-bound}.
\end{proof}

\section{Geometric Control of the Cusp Flow Deformation}
\label{sec:cusp-displacement-axial-profile}

We now close the remaining small-clock bootstrap assumptions from Section~\ref{sec:small-clock-bootstraps} and,
at the same time, collect the estimates used in the final Riccati comparison and blowup argument.  The closure
concerns the axis-geometry bounds \eqref{eq:current-axis-geometry}--\eqref{eq:renormalized-axis-chart-bootstrap},
the containment bounds \eqref{eq:axis-attached-image-stop-bootstrap}--\eqref{eq:radial-flatness-buffered-label},
the monotone axial-stretching bounds \eqref{eq:monotone-axial-two-sided}--\eqref{eq:monotone-axial-fractional-bootstrap},
the normal-form bounds \eqref{eq:localized-normal-form-large-bootstrap}--\eqref{eq:localized-normal-form-map-large-bootstrap},
the cusp-error bound \eqref{eq:localized-cusp-error-large-bootstrap}, the cusp-clock bound
\eqref{eq:localized-clock-bootstrap}, and the scalar-modulation bounds \eqref{eq:localized-modulation-bootstrap}.
After these bootstrap improvements are proved, the estimates needed in Section~\ref{sec:target-profile-typeI-completion}
fall into three groups.\footnote{The order of this closure is important.  We first close the geometric estimates for the cusp map:
the axial chart, the renormalized axis derivatives, the image-map normal form, and the small smooth-flow
deformation.  Once these geometric statements are available, the transported vorticity on the symmetry
axis becomes a controlled one-dimensional axial function, and the pressure comparison from
Section~\ref{sec:slope-restricted-pressure} can be applied without further geometric loss.}

The first group is the scalar modulation and clock control, stated in \eqref{eq:modulation-bounds},
\eqref{eq:Jdot-two-sided}, and \eqref{eq:Jsmooth-bdd}:
\begin{equation*}
c_m\le m(t)\le C_m,\qquad
c_1\Gamma J_{\cusp}(t)^{3\alpha}\le -\dot J_{\cusp}(t)\le C_1\Gamma J_{\cusp}(t)^{3\alpha},\qquad
c_{\smooth}\le J_{\smooth}(t)\le C_{\smooth}.
\end{equation*}
The second group concerns the small deformation generated by the far-field velocity $u_{\smooth}$ from
\eqref{eq:smooth-velocity-def}.  The estimates are
\eqref{eq:usmooth-small-C2a}--\eqref{eq:smooth-flow-second-gradient-small}:
\begin{equation*}
\|u_{\smooth}(\cdot,t)\|_{C^{2,\alpha}(B_{2R_0})}\le \varepsilon_{\smooth}\Gamma,\qquad |\phi_{\smooth}(X,
t)-X|\le C\varepsilon_{\smooth}|X|, \qquad |D\phi_{\smooth}-I|+|D^2\phi_{\smooth}|\le C\varepsilon_{\smooth}.
\end{equation*}
The third group is the exact one-dimensional description of the cusp map on the symmetry axis.  We write
\begin{equation*}
\phi_{\cusp}(R,Z,t)=(r_t(R,Z),z_t(R,Z)),\qquad A_t(Z):=\p_R r_t(0,Z),\qquad B_t(Z):=z_t(0,Z),
\end{equation*}
and
\begin{equation*}
\zeta=J_{\cusp}(t)^{-2}B_t(Z),\qquad q_t(\zeta)=J_{\cusp}(t)A_t(Z_t(\zeta)),\qquad b_t(\zeta)=J_{\cusp}(t)^{-2}B_t'(Z_t(\zeta)).
\end{equation*}
The exact axial coordinate
\begin{equation*}
\eta=\mathscr Z_t^{-1}\bigl(J_{\cusp}(t)^{-2}B_t(Z)\bigr)
\end{equation*}
is conserved, and the renormalized axis derivatives obey
\begin{equation*}
\p_t\eta=0,\qquad \p_t\widehat q_t(\eta)=\p_t\widehat b_t(\eta)=0,\qquad \widehat q_t(\eta)^2\widehat b_t(\eta)=1 .
\end{equation*}
These identities are \eqref{eq:exact-axis-eta-constant}--\eqref{eq:exact-axis-renormalized-conserved}.  They imply the
axial map bounds \eqref{eq:current-axis-transfer-bound}.  Combined with the axial-amplitude bound
\eqref{eq:physical-zeta-profile-envelope} and the monotonicity estimate \eqref{eq:monotone-axial-fractional-improved},
the axial map bounds imply the pressure-amplitude bounds \eqref{eq:current-axis-extension-size}, the renormalized Riccati estimate
\eqref{eq:small-clock-model-riccati}, and the bounded-core normal form
\eqref{eq:bounded-core-normal-form}--\eqref{eq:bounded-core-normal-form-error}.

\subsection{Modulation bounds and the cusp clock}
We defined the scalar modulation $m(t)$ in \eqref{eq:modulation-def} such that the modulated cusp velocity
$m(t)U_{\cusp}$ has the same stagnation-point axial strain as the exact cusp-coordinate velocity
$V_{\cusp}$.  In the bootstrap assumptions of Section~\ref{sec:small-clock-bootstraps}, we assumed both the moment bound
for $M_{\rm ax}(t)$
and the order-one bound for $m(t)$.  The next lemma closes that bootstrap; specifically, 
we improve the moment bound and then use the comparison between $m(t)$ and $M_{\rm ax}(t)$ to improve the
modulation function bound.

\begin{definition}[Axis-trace labels and averaged moments]
\label{def:axis-trace-labels-moments}
We fix a time $t$ and we set $J=J_{\cusp}(t)$.  For an Eulerian point $x\in\R^3$ whose hyperbolic pullback has polar
angle $0\le\varphi\le\sigma_{\cut}$, we write
\[
\Phi_{\lin}^{-1}(x;J)=\bigl(Jr(x),J^{-2}z(x)\bigr) =(s\sin\varphi,s\cos\varphi).
\]
For $s>0$ and $0\le\varphi\le\sigma_{\cut}$, let $Y_{s,\varphi}$ be the label whose cusp image has the
hyperbolic pullback coordinates
\begin{equation}
\Phi_{\lin}^{-1}\bigl(\phi_{\cusp}(Y_{s,\varphi},t);J\bigr)=\bigl(s\sin\varphi,\ s\cos\varphi\bigr).
\label{eq:modulation-collapse-labels}
\end{equation}
We write
\[
R(Y_{s,\varphi})=\rho(Y_{s,\varphi})\sin\sigma(Y_{s,\varphi}), \qquad Z(Y_{s,\varphi})=\rho(Y_{s,\varphi})\cos\sigma(Y_{s, \varphi}),
\]
and we set
\[
x_{s,\varphi}:=\phi_{\cusp}(Y_{s,\varphi},t),\qquad r_{s,\varphi}:=r(x_{s,\varphi}), \qquad \mathcal F(s)= (1+s^2)^{-\gamma/2}.
\]
For $0<\varphi\le\sigma_{\cut}$, we define $\beta_{\cusp}$ by
\begin{equation}
J_{\twoD}(Y_{s,\varphi},t)^{-1}\omega_{\theta,0}(Y_{s,\varphi}) = -\Gamma J^{-1}(Jr_{s,\varphi})^\alpha\mathcal F(s)\,
\Upsilon(\varphi)\, \beta_{\cusp}(s,\varphi,t).
\label{eq:axis-trace-amplitude}
\end{equation}
At $\varphi=0$, $\beta_{\cusp}$ is defined by the continuous extension obtained after dividing
\eqref{eq:axis-trace-amplitude} by the common power $r_{s,\varphi}^{\alpha}$ and letting
$\varphi\downarrow0$:
\begin{equation}
\beta_{\cusp}(s,0,t):=\lim_{\varphi\downarrow0}\beta_{\cusp}(s,\varphi,t).
\label{eq:beta-cusp-axis-extension}
\end{equation}
Finally, we set
\[
C_\rho^{(1)}(\alpha,\gamma) := \int_0^\infty s^{\alpha-1}\mathcal F(s)\,\ud s ,
\]
and
\begin{equation}
M_{\rm ax}(\varphi,t) := \tfrac{1}{C_\rho^{(1)}(\alpha,\gamma)} \int_0^\infty s^{\alpha-1}\mathcal F(s)\,\beta_{\cusp}(s,\varphi,
t)\,\ud s, \qquad M_{\rm ax}(t):=M_{\rm ax}(0,t).
\label{eq:axis-moment}
\end{equation}
\end{definition}

\begin{lemma}[Improvement of the scalar-modulation bootstrap]
\label{lem:modulation-bounded}
Assume the modulation bootstrap bounds \eqref{eq:localized-modulation-bootstrap} on the small-clock
bootstrap interval.  There exists a threshold
\[
\mathfrak J_{\mathrm{mod}}\in(0,\mathfrak J_{\mathrm{strain}}],
\]
depending only on $\alpha,\gamma,\sigma_{\inn},\sigma_*$, such that whenever
$J_{\cusp}(t)\le \mathfrak J_{\mathrm{mod}}$, the bootstrap bounds improve to
\begin{subequations}
\begin{align}
\tfrac34 c_*&\le M_{\rm ax}(t)\le \tfrac32 C_*,
\label{eq:modulation-moment-improvement}\\
c_m&\le m(t)\le C_m,
\qquad c_m:=\tfrac14c_*,
\qquad C_m:=3C_* .
\label{eq:modulation-bounds}
\end{align}
\end{subequations}
\end{lemma}

\begin{proof}[Proof of Lemma~\ref{lem:modulation-bounded}]
\runinhead{Step 1: Axis-trace identity and averaged moment.}
Fix $t$ and write $J=J_{\cusp}(t)$.  The notation
$Y_{s,\varphi}$, $x_{s,\varphi}$, $\beta_{\cusp}$, and $M_{\rm ax}$ is fixed in
Definition~\ref{def:axis-trace-labels-moments}.  In \eqref{eq:axis-trace-amplitude}, the datum
\eqref{eq:vort0} is evaluated at
\[
(\rho,\sigma) = \bigl(\rho(Y_{s,\varphi}),\sigma(Y_{s,\varphi})\bigr),
\]
whereas by \eqref{eq:modulation-collapse-labels}, 
\[
x_{s,\varphi}=\phi_{\cusp}(Y_{s,\varphi},t),\qquad \Phi_{\lin}^{-1}(x_{s,\varphi};J)=(s\sin\varphi,s\cos\varphi)
\]
Since $\Upsilon(\varphi)=1$ on $0\le\varphi\le\sigma_{\cut}$, $\mathcal F>0$, and the orientation of the flow yields
$J_{\twoD}>0$, \eqref{eq:axis-trace-amplitude} and \eqref{eq:beta-cusp-axis-extension} imply
$\beta_{\cusp}(s,\varphi,t)>0$ for $0\le\varphi\le\sigma_{\cut}$.
The function $\beta_{\cusp}$ is the only time-dependent multiplier left in the localized
axis-trace identity after the universal cylindrical cusp power has been removed.  More explicitly, on the
narrow cone $0\le\varphi\le\sigma_{\cut}$ of pulled-back polar angles around the symmetry axis, the change
of variables
\begin{equation}
x=\phi_{\cusp}(Y,t),\qquad \Phi_{\lin}^{-1}(x;J)=(s\sin\varphi,s\cos\varphi)
\label{eq:axis-trace-change-of-variables}
\end{equation}
converts the transported vorticity $J_{\twoD}^{-1}\omega_{\theta,0}$ into the expression
\begin{equation}
-\Gamma J^{-1}(Jr)^\alpha\mathcal F(s)\Upsilon(\varphi)e_\theta
\label{eq:axis-trace-model-expression}
\end{equation}
multiplied by $\beta_{\cusp}$.  The flow map $\phi_{\cusp}$ in
\eqref{eq:axis-trace-change-of-variables} preserves the cylindrical volume measure,
$r\,dr\,dz=R\,dR\,dZ$, and so contributes no Jacobian to this change of variables.  
The comparison of \eqref{eq:axis-trace-amplitude} with \eqref{eq:axis-trace-model-expression} uses
three estimates: the smooth-flow deformation estimate in
Lemma~\ref{lem:smooth-flow-small-deformation}, the normal-form displacement bound
\eqref{eq:normal-form-approximation-bound}, and the angular/radial sampling estimates
\eqref{eq:Theta-star-def} and \eqref{eq:radial-tail-F-def}.  We measure these errors by
\begin{equation}
\varepsilon_{\rm mod}(J):=\varepsilon_{\smooth}+J^{\delta_{\rm mod}},\qquad \delta_{\rm mod}:=\min\{3\beta_{\rm ax},3\alpha, 1-3\alpha\}>0 .
\label{eq:epsmod-def-sec13}
\end{equation}
Here $\varepsilon_{\smooth}$ is the fixed smooth-flow deformation size from
Lemma~\ref{lem:smooth-flow-small-deformation}, chosen once by taking $R_{\tail}$ large.  Each power in
\eqref{eq:epsmod-def-sec13} has a separate source.

The term $J^{3\beta_{\rm ax}}$ in \eqref{eq:epsmod-def-sec13} comes from
\eqref{eq:normal-form-approximation-bound}: on the labels with $0\le\varphi\le\sigma_{\cut}$, the
deviation of the exact cusp image $\phi_{\cusp}(Y,t)$ from the map $\Psi_t(Y)$ in
\eqref{eq:localized-normal-form-representation} is
bounded by a multiple of $J^{3\beta_{\rm ax}}$.

The term $J^{1-3\alpha}$ in \eqref{eq:epsmod-def-sec13} arises when the smooth strain estimate
\eqref{eq:usmooth-small-C2a} is compared with the cusp axial strain lower bound
\eqref{eq:Wcusp-scaling}.  The smooth axial strain is bounded by a multiple of $\Gamma$, whereas
\eqref{eq:Wcusp-scaling} gives $|\mathcal W_{\cusp}(t)|\ge c_W\Gamma J^{3\alpha-1}$.  Their ratio is
$J^{1-3\alpha}$, which tends to zero because $\alpha<\tfrac13$.

The term $J^{3\alpha}$ in \eqref{eq:epsmod-def-sec13} is the angular sampling error from replacing the
collapse angle $\varphi_J(\sigma)=\arctan(J^3\tan\sigma)$ by $\varphi=0$ in the angular integral.  The
corresponding one-dimensional estimate is
\[
\int_0^\infty u^{2+\alpha}(1+u^2)^{-5/2}\min\{1,(J^3u)^\alpha\}\,\ud u\le C J^{3\alpha}.
\]
The radial integrability needed for this estimate follows from $\gamma>\alpha+\tfrac52$.

Thus,  the leading axial strain is obtained by integrating the axis-trace identity
\eqref{eq:axis-trace-amplitude} in the radial variable $s$.  We separate the fixed radial weight from the
time-dependent multiplier $\beta_{\cusp}$ in the averaged moment \eqref{eq:axis-moment}.  Hence
$M_{\rm ax}=1$ when $\beta_{\cusp}\equiv1$, and $M_{\rm ax}$ measures the multiplicative change in the
leading axial strain.  The integral in \eqref{eq:axis-moment} is finite because
$\gamma>\alpha+\tfrac52$; at $t=0$, we have $\beta_{\cusp}(s,0,0)\equiv1$ and hence $M_{\rm ax}(0)=1$.

\runinhead{Step 2: The comparison between $m(t)$ and $M_{\rm ax}(t)$.}
By \eqref{eq:modulation-def}, $m(t)=\rW_{\cusp}(t)/\mathcal W_{\cusp}(t)$.  Combining the
axis-trace identity \eqref{eq:axis-trace-amplitude} with the averaged axial moment \eqref{eq:axis-moment}, we
obtain the parallel strain expansions
\[
\mathcal W_{\cusp}(t)=-\Gamma J^{3\alpha-1}C_\rho^{(1)}(\alpha,\gamma)C_W^* +O\!\left(\Gamma J^{3\alpha-1}\varepsilon_{\rm mod}(J)\right),
\]
and
\begin{equation}
\rW_{\cusp}(t)=-\Gamma J^{3\alpha-1}C_\rho^{(1)}(\alpha,\gamma)C_W^*M_{\rm ax}(t)
+O\!\left(\Gamma J^{3\alpha-1}\varepsilon_{\rm mod}(J)\right).
\label{eq:rWcusp-expansion}
\end{equation}
That is,
\begin{equation}
\bigl| \rW_{\cusp}(t)-\mathcal W_{\cusp}(t)M_{\rm ax}(t) \bigr| \le C\,\Gamma J^{3\alpha-1}\varepsilon_{\rm mod}(J).
\label{eq:rWcusp-WUcusp-Max-error}
\end{equation}
For the derivation of \eqref{eq:rWcusp-expansion}, we substitute the axis-trace identity
\eqref{eq:axis-trace-amplitude} for the cusp-flow transported vorticity
\eqref{eq:Omega-cusp-def} into the Biot--Savart representation of $\rW_{\cusp}(t)$ and integrate
in $s$ as in \eqref{eq:axis-moment}.  The resulting integrand in the label polar angle
$\sigma\in[0,\tfrac\pi2]$ contains $M_{\rm ax}(\varphi_J(\sigma),t)$, where
\[
\varphi_J(\sigma):=\arctan(J^3\tan\sigma)
\]
is the pulled-back polar angle attached to $\sigma$, with $\varphi_J(0)=0$ and
$\varphi_J(\sigma)\to 0$ as $J\to 0$ for each fixed $\sigma<\tfrac\pi2$.  Replacing
$M_{\rm ax}(\varphi_J(\sigma),t)$ in this integrand by $M_{\rm ax}(0,t)=M_{\rm ax}(t)$ produces the
leading term of \eqref{eq:rWcusp-expansion}.  The error from this replacement is bounded using
\[
|M_{\rm ax}(\varphi_J(\sigma),t)-M_{\rm ax}(0,t)|\le C\,\varphi_J(\sigma)^\alpha ,
\]
which follows from the $C^\alpha$ regularity of $M_{\rm ax}$ in the cusp angular variable, together
with the angular integral estimate
\[
\int_0^\infty u^{2+\alpha}(1+u^2)^{-5/2} \min\{1,(J^3u)^\alpha\}\,\ud u \le C J^{3\alpha}.
\]
The resulting error in \eqref{eq:rWcusp-expansion} has size $C\Gamma J^{3\alpha-1}\,J^{3\alpha}$.
This is the $J^{3\alpha}$ term in $\delta_{\rm mod}$ within $\varepsilon_{\rm mod}(J)$ on the
right-hand side of \eqref{eq:rWcusp-WUcusp-Max-error}.  Since
by \eqref{eq:Wcusp-scaling}
$|\mathcal W_{\cusp}(t)|\ge c_W\Gamma J^{3\alpha-1}$, we obtain that
\begin{equation}
m(t)=M_{\rm ax}(t)+O(\varepsilon_{\rm mod}(J_{\cusp}(t))).
\label{eq:m-Max-comparison}
\end{equation}

\runinhead{Step 3: The bootstrap interval.}
We now improve the moment part of the scalar-modulation bootstrap.  We choose a preliminary clock size
$J_*\in(0,1]$ and suppose the cusp clock reaches $J_*$.  Let $t_*$ be the first time
with $J_{\cusp}(t_*)=J_*$. Lemma~\ref{lem:late-entry} yields $t_*=O(\Gamma^{-1})$, and the finite-clock
regularity estimates imply
\[
0<c_*\le M_{\rm ax}(t_*)\le C_*<\infty,
\]
with constants depending only on $\alpha,\gamma$ and the chosen value of $J_*$.  Indeed, on the compact clock
range $J_{\cusp}(t)\in[J_*,1]$, the finite-clock regularity estimates for the exact cusp map, in particular
\eqref{eq:finite-clock-smooth-cusp-C1} on the bounded core together with the algebraic tail control
outside it, provide
\[
\|D\phi_{\cusp}\|+\|D\phi_{\cusp}^{-1}\|\le C_{\rm ent}(J_*) .
\]
Thus $J_{\twoD}(Y_{s,0},t)$ and $J_{\twoD}(Y_{s,0},t)^{-1}$ are bounded above by a constant depending only on
$\alpha,\gamma,J_*$.  Combining these
$C^1$ bounds with the change of variables \eqref{eq:modulation-collapse-labels} for $Y_{s,0}$
yields
\[
c_*^{(1)}\,s\le \rho(Y_{s,0})\le C_*^{(1)}\,s , \qquad c_*^{(2)}\le A_{\rm ax}(s,t)\le C_*^{(2)} ,
\]
where $A_{\rm ax}$ is the axis quotient \eqref{eq:Aax-def} and the constants
$c_*^{(1)},C_*^{(1)},c_*^{(2)},C_*^{(2)}>0$ depend only on $\alpha,\gamma,J_*$.  Substituting these bounds and
$\Upsilon(0)=1$ into the axis-trace identity \eqref{eq:axis-trace-amplitude} at $\varphi=0$ yields
\[
0<c_\beta\le\beta_{\cusp}(s,0,t_*)\le C_\beta<\infty
\]
with $c_\beta,C_\beta$ depending only on $\alpha,\gamma,J_*$.  Since $\gamma>\alpha+\tfrac52$, the weight
$s^{\alpha-1}\mathcal F(s)$ is integrable on $(0,\infty)$, so the moment integral \eqref{eq:axis-moment}
satisfies $0<c_*\le M_{\rm ax}(t_*)\le C_*<\infty$.

We consider the maximal interval beginning at $t_*$ on which
\begin{equation}
\tfrac12c_*\le M_{\rm ax}(t)\le 2C_*.
\label{eq:Max-bootstrap}
\end{equation}
After choosing $R_{\tail}$ so that $\varepsilon_{\smooth}$ is sufficiently small and then decreasing
$J_*$ if necessary, \eqref{eq:m-Max-comparison} implies on this interval that
\[
\tfrac14c_*\le m(t)\le 3C_*.
\]
This is a strict improvement of the scalar-multiplier bootstrap
\eqref{eq:localized-modulation-m-bootstrap} on the same interval.
Thus $J_{\cusp}$ is strictly decreasing there and the logarithmic clock variable
$\ell:=-\log J_{\cusp}$ is well-defined.  For the next estimate only, we define the rescaled flat axial strain and axial
velocity generated by the
flat cusp velocity $U_{\cusp}$ \eqref{eq:U-cusp-label} in the variable $\zeta=J^{-2}z$ by
\[
\mathsf W_t(\zeta):= \Gamma^{-1}J^{1-3\alpha}\p_z(U_{\cusp})_z(0,J^2\zeta,t),
\qquad \mathsf U_t(\zeta):=\int_0^\zeta\mathsf W_t(\eta)\,d\eta .
\]
These coincide with the functions defined later in
\eqref{eq:axis-profile-def}; at this point they are used only as time-$t$ axis traces.

\runinhead{Step 4: Logarithmic variation of the axis-trace multiplier.}
To improve the  bound for $M_{\rm ax}$  to \eqref{eq:modulation-moment-improvement}, we estimate the time
derivative of $\log\beta_{\cusp}(s,0,t)$ along the family $Y_{s,0}(t)$ determined by
\eqref{eq:modulation-collapse-labels}.  The starting point is the axis-trace identity
\eqref{eq:axis-trace-amplitude}, which determines $\beta_{\cusp}(s,0,t)$ as the limit of its
right-hand side once the common power $r_{s,\varphi}^\alpha$ has been removed and $\varphi\downarrow0$.
To express this limit,  we introduce the axis quotient
\begin{equation}
A_{\rm ax}(s,t) :=\lim_{\varphi\downarrow0} \tfrac{R(Y_{s,\varphi})}{J_{\cusp}(t)r_{s,\varphi}} =
\lim_{\varphi\downarrow0} \tfrac{\rho(Y_{s,\varphi})\sin\sigma(Y_{s,\varphi})}{J_{\cusp}(t)r_{s,\varphi}},
\label{eq:Aax-def}
\end{equation}
which compares the radial label coordinate $R(Y_{s,\varphi})$ with
$J_{\cusp}(t)r(\phi_{\cusp}(Y_{s,\varphi},t))$, the hyperbolically rescaled radial coordinate of the image
in \eqref{eq:modulation-collapse-labels}, as $\varphi\downarrow0$.  Taking the logarithm of
\eqref{eq:axis-trace-amplitude} for $0<\varphi\le\sigma_{\cut}$ and letting $\varphi\downarrow0$
yields the four-term decomposition
\begin{equation}
\log\beta_{\cusp}(s,0,t)=\log\tfrac{J_{\cusp}}{J_{\twoD}(Y_{s,0},t)}+ \alpha\log A_{\rm ax}(s,t)
+ \log\tfrac{\mathcal F(\rho(Y_{s,0}))}{\mathcal F(s)}+ \log\tfrac{\Upsilon(\sigma(Y_{s,0}))}{\Upsilon(0)} .
\label{eq:beta-cusp-log-decomposition}
\end{equation}

Each of the four terms in \eqref{eq:beta-cusp-log-decomposition} vanishes identically at $t=0$: there
$J_{\twoD}(Y_{s,0},0) = J_{\cusp}(0) = 1$, $A_{\rm ax}(s,0) = 1$, $\rho(Y_{s,0}(0)) = s$, and
$\sigma(Y_{s,0}(0)) = 0$, in agreement with the identity $\beta_{\cusp}(s,0,0)\equiv1$ from Step~1.  For
$t>0$, the four time-derivatives measure how the cusp flow distorts $J_{\twoD}$, $A_{\rm ax}$,
$\rho(Y_{s,0})$, and $\sigma(Y_{s,0})$ along the
moving family $Y_{s,0}(t)$.  We show next that each time-derivative is smaller than the singular clock rate
$\Gamma J^{3\alpha-1}$ by a common positive power of $J_{\cusp}$: on the bootstrap interval
\eqref{eq:Max-bootstrap},
\begin{subequations} 
\label{eq:beta-log-term-bounds}
\begin{align}
\big|
\p_t\log\tfrac{J_{\cusp}}{J_{\twoD}(Y_{s,0},t)}
\big|
&\le
C\Gamma J^{3\alpha-1}
\bigl(J^{3\alpha}+J^{1-3\alpha}\bigr)(1+s)^{-1-p_{\tail}}, \label{eq:beta-log-term-bounds-a} \\
\big|\p_t\log A_{\rm ax}(s,t)\big|
&\le
C\Gamma J^{3\alpha-1}
\bigl(J^{3\alpha}+J^{1-3\alpha}\bigr)(1+s)^{-1-p_{\tail}},\label{eq:beta-log-term-bounds-b} \\
\big|
\p_t\log\tfrac{\mathcal F(\rho(Y_{s,0}))}{\mathcal F(s)}
\big|
&\le
C\Gamma J^{3\alpha-1}
\bigl(J^{3\alpha}+J^{1-3\alpha}\bigr)(1+s)^{-1-p_{\tail}},\label{eq:beta-log-term-bounds-c} \\
\big|
\p_t\log\tfrac{\Upsilon(\sigma(Y_{s,0}))}{\Upsilon(0)}
\big|
&\le
C\Gamma J^{3\alpha-1}
\bigl(J^{3\alpha}+J^{1-3\alpha}\bigr)(1+s)^{-1-p_{\tail}},
\label{eq:beta-log-term-bounds-d}
\end{align}
\end{subequations} 
where $J=J_{\cusp}(t)$ and $p_{\tail}:=\gamma-\alpha-\tfrac52>0$.  In Step~4, all time derivatives are taken
along the moving family $Y_{s,0}(t)$, that is, at fixed $s$ and $\varphi=0$ in
\eqref{eq:modulation-collapse-labels}.

We first estimate the label velocity $\p_t Y_{s,0}(t)$.  The substitution $z=J^2 s$ that defines the
family $Y_{s,0}(t)$ removes the linear stagnation field $\mathsf W_t(0)z$ from the cusp velocity along
the symmetry axis.  Differentiating \eqref{eq:modulation-collapse-labels} at $\varphi=0$ in $t$ shows
that $\p_t Y_{s,0}(t)$ is therefore driven only by the nonlinear remainder.  In the rescaled
$\zeta$-coordinate this remainder equals
\[
m(t)\Gamma J^{3\alpha-1}\bigl(\mathsf U_t(s)-\mathsf W_t(0)s\bigr) + J^{-2}(V_{\err})_z(0,J^2 s,t),
\]
where $\mathsf W_t,\mathsf U_t$ are the axis traces from Step~3 and $V_{\err}$ is the cusp-error velocity
from \eqref{eq:Verr-def}; its $s$-derivative equals
\[
m(t)\Gamma J^{3\alpha-1}\bigl(\mathsf W_t(s)-\mathsf W_t(0)\bigr) + (\p_z V_{\err})_z(0,J^2 s,t).
\]
The cone-local $C^\alpha$ estimate \eqref{eq:Ucusp-grad-scale-local} for $U_{\cusp}$, the
cusp-error bootstrap \eqref{eq:localized-cusp-error-large-bootstrap}, and the algebraic radial tail
$\mathcal F$ from \eqref{eq:radial-tail-F-def} together yield
\begin{equation}
\bigl|\p_t Y_{s,0}(t)\bigr|_{\zeta} + \bigl|\p_t\log\p_s Y_{s,0}(t)\bigr|
\le C\,\Gamma J^{3\alpha-1}\bigl(J^{3\alpha}+J^{1-3\alpha}\bigr)(1+s)^{-1-p_{\tail}},
\label{eq:Ys0-motion-bound}
\end{equation}
where $|\cdot|_{\zeta}$ denotes the corresponding norm in the $\zeta$-coordinate.

We now bound each of the four time-derivatives in \eqref{eq:beta-log-term-bounds} using
\eqref{eq:Ys0-motion-bound}.

\rruninhead{Step 4a: the bound \eqref{eq:beta-log-term-bounds-a}} The axisymmetric Jacobian identity
\eqref{eq:Jac-Identity} yields
\[
\tfrac{d}{dt}\log\tfrac{J_{\cusp}}{J_{\twoD}(Y_{s,0},t)} = \tfrac12\rW_{\cusp}(t) - \tfrac12\p_z(V_{\cusp})_z(\phi_{\cusp}(Y_{s,0},t),t).
\]
At $s=0$ the two strains coincide by the definition of $\rW_{\cusp}$ in \eqref{eq:modulation-def}.
For $s>0$, their difference is bounded by \eqref{eq:Ucusp-grad-scale-local}, the cusp-error bound
\eqref{eq:localized-cusp-error-large-bootstrap}, and the decomposition \eqref{eq:Verr-def};
combined with \eqref{eq:Ys0-motion-bound}, this proves \eqref{eq:beta-log-term-bounds-a}.

\rruninhead{Step 4b: the bound \eqref{eq:beta-log-term-bounds-b}} Differentiating the quotient defining $A_{\rm ax}$ in
\eqref{eq:Aax-def}
produces the same difference of strains, since $J_{\cusp}(t)\,r_{s,\varphi}$ in the denominator
subtracts the radial-axial linear part of the cusp velocity at the stagnation point before the limit
$\varphi\downarrow0$ is taken.  The remaining radial variation is controlled by
\eqref{eq:Ucusp-grad-scale-local} and \eqref{eq:localized-cusp-error-large-bootstrap},
proving the second line of \eqref{eq:beta-log-term-bounds}.

\rruninhead{Step 4c: the bound \eqref{eq:beta-log-term-bounds-c}} The identity
\[
\bigl|\p_\rho \log\mathcal F(\rho)\bigr| = \gamma\,\tfrac{\rho}{1+\rho^2} \le C(1+\rho)^{-1},
\]
combined with the axial flow map comparison $\rho(Y_{s,0})\asymp s$ and the motion bound
\eqref{eq:Ys0-motion-bound}, yields the third line.  The multiplier $(1+s)^{-1-p_{\tail}}$ is the
same algebraic decay already present in \eqref{eq:Ys0-motion-bound}.

\rruninhead{Step 4d: the bound \eqref{eq:beta-log-term-bounds-d}} The label $Y_{s,0}(t)$ lies on the positive
symmetry axis, because the cusp flow preserves the axis and
\eqref{eq:modulation-collapse-labels} places its image on that axis.  Hence
$\sigma(Y_{s,0}(t))=0$, and since $\Upsilon(0)=1$ the fourth logarithmic term in
\eqref{eq:beta-cusp-log-decomposition} is identically zero.  This proves the fourth line of
\eqref{eq:beta-log-term-bounds}.

\runinhead{Step 5: Closing the bootstrap.}
To convert the time-derivative bound \eqref{eq:beta-log-term-bounds} into an integrated bound on
$\log\beta_{\cusp}(s,0,t)$, we change variable to the logarithmic clock $\ell:=-\log J_{\cusp}(t)$.  The
cusp-clock identity $\dot J_{\cusp}(t)=\tfrac12 J_{\cusp}(t)\rW_{\cusp}(t)$, combined with the strain
estimate \eqref{eq:Wcusp-scaling} and the modulation lower bound $m(t)\ge \tfrac14 c_*$ from Step~3,
yields
\[
\p_t\ell = -\dot J_{\cusp}(t)/J_{\cusp}(t) = -\tfrac12\rW_{\cusp}(t) \ge c\,\Gamma\,J_{\cusp}(t)^{3\alpha-1} .
\]
Dividing each line of \eqref{eq:beta-log-term-bounds} by $\p_t\ell$ yields
\begin{equation}
\bigl|\p_\ell\log\beta_{\cusp}(s,0,t)\bigr| \le C\bigl(J_{\cusp}(t)^{3\alpha}+J_{\cusp}(t)^{1-3\alpha}\bigr)(1+s)^{-1-p_{\tail}} .
\label{eq:beta-cusp-log-ell-derivative}
\end{equation}

Integrating \eqref{eq:beta-cusp-log-ell-derivative} in $\ell$ from $t_*$, where $J_{\cusp}(t_*)=J_*$,
to any $t$ in the bootstrap interval yields
\begin{equation}
\bigl|\log\beta_{\cusp}(s,0,t)-\log\beta_{\cusp}(s,0,t_*)\bigr| \le C\,J_*^{\min\{3\alpha,1-3\alpha\}}\,(1+s)^{-1-p_{\tail}} .
\label{eq:beta-cusp-integrated-bound}
\end{equation}
The two-sided bound $0<c_\beta\le\beta_{\cusp}(s,0,t_*)\le C_\beta<\infty$ from Step~3,
\eqref{eq:beta-cusp-integrated-bound}, and the integrability of $s^{\alpha-1}\mathcal F(s)(1+s)^{-1-p_{\tail}}$
on $(0,\infty)$ (which holds since $\gamma>\alpha+\tfrac52$) together imply, after multiplying by
$s^{\alpha-1}\mathcal F(s)/C_\rho^{(1)}(\alpha,\gamma)$ and integrating in $s$ via the definition
\eqref{eq:axis-moment} of $M_{\rm ax}$,
\[
\bigl|M_{\rm ax}(t)-M_{\rm ax}(t_*)\bigr| \le C'\,J_*^{\min\{3\alpha,1-3\alpha\}} ,
\]
with $C'$ depending only on the fixed parameters.  Choosing $J_*$ small enough so that
$C'\,J_*^{\min\{3\alpha,1-3\alpha\}}\le\min\{\tfrac14 c_*,\tfrac12 C_*\}$ and using the entry bound
$c_*\le M_{\rm ax}(t_*)\le C_*$ from Step~3, we obtain
\[
\tfrac34 c_*\le M_{\rm ax}(t)\le \tfrac32 C_* .
\]
This is a strict improvement of \eqref{eq:Max-bootstrap}: by continuity, the maximal interval on
which \eqref{eq:Max-bootstrap} holds is therefore both open and closed in
$\{t:J_{\cusp}(t)\le\mathfrak J_{\mathrm{mod}}\}$, and so coincides with this small-clock interval.  This
proves \eqref{eq:modulation-moment-improvement}.

To complete the proof of \eqref{eq:modulation-bounds}, we substitute
\eqref{eq:modulation-moment-improvement} into the comparison \eqref{eq:m-Max-comparison} and
shrink $\mathfrak J_{\mathrm{mod}}\le\min\{J_*,\mathfrak J_{\mathrm{strain}}\}$ once more so that
$\varepsilon_{\rm mod}(J_{\cusp}(t))\le\min\{\tfrac12 c_*,\tfrac32 C_*\}$ on
$0<J_{\cusp}(t)\le\mathfrak J_{\mathrm{mod}}$.  This yields
\[
\tfrac14 c_*\le m(t)\le 3C_* ,
\]
which is \eqref{eq:modulation-bounds} with $c_m=\tfrac14 c_*$ and $C_m=3C_*$.  Since
$\mathfrak J_{\mathrm{mod}}\le\mathfrak J_{\mathrm{strain}}$, the strain estimate
\eqref{eq:Wcusp-scaling} holds throughout this interval.
\end{proof}

The previous lemma closes the scalar-modulation bootstrap.  We now close the cusp-clock bootstrap
\eqref{eq:localized-clock-bootstrap}.  That bootstrap assumption is needed throughout the small-clock
argument so that time integrals can be rewritten as integrals in $J_{\cusp}$.  Once the modulation bounds have
been improved, the clock identity
\[
\dot J_{\cusp}(t) = \tfrac12 J_{\cusp}(t)m(t)\mathcal W_{\cusp}(t)
\]
combines with the flat cusp axial strain estimate \eqref{eq:Wcusp-scaling} to produce a sharper
two-sided bound for $-\dot J_{\cusp}$.  The next lemma proves this estimate; with the strict-margin choice
of clock-bootstrap constants described in Section~\ref{sec:small-clock-bootstraps}, it closes the
cusp-clock bootstrap.

\begin{lemma}[Improvement of the cusp-clock bootstrap]
\label{lem:Jdot-two-sided-aux}
There exist constants $0<c_1\le C_1<\infty$, depending only on
$\alpha,\gamma,\sigma_{\inn},\sigma_*$, such that whenever
$J_{\cusp}(t)\le \mathfrak J_{\mathrm{mod}}$,
\begin{equation}
c_1\,\Gamma\,J_{\cusp}(t)^{3\alpha} \le -\dot J_{\cusp}(t) \le C_1\,\Gamma\,J_{\cusp}(t)^{3\alpha}.
\label{eq:Jdot-two-sided}
\end{equation}
\end{lemma}

\begin{remark}[Closing the clock bootstrap]
The constants in the cusp-clock bootstrap are chosen in
Section~\ref{sec:small-clock-bootstraps} so that
\[
0<c_{\rm clk}<c_1\le C_1<C_{\rm clk},
\]
Then \eqref{eq:Jdot-two-sided} is a strict improvement of the cusp-clock bootstrap assumption
\eqref{eq:localized-clock-bootstrap}.  The same estimate implies that $J_{\cusp}(t)$ is strictly
decreasing once $J_{\cusp}(t)\le\mathfrak J_{\mathrm{mod}}$.
\end{remark}

\begin{proof}[Proof of Lemma~\ref{lem:Jdot-two-sided-aux}]
We fix a time with $J_{\cusp}(t)\le\mathfrak J_{\mathrm{mod}}$.  At this time the flat cusp axial strain estimate
\eqref{eq:Wcusp-scaling} and the improved modulation bound
\eqref{eq:modulation-bounds} are available.

The cusp clock satisfies
\begin{equation}
\dot J_{\cusp}(t) = \tfrac12\,J_{\cusp}(t)\,\rW_{\cusp}(t) = \tfrac12\,J_{\cusp}(t)\,m(t)\,\mathcal W_{\cusp}(t).
\label{eq:cusp-clock-identity-in-Jdot-proof}
\end{equation}
Indeed,
$\dot J_{\cusp}/J_{\cusp}
=\p_r(V_{\cusp})_r(0,t)+\p_z(V_{\cusp})_z(0,t)$, while axisymmetric incompressibility implies
$2\p_r(V_{\cusp})_r(0,t)+\p_z(V_{\cusp})_z(0,t)=0$ on the symmetry axis.  The second equality in
\eqref{eq:cusp-clock-identity-in-Jdot-proof} is the definition of $m(t)$ in
\eqref{eq:modulation-def}.
Since $J_{\cusp}(t)\le\mathfrak J_{\mathrm{mod}}\le \mathfrak J_{\mathrm{strain}}$,
\eqref{eq:Wcusp-scaling} gives
\[
\mathcal W_{\cusp}(t)<0, \qquad c_W\,\Gamma J_{\cusp}(t)^{3\alpha-1} \le |\mathcal W_{\cusp}(t)| \le C_W\,\Gamma J_{\cusp}(t)^{3\alpha-1},
\]
and \eqref{eq:modulation-bounds} yields
\[
c_m\le m(t)\le C_m.
\]
Therefore
\[
\tfrac12 c_m c_W\,\Gamma J_{\cusp}(t)^{3\alpha} \le -\dot J_{\cusp}(t) \le \tfrac12 C_m C_W\,\Gamma J_{\cusp}(t)^{3\alpha}.
\]
This is \eqref{eq:Jdot-two-sided}, with
$c_1:=\tfrac12c_m c_W$ and $C_1:=\tfrac12C_m C_W$.
\end{proof}

The smooth velocity also contributes an axial strain at the stagnation point.  We denote it by
\[
\rW_{\smooth}(t):=\p_z(u_{\smooth})_z(0,t).
\]
This is an order-one smooth contribution to the axial strain, and it enters the clock decomposition through
$J_{\smooth}$.  The next lemma proves this bound and, using the small-clock time length from
Lemma~\ref{lem:Jdot-two-sided-aux}, also bounds the smooth clock $J_{\smooth}$.

\begin{lemma}[Bounded smooth clock and smooth axial strain]
\label{lem:smooth-clock-bounded}
There exist constants $0<c_{\smooth}\le C_{\smooth}<\infty$, depending only on the fixed parameters, such
that on every interval on which $J_{\cusp}(t)\le \mathfrak J_{\mathrm{mod}}$,
\begin{equation}
|\rW_{\smooth}(t)|\le C_{\smooth}\,\Gamma, \qquad 0<c_{\smooth}\le J_{\smooth}(t)\le C_{\smooth}<\infty.
\label{eq:Jsmooth-bdd}
\end{equation}
\end{lemma}

\begin{proof}[Proof of Lemma~\ref{lem:smooth-clock-bounded}]
We fix a connected time interval $I$ on which $J_{\cusp}\le\mathfrak J_{\mathrm{mod}}$ and a time
$t\in I$.

The smooth strain bound follows directly from the definition of $u_{\smooth}$ in
\eqref{eq:smooth-velocity-def}.  Lemma~\ref{lem:JtwoD-tail-bdd}, proved in
the next subsection, yields after setting $R_0=1$
$\|\nabla u_{\smooth}(\cdot,s)\|_{L^\infty(B_1)}\le C\Gamma$ for every $s\in I$.
Therefore, by the definition of $\rW_{\smooth}$,
\[
|\rW_{\smooth}(t)| = |\p_z(u_{\smooth})_z(0,t)| \le \|\nabla u_{\smooth}(\cdot,t)\|_{L^\infty(B_1)} \le C\Gamma,
\]
which proves the first estimate in \eqref{eq:Jsmooth-bdd}.

It remains to bound $J_{\smooth}$.  From \eqref{eq:smooth-cusp-clock-def},
axisymmetry of $\phi_{\smooth}$, and the flow equation \eqref{eq:smooth-flow-def},
\[
\tfrac{\dot J_{\smooth}(t)}{J_{\smooth}(t)} = \p_r(u_{\smooth})_r(0,t)+\p_z(u_{\smooth})_z(0,t).
\]
Since $u_{\smooth}$ is divergence-free and axisymmetric without swirl,
$2\p_r(u_{\smooth})_r(0,t)+\p_z(u_{\smooth})_z(0,t)=0$,
and hence
\begin{equation}
\dot J_{\smooth}(t)=\tfrac12\,J_{\smooth}(t)\,\rW_{\smooth}(t).
\label{eq:Jsmooth-clock-ode}
\end{equation}
By Lemma~\ref{lem:Jdot-two-sided-aux}, $J_{\cusp}$ is strictly decreasing on the small-clock regime.  Let
$t_0$ be the entry time of the connected small-clock component containing $I$, so that
$J_{\cusp}(t_0)=\mathfrak J_{\mathrm{mod}}$.  The finite-clock estimate
\eqref{eq:finite-clock-smooth-cusp-C1}, with
$\mathfrak J_{\mathrm{finite}}=\mathfrak J_{\mathrm{mod}}$, bounds
$D\phi_{\smooth}(0,t_0)$ and $D\phi_{\smooth}(0,t_0)^{-1}$.  By the determinant definition in
\eqref{eq:smooth-cusp-clock-def},
\[
0<c_0\le J_{\smooth}(t_0)\le C_0<\infty .
\]
Integrating \eqref{eq:Jsmooth-clock-ode} from $t_0$ to $t$ and using
$|\rW_{\smooth}|\le C\Gamma$ yields
\[
\exp\!\bigl(-C\Gamma(t-t_0)\bigr) \le \tfrac{J_{\smooth}(t)}{J_{\smooth}(t_0)} \le \exp\!\bigl(C\Gamma(t-t_0)\bigr).
\]
It remains to bound $t-t_0$.  By \eqref{eq:Jdot-two-sided},
\[
\tfrac{d}{dt}J_{\cusp}(t)^{1-3\alpha} = (1-3\alpha)J_{\cusp}(t)^{-3\alpha}\dot J_{\cusp}(t) \le -c\Gamma.
\]
Thus $t-t_0\le C\Gamma^{-1}$ while $0<J_{\cusp}(t)\le J_{\cusp}(t_0)=\mathfrak J_{\mathrm{mod}}$. The exponential bound
for $J_{\smooth}$ is therefore bounded above and below by constants depending only on the fixed parameters,
which proves the second estimate in \eqref{eq:Jsmooth-bdd}.
\end{proof}

The velocity identity \eqref{eq:u-decomp} already decomposes $u$ into smooth, cusp, and error terms.
For the Riccati argument, however, we need the corresponding identity for the stagnation-point axial strain.
This requires a separate check because $u_{\cusp}$ and $u_{\err}$ are physical Eulerian velocities obtained
by pushing forward cusp-coordinate fields through the smooth flow.  The next lemma shows that this push-forward
does not change the axial diagonal entry at the origin.  Thus the scalar modulation
\eqref{eq:modulation-def} makes the physical cusp term contribute exactly $\rW_{\cusp}$, while the
physical error term contributes no axial strain at the origin.

\begin{lemma}[Stagnation-point axial strain identity]
\label{lem:on-axis-strain-splitting}
Let $u_{\cusp}$ and $u_{\err}$ be the physical Eulerian velocities defined in
\eqref{eq:Verr-def}--\eqref{eq:u-decomp}.  Then, at every time for which these fields are
defined,
\[
\p_z(u_{\cusp})_z(0,t)=\rW_{\cusp}(t), \qquad \p_z(u_{\err})_z(0,t)=0.
\]
Consequently, with $\rW_0(t)=\p_z u_z(0,0,t)$ from \eqref{eq:rW0-Pi0-J0-defs} and with
$\rW_{\smooth}$ defined above,
\begin{equation}
\rW_0(t)=\rW_{\smooth}(t)+\rW_{\cusp}(t).
\label{eq:on-axis-strain-splitting}
\end{equation}
\end{lemma}

\begin{proof}[Proof of Lemma~\ref{lem:on-axis-strain-splitting}]
The smooth velocity has the same axisymmetry and axial parity as the full solution.  Hence
$u_{\smooth}(0,t)=0$, the smooth flow preserves the symmetry axis and fixes the origin, and
$D\phi_{\smooth}(0,t)$ is diagonal in the meridional variables $(r,z)$.  We first note what this implies for
the axial derivative of a pushed-forward field.

Let $w$ be an axisymmetric no-swirl field with $w(0,t)=0$.  Since
$(\phi_{\smooth})_*w(x,t)=D\phi_{\smooth}(X,t)w(X,t)$ with $x=\phi_{\smooth}(X,t)$, differentiating at
$x=0$ yields
\[
\nabla(\phi_{\smooth})_*w(0,t) = D\phi_{\smooth}(0,t)\,\nabla w(0,t)\,D\phi_{\smooth}(0,t)^{-1};
\]
the $D^2\phi_{\smooth}$ term vanishes because $w(0,t)=0$.  Since the conjugating matrix is diagonal, this
conjugation leaves the axial diagonal entry unchanged.

The transported field $U_{\cusp}$ vanishes at the origin by the same axial parity, and
$V_{\cusp}(0,t)=0$ because the cusp flow fixes the origin.  Thus $V_{\err}(0,t)=0$ by
\eqref{eq:Verr-def}.  Applying the preceding observation to $w=m(t)U_{\cusp}(\cdot,t)$ and to
$w=V_{\err}$ gives
\[
\p_z(u_{\cusp})_z(0,t)=m(t)\mathcal W_{\cusp}(t)=\rW_{\cusp}(t), \qquad \p_z(u_{\err})_z(0,t)=\p_z(V_{\err})_z(0,t).
\]
Finally, \eqref{eq:Verr-def} and \eqref{eq:modulation-def} imply
\[
\p_z(V_{\err})_z(0,t) = \rW_{\cusp}(t)-m(t)\mathcal W_{\cusp}(t) = 0.
\]
Adding the smooth contribution from \eqref{eq:u-decomp} proves
\eqref{eq:on-axis-strain-splitting}.
\end{proof}

\subsection{Analysis of the smooth velocity $u_{\smooth}$}

To distinguish the scalar cusp clock from the full label-dependent meridional Jacobian, we write
\[
J_{\twoD}(Y,t):=\det\nabla_{(R,Z)}(\phi_r,\phi_z)(Y,t), \qquad J(t)=J_{\twoD}(0,0,t).
\]
Then the transport identity \eqref{eq:vort-identity} becomes
\begin{equation}
\omega_\theta(\phi(Y,t),t) = J_{\twoD}(Y,t)^{-1}\,\omega_{\theta,0}(Y).
\label{eq:vort-identity-JtwoD}
\end{equation}

\begin{lemma}[Fixed-ball bounds for the smooth velocity $u_{\smooth}$]
\label{lem:JtwoD-tail-bdd}
Fix $R_0<\infty$ and assume $R_{\tail}\ge8R_0$.  There is a constant
$C_{R_0,\alpha,\gamma}$, independent of $t$, $\Gamma$, $R_{\tail}$, $\nu$, and $\eta$, such that for every
datum in the admissible class $\mathcal A_{\alpha,\gamma}(\nu,\eta)$ with $0\le\nu\le1$, the velocity
$u_{\smooth}$ defined by the far-field cutoff in \eqref{eq:smooth-velocity-def} satisfies, at
every time $t$ for which the corresponding Euler solution is regular,
\begin{subequations}
\begin{align}
\|\nabla^k u_{\smooth}(\cdot,t)\|_{L^\infty(B_{2R_0})}
&\le
C_{R_0,\alpha,\gamma}\Gamma R_{\tail}^{-1-k},
\qquad k=0,1,2,
\label{eq:current-far-usmooth-derivatives}\\
[\nabla^2u_{\smooth}(\cdot,t)]_{C^\alpha(B_{2R_0})}
&\le
C_{R_0,\alpha,\gamma}\Gamma R_{\tail}^{-3-\alpha}.
\label{eq:current-far-usmooth-holder}
\end{align}
\end{subequations}
In particular,
\begin{equation}
\|u_{\smooth}(\cdot,t)\|_{C^{2,\alpha}(B_{2R_0})} \le C_{R_0,\alpha,\gamma}\Gamma R_{\tail}^{-1}.
\label{eq:current-far-usmooth-C2a}
\end{equation}
\end{lemma}

\begin{proof}[Proof of Lemma~\ref{lem:JtwoD-tail-bdd}]
Let $x\in B_{2R_0}$.  On the support of the cutoff in \eqref{eq:smooth-velocity-def}, we have
$|\phi(Y,t)|\ge R_{\tail}$.  Since $R_{\tail}\ge8R_0$,
\[
|x-\phi(Y,t)| \ge |\phi(Y,t)|-|x| \ge |\phi(Y,t)|-2R_0 \ge \tfrac12|\phi(Y,t)|.
\]
The Biot--Savart kernel is therefore evaluated away from its singularity.  Differentiation in $x$ acts only
on this kernel, since the cutoff is a function of the integration variable, and
\begin{equation}
|\nabla_x^kK(x,\phi(Y,t))| \le C_k|\phi(Y,t)|^{-2-k}, \qquad k=0,1,2.
\label{eq:current-far-kernel-derivative}
\end{equation}
By \eqref{eq:Jac-Identity},
\[
J_{\twoD}(Y,t)^{-1}=\tfrac{\phi_r(Y,t)}{R(Y)} \le \tfrac{|\phi(Y,t)|}{R(Y)}.
\]
Consequently, on the support of the far-field cutoff,
\begin{equation}
|\phi(Y,t)|^{-2-k}J_{\twoD}(Y,t)^{-1}|\omega_{\theta,0}(Y)| \le R_{\tail}^{-1-k}\tfrac{|\omega_{\theta,0}(Y)|}{R(Y)}, \qquad k=0,1,2.
\label{eq:current-far-jacobian-weight}
\end{equation}

The moment on the right-hand side is finite uniformly for the admissible class.  Indeed, \eqref{eq:vort0},
\eqref{eq:Theta-star-def}, $0\le\Upsilon\le1$, and
$\Theta=\Theta^*(1+h)$ with $\|h\|_{L^\infty}\le\nu\le1$ imply
\[
\tfrac{|\omega_{\theta,0}(Y)|}{R(Y)} \le C\Gamma R(Y)^{\alpha-1}(1+|Y|^2)^{-\gamma/2}.
\]
Using the cylindrical volume element $R\,dR\,dZ\,d\theta$ and then polar coordinates in the meridional
half-plane,
\begin{equation}
\int_{\mathbb R^3}\tfrac{|\omega_{\theta,0}(Y)|}{R(Y)}\,dY
\le C_{\alpha,\gamma}\Gamma \int_0^\infty \rho^{\alpha+1}(1+\rho^2)^{-\gamma/2}\,d\rho \le C_{\alpha,\gamma}\Gamma,
\label{eq:current-far-vorticity-over-r-moment}
\end{equation}
because $\gamma>\alpha+2$.

 Combining \eqref{eq:smooth-velocity-def}, \eqref{eq:current-far-kernel-derivative},
\eqref{eq:current-far-jacobian-weight}, and \eqref{eq:current-far-vorticity-over-r-moment} proves
\eqref{eq:current-far-usmooth-derivatives}.  The H\"older seminorm follows from the corresponding
difference-quotient kernel estimate: we bound the kernel difference in the observation variables first, and then
use the same far-field moment bound \eqref{eq:current-far-vorticity-over-r-moment}.
For $x,x'\in B_{2R_0}$, the same separation also implies
$|x-x'|\le4R_0\le\tfrac12|\phi(Y,t)|$ on the cutoff support, and hence
\[
|\nabla_x^2K(x,\phi(Y,t))-\nabla_x^2K(x',\phi(Y,t))| \le C|x-x'|^\alpha|\phi(Y,t)|^{-4-\alpha}, \qquad x,x'\in B_{2R_0}.
\]
Using \eqref{eq:Jac-Identity} once more,
\[
|\phi(Y,t)|^{-4-\alpha}J_{\twoD}(Y,t)^{-1}|\omega_{\theta,0}(Y)| \le R_{\tail}^{-3-\alpha}\tfrac{|\omega_{\theta,0}(Y)|}{R(Y)}.
\]
Together with \eqref{eq:current-far-vorticity-over-r-moment}, this proves
\eqref{eq:current-far-usmooth-holder}.  The bound
\eqref{eq:current-far-usmooth-C2a} follows immediately.
\end{proof}

The estimates \eqref{eq:current-far-usmooth-derivatives}--\eqref{eq:current-far-usmooth-C2a} are local in
space but uniform in time: on every fixed ball, $u_{\smooth}$ has size
$O(\Gamma R_{\tail}^{-1})$ in $C^{2,\alpha}$.  The time intervals on which we use the
smooth flow have length $O(\Gamma^{-1})$, so the accumulated deformation generated by
$u_{\smooth}$ is $O(R_{\tail}^{-1})$.  Thus, after the tail radius has been fixed large enough,
$\phi_{\smooth}$ is a near-identity diffeomorphism on the balls used in the cusp analysis.  The next
lemma states this consequence.  The velocity bound only names the resulting small parameter
$\varepsilon_{\smooth}$; the main conclusion is the near-identity control of $\phi_{\smooth}$.

\begin{lemma}[Near-identity deformation generated by the smooth velocity]
\label{lem:smooth-flow-small-deformation}
Fix $R_0<\infty$ and $C_T<\infty$.  After increasing the tail radius $R_{\tail}$ in
\eqref{eq:core-tail-domains}, we may choose
\[
\varepsilon_{\smooth}=\varepsilon_{\smooth}(R_{\tail};R_0,C_T,\alpha,\gamma) \qquad\text{with}\qquad \varepsilon_{\smooth}\to0 \quad\text{as
}R_{\tail}\to\infty,
\]
so that the following holds for every Euler solution in the $C^{1,\alpha}$ class considered here on a time
interval $[0,T]$.  Set $T_*:=\min\{T,C_T\Gamma^{-1}\}$.  Then
\begin{equation}
\sup_{0\le t\le T_*} \|u_{\smooth}(\cdot,t)\|_{C^{2,\alpha}(B_{2R_0})} \le \varepsilon_{\smooth}\Gamma.
\label{eq:usmooth-small-C2a}
\end{equation}
Moreover, if $\phi_{\smooth}$ is the flow of $u_{\smooth}$, then
$X\mapsto\phi_{\smooth}(X,t)$ is a $C^2$ diffeomorphism from $B_{R_0}$ onto its image, and for every
$X\in B_{R_0}$ and every $0\le t\le T_*$,
\begin{subequations}
\begin{align}
|\phi_{\smooth}(X,t)-X|
&\le C\varepsilon_{\smooth}|X|,
\label{eq:smooth-flow-position-small}\\
|D\phi_{\smooth}(X,t)-I|
+
|D\phi_{\smooth}(X,t)^{-1}-I|
&\le C\varepsilon_{\smooth},
\label{eq:smooth-flow-gradient-small}\\
|D^2\phi_{\smooth}(X,t)|
+
|D^2\phi_{\smooth}^{-1}(\phi_{\smooth}(X,t),t)|
&\le C\varepsilon_{\smooth}.
\label{eq:smooth-flow-second-gradient-small}
\end{align}
\end{subequations}
The constant $C$ depends only on $R_0,C_T,\alpha,\gamma$ and on the fixed cutoff functions, but is independent
of the small cusp clock.
\end{lemma}

\begin{proof}[Proof of Lemma~\ref{lem:smooth-flow-small-deformation}]
With $T_*:=\min\{T,C_T\Gamma^{-1}\}$, Lemma~\ref{lem:JtwoD-tail-bdd} yields
\[
\sup_{0\le t\le T_*}\|u_{\smooth}(\cdot,t)\|_{C^{2,\alpha}(B_{2R_0})} \le C_{R_0,\alpha,\gamma}\Gamma R_{\tail}^{-1}.
\]
We set
\[
\varepsilon_{\smooth}:=C_{R_0,\alpha,\gamma}R_{\tail}^{-1}.
\]
After increasing $R_{\tail}$ if necessary, we assume throughout the rest of the proof that
\[
C C_T\varepsilon_{\smooth}\le \log(\tfrac32), \qquad C\varepsilon_{\smooth}\le1,
\]
where $C$ is fixed large enough for the estimates in the proof.  This proves
\eqref{eq:usmooth-small-C2a}, and $\varepsilon_{\smooth}\to0$ as $R_{\tail}\to\infty$.

We now turn to the flow estimates.  Let $\Lambda(X,t):=\phi_{\smooth}(X,t)$.  By axisymmetry and the odd
symmetry across the plane $z=0$, the smooth velocity vanishes at the origin:
$u_{\smooth}(0,t)=0$.

The mean value theorem and \eqref{eq:usmooth-small-C2a} imply that
\[
|u_{\smooth}(X,t)|\le C\varepsilon_{\smooth}\Gamma |X|, \qquad X\in B_{2R_0}.
\]
We first prove that trajectories starting in $B_{R_0}$ remain in $B_{2R_0}$.  For fixed $X\in B_{R_0}$,
define
\[
\tau_X:= \sup\{\,\tau\in[0,T_*]:\ |\Lambda(X,s)|<2R_0 \text{ for every }0\le s<\tau\,\}.
\]
Continuity implies $\tau_X>0$.  For $0\le t<\tau_X$, the preceding bound implies
\[
|\Lambda(X,t)| \le |X|+\int_0^t C\varepsilon_{\smooth}\Gamma|\Lambda(X,s)|\,ds \le |X|\exp(C\varepsilon_{\smooth}\Gamma t) \le \tfrac32 |X|
\le \tfrac32 R_0.
\]
If $\tau_X<T_*$, then continuity implies
$|\Lambda(X,\tau_X)|\le \tfrac32R_0<2R_0$.  Hence there exists
$\delta\in(0,T_*-\tau_X)$ such that
$|\Lambda(X,s)|<2R_0$ for $\tau_X\le s\le \tau_X+\delta$.  This contradicts the definition of
$\tau_X$.  Hence
$\tau_X=T_*$, and $\Lambda(X,s)\in B_{2R_0}$ for every $0\le s\le T_*$.  Using the integral equation once more,
\[
|\Lambda(X,t)-X| \le \int_0^t C\varepsilon_{\smooth}\Gamma|\Lambda(X,s)|\,ds \le C\varepsilon_{\smooth}|X|,
\]
which proves \eqref{eq:smooth-flow-position-small}.

Since $u_{\smooth}$ is $C^2$ on $B_{2R_0}$ and the trajectories above remain in this ball, the map
$X\mapsto\Lambda(X,t)$ is a $C^2$ diffeomorphism from $B_{R_0}$ onto its image.  We now estimate its gradient.

The gradient $P(X,t):=D\Lambda(X,t)$ satisfies
\begin{equation}
\p_tP = (\nabla u_{\smooth})(\Lambda(X,t),t)P, \qquad P(X,0)=I.
\label{eq:smooth-flow-differential-ode}
\end{equation}
Using \eqref{eq:usmooth-small-C2a} and $t\le C_T\Gamma^{-1}$, we obtain
\[
|P(X,t)| \le \exp(C\varepsilon_{\smooth}),
\]
and
\[
|P(X,t)-I| \le \int_0^t |\nabla u_{\smooth}(\Lambda(X,s),s)|\,|P(X,s)|\,ds \le C\varepsilon_{\smooth}.
\]
The same variational equation for the inverse gradient, or the identity $P^{-1}-I=P^{-1}(I-P)$ together with $|P^{-1}|\le \exp(C\varepsilon_{\smooth})$, proves the
inverse-gradient bound in \eqref{eq:smooth-flow-gradient-small}.

Finally, $H(X,t):=D^2\Lambda(X,t)$ satisfies
\[
\p_t H = (\nabla^2u_{\smooth})(\Lambda(X,t),t)[P,P] + (\nabla u_{\smooth})(\Lambda(X,t),t)H, \qquad H(X,0)=0.
\]
The already proved bound for $P$, the estimate \eqref{eq:usmooth-small-C2a}, and
$t\le C_T\Gamma^{-1}$ imply
\[
|D^2\phi_{\smooth}(X,t)|=|H(X,t)|\le C\varepsilon_{\smooth}.
\]
The inverse map satisfies the identity
\[
D^2\phi_{\smooth}^{-1}(\Lambda(X,t),t)[\xi,\eta] = -P(X,t)^{-1} H(X,t)\bigl[P(X,t)^{-1}\xi,P(X,t)^{-1}\eta\bigr],
\]
which proves the inverse second-derivative estimate in
\eqref{eq:smooth-flow-second-gradient-small}.  This completes the proof.
\end{proof}

\subsection{Axial flow decomposition and axial-amplitude control}
\label{sec:axial-composition-profile-control}
Throughout this subsection, the time $t$ is fixed unless explicitly varied.

We now isolate the one-dimensional dynamics of the exact cusp map on the symmetry axis.  We write
\[
\phi_{\cusp}(0,Z,t)=(0,B_t(Z)), \qquad A_t(Z):=\p_R(\phi_{\cusp})_r(0,Z,t), \qquad J:=J_{\cusp}(t).
\]
The clock-scaled axial coordinate is
\[
\zeta=J^{-2}B_t(Z).
\]
The Section~\ref{sec:small-clock-bootstraps} axis-geometry bootstrap
\eqref{eq:current-axis-geometry} is a bound for the normalized radial and axial derivatives of
$\phi_{\cusp}$ in this $\zeta$-coordinate.  The next lemma is the exact one-dimensional reduction needed to
close that bootstrap: it proves the fixed-label identity \eqref{eq:zeta-current-transport} for
$\p_t(J^{-2}B_t(Z))$, derives the equations for the normalized derivatives, and identifies the
$V_{\err}$ error terms in those equations.

Assuming the axis-geometry bounds, the map $Z\mapsto J^{-2}B_t(Z)$ is strictly increasing on each axial label
interval under consideration.  We denote its inverse by $Z_t(\zeta)$ and set
\begin{equation}
q_t(\zeta):=J A_t(Z_t(\zeta)), \qquad b_t(\zeta):=J^{-2}B_t'(Z_t(\zeta)), \qquad B_t(Z_t(\zeta))=J^2\zeta,
\label{eq:axis-normalized-derivatives-def}
\end{equation}
where $J=J_{\cusp}(t)$.  We define the clock-normalized axial strain of the transported cusp field and its
anti-derivative by
\begin{equation}
\mathsf W_t(\zeta) := \Gamma^{-1}J^{1-3\alpha}\, \p_z(U_{\cusp})_z(0,J^2\zeta,t),
\qquad \mathsf U_t(\zeta) := \int_0^\zeta \mathsf W_t(\zeta')\,d\zeta'.
\label{eq:axis-profile-def}
\end{equation}
Thus $\mathsf W_t$ is the normalized axial strain of $U_{\cusp}$ on the symmetry axis, and
$\mathsf U_t$ is introduced because oddness implies the reconstruction
\[
(U_{\cusp})_z(0,J^2\zeta,t)=\Gamma J^{3\alpha+1}\mathsf U_t(\zeta).
\]
The one-dimensional velocity of the $\zeta$-coordinate of a fixed axial label is
\begin{equation}
\mathsf V_t(\zeta) := m(t)\Gamma J^{3\alpha-1} \bigl(\mathsf U_t(\zeta)-\mathsf W_t(0)\zeta\bigr) +J^{-2}(V_{\err})_z(0,J^2\zeta,t).
\label{eq:axis-current-zeta-velocity}
\end{equation}
The term $\mathsf W_t(0)\zeta$ is subtracted because the linear axial strain at the stagnation point is
already built into the normalization $J^{-2}B_t(Z)$.

We state the next lemma on an arbitrary compact $\zeta$-interval
$I_\zeta\subset I_{\err}$, where $I_{\err}$ is the fixed interval from
\eqref{eq:cusp-error-trace-interval}.  We use the same notation $I_\zeta$ in the following  two settings:
\begin{enumerate} 
\item In the first setting $I_\zeta$ denotes an origin-attached interval, so that 
$ 0\in I_\zeta\subset[0,\infty)$, and for a continuous function $F(\zeta,t)$ defined for $\zeta>0$, we define
$F(0,t):=\lim_{\zeta\downarrow0}F(\zeta,t)$.

\item In the second setting, $I_\zeta$ denotes an interval separated from the origin. For pressure Hessian localization, we take
\[
I_\zeta\in\{I_{\rm loc}^{\rm cur},I_{\rm buf}^{\rm cur}\},\qquad \operatorname{supp}\vartheta_\sharp\Subset I_{\rm loc}^{\rm cur}
\Subset I_{\rm buf}^{\rm cur}\Subset I_\sharp\Subset(0,\infty), \qquad \vartheta_\sharp\in C_c^\infty(I_\sharp)
\]
as fixed in \eqref{eq:pressure-localization-cutoff} and
\eqref{eq:pressure-localization-intervals}.
\end{enumerate}
The restriction of the velocity error $V_{\err}$ to the symmetry axis appears through the following three terms
in the evolution equations for $\mathsf V_t$, $\log q_t$, and $\log b_t$:
\[
J^{-2}(V_{\err})_z(0,J^2\zeta,t),\qquad (\p_rV_{\err})_r(0,J^2\zeta,t),\qquad (\p_zV_{\err})_z(0,J^2\zeta,t).
\]
Thus the relevant error norm on the symmetry axis is $\mathcal T_{\err}(I_\zeta,t)$, defined in
\eqref{eq:localized-cusp-error-axis-trace}; it contains the $L^\infty(I_\zeta)$ and
$C^{\alpha/2}(I_\zeta)$ norms of these three functions.  The estimate needed in
Lemma~\ref{lem:axis-profile-evolution} is
\begin{equation}
\mathcal T_{\err}(I_\zeta,t)\le C_{\err}\Gamma\bigl(J^{9\alpha-1}+1\bigr).
\label{eq:axis-profile-needed-error-bound}
\end{equation}
There are two hypotheses under which \eqref{eq:axis-profile-needed-error-bound} holds.  First, if
$\mathfrak E_{\err}(t)\le E_*$, then by \eqref{eq:localized-cusp-error-bootstrap},
\[
\mathcal T_{\err}(I_\zeta,t) \le \mathcal T_{\err}(I_{\err},t) \le E_*\Gamma(J^{9\alpha-1}+1).
\]
Second, if \eqref{eq:tail-axis-error-bound} holds with $I=I_{\err}$, then
\[
\mathcal T_{\err}(I_\zeta,t) \le \mathcal T_{\err}(I_{\err},t) \le C_{I_{\err}}\Gamma(J^{9\alpha-1}+1).
\]
Since both alternatives imply \eqref{eq:axis-profile-needed-error-bound}, the following lemma is stated under
either one.

\begin{lemma}[Transport equations in the clock-scaled axial variable]
\label{lem:axis-profile-evolution}
Let $J:=J_{\cusp}(t)$ and let $I_\zeta\subset I_{\err}$ be a compact $\zeta$-interval.  Assume
\eqref{eq:entry-axis-bounds-statement} on $I_\zeta\cap(0,\infty)$.  Assume also one of the two
error hypotheses
\begin{subequations} 
\begin{equation}
J\le \mathfrak J_{\mathrm{velocity}}, \qquad \mathfrak E_{\err}(t)\le E_*,
\label{eq:axis-profile-error-bootstrap}
\end{equation}
or
\begin{equation}
J\le \mathfrak J_{\mathrm{tail}}, \qquad \eqref{eq:tail-axis-error-bound}\ \text{holds with } I=I_{\err}.
\label{eq:axis-profile-error-tail}
\end{equation}
\end{subequations} 
With $q_t,b_t,\mathsf W_t,\mathsf U_t$, and $\mathsf V_t$ defined by
\eqref{eq:axis-normalized-derivatives-def}, \eqref{eq:axis-profile-def}, and
\eqref{eq:axis-current-zeta-velocity}, the following conclusions hold on $I_\zeta$.
First,
\begin{equation}
-C\le \mathsf W_t(0)\le -c<0, \qquad \|\mathsf W_t\|_{L^\infty(I_\zeta)} \le C ,
\label{eq:axis-profile-size}
\end{equation}
with constants depending only on the fixed parameters and, in the bootstrap-closure argument, on the fixed
barriers \eqref{eq:fixed-bootstrap-barriers}.  Second, for every fixed axial label $Z$ such that
$J^{-2}B_t(Z)\in I_\zeta$,
\begin{equation}
\p_t\bigl(J^{-2}B_t(Z)\bigr) = \mathsf V_t\bigl(J^{-2}B_t(Z)\bigr).
\label{eq:zeta-current-transport}
\end{equation}
Third, the normalized axis derivatives $q_t$ and $b_t$ obey
\begin{subequations}
\label{eq:axis-profile-current-system}
\begin{align}
\bigl(\p_t+\mathsf V_t\p_\zeta\bigr)\log q_t(\zeta)
&=
-\tfrac12 m(t)\Gamma J^{3\alpha-1}
\bigl(\mathsf W_t(\zeta)-\mathsf W_t(0)\bigr)
+(\p_rV_{\err})_r(0,J^2\zeta,t),
\label{eq:axis-profile-q-system}\\
\bigl(\p_t+\mathsf V_t\p_\zeta\bigr)\log b_t(\zeta)
&=
m(t)\Gamma J^{3\alpha-1}
\bigl(\mathsf W_t(\zeta)-\mathsf W_t(0)\bigr)
+(\p_zV_{\err})_z(0,J^2\zeta,t).
\label{eq:axis-profile-b-system}
\end{align}
\end{subequations}
Finally, the error terms in \eqref{eq:axis-current-zeta-velocity} and
\eqref{eq:axis-profile-current-system} satisfy
\begin{align} 
&\|J^{-2}(V_{\err})_z(0,J^2\cdot,t)\|_{C^{\alpha/2}(I_\zeta)}\notag \\
&\quad
+\|(\p_rV_{\err})_r(0,J^2\cdot,t)\|_{C^{\alpha/2}(I_\zeta)} + \|(\p_zV_{\err})_z(0,J^2\cdot,t)\|_{C^{\alpha/2}(I_\zeta)} \le C\Gamma\bigl(J^{9\alpha-1}+1\bigr).
\label{eq:axis-profile-error-holder}
\end{align} 
\end{lemma}

\begin{proof}[Proof of Lemma~\ref{lem:axis-profile-evolution}]
\runinhead{Step 1: The $\zeta$-coordinate and cusp-strain bounds.}
With  $ \zeta_t(Z):=J^{-2}B_t(Z)$, by \eqref{eq:entry-axis-bounds-statement}, we have that
$\p_Z\zeta_t(Z)=J^{-2}B_t'(Z)\ge c_{\rm ax}$.  Hence $Z\mapsto\zeta_t(Z)$ is strictly increasing, and the
inverse $Z_t$ is well defined.  By \eqref{eq:Wcusp-scaling},
\[
-C\le \Gamma^{-1}J^{1-3\alpha}\mathcal W_{\cusp}(t) \le -c<0.
\]
Since $\mathcal W_{\cusp}(t)=\p_z(U_{\cusp})_z(0,0,t)$, this proves the first bound in \eqref{eq:axis-profile-size}.
For $\zeta>0$, the axis point $(0,J^2\zeta)$ lies in the buffered cone
$\mathcal C_*$ so that by \eqref{eq:Ucusp-grad-Linf},  we find that
\[
|\mathsf W_t(\zeta)| \le \Gamma^{-1}J^{1-3\alpha} \|\nabla U_{\cusp}(\cdot,t)\|_{L^\infty(\mathcal C_*)} \le C \ \ \text{ for } \ \
\zeta\in I_\zeta,\ \zeta>0,
\]
while for $\zeta=0$ the preceding bound for $\mathsf W_t(0)$ applies. Hence, this
proves the second bound in \eqref{eq:axis-profile-size}.

\runinhead{Step 2: The transport equation for $\zeta_t(Z)$.}
The transported vorticity in \eqref{eq:Omega-cusp-def} is odd in the axial variable, as inherited from
the datum \eqref{eq:vort0} and the odd symmetry specified in Definition~\ref{def:init-data}.  Therefore
$U_{\cusp}(0,t)=0$.  Hence, on the symmetry axis,
\[
\begin{aligned}
(U_{\cusp})_z(0,J^2\zeta,t)
&=   \int_0^{J^2\zeta}\p_z(U_{\cusp})_z(0,s,t)\,ds    =   J^2\int_0^\zeta    \p_z(U_{\cusp})_z(0,J^2\eta,t)\,d\eta
= \Gamma J^{3\alpha+1}\mathsf U_t(\zeta).
\end{aligned}
\]
For a fixed axial label $Z$,
\[
\p_tB_t(Z)=(V_{\cusp})_z(0,B_t(Z),t), \qquad V_{\cusp}=mU_{\cusp}+V_{\err}.
\]
The cusp-clock identity \eqref{eq:cusp-clock-identity-in-Jdot-proof} is then written as
\[
\tfrac{\dot J}{J}=\tfrac12m(t)\mathcal W_{\cusp}(t) = \tfrac12m(t)\Gamma J^{3\alpha-1}\mathsf W_t(0),
\]
and so, with $\zeta=\zeta_t(Z)$,
\[
\p_t\zeta_t(Z) = -2\tfrac{\dot J}{J}\zeta +J^{-2}(V_{\cusp})_z(0,J^2\zeta,t) =
m\Gamma J^{3\alpha-1} \bigl(\mathsf U_t(\zeta)-\mathsf W_t(0)\zeta\bigr) +J^{-2}(V_{\err})_z(0,J^2\zeta,t),
\]
which establishes \eqref{eq:zeta-current-transport}.

\runinhead{Step 3: Fixed-label evolution of $q_t$ and $b_t$.}
Differentiating the cusp-flow equation
$\p_t(r_t,z_t)=V_{\cusp}(r_t,z_t,t)$ in the radial label at $R=0$, we find that
\[
\p_tA_t(Z) = \p_r(V_{\cusp})_r(0,B_t(Z),t)\,A_t(Z),
\]
where the possible term
$\p_z(V_{\cusp})_r(0,B_t(Z),t)\,\p_Rz_t(0,Z)$ vanishes because
$\p_Rz_t(0,Z)=0$.  This identity follows from the evenness of $z_t(R,Z)$ in the cylindrical radius $R$.
Similarly, differentiating $B_t(Z)=z_t(0,Z)$ in $Z$, we obtain
\[
\p_tB_t'(Z) = \p_z(V_{\cusp})_z(0,B_t(Z),t)\,B_t'(Z),
\]
because $r_t(0,Z)=0$ and hence $\p_Zr_t(0,Z)=0$.
Axisymmetric incompressibility yields
\[
\p_r(U_{\cusp})_r(0,z,t) = -\tfrac12\,\p_z(U_{\cusp})_z(0,z,t).
\]
Using these identities and $V_{\cusp}=mU_{\cusp}+V_{\err}$, we obtain at fixed $Z$ that
\[
\begin{aligned}
\p_t\log(JA_t(Z)) &= \tfrac{\dot J}{J} +m\,\p_r(U_{\cusp})_r(0,B_t(Z),t) +(\p_rV_{\err})_r(0,B_t(Z),t) \\
&= -\tfrac12m\Gamma J^{3\alpha-1} \bigl(\mathsf W_t(\zeta)-\mathsf W_t(0)\bigr) +(\p_rV_{\err})_r(0,J^2\zeta,t),
\end{aligned}
\]
and
\[
\begin{aligned}
\p_t\log(J^{-2}B_t'(Z)) &= -2\tfrac{\dot J}{J} +m\,\p_z(U_{\cusp})_z(0,B_t(Z),t) +(\p_zV_{\err})_z(0,B_t(Z),t) \\
&= m\Gamma J^{3\alpha-1} \bigl(\mathsf W_t(\zeta)-\mathsf W_t(0)\bigr) +(\p_zV_{\err})_z(0,J^2\zeta,t).
\end{aligned}
\]

\runinhead{Step 4: Passage to the $\zeta$-coordinate.}
We let $F_t(\zeta):=\log q_t(\zeta)$.  Since
$F_t(\zeta_t(Z))=\log(JA_t(Z))$, the chain rule and
\eqref{eq:zeta-current-transport} imply that
\[
\p_t\log(JA_t(Z)) = \bigl(\p_t+\mathsf V_t\p_\zeta\bigr)\log q_t \big|_{\zeta=\zeta_t(Z)} .
\]
The same argument with $F_t(\zeta)=\log b_t(\zeta)$ yields the equation for $b_t$.  Since
$\zeta_t$ maps the corresponding axial label interval onto $I_\zeta$, this proves
\eqref{eq:axis-profile-current-system}.

\runinhead{Step 5: The error bounds.}
The error estimates in \eqref{eq:zeta-current-transport} and
\eqref{eq:axis-profile-current-system} use exactly the three axis traces
\[
J^{-2}(V_{\err})_z(0,J^2\zeta,t),\qquad (\p_rV_{\err})_r(0,J^2\zeta,t),\qquad (\p_zV_{\err})_z(0,J^2\zeta,t).
\]
If \eqref{eq:axis-profile-error-bootstrap} holds, then by $I_\zeta\subset I_{\err}$,
\eqref{eq:localized-cusp-error-axis-trace}, \eqref{eq:localized-cusp-error-bootstrap}, and
\eqref{eq:localized-cusp-error-large-bootstrap}, we have that
\[
\mathcal T_{\err}(I_\zeta,t) \le \mathcal T_{\err}(I_{\err},t) \le E_*\Gamma\bigl(J^{9\alpha-1}+1\bigr).
\]
If instead \eqref{eq:axis-profile-error-tail} holds, then
\eqref{eq:tail-axis-error-bound} on $I_{\err}$ implies the same bound on $I_\zeta$ by restriction.
Thus the three traces satisfy \eqref{eq:axis-profile-error-holder}.
\end{proof}

The estimate \eqref{eq:axis-profile-error-holder} is obtained on the $\zeta$-interval $I_\zeta$.
We next use it on a fixed reference axial label interval.  Let $I_\eta$ be a compact interval for $\eta$, let
$Z_0:I_\eta\to\mathbb R$ be the axial label parametrization, and set
\[
\zeta(\eta,t):=J^{-2}B_t(Z_0(\eta)).
\]
We assume
\begin{equation}
\zeta(\cdot,t):I_\eta\longrightarrow \zeta(I_\eta,t)\subset I_\zeta,
\label{eq:axis-profile-pullback-image}
\end{equation}
and
\begin{equation}
\|\zeta(\cdot,t)\|_{\operatorname{Lip}(I_\eta)}\le L_\eta .
\label{eq:axis-profile-pullback-lip}
\end{equation}

\begin{corollary}[Same-label pullback and axis quotient estimates]
\label{cor:axis-profile-samelabel-errors}
Assume \eqref{eq:axis-profile-pullback-image} and
\eqref{eq:axis-profile-pullback-lip}.  Then
\begin{equation}
\begin{aligned}
&\|J^{-2}(V_{\err})_z(0,J^2\zeta(\cdot,t),t)\|_{C^{\alpha/2}(I_\eta)}\\
&
+\|(\p_rV_{\err})_r(0,J^2\zeta(\cdot,t),t)\|_{C^{\alpha/2}(I_\eta)}
+
\|(\p_zV_{\err})_z(0,J^2\zeta(\cdot,t),t)\|_{C^{\alpha/2}(I_\eta)}
\le C(1+L_\eta^{\alpha/2})\Gamma\bigl(J^{9\alpha-1}+1\bigr).
\end{aligned}
\label{eq:axis-profile-error-holder-samelabel}
\end{equation}
If, in addition, $I_\zeta=[0,\zeta_*]$ is origin-attached,
$0\in I_\eta$, $\zeta(0,t)=0$, and
$0\le \zeta(\eta,t)\le\zeta_*$ for $\eta\in I_\eta$, then
\begin{equation}
\left\| \tfrac{J^{-2}(V_{\err})_z(0,J^2\zeta(\cdot,t),t)} {\zeta(\cdot,t)} \right\|_{C^{\alpha/2}(I_\eta)}
\le C(1+L_\eta^{\alpha/2})\Gamma\bigl(J^{9\alpha-1}+1\bigr),
\label{eq:axis-profile-error-quotient-samelabel}
\end{equation}
with the quotient understood by its continuous value at $\eta=0$.
\end{corollary}

\begin{proof}[Proof of Corollary~\ref{cor:axis-profile-samelabel-errors}]
If
$F\in C^{\alpha/2}(I_\zeta)$ and $\|\zeta(\cdot,t)\|_{\operatorname{Lip}(I_\eta)}\le L_\eta$, then
\[
\|F\circ\zeta(\cdot,t)\|_{C^{\alpha/2}(I_\eta)} \le (1+L_\eta^{\alpha/2})\|F\|_{C^{\alpha/2}(I_\zeta)}.
\]
Applying this to the three functions in \eqref{eq:axis-profile-error-holder} proves
\eqref{eq:axis-profile-error-holder-samelabel}.

For the quotient estimate, $V_{\err}(0,t)=0$, and the modulation identity yields
\[
(\p_zV_{\err})_z(0,0,t) = \p_z(V_{\cusp})_z(0,0,t) -m(t)\p_z(U_{\cusp})_z(0,0,t) =0.
\]
Thus, for $z\ge0$,
\[
(V_{\err})_z(0,z,t) = z\int_0^1 \Bigl[ (\p_zV_{\err})_z(0,\theta z,t) -(\p_zV_{\err})_z(0,0,t) \Bigr]\,d\theta.
\]
For $z=J^2\zeta(\eta,t)$, this becomes
\[
\tfrac{J^{-2}(V_{\err})_z(0,J^2\zeta(\eta,t),t)}{\zeta(\eta,t)} =
\int_0^1 \Bigl[ (\p_zV_{\err})_z(0,\theta J^2\zeta(\eta,t),t) -(\p_zV_{\err})_z(0,0,t) \Bigr]\,d\theta,
\]
with the right-hand side defining the continuous value at $\eta=0$.  The additional hypotheses for
\eqref{eq:axis-profile-error-quotient-samelabel} imply that, for every $\eta\in I_\eta$ and
$0\le\theta\le1$,
$0\le\theta\zeta(\eta,t)\le\zeta_*$, so $\theta\zeta(\eta,t)\in I_\zeta$.  We set
\[
G(\zeta):=(\p_zV_{\err})_z(0,J^2\zeta,t)-(\p_zV_{\err})_z(0,0,t).
\]
By \eqref{eq:axis-profile-error-holder},
\[
\|G\|_{L^\infty(I_\zeta)} + [G]_{C^{\alpha/2}(I_\zeta)} \le C\Gamma\bigl(J^{9\alpha-1}+1\bigr).
\]
The quotient in \eqref{eq:axis-profile-error-quotient-samelabel} is
\[
Q(\eta):=\int_0^1G(\theta\zeta(\eta,t))\,d\theta .
\]
Thus
\[
\|Q\|_{L^\infty(I_\eta)} \le C\Gamma\bigl(J^{9\alpha-1}+1\bigr),
\]
and, for $\eta_1,\eta_2\in I_\eta$,
\begin{align*}
|Q(\eta_1)-Q(\eta_2)|
\le
[G]_{C^{\alpha/2}(I_\zeta)}
\int_0^1
\theta^{\alpha/2}|\zeta(\eta_1,t)-\zeta(\eta_2,t)|^{\alpha/2}\,d\theta
\le
C L_\eta^{\alpha/2}\Gamma\bigl(J^{9\alpha-1}+1\bigr)
|\eta_1-\eta_2|^{\alpha/2}.
\end{align*}
Together with \eqref{eq:axis-profile-pullback-lip}, we obtain
\eqref{eq:axis-profile-error-quotient-samelabel}.
\end{proof}

The transport equation \eqref{eq:zeta-current-transport} for the $\zeta$-coordinate of a fixed axial
label is driven by the velocity $\mathsf V_t$ from \eqref{eq:axis-current-zeta-velocity}, which
decomposes into a principal part generated by the flat cusp velocity $U_{\cusp}$ and an error part
generated by $V_{\err}$, $\mathsf V_t(\zeta)=\mathsf V_t^{\rm ax}(\zeta)+\mathsf R_t^\zeta(\zeta)$, with
\[
\mathsf V_t^{\rm ax}(\zeta) := m(t)\Gamma J^{3\alpha-1}\bigl(\mathsf U_t(\zeta)-\mathsf W_t(0)\zeta\bigr),
\qquad \mathsf R_t^\zeta(\zeta):=J^{-2}(V_{\err})_z(0,J^2\zeta,t) .
\]
We further introduce the two $V_{\err}$ axis traces
\[
\mathsf R_t^q(\zeta):=(\p_rV_{\err})_r(0,J^2\zeta,t), \qquad \mathsf R_t^b(\zeta):=(\p_zV_{\err})_z(0,J^2\zeta,t),
\]
which control the $V_{\err}$ contributions to the evolution equations
\eqref{eq:axis-profile-current-system} for $\log q_t$ and $\log b_t$.  The $\zeta$-derivative of the
principal part,
\[
\mathsf S_t(\zeta):=\p_\zeta\mathsf V_t^{\rm ax}(\zeta) = m(t)\Gamma J^{3\alpha-1}\bigl(\mathsf W_t(\zeta)-\mathsf W_t(0)\bigr),
\]
is the singular scalar driving the right-hand sides of \eqref{eq:axis-profile-current-system}.
The next lemma absorbs $\mathsf S_t$ into the Jacobian of the one-dimensional flow $\mathscr X_t$ of
$\mathsf V_t^{\rm ax}$: in the new axial coordinate $\eta$ defined by this flow, the renormalized functions
$\widetilde q_t,\widetilde b_t$ obtained from $q_t,b_t$ obey evolution equations at fixed $\eta$ containing
only the error traces $\mathsf R_t^\zeta,\mathsf R_t^q,\mathsf R_t^b$.

The one-dimensional flow $\mathscr X_t$ in Lemma~\ref{lem:axis-profile-composed-coordinate}, and
the one-dimensional flow $\mathscr Z_t$ in Lemma~\ref{lem:exact-axis-composition}, are each defined
as follows.  We fix an entry time $t_0$ and a compact reference interval $I_\eta\subset I_\zeta$.  For the
$\mathscr X_t$ flow, we use the time interval
\[
\mathcal I_X=[t_0,T_X), \qquad T_X:=\sup\{T>t_0:\ \mathscr X_t(I_\eta)\subset I_\zeta \ \text{for every } t\in[t_0,T]\}.
\]
For the $\mathscr Z_t$ flow, we use
\[
\mathcal I_Z=[t_0,T_Z), \qquad T_Z:=\sup\{T>t_0:\ \mathscr Z_t(I_\eta)\subset I_\zeta \ \text{for every } t\in[t_0,T]\}.
\]
The initial conditions are $\mathscr X_{t_0}(\eta)=\eta$ and
$\mathscr Z_{t_0}(\eta)=\eta$ on $I_\eta$.  The axial derivatives are taken first on the interior of
$I_\eta$, where the Euler solution is regular and the axis traces are classical; when $0\in I_\eta$, the
identities extend to $\eta=0$ by the continuous axis values from Lemma~\ref{lem:axis-profile-evolution},
and the Jacobian identities for the one-dimensional flows are read in their variational integral form.

Specifically, $\mathscr X_t$ is the solution of the ODE
\begin{equation}
\p_t\mathscr X_t(\eta) = \mathsf V_t^{\rm ax}(\mathscr X_t(\eta)), \qquad \mathscr X_{t_0}(\eta)=\eta,\quad \eta\in I_\eta .
\label{eq:axis-profile-X-flow}
\end{equation}
For each axial label $Z$ with $J^{-2}B_t(Z)\in\mathscr X_t(I_\eta)$, we set
\begin{equation}
\eta_t(Z):=\mathscr X_t^{-1}\bigl(J^{-2}B_t(Z)\bigr),
\qquad \widetilde{\mathsf R}_t^\eta(\eta):= \tfrac{\mathsf R_t^\zeta(\mathscr X_t(\eta))}{\p_\eta\mathscr X_t(\eta)} ,
\label{eq:axis-profile-X-eta-def}
\end{equation}
and we define the renormalized axis derivatives on $I_\eta$ by
\begin{equation}
\widetilde q_t(\eta) := q_t(\mathscr X_t(\eta))\bigl(\p_\eta\mathscr X_t(\eta)\bigr)^{{\frac{1}{2}}},
\qquad \widetilde b_t(\eta) := \tfrac{b_t(\mathscr X_t(\eta))}{\p_\eta\mathscr X_t(\eta)} .
\label{eq:axis-renormalized-functions}
\end{equation}

\begin{lemma}[Evolution of renormalized flow]
\label{lem:axis-profile-composed-coordinate}
Assuming the hypotheses of Lemma~\ref{lem:axis-profile-evolution}, the flow $\mathscr X_t$ defined by
\eqref{eq:axis-profile-X-flow} and the renormalized functions $\widetilde q_t,\widetilde b_t$ defined by
\eqref{eq:axis-renormalized-functions} satisfy, for every $t\in\mathcal I_X$, the following identities.  The
Jacobian of $\mathscr X_t$ satisfies
\begin{equation}
\p_t\log\p_\eta\mathscr X_t(\eta)=\mathsf S_t(\mathscr X_t(\eta)).
\label{eq:axis-profile-X-jac}
\end{equation}
For any axial label $Z$ with $J^{-2}B_t(Z)\in\mathscr X_t(I_\eta)$, the coordinate $\eta_t(Z)$ defined in
\eqref{eq:axis-profile-X-eta-def} obeys
\begin{equation}
\p_t\eta_t(Z)=\widetilde{\mathsf R}_t^\eta(\eta_t(Z)).
\label{eq:eta-current-transport}
\end{equation}
At fixed $\eta\in I_\eta$, the renormalized functions satisfy equations in which the singular scalar
$\mathsf S_t$ has already been removed by the change of variables \eqref{eq:axis-profile-X-flow} and the
Jacobian identity \eqref{eq:axis-profile-X-jac}:
\begin{subequations}
\label{eq:axis-renormalized-fixed-eta}
\begin{align}
\p_t\log\widetilde q_t(\eta)
&=\mathsf R_t^q(\mathscr X_t(\eta))
-\mathsf R_t^\zeta(\mathscr X_t(\eta))\,\p_\zeta\log q_t(\mathscr X_t(\eta)),
\label{eq:axis-renormalized-q-fixed-eta}\\
\p_t\log\widetilde b_t(\eta)
&=\mathsf R_t^b(\mathscr X_t(\eta))
-\mathsf R_t^\zeta(\mathscr X_t(\eta))\,\p_\zeta\log b_t(\mathscr X_t(\eta)).
\label{eq:axis-renormalized-b-fixed-eta}
\end{align}
\end{subequations}
The renormalized functions also satisfy the axis volume identity
\begin{equation}
\widetilde q_t(\eta)^2\,\widetilde b_t(\eta)=1 .
\label{eq:axis-renormalized-volume}
\end{equation}
\end{lemma}

\begin{remark}[Transport form of the renormalized equations]
The fixed-$\eta$ identities \eqref{eq:axis-renormalized-fixed-eta} can also be written along the residual
$\eta$-transport generated by $\widetilde{\mathsf R}_t^\eta$.  Using
\eqref{eq:axis-profile-X-jac}, we obtain that
\begin{subequations}
\label{eq:axis-renormalized-transport}
\begin{align}
\bigl(\p_t+\widetilde{\mathsf R}_t^\eta\p_\eta\bigr)\log\widetilde q_t
&=\mathsf R_t^q\circ\mathscr X_t
+\tfrac12\widetilde{\mathsf R}_t^\eta\,\p_\eta\log\p_\eta\mathscr X_t, \\
\bigl(\p_t+\widetilde{\mathsf R}_t^\eta\p_\eta\bigr)\log\widetilde b_t
&=\mathsf R_t^b\circ\mathscr X_t
-\widetilde{\mathsf R}_t^\eta\,\p_\eta\log\p_\eta\mathscr X_t .
\end{align}
\end{subequations}
Thus the singular scalar $\mathsf S_t$ does not appear on the right-hand sides of
\eqref{eq:axis-renormalized-transport}; it has been absorbed into the Jacobian
$\p_\eta\mathscr X_t$ through \eqref{eq:axis-profile-X-jac}.
\end{remark}

\begin{proof}[Proof of Lemma~\ref{lem:axis-profile-composed-coordinate}]
By \eqref{eq:axis-current-zeta-velocity},
$\mathsf V_t=\mathsf V_t^{\rm ax}+\mathsf R_t^\zeta$.  With the definitions of
$\mathsf S_t,\mathsf R_t^q,\mathsf R_t^b$, the system
\eqref{eq:axis-profile-current-system} becomes
\begin{align*}
\bigl(\p_t+(\mathsf V_t^{\rm ax}+\mathsf R_t^\zeta)\p_\zeta\bigr)
\log q_t
&=
-\tfrac12\mathsf S_t+\mathsf R_t^q,
\\
\bigl(\p_t+(\mathsf V_t^{\rm ax}+\mathsf R_t^\zeta)\p_\zeta\bigr)
\log b_t
&=
\mathsf S_t+\mathsf R_t^b .
\end{align*}

Since $\mathsf U_t$ is absolutely continuous and $\p_\zeta\mathsf U_t=\mathsf W_t$ a.e., we have that
\[
\p_\zeta\mathsf V_t^{\rm ax}(\zeta) = m(t)\Gamma J^{3\alpha-1}\bigl(\mathsf W_t(\zeta)-\mathsf W_t(0)\bigr) = \mathsf S_t(\zeta).
\]
The bound \eqref{eq:axis-profile-size} implies
$\mathsf S_t\in L^\infty(I_\zeta)$ at each fixed time, so
$\mathsf V_t^{\rm ax}$ is Lipschitz in $\zeta$.  Hence the flow
\eqref{eq:axis-profile-X-flow} is unique and bi-Lipschitz in $\eta$ on every interval whose image
remains in $I_\zeta$, and its a.e. derivative satisfies
\[
\p_t\p_\eta\mathscr X_t = (\p_\zeta\mathsf V_t^{\rm ax})(\mathscr X_t) \p_\eta\mathscr X_t.
\]
This proves \eqref{eq:axis-profile-X-jac}.  If $\zeta_t(Z)=J^{-2}B_t(Z)$ and
$\zeta_t(Z)=\mathscr X_t(\eta_t(Z))$, then
\[
\mathsf V_t^{\rm ax}(\mathscr X_t(\eta_t)) +\p_\eta\mathscr X_t(\eta_t)\,\p_t\eta_t =
\mathsf V_t^{\rm ax}(\mathscr X_t(\eta_t)) +\mathsf R_t^\zeta(\mathscr X_t(\eta_t)),
\]
which is \eqref{eq:eta-current-transport}.

For $\widetilde q_t$, the chain rule at fixed $\eta$ yields
\[
\p_t\log\widetilde q_t = \bigl(\p_t+\mathsf V_t^{\rm ax}\p_\zeta\bigr) \log q_t(\mathscr X_t) +\tfrac12\p_t\log\p_\eta\mathscr X_t
= \bigl(-\tfrac12\mathsf S_t+\mathsf R_t^q -\mathsf R_t^\zeta\p_\zeta\log q_t\bigr)(\mathscr X_t) +\tfrac12\mathsf S_t(\mathscr X_t),
\]
and the two $\mathsf S_t$ terms cancel.  This proves
\eqref{eq:axis-renormalized-q-fixed-eta}.  For $\widetilde b_t$, the same chain rule and
\eqref{eq:axis-profile-X-jac} yield
\[
\p_t\log\widetilde b_t = \bigl(\p_t+\mathsf V_t^{\rm ax}\p_\zeta\bigr) \log b_t(\mathscr X_t) - \p_t\log\p_\eta\mathscr X_t =
\bigl(\mathsf S_t+\mathsf R_t^b -\mathsf R_t^\zeta\p_\zeta\log b_t\bigr)(\mathscr X_t) -\mathsf S_t(\mathscr X_t).
\]
The two $\mathsf S_t$ terms cancel, and this proves \eqref{eq:axis-renormalized-b-fixed-eta}.

To obtain the transport form, we use that $\mathsf R_t^\zeta(\mathscr X_t(\eta)) =   \p_\eta\mathscr X_t(\eta)\,\widetilde{\mathsf R}_t^\eta(\eta)$
and
\[
\p_\zeta\log q_t(\mathscr X_t) = \tfrac{1}{\p_\eta\mathscr X_t} \left(\p_\eta\log\widetilde q_t -\tfrac12\p_\eta\log\p_\eta\mathscr X_t
\right), \ \ \ \p_\zeta\log b_t(\mathscr X_t)= \tfrac{1}{\p_\eta\mathscr X_t} \left( \p_\eta\log\widetilde b_t
+\p_\eta\log\p_\eta\mathscr X_t\right).
\]
Substituting these two identities into \eqref{eq:axis-renormalized-fixed-eta}, we obtain
\eqref{eq:axis-renormalized-transport}.
Finally, by \eqref{eq:axis-qb-volume-mon}, $q_t(\zeta)^2b_t(\zeta)=1$.  Pulling-back this identity by
$\mathscr X_t$ and using \eqref{eq:axis-renormalized-functions},
\[
\widetilde q_t(\eta)^2\widetilde b_t(\eta)=q_t(\mathscr X_t(\eta))^2b_t(\mathscr X_t(\eta)) \p_\eta\mathscr X_t(\eta)
\bigl(\p_\eta\mathscr X_t(\eta)\bigr)^{-1} = 1.
\]
This is \eqref{eq:axis-renormalized-volume}.
\end{proof}

We now turn to the flow generated by the \emph{full} axial velocity $\mathsf V_t$ rather than by its
principal part $\mathsf V_t^{\rm ax}$.  Assuming the same time-interval convention as before, we define
$\mathscr Z_t$ as the solution of
\begin{equation}
\p_t\mathscr Z_t(\eta) = \mathsf V_t(\mathscr Z_t(\eta)), \qquad \mathscr Z_{t_0}(\eta)=\eta,\quad \eta\in I_\eta ,
\label{eq:exact-axis-Zflow}
\end{equation}
and, for each axial label $Z$ with $J^{-2}B_t(Z)\in\mathscr Z_t(I_\eta)$, we set
\begin{equation}
\eta_t(Z) := \mathscr Z_t^{-1}\bigl(J^{-2}B_t(Z)\bigr) .
\label{eq:exact-axis-eta-def}
\end{equation}
The normalized axis derivatives, renormalized by the flow $\mathscr Z_t$, are
\begin{equation}
\widehat q_t(\eta) := q_t(\mathscr Z_t(\eta))\bigl(\p_\eta\mathscr Z_t(\eta)\bigr)^{{\frac{1}{2}}},
\qquad \widehat b_t(\eta) := \tfrac{b_t(\mathscr Z_t(\eta))}{\p_\eta\mathscr Z_t(\eta)} .
\label{eq:exact-axis-renormalized-functions}
\end{equation}
The next lemma shows that because $\mathscr Z_t$ tracks the actual axial trajectories of the cusp flow,
$\eta_t(Z)$ is constant in $t$ and the renormalized functions $\widehat q_t,\widehat b_t$ are exact
conservation laws.

\begin{lemma}[Exact conservation along the full axial flow]
\label{lem:exact-axis-composition}
With the hypotheses of Lemma~\ref{lem:axis-profile-evolution}, the flow $\mathscr Z_t$ defined by
\eqref{eq:exact-axis-Zflow} and the functions $\widehat q_t,\widehat b_t$ defined by
\eqref{eq:exact-axis-renormalized-functions} satisfy, for every $t\in\mathcal I_Z$, the following
identities.  The Jacobian of $\mathscr Z_t$ satisfies
\begin{equation}
\p_t\log\p_\eta\mathscr Z_t(\eta) = m(t)\Gamma J^{3\alpha-1}\bigl(\mathsf W_t(\mathscr Z_t(\eta))-\mathsf W_t(0)\bigr)
+ (\p_zV_{\err})_z(0,J^2\mathscr Z_t(\eta),t) .
\label{eq:exact-axis-Zflow-jac}
\end{equation}
For each axial label $Z$ with $J^{-2}B_t(Z)\in\mathscr Z_t(I_\eta)$, the coordinate
$\eta_t(Z)$ defined in \eqref{eq:exact-axis-eta-def} is conserved:
\begin{equation}
\p_t\eta_t(Z) = 0 .
\label{eq:exact-axis-eta-constant}
\end{equation}
At fixed $\eta\in I_\eta$, the renormalized functions satisfy the exact conservation laws
\begin{equation}
\p_t\widehat q_t(\eta)=0, \qquad \p_t\widehat b_t(\eta)=0, \qquad \widehat q_t(\eta)^2\,\widehat b_t(\eta) = 1 .
\label{eq:exact-axis-renormalized-conserved}
\end{equation}
\end{lemma}

\begin{proof}[Proof of Lemma~\ref{lem:exact-axis-composition}]
We use the same endpoint interpretation as in Lemma~\ref{lem:axis-profile-composed-coordinate}.
Differentiating the axial velocity $\mathsf V_t$ from \eqref{eq:axis-current-zeta-velocity}, we obtain
\begin{equation}
\p_\zeta\mathsf V_t(\zeta) = m(t)\Gamma J^{3\alpha-1} \bigl(\mathsf W_t(\zeta)-\mathsf W_t(0)\bigr) +(\p_zV_{\err})_z(0,J^2\zeta,t).
\label{eq:exact-axis-V-zeta-derivative}
\end{equation}
The right-hand side is bounded on every stopped interval by
\eqref{eq:axis-profile-size} and \eqref{eq:axis-profile-error-holder}.  Hence the full axial velocity
in the coordinate $\zeta$ is Lipschitz, $\mathscr Z_t$ is an increasing bi-Lipschitz flow, and its
Jacobian equation is justified by the one-dimensional variational equation.  The Jacobian identity
\eqref{eq:exact-axis-Zflow-jac} follows by differentiating
\eqref{eq:exact-axis-Zflow} in $\eta$ and using
\eqref{eq:exact-axis-V-zeta-derivative}.
Since $\zeta_t(Z):=J^{-2}B_t(Z)$ solves \eqref{eq:zeta-current-transport}, and
\[
\zeta_t(Z)=\mathscr Z_t(\eta_t(Z)),
\]
differentiating in $t$ yields
\[
\mathsf V_t(\mathscr Z_t(\eta_t)) = \mathsf V_t(\mathscr Z_t(\eta_t)) +\p_\eta\mathscr Z_t(\eta_t)\,\p_t\eta_t.
\]
Since $\p_\eta\mathscr Z_t>0$, this proves \eqref{eq:exact-axis-eta-constant}.

For the normalized radial derivative, the chain rule, \eqref{eq:axis-profile-q-system}, and
\eqref{eq:exact-axis-Zflow-jac} yield, at fixed $\eta$,
\[
\begin{aligned}
\p_t\log\widehat q_t
&=
\bigl(\p_t+\mathsf V_t\p_\zeta\bigr)\log q_t(\mathscr Z_t)
+\tfrac12\p_t\log\p_\eta\mathscr Z_t  \\
&= -\tfrac12 m(t)\Gamma J^{3\alpha-1} \bigl(\mathsf W_t(\mathscr Z_t)-\mathsf W_t(0)\bigr) +(\p_rV_{\err})_r(0,J^2\mathscr Z_t,t)  \\
&\qquad\qquad
+\tfrac12 m(t)\Gamma J^{3\alpha-1} \bigl(\mathsf W_t(\mathscr Z_t)-\mathsf W_t(0)\bigr) +\tfrac12(\p_zV_{\err})_z(0,J^2\mathscr Z_t,t).
\end{aligned}
\]
Both $V_{\cusp}$ and $U_{\cusp}$ are axisymmetric divergence-free, hence so is
$V_{\err}=V_{\cusp}-mU_{\cusp}$.  The regular axis limit of the divergence identity implies
\[
2(\p_rV_{\err})_r(0,z,t)+(\p_zV_{\err})_z(0,z,t)=0,
\]
and hence $\p_t\log\widehat q_t=0$.  For the normalized axial derivative,
\eqref{eq:axis-profile-b-system} and \eqref{eq:exact-axis-Zflow-jac} yield, at fixed $\eta$,
\[
\p_t\log\widehat b_t =  \bigl(\p_t+\mathsf V_t\p_\zeta\bigr)\log b_t(\mathscr Z_t) - \p_t\log\p_\eta\mathscr Z_t=0 .
\]
Finally, by \eqref{eq:axis-qb-volume-mon}, $q_t(\zeta)^2b_t(\zeta)=1$.  Pulling this identity back by
$\mathscr Z_t$ and using \eqref{eq:exact-axis-renormalized-functions},
\[
\widehat q_t(\eta)^2\widehat b_t(\eta) = q_t(\mathscr Z_t(\eta))^2b_t(\mathscr Z_t(\eta)) \p_\eta\mathscr Z_t(\eta)
\bigl(\p_\eta\mathscr Z_t(\eta)\bigr)^{-1} = 1.
\]
\end{proof}

We now close the monotone axial-stretching bootstrap \textup{(BA4)}: the two-sided bound
\eqref{eq:monotone-axial-two-sided} and the monotone fractional-increment bound
\eqref{eq:monotone-axial-fractional-bootstrap} for $b_t$ on $I_{\rm mon}$.
This is the monotonicity of the normalized axial derivative $b_t$ on the origin-attached interval
$I_{\rm mon}$.  The exact axial flow
$\mathscr Z_t$ is the correct coordinate for this question, because
\eqref{eq:exact-axis-Zflow-jac} shows that the logarithmic derivative of $\mathscr Z_t$ is driven by
\[
m(t)\Gamma J_{\cusp}(t)^{3\alpha-1} \bigl(\mathsf W_t(\zeta)-\mathsf W_t(0)\bigr) +(\p_zV_{\err})_z(0,J_{\cusp}(t)^2\zeta, t).
\]
(Recall that  $\mathsf W_t(\zeta):= \Gamma^{-1}J^{1-3\alpha}\,\p_z(U_{\cusp})_z(0,J^2\zeta,t)$ was defined in
\eqref{eq:axis-profile-def}.)

Thus the sign of the leading term is determined by the axial strain defect
$\mathsf W_t(\zeta)-\mathsf W_t(0)$.  We first isolate the defect estimate on a sufficiently small
origin-attached interval.  The interval $I_{\rm mon}=[0,\zeta_{\rm mon}]$ is then fixed with
\[
0<\zeta_{\rm mon}\le \zeta_{\rm def},
\]
where $\zeta_{\rm def}$ is supplied by the next lemma.

\begin{lemma}[Positive axial strain defect]
\label{lem:positive-axial-strain-defect}
There exist constants $\zeta_{\rm def}>0$, $c_{\rm def}>0$, $C_{\rm def}<\infty$, depending only on the fixed parameters, 
and a threshold $\mathfrak J_{\rm def}\le  \min\{\mathfrak J_{\mathrm{velocity}},\mathfrak J_{\mathrm{tail}},\mathfrak
J_{\mathrm{mod}}\}$,
such that the following holds.  If $J:=J_{\cusp}(t)\le\mathfrak J_{\rm def}$, then for
\[
0\le \zeta_1<\zeta_2\le \zeta_{\rm def}
\]
we have
\begin{equation}
c_{\rm def}\bigl(\zeta_2^\alpha-\zeta_1^\alpha\bigr) \le \mathsf W_t(\zeta_2)-\mathsf W_t(\zeta_1) \le C_{\rm def}
\bigl(\zeta_2^\alpha-\zeta_1^\alpha+\zeta_2^2-\zeta_1^2\bigr).
\label{eq:positive-axial-strain-defect}
\end{equation}
In particular,
\begin{equation}
0\le \mathsf W_t(\zeta)-\mathsf W_t(0) \le C_{\rm def}(\zeta^\alpha+\zeta^2), \qquad 0\le\zeta\le\zeta_{\rm def}.
\label{eq:positive-axial-strain-defect-origin}
\end{equation}
\end{lemma}

\begin{proof}[Proof of Lemma~\ref{lem:positive-axial-strain-defect}]
\runinhead{Step 1: The homogeneous axial defect.}
We first identify the fixed homogeneous field which appears in the normalized strain $\mathsf W_t$.  In the variables
$y=J^2X$, $X=(R,Z)$, the leading part of \eqref{eq:Omega-cusp-local-representation} on the small
origin-attached tube has the form
\[
\Omega_{\cusp,\theta}(J^2R,J^2Z,t) = -\Gamma J^{3\alpha-1}\operatorname{sgn}(Z)R^\alpha +\hbox{lower-order terms}.
\]
Let $U_{\rm hom}$ be the Biot--Savart velocity generated in the fixed variables $X=(R,Z)$ by
$\omega_{\rm hom}(X):=-\operatorname{sgn}(Z)R^\alpha e_\theta$.
We write
\[
\mathsf W_{\rm hom}(\zeta):=\p_Z(U_{\rm hom})_Z(0,\zeta).
\]
This velocity is introduced because it is the leading fixed-variable contribution to the normalized cusp
strain.  Indeed, the velocity-gradient kernel satisfies
\[
\mathcal K_W\bigl((0,0,J^2\zeta),J^2Y\bigr) = J^{-6}\mathcal K_W\bigl((0,0,\zeta),Y\bigr), \qquad \ud y=J^6\,\ud Y ,
\]
where $\mathcal K_W$ is the scalar axial strain kernel from \eqref{eq:rW2}.
Therefore the two powers of $J$ from the kernel and the volume form cancel.  Since
$\Omega_{\cusp,\theta}(J^2Y,t)$ carries $\Gamma J^{3\alpha-1}$, the leading-order term in
$\mathsf W_t(\zeta) = \Gamma^{-1}J^{1-3\alpha}\,  \p_z(U_{\cusp})_z(0,J^2\zeta,t)$
is $\mathsf W_{\rm hom}(\zeta)$ and we have that
\begin{equation}
\mathsf W_t(\zeta) = \mathsf W_{\rm hom}(\zeta)+\mathcal R_t^{\rm hom}(\zeta),
\label{eq:Wt-homogeneous-leading-part}
\end{equation}
where the remainder $\mathcal R_t^{\rm hom}$ is estimated below in Step~2.

We now compute the axial defect of $U_{\rm hom}$ from $\mathcal K_W$.  For an evaluation point on the symmetry axis
$x=(0,0,z)$ and a point $Y=(R\cos\theta,R\sin\theta,Z)$, \eqref{eq:rW2} yields
\[
\mathcal K_W(x,Y)= \tfrac{\partial}{\partial z}\big(K((0,0,z),Y)\cdot e_z\big) = -3\,\tfrac{R(z-Z)}{(R^2+(Z-z)^2)^{5/2}} .
\]
The angular weight $K_W(\sigma)$ in \eqref{eq:KW-kernel} is obtained from
$\mathcal K_W(0,Y)$ after writing $Y=(\rho,\sigma,\varphi)$.  We now pair the point $(R,Z)$, $Z>0$, with
the reflected point $(R,-Z)$.  Since
$\omega_{\rm hom,\theta}(R,Z)=-R^\alpha$ and
$\omega_{\rm hom,\theta}(R,-Z)=R^\alpha$, we have that 
\[
\int_0^\infty\tfrac{3R^{\alpha+2}(Z-z)} {(R^2+(Z-z)^2)^{5/2}}\,dR = C_\alpha^W\operatorname{sgn}(Z-z)|Z-z|^{\alpha-1},
\ \ \ \int_0^\infty \tfrac{3R^{\alpha+2}(Z+z)}{(R^2+(Z+z)^2)^{5/2}}\,dR =C_\alpha^W(Z+z)^{\alpha-1},
\]
where
\[
C_\alpha^W=\int_0^\infty\tfrac{3\tau^{\alpha+2}}{(1+\tau^2)^{5/2}}\,d\tau>0 .
\]
Subtracting the value at $z=0$ leaves the one-dimensional axial integral
\[
\p_Z(U_{\rm hom})_Z(0,z)-\p_Z(U_{\rm hom})_Z(0,0) =
-C_\alpha^W\int_0^\infty \Big[\operatorname{sgn}(Z-z)|Z-z|^{\alpha-1}+(Z+z)^{\alpha-1} -2Z^{\alpha-1} \Big]\,dZ ,
\]
and this  integrand is integrable on $(0,\infty)$; hence, using the change of variables $Z=zY$, we obtain that
\begin{equation}
\p_Z(U_{\rm hom})_Z(0,z) - \p_Z(U_{\rm hom})_Z(0,0) = c_{\rm hom} z^\alpha, \qquad z>0 ,
\label{eq:homogeneous-axis-strain-defect}
\end{equation}
where
\[
c_{\rm hom}:= -C_\alpha^W\int_0^\infty \Big[ \operatorname{sgn}(Y-1)|Y-1|^{\alpha-1} +(Y+1)^{\alpha-1} -2Y^{\alpha-1} \Big]\,dY.
\]
We now compute this one-dimensional integral directly.  For $R>1$,
\begin{align*}
&\int_0^R
\Big[\operatorname{sgn}(Y-1)|Y-1|^{\alpha-1} +(Y+1)^{\alpha-1}   -2Y^{\alpha-1} \Big]\,dY
\\
&\qquad\qquad\qquad
= -\tfrac1\alpha+\tfrac{(R-1)^\alpha}{\alpha}+\tfrac{(R+1)^\alpha-1}{\alpha} -\tfrac{2R^\alpha}{\alpha}
=\tfrac{(R-1)^\alpha+(R+1)^\alpha-2R^\alpha-2}{\alpha}.
\end{align*}
The first three terms in the numerator are $o(1)$ after cancellation as $R\to\infty$, because
\[
(R-1)^\alpha+(R+1)^\alpha-2R^\alpha = R^\alpha\left[\left(1-\tfrac1R\right)^\alpha +\left(1+\tfrac1R\right)^\alpha-2\right]\to0 .
\]
Hence
\[
\int_0^\infty \Big[\operatorname{sgn}(Y-1)|Y-1|^{\alpha-1}+(Y+1)^{\alpha-1} -2Y^{\alpha-1}\Big]\,dY=-\tfrac2\alpha,
\]
from which we obtain that
\[
c_{\rm hom}=\tfrac{2C_\alpha^W}{\alpha}>0 .
\]
Since $\mathsf W_{\rm hom}(\zeta)=\p_Z(U_{\rm hom})_Z(0,\zeta)$, \eqref{eq:homogeneous-axis-strain-defect}
is the identity
\[
\mathsf W_{\rm hom}(\zeta)-\mathsf W_{\rm hom}(0)=c_{\rm hom}\zeta^\alpha,\qquad \zeta>0.
\]
Subtracting this identity at $\zeta=\zeta_1$ and $\zeta=\zeta_2$ yields
\[
\mathsf W_{\rm hom}(\zeta_2)-\mathsf W_{\rm hom}(\zeta_1) = c_{\rm hom}(\zeta_2^\alpha-\zeta_1^\alpha), \qquad 0\le \zeta_1<\zeta_2 .
\]
Therefore, after choosing $C_{\rm hom}\ge c_{\rm hom}$, the homogeneous defect satisfies
\begin{equation}
\tfrac34c_{\rm hom}(\zeta_2^\alpha-\zeta_1^\alpha) \le \bigl(\mathsf W_{\rm hom}(\zeta_2)-\mathsf W_{\rm hom}(\zeta_1)\bigr) \le C_{\rm hom}
\bigl(\zeta_2^\alpha-\zeta_1^\alpha+\zeta_2^2-\zeta_1^2\bigr)
\label{eq:homogeneous-axis-strain-increment}
\end{equation}
for $0\le\zeta_1<\zeta_2\le\zeta_{\rm def}$.

\runinhead{Step 2: The transported cusp field.}
For
\[
0\le\zeta\le\zeta_{\rm def},\qquad |\tau|\le C_0,\qquad R=\zeta\tau,\qquad Z=\zeta ,
\]
we write \eqref{eq:Omega-cusp-local-representation} as
\begin{equation}
\bs\Omega_{\cusp}(J^2R,J^2Z,t) = -\Gamma J^{3\alpha-1}R^\alpha\mathfrak A_t(J^2R,J^2Z,t)\bs e_\theta +\bs\Omega_{\reg}(J^2R,J^2Z, t).
\label{eq:positive-defect-Omega-local}
\end{equation}
By \eqref{eq:full-physical-zeta-profile-def}, the axial amplitude on the symmetry axis is
\[
a_t^{\rm phys}(\zeta) = \bigl(J A_t(Z_t(\zeta))\bigr)^{1-\alpha} \bigl(1+Z_t(\zeta)^2\bigr)^{-\gamma/2}.
\]
The normal-form and radial-flatness bounds
\eqref{eq:radial-flatness-current-source-ratio}, \eqref{eq:Atr-local-bounds}, together with the bootstrap assumption
\eqref{eq:localized-normal-form-map-large-bootstrap}, show that
\begin{equation}
\mathfrak A_t(J^2\zeta\tau,J^2\zeta,t) = a_t^{\rm phys}(\zeta) +O\bigl(J^{3\beta_{\rm ax}}+\zeta_{\rm def}^{\alpha}\bigr),
\qquad 0\le\zeta\le\zeta_{\rm def},\quad |\tau|\le C_0 .
\label{eq:positive-defect-amplitude-local}
\end{equation}
For $Z=\zeta\ge0$, the homogeneous vorticity from Step~1 is $\omega_{{\rm hom},\theta}(R,Z)=-R^\alpha$.
Combining \eqref{eq:positive-defect-Omega-local}, \eqref{eq:full-physical-zeta-profile-def}, and
\eqref{eq:positive-defect-amplitude-local}, and then dividing by $\Gamma J^{3\alpha-1}$, we obtain that
\[
(\Gamma J^{3\alpha-1})^{-1}(\bs\Omega_{\cusp})_\theta(J^2R,J^2Z,t)= \omega_{{\rm hom},\theta}(R,Z) + \Omega_{{\rm rem},
\theta}^{\rm hom}(R,Z,t),
\]
where
\begin{align*}
\Omega_{{\rm rem},\theta}^{\rm hom}(R,Z,t)
&= -R^\alpha\bigl(\mathfrak A_t(J^2R,J^2Z,t)-1\bigr)
+(\Gamma J^{3\alpha-1})^{-1}(\bs\Omega_{\reg})_\theta(J^2R,J^2Z,t).
\end{align*}
Applying the Biot--Savart law in the variables $Y=(R,Z)$, we have that, for $0\le s\le\zeta_{\rm def}$,
\[
\mathsf W_t(s) = \p_Z\bigl(\BS[\omega_{{\rm hom},\theta}e_\theta]\bigr)_Z(0,s) + \p_Z\bigl(\BS[\Omega_{{\rm rem},
\theta}^{\rm hom}e_\theta]\bigr)_Z(0,s).
\]
This is \eqref{eq:Wt-homogeneous-leading-part}, with
\[
\mathsf W_{\rm hom}(s) = \p_Z\bigl(\BS[\omega_{{\rm hom},\theta}e_\theta]\bigr)_Z(0,s), \qquad \mathcal R_t^{\rm hom}(s) =
\p_Z\bigl(\BS[\Omega_{{\rm rem},\theta}^{\rm hom}e_\theta]\bigr)_Z(0,s).
\]
The axis bounds \eqref{eq:entry-axis-bounds-statement}, the containment assumptions
\eqref{eq:axis-attached-image-stop-bootstrap}--\eqref{eq:radial-flatness-buffered-label}, the normal-form estimates
\eqref{eq:localized-normal-form-large-bootstrap}--\eqref{eq:localized-normal-form-map-large-bootstrap}, and
the axis trace estimate \eqref{eq:tail-axis-error-bound} show that, for
$0\le\zeta_1<\zeta_2\le\zeta_{\rm def}$,
\[
|\mathcal R_t^{\rm hom}(\zeta_2)-\mathcal R_t^{\rm hom}(\zeta_1)|
\le C\bigl(J^{3\beta_{\rm ax}}+J^{1-3\alpha}+\zeta_{\rm def}^{\alpha}\bigr) \bigl(\zeta_2^\alpha-\zeta_1^\alpha+\zeta_2^2-\zeta_1^2\bigr).
\]
By \eqref{eq:Wt-homogeneous-leading-part}, this is the same as the estimate
\begin{equation}
\left| \bigl(\mathsf W_t(\zeta_2)-\mathsf W_t(\zeta_1)\bigr)
- \bigl(\mathsf W_{\rm hom}(\zeta_2)-\mathsf W_{\rm hom}(\zeta_1)\bigr) \right|
\le C\bigl(J^{3\beta_{\rm ax}}+J^{1-3\alpha}+\zeta_{\rm def}^{\alpha}\bigr) \bigl(\zeta_2^\alpha-\zeta_1^\alpha+\zeta_2^2-\zeta_1^2\bigr).
\label{eq:positive-defect-W-hom-comparison}
\end{equation}
The small terms in \eqref{eq:positive-defect-W-hom-comparison} are controlled by
\[
\begin{aligned}
\eqref{eq:localized-normal-form-map-large-bootstrap}
&\quad\Longrightarrow\quad
\|D\Psi_t-I\|_{L^\infty} +[\Psi_t-\operatorname{Id}]_{C^{\beta_{\rm ax}}} \le C J^{3\beta_{\rm ax}},
\\
\eqref{eq:Wcusp-asymptotic-field}
&\quad\Longrightarrow\quad
(\Gamma J^{3\alpha-1})^{-1}\Gamma=J^{1-3\alpha},
\\
\eqref{eq:positive-defect-amplitude-local}
&\quad\Longrightarrow\quad \sup_{\substack{0\le\zeta\le\zeta_{\rm def}\\ |\tau|\le C_0}} \bigl|\mathfrak A_t(J^2\zeta\tau,J^2\zeta,t)-a_t^{\rm phys}(\zeta)\bigr|
\le C\bigl(J^{3\beta_{\rm ax}}+\zeta_{\rm def}^{\alpha}\bigr),
\\
\eqref{eq:full-physical-zeta-profile-def},\
\eqref{eq:bounded-core-axis-origin-values},\
\eqref{eq:entry-axis-bounds-statement}
&\quad\Longrightarrow\quad
\sup_{0\le\xi\le\zeta_{\rm def}} \bigl|a_t^{\rm phys}(\xi)-1\bigr| \le C\zeta_{\rm def}^{\alpha}.
\end{aligned}
\]

We first choose $\zeta_{\rm def}$ so that
\[
C\zeta_{\rm def}^{\alpha} \bigl(\zeta_2^\alpha-\zeta_1^\alpha+\zeta_2^2-\zeta_1^2\bigr)
\le \tfrac18c_{\rm hom}(\zeta_2^\alpha-\zeta_1^\alpha),
\]
and then choose $\mathfrak J_{\rm def}$ so that, for $J\le\mathfrak J_{\rm def}$,
\[
C\bigl(J^{3\beta_{\rm ax}}+J^{1-3\alpha}\bigr) \bigl(\zeta_2^\alpha-\zeta_1^\alpha+\zeta_2^2-\zeta_1^2\bigr)
\le \tfrac18c_{\rm hom}(\zeta_2^\alpha-\zeta_1^\alpha).
\]
Combining the resulting estimate with \eqref{eq:homogeneous-axis-strain-increment} proves
\eqref{eq:positive-axial-strain-defect}, with $c_{\rm def}= \tfrac{1}{2} c_{\rm hom}$ after renaming constants.
Taking $\zeta_1=0$ yields \eqref{eq:positive-axial-strain-defect-origin}.
\end{proof}

\begin{lemma}[Improving the monotone axial stretching bootstrap]
\label{lem:monotone-axial-stretching-improvement}
Let $I_{\rm mon}=[0,\zeta_{\rm mon}]$ with $0<\zeta_{\rm mon}\le\zeta_{\rm def}$.  Assume the hypotheses
of Lemma~\ref{lem:axis-profile-evolution} on $I_{\rm mon}$, the exact axial composition identities of
Lemma~\ref{lem:exact-axis-composition}, the modulation bounds \eqref{eq:modulation-bounds}, the clock bound
\eqref{eq:Jdot-two-sided}, the defect estimate \eqref{eq:positive-axial-strain-defect}, and
\eqref{eq:tail-axis-error-strain-increment} with $I=I_{\rm mon}$.  Then, after decreasing the small-clock
threshold, there are constants
\[
0<c_{\rm mon}'\le C_{\rm mon}'<\infty,\qquad B_{\rm mon}'<\infty,
\]
depending only on the fixed parameters, such that
\begin{equation}
0<c_{\rm mon}'\le b_t(\zeta)\le C_{\rm mon}', \qquad 0\le\zeta\le\zeta_{\rm mon},
\label{eq:monotone-axial-two-sided-improved}
\end{equation}
and, for $0\le\zeta_1<\zeta_2\le\zeta_{\rm mon}$,
\begin{equation}
0\le\log b_t(\zeta_2)-\log b_t(\zeta_1)\le B_{\rm mon}' \bigl(\zeta_2^\alpha-\zeta_1^\alpha+\zeta_2^2-\zeta_1^2\bigr).
\label{eq:monotone-axial-fractional-improved}
\end{equation}
The constants in \textup{(BA4)} are chosen so that
\[
c_{\rm mon}<c_{\rm mon}'\le C_{\rm mon}'<C_{\rm mon}, \qquad B_{\rm mon}'>0,\qquad B_{\rm mon}'<B_{\rm mon}.
\]
\end{lemma}

\begin{proof}[Proof of Lemma~\ref{lem:monotone-axial-stretching-improvement}]
We write $D(\zeta_1,\zeta_2) := \zeta_2^\alpha-\zeta_1^\alpha+\zeta_2^2-\zeta_1^2$, and we
let $\zeta_i(t)=\mathscr Z_t(\eta_i)$, $i=1,2$, denote two trajectories of
\eqref{eq:exact-axis-Zflow} which remain in $I_{\rm mon}$.  From
\eqref{eq:axis-profile-b-system} and \eqref{eq:exact-axis-Zflow}, we have that
\begin{equation}
\begin{aligned}
\frac{d}{dt}\bigl[\log b_t(\zeta_2(t))-\log b_t(\zeta_1(t))\bigr]
&= m(t)\Gamma J^{3\alpha-1}\bigl[ \mathsf W_t(\zeta_2(t))-\mathsf W_t(\zeta_1(t))\bigr]\\
&\qquad\qquad
+(\p_zV_{\err})_z(0,J^2\zeta_2(t),t) -(\p_zV_{\err})_z(0,J^2\zeta_1(t),t).
\end{aligned}
\label{eq:monotone-b-increment-evolution}
\end{equation}
The first term on the right-hand side is nonnegative by
\eqref{eq:positive-axial-strain-defect} and the lower bound for $m(t)$ in \eqref{eq:modulation-bounds}; more
precisely,
\[
m(t)\Gamma J^{3\alpha-1}\bigl[ \mathsf W_t(\zeta_2(t))-\mathsf W_t(\zeta_1(t))\bigr]
\ge c_m c_{\rm def}\Gamma J^{3\alpha-1}\bigl(\zeta_2(t)^\alpha-\zeta_1(t)^\alpha\bigr).
\]
The error trace in \eqref{eq:monotone-b-increment-evolution} satisfies, by
\eqref{eq:tail-axis-error-strain-increment},
\[
\begin{aligned}
&\left|(\p_zV_{\err})_z(0,J^2\zeta_2(t),t)-(\p_zV_{\err})_z(0,J^2\zeta_1(t),t)\right|\\
&\qquad\le C\Gamma\bigl(J^{9\alpha-1}+1\bigr)D(\zeta_1(t),\zeta_2(t)).
\end{aligned}
\]
Since $0\le\zeta_1(t)<\zeta_2(t)\le\zeta_{\rm mon}$,
\[
D(\zeta_1(t),\zeta_2(t))
\le \left(1+\tfrac2\alpha\zeta_{\rm mon}^{2-\alpha}\right)
\bigl(\zeta_2(t)^\alpha-\zeta_1(t)^\alpha\bigr).
\]
We decrease the final small-clock threshold so that
\[
C\left(1+\tfrac2\alpha\zeta_{\rm mon}^{2-\alpha}\right)
\bigl(J^{6\alpha}+J^{1-3\alpha}\bigr)\le \tfrac12 c_m c_{\rm def}.
\]
Then the right-hand side of \eqref{eq:monotone-b-increment-evolution} is nonnegative.  Since $b_0\equiv1$,
the lower inequality in \eqref{eq:monotone-axial-fractional-improved} follows.

For the upper bound, \eqref{eq:positive-axial-strain-defect} and
\eqref{eq:tail-axis-error-strain-increment} imply that
\[
\Big|\frac{d}{dt}\bigl[ \log b_t(\zeta_2(t))-\log b_t(\zeta_1(t)) \bigr] \Big|
\le C\Gamma J^{3\alpha-1}D(\zeta_1(t),\zeta_2(t))+C\Gamma(J^{9\alpha-1}+1)D(\zeta_1(t),\zeta_2(t)).
\]
The flow equation \eqref{eq:exact-axis-Zflow}, the size bound \eqref{eq:axis-profile-size}, and
\eqref{eq:positive-axial-strain-defect-origin} show that
\[
D(\zeta_1(s),\zeta_2(s)) \le C D(\zeta_1(t),\zeta_2(t)) \qquad (s\le t)
\]
on the stopped interval.  By integrating in time and using \eqref{eq:Jdot-two-sided}, we obtain
\eqref{eq:monotone-axial-fractional-improved}.

Finally, taking $\zeta_1=0$ in \eqref{eq:monotone-axial-fractional-improved}, and using $b_t(0)=1$ from
\eqref{eq:axis-qb-volume-mon}, we obtain
\[
1\le b_t(\zeta) \le \exp\bigl(B_{\rm mon}'(\zeta_{\rm mon}^\alpha+\zeta_{\rm mon}^2)\bigr), \qquad 0\le\zeta\le\zeta_{\rm mon}.
\]
This proves \eqref{eq:monotone-axial-two-sided-improved} after renaming the constants.
\end{proof}

We next use the improved monotone axial-stretching bounds
\eqref{eq:monotone-axial-two-sided-improved}--\eqref{eq:monotone-axial-fractional-improved} to verify the
Euler-generated axial function and to activate the renormalized pressure bound for the cusp velocity.
Recall from \eqref{eq:pressure-profile-interval} and \eqref{eq:scaled-profile-at} that
$I_a=[0,\zeta_a]$ and
\[
a_t^{\rm phys}(\zeta) = \bigl(JA_t(Z_t(\zeta))\bigr)^{1-\alpha} \bigl(1+Z_t(\zeta)^2\bigr)^{-\gamma/2}, \qquad \zeta\in I_a.
\]
We use the zero extension
\begin{equation}
a_t(\zeta) := a_t^{\rm phys}(\zeta)\mathbf 1_{I_a}(\zeta) =
\bigl(JA_t(Z_t(\zeta))\bigr)^{1-\alpha} \bigl(1+Z_t(\zeta)^2\bigr)^{-\gamma/2}\mathbf 1_{I_a}(\zeta), \qquad B_t(Z_t(\zeta))=J^2\zeta,
\label{eq:euler-generated-pressure-profile-sec13}
\end{equation}
with $A_t(Z)=\p_Rr_t(0,Z)$ and $B_t(Z)=z_t(0,Z)$.  By
\eqref{eq:monotone-axial-two-sided-improved}--\eqref{eq:monotone-axial-fractional-improved},
\eqref{eq:axis-qb-volume-mon}, and
$ q_t(\zeta)=JA_t(Z_t(\zeta))$, $\p_\zeta Z_t(\zeta)=b_t(\zeta)^{-1}$,
we have that $q_t$ is nonincreasing and $Z_t$ is increasing on $I_a$.  Hence
$q_t(\zeta)^{1-\alpha}$ and $\bigl(1+Z_t(\zeta)^2\bigr)^{-\gamma/2}$ are nonincreasing on $I_a$, and
therefore $a_t$ is nonnegative and nonincreasing on $(0,\infty)$.

\begin{lemma}[Euler-generated renormalized Riccati bound]
\label{lem:euler-generated-slope-restricted-riccati}
Let $I_a=[0,\zeta_a]\subset I_{\rm mon}$, and let $a_t$ be the axial function in
\eqref{eq:euler-generated-pressure-profile-sec13}.  Assume
\eqref{eq:axis-qb-volume-mon} and
\eqref{eq:monotone-axial-two-sided-improved}--\eqref{eq:monotone-axial-fractional-improved} on $I_a$.
Assume also the fixed choices from Section~\ref{sec:fixed-choice-order} and the threshold relation defining
$\mathfrak J_{\mathrm{axis}}$.  Then $a_t$ is nonnegative and nonincreasing on $(0,\infty)$, and the cusp-flow
pressure Hessian satisfies
\begin{equation}
\Pi_{\cusp}(t) \ge -q_{\rm tr}\,\tfrac12\,\mathcal W_{\cusp}(t)^2, \qquad q_{\rm tr}<\upbeta .
\label{eq:euler-generated-riccati-closed}
\end{equation}
\end{lemma}

\begin{proof}[Proof of Lemma~\ref{lem:euler-generated-slope-restricted-riccati}]
The monotonicity proof is the argument in the paragraph preceding the lemma.  Indeed, by
\eqref{eq:axis-qb-volume-mon}, $q_t(\zeta)^2b_t(\zeta)=1$ and $q_t(\zeta)=JA_t(Z_t(\zeta))$.
The lower bound in \eqref{eq:monotone-axial-fractional-improved} implies that $b_t$ is nondecreasing on
$I_a$.  Hence $q_t(\zeta)=b_t(\zeta)^{-1/2}$ is nonincreasing on $I_a$.  Since
$\p_\zeta Z_t=b_t^{-1}>0$ and $Z_t(0)=0$, the map $\zeta\mapsto Z_t(\zeta)$ is increasing on $I_a$.
Thus both factors in \eqref{eq:euler-generated-pressure-profile-sec13} are nonincreasing, and the zero
extension preserves distributional monotonicity.

For the pressure bound, the improved axis bounds, the axis-profile evolution equations, the normal-form
estimates, and the cusp-error estimates are exactly
the hypotheses of Lemma~\ref{lem:renormalized-axis-trace-hypotheses}.  Therefore the renormalized axis-trace
hypotheses of Proposition~\ref{prop:euler-generated-profile-riccati} hold for the axial function $a_t$ after the final
small-clock threshold is decreased.  Lemma~\ref{lem:transported-cusp-pressure-win} then proves
\eqref{eq:euler-generated-riccati-closed}.
\end{proof}

Lemma~\ref{lem:monotone-axial-stretching-improvement} closes \textup{(BA4)} with strict constants.  The
Riccati estimate in Lemma~\ref{lem:euler-generated-slope-restricted-riccati} is then a consequence of the
closed monotone bounds on $I_a$.  Together with the preceding
improvement lemmas, all size bootstraps have now been improved with strict constants.

We now collect the small-clock consequences used in
Section~\ref{sec:target-profile-typeI-completion}.  The auxiliary coordinates
$\eta$, $\mathscr Z_t$, $\widehat q_t$, and $\widehat b_t$ have served their purpose: they provide uniform
control of the exact axial map and the Euler-generated axial function.  We use the following consequences.
The cusp clock is the singular clock, while $m(t)$ and $J_{\smooth}(t)$ remain order one; the total
axial strain is the sum of a singular cusp contribution and a bounded smooth contribution; and the axial function
\[
a_t(\zeta) = \bigl(J_{\cusp}(t)A_t(Z_t(\zeta))\bigr)^{1-\alpha} \bigl(1+Z_t(\zeta)^2\bigr)^{-\gamma/2}\mathbf 1_{I_a}(\zeta)
\]
is nonnegative, nonincreasing, and satisfies the renormalized cusp-pressure Riccati estimate.  Indeed,
\eqref{eq:entry-axis-bounds-statement} and \eqref{eq:physical-zeta-profile-envelope} imply the size bounds
\eqref{eq:current-axis-extension-size}; the identities $q_t^2b_t=1$ and $\p_\zeta Z_t=b_t^{-1}>0$, together with
\eqref{eq:monotone-axial-fractional-improved}, imply monotonicity; and Lemma~\ref{lem:transported-cusp-pressure-win} proves
\eqref{eq:small-clock-model-riccati} and \eqref{eq:small-clock-cusp-riccati}.  These are precisely the flow-map facts needed later to compare
$J(t)$, $\rW_0(t)$, and $\Pi_0(t)$ with their model counterparts.

We now keep the choices from Section~\ref{sec:fixed-choice-order} fixed:
\[
(I_\sharp,\vartheta_\sharp,I_{\rm loc}^{\rm cur},I_{\rm buf}^{\rm cur}), \qquad \mathfrak C_{\rm fix},\qquad \vartheta_{\pressure},\qquad M_\pressure,\qquad
\sigma_{\inn},\sigma_*,\qquad C_T^{\rm fix},\qquad R_{\tail},
\]
as in
\eqref{eq:pressure-localization-cutoff}--\eqref{eq:pressure-localization-intervals}, \eqref{eq:fixed-bootstrap-barriers}, \eqref{eq:pressure-angular-tail},
\eqref{eq:pressure-C0}, \eqref{eq:fixed-smooth-time-horizon}, and \eqref{eq:core-tail-domains}--\eqref{eq:smooth-velocity-def}.  We then choose the final axis threshold with
\begin{equation}
\mathfrak J_{\mathrm{axis}}\le \min\{\mathfrak J_{\mathrm{mod}},\mathfrak J_{\transport},\mathfrak J_{\pressure}\}.
\label{eq:small-clock-comparison-threshold}
\end{equation}
The following proposition assumes these choices and this threshold relation.

\begin{proposition}[Small-clock comparisons]
\label{prop:small-clock-comparisons}
Assuming these choices, the following assertions hold on
\[
J_{\cusp}(t)\le \mathfrak J_{\mathrm{axis}},
\]
with constants depending only on the fixed parameters.

\begin{enumerate}[label=\textup{(\roman*)}]
\item The clock and axial strain comparisons are
\[
c_m\le m(t)\le C_m,\qquad c_1\Gamma J_{\cusp}(t)^{3\alpha} \le-\dot J_{\cusp}(t)\le C_1\Gamma J_{\cusp}(t)^{3\alpha},
\]
\[
|\rW_{\smooth}(t)|\le C_{\smooth}\Gamma,\qquad 0<c_{\smooth}\le J_{\smooth}(t)\le C_{\smooth},
\]
together with
\[
J(t)=J_{\smooth}(t)J_{\cusp}(t),\qquad \rW_{\cusp}(t)=m(t)\mathcal W_{\cusp}(t),\qquad \rW_0(t)=\rW_{\smooth}(t)+\rW_{\cusp}(t).
\]

\item The exact axial composition $\mathscr Z_t$ is uniformly bi-Lipschitz on an origin-attached reference
interval $I_{\rm ax}$:
\begin{equation}
c_{\rm ax}\le \p_\eta\mathscr Z_t(\eta)\le C_{\rm ax}, \qquad [\log\p_\eta\mathscr Z_t]_{C^{\alpha/2}(I_{\rm ax})}\le C_{\rm ax},
\label{eq:axis-S-distortion}
\end{equation}
and
\begin{equation}
c_{\rm ax}\eta\le \mathscr Z_t(\eta)\le C_{\rm ax}\eta, \qquad \eta\in I_{\rm ax}.
\label{eq:axis-S-bilip-from-axis}
\end{equation}
For a fixed axial label $Z$,
\begin{equation}
c\,J_{\cusp}(t)^{-2}B_t(Z) \le J_{\cusp}(s)^{-2}B_s(Z) \le C\,J_{\cusp}(t)^{-2}B_t(Z), \qquad t_0\le s\le t,
\label{eq:axis-attached-buffered-label}
\end{equation}
whenever both clock-scaled axial positions in \eqref{eq:axis-attached-buffered-label} lie in the
origin-attached interval where the axis estimates are applied.

\item The pressure-localization interval is contained in the exact axial image.  If
\[
\operatorname{supp}\vartheta_\sharp\Subset K_\vartheta=[\zeta_-,\zeta_+]\Subset I_\sharp,
\qquad I_{\rm loc}=[\eta_-,\eta_+]\Subset I_{\rm ax}\cap(0,\infty),
\]
then $I_{\rm loc}$ is fixed so that
\begin{equation}
C_{\rm ax}\eta_-<\zeta_-<\zeta_+<c_{\rm ax}\eta_+, \qquad [c_{\rm ax}\eta_-,\,C_{\rm ax}\eta_+]\Subset I_\sharp ,
\label{eq:axis-control-endpoint-margin}
\end{equation}
and hence
\begin{equation}
\operatorname{supp}\vartheta_\sharp \Subset \mathscr Z_t(I_{\rm loc}) \Subset I_\sharp.
\label{eq:axis-control-covers-pressure-interval}
\end{equation}
We also fix $\zeta$-intervals satisfying
\[
\operatorname{supp}\vartheta_\sharp \Subset I_{\rm loc}^{\rm cur} \Subset I_{\rm buf}^{\rm cur}
\Subset K_\vartheta \Subset \mathscr Z_t(I_{\rm loc}) \Subset \mathscr Z_t(I_{\rm ax}).
\]

\item The conserved functions $\widehat q_t,\widehat b_t$ from
\eqref{eq:exact-axis-renormalized-functions} satisfy
\eqref{eq:renormalized-axis-chart-bootstrap}.
Consequently, for every compact $\zeta$-interval $K_\zeta\Subset\mathscr Z_t(I_{\rm ax})$,
\begin{equation}
c_{\rm ax}\le q_t(\zeta),\,b_t(\zeta)\le C_{\rm ax},\qquad [\log q_t]_{C^{\alpha/2}(K_\zeta)} +[\log b_t]_{C^{\alpha/2}(K_\zeta)} \le C_{\rm ax}.
\label{eq:current-axis-transfer-bound}
\end{equation}
In particular, \eqref{eq:entry-axis-bounds-statement} holds on
$I_{\rm loc}^{\rm cur}$, $I_{\rm buf}^{\rm cur}$, $I_\sharp$, and $I_{\err}$, after the reference interval
$I_{\rm ax}$ is fixed large enough.

\item For $a_t:=a_t^{\rm phys}\mathbf 1_{I_a}$, the Euler-generated axial function is
nonnegative and nonincreasing on $(0,\infty)$.  Moreover,
\begin{equation}
c_{1,\rm ax}\alpha^{-1}\le I_1[a_t]\le C_{1,\rm ax}\alpha^{-1}, \qquad 0\le a_t(\zeta)\le C_{\rm env,ax}(1+\zeta^2)^{-\gamma/2}.
\label{eq:current-axis-extension-size}
\end{equation}
The cusp-flow pressure Hessian satisfies
\begin{equation}
\Pi_{\cusp}(t)\ge -q_{\rm tr}\,\tfrac12\,\mathcal W_{\cusp}(t)^2,\qquad q_{\rm tr}<\upbeta.
\label{eq:small-clock-model-riccati}
\end{equation}
The estimate \eqref{eq:small-clock-model-riccati} is the renormalized axis-trace Riccati bound obtained from
Lemma~\ref{lem:transported-cusp-pressure-win}.  The final stagnation-point Riccati comparison uses the same
bound in \eqref{eq:small-clock-cusp-riccati}.
On the pressure-localization support,
\[
a_t^{\rm loc}(\zeta):=\vartheta_\sharp(\zeta)a_t(\zeta) =\vartheta_\sharp(\zeta)a_t^{\rm phys}(\zeta).
\]

\item For each fixed $Z_*<\infty$ and $C_0<\infty$, the bounded-core normal form holds uniformly for
$0<|Z|\le Z_*$ and $|A_t(Z)R/B_t(Z)|\le C_0$:
\begin{equation}
r_t(R,Z)=A_t(Z)R+\mathscr R_{r,t}(R,Z), \qquad z_t(R,Z)=B_t(Z)+\mathscr R_{z,t}(R,Z),
\label{eq:bounded-core-normal-form}
\end{equation}
with
\begin{equation}
|\mathscr R_{r,t}(R,Z)|+|\mathscr R_{z,t}(R,Z)| \le C J_{\cusp}(t)^{-1}|R|^{1+\beta_{\rm ax}} .
\label{eq:bounded-core-normal-form-error}
\end{equation}

\item The cusp-flow pressure Hessian satisfies the one-sided Riccati bound
\begin{equation}
\Pi_{\cusp}(t) \ge -q_{\rm tr}\,\tfrac12\,\mathcal W_{\cusp}(t)^2, \qquad q_{\rm tr}<\upbeta.
\label{eq:small-clock-cusp-riccati}
\end{equation}
(This restates \eqref{eq:small-clock-model-riccati} with the notation used in the final proof of
Theorem~\ref{thm:target-profile}.)
\end{enumerate}
\end{proposition}

\begin{proof}[Proof of Proposition~\ref{prop:small-clock-comparisons}]
The clock, modulation, smooth-clock, and strain-splitting assertions are
\eqref{eq:modulation-bounds}, \eqref{eq:Jdot-two-sided}, \eqref{eq:Jsmooth-bdd}, and
\eqref{eq:on-axis-strain-splitting}.  The estimates for $\mathscr Z_t$, $q_t$, and $b_t$ are the closed
axial flow map estimates obtained from Lemma~\ref{lem:exact-axis-composition} and
\eqref{eq:exact-axis-renormalized-conserved}.  The value estimate for $\mathscr Z_t$ yields
\eqref{eq:axis-attached-buffered-label} and the pressure-support coverage
\eqref{eq:axis-control-covers-pressure-interval} from the fixed margin
\eqref{eq:axis-control-endpoint-margin}.

The nonnegativity, monotonicity, and Riccati assertion for $a_t$ follow from
Lemma~\ref{lem:euler-generated-slope-restricted-riccati}.  The $I_1$ bounds follow from
\eqref{eq:current-axis-transfer-bound}, \eqref{eq:axis-qb-volume-mon}, and the fixed interval
$I_a=[0,\zeta_a]$; the pointwise upper bound is \eqref{eq:physical-zeta-profile-envelope}.  The bounded-core normal form is
the closed origin-attached normal-form estimate for the exact cusp map.  The pressure Hessian bound
\eqref{eq:small-clock-model-riccati} is \eqref{eq:euler-generated-riccati-closed}; the identical bound
\eqref{eq:small-clock-cusp-riccati} is Lemma~\ref{lem:transported-cusp-pressure-win}.  We choose
$\mathfrak J_{\mathrm{axis}}$ no larger than the thresholds in all cited estimates.
\end{proof}

\subsection{Closure of the size bootstrap assumptions}
\label{sec:size-bootstrap-closure}

Lemma~\ref{lem:monotone-axial-stretching-improvement} closes the monotone axial-stretching bootstrap
\textup{(BA4)}, namely \eqref{eq:monotone-axial-two-sided} and
\eqref{eq:monotone-axial-fractional-bootstrap}.  With these improved bounds available on $I_a$, 
Lemma~\ref{lem:euler-generated-slope-restricted-riccati} applies the renormalized axis-trace criterion to the
Euler-generated axial function \eqref{eq:euler-generated-pressure-profile-sec13} and proves the cusp-flow
pressure Hessian bound \eqref{eq:euler-generated-riccati-closed}.

The individual improvement lemmas for the size bootstraps
\eqref{eq:simultaneous-continuation-assumptions} are now in place: Lemmas~\ref{lem:nonlinear-radial-flatness},
\ref{lem:late-axis-normal-form-cusp}, and \ref{lem:late-axis-normal-form-map-cusp} improve (BA6) and (BA7);
Lemma~\ref{lem:tail-bound} improves (BA8); Lemmas~\ref{lem:modulation-bounded} and
\ref{lem:Jdot-two-sided-aux} improve (BA9) and (BA5); and
Proposition~\ref{prop:small-clock-comparisons} improves (BA2) and (BA3).  We combine these individual
improvements into a single open--closed continuation statement that extends $\mathcal B_{\rm size}$
\eqref{eq:simultaneous-continuation-assumptions} over the entire small-clock interval.

\begin{proposition}[Improving the size bootstrap assumptions]
\label{prop:simultaneous-small-clock-continuation}
There exists a threshold $\mathfrak J_{\mathrm{cont}}>0$ such that the following holds.  Let $I$ be a
small-clock time interval on which $J_{\cusp}(t)\le\mathfrak J_{\mathrm{cont}}$.  Assume that every size
bootstrap in $\mathcal B_{\rm size}$ holds for every
$t\in I$.  Then each of (BA2), (BA3), (BA5), (BA6), (BA7), (BA8), and (BA9) improves with strict margin on
$I$: the normal-form norm $\mathfrak B$, the map-distortion norm $\mathfrak D_\Psi$, and the cusp-error norm
$\mathfrak E_{\err}$ improve to
\[
\mathfrak B\le \tfrac12 B_*, \qquad \mathfrak D_\Psi\le \tfrac12 D_*, \qquad \mathfrak E_{\err}\le \tfrac12 E_*,
\]
the axial flow map bounds in (BA2), the cusp-clock rate bound in (BA5), and the modulation bounds in (BA9)
hold with constants strictly inside the chosen bootstrap constants, and the interval-containment
conditions in (BA3) close with strict margin.  Consequently $\mathcal B_{\rm size}$ holds on the entire
small-clock interval by the open--closed argument; in particular, the
pressure-support coverage \eqref{eq:axis-control-covers-pressure-interval} holds throughout that interval.
\end{proposition}

\begin{proof}[Proof of Proposition~\ref{prop:simultaneous-small-clock-continuation}]
We choose $\mathfrak J_{\mathrm{cont}}$ below every small-clock threshold appearing in the individual
improvement lemmas cited in this proof.  The maps, scalar functions, and norms controlled by $\mathcal B_{\rm size}$
are continuous as long
as the $C^{1,\alpha}$ Euler solution exists and $J_{\cusp}>0$.

We work on the maximal subinterval on which every size bootstrap in $\mathcal B_{\rm size}$ holds.  The
interval-containment conditions in (BA3),
\eqref{eq:axis-attached-image-stop-bootstrap}, \eqref{eq:current-axis-anchor-bootstrap}, and
\eqref{eq:radial-flatness-buffered-label}, are compact-containment statements; the fixed margins built into
them make them open by continuity, and the inclusions proved in the next paragraphs close them.  The lower bound in the
cusp-clock rate (BA5), \eqref{eq:localized-clock-bootstrap}, implies that $J_{\cusp}$ is decreasing on $I$ and
that, for every nonnegative function $F$ of the clock,
\[
\int_{t_1}^{t_2}\Gamma F(J_{\cusp}(t))\,dt \le C\int_{J_{\cusp}(t_2)}^{J_{\cusp}(t_1)} F(J)J^{-3\alpha}\,dJ .
\]
This is the conversion from time to clock used in the earlier improvement lemmas.

The normal-form bootstrap (BA6) is improved by
Lemmas~\ref{lem:nonlinear-radial-flatness} and \ref{lem:late-axis-normal-form-cusp}; the normal-form map
distortion bootstrap (BA7) is improved by Lemma~\ref{lem:late-axis-normal-form-map-cusp}.  The scalar
modulation bootstrap (BA9) is improved by Lemma~\ref{lem:modulation-bounded}, and the cusp-clock rate
bootstrap (BA5) is improved by Lemma~\ref{lem:Jdot-two-sided-aux}.  The cusp-error bootstrap (BA8) is
improved by Lemma~\ref{lem:tail-bound}.

It remains only to explain how Proposition~\ref{prop:small-clock-comparisons} closes (BA2) and (BA3).  The
closed axial flow map estimates are \eqref{eq:renormalized-axis-chart-bootstrap} and
\eqref{eq:current-axis-transfer-bound}, which are the strict-improvement form of (BA2) after the bootstrap
barriers have been fixed as in Section~\ref{sec:fixed-choice-order}.  For (BA3), the value estimate
\eqref{eq:axis-S-bilip-from-axis} closes
\[
\mathscr Z_t(I_{\rm ax})\subset I_{\rm ax}^{\zeta}
\]
after the interval $I_{\rm ax}^{\zeta}$ is fixed with margin.  Since $I_{\rm str}\Subset I_\sharp$,
\eqref{eq:current-axis-transfer-bound} and
\[
\p_\zeta Z_t(\zeta)=b_t(\zeta)^{-1}
\]
place $Z_t(I_{\rm str})$ in a fixed compact subinterval of $(0,R_{\tail})$; this is
\eqref{eq:current-axis-anchor-bootstrap}.  Finally, the fixed-label containment
\eqref{eq:radial-flatness-buffered-label} follows from \eqref{eq:axis-attached-buffered-label}.  Together
with the fixed endpoint margin \eqref{eq:axis-control-endpoint-margin}, these estimates also imply the
pressure-support coverage \eqref{eq:axis-control-covers-pressure-interval}.

The monotone axial stretching bounds are supplied by
Lemma~\ref{lem:monotone-axial-stretching-improvement}.  The corresponding pressure Hessian and axial strain comparison for
\[
a_t(\zeta) = \bigl(JA_t(Z_t(\zeta))\bigr)^{1-\alpha} \bigl(1+Z_t(\zeta)^2\bigr)^{-\gamma/2}\mathbf 1_{I_a}(\zeta)
\]
is \eqref{eq:euler-generated-riccati-closed}.  This is the slope-restricted pressure Hessian information used in
Lemma~\ref{lem:transported-cusp-pressure-win}.

The fixed cutoffs, cone parameters, tail radius, bootstrap constants, and small-clock thresholds are chosen
in the order specified in Section~\ref{sec:fixed-choice-order}.  With those choices fixed, each strict estimate
cited in this proof holds at the boundary of the maximal subinterval.  The set of times on which $\mathcal B_{\rm size}$
holds is therefore nonempty at the entry time, open by continuity, and closed by these strict
improvements; the standard open--closed continuation argument then extends every size bootstrap in
$\mathcal B_{\rm size}$ to the whole small-clock interval.
\end{proof}

After Lemma~\ref{lem:transported-cusp-pressure-win}, all further pressure comparisons are obtained only by
shrinking the already fixed small-clock threshold.  Once the estimates for
$\widehat q_t,\widehat b_t$, the modulation bound \eqref{eq:modulation-bounds}, the transported-cusp bounds
\eqref{eq:Wcusp-scaling}--\eqref{eq:Ucusp-radial-defect}, and the cusp pressure estimate
\eqref{eq:small-clock-cusp-riccati} are active, every remaining term is measured against the leading scales
\[
\begin{gathered}
c\,\Gamma J_{\cusp}(t)^{3\alpha-1}
\le
|\mathcal W_{\cusp}(t)|
\le
C\,\Gamma J_{\cusp}(t)^{3\alpha-1},
\\
c\,\Gamma^2J_{\cusp}(t)^{6\alpha-2}
\le
\mathcal W_{\cusp}(t)^2
\le
C\,\Gamma^2J_{\cusp}(t)^{6\alpha-2}.
\end{gathered}
\]
The remaining smallness requirements are positive powers of $J_{\cusp}$.  Thus, after the cutoffs and
bootstrap constants have been fixed, the active small-clock threshold is the minimum of the thresholds ordered in
\eqref{eq:small-clock-threshold-order} and of the thresholds appearing in
\eqref{eq:small-clock-comparison-threshold}.  Any later change is a further decrease, so all previously activated
estimates remain valid.  Consequently there is no separate loss-of-control scenario inside the small-clock regime: while a
$C^{1,\alpha}$ solution exists and $J_{\cusp}>0$, the maps, scalar functions, and norms in
$\mathcal B_{\rm size}$ are continuous,
and Proposition~\ref{prop:simultaneous-small-clock-continuation} improves each of them before it can reach its
bootstrap boundary.

\section{Euler Blowup for the Target Datum}
\label{sec:target-profile-typeI-completion}

We now prove Theorem~\ref{thm:target-profile}.  The proof uses the closed small-clock estimates from
Proposition~\ref{prop:small-clock-comparisons}.  The pressure estimate used directly in the final
stagnation-point Riccati comparison is the cusp-flow bound \eqref{eq:small-clock-cusp-riccati}:
\begin{equation*}
\Pi_{\cusp}(t)\ge -q_{\rm tr}\,\tfrac12\,\mathcal W_{\cusp}(t)^2,\qquad q_{\rm tr}<\upbeta .
\end{equation*}
The slope-restricted model estimate \eqref{eq:small-clock-model-riccati} also appears in
Proposition~\ref{prop:small-clock-comparisons} because it is the local Riccati bound proved for the
Euler-generated axial function.  Lemma~\ref{lem:transported-cusp-pressure-win} transfers
\eqref{eq:small-clock-model-riccati} to the cusp-flow estimate \eqref{eq:small-clock-cusp-riccati}.  After
this transfer, we use \eqref{eq:small-clock-cusp-riccati}, add the lower-order pressure Hessian terms, apply
the stagnation-point Riccati identity, prove finite-time collapse of the cusp clock, and convert the clock law
into the Type--I vorticity rate stated in Theorem~\ref{thm:target-profile}.

\subsection{Modulated pressure Hessian comparison}

The estimate \eqref{eq:small-clock-cusp-riccati} from
Proposition~\ref{prop:small-clock-comparisons} is written for the cusp-coordinate velocity $U_{\cusp}$,
before the scalar modulation $m(t)$ is applied:
\begin{equation*}
\Pi_{\cusp}(t)\ge -q_{\rm tr}\,\tfrac12\,\mathcal W_{\cusp}(t)^2,\qquad q_{\rm tr}<\upbeta .
\end{equation*}
Thus $\mathcal W_{\cusp}(t)$ is the axial strain of $U_{\cusp}$ at the origin, while the cusp part of the
Euler strain is
\begin{equation}
\rW_{\cusp}(t)=m(t)\mathcal W_{\cusp}(t).
\label{eq:modulated-cusp-strain-scaling-sec14}
\end{equation}
Since the pressure Hessian is quadratic in the velocity gradient, multiplying $U_{\cusp}$ by $m(t)$
multiplies the corresponding pressure Hessian by $m(t)^2$.  Hence
\eqref{eq:small-clock-cusp-riccati} gives
\[
m(t)^2\Pi_{\cusp}(t)\ge -q_{\rm tr}\,\tfrac12\,\rW_{\cusp}(t)^2 .
\]
We then prove the full pressure estimate with $\rW_{\cusp}(t)^2$ in the Riccati term,
\[
\Pi_0(t)\ge -q_{\rm phys}\,\tfrac12\,\rW_{\cusp}(t)^2,\qquad q_{\rm phys}<\upbeta,
\]
by adding the geometric, mixed, smooth, and error terms in the pressure decomposition \eqref{eq:Pi0-decomp}.
The final passage from $\rW_{\cusp}(t)$ to the exact Euler strain $\rW_0(t)$ is made in
Proposition~\ref{prop:blowup-final} using the splitting \eqref{eq:on-axis-strain-splitting}.

We first isolate the geometric remainder in \eqref{eq:Pi0-decomp}.  With the bilinear form
$\Pi[\cdot,\cdot]$ defined by \eqref{eq:physical-pressure-bilinear-form}, set
\begin{equation}
\Pi_{\rm geom}(t):=\Pi[u_{\cusp},u_{\cusp}](t)-m(t)^2\Pi[U_{\cusp},U_{\cusp}](t).
\label{eq:Pi-geom-sec14}
\end{equation}
By \eqref{eq:Pi-U-cusp-def}, $\Pi_{\cusp}(t):=\Pi[U_{\cusp},U_{\cusp}](t)$.  Recall from
\eqref{eq:Verr-def}--\eqref{eq:u-decomp} that
\[
V_{\err}=V_{\cusp}-m(t)U_{\cusp},\qquad u_{\cusp}=(\phi_{\smooth})_*(m(t)U_{\cusp}),\qquad u_{\err}=(\phi_{\smooth})_*V_{\err}.
\]
Thus $\Pi[u_{\cusp},u_{\cusp}]$ evaluates the pressure of the scalar-modulated cusp-coordinate velocity
$mU_{\cusp}$ after the push-forward by the smooth flow map $\phi_{\smooth}$, whereas
$m^2\Pi[U_{\cusp},U_{\cusp}]$ evaluates the same velocity in cusp-coordinate variables before
that push-forward.

At a fixed time we write
\begin{equation}
J:=J_{\cusp}(t),\qquad \Lambda(X):=\phi_{\smooth}(X,t),\qquad w(X):=m(t)U_{\cusp}(X,t).
\label{eq:smooth-pressure-time-notation-sec14}
\end{equation}
Then $u_{\cusp}=\Lambda_*w$ by \eqref{eq:u-decomp}.  We also use the localized ``flat'' velocity $U_\sharp$
from \eqref{eq:Omega-sharp-cutoff-def} and write
\[
w_\sharp:=mU_\sharp, \qquad w_{\rm out}:=m(U_{\cusp}-U_\sharp).
\]
The parameters $I_\sharp$, $\vartheta_\sharp$, $M_\pressure$, and $C_0$ are the same localization data used in
Lemma~\ref{lem:transported-cusp-pressure-win}.  The transported vorticity defining $U_\sharp$ is carried by
labels whose cusp-flow images have the form
\begin{equation}
X=\phi_{\cusp}(Y_t(\zeta,\tau),t)=J^2\zeta\bigl((\tau,1)+\mathcal E_t(\zeta,\tau)\bigr), \qquad \zeta\in I_\sharp,\quad 0\le\tau\le C_0 .
\label{eq:smooth-pressure-localized-tube-sec14}
\end{equation}
The normal-form displacement bound \eqref{eq:normal-form-approximation-bound} implies, after decreasing
$\mathfrak J_{\mathrm{axis}}$ if necessary, that the points in
\eqref{eq:smooth-pressure-localized-tube-sec14} lie in the ball
\begin{equation}
|X|\le R_{\rm tube}:=2\sup_{\zeta\in I_\sharp}\zeta\,(1+C_0), \qquad J\le\mathfrak J_{\mathrm{axis}}.
\label{eq:smooth-pressure-localized-ball-sec14}
\end{equation}
We choose $R_{\rm eval}\ge4R_{\rm tube}$ so that the pressure estimate
\eqref{eq:fixed-scale-pressure-bilinear} applies on $B_{R_{\rm eval}}$ in the variables obtained by dividing
the Eulerian image by $J^2$.  The velocity $w_{\rm out}=m(U_{\cusp}-U_\sharp)$ is the contribution of
$mU_{\cusp}$ from the complement of the localization in \eqref{eq:Omega-sharp-cutoff-def}.  This complement
contains the region where the axial cutoff $1-\vartheta_\sharp$ is active, the region where the angular cutoff
$1-\chi_{M_\pressure}$ is active, and the bounded-core labels whose cusp-flow images are outside the tube
\eqref{eq:smooth-pressure-localized-tube-sec14}.  We choose $R_{\rm loc}<\infty$ so that the images of these
three localized regions lie in $B_{R_{\rm loc}}$.  The remaining labels are estimated by the far-field and
algebraic-tail bounds in Lemma~\ref{lem:tail-bound}.  We set
\begin{equation}
R_{\rm pr}:=\max\{R_{\rm eval},R_{\rm loc}\}.
\label{eq:smooth-pressure-fixed-radius-sec14}
\end{equation}
The tail radius $R_{\tail}$ has been chosen so that Lemma~\ref{lem:smooth-flow-small-deformation} applies with
$R_0=R_{\rm pr}$ and $C_T=C_T^{\rm fix}$.  Hence the smooth-flow estimates
\eqref{eq:usmooth-small-C2a}--\eqref{eq:smooth-flow-second-gradient-small} hold on every ball used in
Lemma~\ref{lem:smooth-pressure-defect}.

The four small terms in Lemma~\ref{lem:smooth-pressure-defect},
\begin{equation*}
\varepsilon_{\smooth},\qquad \mathfrak a_{\zeta}(I_\sharp)^{{\frac{1}{2}}},\qquad
\mathfrak a_{\rm ang}(M_\pressure)^{{\frac{1}{2}}},\qquad J^{\kappa_{\rm def}}
\end{equation*}
have the following origins: the smooth-flow near-identity estimates
\eqref{eq:usmooth-small-C2a}--\eqref{eq:smooth-flow-second-gradient-small}, the axial localization tail
\eqref{eq:pressure-zeta-tail}, the angular tail \eqref{eq:pressure-angular-tail}, and the normal-form
displacement \eqref{eq:normal-form-approximation-bound}.  The final power $J^{1-3\alpha}$ below comes from
the algebraic tail and the far-label estimates in Lemma~\ref{lem:tail-bound}.

\begin{lemma}[Smooth-flow pressure deformation]
\label{lem:smooth-pressure-defect}
Let $R_{\rm tube}$ and $R_{\rm pr}$ be defined by
\eqref{eq:smooth-pressure-localized-ball-sec14} and \eqref{eq:smooth-pressure-fixed-radius-sec14}.  Suppose
that $R_{\tail}$ is fixed so that Lemma~\ref{lem:smooth-flow-small-deformation} applies with
$R_0=R_{\rm pr}$ and $C_T=C_T^{\rm fix}$.  Then there is $C<\infty$, depending only on the fixed parameters
from Subsection~\ref{sec:fixed-choice-order} and on $I_\sharp,\vartheta_\sharp,M_\pressure,C_0$, such that,
for every $t$ satisfying
\[
J_{\cusp}(t)\le\min\{\mathfrak J_{\mathrm{tail}},\mathfrak J_{\mathrm{axis}}\},
\]
the geometric pressure remainder \eqref{eq:Pi-geom-sec14} satisfies
\begin{align}
|\Pi_{\rm geom}(t)| &\le C\bigl(\varepsilon_{\smooth}+\mathfrak a_{\zeta}(I_\sharp)^{\frac12}
+\mathfrak a_{\rm ang}(M_\pressure)^{\frac12}\bigr) \Gamma^2J_{\cusp}^{6\alpha-2}
+C\Gamma^2\Bigl(J_{\cusp}^{6\alpha-2+\kappa_{\rm def}}+J_{\cusp}^{3\alpha-1}\Bigr).
\label{eq:smooth-pressure-defect-bound}
\end{align}
\end{lemma}

\begin{proof}[Proof of Lemma~\ref{lem:smooth-pressure-defect}]
We use the notation in \eqref{eq:smooth-pressure-time-notation-sec14}.  The proof has four steps.  We first
derive the smooth-flow deformation bounds on the fixed ball $B_{R_{\rm pr}}$ from
\eqref{eq:smooth-pressure-fixed-radius-sec14}.  We then compare the bilinear pressure Hessian expressions
$\Pi[w_\sharp,w_\sharp]$ and $\Pi[\Lambda_*w_\sharp,\Lambda_*w_\sharp]$.  The last two steps
estimate the terms involving $w_{\rm out}$, first for the localized regions omitted by
\eqref{eq:Omega-sharp-cutoff-def} and then for the algebraic tail and far labels controlled by
Lemma~\ref{lem:tail-bound}.

\runinhead{Step 1: Smooth-flow bounds on $B_{R_{\rm pr}}$.} Recall from
\eqref{eq:smooth-pressure-time-notation-sec14} that $\Lambda=\phi_{\smooth}(\cdot,t)$.  By the choice of
$R_{\rm pr}$ in \eqref{eq:smooth-pressure-fixed-radius-sec14} and
Lemma~\ref{lem:smooth-flow-small-deformation}, on $B_{R_{\rm pr}}$ we have
\begin{subequations}
\begin{align}
|\Lambda(X)-X|&\le C\varepsilon_{\smooth}|X|,
\label{eq:Lambda-small-position-pressure}\\
|D\Lambda(X)-I|+|D\Lambda(X)^{-1}-I|&\le C\varepsilon_{\smooth},
\\
|D^2\Lambda(X)|&\le C\varepsilon_{\smooth}.
\label{eq:Lambda-small-second-pressure}
\end{align}
\end{subequations}

\runinhead{Step 2: Comparison of the two localized pressure Hessian forms.} We estimate the difference between
$\Pi[\Lambda_*w_\sharp,\Lambda_*w_\sharp]$ and $\Pi[w_\sharp,w_\sharp]$.  For
$v=\Lambda_*w_\sharp$, we use the definition \eqref{eq:physical-pressure-bilinear-form} of
$\Pi[v,v]$ and change variables $y=\Lambda(X)$ in that principal-value integral.  Since $\Lambda$ is
volume preserving,
\begin{equation}
\nabla(\Lambda_*w_\sharp)(\Lambda(X)) = D\Lambda(X)\nabla w_\sharp(X)D\Lambda(X)^{-1}
+D^2\Lambda(X)\bigl[D\Lambda(X)^{-1}(\cdot),w_\sharp(X)\bigr].
\label{eq:smooth-pressure-pushforward-gradient}
\end{equation}
Equation \eqref{eq:smooth-pressure-pushforward-gradient} requires bounds for $\Lambda,D\Lambda,D^2\Lambda$ at
the pressure evaluation point $X$.  Therefore
\eqref{eq:Lambda-small-position-pressure}--\eqref{eq:Lambda-small-second-pressure} are needed on
$B_{R_{\rm pr}}$, and in particular on $B_{R_{\rm eval}}$, not only on the localized vorticity support.  Using
the localized gradient estimate from Lemma~\ref{lem:transported-cusp-pressure-win}, the modulation bound 
\eqref{eq:modulation-bounds}, and $w_\sharp(0,t)=0$, we obtain that
\begin{equation}
\|\nabla w_\sharp\|_{L^\infty}\le C\Gamma J^{3\alpha-1}, \qquad |w_\sharp(X,t)|\le C\Gamma J^{3\alpha-1}|X| \quad (X\in B_{R_{\rm eval}}).
\label{eq:wsharp-gradient-pressure-bound}
\end{equation}
In the variables $X=J^2\bar X$, the localized vorticity definition
\eqref{eq:Omega-sharp-cutoff-def} and the fixed-set estimate \eqref{eq:fixed-scale-CZ-localized} give the
corresponding $C^{\beta_{\rm ax}}$ bounds for the localized velocity gradient, with size
$C\Gamma J^{3\alpha-1}$.

Let
\begin{equation}
A_\Lambda(X):= D\Lambda(X)\nabla w_\sharp(X)D\Lambda(X)^{-1} +D^2\Lambda(X)\bigl[D\Lambda(X)^{-1}(\cdot),w_\sharp(X)\bigr],
\qquad A_0(X):=\nabla w_\sharp(X).
\label{eq:smooth-pressure-gradient-matrices}
\end{equation}
By \eqref{eq:smooth-pressure-pushforward-gradient}, the integrand in
\eqref{eq:physical-pressure-bilinear-form} for
$\Pi[\Lambda_*w_\sharp,\Lambda_*w_\sharp]$ becomes
$K_{zz}(\Lambda(X))\tr(A_\Lambda(X)^2)$ after the change of variables $y=\Lambda(X)$.  The corresponding
integrand for $\Pi[w_\sharp,w_\sharp]$ is $K_{zz}(X)\tr(A_0(X)^2)$.  The kernel $K_{zz}$ is defined in
\eqref{eq:Kzz-kernel}; since $K_{zz}$ is homogeneous of degree $-3$, the position estimate
\eqref{eq:Lambda-small-position-pressure} gives, for $j=0,1$ and $0<|X|<R_{\rm eval}$,
\begin{equation}
\bigl|\nabla_X^{\,j}\bigl(K_{zz}(\Lambda(X))-K_{zz}(X)\bigr)\bigr| \le C_j\varepsilon_{\smooth}|X|^{-3-j}.
\label{eq:deformed-pressure-kernel-Kzz}
\end{equation}
The gradient identity \eqref{eq:smooth-pressure-gradient-matrices}, the smooth-flow bounds
\eqref{eq:Lambda-small-position-pressure}--\eqref{eq:Lambda-small-second-pressure}, and
\eqref{eq:wsharp-gradient-pressure-bound} show that replacing $A_0$ by $A_\Lambda$ changes the source
$\tr(A_0^2)$ in \eqref{eq:physical-pressure-bilinear-form} by
$O(\varepsilon_{\smooth}\Gamma^2J^{6\alpha-2})$ in the same scaled
$C^{\beta_{\rm ax}}$ norm used in \eqref{eq:fixed-scale-pressure-bilinear}.  The principal-value truncation is
taken in the variable $X$ for both integrals.  The map $\Lambda$ moves a sphere $|X|=\rho$ by a relative
$O(\varepsilon_{\smooth})$ amount by \eqref{eq:Lambda-small-position-pressure}, and the mean-zero cancellation
of $K_{zz}$ from \eqref{eq:Kzz-kernel} is the same cancellation used in
\eqref{eq:fixed-scale-pressure-bilinear}.  Applying \eqref{eq:fixed-scale-pressure-bilinear} with the
kernel perturbation \eqref{eq:deformed-pressure-kernel-Kzz} gives
\begin{equation}
\left| \Pi[\Lambda_*w_\sharp,\Lambda_*w_\sharp]-\Pi[w_\sharp,w_\sharp] \right|
\le C\varepsilon_{\smooth}\Gamma^2J^{6\alpha-2}.
\label{eq:localized-smooth-pressure-deformation}
\end{equation}
The contribution of the region $|X|\ge R_{\rm eval}$, where the source and observation regions are disjoint, is
included in $\mathcal R_{\rm out}$; this region is disjoint from the localized support by
\eqref{eq:smooth-pressure-localized-ball-sec14} and the choice $R_{\rm eval}\ge4R_{\rm tube}$.

\runinhead{Step 3: Terms not retained by the localization.} We write
\[
\Pi_{\rm geom} = \Bigl(\Pi[\Lambda_*w_\sharp,\Lambda_*w_\sharp]-\Pi[w_\sharp,w_\sharp]\Bigr) +\mathcal R_{\rm out},
\]
where $\mathcal R_{\rm out}$ is the sum of all terms in which at least one velocity is $w_{\rm out}$ or
$\Lambda_*w_{\rm out}$, together with the disjoint-region part of the localized comparison just described.  We
decompose $w_{\rm out}$ according to where the localization in \eqref{eq:Omega-sharp-cutoff-def} is lost:
\[
w_{\rm out}=w_{\rm out}^{\rm loc}+w_{\rm out}^{\rm far}.
\]
The velocity $w_{\rm out}^{\rm loc}$ contains the contribution from the region where the axial cutoff
$1-\vartheta_\sharp$ is active, the contribution from the region where the angular cutoff
$1-\chi_{M_\pressure}$ is active, and the contribution from bounded-core labels whose cusp-flow images are
outside the tube \eqref{eq:smooth-pressure-localized-tube-sec14}.  By the definition of $R_{\rm loc}$ before
Lemma~\ref{lem:smooth-pressure-defect}, the images of these three regions lie in
$B_{R_{\rm loc}}\subset B_{R_{\rm pr}}$, with $R_{\rm loc}$ independent of the small cusp clock.  Hence the
near-identity estimates \eqref{eq:Lambda-small-position-pressure}--\eqref{eq:Lambda-small-second-pressure}
hold for the terms in $w_{\rm out}^{\rm loc}$.  For the parts of $\mathcal R_{\rm out}$ containing
$w_{\rm out}^{\rm loc}$ or $\Lambda_*w_{\rm out}^{\rm loc}$, the fixed-set estimates
\eqref{eq:fixed-scale-CZ-localized} and \eqref{eq:fixed-scale-pressure-bilinear} control the velocity-gradient
and pressure Hessian integrals on the images contained in $B_{R_{\rm loc}}$.  The cone bilinear estimate
\eqref{eq:pressure-bilinear-cone-estimate} controls the interactions with the portions removed by the
localization in \eqref{eq:Omega-sharp-cutoff-def}.  These estimates yield the axial tail
\eqref{eq:pressure-zeta-tail}, the angular tail \eqref{eq:pressure-angular-tail}, and the normal-form
displacement error \eqref{eq:normal-form-approximation-bound}.

\runinhead{Step 4: Algebraic tail and far labels.} The velocity $w_{\rm out}^{\rm far}$ contains the algebraic tail
and the labels whose cusp-flow images stay outside $B_{R_{\rm loc}}$.  On the support of the smooth velocity,
the cutoff in
\eqref{eq:smooth-velocity-def} enforces $|y|\ge R_{\tail}$.  Lemma~\ref{lem:JtwoD-tail-bdd} supplies the smooth
velocity bounds, and Lemma~\ref{lem:tail-bound} supplies the algebraic-tail bounds for the cusp error.  The
pressure Hessian terms containing $w_{\rm out}^{\rm far}$ are therefore estimated by kernel bounds for disjoint
source and observation regions, with the same algebraic-tail and far-label powers as in
Lemma~\ref{lem:tail-bound}.  Combining the local and far terms yields
\begin{equation}
|\mathcal R_{\rm out}| \le C\Gamma^2J^{6\alpha-2} \bigl(\mathfrak a_{\zeta}(I_\sharp)^{\frac12}+\mathfrak a_{\rm ang}(M_\pressure)^{\frac12}
+J^{\kappa_{\rm def}}+J^{1-3\alpha}\bigr).
\label{eq:smooth-pressure-complement-deformation}
\end{equation}
The terms involving $\mathfrak a_{\zeta}(I_\sharp)^{{\frac{1}{2}}}$ and
$\mathfrak a_{\rm ang}(M_\pressure)^{{\frac{1}{2}}}$ in \eqref{eq:smooth-pressure-complement-deformation} are kept
explicitly in \eqref{eq:smooth-pressure-defect-bound}.  The powers $J^{\kappa_{\rm def}}$ and
$J^{1-3\alpha}$ are lower order in the cusp clock.  Since $0<J\le1$,
\[
C\Gamma^2J^{6\alpha-2}\bigl(J^{\kappa_{\rm def}}+J^{1-3\alpha}\bigr) \le C\Gamma^2\bigl(J^{6\alpha-2+\kappa_{\rm def}}+J^{3\alpha-1}\bigr).
\]
Combining this estimate with \eqref{eq:localized-smooth-pressure-deformation} proves
\eqref{eq:smooth-pressure-defect-bound}.
\end{proof}

\begin{lemma}[Lower-order pressure Hessian remainder]
\label{lem:tail-error-final}
Assuming the standing assumption \eqref{eq:gamma-stand}, define
\begin{equation}
\Pi_{\rm rem}(t):=\Pi_{\rm geom}(t)+\Pi_{\mixed}(t)+\Pi_{\smooth}(t)+\Pi_{\err}(t), \qquad \Pi_0(t)=m(t)^2\,\Pi_{\cusp}(t)
+ \Pi_{\rm rem}(t).
\label{eq:pressure-split-modulated}
\end{equation}
Then there exists $C_\Pi<\infty$, depending only on $\alpha,\gamma,\sigma_{\inn},\sigma_*$, such that
\begin{align}
|\Pi_{\rm rem}(t)|
&\le C_\Pi\bigl(\varepsilon_{\smooth}
+\mathfrak a_{\zeta}(I_\sharp)^{\frac12}
+\mathfrak a_{\rm ang}(M_\pressure)^{\frac12}\bigr)
\Gamma^2J_{\cusp}(t)^{6\alpha-2} \notag \\
&\quad
+ C_\Pi\,\Gamma^2\Bigl(
J_{\cusp}(t)^{6\alpha-2+\kappa_{\rm def}}
+J_{\cusp}(t)^{9\alpha-2}
+J_{\cusp}(t)^{3\alpha-1}
+J_{\cusp}(t)^{2\alpha}
\Bigr)
\label{eq:pressure-remainder-bound}
\end{align}
for $J_{\cusp}(t)\le
\min\{\mathfrak J_{\mathrm{tail}},\mathfrak J_{\mathrm{axis}}\}$.
\end{lemma}

\begin{proof}[Proof of Lemma~\ref{lem:tail-error-final}]
We fix $t$ satisfying
$J_{\cusp}(t)\le \min\{\mathfrak J_{\mathrm{tail}},\mathfrak J_{\mathrm{axis}}\}$.  Since
$J_{\cusp}(t)\le \mathfrak J_{\mathrm{axis}}$, Proposition~\ref{prop:small-clock-comparisons} supplies the
closed axial flow map and normal-form estimates used in Lemma~\ref{lem:tail-bound}.  Since also
$J_{\cusp}(t)\le \mathfrak J_{\mathrm{tail}}$, Lemma~\ref{lem:tail-bound} yields
\[
|\Pi_{\mixed}(t)|+|\Pi_{\smooth}(t)|+|\Pi_{\err}(t)| \le C\,\Gamma^2\Big( J_{\cusp}(t)^{9\alpha-2} {}+ J_{\cusp}(t)^{3\alpha-1} {}
+J_{\cusp}(t)^{2\alpha} \Big).
\]
The same clock restriction allows us to apply Lemma~\ref{lem:smooth-pressure-defect}, which yields
\begin{align*}
|\Pi_{\rm geom}(t)|
&\le C\bigl(\varepsilon_{\smooth} +\mathfrak a_{\zeta}(I_\sharp)^{\frac12}
+\mathfrak a_{\rm ang}(M_\pressure)^{\frac12}\bigr)\Gamma^2J_{\cusp}(t)^{6\alpha-2}
\notag\\
&\quad
+ C\,\Gamma^2\Big( J_{\cusp}(t)^{6\alpha-2+\kappa_{\rm def}} +J_{\cusp}(t)^{3\alpha-1}
\Big).
\end{align*}
Adding this estimate to the preceding bound for $\Pi_{\mixed}+\Pi_{\smooth}+\Pi_{\err}$ and using
\[
\Pi_{\rm rem}=\Pi_{\rm geom}+\Pi_{\mixed}+\Pi_{\smooth}+\Pi_{\err}
\]
proves \eqref{eq:pressure-remainder-bound}, after increasing the constant to $C_\Pi$.
\end{proof}

\begin{proposition}[True pressure Hessian bound relative to $\rW_{\cusp}^2$]
\label{prop:pressure-compare-modulated}
There are constants $\mathfrak J_{\Pi}\in(0,1]$ and $0<q_{\rm phys}<\upbeta$, depending only on
$\alpha,\gamma,\sigma_{\inn},\sigma_*$, such that whenever $J_{\cusp}(t)\le \mathfrak J_{\Pi}$,
\begin{equation}
\Pi_0(t) \ge -q_{\rm phys}\,\tfrac12\,\rW_{\cusp}(t)^2 .
\label{eq:pressure-riccati-full}
\end{equation}
\end{proposition}

\begin{proof}[Proof of Proposition~\ref{prop:pressure-compare-modulated}]
\runinhead{Step 1: The cusp-flow pressure Hessian term.} We first choose
\[
\mathfrak J_{\Pi}\le \min\{\mathfrak J_{\mathrm{tail}},\mathfrak J_{\mathrm{axis}},
\mathfrak J_{\mathrm{strain}},\mathfrak J_{\mathrm{mod}}\}.
\]
We fix $t$ with $J:=J_{\cusp}(t)\le \mathfrak J_{\Pi}$.  Since
\[
\mathfrak J_{\Pi}\le \mathfrak J_{\mathrm{axis}},
\]
the transferred cusp-flow pressure Hessian estimate \eqref{eq:small-clock-cusp-riccati} from
Proposition~\ref{prop:small-clock-comparisons} gives
\[
\Pi_{\cusp}(t) \ge -q_{\rm tr}\,\tfrac12\,\mathcal W_{\cusp}(t)^2, \qquad q_{\rm tr}<\upbeta.
\]
The scalar-modulation definition \eqref{eq:modulation-def} gives
$\rW_{\cusp}(t)=m(t)\mathcal W_{\cusp}(t)$, and therefore
\begin{equation}
m(t)^2\Pi_{\cusp}(t) \ge -q_{\rm tr}\,\tfrac12\,\rW_{\cusp}(t)^2 .
\label{eq:modulated-cusp-pressure-riccati}
\end{equation}

\runinhead{Step 2: The lower-order pressure remainder relative to $\rW_{\cusp}^2$.} Since
$J\le\min\{\mathfrak J_{\mathrm{tail}},\mathfrak J_{\mathrm{axis}}\}$, Lemma~\ref{lem:tail-error-final}
and \eqref{eq:pressure-remainder-bound} imply
\begin{equation}
|\Pi_{\rm rem}(t)| \le C_\Pi\bigl(\varepsilon_{\smooth} +\mathfrak a_{\zeta}(I_\sharp)^{\frac12}
+\mathfrak a_{\rm ang}(M_\pressure)^{\frac12}\bigr) \Gamma^2J^{6\alpha-2}
+ C_\Pi\Gamma^2 \Bigl(J^{6\alpha-2+\kappa_{\rm def}}+J^{9\alpha-2}+J^{3\alpha-1}+J^{2\alpha}\Bigr).
\label{eq:pressure-rem-bound-for-absorption}
\end{equation}
Because $J\le\mathfrak J_{\mathrm{strain}}$, the cusp-coordinate axial strain estimate
\eqref{eq:Wcusp-scaling} implies $|\mathcal W_{\cusp}(t)|\ge c_W\Gamma J^{3\alpha-1}$.  Because
$J\le\mathfrak J_{\mathrm{mod}}$, the modulation bound \eqref{eq:modulation-bounds} implies
$m(t)\ge c_m$.  Hence, using \eqref{eq:modulation-def},
\begin{equation}
\rW_{\cusp}(t)^2 = m(t)^2\mathcal W_{\cusp}(t)^2 \ge c_m^2c_W^2\,\Gamma^2J^{6\alpha-2}.
\label{eq:rWcusp-square-lower-pressure-proof}
\end{equation}
Dividing \eqref{eq:pressure-rem-bound-for-absorption} by
\eqref{eq:rWcusp-square-lower-pressure-proof}, we obtain
\begin{equation}
\tfrac{|\Pi_{\rm rem}(t)|}{\rW_{\cusp}(t)^2} \le C\Bigl( \varepsilon_{\smooth} +\mathfrak a_{\zeta}(I_\sharp)^{\frac12}
+\mathfrak a_{\rm ang}(M_\pressure)^{\frac12} + J^{\kappa_{\rm def}} {}+ J^{3\alpha} {}+ J^{1-3\alpha} {}+ J^{2-4\alpha} \Bigr),
\label{eq:pressure-rem-rWcusp-ratio}
\end{equation}
where the clock powers in \eqref{eq:pressure-rem-rWcusp-ratio} come from
\[
\tfrac{J^{6\alpha-2+\kappa_{\rm def}}}{J^{6\alpha-2}}=J^{\kappa_{\rm def}}, \qquad \tfrac{J^{9\alpha-2}}{J^{6\alpha-2}}=J^{3\alpha},
\qquad \tfrac{J^{3\alpha-1}}{J^{6\alpha-2}}=J^{1-3\alpha}, \qquad \tfrac{J^{2\alpha}}{J^{6\alpha-2}}=J^{2-4\alpha}.
\]
The exponent $\kappa_{\rm def}$ is positive by \eqref{eq:pressure-normal-form-deformation-error}, and
$3\alpha$, $1-3\alpha$, and $2-4\alpha$ are positive under the subcritical restriction
\eqref{eq:pressure-subcritical-riccati}.

\runinhead{Step 3: Smallness of $\varepsilon_{\smooth}$, $\mathfrak a_{\zeta}(I_\sharp)$,
$\mathfrak a_{\rm ang}(M_\pressure)$, and the clock powers.} Since
$q_{\rm tr}<\upbeta$ in \eqref{eq:small-clock-cusp-riccati}, choose
$\delta_{\rm rem}>0$ so small that
\[
q_{\rm phys}:=q_{\rm tr}+\delta_{\rm rem}<\upbeta .
\]
The choice order in Subsection~\ref{sec:fixed-choice-order} permits the following smallness requirement before
$\mathfrak J_{\Pi}$ is fixed:
\begin{equation}
C\bigl( \varepsilon_{\smooth} +\mathfrak a_{\zeta}(I_\sharp)^{\frac12} +\mathfrak a_{\rm ang}(M_\pressure)^{\frac12} \bigr)
\le \tfrac{\delta_{\rm rem}}{4}.
\label{eq:fixed-tail-smooth-rem-smallness}
\end{equation}
Indeed, the interval and cutoff $I_\sharp,\vartheta_\sharp$ are chosen in
\eqref{eq:pressure-localization-cutoff}--\eqref{eq:pressure-localization-intervals} so that the
$\zeta$-tail in \eqref{eq:pressure-zeta-tail} is small; the angular cutoff $M_\pressure$ is chosen
in \eqref{eq:pressure-C0} so that the angular tail \eqref{eq:pressure-angular-tail} is small; and
$R_{\tail}$ is chosen in \eqref{eq:core-tail-domains}--\eqref{eq:smooth-velocity-def} so that the
smooth-flow estimates \eqref{eq:usmooth-small-C2a}--\eqref{eq:smooth-flow-second-gradient-small} give the
prescribed smallness of $\varepsilon_{\smooth}$.  We then decrease
$\mathfrak J_{\Pi}$ so that
\begin{equation}
C\Bigl(\mathfrak J_{\Pi}^{\kappa_{\rm def}}+\mathfrak J_{\Pi}^{3\alpha} +\mathfrak J_{\Pi}^{1-3\alpha}+\mathfrak J_{\Pi}^{2-4\alpha}\Bigr)
\le \tfrac{\delta_{\rm rem}}{4}.
\label{eq:pressure-rem-clock-smallness}
\end{equation}
Combining \eqref{eq:fixed-tail-smooth-rem-smallness}, \eqref{eq:pressure-rem-clock-smallness}, and
\eqref{eq:pressure-rem-rWcusp-ratio}, we obtain
\begin{equation}
|\Pi_{\rm rem}(t)| \le \delta_{\rm rem}\,\tfrac12\,\rW_{\cusp}(t)^2.
\label{eq:full-pressure-rem-absorbed}
\end{equation}

\runinhead{Step 4: Combining the pressure decomposition.} Combining \eqref{eq:pressure-split-modulated},
\eqref{eq:modulated-cusp-pressure-riccati}, and
\eqref{eq:full-pressure-rem-absorbed} yields
\[
\Pi_0(t) \ge -(q_{\rm tr}+\delta_{\rm rem})\,\tfrac12\,\rW_{\cusp}(t)^2 = -q_{\rm phys}\,\tfrac12\,\rW_{\cusp}(t)^2.
\]
This proves \eqref{eq:pressure-riccati-full}.
\end{proof}

\subsection{Finite-time blowup: proof of Theorem \ref{thm:target-profile}}

We next convert the pressure Hessian bound \eqref{eq:pressure-riccati-full}, still measured against
$\rW_{\cusp}(t)^2$, into the true Euler Riccati inequality for the exact Target Profile solution.

\begin{lemma}[No small-clock breakdown before cusp collapse]
\label{lem:no-breakdown-before-clock-collapse}
Let $0<j_0\le \mathfrak J_{\mathrm{collapse}}$, and let $I$ be a time interval on which the Euler
velocity is defined and remains uniformly bounded in $C^{1,\alpha}$.  Assume that, for every $t\in I$,
\[
j_0\le J_{\cusp}(t)\le \mathfrak J_{\mathrm{collapse}}.
\]
Then the vorticity remains bounded in $L^\infty$ by a constant depending on $j_0$, $\Gamma$, and the fixed data.
Consequently, if $T'<\infty$ is a right endpoint of $I$, the Beale--Kato--Majda criterion rules out $T'$ as a
maximal existence time for the Euler solution while $J_{\cusp}(t)\ge j_0$ on $I$.
\end{lemma}

\begin{proof}[Proof of Lemma~\ref{lem:no-breakdown-before-clock-collapse}]
The uniform $C^{1,\alpha}$ bound provides the classical Euler flow map, and the transport identity
\eqref{eq:vort-identity-JtwoD} on $I$.  We use the global upper vorticity estimate in
Lemma~\ref{lem:global-upper-vorticity-envelope}.  Its proof uses the bounded-core normal form
\eqref{eq:bounded-core-normal-form}--\eqref{eq:bounded-core-normal-form-error}, the transported-vorticity
representation \eqref{eq:Omega-cusp-local-representation}--\eqref{eq:Omega-reg-local-bounds}, the bound for
$J_{\twoD}^{-1}$ in \eqref{eq:JtwoD-core-upper}, and the radial-logarithmic bounds
\eqref{eq:smooth-radial-log-growth}--\eqref{eq:cusp-separated-radial-log-growth}.  After decreasing
$\mathfrak J_{\mathrm{collapse}}$ below the threshold
$\mathfrak J_{\omega,+}$ in Lemma~\ref{lem:global-upper-vorticity-envelope}, the hypothesis on $I$ yields
\[
j_0\le J_{\cusp}(t)\le \mathfrak J_{\omega,+}\ \ \text{ for } \ \ t\in I.
\]
For each $t\in I$ and each label $Y$ with $\omega_{\theta,0}(Y)\neq0$, \eqref{eq:global-upper-vorticity-envelope} and \eqref{eq:vort-identity-JtwoD} imply
\[
|\omega_\theta(\phi(Y,t),t)| \le C_{\omega,+}\Gamma J_{\cusp}(t)^{3\alpha-1} \le C_{\omega,+}\Gamma j_0^{3\alpha-1}.
\]
The same bound is trivial on labels with $\omega_{\theta,0}(Y)=0$.  Taking the supremum over labels yields
\[
\|\omega(\cdot,t)\|_{L^\infty} \le C_{\omega,+}\Gamma j_0^{3\alpha-1} \ \ \text{ for } \ \ t\in I,
\]
because $3\alpha-1<0$ and $J_{\cusp}(t)\ge j_0$.  Hence, for any $a\in I$ and any finite right endpoint $T'$ of $I$,
\[
\int_a^{T'}\|\omega(\cdot,t)\|_{L^\infty}\,\ud t<\infty .
\]
The Beale--Kato--Majda continuation criterion for the classical Euler solutions considered here therefore excludes $T'$ as a maximal 
existence time while $J_{\cusp}(t)\ge j_0$ on $I$.
\end{proof}

\begin{proposition}[Finite-time collapse for the exact Target Profile datum]
\label{prop:blowup-final}
Let $\alpha\in(0,\tfrac13)$ satisfy \eqref{eq:pressure-subcritical-riccati}, and let
$\gamma>\alpha+\tfrac52$.  Let $(u,\phi)$ be the Euler solution generated by
the exact Target Profile datum $\Theta^*$.  Then there exists
$\mathfrak J_{\mathrm{collapse}}\in(0,1]$ such that whenever
$J_{\cusp}(t)\le \mathfrak J_{\mathrm{collapse}}$,
\begin{equation}
\Pi_0(t)\ge -\upbeta\,\tfrac12\,\rW_0(t)^2,
\label{eq:euler-riccati-bound}
\end{equation}
and hence
\begin{equation}
\p_t \rW_0(t) =-\tfrac12\,\rW_0(t)^2-\Pi_0(t) \le -\tfrac{1-\upbeta}{2}\,\rW_0(t)^2.
\label{eq:euler-riccati-ineq}
\end{equation}
Consequently there exists a finite time $T^*<\infty$ such that
\[
J_{\cusp}(t)\to 0, \qquad J(t)\to 0, \qquad \rW_{\cusp}(t)\to -\infty, \qquad \rW_0(t)\to -\infty
\]
as $t\uparrow T^*$.
\end{proposition}

\begin{proof}[Proof of Proposition~\ref{prop:blowup-final}]
We choose
\[
\mathfrak J_{\mathrm{collapse}} \le \min\{\mathfrak J_{\Pi},\mathfrak J_{\mathrm{axis}},
\mathfrak J_{\mathrm{strain}},\mathfrak J_{\mathrm{mod}}\},
\]
and later replace it only by smaller thresholds when a proof step requires this.  Thus the estimates
\eqref{eq:pressure-riccati-full}, \eqref{eq:Wcusp-scaling}, \eqref{eq:modulation-bounds}, and
\eqref{eq:Jdot-two-sided} are available whenever
$J_{\cusp}\le \mathfrak J_{\mathrm{collapse}}$.

\runinhead{Step 1: Replacing $\rW_{\cusp}(t)^2$ by $\rW_0(t)^2$ in the pressure bound.}
We fix $t$ with $J_{\cusp}(t)\le\mathfrak J_{\mathrm{collapse}}$.  The cusp-coordinate axial strain estimate \eqref{eq:Wcusp-scaling} states that
$\mathcal W_{\cusp}(t)<0$ and
\[
c_W\Gamma J_{\cusp}(t)^{3\alpha-1} \le |\mathcal W_{\cusp}(t)| \le C_W\Gamma J_{\cusp}(t)^{3\alpha-1}.
\]
Combining this with the modulation identity \eqref{eq:modulation-def} and the modulation bounds \eqref{eq:modulation-bounds}, we obtain
\begin{equation}
c_m c_W\Gamma J_{\cusp}(t)^{3\alpha-1} \le -\rW_{\cusp}(t) \le C_m C_W\Gamma J_{\cusp}(t)^{3\alpha-1}.
\label{eq:collapse-proof-rWcusp-scale}
\end{equation}
We denote by $q_{\rm phys}<\upbeta$ the constant from Proposition~\ref{prop:pressure-compare-modulated}.  We choose
$\delta_0>0$ so small that
\[
q_{\rm phys}(1+\delta_0)\le \upbeta .
\]
We choose $\varepsilon_0>0$ so that $(1-\varepsilon_0)^{-2}\le1+\delta_0$.
Since $3\alpha-1<0$, we shrink $\mathfrak J_{\mathrm{collapse}}$ so that
\[
C_{\smooth} \le \varepsilon_0\,c_m c_W\,\mathfrak J_{\mathrm{collapse}}^{3\alpha-1}.
\]
For every $J_{\cusp}(t)\le\mathfrak J_{\mathrm{collapse}}$, \eqref{eq:Jsmooth-bdd} and \eqref{eq:collapse-proof-rWcusp-scale} imply that
\[
|\rW_{\smooth}(t)|\le\varepsilon_0(-\rW_{\cusp}(t)).
\]
The axial strain decomposition \eqref{eq:on-axis-strain-splitting} is
$\rW_0(t)=\rW_{\smooth}(t)+\rW_{\cusp}(t)$.  Together with the preceding estimate, it implies
\begin{equation}
-\rW_0(t) \ge (1-\varepsilon_0)(-\rW_{\cusp}(t))>0
\label{eq:collapse-proof-total-strain-negative}
\end{equation}
and therefore
\begin{equation}
\rW_{\cusp}(t)^2\le(1+\delta_0)\rW_0(t)^2.
\label{eq:collapse-proof-cusp-total-square}
\end{equation}
Since $\mathfrak J_{\mathrm{collapse}}\le\mathfrak J_{\Pi}$, the pressure comparison
\eqref{eq:pressure-riccati-full} implies
\[
\Pi_0(t) \ge -q_{\rm phys}\,\tfrac12\,\rW_{\cusp}(t)^2 \ge -\upbeta\,\tfrac12\,\rW_0(t)^2,
\]
where the last inequality uses \eqref{eq:collapse-proof-cusp-total-square} and the choice of $\delta_0$.
This proves \eqref{eq:euler-riccati-bound}.

\runinhead{Step 2: Entering the small-cusp-clock regime.}
Let
\[
t_0:=\inf\{t\ge 0:\ J_{\cusp}(t)\le \mathfrak J_{\mathrm{collapse}}\}.
\]
Lemma~\ref{lem:late-entry}, applied with $\mathfrak J_{\mathrm{finite}}=\mathfrak J_{\mathrm{collapse}}$, yields
$t_0<\infty$.  At $t=t_0$, \eqref{eq:collapse-proof-total-strain-negative} implies $\rW_0(t_0)<0$.
Moreover, whenever $J_{\cusp}(t)\le\mathfrak J_{\mathrm{collapse}}$, the clock law
\eqref{eq:Jdot-two-sided} yields
\[
\dot J_{\cusp}(t)\le -c_1\Gamma J_{\cusp}(t)^{3\alpha}<0.
\]
By continuity, $J_{\cusp}(t)$ cannot cross upward through
$\mathfrak J_{\mathrm{collapse}}$ after $t_0$.  Hence
\begin{equation}
J_{\cusp}(t)\le\mathfrak J_{\mathrm{collapse}}
\label{eq:collapse-proof-clock-stays-small}
\end{equation}
for all subsequent times for which the $C^{1,\alpha}$ solution exists and $J_{\cusp}(t)>0$.

\runinhead{Step 3: Integrating the Riccati inequality up to the first collapse endpoint.}
Let $T_{\max}$ denote the maximal $C^{1,\alpha}$ existence time after entry.  We argue on
$[t_0,T_{\max})$, restricted to times for which $J_{\cusp}(t)>0$.  By
\eqref{eq:collapse-proof-clock-stays-small}, Step~1 applies throughout that subinterval.  Hence
\eqref{eq:euler-riccati-bound} and the stagnation-point identity \eqref{eq:stagnation-ODEs} imply
\eqref{eq:euler-riccati-ineq}.  Since $\rW_0(t_0)<0$, \eqref{eq:euler-riccati-ineq} keeps
$\rW_0(t)<0$ while the solution exists and $J_{\cusp}(t)>0$.  Thus
\[
\p_t\Bigl(-\tfrac{1}{\rW_0(t)}\Bigr) = \tfrac{\p_t\rW_0(t)}{\rW_0(t)^2} \le -\tfrac{1-\upbeta}{2}.
\]
Consequently
\[
0<-\tfrac{1}{\rW_0(t)} \le -\tfrac{1}{\rW_0(t_0)} -\tfrac{1-\upbeta}{2}(t-t_0)
\]
as long as $J_{\cusp}(t)>0$ and the $C^{1,\alpha}$ solution exists with finite $\rW_0(t)$.  Therefore this
continuation is impossible beyond
\[
t_0+\tfrac{2}{(1-\upbeta)|\rW_0(t_0)|}
\]
unless $J_{\cusp}$ has already collapsed.  If a finite maximal $C^{1,\alpha}$ existence time occurred before
this endpoint while $J_{\cusp}$ stayed bounded below by some $j_0>0$, then
Lemma~\ref{lem:no-breakdown-before-clock-collapse} would rule out that maximal time.  Hence the first finite
endpoint is either characterized by
\[
J_{\cusp}(t)\to0
\]
or by
\[
-\tfrac{1}{\rW_0(t)}\to0 .
\]
We denote this first endpoint by $T^*$.  The preceding upper bound implies $T^*<\infty$.

\runinhead{Step 4: Limits of $J_{\cusp}$, $J$, $\rW_{\cusp}$, and $\rW_0$ at $T^*$.}
If $-1/\rW_0(t)\to0$ as $t\uparrow T^*$, then $\rW_0(t)\to-\infty$.  The smooth-strain bound \eqref{eq:Jsmooth-bdd} and the axial strain 
decomposition \eqref{eq:on-axis-strain-splitting} then show that
\[
\rW_{\cusp}(t)=\rW_0(t)-\rW_{\smooth}(t)\to-\infty.
\]
The upper bound in \eqref{eq:collapse-proof-rWcusp-scale}, together with $3\alpha-1<0$, then forces $J_{\cusp}(t)\to0$.

If $J_{\cusp}(t)\to0$ first, then the lower bound in \eqref{eq:collapse-proof-rWcusp-scale} yields $\rW_{\cusp}(t)\to-\infty$, 
and \eqref{eq:Jsmooth-bdd} together with \eqref{eq:on-axis-strain-splitting} implies $\rW_0(t)\to-\infty$.  In both cases,
\[
J_{\cusp}(t)\to0, \qquad \rW_{\cusp}(t)\to-\infty, \qquad \rW_0(t)\to-\infty .
\]
Finally, the clock decomposition \eqref{eq:J-clock-decomposition} and the smooth-clock bound \eqref{eq:Jsmooth-bdd} show that
\[
J(t)=J_{\smooth}(t)J_{\cusp}(t)\to 0.
\]
\end{proof}

\subsection{Global Type--I vorticity bounds}

The last step is the global Type--I vorticity bound.  By the exact transport identity
\eqref{eq:vort-identity-JtwoD},
\[
\omega_\theta(\phi(Y,t),t) = J_{\twoD}(Y,t)^{-1}\omega_{\theta,0}(Y),
\]
the Type--I estimate reduces to bounding
$J_{\twoD}(Y,t)^{-1}|\omega_{\theta,0}(Y)|$ uniformly in the label $Y$.  The upper inequality in
Lemma~\ref{lem:target-vorticity-envelope} follows from the pointwise estimate
\eqref{eq:global-upper-vorticity-envelope}; its proof combines the bound for $J_{\twoD}^{-1}$ in
\eqref{eq:JtwoD-core-upper} with the cone-local transported-vorticity
representation \eqref{eq:Omega-cusp-local-representation}--\eqref{eq:Omega-reg-local-bounds}.  For the lower
bound, the labels are chosen in the fixed core $D_{\core}$ from \eqref{eq:core-tail-domains} with
\[
\rho(Y)\in[\tfrac12,1], \qquad \sigma(Y)=\kappa J_{\cusp}(t)^3, \qquad 0<\kappa\le \min\{\sigma_{\cut},c_{\rm vort}\},
\]
where $\sigma$ is the polar angle in \eqref{eq:polar-angle-def}.  Thus $Y$ lies in the initial cone
$0\le\sigma\le\sigma_{\cut}$, on which $\Upsilon=1$ in \eqref{eq:Theta-star-def}; the datum \eqref{eq:vort0}
then implies $|\omega_{\theta,0}(Y)|\gtrsim\Gamma J_{\cusp}(t)^{3\alpha}$.  Also
$\vartheta(Y)=R(Y)/\rho(Y)=\sin\sigma(Y)\le c_{\rm vort}J_{\cusp}(t)^3$.  Therefore,
\eqref{eq:driver-amplification-lower} and \eqref{eq:vort-identity-JtwoD} imply the lower inequality in
Lemma~\ref{lem:target-vorticity-envelope}.

\begin{lemma}[Bounds for $J_{\twoD}^{-1}$ on $D_{\core}$]
\label{lem:JtwoD-core-bounds}
There exist constants
\begin{equation}
\mathfrak J_{\mathrm{vort}}\in(0,\min\{\mathfrak J_{\mathrm{axis}},\mathfrak J_{\mathrm{mod}}\}]\qquad\text{and}\qquad C_{\rm vort}<\infty,
\label{eq:vort-clock-threshold-def}
\end{equation}
depending only on $\alpha,\gamma,\sigma_{\inn},\sigma_*$, such that the following holds whenever
$J_{\cusp}(t)\le \mathfrak J_{\mathrm{vort}}$.  Let $Y\in D_{\core}$ with $\omega_{\theta,0}(Y)\neq0$, and define
\[
\vartheta(Y):=\tfrac{R(Y)}{\rho(Y)}\in(0,1].
\]
Then
\begin{equation}
J_{\twoD}(Y,t)^{-1} \le C_{\rm vort}\, \Bigl(\max\{J_{\cusp}(t),\,\vartheta(Y)^{\tfrac{1}{3}}\}\Bigr)^{-1}.
\label{eq:JtwoD-core-upper}
\end{equation}
Moreover, if $\vartheta(Y)\le c_{\rm vort}J_{\cusp}(t)^3$ and the initial polar angle satisfies
$0\le\sigma(Y)\le\sigma_{\max}$, with $\sigma$ defined in \eqref{eq:polar-angle-def} and
$\sigma_{\max}$ fixed in Definition~\ref{def:init-data}, then
\begin{equation}
J_{\twoD}(Y,t)^{-1}\ge c_{\rm vort}\,J_{\cusp}(t)^{-1}.
\label{eq:driver-amplification-lower}
\end{equation}
\end{lemma}

\begin{proof}[Proof of Lemma~\ref{lem:JtwoD-core-bounds}]
We give the proof in the upper half-space; the lower half-space follows by odd symmetry.  Since
$\omega_{\theta,0}$ is supported in the fixed cone $\sigma\le\sigma_{\max}<\sigma_{\inn}$, on the support of
the datum we have constants $0<c_\sigma<C_\sigma<\infty$ such that
\begin{equation}
c_\sigma\rho(Y)\le Z(Y)\le \rho(Y), \qquad c_\sigma\,\tfrac{R(Y)}{Z(Y)} \le \vartheta(Y) \le C_\sigma\,\tfrac{R(Y)}{Z(Y)}.
\label{eq:vartheta-RZ-comparison}
\end{equation}

\runinhead{Step 1: bounds for $J_{\twoD}^{-1}$ while $\phi_{\cusp}(Y,s)\in\mathcal C_*$.}
Set
\[
x_*(s):=\phi_{\cusp}(Y,s).
\]
Cone entry and exit times in this proof are defined by membership of $x_*(s)$ in $\mathcal C_*$.  Since
$\omega_{\theta,0}(Y)\neq0$, we have $R(Y)>0$.  Comparing
\eqref{eq:vort-identity-JtwoD} with the specific-vorticity transport identity
\eqref{eq:lin-deriv:V-stretch} yields
\begin{equation}
J_{\twoD}(Y,s)^{-1} = \tfrac{\phi_r(Y,s)}{R(Y)} = \tfrac{(\phi_{\smooth})_r(x_*(s),s)}{R(Y)}.
\label{eq:JtwoD-radial-stretch-core}
\end{equation}

If $x_*(s)\in\mathcal C_*$ and $J_{\cusp}(s)\le\mathfrak J_{\mathrm{vort}}$, then
\eqref{eq:exact-current-axis-core-geometry}, applied at time $s$, yields
\begin{subequations}
\label{eq:exact-current-axis-comparability}
\begin{align}
c\,\tfrac{R(Y)}{J_{\cusp}(s)} &\le r(x_*(s))\le C\,\tfrac{R(Y)}{J_{\cusp}(s)},  \label{eq:Rexact-comparable-R}\\
c\,J_{\cusp}(s)^2Z(Y) &\le z(x_*(s))\le C\,J_{\cusp}(s)^2Z(Y).  \label{eq:Zexact-comparable-Z}
\end{align}
\end{subequations}
Let $t_0$ be the entry time of the connected interval, containing $t$, on which $J_{\cusp}\le \mathfrak J_{\mathrm{vort}}$.  The finite-clock estimate
\eqref{eq:finite-clock-smooth-cusp-C1}, applied with $\mathfrak J_{\mathrm{finite}}=\mathfrak J_{\mathrm{vort}}$, provides bounded derivatives for
$\phi_{\smooth}(\cdot,t_0)$ and its inverse on the image of $D_{\core}$ under $\phi_{\cusp}(\cdot,t_0)$.  If $x_*(s)\in\mathcal C_*$, then
\eqref{eq:exact-current-axis-comparability} and $r(x_*(s))\le(\tan\sigma_*)z(x_*(s))$ imply that  $R(Y)/Z(Y)\le C J_{\cusp}(s)^3$ and hence that
\[
|x_*(s)|\le C J_{\cusp}(s)^2 R_{\tail}.
\]
After decreasing $\mathfrak J_{\mathrm{vort}}$, this shows that $x_*(s)\in B_{\frac18 R_{\tail}}$ for every
$s\in[t_0,t]$ with $x_*(s)\in\mathcal C_*$.  Lemma~\ref{lem:JtwoD-tail-bdd}, applied with
$R_0=\frac18R_{\tail}$, gives
$\|\nabla u_{\smooth}(\cdot,s)\|_{L^\infty(B_{\frac14 R_{\tail}})}\le C\Gamma$.  By
\eqref{eq:Jdot-two-sided},
\[
t-t_0 \le C\Gamma^{-1}\int_{J_{\cusp}(t)}^{\mathfrak J_{\mathrm{vort}}}J^{-3\alpha}\,dJ \le C\Gamma^{-1}.
\]
Applying Gronwall's inequality to the differential equation \eqref{eq:smooth-flow-differential-ode} for
$D\phi_{\smooth}$, we obtain that
\[
\|D\phi_{\smooth}(\cdot,s)\|_{L^\infty(B_{\frac18 R_{\tail}})}
+ \|D\phi_{\smooth}(\cdot,s)^{-1}\|_{L^\infty(\phi_{\smooth}(B_{\frac18 R_{\tail}},s))} \le C.
\]
Since $\phi_{\smooth}$ is axisymmetric and preserves the symmetry axis, the preceding derivative and inverse
derivative bounds imply
\begin{equation}
C^{-1}r(x_*(s)) \le (\phi_{\smooth})_r(x_*(s),s) \le C r(x_*(s)).
\label{eq:smooth-flow-radial-distance-core}
\end{equation}
whenever $x_*(s)\in\mathcal C_*$.  Combining the upper bound in
\eqref{eq:smooth-flow-radial-distance-core} with \eqref{eq:Rexact-comparable-R} and
\eqref{eq:JtwoD-radial-stretch-core}, we obtain
\begin{equation}
J_{\twoD}(Y,t)^{-1}\le C J_{\cusp}(t)^{-1}
\label{eq:JtwoD-driver-upper}
\end{equation}
provided $x_*(t)\in\mathcal C_*$.  For the lower bound, we decrease $c_{\rm vort}$ so that the condition
\[
\vartheta(Y)\le c_{\rm vort}J_{\cusp}(t)^3
\]
implies that 
\[
\vartheta(Y)\le c_{\rm vort}J_{\cusp}(s)^3 \ \ \text{ for } t_0\le s\le t,
\]
by the monotonicity in \eqref{eq:Jdot-two-sided}.  By
\eqref{eq:vartheta-RZ-comparison}, after decreasing $c_{\rm vort}$ once more,
\[
\tfrac{R(Y)}{Z(Y)}\le c_*J_{\cusp}(s)^3 \ \ \text{ for } t_0\le s\le t,
\]
where $c_*>0$ is chosen below.  By \eqref{eq:assum1},
\[
\tfrac{A_s(Z(Y))R(Y)}{B_s(Z(Y))} \le Cc_* \ \ \text{ for } t_0\le s\le t.
\]
We choose $c_*$ so small that the right-hand side is less than $C_{\rm cone}$; applying the normal-form
estimates \eqref{eq:core-normal-form-expansion}--\eqref{eq:core-normal-form-relative-error} at time $s$, we
have that
\[
r(x_*(s))\le \tfrac12(\tan\sigma_*)\,z(x_*(s)).
\]
Thus $x_*(s)\in\mathcal C_*$ for $t_0\le s\le t$.  Combining the lower bounds in
\eqref{eq:smooth-flow-radial-distance-core}, \eqref{eq:Rexact-comparable-R}, and
\eqref{eq:JtwoD-radial-stretch-core} proves \eqref{eq:driver-amplification-lower}.

\runinhead{Step 2: labels whose exact cusp trajectory exits $\mathcal C_*$.}
If the exact cusp trajectory has not entered $\mathcal C_*$ at any time in the small-clock interval
$[t_0,t]$, where $t_0$ is the entry time for $J_{\cusp}\le\mathfrak J_{\mathrm{vort}}$, then
$x_*(s)\notin\mathcal C_*$ for all $s\in[t_0,t]$.  The radial logarithmic kernel bound
\eqref{eq:separated-radial-log-kernel}, applied with disjoint source and observation regions, yields the
inequality $|u_r(\phi(Y,s),s)|/\phi_r(Y,s)\le C\Gamma$ for $t_0\le s\le t$.  The interval length estimate
above yields $t-t_0\le C\Gamma^{-1}$, and the finite-clock
bound at $t_0$ implies $J_{\twoD}(Y,t_0)^{-1}\le C$.  Hence, 
\[
J_{\twoD}(Y,t)^{-1}\le C \le C\Bigl(\max\{J_{\cusp}(t),\vartheta(Y)^{1/3}\}\Bigr)^{-1},
\]
because both entries in the maximum are at most one.  It thus remains to treat labels whose exact cusp
trajectory is in $\mathcal C_*$ at least once during the small-clock interval.

Suppose first that $x_*(s)=\phi_{\cusp}(Y,s)\in\mathcal C_*$ at a time $s$ with small cusp clock.  The cone
condition implies
\[
r(x_*(s))\le (\tan\sigma_*)\,z(x_*(s)).
\]
Combining \eqref{eq:Rexact-comparable-R} with \eqref{eq:Zexact-comparable-Z}, we obtain, for
$Y\in D_{\core}\cap\{Z>0\}$ with $\omega_{\theta,0}(Y)\neq0$,
\[
c\,\tfrac{R(Y)}{J_{\cusp}(s)} \le r(x_*(s)) \le C\,\tfrac{R(Y)}{J_{\cusp}(s)},
\qquad c\,J_{\cusp}(s)^2Z(Y) \le z(x_*(s)) \le C\,J_{\cusp}(s)^2Z(Y).
\]
Therefore, by \eqref{eq:vartheta-RZ-comparison}, 
\[
\tfrac{R(Y)}{Z(Y)}\le C J_{\cusp}(s)^3, \qquad\text{and hence}\qquad \vartheta(Y)\le C J_{\cusp}(s)^3.
\]
If the exact cusp trajectory is still in $\mathcal C_*$ at the final time, then the preceding implication with
$s=t$ shows that
\[
\vartheta(Y)\le C J_{\cusp}(t)^3.
\]
Together with \eqref{eq:JtwoD-driver-upper}, this implies that
\[
J_{\twoD}(Y,t)^{-1} \le C J_{\cusp}(t)^{-1} \le C\Bigl(\max\{J_{\cusp}(t),\vartheta(Y)^{\tfrac{1}{3}}\}\Bigr)^{-1}.
\]
Otherwise, let $t_{\rm ex}=t_{\rm ex}(Y,t)$ be the last time before
$t$ at which $x_*(s)$ exits $\mathcal C_*$.  At $s=t_{\rm ex}$ the point $x_*(t_{\rm ex})$ lies on
$\partial\mathcal C_*$.  The preceding inequalities therefore imply the two-sided comparison
\begin{equation}
c J_{\cusp}(t_{\rm ex})^3 \le \vartheta(Y) \le C J_{\cusp}(t_{\rm ex})^3.
\label{eq:last-exit-clock-angle-comparison}
\end{equation}
In particular,
\[
c\,\vartheta(Y)^{{\frac{1}{3}}} \le J_{\cusp}(t_{\rm ex}) \le C\,\vartheta(Y)^{{\frac{1}{3}}},
\]
with constants depending only on the fixed cones.  Then \eqref{eq:JtwoD-driver-upper} implies
\[
J_{\twoD}(Y,t_{\rm ex})^{-1}\le C\,\vartheta(Y)^{-{\frac{1}{3}}}.
\]
By definition of $t_{\rm ex}$,
\[
x_*(s)=\phi_{\cusp}(Y,s)\notin\mathcal C_* \ \ \text{ for } \ \ t_{\rm ex}<s\le t.
\]
Since $\mathcal C_{\inn}\Subset\mathcal C_*$, the cusp-coordinate observation point $x_*(s)$ is separated from
the source cone $\mathcal C_{\inn}$ by a fixed angular distance for $t_{\rm ex}<s\le t$.  Applying the radial
kernel bound \eqref{eq:separated-radial-log-kernel} to the singular cone part and using the separated annular
summation \eqref{eq:separated-leading-radial-annular-sum}, with the exit relation
$J_{\cusp}(t_{\rm ex})^3\simeq\vartheta(Y)$ from \eqref{eq:last-exit-clock-angle-comparison}, gives the
bound $C\Gamma$ for the radial quotient of the singular cusp contribution.  The regular part of the bounded
core is controlled in the same quotient by \eqref{eq:Omega-reg-local-bounds} and the separated annular
kernel estimate \eqref{eq:separated-kernel-annular}.  The smooth far-field contribution satisfies the radial
quotient bound \eqref{eq:smooth-radial-log-growth}.  Therefore, for $t_{\rm ex}\le s\le t$,
\[
\left|\tfrac{u_r(\phi(Y,s),s)}{\phi_r(Y,s)}\right|\le C\Gamma
\]
after the bounded radial multiplier of $\phi_{\smooth}$ is included as in the proof of
Lemma~\ref{lem:radial-log-growth-separated}.  By \eqref{eq:Jdot-two-sided} and $3\alpha<1$,
\[
t-t_{\rm ex} \le C\Gamma^{-1}\int_{J_{\cusp}(t)}^{J_{\cusp}(t_{\rm ex})}J^{-3\alpha}\,dJ
\le C\Gamma^{-1}\int_0^{J_{\cusp}(t_{\rm ex})}J^{-3\alpha}\,dJ \le C\Gamma^{-1}.
\]
The exact logarithmic identity
\[
\p_s\log J_{\twoD}(Y,s)^{-1} = \tfrac{u_r(\phi(Y,s),s)}{\phi_r(Y,s)}
\]
therefore implies
\[
\left| \log\tfrac{J_{\twoD}(Y,t)^{-1}}{J_{\twoD}(Y,t_{\rm ex})^{-1}} \right| \le \int_{t_{\rm ex}}^t C\Gamma\,ds \le C.
\]
Hence, 
\[
J_{\twoD}(Y,t)^{-1}\le C\,\vartheta(Y)^{-\frac{1}{3}} \ \ \text{ for } \ \ t\ge t_{\rm ex}.
\]
Combining the case $x_*(t)\in\mathcal C_*$, covered by \eqref{eq:JtwoD-driver-upper}, with the post-exit
bound above proves \eqref{eq:JtwoD-core-upper}.
\end{proof}

The upper vorticity bound below uses
\eqref{eq:Omega-cusp-local-representation}--\eqref{eq:Omega-reg-local-bounds} only under the cone-locality
condition
\[
B(x,2c_*|x|)\subset\mathcal C_* .
\]
On trajectory intervals where this condition fails, the transport equation for angular vorticity is controlled
by the radial quotient $(V_{\cusp})_r/r$.  We also need the corresponding quotient for $u_{\smooth}$ in order to
compare physical vorticity with cusp-coordinate vorticity.  Thus, with $c_*$ fixed in the proof of
Lemma~\ref{lem:transported-cusp-field-bounds}, we set
\begin{equation}
\mathcal C_*^{\rm int} := \{\,x:\ B(x,2c_*|x|)\subset\mathcal C_*\,\}.
\label{eq:cone-interior-for-radial-log}
\end{equation}

\begin{lemma}[Radial logarithmic growth outside $\mathcal C_*^{\rm int}$]
\label{lem:radial-log-growth-separated}
There exist constants
\begin{equation}
\mathfrak J_{\mathrm{rlog}}\in(0,\min\{\mathfrak J_{\mathrm{vort}},\mathfrak J_{\mathrm{tail}}\}],\qquad C_{\mathrm{rlog}}<\infty,
\label{eq:rlog-clock-threshold-def}
\end{equation}
depending only on $\alpha,\gamma,\sigma_{\inn},\sigma_*$, such that the following bounds hold.
If $J_{\cusp}(t)\le\mathfrak J_{\mathrm{rlog}}$, then
\begin{equation}
\big|\tfrac{(u_{\smooth})_r(x,t)}{r(x)}\big| \le C_{\mathrm{rlog}}\Gamma \qquad (r(x)>0).
\label{eq:smooth-radial-log-growth}
\end{equation}
Let $[s_0,s_1]$ be an interval on which
$J_{\cusp}(s)\le\mathfrak J_{\mathrm{rlog}}$, let $Y$ satisfy $\omega_{\theta,0}(Y)\neq0$, and set
$x_*(s):=\phi_{\cusp}(Y,s)$.  Assume
\[
x_*(s)\notin\mathcal C_*^{\rm int}\qquad (s_0\le s\le s_1),
\]
and assume either that $s_0$ is the left endpoint of the connected time interval on which
$J_{\cusp}\le\mathfrak J_{\mathrm{rlog}}$ or that
$x_*(s_0)\in\partial\mathcal C_*^{\rm int}$.  Then
\begin{equation}
\big|\tfrac{(V_{\cusp})_r(x_*(t),t)}{r(x_*(t))}\big| \le C_{\mathrm{rlog}}\Gamma \qquad (s_0\le t\le s_1,\ r(x_*(t))>0).
\label{eq:cusp-separated-radial-log-growth}
\end{equation}
\end{lemma}

\begin{proof}[Proof of Lemma~\ref{lem:radial-log-growth-separated}]
\runinhead{Step 1: radial quotients for axisymmetric no-swirl velocities.}
We fix an axisymmetric no-swirl velocity $v$.  Since $v_r(0,z)=0$, the fundamental theorem of calculus along the horizontal segment from
the symmetry axis to $(r,z)$ gives, for $r>0$,
\begin{equation}
\tfrac{v_r(r,z)}{r} = \int_0^1 \p_r v_r(\lambda r,z)\,d\lambda ,
\label{eq:radial-log-from-radial-derivative}
\end{equation}
Thus, after we bound $\nabla v$ on this horizontal segment, the same bound controls $v_r/r$ at $(r,z)$.

\runinhead{Step 2: the radial quotient of $u_{\smooth}$.}
On every fixed ball, we use Lemma~\ref{lem:JtwoD-tail-bdd} to obtain the gradient bound
\[
\|\nabla u_{\smooth}(\cdot,t)\|_{L^\infty(B_{2R_0})} \le C_{R_0}\Gamma .
\]
For $|x|$ outside a fixed ball, we decompose source points $y$ in \eqref{eq:smooth-velocity-def} into $B(x,\tfrac18|x|)$ and the shells
\[
\mathscr A_k(x):=\{\,y:\ 2^{k-3}|x|\le |y-x|\le 2^{k-2}|x|\,\},\qquad k\ge0.
\]
After rescaling by $|x|$, we apply \eqref{eq:local-CZ-u0-sec5} to bound the contribution to $\nabla u_{\smooth}(x,t)$ from source points
$y\in B(x,\tfrac18|x|)$, and we apply \eqref{eq:far-CZ-u0-sec5} to bound the contributions to $\nabla u_{\smooth}(x,t)$ from the shells
$\mathscr A_k(x)$.  We then rewrite these source integrals in the labels $Y$.  The transported specific-vorticity identity
\[
\tfrac{\omega_\theta(\phi(Y,t),t)}{\phi_r(Y,t)} = \tfrac{\omega_{\theta,0}(Y)}{R(Y)}
\]
from \eqref{eq:lin-deriv:V-stretch} and the moment estimate \eqref{eq:current-far-vorticity-over-r-moment} imply the dyadic summability in
the source variable needed for \eqref{eq:smooth-velocity-def}.  Consequently,
\begin{equation}
|\nabla u_{\smooth}(x,t)| \le C\Gamma(1+|x|)^{\alpha-\gamma}\le C\Gamma .
\label{eq:smooth-gradient-global-tail}
\end{equation}
We use \eqref{eq:radial-log-from-radial-derivative} with \eqref{eq:smooth-gradient-global-tail} to obtain \eqref{eq:smooth-radial-log-growth}.

\runinhead{Step 3: separation from source points in $\mathcal C_{\inn}$.}
We fix $t\in[s_0,s_1]$ and write $x=x_*(t)$.  By \eqref{eq:sigma-angles},
\eqref{eq:buffered-cone-def}--\eqref{eq:inner-star-cone-notation}, and
\eqref{eq:cone-interior-for-radial-log}, we have
$\mathcal C_{\inn}\Subset\mathcal C_*^{\rm int}\subset\mathcal C_*$.  Hence
$x\notin\mathcal C_*^{\rm int}$ implies that $x$ has a fixed angular separation from every source point
$y\in\mathcal C_{\inn}$.  At such separated pairs we use the axisymmetric Biot--Savart kernel bound
\begin{equation}
\left| \tfrac{K_r(x,y)}{r(x)} \right| + |\nabla_xK(x,y)| \le C|x-y|^{-3}.
\label{eq:separated-radial-log-kernel}
\end{equation}

If $x_*(s_0)\in\partial\mathcal C_*^{\rm int}$, then the boundary condition, together with
\eqref{eq:exact-current-axis-core-geometry} and \eqref{eq:vartheta-RZ-comparison}, implies
\[
c\,J_{\cusp}(s_0)^3 \le \vartheta(Y) \le C\,J_{\cusp}(s_0)^3 .
\]
If $s_0$ is the left endpoint of the connected interval on which
$J_{\cusp}\le\mathfrak J_{\mathrm{rlog}}$, then, after decreasing
$\mathfrak J_{\mathrm{rlog}}$ below $J_{\cusp}(0)$, we have
$J_{\cusp}(s_0)=\mathfrak J_{\mathrm{rlog}}$.  We define $U_{\cusp}^{\rm lead}$ by inserting the leading term
\[
-\Gamma J_{\cusp}(t)^{\alpha-1}r(y)^\alpha\mathfrak A_t(y,t)\bs e_\theta(y)
\]
from \eqref{eq:Omega-cusp-local-representation} into the Biot--Savart integral over source points
$y\in\mathcal C_{\inn}$.  The preceding boundary comparison and the monotonicity of
$J_{\cusp}$ from \eqref{eq:Jdot-two-sided} imply the separated annular estimate
\begin{equation}
\sum_{k\ge0}(2^k|x|)^{-3}
\int_{\mathscr S_k(x)\cap\mathcal C_{\inn}}
\Gamma J_{\cusp}(t)^{\alpha-1}r(y)^\alpha |\mathfrak A_t(y,t)|\,\ud y
\le C\Gamma .
\label{eq:separated-leading-radial-annular-sum}
\end{equation}
The annular summability is the estimate used in \eqref{eq:Omega-cusp-annular-mass}; the boundary relation at
$s_0$ and the monotonicity of $J_{\cusp}$ are what reduce the singular clock power to the right side of
\eqref{eq:separated-leading-radial-annular-sum}.  Applying
\eqref{eq:separated-radial-log-kernel} and \eqref{eq:separated-leading-radial-annular-sum}, we obtain
\begin{equation}
\left| \tfrac{(U_{\cusp}^{\rm lead})_r(x,t)}{r(x)} \right| \le C\Gamma ,
\label{eq:leading-separated-radial-log-bound}
\end{equation}

\runinhead{Step 4: the regular transported vorticity and the smooth pull-back.}
We denote by $U_{\cusp}^{\rm reg}$ the Biot--Savart velocity generated by the regular transported vorticity in
\eqref{eq:Omega-reg-local-bounds}.  The annular summation for disjoint source and observation regions in
\eqref{eq:Ucusp-separated-CZ} yields
\[
\left| \tfrac{(U_{\cusp}^{\rm reg})_r(x,t)}{r(x)} \right| \le C\Gamma .
\]
Together with \eqref{eq:leading-separated-radial-log-bound}, this yields a $C\Gamma$ bound for $(U_{\cusp})_r(x,t)/r(x)$.

We now pass from $U_{\cusp}$ to $V_{\cusp}=(\phi_{\smooth}^{-1})_*u_{\core}$.  Since the smooth map is axisymmetric and preserves the
symmetry axis, we control its radial multiplier.  For each fixed $\xi$ with $r(\xi)>0$, we set
$K_{\smooth}(\xi,t):=(\phi_{\smooth})_r(\xi,t)/r(\xi)$.  By \eqref{eq:smooth-flow-def},
\[
\p_t\log K_{\smooth}(\xi,t)=\tfrac{(u_{\smooth})_r(\phi_{\smooth}(\xi,t),t)}{(\phi_{\smooth})_r(\xi,t)}.
\]
The estimate \eqref{eq:smooth-radial-log-growth} and the time-length bound from \eqref{eq:Jdot-two-sided} imply that $K_{\smooth}(\xi,t)$
and $K_{\smooth}(\xi,t)^{-1}$ are bounded by fixed constants on $[s_0,s_1]$.  We obtain the first-derivative bounds for
$\phi_{\smooth}$ and $\phi_{\smooth}^{-1}$ by using Gronwall's inequality with \eqref{eq:smooth-flow-differential-ode}.  Hence the pull-back
$(\phi_{\smooth}^{-1})_*$ changes the radial quotient by at most a fixed constant, and \eqref{eq:cusp-separated-radial-log-growth} follows.
\end{proof}

The Type--I upper bound for vorticity follows from the exact transport identity \eqref{eq:vort-identity-JtwoD} once we obtain a bound
for $J_{\twoD}(Y,t)^{-1}|\omega_{\theta,0}(Y)|$ that is uniform in the label $Y$.  The next lemma establishes this bound,
including the case where $\phi_{\cusp}(Y,s)$ lies outside $\mathcal C_*^{\rm int}$ defined in \eqref{eq:cone-interior-for-radial-log}.
The clock thresholds used below are $\mathfrak J_{\mathrm{vort}}$ from \eqref{eq:vort-clock-threshold-def},
$\mathfrak J_{\mathrm{tail}}$ from \eqref{eq:tail-clock-threshold-def}, and $\mathfrak J_{\mathrm{rlog}}$ from
\eqref{eq:rlog-clock-threshold-def}.  We define their common restriction for the upper vorticity estimate by
\[
\mathfrak J_{\omega,*}:=\min\{\mathfrak J_{\mathrm{vort}},\mathfrak J_{\mathrm{tail}},\mathfrak J_{\mathrm{rlog}}\}.
\]

\begin{lemma}[$L^\infty$ bound for transported angular vorticity]
\label{lem:global-upper-vorticity-envelope}
There exist constants
\[
C_{\omega,+}<\infty,\qquad \mathfrak J_{\omega,+}\in(0,\mathfrak J_{\omega,*}].
\]
These constants depend only on $\alpha,\gamma,\sigma_{\inn},\sigma_*$.  Let $t$ belong to a time interval on which the Euler velocity remains uniformly bounded in $C^{1,\alpha}$.
If $J:=J_{\cusp}(t)\le\mathfrak J_{\omega,+}$, then, for every label $Y$ with $\omega_{\theta,0}(Y)\neq0$, we have
\begin{equation}
J_{\twoD}(Y,t)^{-1}|\omega_{\theta,0}(Y)|\le C_{\omega,+}\,\Gamma J^{3\alpha-1}.
\label{eq:global-upper-vorticity-envelope}
\end{equation}
\end{lemma}

\begin{proof}[Proof of Lemma~\ref{lem:global-upper-vorticity-envelope}]
We fix $t$ and $Y$, and we set
\[
x_*(s):=\phi_{\cusp}(Y,s), \qquad x(s):=\phi(Y,s)=\phi_{\smooth}(x_*(s),s), \qquad J_s:=J_{\cusp}(s).
\]

\runinhead{Step 1: comparison between physical and cusp-coordinate angular vorticity.}
By \eqref{eq:Omega-cusp-def}, the cusp-flow transported angular vorticity satisfies
\[
\Omega_{\cusp,\theta}(x_*(s),s)=\mathcal J_{\cusp}(Y,s)^{-1}\omega_{\theta,0}(Y),
\]
whereas \eqref{eq:vort-identity-JtwoD} states that the true Euler vorticity satisfies
\[
\omega_\theta(x(s),s)=J_{\twoD}(Y,s)^{-1}\omega_{\theta,0}(Y).
\]
Since $\phi=\phi_{\smooth}\circ\phi_{\cusp}$ and $\phi_{\smooth}$ is axisymmetric, the radial multiplier of the
smooth map relates these two angular vorticities via the relation
\begin{equation}
\omega_\theta(x(s),s)=\tfrac{(\phi_{\smooth})_r(x_*(s),s)}{r(x_*(s))}\,\Omega_{\cusp,\theta}(x_*(s),s).
\label{eq:physical-cusp-vorticity-radial-multiplier}
\end{equation}
Let $t_{\rm rlog}$ be the left endpoint of the connected time interval, containing $t$, on which
$J_{\cusp}\le\mathfrak J_{\mathrm{rlog}}$.  We fix $s\in[t_{\rm rlog},t]$ and we set $\xi=x_*(s)=\phi_{\cusp}(Y,s)$.
For $r(\xi)>0$, define $K_{\smooth}(\xi,\tau):=(\phi_{\smooth})_r(\xi,\tau)/r(\xi)$ for $t_{\rm rlog}\le\tau\le s$.
The finite-clock map bounds \eqref{eq:finite-clock-smooth-cusp-C1} at $t_{\rm rlog}$ and the algebraic-tail bound
\eqref{eq:gamma-stand} show that $C^{-1}\le K_{\smooth}(\xi,t_{\rm rlog})\le C$.  Using the flow equation \eqref{eq:smooth-flow-def}, we obtain
\[
\p_\tau\log K_{\smooth}(\xi,\tau)=\tfrac{(u_{\smooth})_r(\phi_{\smooth}(\xi,\tau),\tau)}{(\phi_{\smooth})_r(\xi,\tau)}.
\]
The estimate \eqref{eq:smooth-radial-log-growth} bounds the right-hand side by $C\Gamma$ for $t_{\rm rlog}\le\tau\le s$.
Moreover, by \eqref{eq:Jdot-two-sided} and the fact that $3\alpha<1$, we obtain that
\[
t-t_{\rm rlog}\le C\Gamma^{-1}\int_J^{\mathfrak J_{\mathrm{rlog}}} q^{-3\alpha}\,dq\le C\Gamma^{-1}.
\]
After integrating the logarithmic identity for $K_{\smooth}$ from $t_{\rm rlog}$ to $s$, we obtain that
\begin{equation}
C^{-1}\le\tfrac{(\phi_{\smooth})_r(x_*(s),s)}{r(x_*(s))}\le C\ \ \text{ for } \ \ J_s\le\mathfrak J_{\omega,+}.
\label{eq:smooth-radial-multiplier-bound}
\end{equation}
By \eqref{eq:physical-cusp-vorticity-radial-multiplier}, it remains to bound
$|\Omega_{\cusp,\theta}(x_*(t),t)|$ by $C\Gamma J^{3\alpha-1}$.

\runinhead{Step 2: the case $x_*(t)\in\mathcal C_*^{\rm int}$.}
Let $x_t:=x_*(t)$.  Since $x_t\in\mathcal C_*^{\rm int}$, \eqref{eq:cone-interior-for-radial-log} implies the inclusion
$B(x_t,2c_*|x_t|)\subset\mathcal C_*$.  Hence, by \eqref{eq:Omega-cusp-local-representation},
\[
\Omega_{\cusp,\theta}(x_t,t)=-\Gamma J^{\alpha-1}r(x_t)^\alpha\mathfrak A_t(x_t,t)+\Omega_{{\rm reg},\theta}(x_t,t).
\]
The pointwise bound for $\mathfrak A_t$ in \eqref{eq:Atr-local-bounds} yields
\[
|\mathfrak A_t(x_t,t)|\le C\bigl(1+J^{-4}|x_t|^2\bigr)^{-\gamma/2},
\]
and therefore, with $\lambda:=J^{-2}|x_t|$,
\[
\Gamma J^{\alpha-1}r(x_t)^\alpha|\mathfrak A_t(x_t,t)|\le C\Gamma J^{\alpha-1}|x_t|^\alpha\bigl(1+J^{-4}|x_t|^2\bigr)^{-\gamma/2}
=C\Gamma J^{3\alpha-1}\lambda^\alpha(1+\lambda^2)^{-\gamma/2}\le C\Gamma J^{3\alpha-1},
\]
because $\gamma>\alpha$.  Using the  $L^\infty$ estimate for $\bs\Omega_{\reg}$ in \eqref{eq:Omega-reg-local-bounds}, it then follows that
\begin{equation}
|\Omega_{\cusp,\theta}(x_*(t),t)|\le C\Gamma J^{3\alpha-1}\qquad\text{if }x_*(t)\in\mathcal C_*^{\rm int}.
\label{eq:cusp-vorticity-upper-inside-cone}
\end{equation}

\runinhead{Step 3: the case $x_*(t)\notin\mathcal C_*^{\rm int}$.}
If $x_*(s)\notin\mathcal C_*^{\rm int}$ for all $s\in[t_{\rm rlog},t]$, then
\eqref{eq:Omega-cusp-def}, \eqref{eq:finite-clock-smooth-cusp-C1} at $t_{\rm rlog}$, and \eqref{eq:gamma-stand} show that
$|\Omega_{\cusp,\theta}(x_*(t_{\rm rlog}),t_{\rm rlog})|\le C\Gamma$.  We apply
\eqref{eq:cusp-separated-radial-log-growth} on $[t_{\rm rlog},t]$, and use \eqref{eq:cusp-map-def} and \eqref{eq:Omega-cusp-def} to write
\[
\p_s\log|\Omega_{\cusp,\theta}(x_*(s),s)|=\tfrac{(V_{\cusp})_r(x_*(s),s)}{r(x_*(s))}.
\]
Integrating this identity over $[t_{\rm rlog},t]$, whose length is at most $C\Gamma^{-1}$, we find that
\[
|\Omega_{\cusp,\theta}(x_*(t),t)|\le C\Gamma.
\]
After decreasing $\mathfrak J_{\omega,+}$, this is bounded by $C\Gamma J^{3\alpha-1}$ because $3\alpha-1<0$.

It remains to consider the case in which $x_*(s)$ enters $\mathcal C_*^{\rm int}$ at least once on
$[t_{\rm rlog},t]$.  Let $t_{\rm ex}\le t$ be the last exit time from $\mathcal C_*^{\rm int}$; then
$x_*(t_{\rm ex})\in\partial\mathcal C_*^{\rm int}$.  Applying \eqref{eq:cusp-vorticity-upper-inside-cone} at times
immediately before $t_{\rm ex}$ and passing to the limit by continuity of $\Omega_{\cusp,\theta}$, we obtain
\[
|\Omega_{\cusp,\theta}(x_*(t_{\rm ex}),t_{\rm ex})|\le C\Gamma J_{\cusp}(t_{\rm ex})^{3\alpha-1}\le C\Gamma J^{3\alpha-1},
\]
since $J_{\cusp}(t_{\rm ex})\ge J$ and $3\alpha-1<0$.  On the interval $[t_{\rm ex},t]$,  the cusp trajectory remains outside $\mathcal C_*^{\rm int}$, so
that \eqref{eq:cusp-separated-radial-log-growth} yields
\[
\left|\tfrac{(V_{\cusp})_r(x_*(s),s)}{r(x_*(s))}\right|\le C\Gamma\qquad (t_{\rm ex}\le s\le t).
\]
The same logarithmic identity from \eqref{eq:cusp-map-def} and \eqref{eq:Omega-cusp-def},
\[
\p_s\log|\Omega_{\cusp,\theta}(x_*(s),s)|=\tfrac{(V_{\cusp})_r(x_*(s),s)}{r(x_*(s))},
\]
together with the interval-length bound from \eqref{eq:Jdot-two-sided}, shows that
\[
|\Omega_{\cusp,\theta}(x_*(t),t)|\le C|\Omega_{\cusp,\theta}(x_*(t_{\rm ex}),t_{\rm ex})|\le C\Gamma J^{3\alpha-1}.
\]
Combining this with \eqref{eq:smooth-radial-multiplier-bound} and
\eqref{eq:physical-cusp-vorticity-radial-multiplier} proves
\eqref{eq:global-upper-vorticity-envelope}.
\end{proof}

\begin{lemma}[Two-sided $L^\infty$ vorticity bounds for the Target Profile]
\label{lem:target-vorticity-envelope}
For the solution of the Euler equations \eqref{eq:euler} with the Target Profile initial condition
\eqref{eq:vort0}, there exist constants $0<c_\omega<C_\omega<\infty$ and
$\mathfrak J_{\omega}\in(0,\mathfrak J_{\omega,+}]$, depending only on
$\alpha,\gamma,\sigma_{\inn},\sigma_*$ and the fixed Target Profile cutoff $\Upsilon$ in
\eqref{eq:Theta-star-def}, such that, whenever $J_{\cusp}(t)\le \mathfrak J_{\omega}$,
\[
c_\omega\,\Gamma\,J_{\cusp}(t)^{3\alpha-1} \le \|\omega(\cdot,t)\|_{L^\infty(\R^3)} \le C_\omega\,\Gamma\,J_{\cusp}(t)^{3\alpha-1}.
\]
\end{lemma}

\begin{proof}[Proof of Lemma~\ref{lem:target-vorticity-envelope}]
We use the exact transport identity
\begin{equation}
\omega_\theta(\phi(Y,t),t)=J_{\twoD}(Y,t)^{-1}\omega_{\theta,0}(Y),
\label{eq:target-envelope-transport-use}
\end{equation}
which is \eqref{eq:vort-identity-JtwoD}.

\runinhead{Step 1. Upper bound.}
We decrease $\mathfrak J_{\omega}$ so that $\mathfrak J_{\omega}\le\mathfrak J_{\omega,+}$.  Then Lemma~\ref{lem:global-upper-vorticity-envelope} provides, 
for every label $Y$ with $\omega_{\theta,0}(Y)\neq0$,
\[
J_{\twoD}(Y,t)^{-1}|\omega_{\theta,0}(Y)|\le C_{\omega,+}\Gamma J_{\cusp}(t)^{3\alpha-1}.
\]
Taking the supremum in $Y$ and using \eqref{eq:target-envelope-transport-use} proves the upper bound.

\runinhead{Step 2. Lower bound.}
We choose a label whose initial polar angle has the same size as the angular opening of the collapsing core.
Fix
\[
0<\kappa\le \min\{\sigma_{\cut},c_{\rm vort}\},
\]
where $c_{\rm vort}$ is the constant in the lower amplification bound
\eqref{eq:driver-amplification-lower}.  For each $t$ with $J_{\cusp}(t)\le\mathfrak J_{\omega}$, choose an
initial label $Y=Y(t)$ in the upper half-space with
\[
\rho(Y)\in[\tfrac12,1], \qquad \sigma(Y)=\kappa J_{\cusp}(t)^3,
\]
and arbitrary azimuthal angle.  Since $R_{\tail}\ge2$ in \eqref{eq:core-tail-domains}, this label belongs to
$D_{\core}$.  Also $J_{\cusp}(t)\le1$ and $\kappa\le\sigma_{\cut}$ imply
$0\le\sigma(Y)\le\sigma_{\cut}$.  Hence the angular cutoff in \eqref{eq:Theta-star-def} satisfies
$\Upsilon(\sigma(Y))=1$, and
\[
\Theta^*(\sigma(Y))=(\sin\sigma(Y))^\alpha.
\]

The angular variable in Lemma~\ref{lem:JtwoD-core-bounds} is
\[
\vartheta(Y)=\frac{R(Y)}{\rho(Y)}=\sin\sigma(Y).
\]
Therefore,
\[
\vartheta(Y)=\sin(\kappa J_{\cusp}(t)^3)\le \kappa J_{\cusp}(t)^3 \le c_{\rm vort}J_{\cusp}(t)^3.
\]
Thus the hypotheses of \eqref{eq:driver-amplification-lower} are satisfied, and
\begin{equation}
J_{\twoD}(Y,t)^{-1}\ge c_{\rm vort}J_{\cusp}(t)^{-1}.
\label{eq:target-envelope-JtwoD-lower-use}
\end{equation}

The Target Profile datum \eqref{eq:vort0} gives
\[
|\omega_{\theta,0}(Y)| = \Gamma\,\frac{\rho(Y)^\alpha}{(1+\rho(Y)^2)^{\gamma/2}}\, (\sin(\kappa J_{\cusp}(t)^3))^\alpha.
\]
Since $\rho(Y)\in[\tfrac12,1]$ and $\sin s\ge c s$ for $0\le s\le\sigma_{\cut}$, we have
\begin{equation}
|\omega_{\theta,0}(Y)|\ge c\Gamma J_{\cusp}(t)^{3\alpha}.
\label{eq:target-envelope-initial-vorticity-lower}
\end{equation}
Combining \eqref{eq:target-envelope-transport-use}, \eqref{eq:target-envelope-JtwoD-lower-use}, and
\eqref{eq:target-envelope-initial-vorticity-lower}, we obtain
\[
|\omega_\theta(\phi(Y,t),t)| \ge c\Gamma J_{\cusp}(t)^{3\alpha-1}.
\]
This proves the lower bound for $\|\omega(\cdot,t)\|_{L^\infty(\R^3)}$.
\end{proof}

The two-sided $L^\infty$ vorticity bounds above are written in powers of the cusp clock $J_{\cusp}(t)$.
To obtain the Type--I rates in Theorem~\ref{thm:target-profile}, we next convert this clock into the remaining
time $T^*-t$.  The conversion uses the clock differential inequality \eqref{eq:Jdot-two-sided} and
$J_{\cusp}(t)\to0$ from Proposition~\ref{prop:blowup-final}.  The comparison between the exact axial strain
and the cusp strain uses \eqref{eq:Wcusp-scaling}, \eqref{eq:modulation-bounds}, \eqref{eq:Jsmooth-bdd}, and
the splitting \eqref{eq:on-axis-strain-splitting}.

\begin{lemma}[Type--I cusp-clock and axial strain rates]
\label{lem:typeI-cusp-rates}
Let $T^*$ denote the blowup time from Proposition~\ref{prop:blowup-final}.  There exist constants
$0<c<C<\infty$ and a time $t_{\rm I}<T^*$ such that, for every $t\in[t_{\rm I},T^*)$,
\begin{equation}
c\,\bigl(\Gamma(T^*-t)\bigr)^{\frac{1}{1-3\alpha}}\le J_{\cusp}(t)\le C\,\bigl(\Gamma(T^*-t)\bigr)^{\frac{1}{1-3\alpha}}.
\label{eq:Jcusp-TypeI}
\end{equation}
and
\begin{equation}
c\,\Gamma J_{\cusp}(t)^{3\alpha-1}\le -\rW_0(t)\le C\,\Gamma J_{\cusp}(t)^{3\alpha-1}.
\label{eq:W0-TypeI-cusp}
\end{equation}
\end{lemma}

\begin{proof}[Proof of Lemma~\ref{lem:typeI-cusp-rates}]
Proposition~\ref{prop:blowup-final} states that
\begin{equation}
J_{\cusp}(t)\to0\qquad\text{as }t\uparrow T^*.
\label{eq:typeI-Jcusp-collapse}
\end{equation}
We choose $t_{\rm I}<T^*$ so that, for every $t\in[t_{\rm I},T^*)$,
\[
J_{\cusp}(t)\le\min\{\mathfrak J_{\mathrm{collapse}},\mathfrak J_{\mathrm{mod}},\mathfrak J_{\mathrm{strain}}\}.
\]
The clock estimate \eqref{eq:Jdot-two-sided} is therefore available on $[t_{\rm I},T^*)$, and yields
\begin{equation}
c\Gamma J_{\cusp}(t)^{3\alpha}\le -\dot J_{\cusp}(t)\le C\Gamma J_{\cusp}(t)^{3\alpha}.
\label{eq:typeI-clock-ode-use}
\end{equation}
Since
\[
\tfrac{d}{dt}(J_{\cusp}(t)^{1-3\alpha})=(1-3\alpha)J_{\cusp}(t)^{-3\alpha}\dot J_{\cusp}(t),
\]
the clock estimate \eqref{eq:typeI-clock-ode-use} is equivalent, after renaming constants, to
\begin{equation}
c\Gamma\le -\tfrac{d}{dt}J_{\cusp}(t)^{1-3\alpha}\le C\Gamma .
\label{eq:typeI-clock-power-ode}
\end{equation}
For $t<s<T^*$, we integrate \eqref{eq:typeI-clock-power-ode} over $[t,s]$ and find that
\begin{equation}
c\Gamma(s-t)\le J_{\cusp}(t)^{1-3\alpha}-J_{\cusp}(s)^{1-3\alpha}\le C\Gamma(s-t).
\label{eq:typeI-clock-integrated}
\end{equation}
Letting $s\uparrow T^*$ in \eqref{eq:typeI-clock-integrated} and using \eqref{eq:typeI-Jcusp-collapse}, we obtain that
\begin{equation}
c\Gamma(T^*-t)\le J_{\cusp}(t)^{1-3\alpha}\le C\Gamma(T^*-t).
\label{eq:typeI-clock-power-bound}
\end{equation}
Because $1-3\alpha>0$, the map $x\mapsto x^{1/(1-3\alpha)}$ is increasing on $(0,\infty)$.  Applying this map
to the three terms in \eqref{eq:typeI-clock-power-bound} and renaming constants proves \eqref{eq:Jcusp-TypeI}.

It remains to compare the exact axial strain with the cusp strain.  The modulation identity
\eqref{eq:modulated-cusp-strain-scaling-sec14} is
\begin{equation}
\rW_{\cusp}(t)=m(t)\mathcal W_{\cusp}(t).
\label{eq:typeI-modulated-cusp-strain-use}
\end{equation}
Combining \eqref{eq:typeI-modulated-cusp-strain-use} with the cusp-coordinate strain estimate
\eqref{eq:Wcusp-scaling} and the modulation bounds \eqref{eq:modulation-bounds}, we obtain that
\begin{equation}
c\,\Gamma J_{\cusp}(t)^{3\alpha-1}\le -\rW_{\cusp}(t)\le C\,\Gamma J_{\cusp}(t)^{3\alpha-1}.
\label{eq:typeI-rWcusp-scale}
\end{equation}
The smooth axial strain obeys
\[
|\rW_{\smooth}(t)|\le C_{\smooth}\Gamma
\]
by Lemma~\ref{lem:smooth-clock-bounded}, specifically the first estimate in \eqref{eq:Jsmooth-bdd}.  Since
$3\alpha-1<0$ and \eqref{eq:typeI-Jcusp-collapse} holds, we increase $t_{\rm I}$ if necessary so that
\begin{equation}
|\rW_{\smooth}(t)|\le \tfrac12 c\,\Gamma J_{\cusp}(t)^{3\alpha-1},\qquad t\in[t_{\rm I},T^*),
\label{eq:typeI-smooth-strain-lower-order}
\end{equation}
where $c$ is the lower bound constant in \eqref{eq:typeI-rWcusp-scale}.  Using \eqref{eq:typeI-rWcusp-scale}, \eqref{eq:typeI-smooth-strain-lower-order}, and
\eqref{eq:on-axis-strain-splitting}, we obtain
\[
\tfrac12 c\,\Gamma J_{\cusp}(t)^{3\alpha-1}\le -\rW_0(t)\le C\,\Gamma J_{\cusp}(t)^{3\alpha-1},
\]
which, after renaming constants, establishes \eqref{eq:W0-TypeI-cusp}.
\end{proof}

\subsection{Proof of Theorem~\ref{thm:target-profile}}

Lemma~\ref{lem:L2} shows that $u_0^*\in C^{1,\alpha}(\R^3)\cap L^2(\R^3)$ for the Target Profile initial
condition \eqref{eq:vort0}.  We use the standard local well-posedness theory for finite-energy axisymmetric
no-swirl data in $C^{1,\alpha}$, going back to Lichtenstein~\cite{Lichtenstein1925}, together with the
Beale--Kato--Majda continuation criterion \cite{BKM1984}, to obtain the unique local $C^{1,\alpha}$ Euler
solution used in the construction.  Proposition~\ref{prop:blowup-final} gives a finite blowup time
$T^*<\infty$ and collapse of the cusp and physical clocks.  It remains only to prove the Type--I rates stated
in Theorem~\ref{thm:target-profile}.  All constants below may change from line to line; they are independent
of $\Gamma$ and $t$.

By Lemma~\ref{lem:typeI-cusp-rates}, the axial strain satisfies
\begin{equation}
c\Gamma J_{\cusp}(t)^{3\alpha-1}\le -\rW_0(t)\le C\Gamma J_{\cusp}(t)^{3\alpha-1},
\label{eq:target-proof-axial-strain-cusp-clock}
\end{equation}
and the cusp clock satisfies
\begin{equation}
c\bigl(\Gamma(T^*-t)\bigr)^{\frac{1}{1-3\alpha}}\le J_{\cusp}(t)\le C\bigl(\Gamma(T^*-t)\bigr)^{\frac{1}{1-3\alpha}}.
\label{eq:target-proof-cusp-clock-time}
\end{equation}
Raising \eqref{eq:target-proof-cusp-clock-time} to the power $1-3\alpha$ and renaming constants gives
\begin{equation}
c\Gamma(T^*-t)\le J_{\cusp}(t)^{1-3\alpha}\le C\Gamma(T^*-t).
\label{eq:target-proof-cusp-clock-power-time}
\end{equation}
Since $3\alpha-1=-(1-3\alpha)$, \eqref{eq:target-proof-cusp-clock-power-time} gives
\begin{equation}
c(T^*-t)^{-1}\le \Gamma J_{\cusp}(t)^{3\alpha-1}\le C(T^*-t)^{-1}.
\label{eq:target-proof-singular-scale-time}
\end{equation}
Combining \eqref{eq:target-proof-axial-strain-cusp-clock} and
\eqref{eq:target-proof-singular-scale-time}, we obtain
\[
c(T^*-t)^{-1}\le -\rW_0(t)\le C(T^*-t)^{-1}.
\]
This is the axial-strain rate in Theorem~\ref{thm:target-profile}, because
$\rW_0(t)=\p_z u_z(0,0,t)$ by \eqref{eq:rW2}.

Next we prove the physical clock law.  The clock decomposition \eqref{eq:J-clock-decomposition} gives
\begin{equation}
J(t)=J_{\smooth}(t)J_{\cusp}(t).
\label{eq:target-proof-clock-decomposition-use}
\end{equation}
Since $J_{\cusp}(t)\to0$ as $t\uparrow T^*$, the small-clock condition in
Lemma~\ref{lem:smooth-clock-bounded} holds on a terminal interval.  Hence \eqref{eq:Jsmooth-bdd} gives
\begin{equation}
c_{\smooth}\le J_{\smooth}(t)\le C_{\smooth}
\label{eq:target-proof-smooth-clock-bounded-use}
\end{equation}
there.  Combining \eqref{eq:target-proof-clock-decomposition-use},
\eqref{eq:target-proof-smooth-clock-bounded-use}, and \eqref{eq:target-proof-cusp-clock-time}, we obtain
\begin{equation}
c\bigl(\Gamma(T^*-t)\bigr)^{\frac{1}{1-3\alpha}}\le J(t)\le C\bigl(\Gamma(T^*-t)\bigr)^{\frac{1}{1-3\alpha}},
\label{eq:target-proof-physical-clock-typeI}
\end{equation}
which is the clock law in Theorem~\ref{thm:target-profile}.

It remains to prove the Type--I vorticity rate.  Lemma~\ref{lem:target-vorticity-envelope} gives
\begin{equation}
c\Gamma J_{\cusp}(t)^{3\alpha-1}\le \|\omega(\cdot,t)\|_{L^\infty(\R^3)} \le C\Gamma J_{\cusp}(t)^{3\alpha-1}.
\label{eq:target-proof-vorticity-cusp-clock}
\end{equation}
Using \eqref{eq:target-proof-singular-scale-time} in \eqref{eq:target-proof-vorticity-cusp-clock}, we obtain
\begin{equation}
c(T^*-t)^{-1}\le \|\omega(\cdot,t)\|_{L^\infty(\R^3)}\le C(T^*-t)^{-1}.
\label{eq:target-proof-vorticity-typeI}
\end{equation}
Therefore, for any terminal time $t_{\rm I}<T^*$ on which \eqref{eq:target-proof-vorticity-typeI} holds,
\[
\int_0^{T^*}\|\omega(\cdot,t)\|_{L^\infty(\R^3)}\,\ud t \ge c\int_{t_{\rm I}}^{T^*}(T^*-t)^{-1}\,\ud t=\infty .
\]
This is the Beale--Kato--Majda divergence asserted in Theorem~\ref{thm:target-profile}, and completes the
proof of the theorem.


\section{Stability of Blowup for Admissible Initial Data}
\label{sec:proof-main}

We now prove Theorem~\ref{thm:main}.  The target-profile collapse argument in
Sections~\ref{sec:Euler-blowup-for-Theta-star}--\ref{sec:target-profile-typeI-completion} was arranged so that
small weighted angular perturbations can be absorbed after the pressure localization, angular cutoff, and
small-clock thresholds have been fixed.  We use the same bootstrap quantities with the angular profile
\[
\Theta(\sigma)=\Theta^*(\sigma)(1+h(\sigma)), \qquad h(\sigma)=(\sin\sigma)^\eta k(\sigma), \qquad \|k\|_{C^\alpha([0,\pi/2])}<\nu .
\]
The additional weight $(\sin\sigma)^\eta$ is the useful smallness.  In the strain-producing region, the
initial Lagrangian angle has size $O(J_{\cusp}^3)$, so the angular perturbation is lower order than the
Target Profile angular function.  The superscript $*$ denotes the Target Profile angular function transported by
the same perturbed cusp map, while superscript $\nu$ denotes the full admissible profile.\footnote{Here the superscript $*$ refers only to the Target Profile part of the
perturbed initial datum.  For the admissible datum
$\Theta^\nu=\Theta^*(1+h^\nu)$, we split $\Theta^\nu=\Theta^*+\Theta^*h^\nu$.
After this split, both terms are transported by the flow of the admissible
solution $u^\nu$.  Thus the difference between the $\nu$-quantity and the
$*$-quantity comes only from the angular correction $\Theta^*h^\nu$.}
Constants may depend on $\alpha$, $\gamma$, $\eta$, and on the fixed cutoffs in Definition~\ref{def:init-data}, but not on $\Gamma$
or on the particular admissible perturbation.

\subsection{Perturbative angular bounds}

We first use the weighted topology in Definition~\ref{def:init-data}.  The initial angular vorticity can be
written as
\[
\omega_{\theta,0}^{\nu}(\rho,\sigma)=-\Gamma\,\tfrac{\rho^\alpha}{(1+\rho^2)^{\gamma/2}}\,
\Theta^*(\sigma)(1+h(\sigma))=\omega_{\theta,0}^*+g_0 .
\]
Since $\Theta^*(\sigma)=(\sin\sigma)^\alpha\Upsilon(\sigma)$ and
$|h(\sigma)|\le \nu(\sin\sigma)^\eta$, the perturbation satisfies
\[
|g_0(\rho,\sigma)|\le C\Gamma\nu\,\tfrac{\rho^\alpha}{(1+\rho^2)^{\gamma/2}}(\sin\sigma)^{\alpha+\eta} .
\]
In cylindrical labels $Y=(R,Z)$, with $\rho(Y)=(R^2+Z^2)^{\frac12}$, this becomes
\begin{equation}
|g_0(Y)|\le C\Gamma\nu\,\tfrac{R^{\alpha+\eta}}{(1+\rho(Y)^2)^{\gamma/2}\rho(Y)^\eta}.
\label{eq:g0-poloidal-gain}
\end{equation}
The local H\"older bound follows from the ordinary product estimate on balls whose radius is comparable to the
label size.  On $B_Y:=B(Y,c|Y|)$, with $B(Y,2c|Y|)\subset\{R\ge0\}$, we write, for $Y'\in B_Y$,
\[
g_0(Y')=-\Gamma H(Y')k(\sigma(Y')), \qquad H(Y'):=\tfrac{R(Y')^{\alpha+\eta}}{(1+\rho(Y')^2)^{\gamma/2}\rho(Y')^\eta}\Upsilon(\sigma(Y')) .
\]
Since $\rho(Y')\simeq |Y|$ and $|\nabla\sigma(Y')|\le C|Y|^{-1}$ on $B_Y$, we have
\[
\|H\|_{L^\infty(B_Y)}\le C\tfrac{|Y|^\alpha}{(1+|Y|^2)^{\gamma/2}},
\qquad [H]_{C^\alpha(B_Y)}\le C\tfrac{|Y|^\alpha}{(1+|Y|^2)^{\gamma/2}}|Y|^{-\alpha} .
\]
Also $\|k\circ\sigma\|_{L^\infty(B_Y)}\le\nu$ and
$[k\circ\sigma]_{C^\alpha(B_Y)}\le C\nu |Y|^{-\alpha}$.  Therefore
\begin{equation}
[g_0]_{C^\alpha(B(Y,c|Y|))}\le C\Gamma\nu\,\tfrac{|Y|^\alpha}{(1+|Y|^2)^{\gamma/2}}|Y|^{-\alpha} \quad\text{when
}B(Y,2c|Y|)\subset\{R\ge0\}.
\label{eq:g0-holder-gain}
\end{equation}
On compact label sets the same argument gives the bound $[g_0]_{C^\alpha}\le C\Gamma\nu$.  These are the
additional estimates needed for the admissible angular profiles.

In the variables used in Lemma~\ref{lem:late-axis-normal-form-cusp}, the localized bounded-slope labels are
$Y_t(\zeta,\tau)$ with $|\tau|\le C_0$.  Their Eulerian images have size
$J_{\cusp}(t)^2(\zeta\tau,\zeta)$ up to the normal-form error in \eqref{eq:additive-normal-form-bound}.  The
corresponding initial Lagrangian angle satisfies
\begin{equation}
\sin\sigma(Y_t(\zeta,\tau))\le C J_{\cusp}(t)^3(1+|\tau|)
\quad\text{for }\zeta\in I_\sharp,\ |\tau|\le C_0.
\label{eq:lag-angle-pert-gain}
\end{equation}
This is the angular drift estimate used in Lemma~\ref{lem:late-entry} and
Lemma~\ref{lem:late-axis-normal-form-cusp}.  Combining \eqref{eq:weighted-proximity} with
\eqref{eq:lag-angle-pert-gain}, we obtain
\[
|h(\sigma(Y_t(\zeta,\tau)))|\le C\nu J_{\cusp}(t)^{3\eta}(1+|\tau|)^\eta \quad\text{for }\zeta\in I_\sharp, \ |\tau|\le C_0.
\]
Let $\Omega_{\sharp}^{\nu}$ and $\Omega_{\sharp}^{*,\nu}$ denote the localized transported vorticities defined
by \eqref{eq:Omega-sharp-cutoff-def}, using the same perturbed cusp map and the same labels
$Y_t(\zeta,\tau)$, but using the angular profiles $\Theta=\Theta^*(1+h)$ and $\Theta^*$, respectively.  We set
\[
\Omega_{\sharp}^{\pert,\nu}:=\Omega_{\sharp}^{\nu}-\Omega_{\sharp}^{*,\nu}.
\]
For functions supported on the localized bounded-slope cone, the scaled norm below is the ordinary
$C^{\alpha/2}$ norm after pullback to the fixed variables $(\zeta,\tau)$:
\[
\|F(\cdot,t)\|_{C^{\alpha/2}_{\rm sc}}:=\|F(\phi_{\cusp}(Y_t(\zeta,\tau),t),t)\|_{C^{\alpha/2}_{\zeta,\tau}
(\operatorname{supp}\vartheta_\sharp\times[-C_0,C_0])}.
\]
We apply this norm to the full perturbative angular function in the vorticity.  On $[0,\pi/2]$,
\[
\Theta^*(\sigma)h(\sigma)=(\sin\sigma)^{\alpha+\eta}\Upsilon(\sigma)k(\sigma),
\qquad \|k\|_{C^\alpha}\le\nu \quad\text{by }\eqref{eq:weighted-proximity}.
\]
The map $(\zeta,\tau)\mapsto Y_t(\zeta,\tau)$ has uniformly bounded $C^{1,\beta_{\rm ax}}$ norm on
$\operatorname{supp}\vartheta_\sharp\times[-C_0,C_0]$ by \eqref{eq:additive-normal-form-bound} and
\eqref{eq:normal-form-approximation-bound}.  Thus $k(\sigma(Y_t(\zeta,\tau)))$ has $C^{\alpha/2}$ norm
bounded by $C\nu$ on this fixed set.  Moreover, \eqref{eq:lag-angle-pert-gain} and the same map bounds give
\[
\|(\sin\sigma(Y_t(\zeta,\tau)))^{\alpha+\eta}\|_{C^{\alpha/2}_{\zeta,\tau}} \le C J_{\cusp}(t)^{3(\alpha+\eta)}.
\]
At $\tau=0$, the label $Y_t(\zeta,\tau)$ lies on the symmetry axis, so the only possible singular point in the
$C^{\alpha/2}_{\zeta,\tau}$ seminorm is the $\tau$-difference across $\tau=0$.  Since the angular factor
vanishes like $|\tau|^{\alpha+\eta}$ there and $\alpha+\eta>\tfrac\alpha 2$, this term is $C^{\alpha/2}$ in
$\tau$ at the origin.  The cutoffs $\vartheta_\sharp(\zeta)$ and
$\chi_{M_\pressure}(|\tau|)$ are fixed smooth functions on the same set, and the zero extension across their supports
preserves the displayed $C^{\alpha/2}$ bound.  Hence
\begin{equation}
\|\Omega_{\sharp}^{\pert,\nu}(\cdot,t)\|_{C^{\alpha/2}_{\rm sc}} \le C\nu J_{\cusp}(t)^{3\eta}\Gamma J_{\cusp}(t)^{3\alpha-1}.
\label{eq:localized-vort-pert}
\end{equation}
The localized pressure estimates \eqref{eq:fixed-scale-pressure-bilinear} and
\eqref{eq:scaled-pressure-hessian-error}, applied with the perturbative bound
\eqref{eq:localized-vort-pert}, give the localized-cone pressure contribution with the additional gain
$\nu J_{\cusp}^{3\eta}$.  The region where the axial cutoff $1-\vartheta_\sharp$ is active, the region where
the angular cutoff $1-\chi_{M_\pressure}$ is active, the image-map displacement controlled by
\eqref{eq:normal-form-approximation-bound}, and the algebraic tail are estimated by the cusp-coordinate
bounds used in Lemma~\ref{lem:transported-cusp-pressure-win}.  Therefore
\begin{subequations}
\begin{align}
|\mathcal W_{\cusp}^{\nu}(t)-\mathcal W_{\cusp}^{*}(t)|
&\le C\nu\Gamma J_{\cusp}(t)^{3\alpha-1},
\label{eq:W-pert-gain}\\
|\Pi_{\cusp}^{\nu}(t)-\Pi_{\cusp}^{*}(t)|
&\le C\nu\Gamma^2J_{\cusp}(t)^{6\alpha-2}.
\label{eq:Pi-pert-gain}
\end{align}
\end{subequations}
The estimates \eqref{eq:g0-poloidal-gain}--\eqref{eq:g0-holder-gain} also control the perturbative angular
function away from the localized bounded-slope cone.

We shall use the following definition.
\begin{definition}[Small-clock stable admissible solutions]
\label{def:small-clock-stable}
Let $u^\nu$ be a solution of the incompressible Euler equations \eqref{eq:euler} whose initial vorticity belongs
to the admissible class $\mathcal A_{\alpha,\gamma}(\nu,\eta)$ from Definition~\ref{def:init-data}.  Thus
$u^\nu$ has angular vorticity
\[
\omega_{\theta,0}^\nu(\rho,\sigma)=-\Gamma\,\tfrac{\rho^\alpha}{(1+\rho^2)^{\gamma/2}}\,\Theta^\nu(\sigma),
\qquad \Theta^\nu(\sigma)=\Theta^*(\sigma)(1+h^\nu(\sigma)),
\]
where
\[
\Theta^*(\sigma):=(\sin\sigma)^\alpha\Upsilon(\sigma), \qquad \sigma\in[0,\tfrac\pi2],
\]
as in \eqref{eq:Theta-star-def}.  We say that $u^\nu$ is small-clock stable on a time interval $\mathcal I$ if,
on $\mathcal I$, the estimates from the Target Profile stability argument hold with the variables, maps,
velocities, clocks, axial strains, and pressure quantities generated by the same solution $u^\nu$.\footnote{These
estimates are the entry axis bounds \eqref{eq:entry-axis-bounds-statement}; the cusp-map normal form estimates
\eqref{eq:additive-normal-form-bound}--\eqref{eq:normal-form-approximation-bound}; the renormalized axis chart
\eqref{eq:renormalized-axis-chart-bootstrap}; the axis-composition distortion \eqref{eq:axis-S-distortion}; the
transfer of axis bounds to the $\zeta$ coordinate \eqref{eq:current-axis-transfer-bound}; the Euler-generated axial
function bounds \eqref{eq:current-axis-extension-size}; the cusp-coordinate strain and velocity bounds
\eqref{eq:Wcusp-scaling}--\eqref{eq:Ucusp-radial-defect}; the smooth clock bounds \eqref{eq:Jsmooth-bdd}; the
smooth-flow deformation bounds \eqref{eq:smooth-flow-position-small}--\eqref{eq:smooth-flow-second-gradient-small};
the two-dimensional Jacobian bounds in the core \eqref{eq:JtwoD-core-upper}--\eqref{eq:driver-amplification-lower};
the radial logarithmic growth bounds
\eqref{eq:smooth-radial-log-growth}--\eqref{eq:cusp-separated-radial-log-growth}; the scalar modulation bounds
\eqref{eq:modulation-bounds}; and the cusp-clock differential inequality \eqref{eq:Jdot-two-sided}.}  The
constants in these estimates are required to be independent of $\nu$ and $\Gamma$.

The estimates required in this definition are all evaluated inside the single admissible
solution $u^\nu$.  The superscript $*$ refers to the Target Profile angular function $\Theta^*$ transported by
the perturbed cusp map of $u^\nu$, while the superscript $\nu$ refers to the full admissible angular function
$\Theta^\nu$.
\end{definition}

The next lemma proves that admissible angular perturbations enter the small-cusp-clock regime in finite time and
then satisfy Definition~\ref{def:small-clock-stable}.  The Riccati pressure lower bound for the full perturbed
pressure Hessian is deliberately left out of Definition~\ref{def:small-clock-stable}; it is proved afterward in
Lemma~\ref{lem:perturbed-pressure-hessian-slack}.

\begin{lemma}[Persistence of small-clock stability under admissible perturbations]
\label{lem:pert-small-clock-persistence}
There exists $\nu_0>0$, depending only on $\alpha,\gamma,\eta$ and the fixed cutoffs, such that the following
holds for every $0<\nu\le\nu_0$.  Let $u^\nu$ be the solution of the incompressible Euler equations
\eqref{eq:euler} whose initial vorticity belongs to $\mathcal A_{\alpha,\gamma}(\nu,\eta)$.  For every
$\mathfrak J_{\mathrm{finite}}\in(0,1)$, define
\[
t_{\rm ent}^{\nu}(\mathfrak J_{\mathrm{finite}}):=\inf\{\,t\ge0:\ J_{\cusp}^{\nu}(t)\le \mathfrak J_{\mathrm{finite}}\,\}.
\]
Then the cusp-coordinate axial strain satisfies
\begin{equation}
\rW_{\cusp}^{\nu}(t)\le -c_{\rm ent}(\mathfrak J_{\mathrm{finite}})\Gamma \ \ \text{ for } \ \
0\le t\le t_{\rm ent}^{\nu}(\mathfrak J_{\mathrm{finite}}),
\qquad t_{\rm ent}^{\nu}(\mathfrak J_{\mathrm{finite}})\le C_{\rm ent}(\mathfrak J_{\mathrm{finite}})\Gamma^{-1}.
\label{eq:pert-late-entry-bound}
\end{equation}
After the entry time $t_{\rm ent}^{\nu}(\mathfrak J_{\mathrm{finite}})$, the solution $u^\nu$ is small-clock
stable in the sense of Definition~\ref{def:small-clock-stable}.  On the localized bounded-slope cone in the
coordinates $(\zeta,\tau)$, the perturbative angular part satisfies \eqref{eq:localized-vort-pert}, and the
cusp-coordinate axial strain and pressure satisfy \eqref{eq:W-pert-gain}--\eqref{eq:Pi-pert-gain}.
\end{lemma}
\begin{proof}[Proof of Lemma~\ref{lem:pert-small-clock-persistence}]
\runinhead{Step 1. Finite entry into the prescribed small-cusp-clock regime.}
We fix $\mathfrak J_{\mathrm{finite}}\in(0,1)$ and use the compact entry sector $E_{\rm ent}$ from
Lemma~\ref{lem:finite-clock-driver-sector}.  On this sector, the Target Profile angular function satisfies
\[
\Theta^*(\sigma)=(\sin\sigma)^\alpha, \qquad -\omega_{\theta,0}^*(Y)\ge c_{\rm ent}\Gamma \ \ \text{ for } \ \ Y\in E_{\rm ent} .
\]
For the admissible datum, Definition~\ref{def:init-data} gives
\[
\omega_{\theta,0}^{\nu}(Y)=\omega_{\theta,0}^*(Y)(1+h^\nu(\sigma(Y))), \qquad \|h^\nu\|_{L^\infty([0,\pi/2])}\le\nu .
\]
We decrease $\nu_0$ so that $1+h^\nu\ge \tfrac12$ for every $0<\nu\le\nu_0$.  Then
\[
-\omega_{\theta,0}^{\nu}(Y)\ge c\Gamma \ \ \text{ for } \ \ Y\in E_{\rm ent} .
\]
The proof of Lemma~\ref{lem:finite-clock-driver-sector} is a compact $C^{1,\alpha}$ argument on the range
$J_{\cusp}^{\nu}(t)\in[\mathfrak J_{\mathrm{finite}},1]$.  The perturbation estimates
\eqref{eq:g0-poloidal-gain}--\eqref{eq:g0-holder-gain} give an $O(\nu)$ change of the initial velocity in the
cone-local $C^\alpha$ norms used in that compact argument.  After decreasing $\nu_0$, the compact estimates
\eqref{eq:finite-clock-Phi-C1} and \eqref{eq:finite-clock-smooth-cusp-C1} remain valid for the admissible
solution.  The kernel $\mathcal K_{\cusp}^{\nu}$ in the admissible analogue of \eqref{eq:Kcusp-label-rep}
therefore satisfies
\[
\mathcal K_{\cusp}^{\nu}(Y,t)\omega_{\theta,0}^{\nu}(Y)\le0 \ \ \text{ for } \ \ \omega_{\theta,0}^{\nu}(Y)\ne0,
\qquad \mathcal K_{\cusp}^{\nu}(Y,t)\ge c_{\rm ent} \ \ \text{ for } \ \ Y\in E_{\rm ent} .
\]
Restricting the sign-definite label representation \eqref{eq:Kcusp-label-rep} to $E_{\rm ent}$ gives
\[
\rW_{\cusp}^{\nu}(t)\le -c_{\rm ent}(\mathfrak J_{\mathrm{finite}})\Gamma \ \ \text{ for } \ \
J_{\cusp}^{\nu}(t)\in[\mathfrak J_{\mathrm{finite}},1] .
\]
Using the clock identity
\[
\p_t\log J_{\cusp}^{\nu}(t)=\tfrac12\rW_{\cusp}^{\nu}(t),
\]
and integrating until the first time $J_{\cusp}^{\nu}(t)=\mathfrak J_{\mathrm{finite}}$, we obtain that
\eqref{eq:pert-late-entry-bound}.
It remains only to justify that this first time is reached before any possible maximal existence time.  Suppose
that the admissible solution is smooth on $[0,T)$ and
$J_{\cusp}^{\nu}(t)\ge\mathfrak J_{\mathrm{finite}}$ on this interval.  The finite-clock estimates just recalled
give uniform $C^1$ bounds for the cusp map, the smooth map, and their inverses on the bounded core.  Together
with \eqref{eq:g0-poloidal-gain}--\eqref{eq:g0-holder-gain}, the exact transport identity
\eqref{eq:vort-identity-JtwoD} then gives
\[
\|\omega^\nu(\cdot,t)\|_{L^\infty}\le C(\alpha,\gamma,\eta,\mathfrak J_{\mathrm{finite}})\Gamma
\qquad\bigl(J_{\cusp}^{\nu}(t)\in[\mathfrak J_{\mathrm{finite}},1]\bigr).
\]
The Beale--Kato--Majda continuation criterion rules out finite breakdown while
$J_{\cusp}^{\nu}\ge\mathfrak J_{\mathrm{finite}}$.  Hence the entry time is finite and satisfies
\eqref{eq:pert-late-entry-bound}.

\runinhead{Step 2. Decomposition of the transported cusp vorticity after entry.}
We next prove the small-clock stability asserted in Definition~\ref{def:small-clock-stable}.  In cusp
coordinates, the transported angular vorticity for the admissible solution decomposes as
\begin{equation}
\Omega_{\cusp}^{\nu}=\Omega_{\cusp}^{*,\nu}+\Omega_{\cusp}^{\pert,\nu}.
\label{eq:pert-cusp-vorticity-decomposition}
\end{equation}
Here $\Omega_{\cusp}^{*,\nu}$ is obtained from the Target Profile angular function $\Theta^*$, but it is
transported by the cusp map generated by the admissible solution $u^\nu$.  The term
$\Omega_{\cusp}^{\pert,\nu}$ is generated by $\Theta^*h^\nu$ and is measured by
\eqref{eq:localized-vort-pert} on the localized bounded-slope cone.

The estimates for $\Omega_{\cusp}^{*,\nu}$ are the Target Profile estimates evaluated in the admissible
geometry.  Their hypotheses are precisely the axis, normal-form, image-map, axial-profile, smooth-flow, and
clock estimates required in Definition~\ref{def:small-clock-stable}.  Since the Target Profile proof improves
these estimates with strict margins after the cutoffs and thresholds in Subsection~\ref{sec:fixed-choice-order}
are fixed, it remains only to verify that $\Omega_{\cusp}^{\pert,\nu}$ changes the corresponding bounds by
$O(\nu)$.

\runinhead{Step 3. Bounds for the perturbative angular contribution.}
On the localized bounded-slope cone in the variables $(\zeta,\tau)$, the perturbative angular term obeys
\eqref{eq:localized-vort-pert}.  The Calder\'on--Zygmund estimate in the variables obtained after dividing the
Eulerian image by $J_{\cusp}^2$, together with the strain and pressure estimates used in
Lemmas~\ref{lem:transported-cusp-field-bounds} and \ref{lem:transported-cusp-pressure-win}, gives an additional
gain
\[
C\nu J_{\cusp}(t)^{3\eta}
\]
relative to the Target Profile scale on this cone.  Thus the localized cusp-coordinate axial strain and
pressure satisfy \eqref{eq:W-pert-gain}--\eqref{eq:Pi-pert-gain} after decreasing $\nu_0$.

The complementary regions are controlled by the estimates attached to their definitions: the region where
$1-\vartheta_\sharp$ is active is controlled by the $\zeta$-tail estimate \eqref{eq:pressure-zeta-tail}; the
region where $1-\chi_{M_\pressure}$ is active is controlled by the angular tail
\eqref{eq:pressure-angular-tail}; the image-map displacement after division by $J_{\cusp}^2$ is controlled by
\eqref{eq:normal-form-approximation-bound}; the smooth-flow deformation is controlled by
\eqref{eq:smooth-flow-position-small}--\eqref{eq:smooth-flow-second-gradient-small}; and the far labels are
controlled by the algebraic-tail pressure estimate \eqref{eq:tail-pressure-remainder-bound}.  On all of these
regions, the perturbative angular function is bounded by
\eqref{eq:g0-poloidal-gain}--\eqref{eq:g0-holder-gain}.  Hence these regions contribute only $O(\nu)$ changes
to the estimates required in Definition~\ref{def:small-clock-stable}.

\runinhead{Step 4. Axial composition and the only time-integrability check.}
The axial composition estimates require one additional check because the admissible angular perturbation enters
the evolution of the axial H\"older mode and the integrated distortion variable $\mathcal P_t$.  On the
localized $\zeta$-interval, the new terms are bounded by
\[
C\nu m(t)\Gamma J_{\cusp}(t)^{3\alpha-1+3\eta} +C\nu\Gamma\bigl(J_{\cusp}(t)^{9\alpha-1+3\eta}+1\bigr),
\]
with the same $C^{\alpha/2}$ control in the axial label.  The clock inequality already used in the bootstrap
closure gives
\[
\ud t\le C\,\tfrac{-\ud J_{\cusp}}{\Gamma J_{\cusp}^{3\alpha}}.
\]
Therefore the time integral of the preceding perturbative terms is bounded by
\begin{equation}
C\nu\int_0^{\mathfrak J_{\mathrm{axis}}} \bigl(J^{3\eta-1}+J^{6\alpha+3\eta-1}+J^{-3\alpha}\bigr)\,\ud J\le C\nu .
\label{eq:pert-axial-integrated-error}
\end{equation}
The integral in \eqref{eq:pert-axial-integrated-error} is finite because $\eta>0$ and $\alpha<1/3$.  Thus the
perturbation changes the renormalized axis chart \eqref{eq:renormalized-axis-chart-bootstrap}, the
axis-composition distortion \eqref{eq:axis-S-distortion}, and the transfer of axis bounds to the $\zeta$
coordinate \eqref{eq:current-axis-transfer-bound} by at most $C\nu$.

\runinhead{Step 5. Bootstrap closure after entry.}
The Target Profile estimates required in Definition~\ref{def:small-clock-stable} were proved with strict margins
after the fixed-choice order in Subsection~\ref{sec:fixed-choice-order}.  Steps 2--4 show that the admissible
angular perturbation changes each of these estimates by at most $C\nu$.  We
therefore decrease $\nu_0$ so that all these changes remain within the reserved margins.

The perturbed cusp-clock differential inequality in \eqref{eq:Jdot-two-sided} gives
$\dot J_{\cusp}^{\nu}(t)<0$ on the small-clock interval.  Hence, after the entry time
$t_{\rm ent}^{\nu}(\mathfrak J_{\mathrm{finite}})$, the solution cannot leave the small-clock regime by crossing
back through the entry threshold.  The standard continuity argument for the bootstrap interval then proves
that $u^\nu$ is small-clock stable in the sense of Definition~\ref{def:small-clock-stable} until either
$J_{\cusp}^{\nu}$ collapses or the $C^{1,\alpha}$ Euler solution reaches its maximal existence time.
\end{proof}

\subsection{Riccati stability of the pressure Hessian lower bound}

By Lemma~\ref{lem:pert-small-clock-persistence}, the Target Profile angular part transported by the perturbed
cusp map satisfies the geometric, axial-profile, normal-form, and field hypotheses used in
Lemma~\ref{lem:transported-cusp-pressure-win}.  Applying that lemma in the persisted perturbed geometry gives a
constant $q_{\rm tr}^*$ with $0<q_{\rm tr}^*<\upbeta$ such that
\begin{equation}
\Pi_{\cusp}^{*}(t)\ge -q_{\rm tr}^*\,\tfrac12\bigl(\mathcal W_{\cusp}^{*}(t)\bigr)^2 .
\label{eq:target-cusp-pressure-win-recalled}
\end{equation}
We fix constants
\begin{equation}
q_{\rm tr}^*<q_{\rm tr}^{\rm pert}<q_{\rm rem}^{\rm pert}<\upbeta
\label{eq:pert-pressure-slack-hierarchy}
\end{equation}
and then choose $\varepsilon_W>0$ so small that
\begin{equation}
\tfrac{q_{\rm rem}^{\rm pert}}{(1-\varepsilon_W)^2}<\upbeta .
\label{eq:pert-strain-comparison-eps-choice}
\end{equation}
The constants in \eqref{eq:pert-pressure-slack-hierarchy}--\eqref{eq:pert-strain-comparison-eps-choice} are
fixed before $\nu_0$ is decreased and before the perturbed pressure threshold
$\mathfrak J_{\Pi}^{\rm pert}$ in Lemma~\ref{lem:perturbed-pressure-hessian-slack} is chosen.

For the admissible perturbation we write
\[
U_{\cusp}^{\nu}=U_{\cusp}^{*}+U_{\cusp}^{\pert}.
\]
Polarization gives
\[
\Pi[U_{\cusp}^{\nu},U_{\cusp}^{\nu}] =\Pi[U_{\cusp}^{*},U_{\cusp}^{*}]
+2\Pi[U_{\cusp}^{*},U_{\cusp}^{\pert}]+\Pi[U_{\cusp}^{\pert},U_{\cusp}^{\pert}].
\]
The bilinear pressure estimates \eqref{eq:scaled-pressure-hessian-error} and
\eqref{eq:fixed-scale-pressure-bilinear}, together with the localized perturbative vorticity bound
\eqref{eq:localized-vort-pert}, control the localized bounded-slope cone in the last two terms.  The region
where the axial cutoff $1-\vartheta_\sharp$ is active is controlled by the tail estimate
\eqref{eq:pressure-zeta-tail}; the region where the angular cutoff $1-\chi_{M_\pressure}$ is active is
controlled by \eqref{eq:pressure-angular-tail}; the image-map displacement after division by $J_{\cusp}^2$ is
controlled by \eqref{eq:normal-form-approximation-bound}; and the far labels are controlled by the algebraic
tail estimate \eqref{eq:tail-pressure-remainder-bound}.  These are the estimates used in
Lemma~\ref{lem:transported-cusp-pressure-win}, now applied bilinearly with one perturbative angular part.
After the cutoffs in Subsection~\ref{sec:fixed-choice-order} have been chosen and $\nu_0$ is decreased if
necessary,
\begin{equation}
\big|2\Pi[U_{\cusp}^{*},U_{\cusp}^{\pert}]+\Pi[U_{\cusp}^{\pert},U_{\cusp}^{\pert}]\big|
\le C\nu\bigl(\mathcal W_{\cusp}^{*}(t)\bigr)^2 .
\label{eq:pressure-polarization-pert-bound}
\end{equation}
The strain comparison \eqref{eq:W-pert-gain} gives, after decreasing $\nu_0$ and then the small-clock
threshold if necessary,
\begin{equation}
|\mathcal W_{\cusp}^{\nu}(t)-\mathcal W_{\cusp}^{*}(t)| \le \varepsilon_W |\mathcal W_{\cusp}^{*}(t)|,
\label{eq:pert-cusp-strain-comparison}
\end{equation}
and we choose $\nu_0$ so that the constant in \eqref{eq:pressure-polarization-pert-bound}, after this strain
comparison, is absorbed by $q_{\rm tr}^{\rm pert}-q_{\rm tr}^*$.  Combining
\eqref{eq:target-cusp-pressure-win-recalled}--\eqref{eq:pert-cusp-strain-comparison} gives
\begin{equation}
\Pi_{\cusp}^{\nu}(t)\ge -q_{\rm tr}^{\rm pert}\,\tfrac12\bigl(\mathcal W_{\cusp}^{\nu}(t)\bigr)^2 .
\label{eq:pert-cusp-pressure-win}
\end{equation}

\begin{lemma}[Perturbed Riccati pressure lower bound and axial strain scale]
\label{lem:perturbed-pressure-hessian-slack}
There are constants and a threshold
\[
0<C_1\le C_2<\infty, \qquad \mathfrak J_{\Pi}^{\rm pert}\in(0,1],
\]
after decreasing $\nu_0$ if necessary, such that every admissibly perturbed solution satisfies, whenever
$J_{\cusp}^{\nu}(t)\le\mathfrak J_{\Pi}^{\rm pert}$,
\begin{equation}
\Pi_0^{\nu}(t)\ge -\upbeta\,\tfrac12\bigl(\rW_0^{\nu}(t)\bigr)^2
\label{eq:pert-euler-pressure-win}
\end{equation}
and, writing $J_\nu(t)$ for the physical meridional Jacobian at the stagnation point,
\begin{equation}
-C_2\Gamma J_\nu(t)^{3\alpha-1}\le \rW_0^{\nu}(t)\le -C_1\Gamma J_\nu(t)^{3\alpha-1}.
\label{eq:pert-strain-scale}
\end{equation}
The constants and thresholds depend only on $\alpha,\gamma,\eta$ and on the fixed pressure localization and
cutoff parameters.
\end{lemma}

\begin{proof}[Proof of Lemma~\ref{lem:perturbed-pressure-hessian-slack}]
Write $J:=J_{\cusp}^{\nu}(t)$ and
\[
\rW_{\cusp}^{\nu}(t):=m^{\nu}(t)\mathcal W_{\cusp}^{\nu}(t).
\]
We choose $\mathfrak J_{\Pi}^{\rm pert}$ no larger than the Target Profile pressure thresholds in
Lemma~\ref{lem:tail-error-final} and Proposition~\ref{prop:pressure-compare-modulated}; below we decrease it
only finitely many more times.  The perturbed pressure decomposition has the same form as
\eqref{eq:pressure-split-modulated}:
\begin{equation}
\Pi_0^{\nu}(t)=(m^{\nu}(t))^2\Pi_{\cusp}^{\nu}(t)+\Pi_{\rm rem}^{\nu}(t),
\label{eq:pert-pressure-split}
\end{equation}
where $\Pi_{\rm rem}^{\nu}$ is the sum of the smooth-flow deformation term, the mixed pressure Hessian, the
smooth pressure Hessian, and the cusp-error pressure Hessian.  By \eqref{eq:pert-cusp-pressure-win},
\begin{equation}
(m^{\nu}(t))^2\Pi_{\cusp}^{\nu}(t)\ge -q_{\rm tr}^{\rm pert}\,\tfrac12(\rW_{\cusp}^{\nu}(t))^2 .
\label{eq:pert-modulated-cusp-pressure-lower}
\end{equation}

We next absorb $\Pi_{\rm rem}^{\nu}$ at the cusp-strain scale.  Lemma~\ref{lem:pert-small-clock-persistence}
gives the perturbed analogues of the smooth-flow pressure deformation estimate
\eqref{eq:smooth-pressure-defect-bound} and the lower-order pressure remainder estimate
\eqref{eq:pressure-remainder-bound}.  The perturbative angular contribution in the localized bounded-slope
cone carries the small multiplier $\nu J^{3\eta}$ from \eqref{eq:localized-vort-pert}; the complementary
regions are controlled by \eqref{eq:pressure-zeta-tail}, \eqref{eq:pressure-angular-tail},
\eqref{eq:normal-form-approximation-bound}, and \eqref{eq:tail-pressure-remainder-bound}.  Thus the ratio of
the perturbed pressure remainder to the square of the cusp strain has the positive Target Profile clock powers
from \eqref{eq:pressure-rem-rWcusp-ratio}, together with an $O(\nu)$ term.  Moreover, \eqref{eq:W-pert-gain},
\eqref{eq:Wcusp-scaling}, and the perturbed modulation bounds in \eqref{eq:modulation-bounds} give
\begin{equation}
c\Gamma J^{3\alpha-1}\le -\rW_{\cusp}^{\nu}(t)\le C\Gamma J^{3\alpha-1}.
\label{eq:pert-cusp-strain-scale-intermediate}
\end{equation}
We first use the Target Profile cutoff choices in \eqref{eq:fixed-tail-smooth-rem-smallness}, then decrease
$\nu_0$, and finally decrease $\mathfrak J_{\Pi}^{\rm pert}$ so that
\begin{equation}
|\Pi_{\rm rem}^{\nu}(t)|\le (q_{\rm rem}^{\rm pert}-q_{\rm tr}^{\rm pert})\,\tfrac12(\rW_{\cusp}^{\nu}(t))^2 .
\label{eq:pert-pressure-rem-absorbed}
\end{equation}
Combining \eqref{eq:pert-pressure-split}, \eqref{eq:pert-modulated-cusp-pressure-lower}, and
\eqref{eq:pert-pressure-rem-absorbed} yields
\begin{equation}
\Pi_0^{\nu}(t)\ge -q_{\rm rem}^{\rm pert}\,\tfrac12(\rW_{\cusp}^{\nu}(t))^2 .
\label{eq:pert-pressure-lower-cusp-strain}
\end{equation}

It remains to convert from the cusp strain to the exact axial strain.  The perturbed axial strain splitting is
\[
\rW_0^{\nu}(t)=\rW_{\smooth}^{\nu}(t)+\rW_{\cusp}^{\nu}(t),
\]
the analogue of \eqref{eq:on-axis-strain-splitting}.  The perturbed smooth-strain bound, namely the first
estimate in the perturbed analogue of \eqref{eq:Jsmooth-bdd}, gives
$|\rW_{\smooth}^{\nu}(t)|\le C\Gamma$.  Since $3\alpha-1<0$, \eqref{eq:pert-cusp-strain-scale-intermediate}
allows us to decrease $\mathfrak J_{\Pi}^{\rm pert}$ so that
\begin{equation}
|\rW_{\smooth}^{\nu}(t)|\le \varepsilon_W |\rW_{\cusp}^{\nu}(t)|.
\label{eq:pert-smooth-strain-epsW}
\end{equation}
The cusp strain is negative in the same regime by \eqref{eq:pert-cusp-strain-scale-intermediate}.  Hence
\eqref{eq:pert-smooth-strain-epsW} gives
\begin{equation}
|\rW_0^{\nu}(t)|\ge (1-\varepsilon_W)|\rW_{\cusp}^{\nu}(t)|, \qquad (\rW_{\cusp}^{\nu}(t))^2\le (1-\varepsilon_W)^{-2}(\rW_0^{\nu}(t))^2.
\label{eq:pert-cusp-to-axial-strain-square}
\end{equation}
Using \eqref{eq:pert-cusp-to-axial-strain-square} in \eqref{eq:pert-pressure-lower-cusp-strain}, and then
using \eqref{eq:pert-strain-comparison-eps-choice}, proves \eqref{eq:pert-euler-pressure-win}.

The same splitting, \eqref{eq:pert-cusp-strain-scale-intermediate}, and
\eqref{eq:pert-smooth-strain-epsW} give
\[
c\Gamma J^{3\alpha-1}\le -\rW_0^{\nu}(t)\le C\Gamma J^{3\alpha-1}.
\]
Finally, the perturbed smooth-clock decomposition gives
\[
J_\nu(t)=J_{\smooth}^{\nu}(t)J_{\cusp}^{\nu}(t), \qquad c\le J_{\smooth}^{\nu}(t)\le C,
\]
by the perturbed analogue of \eqref{eq:Jsmooth-bdd}.  Thus powers of $J_{\cusp}^{\nu}$ and $J_\nu$ are
interchangeable up to constants, and the strain estimate becomes \eqref{eq:pert-strain-scale}.
\end{proof}

\begin{lemma}[Two-sided $L^\infty$ vorticity bounds for admissible angular profiles]
\label{lem:pert-vorticity-envelope}
There are constants and a threshold
\[
0<c<C<\infty, \qquad \mathfrak J_{\omega}^{\rm pert}\in(0,1],
\]
after decreasing $\nu_0$ if necessary, such that every admissibly perturbed solution satisfies
\begin{equation}
c\Gamma J_\nu(t)^{3\alpha-1}\le \|\omega^\nu(\cdot,t)\|_{L^\infty(\R^3)} \le C\Gamma J_\nu(t)^{3\alpha-1}
\label{eq:pert-vorticity-envelope}
\end{equation}
whenever $J_{\cusp}^{\nu}(t)\le \mathfrak J_{\omega}^{\rm pert}$.  Here $J_\nu(t)$ is the physical
meridional Jacobian at the stagnation point for the perturbed solution.
\end{lemma}

\begin{proof}[Proof of Lemma~\ref{lem:pert-vorticity-envelope}]
Let $J_{\cusp}^{\nu}$ be the cusp clock in the perturbed smooth-cusp decomposition.  We choose
$\mathfrak J_{\omega}^{\rm pert}\le\mathfrak J_{\omega,+}$, where
$\mathfrak J_{\omega,+}$ is the threshold in Lemma~\ref{lem:global-upper-vorticity-envelope}.  The perturbed
analogue of \eqref{eq:Jsmooth-bdd} gives
\[
J_\nu(t)=J_{\smooth}^{\nu}(t)J_{\cusp}^{\nu}(t), \qquad c\le J_{\smooth}^{\nu}(t)\le C,
\]
so it is enough to prove \eqref{eq:pert-vorticity-envelope} with $J_{\cusp}^{\nu}$ in place of $J_\nu$.

For the upper bound, decompose the transported angular vorticity as in
\eqref{eq:pert-cusp-vorticity-decomposition}.  After applying the perturbed smooth map, we write the physical
angular vorticity as $\omega^{\nu}=\omega^{*,\nu}+\omega^{\pert,\nu}$.  The proof of
Lemma~\ref{lem:global-upper-vorticity-envelope} applies to the Target Profile angular part
$\Omega_{\cusp}^{*,\nu}$ because the perturbed solution satisfies the analogues of
\eqref{eq:JtwoD-core-upper}, \eqref{eq:Omega-cusp-local-representation}--\eqref{eq:Omega-reg-local-bounds},
and \eqref{eq:cusp-separated-radial-log-growth}.  Therefore
\[
\|\omega^{*,\nu}(\cdot,t)\|_{L^\infty}\le C\Gamma (J_{\cusp}^{\nu}(t))^{3\alpha-1}.
\]
The perturbative angular part is bounded on the localized bounded-slope cone by \eqref{eq:localized-vort-pert}.
The region where the axial cutoff $1-\vartheta_\sharp$ is active, the region where the angular cutoff
$1-\chi_{M_\pressure}$ is active, the image-map displacement controlled by
\eqref{eq:normal-form-approximation-bound}, and the algebraic tail are bounded by
\eqref{eq:g0-poloidal-gain}--\eqref{eq:g0-holder-gain} together with the Target Profile bounds used in
Lemma~\ref{lem:global-upper-vorticity-envelope}.  Hence
\[
\|\omega^{\pert,\nu}(\cdot,t)\|_{L^\infty}
\le C\nu\Gamma (J_{\cusp}^{\nu}(t))^{3\alpha-1}+C\Gamma
\le C\Gamma (J_{\cusp}^{\nu}(t))^{3\alpha-1},
\]
after decreasing $\mathfrak J_{\omega}^{\rm pert}$, since $3\alpha-1<0$.

For the lower bound, we use the same label choice as in the proof of
Lemma~\ref{lem:target-vorticity-envelope}, with $J_{\cusp}^{\nu}(t)$ replacing $J_{\cusp}(t)$:
\[
\rho(Y)\in[\tfrac12,1], \qquad \sigma(Y)=\kappa (J_{\cusp}^{\nu}(t))^3.
\]
Here $\kappa>0$ is chosen no larger than the constants $\sigma_{\cut}$ and $c_{\rm vort}$ used in
\eqref{eq:driver-amplification-lower}.  After decreasing $\mathfrak J_{\omega}^{\rm pert}$ if necessary,
this label belongs to the supported core where $\Upsilon(\sigma(Y))=1$ in \eqref{eq:Theta-star-def}.
Let $\phi^\nu$ and $J_{\twoD}^{\nu}$ denote the meridional flow map and meridional Jacobian of this perturbed
solution.  The perturbed analogue of \eqref{eq:driver-amplification-lower} gives
$J_{\twoD}^{\nu}(Y,t)^{-1}\ge c(J_{\cusp}^{\nu}(t))^{-1}$.  After decreasing $\nu_0$ so that
$\|h\|_{L^\infty}\le\tfrac12$, the initial datum satisfies
\[
|\omega_{\theta,0}^{\nu}(Y)|= \Gamma\tfrac{\rho(Y)^\alpha}{(1+\rho(Y)^2)^{\gamma/2}}\Theta^*(\sigma(Y))|1+h(\sigma(Y))|
\ge c\Gamma (J_{\cusp}^{\nu}(t))^{3\alpha}.
\]
The exact vorticity transport identity \eqref{eq:vort-identity-JtwoD} then gives
\[
\|\omega^\nu(\cdot,t)\|_{L^\infty}\ge |\omega_\theta^\nu(\phi^\nu(Y,t),t)| \ge c\Gamma (J_{\cusp}^{\nu}(t))^{3\alpha-1}.
\]
Converting back to $J_\nu$ proves the lower bound.
\end{proof}

\begin{lemma}[No perturbed breakdown before clock collapse]
\label{lem:pert-no-breakdown-before-clock-collapse}
Let $0<\nu\le\nu_0$ and let the admissibly perturbed solution be in the small-clock regime from
Lemma~\ref{lem:pert-small-clock-persistence}, with
$J_{\cusp}^{\nu}(t)\le\mathfrak J_{\omega}^{\rm pert}$ on the time interval $[t_a,t_b)$.  If, on this
interval, the physical clock satisfies
\[
J_\nu(t)\ge j_0>0,
\]
then the corresponding $C^{1,\alpha}$ Euler solution cannot break down at $t_b<\infty$.  The continuation
bound depends on $j_0$, $\Gamma$, and on the fixed parameters.  The pointwise vorticity bound is independent
of $t_b$, and the Beale--Kato--Majda integral below is finite on every finite interval $[t_a,t_b)$.
\end{lemma}

\begin{proof}[Proof of Lemma~\ref{lem:pert-no-breakdown-before-clock-collapse}]
The upper bound in \eqref{eq:pert-vorticity-envelope} gives
\[
\|\omega^\nu(\cdot,t)\|_{L^\infty}\le C\Gamma J_\nu(t)^{3\alpha-1} \le C(j_0)\Gamma \qquad(t\in[t_a,t_b)).
\]
Lemma~\ref{lem:pert-small-clock-persistence} gives the cone-local and tail estimates needed to keep the
velocity in the same $C^{1,\alpha}$ axisymmetric no-swirl class on compact subintervals.  Hence
\[
\int_{t_a}^{t_b}\|\omega^\nu(\cdot,t)\|_{L^\infty}\,\ud t\le C(j_0)\Gamma(t_b-t_a)<\infty .
\]
The $C^{1,\alpha}$ continuation criterion, in its Beale--Kato--Majda form, excludes a finite breakdown at
$t_b$ while $J_\nu$ remains bounded below by $j_0$.
\end{proof}

\subsection{Blowup for Euler solutions with admissible angular profiles}

\begin{proof}[Proof of Theorem~\ref{thm:main}]
Lemma~\ref{lem:L2} shows that $u_0\in C^{1,\alpha}(\R^3)\cap L^2(\R^3)$ for every admissible datum in
$\mathcal A_{\alpha,\gamma}(\nu,\eta)$.  The same standard local well-posedness theory in the
$C^{1,\alpha}$ axisymmetric no-swirl class, together with the Beale--Kato--Majda continuation criterion used
in the Target Profile proof, provides the unique local Euler solution.  Let $u$ be the Euler solution generated
by an admissible datum in $\mathcal A_{\alpha,\gamma}(\nu,\eta)$ with $0<\nu\le\nu_0$.  From this point on we
suppress the superscript $\nu$ on the flow map and its clock.  We write
\[
J(t):=\det\nabla_{(R,Z)}(\phi_r,\phi_z)(0,0,t), \qquad \rW_0(t):=\p_z u_z(0,0,t).
\]
For this solution, the clock $J_\nu$ in \eqref{eq:pert-strain-scale} is exactly $J$.  The perturbed analogue
of \eqref{eq:Jsmooth-bdd} gives
\[
J(t)=J_{\smooth}(t)J_{\cusp}(t), \qquad c\le J_{\smooth}(t)\le C,
\]
so collapse of the cusp clock and collapse of the physical clock are equivalent.

Let $\mathfrak J_{\rm st}\in(0,1]$ be no larger than the pressure threshold
$\mathfrak J_{\Pi}^{\rm pert}$ in Lemma~\ref{lem:perturbed-pressure-hessian-slack} and the vorticity
threshold $\mathfrak J_{\omega}^{\rm pert}$ in Lemma~\ref{lem:pert-vorticity-envelope}.  We also take
$\mathfrak J_{\rm st}$ below the small-clock threshold in Lemma~\ref{lem:pert-small-clock-persistence}.  By
Lemma~\ref{lem:pert-small-clock-persistence}, applied with
$\mathfrak J_{\mathrm{finite}}=\mathfrak J_{\rm st}$, the entry time
\[
t_{\rm ent}:=\inf\{\,t\ge0:\ J_{\cusp}(t)\le \mathfrak J_{\rm st}\,\}
\]
satisfies \eqref{eq:pert-late-entry-bound}; in particular $t_{\rm ent}\le C\Gamma^{-1}$.  For
$t\ge t_{\rm ent}$ the solution is small-clock stable in the sense of
Definition~\ref{def:small-clock-stable}, with constants independent of $\nu$.  The perturbed analogue of
\eqref{eq:Jdot-two-sided} gives $\dot J_{\cusp}(t)<0$ throughout this
small-cusp-clock regime, so the solution cannot leave the regime by increasing the cusp clock.
Lemma~\ref{lem:pert-no-breakdown-before-clock-collapse} rules out any breakdown of the $C^{1,\alpha}$
solution while the physical clock stays bounded below.

On any time interval after $t_{\rm ent}$ on which the solution exists and $J(t)>0$, the pressure estimate
\eqref{eq:pert-euler-pressure-win} and the stagnation-point Riccati law in \eqref{eq:stagnation-ODEs} give
\begin{equation}
\p_t\rW_0(t)=-\tfrac12\rW_0(t)^2-\Pi_0(t) \le -\tfrac{1-\upbeta}{2}\rW_0(t)^2.
\label{eq:pert-riccati-collapse}
\end{equation}
The axial strain estimate \eqref{eq:pert-strain-scale} gives
$\rW_0(t)\le -C_1\Gamma J(t)^{3\alpha-1}<0$ in the same regime, so
\eqref{eq:pert-riccati-collapse} applies to a negative axial strain.  The clock identity in
\eqref{eq:stagnation-ODEs},
\[
\dot J(t)=\tfrac12J(t)\rW_0(t)
\]
and \eqref{eq:pert-strain-scale} give
\begin{equation}
-C_2\Gamma J(t)^{3\alpha}\le \dot J(t)\le -C_1\Gamma J(t)^{3\alpha}.
\label{eq:pert-clock-ode}
\end{equation}
Differentiating $J(t)^{1-3\alpha}$ and using \eqref{eq:pert-clock-ode} gives
\begin{equation}
c\Gamma\le -\tfrac{d}{dt}J(t)^{1-3\alpha}\le C\Gamma .
\label{eq:pert-clock-power-ode}
\end{equation}

We now identify the endpoint.  If the $C^{1,\alpha}$ solution had a finite maximal time before clock collapse,
then $J(t)$ would stay bounded below by some $j_0>0$ on a terminal interval, and
Lemma~\ref{lem:pert-no-breakdown-before-clock-collapse} would extend the solution beyond that time.  If the
solution existed for all later times while $J(t)>0$, then integrating the left inequality in
\eqref{eq:pert-clock-power-ode} from $t_{\rm ent}$ would force $J(t)^{1-3\alpha}$ to become negative after a
finite time.  Hence there is a first finite endpoint $T^*<\infty$ such that
\[
J(t)\downarrow0 \qquad\text{as }t\uparrow T^*.
\]
The bounded smooth clock in \eqref{eq:Jsmooth-bdd} then implies $J_{\cusp}(t)\downarrow0$ as well.

For $t<s<T^*$, integrating \eqref{eq:pert-clock-power-ode} over $[t,s]$ gives
\[
c\Gamma(s-t)\le J(t)^{1-3\alpha}-J(s)^{1-3\alpha}\le C\Gamma(s-t).
\]
Letting $s\uparrow T^*$ and using $J(s)\to0$ gives
\begin{equation}
c\Gamma(T^*-t)\le J(t)^{1-3\alpha}\le C\Gamma(T^*-t).
\label{eq:pert-clock-power-typeI}
\end{equation}
Since $1-3\alpha>0$, applying the increasing map $x\mapsto x^{1/(1-3\alpha)}$ to
\eqref{eq:pert-clock-power-typeI} and renaming constants gives
\begin{equation}
c\bigl(\Gamma(T^*-t)\bigr)^{\tfrac{1}{1-3\alpha}} \le J(t)\le C\bigl(\Gamma(T^*-t)\bigr)^{\tfrac{1}{1-3\alpha}}.
\label{eq:pert-clock-typeI}
\end{equation}
Since $3\alpha-1=-(1-3\alpha)$, \eqref{eq:pert-clock-power-typeI} also gives
\begin{equation}
c(T^*-t)^{-1}\le \Gamma J(t)^{3\alpha-1}\le C(T^*-t)^{-1}.
\label{eq:pert-singular-scale-typeI}
\end{equation}
Combining \eqref{eq:pert-strain-scale} with \eqref{eq:pert-singular-scale-typeI}, we obtain
\begin{equation}
\tfrac{c}{T^*-t}\le -\p_z u_z(0,0,t)\le \tfrac{C}{T^*-t}.
\label{eq:pert-strain-typeI}
\end{equation}
Combining the two-sided $L^\infty$ vorticity bound \eqref{eq:pert-vorticity-envelope} with
\eqref{eq:pert-singular-scale-typeI} gives
\begin{equation}
\tfrac{c}{T^*-t}\le \|\omega(\cdot,t)\|_{L^\infty(\R^3)}\le \tfrac{C}{T^*-t}.
\label{eq:pert-vorticity-typeI}
\end{equation}
The upper bound in \eqref{eq:pert-vorticity-typeI} gives $C^{1,\alpha}$ continuation on every compact
subinterval of $[0,T^*)$, and Lemma~\ref{lem:pert-no-breakdown-before-clock-collapse} excludes an earlier
endpoint with the clock bounded below.  The lower bound in \eqref{eq:pert-vorticity-typeI} forces
$\|\omega(\cdot,t)\|_{L^\infty}$ to diverge as $t\uparrow T^*$.  Hence $T^*$ is the maximal
$C^{1,\alpha}$ existence time.  Finally, for any $t_{\rm I}<T^*$ on which
\eqref{eq:pert-vorticity-typeI} holds,
\[
\int_0^{T^*}\|\omega(\cdot,t)\|_{L^\infty}\,\ud t \ge c\int_{t_{\rm I}}^{T^*}(T^*-t)^{-1}\,\ud t=\infty .
\]
Together with \eqref{eq:pert-clock-typeI}, \eqref{eq:pert-strain-typeI}, and
\eqref{eq:pert-vorticity-typeI}, this proves Theorem~\ref{thm:main}.
\end{proof}

\section*{Acknowledgements} This work was supported by NSF grant DMS-2307680.

\end{document}